\documentclass[11pt,reqno]{amsart}

\usepackage[leqno]{amsmath}
\usepackage{amsthm}
\usepackage{amsfonts, tikz-cd}
\usepackage{amssymb}
\usepackage{eucal}
\usepackage[all]{xy}
\usepackage{mathtools}
\usepackage{stmaryrd}
\usepackage{csquotes}
\usepackage{enumitem}
\usepackage{upgreek}

\setenumerate{itemsep=2pt,parsep=2pt,before={\parskip=2pt}}

\usepackage[colorlinks=true,hyperindex, linkcolor=magenta, pagebackref=false, citecolor=cyan,pdfpagelabels]{hyperref}

\usepackage{xspace}
\usepackage{epsfig,epic,eepic,latexsym,color}
\usepackage[all]{xy}
\usepackage{comment} 
\numberwithin{equation}{section}

\usepackage[english]{babel}
\usepackage{latexsym}

\newcommand{\nc}{\newcommand}
\nc{\rnc}{\renewcommand}
\rnc{\P}{\mathbf P}
\nc{\R}{\mathbf R}
\rnc{\rm}{\mathrm}
\rnc{\bf}{\mathbf}
\nc{\cal}{\mathcal}

\nc{\C}{\mathbf C}
\nc{\Q}{\mathbf Q}
\nc{\Z}{\mathbf Z}
\nc{\N}{\mathbf N}
\nc{\A}{\mathbf A}

\nc{\an}{\operatorname{an}}
\nc{\perfd}{\operatorname{perfd}}
\nc{\perf}{\operatorname{new, perf}}
\nc{\diam}{\diamondsuit}

\nc{\htt}{\operatorname{ht}}
\nc{\Nm}{\operatorname{Nm}}
\nc{\Ker}{\operatorname{Ker}}
\nc{\mmod}{\operatorname{mod}}
\nc{\End}{\operatorname{End}}
\nc{\Tor}{\operatorname{Tor}}
\nc{\coker}{\operatorname{Coker}}

\nc{\Aut}{\operatorname{Aut}}
\nc{\cont}{\text{cont}}
\nc{\sep}{\text{sep}}
\nc{\Hom}{\mathrm{Hom}}
\nc{\Gal}{\mathrm{Gal}}
\nc{\Spec}{\text{Spec}\,}
\nc{\Spd}{\text{Spd}\,}
\nc{\RZ}{\operatorname{RZ}}
\nc{\hocolim}{\operatorname{hocolim}}

\rnc{\t}{\tau}
\nc{\mm}{\pmb{\mu}}
\rnc{\a}{\alpha}
\nc{\n}{\mathfrak n}
\nc{\m}{\mathfrak m}
\nc{\mfs}{\mathfrak s}
\nc{\mf}{\mathfrak f}
\nc{\e}{\varepsilon}
\nc{\dd}{\delta}
\nc{\Imm}{\operatorname{Im}}
\rnc{\sp}{\operatorname{sp}}

\nc{\p}{\mathfrak p}
\nc{\q}{\mathfrak q}

\nc{\Sym}{\operatorname{Sym}}
\nc{\codim}{\operatorname{codim}}
\nc{\rk}{\operatorname{rk}}
\nc{\GL}{\operatorname{GL}}
\nc{\SL}{\operatorname{SL}}
\nc{\Lie}{\operatorname{Lie}}
\nc{\Ind}{\operatorname{Ind}}
\nc{\Div}{\underline{Div}}
\nc{\Pic}{\mathbf{Pic}}
\nc{\uPic}{\underline{ \mathbf{Pic}}}
\nc{\rH}{\mathrm{H}}
\nc{\Spf}{\operatorname{Spf}\,}
\nc{\Frac}{\operatorname{Frac}}
\nc{\colim}{\operatorname{colim}}
\nc{\Spa}{\operatorname{Spa}\,}

\rnc{\an}{\operatorname{an}}
\nc{\et}{\text{\'et}}
\nc{\set}{\text{s\'et}}
\nc{\fet}{\text{f\'et}}
\nc{\sfet}{\text{sf\'et}}

\newcommand{\adj}[4]{#1\negmedspace: #2\rightleftarrows #3:\negmedspace #4}

\rnc{\et}{\text{\'et}}
\nc{\etaff}{\text{\'etaff}}
\nc{\proet}{\text{pro\'et}}
\nc{\qproet}{\text{qpro\'et}}
\nc{\qp}{\rm{qp}}
\nc{\coh}{\text{coh}}
\nc{\acoh}{\text{acoh}}
\nc{\aft}{\text{aft}}
\nc{\afp}{\text{afp}}
\nc{\aqc}{\text{aqc}}
\nc{\qc}{\text{qc}}

\nc{\xr}{\xrightarrow}

\nc{\eps}{\epsilon}
\nc{\ov}{\overline}
\nc{\ud}{\underline}
\nc{\wdh}{\widehat}

\nc{\F}{\mathcal F}
\nc{\G}{\mathcal G}
\nc{\E}{\mathcal E}
\nc{\K}{\mathcal K}
\nc{\I}{\mathcal I}
\nc{\sQ}{\mathcal Q}

\nc{\X}{\mathfrak X}
\nc{\Y}{\mathfrak Y}
\nc{\T}{\mathfrak T}

\nc{\LL}{\mathcal{L}}
\rnc{\S}{\mathcal S}
\nc{\M}{\mathcal M}
\nc{\sU}{\mathfrak U}
\nc{\V}{\mathfrak V}

\nc{\ra}{\rangle}

\nc{\os}{\overset}

\rnc{\O}{\mathcal O}
\nc{\J}{\mathcal J}
\theoremstyle{definition}

\newtheorem{thm1}{Theorem}[section]
\newtheorem{lemma1}[thm1]{Lemma}
\newtheorem{defn1}[thm1]{Definition}

\newtheorem{rmk1}[thm1]{Remark}

\newtheorem{cor1}[thm1]{Corollary}
\newtheorem{warning1}[thm1]{Warning}

\newtheorem{thm}{Theorem}[subsection]
\newtheorem{lemma}[thm]{Lemma}
\newtheorem{defn}[thm]{Definition}

\newtheorem{prop}[thm]{Proposition}
\newtheorem{Cor}[thm]{Corollary}
\newtheorem{rmk}[thm]{Remark}
\newtheorem{warning}[thm]{Warning}

\newtheorem{question}[thm]{Question}
\newtheorem{example}[thm]{Example}

\newtheorem{setup}[thm]{Set-up}
\newtheorem{cor}[thm]{Corollary}

\begin{document}

\title{Almost Coherent Modules and Almost Coherent Sheaves}
\bibliographystyle{halpha-abbrv}
\author{Bogdan Zavyalov}
\maketitle

\begin{abstract}
We extend the theory of almost coherent modules that was introduced in ``Almost Ring Theory'' \cite{GR} by Gabber and Ramero.
Then we globalize it by developing a new theory of almost coherent sheaves on schemes and on a class of ``nice'' formal
schemes. We show that these sheaves satisfy many properties similar to usual coherent sheaves, i.e., the Almost Proper 
Mapping Theorem, the Formal GAGA, etc. We also construct an almost version of the Grothendieck twisted image functor $f^!$ 
and verify its properties. Lastly, we study sheaves of $p$-adic nearby cycles on admissible formal models of rigid-analytic varieties 
and show that these sheaves provide examples of almost coherent sheaves. This gives a new proof of the finiteness result for 
\'etale cohomology of proper rigid-analytic varieties obtained before in the work of Peter Scholze ``$p$-adic Hodge Theory For 
Rigid-Analytic Varieties'' \cite{Sch1}.
\end{abstract}

\tableofcontents
\newpage

\section{Introduction}
\subsection{Motivation}
The purpose of this paper is two-fold. The first goal is to develop a sufficiently rich theory of almost coherent sheaves on schemes and a class of formal schemes. The second goal is to provide the reader with one interesting source of examples of almost coherent sheaves. Namely, we show that the complex of $p$-adic nearby cycles $\bf{R}\nu_* \left(\cal{E}\right)$ has quasi-coherent, almost coherent cohomology sheaves for any admissible formal $\O_C$-scheme $\X$ and $\O_{X^\diam}^+/p$-vector bundle $\cal{E}$ (see Definition~\ref{defn:vect-bundles}). \smallskip

Before we discuss the content of each chapter in detail, we explain the motivation behind the work done in this manuscript. \smallskip 

The first source motivation comes from the work of P.\,Scholze on the finiteness of $\bf{F}_p$-cohomology groups of proper rigid-analytic varieties over $p$-adic fields (see \cite{Sch1}). The second source of motivation (clearly related to the first one) is the desire to set up a robust enough theory of almost coherent sheaves that is crucially used in our proof of Poincar\'e Duality for $\bf{F}_p$-local systems on smooth and proper rigid-analytic varieties over $p$-adic fields in \cite{Zav4}. \smallskip

We start with the work of P.\,Scholze. In \cite{Sch1}, he showed that there is an almost isomorphism \[
\rm{H}^i(X, \bf{F}_p) \otimes \O_C/p \simeq^a \rm{H}^i(X, \O_{X_\et}^+/p)
\]
for any proper rigid-analytic variety $X$ over a $p$-adic algebraically closed field $C$. This almost isomorphism allows us to reduce studying certain properties of  $\rm{H}^i(X, \bf{F}_p)$ for a $p$-adic proper rigid-analytic space $X$ to studying almost properties of the cohomology groups $\rm{H}^i(X, \O_{X_\et}^+/p)$, or the full complex $\bf{R}\Gamma(X, \O_{X_\et}^+/p)$. For instance, Scholze shows that $\rm{H}^i(X, \bf{F}_p)$ are finite groups by deducing it from almost coherence of $\rm{H}^i(X, \O_{X_\et}^+/p)$ over $\O_C/p$. \smallskip

Scholze's argument does not involve any choice of an admissible formal model for $X$ and is performed entirely on the generic fiber via an elaborate study of cancellations in certain spectral sequences. A different natural approach to studying $\bf{R}\Gamma(X, \O_{X_\et}^+/p)$ is to rewrite this complex as 
\[
\bf{R}\Gamma\left(X, \O_{X_\et}^+/p\right) \simeq \bf{R}\Gamma\left(\X_0, \bf{R}t_*\O_{X_\et}^+/p\right)
\]
for a choice of an admissible formal $\O_C$-model $\X$ and the natural morphism of ringed sites 
\[
t \colon (X_\et, \O_{X_\et}^+/p) \to (\X_{0, \rm{Zar}}, \O_{\X_0})
\]
with $\X_0$ the mod-$p$ fiber of $\X$. Then we can separately study the complex $\bf{R}t_*\left(\O_{X_\et}^+/p\right)$ and the functor $\bf{R}\Gamma(\X,-)$. In order to make this strategy work, we develop the notion of almost coherent sheaves on $\X$ and $\X_0$ and show its various properties similar to the properties of coherent sheaves. This occupies Chapters~\ref{almost-commutative-algebra}-\ref{section:cohomological-properties}. While Chapters~\ref{section:proetale-v-sites}~and~\ref{section:almost-coherent-nearby-cycles} are devoted to showing that the complex $\bf{R}t_*\left(\O_{X_\et}^+/p\right)$ has almost coherent cohomology groups, and to generalizing these finiteness results to all $\O^+/p$-vector bundles. Combining it with the Almost Proper Mapping Theorem~\ref{thm:intro-almost-proper-schemes}, we reprove \cite[Lemma 5.8 and Theorem 5.1]{Sch1} in a slightly greater generality (allowing arbitrary Zariski-constructible coefficients as opposed to local systems).

\begin{thm}\label{thm:intro-1}(Lemma~\ref{almost-scholze},  Lemma~\ref{lemma:relative-primitive}, and Lemma~\ref{lemma:et-qp-v-overconv-coeff})  Let $C$ be an algebraically closed $p$-adic non-archimedean field, let $X$ be a proper rigid-analytic variety over $C$, and let $\F$ be a Zariski-constructible sheaf of $\bf{F}_p$-modules (see Definition~\ref{defn:zc}). Then 
\begin{enumerate}[label=\textbf{(\arabic*)}]
    \item $\rm{H}^i(X, \F \otimes_{\bf{F}_p}\O_{X_\et}^+/p)$ is an almost finitely generated $\O_C/p$-module for $i\geq 0$;
    \item the natural morphism
    \[
    \rm{H}^i\left(X, \F\right) \otimes_{\bf{F}_p} \O_C/p \to \rm{H}^i\left(X, \F\otimes_{\bf{F}_p} \O_{X_\et}^+/p\right)
    \]
    is an almost isomorphism for $i\geq 0$;
    \item $\rm{H}^i(X, \F \otimes_{\bf{F}_p}\O_{X_\et}^+/p)$ is almost zero for $i>2\dim X$.
\end{enumerate}
\end{thm}

\begin{thm}\label{thm:intro-2}(Lemma~\ref{lemma:primitive-zc})\footnote{Theorem~\ref{thm:intro-2} can also be easily deduced from the results of \cite{Bhatt-Hansen}.}  In the notation of Theorem~\ref{thm:intro-1}, we have
\begin{enumerate}[label=\textbf{(\arabic*)}]
    \item $\rm{H}^i(X, \F)$ is a finite group for $i\geq 0$;
    \item $\rm{H}^i(X, \F)\simeq 0$ for $i > 2\dim X$.
\end{enumerate}
\end{thm}

Now we discuss the role of this paper in our proof of Poincar\'e Duality in \cite{Zav4}. We start with the precise formulation of this result. 

\begin{thm}\label{thm:later}\cite{Zav4} Let $C$ be an algebraically closed $p$-adic non-archimedean field, let $X$ be a rigid-analytic variety over $C$ of pure dimension $d$, and let $\bf{L}$ be an $\bf{F}_p$-local system on $X_\et$. Then there is a canonical trace map 
\[
t_X\colon\rm{H}^{2d}\left(X, \bf{F}_p(d)\right) \to \bf{F}_p
\]
such that the induced pairing
\[
\rm{H}^i\left(X, \bf{L}\right)\otimes \rm{H}^{2d-i}(X, \bf{L}^{\vee}(d))\xr{-\cup-} \rm{H}^{2d}(X, \bf{F}_p(d)) \xr{t_X} \bf{F}_p
\]
is perfect. 
\end{thm}

The essential idea of the proof (at least for $\bf{L}=\ud{\bf{F}}_p$) is to use Theorem~\ref{thm:intro-1} to reduce Poincar\'e Duality to almost duality for the complex $\bf{R}\Gamma(X, \O_{X_\et}^+/p)$. We study this complex via the isomorphism
\[
\bf{R}\Gamma(X, \O_{X_\et}^+/p) \simeq \bf{R}\Gamma(\X_0, \bf{R}t_*\O_{X_\et}^+/p).
\]
Roughly, we separately show almost duality for the ``nearby cycles functor'' $\bf{R}t_*$, and then establish an almost version of Grothendieck Duality for the $\O_C/p$-scheme $\X_0$. Even to formulate these things precisely, one needs a good theory of almost (coherent) sheaves that globalizes the theory of almost (coherent) modules. For this theory to be useful, we have to establish that almost coherent sheaves share many properties similar to classical coherent sheaves {\it and} the ``nearby cycles''  $\bf{R}t_*\left(\O_{X_\et}^+/p\right)$ are almost coherent. \smallskip 

The main goal of Sections~\ref{almost-commutative-algebra}-\ref{section:cohomological-properties} is to develop this general theory of almost (coherent) sheaves. In Section~\ref{section:proetale-v-sites}, we study $\O^+/p$, $\O^+$, and $\O$-vector bundles in different topologies. Section~\ref{section:almost-coherent-nearby-cycles} is devoted to verifying that ``nearby cycles'' are almost coherent. That being said, we now discuss the content and main results of each section in more detail. 

\subsection{Foundations of almost mathematics (Sections~\ref{almost-commutative-algebra}-\ref{section:cohomological-properties})}
Section~\ref{almost-mathematics} defines the category of almost modules and studies its main properties. This section is very motivated by \cite{GR}. However, it seems that some results that we need later in the paper are not present in \cite{GR}, so we give an (almost) self-contained introduction to almost commutative algebra. We define the category of almost modules (see the discussion after Corollary~\ref{cor:sigma-serre-subcategory}), the almost tensor product functor $-\otimes_{R^a} -$ (see Proposition~\ref{many-functors}(1)), the almost Hom functor $\rm{alHom}_{R^a}(-, -)$ (see Proposition~\ref{many-functors}(3)), and the notion of almost finitely generated (see Definitions~\ref{defn:almost-finitely-generated}), almost finitely presented (see Definition~\ref{defn:almost-finitely-presented}), and almost coherent modules (see Definition~\ref{acoh}). We show that almost coherent modules satisfy many natural properties similar to the properties of classical coherent modules. We summarize some of them in the following theorem:

\begin{thm}\label{thm:intro-section-2}(Lemma~\ref{main-coh}, Propositions~\ref{tensor-product-coh}, \ref{alHom-derived-coh}, \ref{base-change-hom-2var-derived}, Theorem~\ref{thm:faithuflly-flat descent}, and Lemma~\ref{almost-flat-descent}) Let $R$ be a ring with an ideal $\m$ such that $\widetilde{\m}\coloneqq \m\otimes_R \m$ is $R$-flat and $\m^2=\m$. 
\begin{enumerate}[label=\textbf{(\arabic*)}]
    \item Almost coherent $R^a$-modules form a weak Serre subcategory of $\bf{Mod}_R^a$; 
    \item If $R$ is an almost coherent ring (i.e. free rank-$1$ $R$-module is almost coherent), and $M^a, N^a$ two objects in $\bf{D}^-_{acoh}(R)^a$. Then $M^a\otimes_{R^a}^L N^a\in \bf{D}^-_{acoh}(R)^a$;
    \item If $R$ is an almost coherent ring, and $M^a \in \bf{D}^-_{acoh}(R)^a$, $N^a \in \bf{D}^+_{acoh}(R)^a$. Then \[
    \bf{R}\rm{alHom}_{R^a}(M^a, N^a)\in \bf{D}^+_{acoh}(R)^a;
    \]
    \item If $R$ is an almost coherent ring, $M^a\in \bf{D}_{acoh}^-(R)^a$, $N^a\in \bf{D}^+(R)^a$, and $P^a$ an almost flat $R^a$-module. Then the natural map $\bf{R}\rm{Hom}_{R^a}(M^a, N^a) \otimes_{R^a} P^a \to \bf{R}\rm{Hom}_{R^a}(M^a, N^a\otimes_{R^a} P^a)$ is an almost isomorphism;
    \item Descent of almost modules along an almost faithfully flat morphism $R\to S$ is always effective;
    \item Let $R \to S$ be an almost faithfully flat map, and let $M^a$ be an $R^a$-module. Suppose that $M^a\otimes_{R^a} S^a$ is almost finitely generated (resp. almost finitely presented, resp. almost coherent) $S^a$-module. Then so is $M^a$.
\end{enumerate}
\end{thm}

If $R$ is $I$-adically adhesive for some finitely generated ideal $I$ (see Definition~\ref{defn:adhesive}), we can show that almost finitely generated $R$-modules satisfy a (weak) version of the Artin--Rees Lemma, and behave nicely with respect to the completion functor. These results will be crucial for globalizing the theory of almost coherent modules on formal schemes. 

\begin{lemma}(Lemma~\ref{Artin-Rees} and Lemma~\ref{completion-finitely-generated}) Let $R$ be an $I$-adically adhesive ring with an ideal $\m$ such that $I\subset \m$, $\m^2=\m$, and $\m\otimes_R \m$ is $R$-flat (see Set-up~\ref{set-up1}). Let $M$ be an almost finitely generated $R$-module. Then
\begin{enumerate}[label=\textbf{(\arabic*)}]
    \item for any $R$-submodule $N\subset M$, the induced topology on $N$ coincides with the $I$-adic topology;
    \item The natural morphism $M\otimes_R \wdh{R} \to \wdh{M}$ is an isomorphism. In particular, if $R$ is $I$-adically complete, then any almost finitely generated $R$-module is also $I$-adically complete.
\end{enumerate}
\end{lemma}

If $R$ is a (topologically) finitely generated algebra over a perfectoid valuation ring $K^+$ (see Definition~\ref{defn:valuation-perfectoid}), we can say even more. In this case, it turns out that $R$ is almost noetherian (see Definition~\ref{defn:almost-noetherian}), so the theory simplifies significantly. Another useful result is that it suffices to check that a derived complete complex is almost coherent after taking the derived quotient by a pseudo-uniformizer $\varpi$. This is very handy in practice because it reduces many (subtle) integral questions to the torsion case, where there are no topological subtleties. 

\begin{thm}\label{thm:intro-mod-p}(Theorem~\ref{thm:top-ft-almost-noetherian}, Theorem~\ref{thm:finite-type-almost-noetherian}, Theorem~\ref{thm:check-almost-coh-mod-p}) Let $K^+$ be a perfectoid valuation ring with a pseudo-uniformizer $\varpi$ as in Lemma~\ref{lemma:roots-of-pseudounformizer}, and let $R$ be a $K^+$-algebra. Then
\begin{enumerate}[label=\textbf{(\arabic*)}]
    \item $R$ is almost noetherian if $R$ is (topologically) finite type over $K^+$;
    \item Suppose $R$ is a topologically finite type $K^+$-algebra and $M$ is a derived $\varpi$-adically complete object in $\bf{D}(R)$ such that $[M/\varpi]\in \bf{D}^{[c, d]}_{acoh}(R/\varpi)$. Then $M\in \bf{D}^{[c,d]}_{acoh}(R)$. 
\end{enumerate}
\end{thm}

\smallskip

We discuss the extension of almost mathematics to ringed sites in Section~$3$. The main goal is to generalize all constructions from almost mathematics to a general ringed site. We define the notion of almost $\O_X$-modules on a ringed site $(X, \O_X)$ (see Definition~\ref{defn:almost-sheaves-1}) and of $\O_X^a$-modules (see Definition~\ref{defn:almost-sheaves-2}) and show that they are equivalent:

\begin{thm}\label{thm:intro-section-equivalent-almost sheaves} (Theorem~\ref{almost-two-different-variants}) Let $R$ be as in Theorem~\ref{thm:intro-section-2} and $(X, \O_X)$ a ringed $R$-site. Then the functor
\[
(-)^a\colon \bf{Mod}_{\O_X}^a \to \bf{Mod}_{\O_X^a}
\]
is an equivalence of categories.
\end{thm}

We also define the functors $- \otimes -$, $\rm{Hom}_{\O_X^a}(-, -)$, $\rm{alHom}_{\O_X^a}(-, -)$, $\ud{\cal{H}om}_{\O_X^a}(-, -)$, $\ud{al\cal{H}om}_{\O_X^a}(-, -)$, $f_*$, and $f^*$ on the category of $\O_X^a$-modules. We refer to Section~\ref{basic-functors-sheaves} for an extensive discussion of these functors. Then we study the derived category of $\O_X^a$-modules and derived analogues of the functors mentioned above. This is done in Sections~\ref{section:derived-sheaves} and~\ref{section:derived-sheaves-functors}. \smallskip{}

We develop the theory of almost finitely presented and almost (quasi-)coherent sheaves on schemes and on a class of formal schemes in Section~\ref{acoh-sheaves}. The main goal is to show that these sheaves behave similarly to classical coherent sheaves in many aspects. \smallskip{}

We roughly define almost finitely presented $\O_X^a$-modules as modules such that, for any finitely generated sub-ideal $\m_0\subset \m$, can be locally approximated by finitely presented $\O_X$-modules up to modules annihilated by $\m_0$ (see Definition~\ref{defn:almost-finitely-resented-schemes} for a precise definition). Sections~\ref{acoh-sheaves}-\ref{section:derived-category-schemes} are mostly concerned with local properties of these sheaves. We summarize some of the main results below: 

\begin{thm}\label{intro:thm-local-properties-schemes}(Corollary~\ref{main-coh-glob}, Theorem~\ref{derived-schemes-2}, Lemmas~\ref{tensor-almost-coh-derived}, \ref{pullback-almost-coh-derived}, \ref{pushforward-almost-coh-derived}, and \ref{hom-alcoh-derived}) Let $R$ be a ring with an ideal $\m$ such that $\widetilde{\m}\coloneqq \m\otimes_R \m$ is $R$-flat and $\m^2=\m$. 
\begin{enumerate}[label=\textbf{(\arabic*)}]
    \item For any $R$-scheme $X$, almost coherent $\O_X^a$-modules form a weak Serre subcategory of $\bf{Mod}_{\O_X}^a$;
    
    \item The functor 
    \[
        \widetilde{(-)} \colon \bf{D}_*(R)^a \to \bf{D}_{aqc, *}(\Spec R)^a
    \] 
    is a $t$-exact equivalence of triangulated categories for $* \in \{``\ ", \rm{acoh} \}$. Its quasi-inverse is given by $\bf{R}\Gamma(\Spec R, -)$. In particular, an almost quasi-coherent $\O_{\Spec R}^a$-module $\F^a$ is almost coherent if and only if $\F^a(\Spec R)$ is an almost coherent $R^a$-module;
    \item The natural morphism $\widetilde{M^a\otimes^L_{R^a} N^a} \to \widetilde{M^a}\otimes^L_{\O_{\Spec R}^a}\widetilde{N^a}$
    is an isomorphism for any $M^a, N^a\in \bf{D}(R)^a$;
    \item Let $f\colon \Spec B \to \Spec A$ be an $R$-morphism of affine schemes. Then $\bf{L}f^*(\widetilde{M^a})$ is functorially isomorphic to $\widetilde{M^a\otimes^L_{A^a} B^a}$ for any $M^a\in \bf{D}(A)^a$;
    \item Let $f\colon X\to Y$ be a quasi-compact and quasi-separated morphism of $R$-schemes. Suppose that $Y$ is quasi-compact. Then $\bf{R}f_*$ carries $\bf{D}^*_{aqc}(X)^a$ to $\bf{D}^*_{aqc}(Y)^a$ for any $*\in \{``\text{ ''}, -, +, b\}$;
    \item Suppose that $R$ is almost coherent. Then the natural maps
    \[
        \widetilde{\bf{R}\rm{alHom}_{R^a}(M^a,N^a)} \to \bf{R}\ud{al\cal{H}om}_{\O_{\Spec R}^a}(\widetilde{M^a}, \widetilde{N^a}),
    \]
    \[
        \widetilde{\bf{R}\rm{Hom}_{R^a}(M^a,N^a)} \to \bf{R}\ud{\cal{H}om}_{\O_{\Spec R}^a}(\widetilde{M^a}, \widetilde{N^a})
    \]
are almost isomorphisms for $M^a\in \bf{D}^-_{acoh}(R)^a$, $N^a\in \bf{D}^+(R)^a$.
\end{enumerate}
\end{thm}

We also show that, for a quasi-compact and quasi-separated scheme $X$, any almost finitely presented $\O_X^a$-module admits a global approximation by finitely presented $\O_X$-modules. This result is crucial for establishing {\it global} properties of almost finitely presented $\O_X^a$-modules, and it will be systematically used in Chapter~\ref{section:cohomological-properties}. \smallskip

\begin{thm}\label{thm:intro-global-approximation} (Corollary~\ref{cor:approximate-afpr}) Let $X$ be a quasi-compact and quasi-separated $R$-scheme, and let $\F$ be an almost quasi-coherent $\O_X$-module. Then $\F$ is almost finitely presented (resp. almost finitely generated) if and only if for any finitely generated ideal $\m_0\subset \m$ there is a morphism $f\colon \G \to \F$ such that $\G$ is a quasi-coherent finitely presented (resp. finitely generated) $\O_X$-module , $\m_0(\ker f)=0$ and $\m_0(\coker f)=0$.
\end{thm}

We now discuss the content of Sections~\ref{acoh-sheaves-formal}-\ref{section:derived-functors-formal-schemes}. The main goal of these sections is to show an analogue of Theorem~\ref{intro:thm-local-properties-schemes} for a class of formal schemes. To achieve this, we restrict our attention to the class of topologically finitely presented schemes over a topologically universally adhesive ring $R$ (see Setup~\ref{set-up3}). This, in particular, includes admissible formal schemes over a mixed characteristic, $p$-adically complete rank-$1$ valuation ring $\O_C$ with algebraically closed fraction field $C$. \smallskip

One of the main difficulties in developing a good theory of almost coherent $\O_\X^a$-modules on a formal scheme $\X$ is that there is no good abelian theory of ``quasi-coherent'' on $\X$. The theory of quasi-coherent sheaves is an important tool used in developing the theory of almost coherent sheaves on schemes that does not have an immediate counterpart in the world of formal schemes. \smallskip

We overcome this issue in two different ways: we use the notion of adically quasi-coherent $\O_\X$-modules introduced in \cite{FujKato} (see Definition~\ref{defn:adically-quasi-coherent-sheaves}) and the notion of derived quasi-coherent $\O_\X$-modules introduced in \cite{Lurie-spectral} (see Definition~\ref{defn:derived-quasi-coherent-sheaves}). The first notion has the advantage that every adically quasi-coherent $\O_\X$-module is an actual $\O_\X$-module. However, these modules do not form a weak Serre subcategory inside $\bf{Mod}_{\O_\X}$, so they are not always very useful in practice. The latter definition has the advantage that derived quasi-coherent $\O_\X$-modules form a triangulated subcategory inside $\bf{D}(\X)$; it is quite convenient for certain purposes. However, derived quasi-coherent $\O_\X$-modules are merely objects of $\bf{D}(\X)$ and not actual $\O_\X$-modules in the classical sense. Therefore, we usually use adically quasi-coherent $\O_\X$-modules when needed except for Section~\ref{section:derived-category-formal-schemes}, where the notion of derived quasi-coherent $\O_\X$-modules seems to be more useful for our purposes. In particular, it allows us to define the functor 
\[
(-)^{L\Updelta}\colon \bf{D}_{acoh}(A)^a \to \bf{D}_{acoh}(\Spf A)^a
\]
for any topologically finitely presented $R$-algebra $A$ in a way that ``extends'' the classical functor $(-)^\Updelta\colon \bf{Mod}_A^{\rm{acoh}} \to \bf{Mod}_{\O_\X}$ (see Definition~\ref{defn:updelta-derived} and Lemma~\ref{updelta-derived}). 

\begin{thm}\label{intro:thm-local-properties-formal-schemes}(Lemma~\ref{main-coh-glob-formal}, Corollary~\ref{cor:derived-almost-coherent-sheaves-formal-affine-schemes}, Lemmas~\ref{tensor-almost-coh-formal-derived}, \ref{pullback-almost-coh-formal-derived}, \ref{hom-alcoh-formal-derived})  Let $R$ be a ring with a finitely generated ideal $I$ such that $R$ is $I$-adically complete, $I$-adically topologically universally adhesive, $I$-torsion free with an ideal $\m$ such that $I\subset \m$, $\m^2=\m$ and $\widetilde{\m}\coloneqq \m\otimes_R \m$ is $R$-flat. 
\begin{enumerate}[label=\textbf{(\arabic*)}]
    \item For any topologically finitely presented formal $R$-scheme $\X$, almost coherent $\O_\X^a$-modules form a weak Serre subcategory of $\bf{Mod}_{\O_\X}^a$;
    \item The functor 
    \[
        \bf{R}\Gamma(\Spf R, -) \colon \bf{D}_{acoh}(\Spf R)^a \to \bf{D}_{acoh}(R)^a
    \]
    is a $t$-exact equivalence of triangulated categories;
    \item The natural morphism $(M^a\otimes^L_{R^a} N^a)^{L\Updelta} \to (M^a)^{L\Updelta}\otimes^L_{\O_{\Spf R}^a}(N^a)^{L\Updelta}$
    is an isomorphism for any $M^a, N^a\in \bf{D}_{acoh}(R)^a$;
    \item Let $\mf\colon \Spf B \to \Spf A$ be a morphism of topologically finitely presented affine formal $R$-schemes. Then $\bf{L}\mf^*\left(\left(M^a\right)^{L\Updelta}\right)$ is functorially isomorphic to $(M^a\otimes^L_{A^a} B^a)^{L\Updelta}$ for any $M^a\in \bf{D}_{acoh}(A)^a$;
    \item The natural morphisms
    \[
        \left(\bf{R}\rm{alHom}_{R^a}\left(M^a,N^a\right)\right)^{L\Updelta} \to \bf{R}\ud{al\cal{H}om}_{\O_{\Spf R}^a}\left(\left(M^a\right)^{L\Updelta}, \left(N^a\right)^{L\Updelta}\right),
    \]
    \[
        \left(\bf{R}\rm{Hom}_{R^a}\left(M^a,N^a\right)\right)^{L\Updelta} \to \bf{R}\ud{\cal{H}om}_{\O_{\Spf R}^a}\left(\left(M^a\right)^{L\Updelta}, \left(N^a\right)^{L\Updelta}\right)
    \]
are almost isomorphisms for $M^a\in \bf{D}^-_{acoh}(R)^a$, $N^a\in \bf{D}^+_{acoh}(R)^a$.
\end{enumerate}
\end{thm}

Similarly to the case of schemes, almost coherent sheaves on formal schemes satisfy the global approximation property: 

\begin{thm}\label{intro:thm-approximate-formal}(Theorem~\ref{thm:approximate-formal}) Let $R$ be as in Theorem~\ref{intro:thm-local-properties-formal-schemes}, let $\X$ be a finitely presented formal $R$-scheme, and let $\F$ be an almost finitely generated (resp. almost finitely presented) $\O_\X$-module. Then, for any finitely generated ideal $\m_0\subset \m$, there is an adically quasi-coherent, finitely generated (resp. finitely presented) $\O_\X$-module $\G$ and a map $\phi\colon \G \to \F$ such that $\m_0(\coker\phi)=0$ and $\m_0(\ker\phi)=0$.
\end{thm}

We discuss the global properties of almost coherent sheaves in Chapter~\ref{section:cohomological-properties}. Namely, we generalize certain cohomological properties of classical coherent sheaves to the case of almost coherent sheaves. We start with the almost version of the Proper Mapping Theorem:

\begin{thm}\label{thm:intro-almost-proper-schemes}(Theorem~\ref{almost-proper-mapping}) Let $R$ be a universally coherent\footnote{Any finitely presented $R$-algebra $A$ is coherent} ring with an ideal $\m$ such that $\widetilde{\m}\coloneqq \m\otimes_R \m$ is $R$-flat and $\m^2=\m$. And let $f\colon X \to Y$ be a proper morphism of finitely presented $R$-schemes with quasi-compact $Y$. Then $\bf{R}f_*$ carries $\bf{D}^*_{acoh}(X)^a$ to $\bf{D}^*_{acoh}(Y)^a$ for $*\in \{``\ \text{''}, -, +, b\}$.
\end{thm}

The essential idea of the proof is to reduce Theorem~\ref{thm:intro-almost-proper-schemes} to the classical Proper Mapping Theorem over a universally coherent base \cite[Theorem I.8.1.3]{FujKato}. The key input to make this reduction work is Theorem~\ref{thm:intro-global-approximation}. 

We also prove a version of the Almost Proper Mapping Theorem for a morphism of formal schemes:

\begin{thm}\label{thm:intro-almost-proper-formal-schemes} (Theorem~\ref{almost-proper-mapping-formal}) Let $\Y$ be a topologically finitely presented formal $R$-scheme for $R$ as in Setup~\ref{set-up3}. And let $\mf\colon \X \to \Y$ be a proper, topologically finitely presented morphism. Then $\bf{R}\mf_*$ carries $\bf{D}^*_{acoh}(\X)^a$ to $\bf{D}^*_{acoh}(\Y)^a$ for $*\in \{``\ \text{''}, -, +, b\}$.
\end{thm}

Then we characterize quasi-coherent, almost coherent complexes on finitely presented, separated schemes over a universally coherent base ring $R$. This is an almost analogue of \cite[\href{https://stacks.math.columbia.edu/tag/0CSI}{Tag 0CSI}]{stacks-project}. We follow the same proof strategy but adjust it in certain places to make it work in the almost setting. This result is important for us as it will later play a crucial role in the proof of the Formal GAGA Theorem for almost coherent sheaves. 

\begin{thm}\label{thm:intro-check-perfect}(Theorem~\ref{check-perfect}) Let $R$ be a universally coherent ring with an ideal $\m$ such that $\widetilde{\m}\coloneqq \m\otimes_R \m$ is $R$-flat and $\m^2=\m$. Let $X$ be a separated, finitely presented $R$-scheme, and let $\F \in \mathbf{D}^-_{qc}(X)$ be an object such that 
\[
\mathbf{R}\Hom_X(\mathcal P, \F) \in \mathbf{D}_{acoh}^-(R)
\]
for any $\mathcal P \in \rm{Perf}(X)$. Then $\F\in \mathbf{D}^-_{qc, acoh}(X)$. 
\end{thm}

\begin{thm}\label{thm:intro-formal-GAGA}(Corollary~\ref{cor:almost-GAGA}) Let $R$ be as in Theorem~\ref{thm:intro-almost-proper-formal-schemes}, and let $X$ be a finitely presented $R$-scheme. Then the functor 
\[
\mathbf{L}c^*\colon \mathbf{D}_{acoh}^*(X)^a \to \mathbf{D}_{acoh}^*(\X)^a
\]
induces an equivalence of categories for $*\in  \{\text{`` ''}, +, -, b\}$. 
\end{thm}

We note that the standard proof of the classical formal GAGA theorem via projective methods has no chance of working in the almost coherent situation (due to the lack of ``finiteness'' for almost coherent sheaves). Instead, we ``explicitly'' construct a pseudo-inverse to $\mathbf{L}c^*$ in the derived world by adapting an argument from the paper of J.\,Hall \cite{JH}.  \smallskip{}

The last thing we discuss in Section~\ref{section:cohomological-properties} is the almost version of the Grothendieck Duality. This is a crucial technical tool in our proof of Poincar\'e Duality in \cite{Zav4}. So we develop some foundations of the $f^!$ functor in the almost world in this manuscript. We summarize the main properties of this functor below:  

\begin{thm}\label{thm:intro-almost-grothendieck-duality} (Theorem~\ref{thm:almost-grothendieck-duality}) Let $R$ be as in Theorem~\ref{thm:intro-almost-proper-schemes}, and $\rm{FPS}_R$ be the category of finitely presented, separated $R$-schemes. Then there is a well-defined pseudo-functor $(-)^!$ from $\rm{FPS}_R$ into the $2$-category of categories such that
\begin{enumerate}[label=\textbf{(\arabic*)}]
    \item $(X)^!=\bf{D}^+_{aqc}(X)^a$,
    \item for a smooth morphism $f\colon X \to Y$ of pure relative dimension $d$, $f^! \simeq \bf{L}f^*(-)\otimes^L_{\O_X^a} \Omega^d_{X/Y}[d]$.
    \item for a proper morphism $f\colon X \to Y$, $f^!$ is right adjoint to $\bf{R}f_*\colon \bf{D}^+_{acoh}(X)^a\to \bf{D}^+_{acoh}(Y)^a$.
\end{enumerate}
\end{thm}

\subsection{$\O^+/p$, $\O^+$, and $\O$-vector bundles (Section~\ref{section:proetale-v-sites})}

The main goal of Section~\ref{section:almost-coherent-nearby-cycles} is to study the categories of $\O^+/p$-vector bundles in the \'etale, quasi-pro\'etale, and $v$-topologies. We also show that $\O^+/p$-vector bundles can be trivialized by some particular \'etale covers. These results will play a crucial role in Section~\ref{section:almost-coherent-nearby-cycles}. Also, as an application of our results, we give a new proof of the theorem of Kedlaya--Liu saying that, for a perfectoid space $X$, the categories of $\O$-vector bundles in the analytic and $v$-topologies are equivalent. \smallskip

We formulate the results of this section more precisely below:

\begin{thm}[Corollary~\ref{cor:different-vector-bundles-equivalent}]\label{thm:intro-different-vector-bundles-equivalent} Let $X$ be a strongly sheafy adic space over $\Spa(\Q_p, \Z_p)$ (see Definition~\ref{defn:strongly-sheafy}). Then 
\begin{enumerate}[label=\textbf{(\arabic*)}]
    \item\label{thm:intro-different-vector-bundles-equivalent-1} the categories $\rm{Vect}(X_{\et}; \O_{X_\et}^+/p)$, $\rm{Vect}(X^\diam_\qp; \O_{X^\diam_\qp}^+/p)$, and $\rm{Vect}(X^\diam_v; \O_{X^\diam}^+/p)$ are equivalent;
    \item\label{thm:intro-different-vector-bundles-equivalent-2} These equivalences preserve cohomology groups;
    \item\label{thm:intro-different-vector-bundles-equivalent-3} for any $\O_{X^\diam}^+/p$-vector bundle $\cal{E}$ and a point $x\in X$, there exists an open affinoid subspace $x\in U_x\subset X$ and a finite \'etale surjective morphism $\widetilde{U}_x \to U_x$ such that $\cal{E}|_{\widetilde{U}_x^\diam}$ is a free vector bundle.
\end{enumerate}
\end{thm}

Theorem~\ref{thm:intro-different-vector-bundles-equivalent}\ref{thm:intro-different-vector-bundles-equivalent-1},\ref{thm:intro-different-vector-bundles-equivalent-2} is essentially due to B.\,Heuer (see \cite[\textsection 2]{Heuer-G-torsor} for a similar result in a slightly different level of generality). However, Theorem~\ref{thm:intro-different-vector-bundles-equivalent}\ref{thm:intro-different-vector-bundles-equivalent-3} does not seem to follow from \cite{Heuer-G-torsor} and is crucial for our arguments in Chapter~\ref{section:almost-coherent-nearby-cycles}. \smallskip

We also prove a version of Theorem~\ref{thm:intro-different-vector-bundles-equivalent} for $\O^+$-vector bundles: 

\begin{thm}[Theorem~\ref{thm:integral-Kedlaya-Liu}, Corollary~\ref{cor:etale-locally-trivial}]\label{thm:intro-integral-Kedlaya-Liu} Let $X$ be a perfectoid space over $\Spa(\Q_p, \Z_p)$. Then 
\begin{enumerate}[label=\textbf{(\arabic*)}]
    \item\label{thm:intro-integral-Kedlaya-Liu-1} the categories $\rm{Vect}(X_{\et}; \O_{X_\et}^+)$, $\rm{Vect}(X^\diam_\qp; \O_{X^\diam_\qp}^+)$, and $\rm{Vect}(X^\diam_v; \O_{X^\diam}^+)$ are equivalent;
    \item\label{thm:intro-integral-Kedlaya-Liu-2} These equivalences preserve cohomology groups;
    \item\label{thm:intro-integral-Kedlaya-Liu-3} for any $\O_{X^\diam}^+$-vector bundle $\cal{E}$ and a point $x\in X$, there exists an open affinoid subspace $x\in U_x\subset X$ and a finite \'etale surjective morphism $\widetilde{U}_x \to U_x$ such that $\cal{E}|_{\widetilde{U}_x^\diam}$ is a free vector bundle.
\end{enumerate}
\end{thm}

We also refer to Theorem~\ref{thm:integral-Kedlaya-Liu} for a slightly more precise statement. As an application of our methods, we can also deduce the following theorem of Kedlaya--Liu: 

\begin{thm}[{\cite[Theorem 3.5.8]{KedLiu2}, \cite[Lemma 17.1.8]{Berkeley-notes}, \cite[Theorem 4.27]{Heuer-G-torsor}}, Thereom~\ref{thm:rational-vector-bundles-on-perfectoids}]\label{thm:intro-KL} Let $X$ be a perfectoid space over $\Spa(\Q_p, \Z_p)$.
\begin{enumerate}[label=\textbf{(\arabic*)}]
    \item\label{thm:intro-KL-1} the categories $\rm{Vect}(X_{\rm{an}}, \O_X)$, $\rm{Vect}(X_{\et}; \O_{X_\et})$, $\rm{Vect}(X^\diam_\qp; \O_{X^\diam_\qp})$, and $\rm{Vect}(X^\diam_v; \O_{X^\diam})$ are equivalent. Furthermore, if $X=\Spa(R, R^+)$ is an affinoid perfectoid, all these categories are equivalent to the category of finite projective $R$-modules;
    \item\label{thm:intro-KL-2} these equivalences preserve cohomology groups.
\end{enumerate}
\end{thm}

We note that the proof of Theorem~\ref{thm:intro-KL} is quite different from the proofs of \cite[Theorem 3.5.8]{KedLiu2} and \cite[Lemma 17.1.8]{Berkeley-notes}. However, it is quite similar to the proof of \cite[Theorem 4.27]{Heuer-G-torsor} (with appropriate simplifications). We also note that \cite[Theorem 4.27]{Heuer-G-torsor} proves a stronger result that applies to $G$-torsors for any rigid group $G$. We also show that any $\O$-vector bundle in the $v$-topology admits an $\O^+$-lattice after a very explicit \'etale cover:

\begin{thm}[Corollary~\ref{cor:etale-locally-trivial-precise}]\label{thm:intro-etale-locally-trivial-precise} Let $X$ be a strongly sheafy adic space over $\Spa(\Q_p, \Z_p)$, and let $\cal{E}$ be an $\O_{X^\diam}$-vector bundle. Then, for each $x\in X$, there is an open subspace $x\in U_x \subset X$, a finite \'etale surjective morphism $\widetilde{U}_x \to U_x$, and an $\O_{\widetilde{U}^\diam_x}^+$-vector bundle $\cal{E}_x^+$ such that $\cal{E}_x^+\big[\frac{1}{p}\big] \simeq \cal{E}|_{\widetilde{U}_x}$.
\end{thm}

\subsection{$p$-adic nearby cycles sheaves (Section \ref{section:almost-coherent-nearby-cycles})}

The main goal of Section~\ref{section:almost-coherent-nearby-cycles} is to give the main non-trivial example of almost coherent sheaves: the $p$-adic nearby cycles sheaves. \smallskip

We fix a $p$-adic perfectoid field $K$ and a rigid-analytic variety $X$ over $K$ with an admissible formal $\O_K$-model $\X$. \smallskip

The rigid-analytic variety $X$ comes with a morphism of ringed sites 
\[
\nu\colon (X^\diam_v, \O_{X^\diam}^+) \to (\X_{\rm{Zar}}, \O_\X)
\]
and a morphism 
\[
\nu\colon (X^\diam_v, \O_{X^\diam}^+/p) \to (\X_{0, \rm{Zar}}, \O_{\X_0})
\]
where $\X_0$ is the mod-$p$ fiber of $\X$, $X^\diam_v$ is the $v$-site of the associated diamond (see Appendix~\ref{section:v-topology}), and $\O_{X^\diam}^+$ is its integral ``untilted'' structure sheaf (see Definition~\ref{defn:v-structure-sheaf}). \smallskip

The main goal of Section~\ref{section:almost-coherent-nearby-cycles} is to show that the nearby cycles functor $\bf{R}\nu_*$ sends some class of $\O_{X^\diam}^+/p$-sheaves to complexes of almost coherent $\O_{\X_0}$-modules. More precisely, we show that, for any $\O_{X^\diam}^+/p$-vector bundle $\cal{E}$, the complex $\bf{R}\nu_*\cal{E}$ has quasi-coherent and almost coherent cohomology sheaves. We also give a bound on its almost cohomological dimension.

\begin{thm}\label{intro:thm-main-thm-small}(Theorem~\ref{thm:main-thm-small}) Let $\X$ be an admissible formal $\O_K$-scheme with adic generic fiber $X$ of dimension $d$ and mod-$p$ fiber $\X_0$, and let $\cal{E}$ be an $\O^+_{X^{\diam}}/p$-vector bundle. Then 
\begin{enumerate}[label=\textbf{(\arabic*)}]
    \item $\bf{R}\nu_*\cal{E}\in \bf{D}^+_{qc, acoh}(\X_0)$ and $(\bf{R}\nu_*\cal{E})^a\in \bf{D}^{[0, 2d]}_{acoh}(\X_0)^a$;
    \item if $\X=\Spf A$ is affine, then the natural map 
    \[
    \widetilde{\rm{H}^i\left(X^\diam_v, \cal{E} \right)} \to \rm{R}^i\nu_*\left(\cal{E}\right)
    \]
    is an isomorphism for every $i\geq 0$;
    \item the formation of $\rm{R}^i\nu_*(\cal{E})$ commutes with \'etale base change, i.e., for any \'etale morphism $\mf \colon \Y \to \X$ with adic generic fiber $f\colon Y\to X$, the natural morphism 
    \[
    \mf^*_0 \left(\rm{R}^i\nu_{\X, *}(\cal{E}) \right)\to \rm{R}^i\nu_{\Y, *}\left(\cal{E}|_{Y^\diam}\right)
    \]
    is an isomorphism for any $i\geq 0$;
    \item if $\X$ has an open affine covering $\X=\bigcup_{i\in I} \sU_i$ such that $\cal{E}|_{(\sU_{i, K})^\diam}$ is small (see Definition~\ref{defn:small-mod-p}), then
    \[
    \left(\bf{R}\nu_{*}\cal{E}\right)^a \in \bf{D}^{[0, d]}_{acoh}(\X_0)^a;
    \]
    \item there is an admissible blow-up $\X'\to \X$ such that $\X'$ has an open affine covering $\X'=\bigcup_{i\in I} \sU_i$ such that $\cal{E}|_{(\sU_{i, K})^\diam}$ is small. 
    
    In particular, there is a cofinal family of admissible formal models $\{\X'_i
    \}_{i\in I}$ of $X$ such that 
    \[
    \left(\bf{R}\nu_{\X'_i, *}\cal{E}\right)^a\in \bf{D}^{[0, d]}_{acoh}(\X'_{i, 0})^a.
    \]
    for each $i\in I$. 
\end{enumerate}
\end{thm}

\begin{rmk} We note that Theorem~\ref{intro:thm-main-thm-small} implies that the nearby cycles complex $\bf{R}\nu_*\cal{E}$ is quasi-coherent on the nose (as opposed to being almost quasi-coherent). This is quite unexpected to the author since all previous results on the cohomology groups of $\O^+/p$ were only available in the almost category. 
\end{rmk}

\begin{rmk} We do not know if an admissible blow-up $\X' \to \X$ in the formulation of Theorem~\ref{intro:thm-main-thm-small} is really necessary or just an artifact of the proof. More importantly, we do not know if, for every $\O_{X^\diam}^+/p$-vector bundle $\cal{E}$, there is an admissible formal model $\X$ such that the ``nearby cycles'' sheaf $\bf{R}\nu_{\X, *}\cal{E}$ lies in $\bf{D}^{[0, d]}_{acoh}(\X_0)^a$. 
\end{rmk}

The proof of Theorem~\ref{intro:thm-main-thm-small} crucially uses Theorem~\ref{thm:intro-different-vector-bundles-equivalent}, and especially Theorem~\ref{thm:intro-different-vector-bundles-equivalent}\ref{thm:intro-different-vector-bundles-equivalent-3}. \smallskip

Another family of sheaves for which we can establish a good behavior of ``nearby cycles'' is given by sheaves of the form $\F\otimes \O_{X^\diam}^+/p$ for a Zariski-constructible \'etale sheaf of $\bf{F}_p$-modules (see Definition~\ref{defn:zc}). Namely, in this case, we can get a better cohomological bound and show that nearby cycles almost commute with proper base change, as this happens in algebraic geometry.

\begin{thm}\label{intro:main-thm}(Theorem~\ref{main-thm} and Lemma~\ref{lemma:nearby-pushforward}) Let $\X$ be an admissible formal $\O_K$-scheme with adic generic fiber $X$ of dimension $d$ and mod-$p$ fiber $\X_0$, and let $\F \in \bf{D}^{[r, s]}_{zc}(X; \bf{F}_p)$. 
Then 
\begin{enumerate}[label=\textbf{(\arabic*)}]
    \item there is an isomorphism $\bf{R}t_*\left(\F\otimes \O_{X_\et}^+/p\right) \simeq \bf{R}\nu_*\left(\F\otimes \O_{X^\diam}^+/p\right)$;
    \item $\bf{R}\nu_*\left(\F\otimes \O_{X^\diam}^+/p\right) \in \bf{D}^+_{qc, acoh}(\X_0)$, and $\bf{R}\nu_*\left(\F\otimes \O_{X^\diam}^+/p\right)^a \in \bf{D}^{[r, s+d]}_{acoh}(\X_0)^a$;
    \item if $\X=\Spf A$ is affine, then the natural map 
    \[
    \widetilde{\rm{H}^i\left(X^\diam_v, \F\otimes \O_{X^\diam}^+/p \right)} \to \rm{R}^i\nu_*\left(\F\otimes \O_{X^\diam}^+/p\right)
    \]
    is an isomorphism for every $i\geq 0$;
    \item the formation of $\rm{R}^i\nu_*\left(\F\otimes \O_{X^\diam}^+/p\right)$ commutes with \'etale base change, i.e., for any \'etale morphism $\mf \colon \Y \to \X$ with adic generic fiber $f\colon Y\to X$, the natural morphism 
    \[
    \mf^*_0 \left(\rm{R}^i\nu_{\X, *}\left(\F\otimes \O_{X^\diam}^+/p\right) \right)\to \rm{R}^i\nu_{\Y, *}\left(f^{-1}\F\otimes \O_{Y^\diam}^+/p\right)
    \]
    is an isomorphism for any $i\geq 0$;
    \item if $\mf\colon \X \to \Y$ is a proper morphism of admissible formal $\O_K$-schemes with adic generic fiber $f\colon X \to Y$, then the natural morphism
    \[
    \bf{R}\nu_{\Y, *}\left(\bf{R}f_{*} \F\otimes \O^+_{Y^\diam}/p\right) \to \bf{R}\mf_{0, *}\left(\bf{R}\nu_{\X, *}\left(\F\otimes \O_{X^\diam}^+/p\right)\right)
    \]
is an almost isomorphism. 
\end{enumerate}
\end{thm}

We also show an integral version of Theorem~\ref{intro:thm-main-thm-small}:

\begin{thm}\label{intro:thm-main-thm-integral}(Theorem~\ref{thm:main-thm-integral}) Let $\X$ be an admissible formal $\O_K$-scheme with adic generic fiber $X$ of dimension $d$, and let $\cal{E}$ be an $\O^+_{X^{\diam}}$-vector bundle. 
Then 
\begin{enumerate}[label=\textbf{(\arabic*)}]
    \item $\bf{R}\nu_*\cal{E}\in \bf{D}^+_{qc, acoh}(\X)$ and $(\bf{R}\nu_*\cal{E})^a\in \bf{D}^{[0, 2d]}_{acoh}(\X)^a$;
    \item if $\X=\Spf A$ is affine, then the natural map 
    \[
    \rm{H}^i\left(X^\diam_v, \cal{E} \right)^{\Updelta} \to \rm{R}^i\nu_*\left(\cal{E}\right)
    \]
    is an isomorphism for every $i\geq 0$;
    \item the formation of $\rm{R}^i\nu_*(\cal{E})$ commutes with \'etale base change, i.e., for any \'etale morphism $\mf \colon \Y \to \X$ with adic generic fiber $f\colon Y\to X$, the natural morphism 
    \[
    \mf^* \left(\rm{R}^i\nu_{\X, *}(\cal{E}) \right)\to \rm{R}^i\nu_{\Y, *}\left(\cal{E}|_{Y^\diam}\right)
    \]
    is an isomorphism for any $i\geq 0$;
    \item if $\X$ has an open affine covering $\X=\bigcup_{i\in I} \sU_i$ such that $\cal{E}|_{(\sU_{i, K})^\diam}$ is small (see Definition~\ref{defn:small-integrally}), then
    \[
    \left(\bf{R}\nu_{*}\cal{E}\right)^a \in \bf{D}^{[0, d]}_{acoh}(\X)^a;
    \]
    \item there is an admissible blow-up $\X'\to \X$ such that $\X'$ has an open affine covering $\X'=\bigcup_{i\in I} \sU_i$ such that $\cal{E}|_{(\sU_{i, K})^\diam}$ is small. 
    
    In particular, there is a cofinal family of admissible formal models $\{\X'_i
    \}_{i\in I}$ of $X$ such that 
    \[
    (\bf{R}\nu_{\X'_i, *}\cal{E})^a\in \bf{D}^{[0, d]}_{acoh}(\X'_{i})^a.
    \]
    for each $i\in I$. 
\end{enumerate}
\end{thm}

Theorem~\ref{intro:thm-main-thm-integral} has an interesting consequence saying that $v$-cohomology groups of any $\O_{X^\diam}^+$-vector bundle are almost coherent and almost vanish in degrees larger than $2\dim X$. This (together with Theorem~\ref{thm:intro-1}) indicates that there should probably be stronger (almost) finiteness results for some bigger class $\cal{O}_{X^\diam}^+$-modules.

\begin{thm}\label{intro:thm-finiteness-bundles}(Theorem~\ref{thm:finiteness-bundles}) Let $K$ be a $p$-adic perfectoid field, let $X$ be a proper rigid-analytic $K$-variety of dimension $d$, and let be $\cal{E}$ an $\O_{X^\diam}^+$-vector bundle (resp. $\O_{X^\diam}^+/p$-vector bundle). Then 
\[
\bf{R}\Gamma(X_v^\diam, \cal{E}) \in \bf{D}^{[0, 2d]}_{acoh}(\O_K)^a.
\]
\end{thm}

We now explain the main steps of our proofs of Theorems~\ref{intro:thm-main-thm-small}~and~\ref{intro:thm-main-thm-integral} for $\cal{E}=\O_{X^\diam}^+/p$ and $\cal{E}=\O_{X^\diam}^+$ respectively:

\begin{proof}[Proof Sketch]
\begin{enumerate}[label=\textbf{(\arabic*)}]
    \item We first show that the sheaves $\rm{R}^i\nu_*(\O_{X^\diam}^+/p)$ are  quasi-coherent. The main key input is that cohomology of $\O_{X^\diam}^+/p$-vector bundles vanish on strictly totally disconnected spaces (see Definition~\ref{defn:totally-disconnected}), and that each affinoid rigid-analytic variety admits a $v$-covering such that all terms of its \v{C}ech nerve are strictly totally disconnected.
    \item The same ideas can be used to show that the formation of $\rm{R}^i\nu_*(\O_{X^\diam}^+/p)$ commutes with \'etale base change. 
    \item We show that the $\O_{\X_0}$-modules $\rm{R}^i\nu_*\left(\O_{X^\diam}^+/p\right)$ are almost coherent for smooth $X$. This is done in three steps: first, we find an admissible blow-up $\X' \to \X$ such that $\X'$ has an open affine covering $\X' = \bigcup_{i\in I} \sU_i$ such that each $\sU_i=\Spf A_i$ admits a finite rig-\'etale morphism to $\wdh{\bf{A}}^d_{\O_K}$, then we show that the cohomology groups $\rm{H}^i(\sU^\diam_{i, K, v}, \O_{X^\diam}^+/p)$ are almost coherent over $A_i/pA_i$, then we conclude almost coherence of $\rm{R}^i\nu_*\left(\O_{X^\diam}^+/p\right)$.
    
    The first step is the combination of \cite[Proposition 3.7]{BLR3} and Theorem~\ref{etale-finite-etale-formal}. The first result allows us to choose an admissible blow-up $\X' \to \X$ with an open affine covering $\X'=\bigcup_{i\in I} \sU_i$ such that each $\sU_i$ admits a rig-\'etale morphism $\sU_i \to \wdh{\bf{A}}^d_{\O_K}$. Then Theorem~\ref{etale-finite-etale-formal} guarantees that we can change these morphisms so that $\sU_i \to \wdh{\bf{A}}^d_{\O_K}$ are {\it finite} and rig-\'etale.  
    
    The second step follows the strategy presented in \cite{Sch1}. We construct an explicit affinoid perfectoid cover of $\sU_i$ that is a $\bf{Z}_p(1)^d$-torsor. So we reduce studying of $\rm{H}^i(\sU^\diam_{i, K, v}, \O_{X^\diam}^+/p)$ to studying cohomology groups of $\bf{Z}_p(1)^d$ that can be explicitly understood via the Koszul complex. 
    
    The last step is the consequence of the Almost Proper Mapping Theorem~\ref{thm:intro-almost-proper-schemes} and the already obtained results.
    
   \item The next step is to show that $\rm{R}^i\nu_*\left(\O_{X^\diam}^+/p\right)$ is almost coherent for a general $X$. This is done by choosing a proper hypercovering by smooth spaces $X_{\bullet}$ and then using a version of cohomological $v$-descent to conclude almost coherence of the $p$-adic nearby cycles sheaves. As an important technical tool, we use the theory of diamonds developed in \cite{Sch2}.
   
   \item Next we show that $\bf{R}\nu_* \left(\O_{X^\diam}^+/p\right)$ is almost concentrated in degrees $[0,d]$. This claim is quite subtle. The key input is the version of the purity theorem \cite[Theorem 10.11]{BS3} that implies that any {\it finite} (but not necessarily \'etale) adic space over an affinoid perfectoid space has a diamond that is isomorphic to a diamond of an affinoid perfectoid space. This allows us to reduce the question of cohomological bounds of $\bf{R}\nu_* \left(\O_{X^\diam}^+/p\right)^a$ to the question about the cohomological dimension of the pro-finite group $\bf{Z}_p(1)^{d}$. This can be explicitly understood via the Koszul complex again.
   
   \item Finally, we show Theorem~\ref{intro:thm-main-thm-integral} by reducing it to Theorem~\ref{intro:thm-main-thm-small}. The key input is Theorem~\ref{thm:intro-mod-p} that allows us to check finiteness mod-$p$. 
\end{enumerate}
\end{proof}

\subsection{Acknowledgements}
We are very grateful to B.\,Bhatt, B.\,Conrad, S.\,Petrov, and D.\,B.\,Lim for many fruitful discussions. We express additional gratitude to B. Bhatt for bringing \cite[Theorem 10.11]{BS3} and \cite{Guo} to our attention. We are thankful to B.\,Conrad for reading the first draft of this paper and making useful suggestions on how to improve the exposition of this paper. Part of this work was carried out at the mathematics department of the University of Michigan. We thank them for their hospitality. We heartfully thank O.\,Gabber, D.\,Hansen, B.\,Heuer, K.\,Kedlaya, P.\,Scholze, and K.\,Shimomoto for their valuable comments on the previous version of this draft. The later stages of this work were done at the Max Planck Institute for Mathematics in Bonn, Germany. We thank them for the excellent working conditions and funding. 
\subsection{Notation}

A {\it non-archimedean field} $K$ is always assumed to be complete. A non-archimedean field $K$ is called {\it $p$-adic} if its ring of powerbounded elements $\O_K=K^\circ$ is a ring of mixed characteristic $(0, p)$.

We follow \cite[\href{https://stacks.math.columbia.edu/tag/02MN}{Tag 02MN}]{stacks-project} for the definition of a (weak) Serre subcategory of an abelian category $\cal{A}$.  

For an $R$-ringed site $(X, \O_X)$, an element of the derived category $\F\in \bf{D}(X)$, and an element $\varpi\in R$, we denote by $[\F/\varpi]$ the cone of the multiplication by $\varpi$-morphism, i.e.
\[
[\F/\varpi] \coloneqq \rm{cone}(\F \xr{\varpi} \F).
\]

Namely, we say that a non-empty full subcategory $\mathcal{C}$ of an abelian category $\cal{A}$ is a {\it Serre subcategory} if, for any exact sequence $A \to B \to C$ with $A, C \in \mathcal{C}$, we have $B \in \mathcal{C}$. We say that $\cal{C}$ is a {\it weak Serre subcategory} if, for any exact sequence 
\[
A_0 \to A_1 \to A_2 \to A_3 \to A_4
\]
with $A_0, A_1, A_3, A_4 \in \mathcal{C}$, we have $A_2 \in \mathcal{C}$. Look at \cite[\href{https://stacks.math.columbia.edu/tag/02MP}{Tag 02MP}]{stacks-project} and \cite[\href{https://stacks.math.columbia.edu/tag/0754}{Tag 0754}]{stacks-project} for an alternative way to describe (weak) Serre subcategories. \smallskip

If $\cal{C}$ is a Serre subcategory of an abelian category $\cal{A}$ we define the {\it quotient category} as a pair $(\cal{A}/\cal{C}, F)$ of an abelian category $\cal{A}/\cal{C}$ and an exact functor
\[
F\colon \cal{A} \to \cal{A}/\cal{C}
\]
such that, for any exact functor $G\colon \cal{A} \to \cal{B}$ to an abelian category $\cal{B}$ with $\cal{C} \subset \ker G$, there is a factorization $G=H\circ F$ for a unique exact functor $H\colon \cal{A}/\cal{C} \to \cal{B}$. The quotient category always exists by \cite[\href{https://stacks.math.columbia.edu/tag/02MS}{Tag 02MS}]{stacks-project}. \smallskip

If $\cal{B}$ is a full triangulated subcategory of a triangulated category $\cal{D}$ we define the {\it Verdier quotient} as a pair $(\cal{D}/\cal{B}, F)$ of a triangulated category $\cal{D}/\cal{B}$ and an exact functor
\[
F\colon \cal{D} \to \cal{D}/\cal{B}
\]
such that, for any exact functor $G\colon \cal{D} \to \cal{D'}$ to a pre-triangulated category $\cal{D'}$ with $\cal{B} \subset \ker G$, there is a factorization $G=H\circ F$ for a unique exact functor $H\colon \cal{D}/\cal{B} \to \cal{D'}$. The Verdier quotient always exists by \cite[\href{https://stacks.math.columbia.edu/tag/05RJ}{Tag 05RJ}]{stacks-project}. \smallskip

We say that a diagram of categories 
\[
\begin{tikzcd}
\cal{A} \arrow{r}{f} \arrow[d, swap, "h"]& \cal{B} \arrow{d}{g} \\
\cal{C} \arrow{r}{k} \arrow[ru, Rightarrow, "\alpha"]& \cal{D}
\end{tikzcd}
\]
is {\it $(2, 1)$-commutative} if $\alpha\colon k\circ h \Rightarrow g\circ f$ is a natural isomorphism of functors. \smallskip{}

For an abelian group $M$ and commuting endomorphisms $f_1, \dots, f_n$, we define the {\it Koszul complex}
\[
K(M; f_1,\dots, f_n) \coloneqq M \to M \otimes_{\Z} \Z^n \to M\otimes_{\Z} \wedge^2 \left(\Z^n\right) \to \dots \to M \otimes_{\Z} \wedge^n \left(\Z^n\right) 
\]
viewed as a chain complex in cohomological degrees $0, \dots, n$. The differential 
\[
d^k\colon M\otimes_{\Z} \wedge^k\left(\Z^n\right) \simeq \bigoplus_{1\leq i_1< \dots <i_k \leq n} M \to M\otimes_{\Z} \wedge^{i+1}\left(\Z^n\right) \simeq \bigoplus_{1\leq j_1< \dots <j_{k+1} \leq n} M
\]
from $M$ in spot $i_1 < \dots < i_k$ to $M$ in spot $j_1 < \dots < j_{k+1}$ is nonzero only if $\{i_1, \dots, i_k\} \subset \{j_1, \dots, j_{k+1}\}$, in which case it is given by $(-1)^{m-1}f_{j_m}$, where $m \in \{1, \dots , k + 1\}$ is the unique integer such that $j_m \notin \{i_1, \dots , i_k\}$. 

If $M$ is an $R$-module and $f_i$ are elements of $R$ the complex $K(M; f_1, \dots, f_n)$ is a complex of $R$-modules and can be identified with 
\[
M \to M \otimes_{R} R^n \to M\otimes_{R} \wedge^2 \left(R^n\right) \to \dots \to M \otimes_{R} \wedge^n \left(R^n\right).
\]

\section{Almost commutative algebra}\label{almost-commutative-algebra}

This chapter is devoted to the study of almost coherent modules. We recall some basic definitions of almost mathematics in Section~\ref{almost-mathematics}. Then we discuss the main properties of almost finitely generated and almost finitely presented modules in Section~\ref{section-almost-finitely-presented}. These two sections closely follow the discussion of almost mathematics in \cite{GR}. Section~\ref{section-almost-coherent} is dedicated to almost coherent modules and almost coherent rings. We show that almost coherent modules form a weak Serre subcategory of $R$-modules, and they coincide with almost finitely presented ones in the case of almost coherent rings. We discuss base change results in Section~\ref{base-change-section}. Finally, we develop some topological aspects of almost finitely generated modules over ``topologically universally adhesive rings'' in Section~\ref{adhesive-almost}.

\subsection{The category of almost modules}\label{almost-mathematics}
We begin this section by recalling basic definitions of almost mathematics from \cite{GR}. We fix a ``base'' ring $R$ with an ideal $\m$ such that $\m^2=\m$ and $\widetilde{\m}=\m\otimes_R \m$ is flat. We always do almost mathematics with respect to $\m$.

\begin{lemma}\label{almost-zero} Let $M$ be an $R$-module. Then the following are equivalent:
\begin{enumerate}[label=\textbf{(\arabic*)}]
	\item The module $\m M$ is the zero module.
	\item The module $\m \otimes_R M$ is the zero module.
	\item The module $\widetilde{\m}\otimes_R M$ is the zero module.
	\item The module $M$ is annihilated by $\e$ for every $\e\in \m$.
\end{enumerate}
\end{lemma}
\begin{proof}
Note that the multiplication map $\m\otimes_R \m \to \m$ is surjective as $\m^2=\m$. This implies that we have surjections 
\[
\widetilde{\m}\otimes_R M \twoheadrightarrow \m\otimes_R M \twoheadrightarrow \m M.
\] 
This shows that \textbf{(3)} implies \textbf{(2)}, and \textbf{(2)} implies \textbf{(1)}. It is clear that \textbf{(2)} implies \textbf{(3)}, and \textbf{(1)} is equivalent to \textbf{(4)}. So the only thing we are left to show is that \textbf{(1)} implies \textbf{(2)}. 

Suppose that $\m M \simeq 0$. Pick an arbitrary basic element $a\otimes m \in \m\otimes_R M$ with $a\in \m$, $m\in M$. Since $\m^2=\m$ there is a finite number of elements $y_1, \dots, y_k, x_1, \dots, x_k \in \m$ such that \[a=\sum_{i=1}^k x_iy_i.\]
Then we have an equality
\[
a\otimes m = \sum_{i=1}^k x_iy_i\otimes m= \sum_{i=1}^k x_i\otimes y_im=0.
\]\end{proof}

\begin{defn} An $R$-module $M$ is  {\it almost zero}, if any of the equivalent conditions of Lemma \ref{almost-zero} is satisfied for $M$.
\end{defn}

\begin{lemma}\label{tilde-m} Under the assumption as above, the ``multiplication'' morphism  $\widetilde{\m}\otimes_R\widetilde{\m} \to \widetilde{\m}$ is an isomorphism.
\end{lemma}
\begin{proof}
We consider a short exact sequence
\[
0 \to \m \to R \to R/\m \to 0.
\]
Note that $ \left(R/\m\right) \otimes_R \m = \m/\m^2=0$, so we get a short exact sequence
\[
0 \to \Tor^R_1(R/\m, \m) \to \widetilde{\m} \to \m\to 0.
\]

Since $\Tor^R_1(R/\m, \m)$ is almost zero, Lemma~\ref{almost-zero} says that after applying the functor $-\otimes_R \widetilde{\m}$ we get an isomorphism
\[
\widetilde{\m}\otimes_R \widetilde{\m} \simeq \m\otimes_R\widetilde{\m}.
\]
Since $\widetilde{\m}$ is $R$-flat, we also see that $\m\otimes_R\widetilde{\m}$ injects into $\widetilde{\m}$. Moreover, it maps isomorphically onto its image $\m\widetilde{\m}=\widetilde{\m}$ as $\m^2=\m$. Taken together, it shows that 
\[
\widetilde{\m}\otimes_R \widetilde{\m} \simeq \widetilde{\m}.
\]
It is straightforward to see that the constructed isomorphism is the ``multiplication'' map.
\end{proof}

We denote by $\Sigma_R$ the category of almost zero $R$-modules considered as a full subcategory of $\textbf{Mod}_R$.

\begin{Cor}\label{cor:sigma-serre-subcategory} The category $\Sigma_R$ is a Serre subcategory of $\textbf{Mod}_R$.\footnote{We refer to   \cite[\href{https://stacks.math.columbia.edu/tag/02MN}{Tag 02MN}]{stacks-project} for the discussion of (weak) Serre categories.} 
\end{Cor}
\begin{proof}
This follows directly from criterion \textbf{(3)} from Lemma \ref{almost-zero}, flatness of $\widetilde{\m}$ and  \cite[\href{https://stacks.math.columbia.edu/tag/02MP}{Tag 02MP}]{stacks-project}.\end{proof}

This corollary allows us to define the quotient category\footnote{We refer to \cite[\href{https://stacks.math.columbia.edu/tag/02MS}{Tag 02MS}]{stacks-project} for the discussion of quotient categories.} $\textbf{Mod}_R^a\coloneqq \textbf{Mod}_R$/$\Sigma_R$ that we call as the category of almost $R$-modules. Note that the localization functor
\[
(-)^a\colon \textbf{Mod}_R \to \textbf{Mod}_R^a
\]
is an exact and essentially surjective functor. We refer to elements of $\bf{Mod}_R^a$ as almost $R$-modules or $R^a$-modules. We will usually denote them by $M^a$ to distinguish almost $R$-modules from $R$-modules. \smallskip

To simplify the exposition, we will use the notation $\bf{Mod}_R^a$ and $\bf{Mod}_{R^a}$ interchangeably. 

\begin{defn}\label{almost-iso-def} A morphism $f: M \to N$ is called {\it an almost isomorphism} (resp. {\it almost injection}, resp. {\it almost surjection}) if the corresponding morphism $f^a:M^a \to N^a$ is an isomorphism (resp. injection, resp. surjection) in $\textbf{Mod}_R^a$.
\end{defn}

\begin{rmk}\label{rmk:main-example} For any $R$-module $M$, the natural morphism $\pi\colon \widetilde{\m}\otimes_R M \to M$ is an almost isomorphism. Indeed, it suffices to show that 
\[
\widetilde{\m} \otimes_R \ker \pi \simeq 0 \text{ and } \widetilde{\m} \otimes_R \coker \pi \simeq 0. 
\]
Using $R$-flatness of $\widetilde{\m}$, we can reduce the question to showing that the map 
\[
\widetilde{\m}\otimes_R \pi \colon \widetilde{\m}\otimes_R \widetilde{\m}\otimes_R M \to \widetilde{\m}\otimes_R M
\]
is an isomorphism. This follows from Lemma~\ref{tilde-m}.
\end{rmk}

\begin{defn}\label{almost-isomorphic} Two $R$-modules $M$ and $N$ are called {\it almost isomorphic} if  $M^a$ is isomorphic to $N^a$ in $\textbf{Mod}_R^a$.
\end{defn}

\begin{lemma}\label{prop-almost} Let $f\colon M \to N$ be a morphism of $R$-modules, then
\begin{enumerate}[label=\textbf{(\arabic*)}]
	\item\label{prop-almost-1} The morphism $f$ is an almost injection (resp. almost surjection, resp. almost isomorphism) if and only if $\ker(f)$ (resp. $\coker(f)$, resp. $\ker(f)$ and $\coker(f)$) is an almost zero module.
        \item\label{prop-almost-2} We have a functorial bijection $\Hom_{R}(\widetilde{\m}\otimes_R M, N) \simeq \Hom_{\textbf{Mod}_R^a}(M^a, N^a)$.
        \item\label{prop-almost-3} Modules $M$ and $N$ are almost isomorphic (not necessary via the morphism $f$) if and only if $\widetilde{\m}\otimes_R M \simeq \widetilde{\m}\otimes_R N$.
\end{enumerate}
\end{lemma}
\begin{proof} \textbf{(1)} just follows from definition of the quotient category. \textbf{(2)} is discussed in detail in \cite[page 12, (2.2.4)]{GR}. \smallskip

Now we show that \textbf{(3)} follows from \textbf{(1)} and \textbf{(2)}. Remark~\ref{rmk:main-example} implies that $M$ and $N$ are almost isomorphic if $\widetilde{\m}\otimes_R M \simeq \widetilde{\m}\otimes_R N$. \smallskip

Now suppose that there is an almost isomorphism $\varphi\colon M^a \to N^a$. It has a representative $f\colon \widetilde{\m}\otimes_R M \to N$ by \textbf{(2)}. Now \textbf{(1)} and $R$-flatness of $\widetilde{\m}$ imply that $\widetilde{\m} \otimes_R f\colon \widetilde{\m}\otimes_R \widetilde{\m} \otimes_R M\to \widetilde{\m}\otimes_R N$ is an isomorphism. Lemma~\ref{tilde-m} ensures that $\widetilde{\m}\otimes_R \widetilde{\m} \simeq \widetilde{\m}$, so $\widetilde{\m} \otimes_R f$ gives an isomorphism
\[
\widetilde{\m} \otimes_R f:\widetilde{\m} \otimes_R M\to \widetilde{\m}\otimes_R N. \qedhere
\]
\end{proof}

We now define the {\it functor of almost sections} 
\[
(-)_*\colon \textbf{Mod}_R^a \to \textbf{Mod}_R
\]
via the formula
\[
(M^a)_*\coloneqq \Hom_{\textbf{Mod}_R^a}(R^a, M^a)=\Hom_{R}(\widetilde{\m}, M)
\]
for any $R^a$-module $M^a$ with an $R$-module representative $M$. The construction is clearly functorial in $M^a$, so it defines the functor $(-)_*\colon \textbf{Mod}_R^a \to \textbf{Mod}_R$. \smallskip

The functor of almost sections will be the right adjoint to the almostification functor $(\text{-})^a$. Before we discuss why this is the case, we need to define the unit and counit transformations. \smallskip

We start with the unit of the adjunction. For any $R$-module $M$, there is a functorial morphism
\[
\eta_{M, *}\colon M \to \rm{Hom}_R(\widetilde{\m}, M)= M^a_*
\]
that can easily be seen to be an almost isomorphism. 

This allows us to define a functorial morphism
\[
\e_{N^a, *}\colon (N^a_*)^a \to N^a
\]
for any $R^a$-module $N^a$. Namely, the map $\eta_{N, *}\colon N \to N^a_*$ is an almost isomorphism, so we can invert it in the almost category and define 
\[
\e_{N^a, *}\coloneqq (\eta_{N, *}^a)^{-1}\colon (N^a_*)^a \to N^a
\]

Now we define another functor 
\[
(-)_!\colon \textbf{Mod}_R^a \to \textbf{Mod}_R
\] 
that will be a left adjoint to the almostification functor $(-)^a$. Namely, we put
\[
(M^a)_!\coloneqq (M^a)_*\otimes_R \widetilde{\m}\xleftarrow{\sim} M\otimes_R \widetilde{\m}
\]
for any $R^a$-module $M^a$ with an $R$-module representative $M$. This construction is clearly functorial in $M^a$, so it does define a functor. Similarly to the discussion above, for any $R$-module $M$, we define the transformation
\[
\e_{M, !}\colon (M^a)_!=\widetilde{\m}\otimes_R M \to M
\]
as the map induces by the natural morphism $\widetilde{\m} \to R$. Clearly, $\e_{M, !}$ is an almost isomorphism for any $M$. Therefore, this actually allows us to define the morphism 
\[
\eta_{N^a, !}\colon N^a \to (\widetilde{\m} \otimes_R N)^a\simeq (N_!^a)^a
\]
as $\eta_{N^a, !}=(\e_{N, !}^a)^{-1}$. We summarize the main properties of these functors in the following lemma:

\begin{lemma}\label{lemma:adjoint-almost} Let $R$ and $\m$ be  as above. Then \begin{enumerate}[label=\textbf{(\arabic*)}]
	\item\label{adjoint-almost-1} The functor $(-)_*$ is the right adjoint to $(-)^a$. In particular, it is left exact. 
	\item\label{adjoint-almost-2} The unit of the adjunction is equal to $\eta_{M, *}$, the counit of the adjunction is equal to $\e_{N^a, *}$. In particular, both are isomorphisms. 
	\item\label{adjoint-almost-3} The functor $(-)_!$ is the left adjoint to the localization functor $(-)^a$. 
    \item\label{adjoint-almost-4} The functor $(-)_!\colon \textbf{Mod}_R^a \to \textbf{Mod}_R$ is exact.
	\item\label{adjoint-almost-5} The unit of the adjunction is equal to $\eta_{N^a, !}$, the counit of the adjunction is equal to $\e_{M, !}$. In particular, both are almost isomorphisms. 
\end{enumerate}
\end{lemma}
\begin{proof}
This is explained in \cite[Proposition 2.2.13 and Proposition 2.2.21]{GR}.
\end{proof}

\begin{cor}\label{cor:limits-colimits-almost} Let $R$ and $\m$ be  as above. Then $(-)^a\colon \bf{Mod}_R\to \bf{Mod}_R^a$ commutes with limits and colimits. In particular, $\bf{Mod}_R^a$ is complete and cocomplete, and filtered colimits and (arbitrary) products are exact in $\bf{Mod}_R^a$.
\end{cor}
\begin{proof}
    The first claim follows from the fact that $(-)^a$ admits left and right adjoints. The second claim follows the first claim, exactness of $(-)^a$, and analogous exactness properties in $\bf{Mod}_R$. 
\end{proof}

The last thing we need to address in this section is how almost mathematics interacts with base change. We want to be able to talk about preservation of various properties of modules under a base change along a map $R \to S$. The issue here is to define the corresponding ideal $\m_S$ as in the definition of almost mathematics. It turns out that the most naive ideal $\m_S\coloneqq \m S$ does define an ideal of almost mathematics in $S$, but this is not entirely formal and crucially uses our choice of definition for an ideal of almost mathematics. \smallskip

More precisely, if one starts with a flat ideal $\m\subset R$, then the ideal $\m_S\subset S$ is not necessarily flat. However, we show that flatness of $\widetilde{\m}$ implies flatness of $\widetilde{\m_S}$. For this reason, it is essential to not impose a stronger condition on $\m$ to be $R$-flat in the foundations of almost mathematics.

\begin{lemma}\label{base-change} Let $f\colon R \to S$ be a ring homomorphism, and let $\m_S$ be the ideal $\m S\subset S$. Then we have the equality $\m_S^2=\m_S$ and the $S$-module $\widetilde{\m_S}\coloneqq \m_S \otimes_S \m_S$ is $S$-flat.
\end{lemma}
\begin{proof}
The equality $\m_S^2=\m_S$ follows from the analogous assumption on $\m$ and the construction of $\m_S$. Regarding the flatness issue, we claim that $\m_S \otimes_S \m_S \simeq (\m\otimes_R S)\otimes_S(\m \otimes_R S)$. That would certainly imply the desired flatness statement. To prove this claim, we look at the following short exact sequence
\[
0 \to \m \to R \to R/\m \to 0.
\]
We apply $ - \otimes_R S$ to get a short exact sequence
\[
0 \to \Tor_1^R(R/\m, S) \to \m\otimes_R S \to \m S \to 0.
\]
We observe that $\Tor_1^R(R/\m, S)$ is almost zero, so both $\Tor_1^R(R/\m, S) \otimes_S \m S$ and $\Tor_1^R(R/\m, S)\otimes_S (\m\otimes_R S)$ are zero modules due to Lemma \ref{almost-zero}. So we use functors $-\otimes_S (\m\otimes_R S)$ and $-\otimes_S \m S$ to obtain isomorphisms
\[
(\m\otimes_R S) \otimes_S (\m\otimes_R S) \simeq \m S\otimes_R (\m\otimes_R S) \simeq (\m S)\otimes_S (\m S).
\]
Thus we get the desired equality.
\end{proof}

\begin{lemma}\label{lemma:almost-zero-preserve} Let $f\colon R \to S$ be a ring homomorphism, and let $F \colon \bf{Mod}_R \to \bf{Mod}_S$ be an $R$-linear functor (resp. let $F \colon \bf{Mod}^{op}_R \to \bf{Mod}_S$ be an $R$-linear functor). Then $F$ sends almost zero $R$-modules to almost zero $S$-modules.
\end{lemma}
\begin{proof}
Suppose that $M$ is an almost zero $R$-module, so $\e M=0$ for any $\e \in \m$. Then $\e F(M)=0$ because $F$ is $R$-linear, so $F(M)$ is almost zero by Lemma~\ref{almost-zero}.
\end{proof}

\begin{cor}\label{cor:preserve-almost-iso} Let $f\colon R \to S$ be a ring homomorphism, and let $F \colon \bf{Mod}_R \to \bf{Mod}_S$ be a left or right exact $R$-linear functor (resp. let $F \colon \bf{Mod}^{op}_R \to \bf{Mod}_S$ be a left or right exact $R$-linear functor). Then $F$ preserves almost isomorphisms.
\end{cor}
\begin{proof}
We only show the case of a left exact functor $F\colon \bf{Mod}_R \to \bf{Mod}_S$, all other cases are analogous to this one. We choose any almost isomorphism $f\colon M' \to M''$ and wish to show that $F(f)$ is an almost isomorphism. For this, we consider the following exact sequences:
\[
0 \to K \to M' \to M \to 0,
\]
\[
0 \to M \to M'' \to Q \to 0.
\]
We know that $K$ and $Q$ are almost zero by our assumption on $f$. Now, the above short exact sequences induce the following exact sequences:
\[
0 \to F(K) \to F(M') \to F(M) \to \rm{R}^1F(K),
\]
\[
0 \to F(M) \to F(M'') \to F(Q).
\]
Lemma~\ref{lemma:almost-zero-preserve} guarantees that $F(K)$, $\rm{R}^1F(K)$, and $F(Q)$ are almost zero $S$-modules. Therefore, the morphisms $F(M') \to F(M)$ and $F(M) \to F(M'')$ are both almost isomorphisms. In particular, the composition $F(M') \to F(M'')$ is an almost isomorphism as well. 
\end{proof}

\subsection{Basic functors on categories of almost modules}

The category of almost modules admits certain natural functors induced from the category of $R$-modules. It has two versions of the Hom-functor and the tensor product functor. We summarize the properties of these functors in the following proposition:

\begin{prop}\label{many-functors} Let $R, \m$ be as above. Then
\begin{enumerate}[label=\textbf{(\arabic*)}]
\item\label{many-functors-1} We define {\it tensor product} functor $-\otimes_{R^a}-\colon \bf{Mod}_R^a\times \bf{Mod}_R^a \to \bf{Mod}_R^a$ as 
\[
(M^a, N^a) \mapsto (M^a_!\otimes_R N^a_!)^a.
\]
Then there is a natural transformation of functors 
\[
\begin{tikzcd}[column sep = 5em]
\bf{Mod}_R\times \bf{Mod}_R \arrow{r}{-\otimes_R -}\arrow{d}{(-)^a\times (-)^a} &  \bf{Mod}_R  \arrow{d}{(-)^a}\\
\bf{Mod}_R^a\times \bf{Mod}_R^a  \arrow[ru, Rightarrow, "\rho"] \arrow{r}{-\otimes_{R^a}-}& \bf{Mod}_R^a
\end{tikzcd}
\]
that makes the diagram $(2, 1)$-commutative. In particular, there is a functorial isomorphism $(M\otimes_R N)^a \simeq M^a\otimes_{R^a} N^a$ for any $M, N\in \bf{Mod}_R$. 
\item\label{many-functors-2} There is a functorial isomorphism 
\[
\rm{Hom}_{R^a}(M^a, N^a)\simeq \rm{Hom}_R(\widetilde{\m}\otimes M, N).
\]
for any $M,N \in \bf{Mod}_R$. In particular, there is a canonical structure of an $R$-module on the group $\rm{Hom}_{R^a}(M^a, N^a)$; thus defines the functor
\[
\rm{Hom}_{R^a}(-, -) \colon \bf{Mod}_{R^a}^{op}\times \bf{Mod}_{R^a} \to \bf{Mod}_R
\]
\item\label{many-functors-3} We define the functor $\rm{alHom}_{R^a}(-, -)\colon \bf{Mod}_{R^a}^{op}\times \bf{Mod}_{R^a} \to \bf{Mod}_{R^a}$ of  {\it almost homomorphisms} as 
\[
(M^a, N^a) \mapsto \rm{Hom}_{R^a}(M^a, N^a)^a.
\]
Then there is a natural transformation of functors 
\[
\begin{tikzcd}[column sep = 5em]
\bf{Mod}_R^{op}\times \bf{Mod}_R \arrow{r}{\rm{Hom}_R(-, -)}\arrow{d}{(-)^a\times (-)^a} &  \bf{Mod}_R  \arrow[dl, Rightarrow, "\rho"] \arrow{d}{(-)^a}\\
\bf{Mod}_{R^a}^{op}\times \bf{Mod}_{R^a}   \arrow{r}{\rm{alHom}_{R^a}(-, -)}& \bf{Mod}_{R^a}
\end{tikzcd}
\]
that makes the diagram $(2, 1)$-commutative. In particular, $\rm{alHom}_{R^a}(M^a, N^a) \cong^a \rm{Hom}_R(M, N)^a$ for any $M, N \in \bf{Mod}_R$.
\end{enumerate}
\end{prop}
\begin{proof}
\ref{many-functors-1}. We define 
\[
\rho_{M, N}\colon (M^a_!\otimes_R N^a_!)^a \to (M\otimes_R N)^a
\]
to be the morphism induced by 
\[
M^a_!\simeq \widetilde{\m}\otimes_R M \to M \text{ and } N^a_!\simeq \widetilde{\m}\otimes_R N \to N. 
\]
It is clear that $\rho_{M,N}$ is functorial in both variables, so it defines a natural transformation of functors $\rho$. We also need to check that $\rho_{M,N}$ is an isomorphism for any $M$ and $N$. This follows from the following two observations: $\rho_{M,N}$ is an isomorphism if and only if $\rho_{M,N}\otimes_R \widetilde{\m}$ is an isomorphism; and $\rho_{M,N}\otimes_R \widetilde{\m}$ is easily seen to be an isomorphism as $\widetilde{\m}\otimes_{R} \widetilde{\m} \to \widetilde{\m}$ is an isomorphism. \medskip

\ref{many-functors-2} is just a reformulation of Lemma~\ref{prop-almost}\ref{prop-almost-2}. \medskip

In order to show \ref{many-functors-3}, we need to define a functorial morphism 
\[
\rho_{M, N}\colon \rm{Hom}_R(M, N)^a \to \rm{alHom}_{R^a}(M^a, N^a).
\]
We start by using the functorial identification 
\[
\rm{alHom}_{R^a}(M^a, N^a) \cong^a \rm{Hom}_R(\widetilde{\m}\otimes M, N)^a
\]
from \ref{many-functors-2}. Namely, we define $\rho_{M, N}$ as the morphism $\rm{Hom}_R(M, N)^a \to \rm{Hom}_R(\widetilde{\m}\otimes M, N)^a$ induced by the map $\widetilde{\m}\otimes M \to M$. This is clearly functorial in both variables, so it defines the natural transformation $\rho$. \smallskip

We also need to check that $\rho_{M,N}$ is an isomorphism for any $M$ and $N$. This boils down to the fact that $\rm{Hom}_R(-, N)$ sends almost isomorphisms to almost isomorphisms. This, in turn, follows from Corollary~\ref{cor:preserve-almost-iso}.
\end{proof}

\begin{rmk} It is straightforward to check that if $N$ has a structure of an $S^a$-module for some $R$-algebra $S$, then the $R^a$-modules $M^a\otimes_{R^a} N^a$, $\rm{alHom}_{R^a}(M^a, N^a)$ have functorial-in-$M^a$ structures of $S^a$-modules. This implies that the functors  $-\otimes_{R^a} N^a$, $\rm{alHom}_{R^a}(-, N^a)$ naturally land in $\bf{Mod}_S^a$, i.e. define functors
\[
-\otimes_{R^a} N^a \colon \bf{Mod}_R^a \to \bf{Mod}_S^a, \text{ and } \rm{alHom}_{R^a}(-, N^a) \colon \bf{Mod}_R^{a, op} \to \bf{Mod}_S^a
\] 
Similarly, $\rm{Hom}_{R^a}(-, N^a)$ defines a functor $\bf{Mod}_{R}^a \to \bf{Mod}_{S}$.
\end{rmk}

The functor of almost homomorphisms is quite important, as it turns out to be the {\it inner Hom functor}, i.e. it is right adjoint to the tensor product. 

\begin{lemma}\label{inner-hom-strong} Let $f\colon R \to S$ be a ring homomorphism, and let $M^a$ be an $R^a$-module and $N^a, K^a$ be $S^a$-modules. Then there is a functorial $S$-linear isomorphism
\[
\rm{Hom}_{S^a}(M^a\otimes_{R^a} N^a, K^a) \simeq \rm{Hom}_{R^a}(M^a, \rm{alHom}_{S^a}(N^a, K^a)).
\]
\end{lemma}
\begin{proof}
This is a consequence of the usual $\otimes$-$\rm{Hom}$-adjunction, Proposition~\ref{many-functors}, and the fact that $\widetilde{\m}^{\otimes 2}\simeq \widetilde{\m}$. Indeed, we have the following sequence of functorial isomorphisms
\begin{align*}
\rm{Hom}_{S^a}(M^a\otimes_{R^a} N^a, K^a) & \simeq \rm{Hom}_{S}(\widetilde{\m}\otimes_R M\otimes_{R} N, K)\\
& \simeq \rm{Hom}_{S}((\widetilde{\m}\otimes_R M)\otimes_{R} (\widetilde{\m}\otimes_R N), K)\\
& \simeq \rm{Hom}_R(\widetilde{\m}\otimes_R M, \rm{Hom}_S(\widetilde{\m}\otimes_R N, K)) \\
&\simeq \rm{Hom}_{R^a}(M, \rm{alHom}_{S^a}(N^a, K^a)).
\end{align*}
The first isomomorphism follows from Proposition~\ref{many-functors}\ref{many-functors-1}, \ref{many-functors-2}, the second isomorphism follows from the observation $\widetilde{\m}^{\otimes 2}\simeq \widetilde{\m}$, the third isomorphism is just the classical $\otimes$-$\rm{Hom}$-adjunction, and the last isomorphism is a consequence of Proposition~\ref{many-functors}\ref{many-functors-2}, \ref{many-functors-3}.
\end{proof}

\begin{cor}\label{inner-hom-algebra} 
\begin{enumerate}[label=\textbf{(\arabic*)}]
\item\label{inner-hom-algebra-1}  Let $N$ be an $R^a$-module, then the functor $-\otimes_{R^a} N^a$ is left adjoint to the functor $\rm{alHom}_{R^a}(N^a, -)$.
\item\label{inner-hom-algebra-2} Let $R \to S$ be a ring homomorphism. Then the functor $-\otimes_{R^a} S^a\colon \bf{Mod}_R^a \to \bf{Mod}_S^a$ is left adjoint to the forgetful functor. 
\end{enumerate}
\end{cor}
\begin{proof}
Part~\ref{inner-hom-algebra-1} follows from Lemma~\ref{inner-hom-strong} by taking $S$ to be equal to $R$. Part~\ref{inner-hom-algebra-2} follows from Lemma~\ref{inner-hom-strong} by taking $N^a$ to be equal to $S^a$.
\end{proof}

We finish the section by introducing certain types of $R^a$-modules that will be used throughout the paper. 

\begin{defn}\label{defn:almost-flat}
\begin{itemize}
\item  An $R^a$-module $M^a$ is {\it flat} if the functor $M^a\otimes_{R^a} - \colon \bf{Mod}_R^a \to \bf{Mod}_R^a$ is exact.
\item  An $R^a$-module $M^a$ is {\it faithfully flat} if it is flat and $N^a\otimes_{R^a} M^a\simeq 0$ if and only if $N^a\simeq 0$.
\item An $R$-module $M$ is {\it almost flat} (resp. {\it almost faithfully flat}) if an $R^a$-module $M^a$ is {\it flat} (resp. {\it faithfully flat}) 
\item An $R^a$-module $I^a$ is {\it injective} if the functor $\rm{Hom}_{R^a}(-, I^a) \colon \bf{Mod}_R^{a, op} \to \bf{Mod}_R$ is exact.
\item An $R^a$-module $P^a$ is {\it almost projective} if the functor $\rm{alHom}_{R^a}(P^a, -) \colon \bf{Mod}_R^{a} \to \bf{Mod}_R^a$ is exact.
\end{itemize}
\end{defn}
 
\begin{lemma}\label{where-sends} The functor $(-)^a\colon \bf{Mod}_R \to \bf{Mod}_R^a$ sends flat (resp. faithfully flat, resp. injective, resp. projective) $R$-modules to flat (resp. faithfully flat, resp. injective, resp. almost projective) $R^a$-modules.
\end{lemma}
\begin{proof}
The case of flat modules is clear from Lemma~\ref{many-functors}$(\bf{1})$. Now suppose that $M$ is a faithfully flat $R$-module. Recall that $M\otimes_R -\colon \bf{Mod}_R \to \bf{Mod}_R$ is an exact and faithful functor. Therefore, if $M\otimes_R N$ is almost zero, it implies that so is $N$. Thus Lemma~\ref{many-functors}$(\bf{1})$ ensures that $M^a$ is almost faithfully flat. \smallskip

The case of injective modules follows from the fact that $(-)^a$ admits an exact left adjoint functor $(-)_!$. The case of projective modules is clear from the definition. 
\end{proof}

\begin{lemma}\label{lemma:where-sends-2} The functor $(-)_!\colon \bf{Mod}^a_R \to \bf{Mod}_R$ sends flat $R^a$-modules to flat $R$-modules.
\end{lemma}
\begin{proof}
    This follows from the formula $M^a_!\otimes_R N\simeq (M^a\otimes_{R^a} N^a)_!$ for any $R^a$-module $M^a$ and an $R$-module $N$.
\end{proof}

\begin{warning} If $M^a$ is a faithfully flat $R^a$-module, the $R$-module $M^a_!$ may not be faithfully flat. For instance, $R^a$ is a faithfully flat $R^a$-module, but $R^a_!=\widetilde{\m}$ is not a faithfully flat $R$-module. For example, $\widetilde{\m}\otimes_R R/\m \simeq 0$. 
\end{warning}

\begin{cor}\label{enough-almost-projective} Any bounded above complex $C^{\bullet, a}\in \bf{Comp}^-(R^a)$ admits a resolution $P^{\bullet, a} \to C^{\bullet, a}$ by a bounded above complex of almost projective modules.
\end{cor}
\begin{proof}
We consider the complex $C^{\bullet, a}_!\in \bf{Comp}^-(R)$; it admits a resolution by a complex of free modules $p\colon P^{\bullet} \to C^{\bullet, a}_!$. Now we apply $(-)^a$ to this morphism to obtain the maps
\[ 
P^{\bullet, a} \xrightarrow{p^a} (C^{\bullet, a}_!)^a \xleftarrow{\e} C^{\bullet, a}.
\]
The map $\e$ is an {\it isomorphism} in $\bf{Comp}(R^a)$ by Lemma~\ref{lemma:adjoint-almost}, and $p^a$ is a quasi-isomorphism by construction. Thus,  $\e^{-1} \circ p^a\colon P^{\bullet, a} \to C^{\bullet, a}$ is a quasi-isomorphism in $\bf{Comp}(R^a)$. Now note that each term of $P^{\bullet, a}$ is almost projective by Lemma~\ref{where-sends}.
\end{proof}

\subsection{Derived category of almost modules}\label{der-category-modules}

We define the derived category of almost modules in two different ways and show that these definitions coincide. Later we define certain derived functors on the derived category of almost modules. We pay some extra attention to show that the functors in this section are well-defined on unbounded derived categories. \smallskip

We start the section by introducing two different notions of the derived category of almost modules and then show that they are actually the same. 

\begin{defn} We define the {\it derived category of almost $R$-modules} as $\bf{D}(R^a)\coloneqq \bf{D}(\bf{Mod}_R^a)$. \smallskip 
\end{defn}

We define the bounded version of the derived category of almost $R$-modules $\bf{D}^*(R^a)$ for $*\in \{+, -, b\}$ as the full subcategory consisting of bounded below (resp. bounded above, resp. bounded) complexes. 

\begin{defn} We define the {\it almost derived category of $R$-modules} as the Verdier quotient $\bf{D}(R)^a\coloneqq \bf{D}(\bf{Mod}_R)/\bf{D}_{\Sigma_R}(\bf{Mod}_R)$. 
\end{defn}

We recall that $\Sigma_R$ is the Serre subcategory of $\bf{Mod}_R$ that consists of almost zero modules, and $\bf{D}_{\Sigma_R}(\bf{Mod}_R)$ is the full triangulated category of elements in $\bf{D}(\bf{Mod}_R)$ with almost zero cohomology modules. \smallskip 

We note that the functor $(-)^a\colon \bf{Mod}_R \to \bf{Mod}_R^a$ is exact and additive. Thus, it can be derived to the functor $(-)^a\colon \bf{D}(R) \to \bf{D}(R^a)$. Similarly, the functor $(-)_!\colon \bf{Mod}_R^a \to \bf{Mod}_R$ is additive and exact, so it can be derived to the functor $(-)_!\colon \bf{D}(R^a) \to \bf{D}(R)$. The standard argument shows that $(-)_!$ is a left adjoint functor to the functor $(-)^a$ since this already happens on the level of abelian categories. Now we also want to derive the functor $(-)_*\colon \bf{Mod}_{R}^a \to \bf{Mod}_R$. In order to do this on the level of unbounded derived categories, we need to show that $\bf{D}(R^a)$ has enough $K$-injective objects. \smallskip

\begin{defn} We say that a complex of $R^a$-module $I^{\bullet, a}$ is {\it K-injective} if $\rm{Hom}_{K(R^a)}(C^{\bullet, a}, I^{\bullet, a})=0$ for any acyclic complex $C^{\bullet, a}$ of $R^a$-modules.
\end{defn}

\begin{rmk} We remind the reader that $K(R^a)$ stands for the homotopy category of $R^a$-modules.
\end{rmk}

\begin{lemma}\label{K-inj}The functor $(-)^a\colon \bf{Comp}(R) \to \bf{Comp}(R^a)$ sends $K$-injective $R$-complexes to $K$-injective $R^a$-complexes. 
\end{lemma}
\begin{proof}
We note that $(-)^a$ admits an exact left adjoint $(-)_!$ thus \cite[\href{https://stacks.math.columbia.edu/tag/08BJ}{Tag 08BJ}]{stacks-project} ensures that $(-)^a$ preserves $K$-injective complexes. 
\end{proof}

\begin{cor}\label{enough-K-inj} Every object $M^{\bullet, a}\in \bf{Comp}(R^a)$ is quasi-isomorphic to a $K$-injective complex. 
\end{cor}
\begin{proof}
We know that the complex $M^{\bullet}\in \bf{Comp}(R)$ is quasi-isomorphic to a $K$-injective complex $I^{\bullet}$ by \cite[\href{https://stacks.math.columbia.edu/tag/090Y}{Tag 090Y}]{stacks-project} (or \cite[\href{https://stacks.math.columbia.edu/tag/079P}{Tag 079P}]{stacks-project}). Now we use Lemma~\ref{K-inj} to say that $I^{\bullet, a}$ is a $K$-injective complex that is quasi-isomorphic to $M^{\bullet, a}$.
\end{proof}

Now, as the first application of Corollary~\ref{enough-K-inj}, we define the functor $(-)_*\colon \bf{D}(R^a) \to \bf{D}(R)$ as the derived functor of $(-)_*\colon \bf{Mod}_{R}^a \to \bf{Mod}_R$. This functor exists by \cite[\href{https://stacks.math.columbia.edu/tag/070K}{Tag 070K}]{stacks-project}.
 \smallskip

\begin{lemma}\label{adjoint-derived} 
\begin{enumerate}[label=\textbf{(\arabic*)}]
	\item\label{adjoint-derived-1} The functors
    $
        \begin{tikzcd}
            \bf{D}(R)\arrow[r, swap, shift right=.75ex, "(-)^a"] & \bf{D}(R^a)\arrow[l, swap, shift right=.75ex, "(-)_!"] 
        \end{tikzcd}
    $ are adjoint. Moreover, the unit (resp. counit) morphism 
    \[
   (M^a)_! \to M \text{ (resp. } N \to (N_!)^a)
    \] is an almost isomorphism (resp. isomorphism) for any $M\in \bf{D}(R), N\in \bf{D}(R^a)$. In particular, the functor $(-)^a$ is essentially surjective.
    \item\label{adjoint-derived-2} The functors
    $
        \begin{tikzcd}
            \bf{D}(R)\arrow[r, shift left=.75ex, "(-)^a"] & \bf{D}(R^a)\arrow[l, shift left=.75ex, "(-)_*"] 
        \end{tikzcd}
    $ are adjoint. Moreover, the unit (resp. counit) morphism
    \[
    M \to (M^a)_* \text{ (resp. }(N_*)^a \to N)
    \] 
    is an almost isomorphism (resp. isomorphism) for any $M\in \bf{D}(R), N\in \bf{D}(R^a)$.
\end{enumerate}
\end{lemma}
\begin{proof}
We start the proof by showing \ref{adjoint-derived-1}. First, we note that the functors $(-)_!$ and $(-)^a$ are adjoint by the discussion above. Now we show that the cone of the counit map is always in $\bf{D}_{\Sigma_R}(R)$. As both functors $(-)^a$ and $(-)_!$ are exact on the level of abelian categories, it suffices to show the claim for $M\in \bf{Mod}_R^a$. But then the statement follows from Lemma~\ref{lemma:adjoint-almost}\ref{adjoint-almost-5}. The same argument shows that the unit map $N \to (N_!)^a$ is an isomorphism for any $N\in \bf{D}(R^a)$.\medskip

Now we go to \ref{adjoint-derived-2}. We define the functor $(-)_*\colon \bf{D}(R^a) \to \bf{D}(R)$ as the right derived functor of the left exact additive functor $(-)_*\colon \bf{Mod}_R^a \to \bf{Mod}_R$. This functor exists by \cite[\href{https://stacks.math.columbia.edu/tag/070K}{Tag 070K}]{stacks-project} and Corollary~\ref{enough-K-inj}. The functor $(-)_*$ is right adjoint to $(-)^a$ by \cite[\href{https://stacks.math.columbia.edu/tag/0DVC}{Tag 0DVC}]{stacks-project}. \smallskip

We check that the natural map $M \to (M^a)_*$ is an almost isomorphism for any $M\in \bf{D}(R)$. We choose some $K$-injective resolution $M \xrightarrow{\sim} I^{\bullet}$. Then Lemma~\ref{K-inj} guarantees that $M^a \to I^{\bullet, a}$ is a $K$-injective resolution of the complex $M^a$. The map $M \to (M^a)_*$ has a representative
\[
I^{\bullet} \to (I^{\bullet, a})_*.
\]
This map is an almost isomorphism of complexes by Lemma~\ref{lemma:adjoint-almost}\ref{adjoint-almost-2}. Thus, the map $M \to (M^a)_*$ is an almost isomorphism. A similar argument shows that the counit map $(N_*)^a \to N$ is an (almost) isomorphism for any $N\in \bf{D}(R^a)$. 
\end{proof}

\begin{thm}\label{derived-the-same} The functor $(-)^a\colon \bf{D}(R) \to \bf{D}(R^a)$ induces an equivalence of triangulated categories $(-)^a\colon \bf{D}(R)^a \to \bf{D}(R^a)$.
\end{thm}
\begin{proof}
We recall that the Verdier quotient is constructed as the localization of $\bf{D}(R)$ along the morphisms $f\colon C \to C'$ such that $\rm{cone}(f)\in \bf{D}_{\Sigma_R}(R)$. For instance, this is the definition of Verdier quotient in \cite[\href{https://stacks.math.columbia.edu/tag/05RI}{Tag 05RI}]{stacks-project}. Now we see that a morphism $f^a\colon C^a \to C'^a$ is invertible in $\bf{D}(R^a)$ if and only if $\rm{cone}(f)\in \bf{D}_{\Sigma_R}(R)$ by the definition of $\Sigma_R$ and the exactness of $(-)^a$. Moreover, $(-)^a$ admits a right adjoint such that $(-)^a\circ (-)_* \to \rm{Id}$ is an isomorphism of functors. Thus, we can apply \cite[Proposition 1.3]{GZ} to say that the induced functor $(-)^a\colon \bf{D}(R)^a \to \bf{D}(R^a)$ must be an equivalence. 
\end{proof}




\begin{rmk} Theorem~\ref{derived-the-same} shows that the two notions of the derived category of almost modules are the same. In what follows, we do not distinguish $\bf{D}(R^a)$ and $\bf{D}(R)^a$ anymore. 
\end{rmk}

\subsection{Basic functors on derived categories of almost modules}\label{basic-functors-derived-modules}

Now we can ``derive'' certain functors constructed in the previous section. We start by defining the derived versions of different Hom functors, after that we move to the case of the derived tensor product functor. 

\begin{defn} We define the {\it derived Hom} functor 
\[
\bf{R}\rm{Hom}_{R^a}(-, -) \colon \bf{D}(R^a)^{op}\times \bf{D}(R^a) \to \bf{D}(R)
\] 
as it is done in \cite[\href{https://stacks.math.columbia.edu/tag/0A5W}{Tag 0A5W}]{stacks-project} using the fact that $\bf{Comp}(R^a)$ has enough $K$-injective complexes. 

We define {\it Ext modules} via the following formula
\[
\rm{Ext}^i_{R^a}(M^a, N^a)\coloneqq \rm{H}^i(\bf{R}\rm{Hom}_{R^a}(M^a, N^a))\in \bf{Mod}_R
\] 
for $M^a, N^a\in\bf{Mod}_R^a$.
\end{defn}

Explicitly, for any $M^a, N^a \in \bf{D}(R^a)$, the construction of the complex $\bf{R}\rm{Hom}_{R^a}(M^a, N^a)$ goes as follows. We choose a representative $C^{\bullet, a} \to M^a$ and a $K$-injective resolution $N^a \to I^{\bullet, a}$. Then we set $\bf{R}\rm{Hom}_{R^a}(M^a, N^a)=\rm{Hom}^\bullet_{R^a}(C^{\bullet, a}, I^{\bullet, a})$. This construction is independent of the choices and is functional in both variables. We refer to \cite[\href{https://stacks.math.columbia.edu/tag/0A5W}{Tag 0A5W}]{stacks-project} for the details. 

\begin{rmk} We see that \cite[\href{https://stacks.math.columbia.edu/tag/0A64}{Tag 0A64}]{stacks-project} implies that there is a functorial isomorphism 
\[
\rm{H}^i\left(\bf{R}\rm{Hom}_{R^a}\left(M^a, N^a\right)\right) \simeq \rm{Hom}_{\bf{D}(R)^a}\left(M^a, N^a[i]\right).
\]
\end{rmk}

\begin{lemma}\label{derived-hom-alg}
\begin{enumerate}[label=\textbf{(\arabic*)}]
\item\label{derived-hom-alg-1} There are functorial isomorphisms 
\[
\rm{Hom}_{\bf{D}(R)^a}(M^a, N^a)\simeq \rm{Hom}_{\bf{D}(R)}(M^a_!, N) \text{ and } \bf{R}\rm{Hom}_{R^a}(M^a, N^a)\simeq \bf{R}\rm{Hom}_{R}(M^a_!, N)
\]
for any $M,N \in \bf{D}(R)$. 
\item\label{derived-hom-alg-2} For any chosen $M^a\in \bf{Mod}_R^a$, the functor $\bf{R}\rm{Hom}_{R^a}(M^a, -)\colon \bf{D}(R)^a \to \bf{D}(R)$ is isomorphic to the (right) derived functor of $\rm{Hom}_{R^a}(M^a, -)$.
\end{enumerate}
\end{lemma}
\begin{proof}
The first claim easily follows from the fact that $(-)^a$ is a right adjoint to the exact functor $(-)_!$. We leave the details to the reader. \smallskip

The second claim follows from \cite[\href{https://stacks.math.columbia.edu/tag/070K}{Tag 070K}]{stacks-project} and Corollary~\ref{enough-K-inj}.
\end{proof}

\begin{defn}\label{defn-almost-hom-sheaf} We define the {\it derived functor of almost homomorphisms} 
\[
\bf{R}\rm{alHom}_{R^a}(-, -) \colon \bf{D}(R^a)^{op}\times \bf{D}(R^a) \to \bf{D}(R^a)
\] 
as 
\[
\bf{R}\rm{alHom}_{R^a}(M^a, N^a)\coloneqq \bf{R}\rm{Hom}_{R^a}(M^a, N^a)^a=\bf{R}\rm{Hom}_R(M^a_!, N)^a. 
\]
We define the {\it almost Ext modules} as $R^a$-modules defined by \[\rm{alExt}^i_{R^a}(M^a, N^a)\coloneqq \rm{H}^i(\bf{R}\rm{alHom}_{R^a}(M^a, N^a))\] for $M^a, N^a\in\bf{Mod}_R^a$.
\end{defn}

\begin{defn} We define the {\it the complex of almost homomorphisms} $\rm{alHom}^{\bullet}_{R^a}(K^{\bullet, a}, L^{\bullet, a})$ for $K^{\bullet, a}, L^{\bullet, a}\in \bf{Comp}(R^a)$ as follows:
\[
\rm{alHom}_{R^a}^n(K^{\bullet, a}, L^{\bullet, a})\coloneqq \prod_{n=p+q} \rm{alHom}_{R^a}(K^{-q, a}, L^{p, a})
\]
with the differentials
\[
\rm{d}(f)=\rm{d}_{L^{\bullet, a}}\circ f - (-1)^n f\circ \rm{d}_{K^{\bullet, a}}.
\]
\end{defn}

\begin{lemma}\label{almost-proj-good-da} Let $P^{\bullet, a}$ be a bounded above complex of $R^a$-modules with almost projective cohomology modules. Suppose that $M^{\bullet, a}\to N^{\bullet, a}$ is an almost quasi-isomorphism of bounded below complex of $R^a$-modules. Then the natural morphism
\[
\rm{alHom}^{\bullet}_{R^a}(P^{\bullet, a}, M^{\bullet, a}) \to \rm{alHom}^{\bullet}_{R^a}(P^{\bullet, a}, N^{\bullet, a})
\]
is an almost quasi-isomorphism.
\end{lemma}
\begin{proof}
We note that as in the case of the usual Hom-complexes, there are convergent\footnote{Here we use that $P^{\bullet, a}$ is bounded above, $M^{\bullet, a}$ and $N^{\bullet, a}$ are bounded below} spectral sequences
\[
\rm{E}_1^{i,j}=\rm{H}^j\left(\rm{alHom}^{\bullet}_{R^a}(P^{-i, a}, M^{\bullet, a})\right) \Rightarrow \rm{H}^{i+j}\left(\rm{alHom}^{\bullet}_{R^a}\left(P^{\bullet, a}, M^{\bullet, a}\right)\right)
\] 
\[
\rm{E'}_1^{i,j}=\rm{H}^j\left(\rm{alHom}^{\bullet}_{R^a}(P^{-i, a}, N^{\bullet, a})\right) \Rightarrow \rm{H}^{i+j}\left(\rm{alHom}^{\bullet}_{R^a}\left(P^{\bullet, a}, N^{\bullet, a}\right)\right)
\] 
Moreover, there is a natural morphism of spectral sequences $\rm{E}_1^{i,j} \to \rm{E'}_1^{i,j}$. Thus, it suffices to show that the associated map on the first page is an almost isomorphism at each entry. For this, we use the fact that $\rm{alHom}_{R^a}(P^{-i,a}, -)$ is exact to rewrite the first page of this spectral sequence as 
\[
\rm{E}_1^{i,j}=\rm{alHom}_{R^a}\left(P^{-i, a}, \rm{H}^j(M^{\bullet, a})\right)
\]
and the same for $\rm{E'}_1^{i,j}$. So the question boils down to showing that the natural morphisms 
\[
\rm{alHom}_{R^a}\left(P^{-i, a}, \rm{H}^j(M^{\bullet, a})\right) \to \rm{alHom}_{R^a}\left(P^{-i, a}, \rm{H}^j(N^{\bullet, a})\right)
\] 
are almost isomorphisms. But this is clear as $M^{\bullet, a} \to N^{\bullet, a}$ is an almost quasi-isomorphism.
\end{proof}

\begin{lemma}\label{almost-project-for-alhom} Let $P_1^{\bullet, a} \to P_2^{\bullet, a}$ be an almost quasi-isomorphism of bounded above complexes with almost projective cohomology modules. Suppose that $M^{\bullet, a}$ is a bounded below complex of $R^a$-modules. Then the natural morphism
\[
\rm{alHom}^{\bullet}_{R^a}(P_2^{\bullet, a}, M^{\bullet, a}) \to \rm{alHom}^{\bullet}_{R^a}(P_1^{\bullet, a}, M^{\bullet, a})
\]
is an almost quasi-isomorphism.
\end{lemma}
\begin{proof}
We choose some injective resolution $M^{\bullet, a} \to I^{\bullet, a}$ of the bounded below complex $M^{\bullet, a}$. Then we have a commutative diagram
\[
\begin{tikzcd}
\rm{alHom}^{\bullet}_{R^a}(P_2^{\bullet, a}, M^{\bullet, a}) \arrow{d} \arrow{r} & \rm{alHom}^{\bullet}_{R^a}(P_1^{\bullet, a}, M^{\bullet, a}) \arrow{d} \\
\rm{alHom}^{\bullet}_{R^a}(P_2^{\bullet, a}, I^{\bullet, a}) \arrow{r} & \rm{alHom}^{\bullet}_{R^a}(P_1^{\bullet, a}, M^{\bullet, a}).
\end{tikzcd}
\]
The bottom horizontal arrow is an almost quasi-isomorphism by the standard categorical argument with injective resolutions. The vertical maps are almost quasi-isomorphism by Lemma~\ref{almost-proj-good-da}.
\end{proof}

\begin{prop}\label{derived-al-hom}
\begin{enumerate}[label=\textbf{(\arabic*)}]
\item\label{derived-al-hom-1} There is a natural transformation of functors   
\[
\begin{tikzcd}[column sep = 5em]
\bf{D}(R)^{op}\times \bf{D}(R) \arrow{r}{\bf{R}\rm{Hom}_{R}(-, -)} \arrow{d}{(-)^a\times (-)^a}& \bf{D}(R)  \arrow{d}{(-)^a}\arrow[dl, Rightarrow, "\rho"]\\
\bf{D}(R^a)^{op}\times \bf{D}(R^a)  \arrow{r}{\bf{R}\rm{alHom}_{R^a}(-, -)}& \bf{D}(R^a)
\end{tikzcd}
\]
that makes the diagram $(2, 1)$-commutative. In particular, \[
\bf{R}\rm{alHom}_{R^a}(M^a, N^a)\cong^a \bf{R}\rm{Hom}_R(M, N)^a
\]
for any $M,N\in \bf{D}(R)$. 
\item\label{derived-al-hom-2} For any chosen $M^a\in \bf{Mod}_R^a$, the functor $\bf{R}\rm{alHom}_{R^a}(M^a, -)\colon \bf{D}(R^a) \to \bf{D}(R^a)$ is isomorphic to the (right) derived functor of $\rm{alHom}_{R^a}(M^a, -)$.
\item\label{derived-al-hom-3} For any chosen $N^a\in \bf{Mod}_R^a$, the functor $\bf{R}\rm{alHom}_{R^a}(-, N^a)\colon \bf{D}^-(R^a)^{op} \to \bf{D}(R^a)$ is isomorphic to the (right) derived functor of $\rm{alHom}_{R^a}(-, N^a)$.
\end{enumerate}
\end{prop}
\begin{proof}
In order to show Part~\ref{derived-al-hom-1}, we construct functorial morphisms 
\[
\rho_{M, N}\colon \bf{R}\rm{Hom}_R(M, N)^a \to \bf{R}\rm{alHom}_{R^a}(M^a, N^a). 
\]
for any $M, N\in \bf{D}(R)$. We recall that there is a functorial identification 
\[
\bf{R}\rm{alHom}_{R^a}(M^a, N^a) \cong^a \bf{R}\rm{Hom}_R(M^a_!, N)^a\cong^a\bf{R}\rm{Hom}_R(\widetilde{\m}\otimes_R M, N)^a.
\]
So we define 
\[
\rho_{M, N}\colon \bf{R}\rm{Hom}_R(M, N)^a \to \bf{R}\rm{Hom}_R(\widetilde{\m}\otimes_R M, N)^a 
\]
as the morphism induced by the canonical map $\widetilde{\m}\otimes_R M \to M$. This is clearly functorial, so it defines the natural transformation of functors. The only thing we are left to show is that $\rho_{M, N}$ is an almost isomorphism for any $M, N\in \bf{D}(R)$. \smallskip

We recall that $\bf{R}\rm{Hom}_R(M, N)$ is isomorphic to $\rm{Hom}_R^{\bullet}(C^\bullet, I^\bullet)$ for any choice of a $K$-injective resolution of $N\xrightarrow{\sim} I^{\bullet}$ and any resolution $M\xrightarrow{\sim} C^{\bullet}$. Since $\widetilde{\m}\otimes_R C^{\bullet}$ is a resolution of $\widetilde{\m}\otimes_R M$ due to $R$-flatness of $\widetilde{\m}$, we reduce the question to showing that the natural map
\[
\alpha \colon \rm{Hom}_R^{\bullet}(C^\bullet, I^\bullet) \to \rm{Hom}_R^{\bullet}(\widetilde{\m}\otimes_R C^\bullet, I^\bullet)
\]
is an almost quasi-isomorphism of complexes. For this, it suffices to show that $\alpha$ is an isomorphism of complexes. Now note that the degree-$n$ part of $\alpha$ is the map
\[
\prod_{p+q=n} \rm{Hom}_R(C^{-q}, I^p) \to \prod_{p+q=n} \rm{Hom}_R(\widetilde{\m}\otimes_R C^{-q}, I^p).
\]
Since (infinite) products are exact in $\bf{Mod}_R^a$, and any (infinite) product of almost zero modules is almost zero, it suffices to show that each particular map $\rm{Hom}_R(C^{-q}, I^p) \to \rm{Hom}_R(\widetilde{\m}\otimes_R C^{-q}, I^p)$ is an almost isomorphism. This follows from Proposition~\ref{many-functors}\ref{derived-al-hom-3}. \medskip

Part~\ref{derived-al-hom-2} is similar to that of Proposition~\ref{derived-hom-alg}. \medskip

Part~\ref{derived-al-hom-3} is also similar to Part~\ref{derived-hom-alg-2} of Proposition~\ref{derived-hom-alg}, but there are some subtleties due to the fact that $\bf{Mod}_R^a$ does not have enough projective objects. We fix this issue by using \cite[\href{https://stacks.math.columbia.edu/tag/06XN}{Tag 06XN}]{stacks-project} instead of \cite[\href{https://stacks.math.columbia.edu/tag/070K}{Tag 070K}]{stacks-project}. We apply it to the subset $\mathcal P$ of bounded above complexes with almost projective terms. This result is indeed applicable in our situation due to Corollary~\ref{enough-almost-projective} and Lemma~\ref{almost-project-for-alhom}.  
\end{proof}

Now we deal with the case of the derived tensor product functor.

\begin{defn} We say that a complex of $R^a$-module $K^{\bullet, a}$ is {\it almost K-flat} if the naive tensor product complex $C^{\bullet, a}\otimes^{\bullet}_{R^a} K^{\bullet, a}$ is acyclic for any acyclic complex $C^{\bullet, a}$ of $R^a$-modules
\end{defn}


\begin{lemma}\label{almost-K-flat} The functor $(-)^a\colon \bf{Comp}(R) \to \bf{Comp}(R^a)$ sends $K$-flat $R$-complexes to almost $K$-flat $R^a$-complexes. 
\end{lemma}
\begin{proof}
Suppose that $C^{\bullet, a}$ is an acyclic complex of $R^a$-modules and $K^{\bullet}$ is a $K$-flat compelx. Then we see that
\[
C^{\bullet, a}\otimes^{\bullet}_{R^a} K^{\bullet, a} \cong^a (C^{\bullet}\otimes^{\bullet}_{R} K^{\bullet})^a \cong^a (\widetilde{\m}\otimes_R C^{\bullet}\otimes^{\bullet}_{R} K^{\bullet})^a \cong^a  \left(\left(\widetilde{\m}\otimes_R C^{\bullet}\right)\otimes^{\bullet}_{R} K^{\bullet}\right)^a.
\]
The latter complex is acyclic as $\widetilde{\m}\otimes C^{\bullet}$ is acyclic and $K^{\bullet}$ is $K$-flat.
\end{proof}

\begin{cor}\label{enough-K} Every object $M^{\bullet, a}\in \bf{Comp}(R^a)$ is quasi-isomorphic to an almost $K$-flat complex. 
\end{cor}
\begin{proof}
We know that the complex $M^{\bullet}\in \bf{Comp}(R)$ is quasi-isomorphic to a $K$-flat complex $K^{\bullet}$ by \cite[\href{https://stacks.math.columbia.edu/tag/06Y4}{Tag 06Y4}]{stacks-project}. Now we use Lemma~\ref{almost-K-flat} to say that $K^{\bullet, a}$ is an almost $K$-flat complex that is quasi-isomorphic to $M^{\bullet, a}$.
\end{proof}

\begin{defn} We define the {\it derived tensor product functor} 
\[
-\otimes^L_{R^a}-\colon \bf{D}(R)^a\times \bf{D}(R)^a \to \bf{D}(R)^a
\] by the rule 
$
(M^a, N^a) \mapsto (M_!\otimes^L_R N_!)^a
$
for any $M^a, N^a \in \bf{D}(R)^a$.
\end{defn}

\begin{prop}\label{derived-tensor-product}
\begin{enumerate}[label=\textbf{(\arabic*)}]
\item\label{derived-tensor-product-1} There is a natural transformation of functors 
\[
\begin{tikzcd}
\bf{D}(R)\times \bf{D}(R) \arrow{r}{-\otimes_R^L -} \arrow[d, swap, "(-)^a\times (-)^a"]& \bf{D}(R) \arrow{d}{(-)^a} \\
\bf{D}(R)^a\times \bf{D}(R)^a \arrow{r}{-\otimes_{R^a}^L -} \arrow[ru, Rightarrow, "\rho"]& \bf{D}(R)^a
\end{tikzcd}
\]
that makes the diagram $(2, 1)$-commutative. In particular, there is a functorial isomorphism $(M\otimes^L_R N)^a\simeq M^a\otimes^L_{R^a} N^a$ for any $M, N\in \bf{D}(R)$. 
\item\label{derived-tensor-product-2} For any chosen $M^a\in \bf{Mod}_R^a$, the functor $M^a\otimes^L_{R^a}-\colon \bf{D}(R)^a \to \bf{D}(R)^a$ is isomorphic to the (left) derived functor of $M^a\otimes_{R^a}-$.
\end{enumerate}
\end{prop}
\begin{proof}
The proof of Part~\ref{derived-tensor-product-1} is similar to that of Lemma~\ref{many-functors}\ref{many-functors-1}. We leave the details to the reader. \medskip

The proof of Part~\ref{derived-tensor-product-2} is similar to that of Proposition~\ref{derived-al-hom}\ref{derived-al-hom-2}. The claim follows by applying \cite[\href{https://stacks.math.columbia.edu/tag/06XN}{Tag 06XN}]{stacks-project} with $\mathcal P$ being the subset of almost $K$-flat complexes. This result is indeed applicable in our situation due to Corollary~\ref{enough-K} and the almost version of \cite[\href{https://stacks.math.columbia.edu/tag/064L}{Tag 064L}]{stacks-project}.  
\end{proof}

\begin{lemma}\label{o-hom-adj} Let $M^a, N^a, K^a \in \bf{D}(R)^a$, then we have a functorial isomorphism
\[
\bf{R}\rm{Hom}_{R^a}(M^a\otimes^L_{R^a} N^a, K^a) \simeq \bf{R}\rm{Hom}_{R^a}(M^a, \bf{R}\rm{alHom}_{R^a}(N^a, K^a)).
\]
In particular, the functors $
        \begin{tikzcd}
          \bf{R}\rm{alHom}_{R^a}(N^a, -)\colon   \bf{D}(R)^a\arrow[r, swap, shift right=.75ex] & \bf{D}(R)^a\colon -\otimes^L_{R^a} N^a \arrow[l, swap, shift right=.75ex] 
        \end{tikzcd}
    $ are adjoint. 
\end{lemma}
\begin{proof} The claim follows from the following sequence of canonical identifications:
\begin{align*} 
\bf{R}\rm{Hom}_{R^a}(M^a\otimes^L_{R^a} N^a, K^a) &\simeq \bf{R}\rm{Hom}_R((\widetilde{\m}\otimes_R M)\otimes^{L}_{R} (\widetilde{\m}\otimes_R N), K) & \rm{Lemma}~\ref{derived-hom-alg}\ref{derived-hom-alg-1}\\
& \simeq \bf{R}\rm{Hom}_R(\widetilde{\m}\otimes_R M, \bf{R}\rm{Hom}_{R}(\widetilde{\m}\otimes_R N, K)) & \text{\cite[\href{https://stacks.math.columbia.edu/tag/0A5W}{Tag 0A5W}]{stacks-project}} \\
& \simeq \bf{R}\rm{Hom}_{R^a}(M^{a}, \bf{R}\rm{Hom}_{R}(\widetilde{\m}\otimes_R N, K)^a) & \rm{Lemma}~\ref{derived-hom-alg}\ref{derived-hom-alg-1} \\
&\simeq \bf{R}\rm{Hom}_{R^a}(M^{a}, \bf{R}\rm{alHom}_{R^a}(N^{a}, K^a)).& \rm{Definition}~\ref{defn-almost-hom-sheaf}  
\end{align*}
\end{proof}

\begin{defn} Let $f\colon R \to S$ be a ring homomorphism. We define the {\it base change functor} 
\[
-\otimes^L_{R^a}S^a \colon \bf{D}(R)^a \to \bf{D}(S)^a
\] by the rule 
$
M^a \mapsto (M_!\otimes^L_R S)^a
$
for any $M^a\in \bf{D}(R)^a$.
\end{defn}

\begin{prop}\label{derived-base-change}
\begin{enumerate}[label=\textbf{(\arabic*)}]
\item\label{derived-base-change-1} There is a natural transformation of functors 
\[
\begin{tikzcd}
\bf{D}(R) \arrow{r}{-\otimes_R^L S} \arrow[d, swap, "(-)^a"]& \bf{D}(S) \arrow{d}{(-)^a} \\
\bf{D}(R)^a \arrow{r}{-\otimes_{R^a}^L S^a} \arrow[ru, Rightarrow, "\rho"]& \bf{D}(S)^a
\end{tikzcd}
\]
that makes the diagram $(2, 1)$-commutative. In particular, there is a functorial isomorphism $(M\otimes^L_R S)^a\simeq M^a\otimes^L_{R^a} S^a$ for any $M\in \bf{D}(R)$. 
\item\label{derived-base-change-2} The functor $-\otimes^L_{R^a}S^a\colon \bf{D}(R)^a \to \bf{D}(S)^a$ is isomorphic to the (left) derived functor of $-\otimes^L_{R^a}S^a$.
\end{enumerate}
\end{prop}
\begin{proof}
The proof is identical to Proposition~\ref{derived-tensor-product}.  
\end{proof}

\begin{lemma} Let $R\to S$ be a ring homomorphism, and let $M^a \in \bf{D}(R)^a$, $N^a\in \bf{D}(S)^a$. Then we have a functorial isomorphism
\[
\bf{R}\rm{Hom}_{S^a}(M^a\otimes^L_{R^a} S^a, N^a) \simeq \bf{R}\rm{Hom}_{R^a}(M^a, N^a).
\]
In particular, the functors $
        \begin{tikzcd}
          \rm{Forget}\colon   \bf{D}(S)^a\arrow[r, swap, shift right=.75ex] & \bf{D}(R)^a\colon -\otimes^L_{R^a} S^a \arrow[l, swap, shift right=.75ex] 
        \end{tikzcd}
    $ are adjoint. 
\end{lemma}
\begin{proof}
The proof is similar to that of Lemma~\ref{o-hom-adj}.
\end{proof}

\subsection{Almost finitely generated and almost finitely presented modules}\label{section-almost-finitely-presented}

In this section, we discuss the notions of almost finitely generated and almost finitely presented modules. Our discussion closely follows \cite{GR}. The main difference is that we avoid any use of ``uniform structures'' in our treatment; we think that it simplifies the exposition. We recall that we fixed some ``base'' ring $R$ with an ideal $\m$ such that $\m^2=\m$ and $\widetilde{\m}=\m\otimes_R \m$ is flat, and we always do almost mathematics with respect to this ideal. 

\begin{defn}\label{defn:almost-finitely-generated} An $R$-module $M$ is called {\it almost finitely generated}, if for any $\e \in \m$ there is an integer $n_{\e}$ and an $R$-homomorphism 
\[
R^{n_\e} \xr{f} M
\]
such that $\coker(f)$ is killed by $\e$.
\end{defn}

\begin{defn}\label{defn:almost-finitely-presented} An $R$-module $M$ is called {\it almost finitely presented}, if for any $\e, \dd \in \m$ there are integers $n_{\e, \dd}$, $m_{\e, \dd}$ and a complex 
\[
R^{m_{\e, \dd}} \xr{g} R^{n_{\e, \dd}} \xr{f} M
\]
such that $\coker(f)$ is killed by $\e$ and $\dd(\ker f) \subset \operatorname{Im} g$.
\end{defn}

\begin{rmk} Clearly, any almost finitely presented $R$-module is almost finitely generated.
\end{rmk}

\begin{rmk} A typical example of an almost finitely presented module that is not finitely generated is $M=\oplus_{n\geq 1} \O_C/p^{1/n}\O_C$ for an algebraically closed non-archimedean field $C$ of mixed characteristic $(0, p)$.
\end{rmk}

The next few lemmas discuss basic properties of almost finitely generated and almost finitely presented modules. For example, it is not entirely obvious that these notions transfer across almost isomorphisms. We show that this is actually the case, so these notions descend to $\mathbf{Mod}_R^a$. We also show that almost finitely generated and almost finitely presented modules have many good properties that are similar to those of usual finitely generated and finitely presented modules. \smallskip 

Our first main goal is to get alternative criteria for a module to be almost finitely generated (resp. almost finitely presented) and show that this notion descends to the category of almost modules.

\begin{lemma}\label{almost-finite} Let $M$ be an $R$-module, then $M$ is almost finitely generated if and only if for any finitely generated ideal $\m_0 \subset \m$ there a morphism $R^n \xr{f} M$ such that $\m_0(\coker f)=0$.
\end{lemma}
\begin{proof}
The ``if'' part is clear, so we only need to deal with the ``only if'' part. We choose a set of generators $(\e_0, \dots, \e_n)$ for an ideal $\m_0$. By assumption, we have $R$-morphisms
\[
f_{i}\colon R^{n_{\e_i}} \to M
\]
such that $\e_i(\coker f_i)=0$ for all $i$. Then the sum of these morphisms
\[
f\coloneqq \bigoplus_{i=1}^n f_{i}\colon R^{\sum n_{\e_i}} \to M
\]
defines a map such that $\m_0(\coker f)=0$. Since $\m_0$ was an arbitrary morphism, this finishes the proof.
\end{proof}

\begin{lemma}\label{claim} Let $M$ be an almost finitely presented $R$-module, and let $\varphi\colon R^n \to M$ be an $R$-homomorphism such that $\m_1(\coker \varphi)=0$ for some ideal $\m_1\subset \m$. Then for every finitely generated ideal $\m_0\subset \m_1\m$ there is morphism $\psi\colon R^m \to M$ such that 
\[
R^m \xr{\psi} R^n \xr{\varphi} M 
\]
is a three-term complex and $\m_0(\Ker\varphi)\subset \operatorname{Im}(\psi)$.
\end{lemma}
\begin{proof}

Since $M$ is almost finitely presented, for any two elements $\e_1, \e_2\in \m$, we can find a complex 
\[
R^{m_2} \xr{g} R^{m_1} \xr{f} M
\]
such that $\e_1(\coker f)=0$ and $\e_2(\ker f) \subset \operatorname{Im} g$. Now we choose an element $\delta \in \m_1$ and wish define morphisms
\[
\a\colon R^{m_1} \to R^n \text{ and } \beta\colon R^{n} \to R^{m_1}
\]
such that $\varphi \circ \a = \dd f$ and $f \circ \beta = \e_1 \varphi$. 
\[
\begin{tikzcd}
R^{m_2} \arrow[r, "g"] & R^{m_1} \arrow[r,"f"] \arrow[d, xshift=0.7ex, "{\a}"]  & M \\
& R^n \arrow[u, xshift=-0.7ex, "{\beta}"] \arrow[ru, "{\varphi}"] & 
\end{tikzcd}
\]

To achieve this goal, we define $\a$ and $\beta$ in the following way: we fix a basis $e_1, \dots, e_{m_1}$ of $R^{m_1}$, a basis $e'_1, \dots, e'_n$ of $R^n$, and then put $\a$ and $\beta$ to be the unique $R$-linear morphisms such that  
\[
\a(e_i)=y_i \in R^{n} \text{ for some $y_i$ such that } \varphi(y_i)=\delta f(e_i),
\]
\[
\beta(e'_j)=x_j \in R^{m_1} \text{ for some $x_j$ such that } f(x_j)=\e_1 \varphi(e'_j).
\]
It is clear that $\varphi \circ \a=\delta f$ and $f\circ \beta = \e_1\varphi$ as it holds on the basis elements. \smallskip

We define the morphism $\psi \colon R^n\oplus R^{m_2} \to R^n$
by the rule
\[
\psi(x,y)=\a\circ\beta(x)-(\e_1\dd) x+\a\circ g(y).
\]

Now we show that 
\[
\varphi \circ \psi =0 \text{ and } \e_1\e_2\dd \Ker \varphi \subset \operatorname{Im}\psi. \\
\]

{\it We start by showing that $\varphi \circ \psi=0$}: it suffices to prove that 
\[
(\a\circ g)(y) \in \Ker \varphi \text{ for } y\in R^{m_2}, \text{ and } (\a\circ\beta)(x)-(\e_1\dd)x \in \Ker \varphi \text{ for } x\in R^n
\]
We note that we have an equality
\[
(\varphi \circ \a\circ g)(y)=\dd (f \circ g)(y)=\dd 0 = 0,
\]
so $(\a\circ g)(y) \in \Ker (\varphi)$. We also have an equality
\begin{align*}
\left(\varphi \circ \left(\a\circ\beta - \e_1\dd\right)\right)\left(x\right)&=(\varphi \circ \a \circ \beta )(x)- \e_1\dd \varphi(x)\\ 
&= \dd (f\circ \beta )(x)-\e_1\dd \varphi(x) \\
&= \dd \e_1 \varphi(x) -\e_1\dd \varphi(x)\\
&= 0.
\end{align*}
this shows that $(\a\circ\beta)(x)-(\e_1\dd)x \in \Ker(\varphi)$ as well. \smallskip

{\it We show that $(\e_1\e_2\dd)\Ker \varphi \subset \operatorname{Im}(\psi)$}: we observe that for any $x\in \Ker\varphi$ we have $\beta(x)\subset \Ker f $ as $f\circ \beta = \e_1\varphi$. This implies that $\e_2\beta(x)\in \operatorname{Im} g $ since $\e_2\Ker f \subset \Imm g$. Thus, there is $y\in R^{m_2}$ such that $g(y)=\e_2 \beta(x)$, so $(\a\circ g)(y)=\e_2 \a\circ \beta(x)$. This shows that 
\begin{align*}
\psi(-\e_2 x, y)=-\e_2(\a\circ \beta)(x)+\e_1\e_2\dd x +(\a\circ g)(y)=\\ 
-\e_2(\a\circ \beta)(x)+\e_1\e_2\dd x + \e_2(\a\circ \beta)(x)= \e_1\e_2\dd x
\end{align*}
We conclude that $\e_1\e_2\dd x \in \operatorname{Im}(\psi)$ for any $x\in \Ker(\varphi)$. \smallskip

Finally, we recall that $\m_0$ is a finitely generated ideal, and that $\m_0 \subset \m_1\m=\m_1\m^2\subset \m_1$. This means that we can find a finite set $I$, and a finite set of elements $\e_{i,1}, \e_{i,2}\in \m, \dd_i \in \m_1$ such that $\m_0$ is contained in the ideal $J:=(\e_{i,1}\e_{i,2}\dd_{i})_{i\in I}$ (the ideal generated by all the products $\e_{i,1}\e_{i,2}\dd_{i}$). The previous discussion implies that for each $i\in I$, we have a map $\psi_i:R^{k_i} \to R^n$ such that $\varphi\circ \psi_i=0$ and $(\e_{i,1}\e_{i,2}\dd_{i})(\Ker\varphi) \subset \operatorname{Im}\psi_i$. By passing to the homomorphism 
\[
\psi\coloneqq \bigoplus_{i\in I} \psi_i\colon R^{\sum k_i} \to R^n,
\]
we get a map $\psi$ such that $\varphi \circ \psi=0$ and $\m_0(\Ker \varphi) \subset \operatorname{Im}(\psi)$. Therefore, $\psi$ does the job.
\end{proof}

\begin{lemma}\label{almost-finitely-presented} Let $M$ be an $R$-module. Then the following conditions are equivalent:
\begin{enumerate}[label=\textbf{(\arabic*)}]
\item\label{almost-finitely-presented-1} The $R$-module $M$ is almost finitely presented.
\item\label{almost-finitely-presented-2} For any finitely generated ideal $\m_0 \subset \m$ there exist a finitely presented $R$-module $N$ and a homomorphism $f\colon N \to M$ such that $\m_0(\ker f)=0$ and $\m_0(\coker f)=0$.
\item\label{almost-finitely-presented-3} For any finitely generated ideal $\m_0 \subset \m$ there exist integers $n, m$ and a three-term complex 
\[
R^m \xr{g} R^n \xr{f} M
\] 
such that $\m_0(\coker f)=0$ and $\m_0(\Ker f)\subset \operatorname{Im}g$.
\end{enumerate}
\end{lemma}
\begin{proof}
It is clear that \ref{almost-finitely-presented-3} implies both \ref{almost-finitely-presented-1} and \ref{almost-finitely-presented-2}. \smallskip

We show that \ref{almost-finitely-presented-1} implies \ref{almost-finitely-presented-3}. Since $M$ is an almost finitely generated $R$-module, Lemma~\ref{almost-finite} guarantees that, for any finitely generated ideal $\m' \subset \m$, there is a morphism $R^n \xr{f} M $ such that $\m'(\coker f)=0$. \smallskip

We know that $\m_0 \subset \m=\m^2$; this easily implies that there is a finitely generated ideal $\m_1\subset \m$ such that $\m_0\subset \m_1\m \subset \m_1$. So, using $\m'=\m_1$, we can find a homomorphism $R^n \xr{\varphi} M$ such that $\m_1(\coker \varphi)=0$. Lemma \ref{claim} claims that we can also find a homomorphism $\psi\colon R^m \to R^n$ such that 
\[
R^m \xr{\psi} R^n \xr{\varphi} M
\]
is a three-term complex and $\m_0(\ker \varphi) \subset \operatorname{Im}\psi$. Since $\m_0 \subset \m_1$ and $\m_1(\coker \varphi)=0$, we get that $\m_0(\coker \varphi)=0$ as well. This finishes the proof since $\m_0$ was an arbitrary finitely generated sub-ideal of $\m$. \\

Now we show that \ref{almost-finitely-presented-2} implies \ref{almost-finitely-presented-3}. We pick an arbitrary finitely generated ideal $\m_0\subset \m$, and we try to find a three-term complex 
\[
R^m \xr{g} R^n \xr{f} M
\]
such that $\m_0(\coker f)=0$ and $\m_0(\ker f)\subset \Imm(g)$. To achieve this, we use the assumption in \ref{almost-finitely-presented-2} to find a morphism $h\colon N \to M$ such that $N$ is a finitely presented $R$-module, $\m_0(\coker h)=0$, and $\m_0(\ker h)=0$. Since $N$ is finitely presented we can find a short exact sequence
\[
R^m\xr{g} R^n \xr{f'} N \to 0
\]
It is straightforward to see that a three-term complex 
\[
R^m \xr{g} R^n \xr{f\coloneqq h\circ f'} M
\]
satisfies the condition that $\m_0(\coker f)=0$ and $\m_0(\ker f)\subset \Imm(g)$. 
\end{proof}

\begin{lemma}\label{close-almost} Let $M$ be an $R$-module, and suppose that for any finitely generated ideal $\m_0 \subset \m$ there exists a morphism $f\colon N \to M$ such that $\m_0(\ker f)=0$, $\m_0(\coker f)=0$ and $N$ is almost finitely generated (resp. almost finitely presented). Then $M$ is also almost finitely generated (resp. almost finitely presented).
\end{lemma}
\begin{proof}
We give a proof only in the almost finitely presented case; the other case is easier. We pick an arbitrary finitely generated ideal $\m_0 \subset \m$ and another finitely generated ideal $\m_1\subset \m$ such that $\m_0 \subset \m_1^2$. Then we use the assumption to get a morphism 
\[
f\colon N \to M
\]
such that $\m_1(\Ker f)=0, \m_1(\coker f)=0$ and $N$ is an almost finitely presented $R$-module. Lemma~\ref{almost-finitely-presented} guarantees that there is a three-term complex
\[
R^m \xr{h} R^n \xr{g} N
\]
such that $\m_1(\coker g)=0$ and $\m_1(\Ker g) \subset \Imm h$. Then we can consider a three-term complex
\[
R^m \xr{h} R^n \xr{f'\coloneqq f\circ g} M,
\]
it is easily seen that $\m_1^2(\coker f')=0$ and $\m_1^2(\ker f') \subset \Imm(h)$. Since $\m_0 \subset \m_1^2$ we conclude that $\m_0(\coker f')=0$ and $\m_0(\ker f') \subset \Imm(h)$. This shows that $M$ is almost finitely presented.
\end{proof}

\begin{lemma}\label{lemma:afp-compact-1} Let $M$ be an $R$-module, and $\{N_i\}_{i\in I}$ is a filtered diagram of $R$-modules. Then
    \begin{enumerate}[label=\textbf{(\arabic*)}]
        \item\label{lemma:afp-compact-1-1} The natural morphism 
        \[
        \gamma^0_M\colon \colim_I \rm{Hom}_R(M, N_i) \to \rm{Hom}_R(M, \colim_I N_i)
        \]
        is almost injective for an almost finitely generated $M$;
        \item\label{lemma:afp-compact-1-2} The natural morphism 
        \[
        \gamma^0_M \colon \colim_I \rm{Hom}_R(M, N_i) \to \rm{Hom}_R(M, \colim_I N_i)
        \]
        is an almost isomorphism and 
        \[
        \gamma^1_M \colon \colim \rm{Ext}^1_R(M, N_i) \to \rm{Ext}^1_R(M, \colim N_i)
        \]
        is almost injective for an almost finitely presented $M$.
    \end{enumerate}
\end{lemma}
\begin{proof}
    We give a proof for an almost finitely presented $M$; the case of an almost finitely generated $M$ is similar. \smallskip
    
    {\it Step $1$: The case of a finitely presented $M$.} In this case, $\gamma^0_M$ is an isomorphism and $\gamma^1_M$ is injective due to \cite[\href{https://stacks.math.columbia.edu/tag/064T}{Tag 064T}]{stacks-project} and \cite[\href{https://stacks.math.columbia.edu/tag/0G8W}{Tag 0G8W}]{stacks-project}. \smallskip
    
    {\it Step $2$: General case.} We fix a finitely generated ideal $\m_0\subset \m$. Since $\m_0 \subset \m=\m^4$, there is a finitely generated ideal $\m_1$ such that $\m_0 \subset \m_1^4$. Now we use Lemma~\ref{almost-finitely-presented}\ref{almost-finitely-presented-2} to find a finitely presented module $M'$ and a morphism $f\colon M' \to M$ such that $\ker(f)$ and $\coker(f)$ are annihilated by $\m_1$. We denote the image of $f$ by $M''$ and consider the short exact sequences
\[
0 \to K \to M' \to M'' \to 0 \ ,
\]
\[
0 \to M'' \to M \to Q\to 0
\]
with $K$ and $Q$ being annihilated by $\m_1$. After applying the functors $\colim_I \rm{Hom}_R(-, N_i)$ and $\rm{Hom}_R(-,\colim_I N_i)$ and considering the associated long exact sequences, we see that 
\[
b_i\colon \colim_I \rm{Ext}^i_R(M, N_i) \to \colim_I \rm{Ext}^i_R(M', N_i)
\]
and
\[
c_i\colon \rm{Ext}^i_R(M, \colim_I N_i) \to \rm{Ext}^i_R(M', \colim_I N_i)
\]
have kernels and cokernels annihilated by $\m^2_1$ for any $i\geq 0$. Now we consider a commutative diagram
\[
\begin{tikzcd}
    \colim_I \rm{Ext}^i_R(M', N_i) \arrow{r}{\gamma^i_{M'}} & \rm{Ext}^i_R(M', \colim_I N_i)\\
    \colim_I \rm{Ext}^i_R(M, N_i) \arrow{u}{b_i} \arrow{r}{\gamma^i_M}& \rm{Ext}^i_R(M, \colim_I N_i) \arrow{u}{c_i}
\end{tikzcd}
\]
By Step~$1$, we know that $\gamma^i_{M'}$ is an isomorphism for $i=0$ and injective for $i=1$. Moreover, we know that $b_i$ and $c_i$ have kernels and cokernels annihilated by $\m_1^2$. Then it is easy to see that $\coker(\gamma_M^0)$, $\ker(\gamma_M^0)$, and $\ker(\gamma_M^1)$ are annihilated by $\m_1^4$. In particular, they are annihilated by $\m_0\subset \m_1^4$. Since $\m_0$ was arbitrary finitely generated sub-ideal $\m_0\subset \m$, we conclude that $\gamma^0_M$ is an almost isomorphism and $\gamma^1_M$ is almost injective. 
\end{proof}

\begin{lemma}\label{lemma:afp-compact-2} Let $M$ be an $R$-module.
    \begin{enumerate}[label=\textbf{(\arabic*)}]
        \item\label{lemma:afp-compact-2-1} If, for any filtered diagram of $R$-modules $\{N_i\}_{i\in I}$, the natural morphism 
        \[
        \colim_I \rm{Hom}_R(M, N_i) \to \rm{Hom}_R(M, \colim_I N_i)
        \]
        is almost injective, then $M$ is almost finitely generated.
        \item\label{lemma:afp-compact-2-2} If, for any filtered system of $R$-modules $\{N_i\}$, the natural morphism
        \[
        \colim_I \rm{Hom}_R(M, N_i) \to \rm{Hom}_R(M, \colim_I N_i)
        \]
        is an almost isomorphism, then $M$ is almost finitely presented.
    \end{enumerate}
\end{lemma}
\begin{proof}
    \ref{lemma:afp-compact-2-1}: Note that $M \simeq \colim_I M_i$ is a filtered colimit of its finitely generated {\it submodules}. Therefore, we see that 
    \[
    \colim_I \rm{Hom}_R(M, M/M_i) \simeq^a \rm{Hom}_R\left(M, \colim_I (M/M_i)\right) \simeq 0.
    \]
    Consider an element $\alpha$ of $\colim_I \rm{Hom}_R(M, M/M_i)$ that has a representative the quotient morphism $M \to M/M_i$  (for some choice of $i\in I$). Then, for every $\e\in \m$, $\e \alpha =0$ in $\colim_I \rm{Hom}_R(M, M/M_i)$. Explicitly, this means that there is $j\geq i$ such that $\e M \subset M_j$. Now we choose a surjection $R^{n_j} \to M_j$ to see that the composition $f \colon R^{n_j} \to M$ gives a map with $\e(\coker f)=0$. Now note that this property is preserved by replacing $j$ with any $j'>j$. Therefore, for any $\m_0=(\e_1, \dots, \e_n)$, we can find a finitely generated submodule $M_i \subset M$ such that $\m_0M \subset M_i$. Therefore, $M$ is almost finitely generated. \smallskip
    
    \ref{lemma:afp-compact-2-2}: Fix any finitely generated sub-ideal $\m_0 = (\e_1, \dots, \e_n) \subset \m$. We use \cite[\href{https://stacks.math.columbia.edu/tag/00HA}{Tag 00HA}]{stacks-project} to write $M\simeq \colim_{\Lambda} M_\lambda$ as a filtered colimit of {\it finitely presented} $R$-modules. By assumption, the natural morphism
    \[
    \colim_{\Lambda} \rm{Hom}_R(M, M_\lambda) \to \rm{Hom}_R(M, \colim_{\Lambda} M_\lambda)=\rm{Hom}_R(M, M)
    \]
    is an almost isomorphism. In particular, $\e_i \rm{Id}_M$ is in the image of this map for every $i=1,\dots, n$. This means that, for every $\e_i$, there is $\lambda_i\in \Lambda$ and a morphism $g_i \colon M \to M_{\lambda_i}$ such that the composition 
    \[
    f_{\lambda_i} \circ g_i=\e_i \rm{Id}_M,
    \]
    where $f_{\lambda_i}\colon M_{\lambda_i} \to M$ is the natural morphism to the colimit. Note that the existence of such $g_i$ is preserved by replacing $\lambda_i$ with any $\lambda'_i \geq \lambda_i$. Therefore, using that $\{M_\lambda\}$ is a filtered diagram, we can find an index $\lambda$ with maps
    \[
    g_i \colon M \to M_\lambda
    \]
    such that $f_\lambda \circ g_i = \e_i \rm{Id}_M$. We consider the morphism 
    \[
    F_i\coloneqq g_i\circ f_\lambda - \e_i \rm{Id}_{M_\lambda} \colon M_\lambda \to M_\lambda.
    \]
    We note that $\rm{Im}(F_i) \subset \ker(f_\lambda)$ because
    \[
    f_\lambda \circ g_i \circ f_\lambda - f_\lambda \e_i \rm{Id}_{M_i} = \e_i f_\lambda - \e_i f_\lambda = 0.
    \]
    We also have that $\e_i \ker(f_\lambda) \subset \rm{Im}(F_i)$ because $F_i|_{\ker(f_\lambda)} = \e_i \rm{Id}$. Therefore, $\sum_i \rm{Im}(F_i)$ is a finite $R$-module such that 
    \[
    \m_0(\ker f_\lambda) \subset \sum_i \rm{Im}(F_i) \subset \ker(f_\lambda).
    \]
    Therefore, $f\colon M'\coloneqq M_\lambda/(\sum_i \rm{Im}(F_i)) \to M$ is a morphism such that its source $M'$ is finitely presented, $\m_0(\ker f)=0$, and $\m_0(\coker f)=0$. Since $\m_0 \subset \m$ was an arbitrary finitely generated sub-ideal, we conclude that $M$ is almost finitely presented. 
\end{proof}

\begin{cor}\label{cor:compact-criterion} Let $M$ be an $R$-module. Then
    \begin{enumerate}[label=\textbf{(\arabic*)}]
        \item $M$ is almost finitely generated if and only if, for every filtered diagram $\{N_i^a\}_{i\in I}$ of $R^a$-modules, the natural morphism
        \[
        \colim_I \rm{alHom}_R(M^a, N^a_i) \to \rm{alHom}_R(M^a, \colim_I N^a_i)
        \]
        is injective in $\bf{Mod}_R^a$;
        \item $M$ is almost finitely presented if and only if, for every filtered diagram $\{N_i^a\}_{i\in I}$ of $R^a$-modules, the natural morphism
        \[
        \colim_I \rm{alHom}_R(M^a, N^a_i) \to \rm{alHom}_R(M^a, \colim_I N^a_i)
        \]
        is an isomorphism in $\bf{Mod}_R^a$;
    \end{enumerate}
\end{cor}
\begin{proof}
    It formally follows from Lemma~\ref{lemma:afp-compact-1}, Lemma~\ref{lemma:afp-compact-2}, Proposition~\ref{many-functors}\ref{many-functors-3}, and Corollary~\ref{cor:limits-colimits-almost}.
\end{proof}

\begin{Cor}\label{almost-almost-finitely-presented} Let $M$ and $N$ be two almost isomorphic $R$-modules (see Definition \ref{almost-isomorphic}). Then $M$ is almost finitely generated (resp. almost finitely presented) if and only if so is $N$.
\end{Cor}
\begin{proof}
    Corollary~\ref{cor:compact-criterion} implies that $M$ is almost finitely generated (resp. almost finitely presented) if and only if $M^a_!$ is. Since $M^a_! \simeq N^a_!$, we get the desired result. 
\end{proof}

\begin{cor}\label{cor:almost-iso-almost-equi} Let $R\to S$ be an almost isomorphism of rings. Then the forgetful functor $\bf{Mod}_{S^a}^* \to \bf{Mod}_{R^a}^*$ is an equivalence for $*\in \{``\text{ ''}, \rm{aft}, \rm{afp}\}$.
\end{cor}
\begin{proof}
    Corollary~\ref{cor:compact-criterion} ensures that it suffices to prove the claim for $*=``\text{ ''}$ as the property of being almost finitely generated (resp.\,almost finitely presented) depends only on the category $\bf{Mod}_{R^a}$ and not on the ring $R$ itself.\smallskip

    Corollary~\ref{inner-hom-algebra}~\ref{inner-hom-algebra-2} guarantees that the forgetful functor admits a right adjoint $-\otimes_{R^a} S^a \colon \bf{Mod}_R^a \to \bf{Mod}_S^a$. Therefore, it suffices to show that the natural morphisms
    \[
    M^a\to M^a\otimes_{R^a} S^a
    \]
    and
    \[
    N^a\otimes_{R^a} S^a \to N^a
    \]
    are isomorphisms for any $M\in \bf{Mod}_R^a$ and $N\in\bf{Mod}_S^a$. This is obvious from the fact that $R^a \to S^a$ is an isomorphism of $R^a$-modules. 
\end{proof}

\begin{defn} We say that an $R^a$-module $M^a\in \mathbf{Mod}_R^a$ is {\it almost finitely generated} (resp. {\it almost finitely presented}) if its representative $M\in \mathbf{Mod}_R$ is almost finitely generated (resp. almost finitely presented). This definition does not depend on the choice of a 
 representative due to Lemma~\ref{almost-almost-finitely-presented}
\end{defn}

We now want to establish certain good properties of almost finitely presented modules in short exact sequences. This will be crucial later to develop a good theory of almost coherent modules.

\begin{lemma}\label{main-almost} Let $0 \to M' \xr{\varphi} M \xr{\psi} M'' \to 0$ be an exact sequence of $R$-modules, then
\begin{enumerate}[label=\textbf{(\arabic*)}]
\itemsep0.5ex
\item\label{main-almost-1} If $M$ is almost finitely generated, then so is $M''$.
\item\label{main-almost-2} If $M'$ and $M''$ are almost finitely generated (resp. finitely presented), then so is $M$.
\item\label{main-almost-3} If $M$ is almost finitely generated and $M''$ is almost finitely presented, then $M'$ is almost finitely generated.
\item\label{main-almost-4} If $M$ is almost finitely presented and $M'$ is almost finitely generated, then $M''$ is almost finitely presented.
\end{enumerate}
\end{lemma}
\begin{proof}
This can be easily deduced from Lemma~\ref{lemma:afp-compact-1} and Lemma~\ref{lemma:afp-compact-2} via the five lemma (or diagram chase). We only note that the $\rm{Ext}^1$ part of Lemma~\ref{lemma:afp-compact-1}~\ref{lemma:afp-compact-1-2} is crucial to make the argument work. 
\end{proof}

\begin{Cor}\label{main-almost-almost} Let $0 \to M'^a \xr{\varphi} M^a \xr{\psi} M''^a \to 0$ be an exact sequence of  $R^a$-modules. Then all the conclusions of Lemma \ref{main-almost} still hold.
\end{Cor}
\begin{proof}
We use Lemma~\ref{lemma:adjoint-almost}\ref{adjoint-almost-4},\ref{adjoint-almost-5} to see that the sequence 
\[
0 \to (M'^a)_! \xr{\varphi_!} (M^a)_! \xr{\psi_!} (M''^a)_! \to 0
\]
is exact and almost isomorphic to the original one. Moreover, Corollary~\ref{almost-almost-finitely-presented} says that each of those modules $N^a_!$ is almost finitely generated (resp. almost finitely presented) if and only if so is the corresponding $N^a$. Thus, the problem is reduced to Lemma~\ref{main-almost}.
\end{proof}

\begin{lemma}\label{tensor-of-fp} Let $M^a ,N^a$ be two almost finitely generated (resp. almost finitely presented) $R^a$-modules, then so is $M^a\otimes_{R^a} N^a$. Similarly, $M\otimes_R N$ is almost finitely generated (resp. almost finitely presented) for any almost finitely generated (resp. almost finitely presented) $R$-modules $M$ and $N$.
\end{lemma}
\begin{proof}
We show the claim only in the case of almost finitely presented modules; the case of almost finitely generated modules is significantly easier. Moreover, we use Proposition~\ref{many-functors}\ref{many-functors-1} to reduce the question to showing that the tensor product of two almost finitely presented $R$-modules is almost finitely presented. \smallskip

{\it Step 1. The case of finitely presented modules}: If both $M$ and $N$ are finitely presented, then this is a standard fact proven in \cite[II, \textsection 3.6, Proposition 6]{Bou}. \smallskip

{\it Step 2. The case of $M$ being finitely presented}: Now we deal with the case of a finitely presented $R$-module $M$ and an almost finitely presented $N$. We fix a finitely generated ideal $\m_0\subset \m$ and a finitely generated ideal $\m_1$ such that $\m_0\subset \m_1^2$. Now we use Lemma~\ref{almost-finitely-presented}\ref{almost-finitely-presented-2} to find a finitely presented module $N'$ and a morphism $f\colon N' \to N$ such that $\ker(f)$ and $\coker(f)$ are annihilated by $\m_0$. We denote the image of $f$ by $N''$ and consider the short exact sequences
\[
0 \to K \to N' \to N'' \to 0 \ ,
\]
\[
0 \to N'' \to N \to Q\to 0
\]
with $K$ and $Q$ being annihilated by $\m_0$. After applying the functor $M\otimes_R -$, we get the following exact sequences:
\[
M \otimes_R K \to M\otimes_R N' \to M\otimes_R N'' \to 0 \ ,
\]
\[
\rm{Tor}_1^R(M, Q) \to M \otimes_R N'' \to M\otimes_R N \to M\otimes_R Q \to 0.
\]
We note that $M \otimes_R K, \rm{Tor}_1^R(M, Q)$, and $M\otimes_R Q$ are annihilated by $\m_0$. Now it is straightforward to conclude that the map 
\[
M\otimes_R f\colon M\otimes N' \to M\otimes N
\] 
has kernel and cokernel annihilated by $\m_1\subset \m_0^2$. Moreover, $M\otimes N'$ is a finitely presented module by Step $1$. Since $\m_1$ was an arbitrary finitely generated subideal of $\m$, we conclude that $M\otimes N$ is almost finitely presented due to Lemma~\ref{almost-finitely-presented}\ref{almost-finitely-presented-2}. \smallskip 

{\it Step 3. The general case}: Repeat the argument of Step~$2$ once again using Step~$2$ in place of Step~$1$ at the end, and Lemma~\ref{close-almost} in place of Lemma~\ref{almost-finitely-presented}\ref{almost-finitely-presented-2}.
\end{proof}

\begin{lemma}\label{base-change-hom-2var} Let $M$ be an almost finitely presented $R$-module, let $N$ be any $R$-module, and let $P$ be an almost flat $R$-module. Then the natural map $\rm{Hom}_R(M, N) \otimes_R P \to \rm{Hom}_R(M, N\otimes_R P)$ is an almost isomorphism. \smallskip

Similarly, $\rm{Hom}_{R^a}(M^a, N^a) \otimes_{R^a} P^a \to \rm{Hom}_{R^a}(M^a, N^a\otimes_{R^a} P^a)$ is an almost isomorphism for any almost finitely presented $R^a$-module $M^a$, any $R^a$-module $N^a$, and an almost flat $R^a$-module $P^a$.
\end{lemma}
\begin{proof}
Proposition~\ref{many-functors}\ref{many-functors-1} and~\ref{many-functors-3} ensure that it suffices to prove the claim for the case of honest $R$-modules $M$, $N$, and $P$. \smallskip

{\it Step 1. The case of a finitely presented module $M$}: We choose a presentation of $M$:
\[
R^n \to R^m \to M \to 0
\]
Then we use that $P$ is almost flat to get a morphism of almost exact sequences:
\[
\begin{tikzcd}
0 \arrow{r} & \rm{Hom}_R(M, N) \otimes_R P \arrow{r} \arrow{d} & \rm{Hom}_R(R^m, N) \otimes_R P \arrow{r}\arrow{d} & \rm{Hom}_R(R^n, N) \otimes_R P \arrow{d} \\ 
0 \arrow{r} & \rm{Hom}_R(M, N \otimes_R P) \arrow{r} & \rm{Hom}_R(R^m, \otimes_R P)  \arrow{r} & \rm{Hom}_R(R^n, N\otimes_R P). \\ 
\end{tikzcd}
\]
Clearly, the second and third vertical arrows are (almost) isomorphisms, so the first vertical arrow is an almost isomorphism as well. \medskip

{\it Step 2. The General Case}: The case of an almost finitely presented module $M$ follows from the finitely presented case by approximating $M$ by finitely presented $R$-modules. This is similar to the strategy used in Lemma~\ref{tensor-of-fp}; we leave the details to the reader. 
\end{proof}

The last thing we will need is the interaction between properties of an $R$-module $M$ and its ``reduction'' $M/I$ for some finitely generated ideal $I\subset \m$. For example, we know that for an ideal $I\subset \rm{rad}(R)$ and a finite module $M$, Nakayama's lemma states that $M/I=0$ if and only if $M=0$. Another thing is that an $I$-adically complete module $M$ is $R$-finite if and only if $M/I$ is $R/I$-finite. It turns out that both facts have their ``almost'' analogues.

\begin{lemma}\label{Nakayama} Let $I\subset \m \cap \rm{rad}(R)$ be a finitely generated ideal. If $M$ is an almost finitely generated $R$-module such that $M/IM\simeq 0$. Then $M \simeq 0$. If $M/IM\cong^a 0$, then $M\cong^a 0$.
\end{lemma}
\begin{proof}
We use the definition of an almost finitely generated module to find a finite submodule $N$ containing $IM$. If $M/IM$ is isomorphic to the zero module, then inclusion $IM \subset N \subset M$ implies that $N=M$. Thus $M$ is actually finitely generated, now we use the usual Nakayama's Lemma to finish the proof.

If $M/IM$ is merely almost isomorphic to the zero module, then we see that inclusion $IM \subset M$ is an almost isomorphism. In particular, $\m M$ is almost isomorphic to $IM$. Using that $\m^2=\m$, we obtain an {\it equality} 
\[
\m M=\m^2 M=\m(IM)=I(\m M)
\]
Thus we can apply the argument from above to conclude that $\m M=0$. This finishes the proof as $\m M \cong^a M$. 
\end{proof}

\begin{lemma}\label{finite-mod-ideal} Let $R$ be $I$-adically complete for some finitely generated $I\subset \m$. Then an $I$-adically complete $R$-module $M$ is almost finitely generated  if and only if $M/IM$ is almost finitely generated. 
\end{lemma}
\begin{proof}
\cite[Lemma 5.3.18]{GR}
\end{proof}

\subsection{Almost coherent modules and almost coherent rings}\label{section-almost-coherent}

This section is devoted to the study of almost coherent modules which are ``almost'' analogues of classical coherent modules. We show that these modules form a weak Serre subcategory in $\mathbf{Mod}_R$. Then we study the special case of almost coherent modules over an almost coherent ring. In this case, we show that almost coherent modules are equivalent to almost finitely presented modules.\smallskip

We recall that we fixed some ``base'' ring $R$ with an ideal $\m$ such that $\m^2=\m$ and $\widetilde{\m}=\m\otimes_R \m$ is flat, and always we do almost mathematics with respect to this ideal.

\begin{defn}\label{acoh} We say that an (almost) $R$-module $M$ is {\it almost coherent} if it is almost finitely generated and every almost finitely generated almost submodule $N^a \subset M^a$ is almost finitely presented.
\end{defn}

\begin{rmk} An almost submodule $f\colon N^a \hookrightarrow
 M^a$ does not necessarily give rise to a submodule $N'\subset M$ for some $(N')^a \simeq N^a$. The most we can say is that there is an injection $f_!\colon (N^a)_! \hookrightarrow
 (M^a)_!$ whose almostification is equal to the morphism $f$ (this follows from Lemma~\ref{prop-almost}\ref{prop-almost-2}). 
\end{rmk}

\begin{lemma}\label{lemma:almost-iso-almost-equi-coh} Let $R\to S$ be an almost isomorphism of rings. Then the forgetful functor $\bf{Mod}_{S^a}^{\rm{acoh}} \to \bf{Mod}_{R^a}^{\rm{acoh}}$ is an equivalence.
\end{lemma}
\begin{proof}
    This follows directly from Corollary~\ref{cor:almost-iso-almost-equi} and Definition~\ref{acoh}.
\end{proof}

\begin{lemma}\label{almost-coh-module} Let $M^a$ be an almost $R$-module with a representative $M\in \mathbf{Mod}_R$. Then the following are equivalent
\begin{enumerate}[label=\textbf{(\arabic*)}]
	\item The almost module $M^a$ is almost coherent.
	\item The $R$-module $(M^a)_*$ is almost finitely generated, and any almost finitely generated $R$-submodule of $(M^a)_*$ is almost finitely presented.
	\item The $R$-module $(M^a)_!$ is almost finitely generated, and any almost finitely generated $R$-submodule of $(M^a)_!$ is almost finitely presented.
\end{enumerate}
\end{lemma}
\begin{proof}
First of all, we note that Corollary~\ref{almost-almost-finitely-presented} guarantees that $M$ is almost finitely generated if and only if so is $(M^a)_*$. Second, Lemma~\ref{lemma:adjoint-almost} implies that the functor $(-)_*$ is left exact. Therefore, any almost submodule $N^a\subset M^a$ gives rise to an actual submodule $(N^a)_*\subset (M^a)_*$ that is almost isomorphic to $N$. In reverse, any submodule $N\subset (M^a)_*$ gives rise to an almost submodule of $M^a$. Hence, we see that all almost finitely generated almost submodules of $M^a$ are almost finitely presented if and only if all actual almost finitely generated submodules of $M_*$ are almost finitely presented (here we again use Corollary~\ref{almost-almost-finitely-presented}). This shows the equivalence of \textbf{(1)} and \textbf{(2)}. The same argument shows that \textbf{(1)} is equivalent to \textbf{(3)}.
\end{proof}

Note that it is not that clear whether a coherent $R$-module is almost coherent. The issue is that in the definition of almost coherent modules we need to be able to handle all almost finitely generated almost submodules and not only finitely generated. The lemma below is a useful tool to deal with such problems; in particular, it turns out (Corollary \ref{coh-acoh}) that all coherent modules are indeed almost coherent, but we do not know a direct way to see it.

\begin{lemma}\label{almost-coherent} Let $M$ be an $R$-module. Then $M$ is an almost coherent module if one of the following holds:
\begin{enumerate}[label=\textbf{(\arabic*)}]
\item\label{almost-coherent-1} For any finitely generated ideal $\m_0\subset \m$ there exists a coherent $R$-module $N$ and morphism $f\colon N \to M$ such that $\m_0(\ker f)=0$ and $\m_0(\coker f)=0$.
\item\label{almost-coherent-2} For any finitely generated ideal $\m_0\subset \m$ there exist an {\it almost} coherent $R$-module $N$ and a morphism $f\colon N \to M$ such that $\m_0(\ker f)=0$ and $\m_0(\coker f)=0$.
\end{enumerate}
\end{lemma}
\begin{proof}
We start the proof by noting that $M$ comes with the natural almost isomorphism $M \to M^a_*$. Since both assumptions on $M$ pass through this almost isomorphism,  Lemma~\ref{almost-coh-module} implies that it suffices to show that $M_*\coloneqq M^a_*$ is almost coherent. \smallskip

Lemma~\ref{almost-finitely-presented} guarantees that $M_*$ is almost finitely generated. Thus, we only need to check the second condition from Definition~\ref{acoh}. So we pick an arbitrary almost finitely generated $R$-submodule $M_1 \subset M_*$ and wish to show that it is almost finitely presented. We choose an arbitrary finitely generated ideal $\m_0\subset \m$ and another finitely generated ideal $\m_1 \subset \m$ such that $\m_0 \subset \m_1^2$. \smallskip

We use Lemma~\ref{close-almost} to find a morphism $\varphi\colon R^n \to M_1$ such that $\m_1(\coker \varphi)=0$. We denote by $e_1, \dots, e_n$ the standard basis of $R^n$ and by $x_i\coloneqq \varphi(e_i)$ the image of $e_i$ in $M_1$. We also choose a set of generators $(\e_1, \dots, \e_m)$ of the ideal $\m_1$. \smallskip

By assumption, there is a morphism $f\colon N \to M_*$ with a(n) (almost) coherent $R$-module $N$ such that $\m_1(\coker f)=0$ and $\m_1(\ker f)=0$. This implies that $\e_ix_j$ is in the image of $f$ for any $i=1, \dots, m, j=1,\dotsm n$. Let us choose some $y_{i,j} \in N$ such that $f(y_{i,j})=\e_ix_j$, and define an $R$-module $N'$ as the submodule of $N$ generated by all $y_{i,j}$. By construction, $N'$ is a finite $R$-module. Since $N$ is a (almost) coherent module, we conclude that $N'$ is (almost) finitely presented. \smallskip

We observe that $f'\coloneqq  f|_{N'}$ naturally lands in $M_1$, and we have $\m_1(\ker f')=0$ and $\m_1^2(\coker f')=0$. Since $\m_0\subset \m_1^2$, this shows that the morphism 
\[
N' \xr{f'} M_1
\]
has kernel and cokernel killed by $\m_0$. Lemma~\ref{close-almost} shows that $M_1$ is almost finitely presented.
\end{proof}

\begin{question} Does the converse of this Lemma hold?
\end{question} 

\begin{Cor}\label{coh-acoh} Any coherent $R$-module $M$ is almost coherent.
\end{Cor}

The next thing we want to show is that almost coherent modules form a weak Serre subcategory of $\mathbf{Mod}_R$. This is an almost analogue of the corresponding statement in the classical case.

\begin{lemma}\label{main-coh} Let $R$ and $\m$ as above. Then  
\begin{enumerate}[label=\textbf{(\arabic*)}]
\item\label{main-coh-1} An almost finitely generated almost submodule of an almost coherent module is almost coherent;
\item\label{main-coh-2} Let $\varphi\colon N^a \to M^a$ be an almost homomorphism from an almost finitely generated $R^a$-module to an almost coherent $R^a$-module, then $\ker \varphi$ is an almost finitely generated $R^a$-module;
\item\label{main-coh-3} Let $\varphi\colon N^a \to M^a$ be an injective almost homomorphism of almost coherent $R^a$-modules, then $\coker \varphi$ is an almost coherent $R^a$-module;
\item\label{main-coh-4} Let $\varphi\colon N^a \to M^a$ be an almost homomorphism of almost coherent $R^a$-modules, then $\ker \varphi$ and $\coker \varphi$ are almost coherent $R^a$-modules;
\item\label{main-coh-5} Given a short exact sequence of $R^a$-modules $0 \to M'^a \to M^a \to M''^a \to 0$ if two out of three are almost coherent, so is the third.
\end{enumerate}
\end{lemma}
\begin{proof}

\ref{main-coh-1}: This is evident from the definition of an almost coherent almost module. \smallskip

\ref{main-coh-2}: Let us define $N''^a:=\Imm \varphi$ and $N'^a:=\ker \varphi$, then Corollary~\ref{main-almost-almost} implies that $N''^a$ is an almost finitely generated almost submodule of $M^a$. Furthermore, it is almost finitely presented since $M^a$ is almost coherent. Thus, Corollary~\ref{main-almost-almost} implies that $N'$ is almost finitely generated as well. \smallskip

\ref{main-coh-3}: We denote $\coker \varphi$ by $M''^a$, then we have a short exact sequence
\[
0 \to N^a \to M^a \to M''^a \to 0.
\]
Corollary~\ref{main-almost-almost} implies that $M''^a$ is almost finitely generated. Let us choose any almost finitely generated almost submodule $M_1''^a \subset M''^a$ and denote its pre-image in $M^a$ by $M_1^a$. Then we have a short exact sequence
\[
0 \to N^a \to M_1^a \to M_1''^a \to 0.
\]
Corollary~\ref{main-almost-almost} guarantees that $M_1^a$ is an almost finitely generated almost submodule of $M^a$. Since $M^a$ is almost coherent, we see that $M_1^a$ is an almost finitely presented $R^a$-module. Therefore, Corollary~\ref{main-almost-almost} implies that $M_1''^a$ is also almost finitely presented. Hence, the $R^a$-module $M''^a$ is almost coherent. \smallskip

\ref{main-coh-4}: We know that $N'^a\coloneqq \ker \varphi$ is almost finitely generated by \ref{main-coh-2}. Since $N^a$ is almost coherent, we conclude that $N'^a$ is almost coherent by \ref{main-coh-1}. We define $N''^a\coloneqq \Imm \varphi$ and $M''^a\coloneqq \coker \varphi$, then we note that we have two short exact sequences
\begin{align*}
0 \to N'^a \to N^a \to N''^a \to 0, \\
0 \to N''^a \to M^a \to M''^a \to 0.
\end{align*}
We observe that \ref{main-coh-3} shows that $N''^a$ is almost coherent, then we use \ref{main-coh-3} once more to conclude that $M''^a$ is also almost coherent. \smallskip

\ref{main-coh-5}: The only thing that we are left to show is that if $M'^a$ and $M''^a$ are almost coherent, so is $M^a$. It is almost finitely generated by Corollary~\ref{main-almost-almost}. In order to check the second condition from Definition~\ref{acoh}, we choose an almost finitely generated almost submodule $M_1^a\subset M^a$. Let us denote by $M_1''^a$ its image in $M''^a$, and by $M_1'^a$ the kernel of this map. So we have a short exact sequence
\[
0 \to M_1'^a \to M_1^a \to M_1''^a \to 0.
\]
Corollary~\ref{main-almost-almost} guarantees that $M_1''^a$ is an almost finitely generated almost submodule of the almost coherent $R^a$-module $M''^a$. Hence, \ref{main-coh-1} implies that $M_1''^a$ is almost coherent, in particular, it is almost finitely presented. Moreover, we use \ref{main-coh-2} to see that $M_1'^a$ is an almost finitely generated almost submodule of $M'^a$. Since $M'^a$ is almost coherent, we conclude that $M_1'^a$ is almost finitely presented. Finally, Corollary~\ref{main-almost-almost} shows that $M_1^a$ is almost finitely presented as well. This finishes the proof of almost coherence of the $R^a$-module $M^a$.
\end{proof}

\begin{cor}\label{tensor-hom-coh} Let $M^a$ be an almost finitely presented $R^a$-modules and let $N^a$ be an almost coherent $R^a$-module. Then $M^a\otimes_{R^a} N^a$ and $\rm{alHom}_{R^a}(M^a, N^a)$ are almost coherent.  
\end{cor}
\begin{proof} We use Proposition~\ref{many-functors}\ref{many-functors-1},\ref{many-functors-3} to reduce the question to show that $M\otimes_R N$ and $\rm{Hom}_R(M, N)$ are almost coherent $R$-modules for any almost finitely presented $R$-module $M$ and almost coherent $R$-module $N$. \smallskip

{\it Step 1. The case of a finitely presented module $M$}: In this case, we choose a presentation of $M$ as the quotient
\[
R^n \to R^m \to M \to 0.
\]
Then we have short exact sequences
\[
N^n \to N^m \to M\otimes_R N \to 0 
\]
and 
\[
0 \to \rm{Hom}_R(M , N) \to N^m \to N^n.
\]
We note that Lemma~\ref{main-coh}\ref{main-coh-5} implies that $N^m$ and $N^n$ are almost coherent. Thus, Lemma~\ref{main-coh}\ref{main-coh-5} guarantees that both $M\otimes_R N$ and $\rm{Hom}_R(M , N)$ are almost coherent as well. \smallskip

{\it Step 2. The General Case}: The argument is similar to the one used in Step~2 of the proof of Lemma~\ref{tensor-of-fp}. We approximate $M$ by finitely presented $R$-modules. This gives us approximations of $M^a\otimes_{R^a} N^a$ and $\rm{alHom}_{R^a}(M^a, N^a)$ by almost coherent modules. Now Lemma~\ref{almost-coherent} guarantees that these modules are almost coherent. We leave the details to the interested reader. 
\end{proof}

We define $\bf{Mod}^{acoh}_R$ (resp. $\bf{Mod}_{R^a}^{acoh}$) to be the strictly full\footnote{i.e. full subcategory that is closed under isomorphisms.} subcategory of $\bf{Mod}_R$ (resp. $\bf{Mod}_{R^a}$) consisting of almost coherent $R$-modules (resp. $R^a$-modules).

\begin{Cor}\label{8} The category $\bf{Mod}^{\acoh}_R$ (resp. $\bf{Mod}_{R^a}^{\acoh}$) is a weak Serre subcategory of $\bf{Mod}_R$ (resp. $\bf{Mod}_{R^a}$).
\end{Cor}

Corollary~\ref{8} and the discussion in \cite[\href{https://stacks.math.columbia.edu/tag/06UP}{Tag 06UP}]{stacks-project} ensure that $\bf{D}_{acoh}(R)$ and\footnote{These are respectively full subcategories of $\bf{D}(R)$ and $\bf{D}(R)^a$ of complexes with almost coherent cohomology modules.} $\bf{D}_{acoh}(R)^a$ are strictly full saturated\footnote{A strictly full subcategory $\mathcal{D'}$ of a triangulated category $\mathcal{D}$ is saturated if $X\oplus Y \in \mathcal{D'}$ implies $X, Y \in \mathcal{D'}$.} triangulated subcategories of $\bf{D}(R)$ and $\bf{D}(R)^a$ respectively. We define $\mathbf{D}^+_{acoh}(R)\coloneqq \mathbf{D}_{acoh}(R)\cap  \mathbf{D}^+(R)$ and similarly for all other bounded versions. 

\begin{lemma}\label{lemma:coherent-complex-if-approximated} Let $M\in \bf{D}(R)$ be a complex of $R$-modules. Then $M\in \bf{D}_{acoh}(R)$ if one of the following holds:
\begin{enumerate}[label=\textbf{(\arabic*)}]
    \item For every finitely generated ideal $\m_0\subset \m$, there is $N\in \bf{D}_{coh}(R)$ and a morphism $f\colon N \to M$ such that $\m_0\left(\rm{H}^i\left(\rm{cone}\left(f\right)\right)\right)=0$ for every $i\in \Z$,
    \item For every finitely generated ideal $\m_0\subset \m$, there is $N\in \bf{D}_{acoh}(R)$ and a morphism $f\colon N \to M$ such that $\m_0\left(\rm{H}^i\left(\rm{cone}\left(f\right)\right)\right)=0$ for every $i\in \Z$.
\end{enumerate}
\end{lemma}
\begin{proof}
    This is an easy consequence of Lemma~\ref{almost-coherent} and the definition of $\bf{D}_{acoh}(R)$. 
\end{proof}

The last part of this subsection is dedicated to the study of almost coherent rings and almost coherent modules over almost coherent rings. Recall that coherent modules over a coherent ring coincide with finitely presented ones. Similarly, we will show that almost coherent modules over an almost coherent ring turn out to be the same as almost finitely presented ones.

\begin{defn} We say that a ring $R$ is almost coherent if the rank-$1$ free module $R$ is almost coherent as an $R$-module.
\end{defn}

\begin{lemma}\label{coh-acoh-2} A coherent ring $R$ is almost coherent.
\end{lemma}
\begin{proof}
Apply Corollary~\ref{coh-acoh} to a rank-$1$ free module $R$.
\end{proof}

\begin{lemma}\label{acoh-fp-coh} If $R$ is an almost coherent ring, then any almost finitely presented $R$-module $M$ is almost coherent.
\end{lemma}
\begin{proof}

{\it Step 1:} If $M$ is finitely presented over $R$, then we can write it as a cokernel of a map between free finite rank modules. A free finite rank module over an almost coherent ring is almost coherent due to Lemma~\ref{main-coh}\ref{main-coh-5}. A cokernel of a map of almost coherent modules is almost coherent due to Lemma~\ref{main-coh}\ref{main-coh-4}. Therefore, any finitely presented $M$ is almost coherent. \smallskip

{\it Step 2:} Suppose that $M$ is merely almost finitely presented. Lemma~\ref{almost-finitely-presented} guarantees that, for any finitely generated $\m_0\subset \m$, we can find a finitely presented module $N$ and a map $f\colon N \to M$ such that $\ker f$ and $\coker f$ are annihilated by $\m_0$. We know that $N$ is almost coherent by Step~1. Therefore, Lemma~\ref{almost-coherent}\ref{almost-coherent-2} implies that $M$ is almost coherent as well.
\end{proof}

\begin{Cor}\label{acoh-fp-coh-2} Let $R$ be an almost coherent ring. Then an $R$-module $M$ is almost coherent if and only if it is almost finitely presented.
\end{Cor}
\begin{proof}
The ``only if'' part is clear from the definition, the ``if'' part follows from Lemma \ref{acoh-fp-coh}.
\end{proof}

Our next big goal is to show that bounded above almost coherent complexes over an almost coherent ring are exactly almost pseudo-coherent complexes" in some precise way. More precisely, any element $M\in \bf{D}^-_{acoh}(R)$ can be ``approximated'' up to any small torsion by complexes of finite free modules.

\begin{prop}\label{prop:approximate-almost-coherent} Let $R$ be an almost coherent ring and $M\in \bf{D}^-(R)$. Then $M\in \bf{D}^-_{acoh}(R)$ if and only if, for every finitely generated ideal $\m_0\subset \m$, there is a complex $F^\bullet$ of finite free $R$-modules, and a morphism 
\[
f\colon F^\bullet \to M
\]
such that $\m_0\left(\rm{H}^i(\rm{cone}(f))\right)=0$ for every $i\in \Z$. Moreover, if $M\in \bf{D}^{\leq 0}_{coh}(R)$ one can choose $F^\bullet \in \bf{Comp}^{\leq 0}(R)$.
\end{prop}
\begin{proof}
    The ``if'' direction is Lemma~\ref{lemma:coherent-complex-if-approximated}. So we only need to prove the ``only if'' direction. For this direction, we fix a finitely generated ideal $\m_0\subset \m$ and another finitely generated ideal $\m_1\subset \m$ such that $\m_0\subset \m_1^2$. \smallskip
    
    Without loss of generality, we may and do assume that $M\in \bf{D}^{\leq 0}(R)$, and then we choose a complex $M^\bullet\in \bf{Comp}^{\leq 0}(R)$ that represents $M$. Now we prove a slightly more precise claim:
    
    {\it Claim: For every $n\in \Z$, there is a complex of finite free modules $F^\bullet_n$ with a morphism $f_n\colon F^\bullet_n \to M^\bullet$ such that
    \begin{enumerate}[label=\textbf{(\arabic*)}]
        \item $F^\bullet_n\in \bf{Comp}^{[-n, 0]}(R)$;
        \item $\sigma^{\geq -n+1} F^\bullet_n = F^\bullet_{n-1}$ and $\sigma^{\geq -n+1} f_n =f_{n-1}$, where $\sigma^{\geq n-1}$ is the naive truncation;
        \item kernels and cokernels of $\rm{H}^i(f_n)$ are annihilated by $\m_1$ for $i\geq n+1$;
        \item the cokernel of $\rm{H}^n(f_n)$ is annihilated by $\m_1$;
    \end{enumerate}}
    
    {\it Proof of the claim:} We argue by descending induction on $n$. If $n\geq 1$, $F^\bullet =0$ works. Now we suppose that we can construct $F^\bullet_n$, and wish to construct $F^\bullet_{n-1}$. Consider the morphism $f_n$ presented as a commutative diagram
    \[
    \begin{tikzcd}
    0\arrow{r}\arrow{d} & 0 \arrow{r}\arrow{d}& F^n_n \arrow{r}{\rm{d}^n_F} \arrow{d}{f^n_n} & F^{n+1}_n \arrow{r}{\rm{d}^{n+1}_F}\arrow{d}{f^{n+1}_n} \arrow {r} & \dots \arrow{d} \\
    M^{n-2} \arrow{r}{\rm{d}_M^{n-2}}& M^{n-1} \arrow{r}{\rm{d}_M^{n-1}} & M^{n}\arrow{r}{\rm{d}_M^n} & M^{n+1}\arrow{r}{\rm{d}_M^{n+1}} &\dots
    \end{tikzcd}
    \]
    Firstly, $\ker(\rm{d}^n_F)$ is almost coherent as a kernel between finitely presented modules over an almost coherent ring. Secondly, the $R$-module 
    \[
    B^n\coloneqq \ker\left(\ker\left(\rm{d}^n_F\right) \to \rm{H}^n\left(M\right)\right),
    \]
    is also almost coherent as a kernel between almost coherent modules. Therefore, there is a finite free $R$-module $F'^{n-1}$ and a morphism 
    \[
    \rm{d}'\colon F'^{n-1}\to B^n
    \]
    such that $\m_1(\coker \rm{d}')=0$. Since $\rm{H}^{n-1}(M)$ is almost coherent, we can find a finite free $R$-module $F''^{n-1}$ and a morphism 
    \[
    \lambda\colon F''^{n-1} \to \rm{H}^{n-1}(M)
    \]
    such that $\m_1(\coker \lambda)=0$. Let $\nu\colon F''^{n-1} \to Z^{n-1}(M^\bullet)$ be any lift of $\lambda$ to the module of closed elements $Z^{n-1}(M^\bullet)=\ker(\rm{d}_M^{n-1})$. We define
    \[
    f''^{n-1}\colon F''^{n-1} \to M^{n-1}
    \]
    be the composition of $\nu$ with the inclusion $Z^{n-1}(M^\bullet) \to M^{n-1}$. \smallskip
    
    Now we wish to define $F^\bullet_{n-1}$ and $f_{n-1}$. We start with $F^\bullet_{n-1}$; we put $F^{m}_{n-1}=F^m_n$ if $m\geq n$, $F^m_{n-1}=0$ if $m<n-1$, $F^{n-1}_{n-1}=F'^{n-1}\oplus F''^{n-1}$, and define the only non-evident differential 
    \[
    \rm{d}^{n-1}_F\colon F^{n-1}_{n-1}=F'^{n-1}\oplus F''^{n-1}\to F^{n}_n
    \]
    to be zero on $F''^{n-1}$ and equal to $\rm{d}'$ on $F'^{n-1}$. It is evident that $\rm{d}^n_F \circ \rm{d}^{n-1}_F=0$, so this structure defines a complex $F^\bullet_{n-1}$ of finite free $R$-modules. \smallskip
    
    We are only left to define $f_{n-1}$. We must put $f^m_{n-1}=f^m_n$ if $m>n-1$ and $f^m_{n-1}=0$ if $m<n-1$, so the only question is to define $f^{n-1}_{n-1}$. By construction, we have $f_n^n(\rm{d}'F'^{n-1})\subset \rm{d}^{n-1}_MM^{n-1}$, so we can find 
    \[
    f'_{n-1}\colon F'^{n-1}\to M^{n-1}
    \]
    such that $\rm{d}^{n-1}\circ f'_{n-1}=f_n^n \circ \rm{d}'$. Thus we define 
    \[
    f^{n-1}_{n-1}\colon F^{n-1}_{n-1}=F'^{n-1}\oplus F''^{n-1} \to M^{n-1}
    \]
    to be $f'_{n-1}$ on $F'^{n-1}$ and $f''_{n-1}$ on $F''^{n-1}$. Then it is evident from the construction that $f^{\bullet}_{n-1}$ is a morphism of complexes, i.e. the diagram
    \[
    \begin{tikzcd}[column sep=4em, row sep = 3em]
    0\arrow{r}\arrow{d} & F^{n-1}_{n-1} \arrow{r}{\rm{d}^{n-1}_F}\arrow{d}{f^{n-1}_{n-1}}& F^n_{n-1} \arrow{r}{\rm{d}^n_F} \arrow{d}{f^n_{n-1}} & F^{n+1}_{n-1} \arrow{r}{\rm{d}^{n+1}_F}\arrow{d}{f^{n+1}_{n-1}} \arrow {r} & \dots \arrow{d} \\
    M^{n-2} \arrow{r}{\rm{d}^{n-2}_M}& M^{n-1} \arrow{r}{\rm{d}^{n-1}_M} & M^{n}\arrow{r}{\rm{d}^n_M} & M^{n+1}\arrow{r}{\rm{d}^{n+1}_M} &\dots
    \end{tikzcd}
    \]
    By construction, the kernel and cokernel of $\rm{H}^{n}(f_{n-1})$ are annihilated by $\m_1$, and the cokernel of $\rm{H}^{n-1}(f_{n-1})$ is annihilated by $\m_1$. So this finishes the proof of the claim. \smallskip
    
    Now the morphism $f\colon F^\bullet \to M^\bullet$ simply comes as the colimit of $f_n$, i.e.  
    \[
    f=\colim f_n \colon F^\bullet \coloneqq \colim F^\bullet_n \to M^\bullet. 
    \]
    It is easy to see that the cohomology groups of $\rm{cone}(f)$ are annihilated by $\m_0\subset \m_1^2$. 
\end{proof}

\begin{cor} Let $R$ be a coherent ring and $M\in \bf{D}^b(R)$. Then $M\in \bf{D}^b_{acoh}(R)$ if and only if, for every finitely generated ideal $\m_0\subset \m$, there is a complex $N\in \bf{D}^b_{coh}(R)$ and a morphism $f\colon N \to M$ such that $\m_0(\rm{H}^i(\rm{cone}(f)))=0$ for all $i$.
\end{cor}
\begin{proof}
    The ``if'' direction is Lemma~\ref{lemma:coherent-complex-if-approximated}. So we only need to deal with the ``only if'' direction. Assume that $M\in \bf{D}^b(R)$. Then Proposition~\ref{prop:approximate-almost-coherent} implies that there is $F\in \bf{D}^-_{coh}(R)$ and a morphism $f\colon F\to M$ such that $\m_0(\rm{H}^i(\rm{cone}(f)))=0$ for all $i$. Now replace $F$ by $F'\coloneqq \tau^{\geq a}F$ to get the desired approximation with $F'\in \bf{D}^b_{coh}(R)$.
\end{proof}

\begin{prop}\label{tensor-product-coh} Let $R$ be an almost coherent ring, and let $M^a, N^a$ be objects in $\bf{D}^-_{acoh}(R)^a$. Then $M^a\otimes_{R^a}^L N^a\in \bf{D}^-_{acoh}(R)^a$. 
\end{prop}
\begin{proof}
Proposition~\ref{derived-tensor-product} ensures that it suffices to show that $M\otimes^L_R N\in \bf{D}^-_{acoh}(R)$ for $M$, $N\in \bf{D}^-_{coh}(R)$. Clearly, we can cohomologically shift both $M$ and $N$ to assume that they lie $\bf{D}^{\leq 0}_{coh}(R)$.\smallskip 

Now we fix a finitely generated ideal $\m_1 \subset \m$ and use Proposition~\ref{prop:approximate-almost-coherent} to find an exact triangle 
\[
F^\bullet \to M \to Q
\]
where $F^\bullet \in \bf{D}^{\leq 0}(R)$ a complex of finite free modules and $\rm{H}^i(Q)$ are all annihilated by $\m_1$. Then it is easy to see that the kernel and cokernel of the map
\[
\rm{H}^{-i}(F^\bullet \otimes^L_R N) \to \rm{H}^{-i}(M\otimes^L_R N)
\]
are annihilated by $\m_1^{i+1}$. Now we note that, clearly, 
\[
F^\bullet \otimes^L_R N \simeq F^\bullet \otimes^\bullet_R N
\]
lies in $\bf{D}^-_{coh}(R)$ because $F^\bullet$ is a complex of finite free modules. For each pair of an integer $i\geq 0$ and a finitely generated ideal $\m_0 \subset \m=\m^{i+1}$, we can find another finitely generated ideal $\m_1$ such that $\m_0 \subset \m_1^{i+1}$. Therefore, the map
\[
\rm{H}^{-i}(F^\bullet \otimes^L_R N) \to \rm{H}^{-i}(M\otimes^L_R N)
\]
is a morphism with an almost coherent source and $\m_0$-torsion kernel and cokernel. Therefore, Lemma~\ref{almost-coherent}~\ref{almost-coherent-2} implies the claim. 
\end{proof}

\begin{prop}\label{alHom-derived-coh} Let $R$ be an almost coherent ring, and let $M^a \in \bf{D}^-_{acoh}(R)^a, N^a \in \bf{D}^+_{acoh}(R)^a$. Then $\bf{R}\rm{alHom}_{R^a}(M^a, N^a)\in \bf{D}^+_{acoh}(R)^a$. 
\end{prop}
\begin{proof}
The proof is similar to that of Proposition~\ref{tensor-product-coh}. We use Proposition~\ref{derived-al-hom} and the same approximation argument to reduce to the case $M=F^\bullet$ is a bounded above complex of finite free modules. In this case, the claim is essentially obvious due to the explicit construction of the Hom-complex $\rm{Hom}^\bullet_R(F^\bullet, N)$. 
\end{proof}

\begin{prop}\label{base-change-hom-2var-derived} Let $R$ be an almost coherent ring, let $M\in \bf{D}_{acoh}^-(R)$, $N\in \bf{D}^+(R)$, and let $P$ be an almost flat $R$-module. Then the natural map $\bf{R}\rm{Hom}_R(M, N) \otimes_R P \to \bf{R}\rm{Hom}_R(M, N\otimes_R P)$ is an almost isomorphism. \smallskip

Similarly, $\bf{R}\rm{Hom}_{R^a}(M^a, N^a) \otimes^L_{R^a} P^a \to \bf{R}\rm{Hom}_{R^a}(M^a, N^a\otimes^L_{R^a} P^a)$ is an almost isomorphism for any $M^a\in \bf{D}^-_{acoh}(R)^a$, $N^a\in \bf{D}^+(R)^a$, and let $P^a$ an almost flat $R^a$-module.
\end{prop}
\begin{proof}
The proof is similar to that of the above lemmas.
\end{proof}

\begin{cor} Let $R$ be an almost coherent ring, let $M^a\in \bf{D}_{acoh}^-(R)^a$, $N\in \bf{D}^+(R)^a$, and let $P^a$ be an almost flat $R^a$-module. Then the natural map 
\[
\bf{R}\rm{alHom}_{R^a}(M^a, N^a) \otimes^L_{R^a} P^a \to \bf{R}\rm{alHom}_{R^a}(M^a, N^a\otimes_{R^a} P^a)
\]
is an  isomorphism in $\bf{D}(R^a)$.
\end{cor}

\subsection{Almost noetherian rings}

The main goal of this section is to define the almost analogue of the noetherian property. We also verify some of its basic properties. Even though most of the basic facts about noetherian rings carry over to the almost world, we warn the reader that Hilbert's Nullstellensatz seems to be more subtle in the almost world (see Warning~\ref{warn:no-hilbert}); we are able to establish it only in some very particular situations in Section~\ref{section:Examples}. \smallskip

As in the previous sections, we fix a ring $R$ with an ideal $\m$ such that $\m^2=\m$ and $\widetilde{\m}=\m\otimes_R \m$ is flat, and we always do almost mathematics with respect to this ideal. 

\begin{defn}\label{defn:almost-noetherian} A ring $R$ is {\it almost noetherian} if every ideal $I\subset R$ is almost finitely generated. 
\end{defn}

The main goal is to show that every almost finitely generated module over an almost noetherian ring is almost finitely presented. In particular, an almost noetherian ring is almost coherent.  

\begin{lemma}\label{lemma:almost-noetherian-Rn} Let $R$ be an almost noetherian ring, and $M\subset R^n$ an $R$-submodule. Then $M$ is almost finitely generated. 
\end{lemma}
\begin{proof}
    We argue by induction on $n$. The base of induction is $n=1$, where the claim follows from the definition of an almost noetherian ring. \smallskip
    
    Suppose we know the claim for $n-1$, so we deduce the claim for $n$. Denote by $R^{n-1} \subset R^n$ a free $R$-module spanned by first $n-1$ standard basis elements of $R^n$, and denote by $M'\coloneqq M \cap R^{n-1}$ the intersection of $M$ with $R^{n-1}$. Then we have a short exact sequence
    \[
    0 \to M' \to M \to M''\to 0,
    \]
    where $M''$ is naturally an $R$-submodule of $R\simeq R^{n}/R^{n-1}$. By the induction hypothesis, $M'$ is almost finitely generated. $M''$ is almost finitely generated by almost noetherianness of $R$. Therefore, $M$ is almost finitely generated by Lemma~\ref{main-almost}~\ref{main-almost-2}.
\end{proof}

\begin{lemma}\label{lemma:almost-noetherian-almost-fg-almost-fp} Let $R$ be an almost noetherian ring. Then any almost finitely generated $R$-module $M$ is almost finitely presented.
\end{lemma}
\begin{proof}
    Pick any finitely generated sub-ideal $\m_0 \subset \m$. By Lemma~\ref{almost-finite}, there is an $R$-linear homomorphism 
    \[
        f\colon R^n \to M
    \]
    such that $\m_0(\coker f)=0$. Consider $N \coloneqq \ker(f)$. Lemma~\ref{lemma:almost-noetherian-Rn} ensures that $N$ is also almost finitely generated, so there is an $R$-linear homomorphism
    \[
        g'\colon R^m \to N
    \]
    such that $\m_0(\coker g')=0$. Therefore, the composition
    \[
    R^m \xr{g} R^n \xr{f} M
    \]
    is a three-term complex with $\m_0(\coker f)=0$ and $\m_0(\ker f) \subset \rm{Im}(g)$. Since $\m_0$ was an arbitrary finitely generated sub-ideal in $\m$, we conclude that $M$ is almost finitely presented by Lemma~\ref{almost-finitely-presented}~\ref{almost-finitely-presented-3}.
\end{proof}

\begin{cor}\label{cor:almost-noetherian-aft=afpr} A ring $R$ is almost noetherian if and only if any almost finitely generated $R$-module $M$ is almost finitely presented.
\end{cor}
\begin{proof}
    If $R$ is almost noetherian, then any almost finitely generated $R$-module is almost finitely presented due to Lemma~\ref{lemma:almost-noetherian-almost-fg-almost-fp}. \smallskip
    
    Now we suppose that every almost finitely generated $R$-module is almost finitely presented, and we wish to show that $R$ is almost noetherian. Consider an ideal $I \subset R$. Then $R/I$ is clearly a finitely generated $R$-module, in particular, it is almost finitely generated. Therefore, it is almost finitely presented by our assumption on $R$. Now the short exact sequence
    \[
    0\to I \to R \to R/I \to 0
    \]
    and Lemma~\ref{main-almost}~\ref{main-almost-3} imply that $I$ is almost finitely generated. 
\end{proof}

\begin{cor}\label{cor:almost-noetherian-iso} Let $R\to R'$ be an almost isomorphism of rings. Then $R$ is almost noetherian if and only if $R'$ is.
\end{cor}

\begin{cor}\label{cor:submodule-of-aft-noetherian} Let $R$ be an almost noetherian ring, and $M$ an almost finitely generated $R$-module. Then any submodule $N \subset M$ is almost finitely generated.
\end{cor}
\begin{proof}
    Consider the short exact sequence
    \[
    0 \to N \to M \to M/N \to 0.
    \]
    By construction, $M/N$ is almost finitely generated and, therefore, almost finitely presented by Lemma~\ref{lemma:almost-noetherian-almost-fg-almost-fp}. So Lemma~\ref{main-almost}~\ref{main-almost-3} implies that $N$ is almost finitely generated. 
\end{proof}

\begin{cor}\label{cor:almost-noetherian-almost-coherent} Let $R$ be an almost noetherian ring. Then $R$ is almost coherent.
\end{cor}
\begin{proof}
    Lemma~\ref{almost-coh-module} guarantees that it suffices to show that $R_!\simeq \widetilde{\m}$ is almost finitely generated and every finitely generated sub-module of $R_!$ is almost finitely presented. The first property is trivial since $R_!$ is almost isomorphic to $R$, and the second one follows from Lemma~\ref{lemma:almost-noetherian-almost-fg-almost-fp}.
\end{proof}

\begin{cor}\label{cor:almost-noetherian-almost-fg=almost-coh} Let $R$ be an almost noetherian ring. Then an $R$-module $M$ (resp. an $R^a$-module $M^a$) is almost coherent if and only if it is almost finitely generated. 
\end{cor}
\begin{proof}
    It suffices to prove the claim for an honest $R$-module $M$. Corollary~\ref{cor:almost-noetherian-almost-coherent} and Corollary~\ref{acoh-fp-coh-2} imply that $M$ is almost coherent if and only if it is almost finitely presented. Now Lemma~\ref{lemma:almost-noetherian-almost-fg-almost-fp} says that $M$ is almost finitely presented if and only if it is almost finitely generated. This finishes the proof. 
\end{proof}

\begin{warning}\label{warn:no-hilbert} Unlike the case of usual noetherian rings, Hilbert's Nullstellensatz is more subtle in the almost world. In particular, we do not know if a polynomial algebra in a finite number of variables over an almost noetherian ring is almost noetherian. However, we show will that Hilbert's Nullstellensatz holds for perfectoid valuation rings in Section~\ref{section:Examples}.
\end{warning}

\begin{example} Let $\bf{B}_I$ be the period ring from \cite[Definition 1.6.2]{FF}. Then \cite[Corollary 8.16]{Wear} implies that the rings $\bf{B}_I^+$ are almost noetherian for any closed interval $I\subset (0, \infty)$. Another family of examples of almost noetherian rings will be constructed in Section~\ref{section:Examples}. 
\end{example}

\subsection{Base change for almost modules}\label{base-change-section}

In this section, we discuss the behavior of almost modules with respect to base change. Recall that, for a ring homomorphism $\varphi \colon  R\to S$, we always do almost mathematics on $S$-modules with respect to the ideal $\m_S\coloneqq \m S$; look at Lemma~\ref{base-change} for details.

\begin{lemma}\label{trivial-base-change} Let $\varphi \colon R \to S$ be a ring homomorphism, and let $M^a$ be an almost finitely generated (resp. almost finitely presented) $R^a$-module. Then the module $M^a_S\coloneqq M^a\otimes_{R^a} S^a$ is almost finitely generated (resp. almost finitely presented).
\end{lemma}
\begin{proof}
The claim follows from Lemma~\ref{almost-finitely-presented}\ref{almost-finitely-presented-2} and the fact that, for any finitely generated ideal $\m'_0\subset \m_S$, there is a finitely generated ideal $\m_0 \subset \m$ such that $\m'_0 \subset \m_0S$. We give a complete proof only in the case of finitely presented modules because the other case is an easier version of the same argument. \smallskip

First, we note that it suffices to show that $M\otimes_R S$ is almost finitely presented. Now we note that, for any finitely generated ideal $\m'_0\subset \m_S$, there is a finitely generated ideal $\m_0 \subset \m$ such that $\m'_0 \subset \m_0S$. Therefore, it suffices to check the condition of Lemma~\ref{almost-finitely-presented}\ref{almost-finitely-presented-2} only for ideals of the form $\m_0S$, where $\m_0\subset \m$ is a finitely generated sub-ideal. Then we choose some finitely generated ideal $\m_1\subset \m$ such that $\m_0\subset \m_1^2$ and use Lemma~\ref{almost-finitely-presented}\ref{almost-finitely-presented-2} to find a finitely presented module $N$ and a map $f\colon N \to M$ such that $\m_1(\Ker f)=\m_1(\coker f)=0$. Consider an exact sequence 
\[
0 \to K \to N \xr{f} M \to Q \to 0
\]
and denote the image $f$ by $M'$. Then we have the following exact sequences:
\[
K\otimes_R S \to N\otimes_R S \to M'\otimes_R S \to 0
\]
\[
\Tor_1^R(Q, S) \to M'\otimes_R S \to M\otimes_R S \to Q\otimes_R S
\]
Since $K\otimes_R S$, $\Tor_1^R(Q, S)$ and $Q\otimes_R S$ are killed by $\m_1S$, we conclude that $\coker(f\otimes_R S)$ and $\ker(f\otimes_R S)$ are annihilated by $\m_1^2S$. In particular, they are killed by $\m_0S$. Since $N\otimes_R S$ is finitely presented over $S$, Lemma~\ref{almost-finitely-presented} finishes the proof.
\end{proof}

\begin{cor}\label{derived-trivial-base-change} Let $R\to S$ be a ring homomorphism of almost coherent rings, and let $M^a$ be an object of $\bf{D}^-_{acoh}(R)^a$. Then $M^a\otimes_{R^a}^L S^a\in \bf{D}^-_{acoh}(S)^a$. 
\end{cor}
\begin{proof}
The proof is similar to that of Proposition~\ref{tensor-product-coh}. We use Proposition~\ref{derived-base-change} and a similar approximation argument based on Proposition~\ref{prop:approximate-almost-coherent} to reduce to the case $M\simeq F^\bullet$, where $F^\bullet$ is a bounded above complex of finite free modules. In this case, the claim is essentially obvious. 
\end{proof}

\begin{lemma}\label{extension} Let $S$ be a $R$-algebra that is finite (resp. finitely presented) as an $R$-module, and let $M^a$ be an $S^a$-module. Then $M^a$ is almost finitely generated (resp. almost finitely presented) over $R^a$ if and only if it is almost finitely generated (resp. almost finitely presented) over $S^a$.
\end{lemma}
\begin{proof}
As always, we first reduce the question to the case of an honest $S$-module $M$. Now we use the observation that it suffices to check the condition of Lemma~\ref{almost-finitely-presented}\ref{almost-finitely-presented-2} only for the ideals of the form $\m_0S$ for some finitely generated ideal $\m_0\subset \m \subset R$. Then the only non-trivial direction is to show that $M$ is almost finitely presented over $S$ if it is almost finitely presented over $R$. This is proven in a more general situation in Lemma~\ref{extension-noncommutative}
\end{proof} 

\begin{lemma}\label{extension-noncommutative} Let $S$ be a possibly non-commutative $R$-algebra that is finite as a left (resp. right) $R$-module, and let $M$ be a left (resp. right) $S$-module that is almost finitely presented over $R$. Then $M$ is almost finitely presented over $S$ (i.e. for every finitely generated ideal $\m_0\subset \m$, there exists a finitely presented left (resp. right) $S$-module $N$ and a map $N \to M$ such that $\ker f$ and $\coker f$ are annihilated by $\m_0$).
\end{lemma}

\begin{rmk} This lemma will actually be used for a non-commutative ring $S$ in the proof of Theorem~\ref{check-perfect-pn} that, in turn, will be used in the proof of formal GAGA for almost coherent sheaves Theorem~\ref{GAGA}. Namely, we will apply to result to $S=\rm{End}_{\bf{P}^N}(\O\oplus \O(1)\oplus \dots \O(N))$. \smallskip

Besides this application, we will usually use Lemma~\ref{extension-noncommutative} when $R$ and $S$ are almost coherent commutative rings. In this case, the proof of Lemma~\ref{extension-noncommutative} can be significantly simplified. 
\end{rmk}

\begin{proof}
We give a proof for left $S$-modules; the proof for right $S$-modules is the same. We start the proof by choosing some generators $x_1, \dots, x_n$ of $S$ as an $R$-module. Then we pick a finitely generated ideal $\m_0 \subset \m$ and another finitely generated ideal $\m_1$ such that $\m_0 \subset \m_1^2$. And we also choose some generators $(\e_1, \dots, \e_k)=\m_1$ and find a three-term complex
\[
R^t \xr{g} R^m \xr{f} M
\]
such that $\m_1(\coker f)=0$ and $\m_1(\ker f) \subset \Imm g$. We consider the images $y_i\coloneqq f(e_i)\in M$ of the standard basis elements in $R^m$. Then we can find some $\beta_{i,j,s,r}\in R$ such that
\[
\e_sx_iy_j=\sum_{r=1}^m \beta_{i,j,s,r} \cdot y_r \text{ with } \beta_{i,j,s,r} \in R
\]
for any $s=1, \dots k; \ i=1, \dots, n; \ j = 1, \dots, m$. Furthermore, we have $t$ ``relations''
\[
\sum_{j=1}^m \a_{i,j} y_j = 0 \text{ with } \a_{i,j} \in R
\]
such that for any relation $\sum_{i=1}^m b_{i} y_i = 0$ with $b_i \in R$ and any $\e \in \m_1$, we have that the vector $\{\e b_i\}_{i=1}^m \in R^m$ lives in the $R$-subspace generated by vectors $\{\a_{i,j}\}_{i=1}^m$ for $j=1,\dots, t$. Or, in other words, if $\sum_{j=1}^m \a_{i,j} y_j = 0$ then $\e(\sum_{j=1}^m \a_{i,j} e_j) \in \rm{Im}(g)$ for any $\e\in \m_1$. \smallskip

Now we are finally ready to define a three-term complex
\[
S^{nmk+t } \xr{\psi} S^{m} \xr{\varphi} M
\]
We define the map $\varphi$ to be the unique $S$-linear homomorphism such that $\varphi(e_i)=y_i$ for the standard basis in $S^m$. We define $\psi$ as the unique $S$-linear homomorphism such that 
\[
\psi(f_{i,j,s})=\e_sx_ie_j - \sum_{r=1}^m \beta_{i,j,s,r}\cdot e_r \text { and } \psi(f'_l)= \sum_{j=1}^m \a_{l,j}e_j
\]
for the standard basis 
\[
\left\{f_{i,j,s}, f'_l\right\}_{i\leq n, j\leq m, s\leq k, l\leq t} \in S^{nmk+t}
\]  
Then we clearly have that $\varphi\circ \psi=0$ and that $\m_1(\coker \varphi)=0$. We claim that $\m_1^2(\ker\varphi)\subset \Imm \psi$. \smallskip

Let $\varphi(\sum_{i=1}^m c_ie_i)=0$ for some elements $c_i\in S$. We can write each 
\begin{equation}\label{eqn:3}
    c_i=\sum_{j=1}^n r_{i,j}x_j \text{ with } r_{i,j} \in R
\end{equation}

because $x_1, \dots, x_n$ are $R$-module generators of $S$. Thus, the condition that $\varphi(\sum_{i=1}^m c_ie_i)=0$ is equivalent to $\sum_{i,j} r_{i,j}x_jy_i=0$. Now recall that for any $s=1, \dots k$ we have
\[
\e_sx_jy_i=\sum_{r=1}^m \beta_{j,i,s,r}\cdot y_r.
\] 
Therefore, multiplying equation~\ref{eqn:3} by $\e_s$, we get an equality 
\[
0=\e_s\left(\sum_{i,j} r_{i,j}x_jy_i\right) =\sum_{i,j}r_{i,j} \left(\sum_{r=1}^m \beta_{j,i,s,r}\cdot y_r\right)=\sum_{r=1}^m \left(\sum_{i,j} r_{i,j} \beta_{j,i,s,r}\right)y_r
\]
This means that for any $s'=1, \dots, k$ the vector $\{\e_{s'}(\sum_{i,j} r_{i,j}\beta_{j,i,s,r})\}_{r=1}^m\in R^m$ lives in an $R$-subspace generated by vectors $\{\a_{i,j}\}_{i=1}^m$. In particular, for any $r$ and $s'$, $\e_{s'}(\sum_{i,j} r_{i,j}\beta_{j,i,s,r}e_r)$ is equal to $\psi\left(\text{some sum of }f'_l\right)$ by the definition of $\psi$. \smallskip

After unwinding all the definitions we get the following:
\begin{align*}
\e_{s'}\e_{s}\left(\sum_{i=1}^m c_i e_i\right)&=\e_{s'}\e_{s} \left( \sum_{i,j} r_{i,j}x_je_i \right) \\
&=\e_{s'} \left(\sum_{i,j} r_{i,j}\left(\e_sx_je_i-\sum_r \beta_{j,i,s,r}e_r+\sum_r \beta_{j,i,s,r}e_r\right)\right) \\ 
&=\e_{s'}\left(\sum_{i,j} r_{i,j}\left(\e_sx_je_i-\sum_r \beta_{j,i,s,r}e_r\right)\right) +\e_{s'}\left(\sum_r\left(\sum_{i,j} r_{i,j}\beta_{j,i,s,r}\right)e_r \right) \\
&= \psi\left(\e_{s'}\sum_{i,j}r_{i,j}f_{j,i,s}\right)+\psi\left(\text{some sum of }f'_l\right)
\end{align*}

So we see that $\m_1^2\ker(\varphi) \subset \Imm \psi$. In particular, we have $\m_0\ker (\varphi) \subset \Imm \psi$. Now we replace the map $\varphi \colon S^n \to M$ with the induced map
\[
\ov{\varphi}\colon \coker (\psi) \to M
\]
to get a map from a finitely presented left $S$-module such that $\ker(\ov{\varphi})$ and $\coker(\ov{\varphi})$ are annihillated by $\m_0$.
\end{proof}

\subsection{Almost faithfully flat algebras}

In this section, we study the almost faithfully flat morphisms of {\it algebras}. This notion turns out to be quite subtle in the almost world due to the following two observations. The first observation is that, for an almost faithfully flat morphism $R\to S$, the $R^a_!$-module $S^a_!$ is always flat, but not necessarily faithfully flat (see Warning~\ref{warning:not-fully-faithful}). Another observation is that $S^a_!$ usually does not have a structure of an $R$-algebra. \smallskip

For these reasons, it is not evident how to relate almost faithful flatness of an $R$-algebra $S$ to some classical faithful flatness. In order to make this possible, we replace the $(-)_!$-functor with another functor $(-)_{!!}$ that takes into account the $R$-algebra structure on $S$. This functor will send almost faithfully flat $R$-algebras into faithfully flat $R$-algebra, however, it will not, in general, send flat $R$-algebras into flat $R$-algebras. However, this functor will suffice for the purpose of studying almost faithfully flat morphisms. \smallskip

In this section, we follow \cite{GR} pretty closely. \smallskip

For the rest of the section, we fix a ring $R$ with an ideal of almost mathematics $\m$.  

\begin{defn} A homomorphism of $R$-algebras $A\to B$ is {\it almost flat} (resp. {\it almost faithfully flat}) if $B^a$ is a flat (resp. faithfully flat) $A^a$-module (see Definition~\ref{defn:almost-flat}).
\end{defn}

\begin{lemma} Any (faithfully) flat $A$-algebra $B$ is almost (faithfully) flat. 
\end{lemma}
\begin{proof}
    This follows directly from Lemma~\ref{where-sends}.
\end{proof}

\begin{lemma}\label{lemma:almost-ff-almost-injective} Let $A$ be an $R$-algebra and $f\colon A \to B$ a morphism of $R$-algebras. Then $B$ is almost faithfully flat over $A$ if and only if $B^a$ is a flat $A^a$-module and $A^a \to B^a$ is universally injective, i.e., for any $A^a$-module $M^a$, the natural morphism $M^a\to M^a\otimes_{A^a} B^a$ is injective in $\bf{Mod}_A^a$.
\end{lemma}
\begin{proof}
    Suppose that $B$ is almost faithfully flat. Then $B^a$ is a flat $A^a$-module by definition. So we only need to show that $A^a\to B^a$ is universally injective. Pick any $M^a\in \bf{Mod}_A^a$ and consider an $A^a$-module
    \[
    N^a\coloneqq \ker(M^a\to M^a\otimes_{A^a} B^a).
    \]
    Flatness of $B^a$ implies that the morphism
    \[
    N^a\otimes_{A^a} B^a \to M^a\otimes_{A^a} B^a
    \]
    is injective.  Now we also note that the morphism
    \[
    N^a\otimes_{A^a} B^a \to M^a\otimes_{A^a} B^a  
    \]
    is equal to zero by our choice of $N^a$. This implies that $N^a\otimes_{A^a} B^a \simeq 0$. Since $B^a$ is faithfully flat over $A^a$, we conclude that $N^a\simeq 0$.\smallskip 
    
    Now we suppose that $B^a$ is a flat $A^a$-module and $A^a\to B^a$ is universally injective. Thus, for any $A^a$-module $M^a$, we have an injection $M^a\to M^a\otimes_{A^a} B^a$. So if $M^a\otimes_{A^a} B^a\simeq 0$, we conclude that $M^a\simeq 0$. Thus $B^a$ is faithfully flat over $A^a$. 
\end{proof}

\begin{cor}\label{cor:almost-ff-cokernel-flat} Let $A$ be an $R$-algebra and $f\colon A \to B$ is a morphism of $R$-algebras. Then $B$ is almost faithfully flat over $A$ if and only if $B^a$ and $\coker(f^a)$ are flat $A^a$-modules.
\end{cor}
\begin{proof}
    By Lemma~\ref{lemma:almost-ff-almost-injective}, it suffices to show that $f^a$ is universally injective if and only if $\coker(f^a)$ is $A^a$-flat. We note that, for any $A^a$-module $M^a$, $\ker(M^a\to M^a\otimes_{A^a}B^a) \simeq \rm{H}^{-1}(M^a\otimes^L_{A^a} \coker(f^a))$. In particular, 
    \[
    \rm{H}^{-1}\left(M^a\otimes^L_{A^a} \coker(f^a)\right)\simeq 0
    \]
    for any $A^a$-module $M^a$ if and only if the functor $-\otimes_{A^a} \coker(f^a)\colon\bf{Mod}_A^a \to \bf{Mod}_A^a$ is exact. In other words, $A^a\to B^a$ is universally injective if and only if $\coker(f^a)$ is flat over $A^a$. 
\end{proof}

Now we define the functor $(-)_{!!}\colon \bf{Alg}_R\to\bf{Alg}_R$. We start by constructing an $R$-algebra structure on $R\oplus A^a_! = R\oplus (\widetilde{\m}\otimes_R A)$ by defining the multiplication as
\[
(r\oplus a)\cdot (r'\oplus a') = (rr') \oplus(ra' + r'a + aa')
\]
and the summation law coordinate-wise. One easily checks that this is a well-defined (unital, commutative) $R$-algebra structure on $R\oplus A^a_!$. We consider the $R$-submodule $I_A$ of $R\oplus A_!$ generated by elements of the form $(mn, -m\otimes n\otimes 1_A)$ for $m$, $n\in \m$. 

\begin{lemma} The $R$-module $I_A \subset R\oplus A^a_!$ is an ideal.
\end{lemma}
\begin{proof}
It suffices to show that, for any element $(r, x\otimes y \otimes a)$ in $R\oplus A^a_!$, the product 
\[
(r\oplus x\otimes y \otimes a)\cdot (mn \oplus -m\otimes n \otimes 1_A)
\]
lies in $I_A$ for any $m$, $n\in \m$. By definition,
\begin{align*}
(r \oplus x\otimes y \otimes a)\cdot (mn \oplus -m\otimes n \otimes 1_A) & = (rmn) \oplus(-rm\otimes n\otimes 1_A + xm\otimes yn \otimes a -xm\otimes yn \otimes a) \\
& = r(mn \oplus -m\otimes n\otimes 1_A) \in I_A. \qedhere
\end{align*}
\end{proof}

\begin{defn} The functor $(-)_{!!}\colon \bf{Alg}_R\to \bf{Alg}_R$ is defined as
\[
A \mapsto (R\oplus A^a_!)/I_A
\]
with the induced $R$-algebra structure. 
\end{defn}

For any $R$-algebra $A$, there is a functorial $R$-algebra homomorphism $R\oplus A^a_! \to A$ defined by
\[
r\oplus(m\otimes n \otimes a)\mapsto r+mna.
\]
Clearly, this homomorphism is zero on $I_A$, so it descends to an $R$-algebra homomorphism $\eta\colon A_{!!} \to A$. 

\begin{lemma}\label{lemma:properties-of-!!}
\begin{enumerate}[label=\textbf{(\arabic*)}]
    \item\label{lemma:properties-of-!!-1} For any $R$-algebra $A$, the natural morphism $\eta\colon A_{!!} \to A$ is an almost isomorphism. 
    \item\label{lemma:properties-of-!!-2} A morphism of $R$-algebras $f\colon A\to B$ is almost injective (as a morphism of $R$-modules) if and only if $f_{!!}\colon A_{!!}\to B_{!!}$ is injective.
    \item\label{lemma:properties-of-!!-3} For any morphism of $R$-algebras $f\colon A\to B$, there is a canonical isomorphism of $A_{!!}$-modules $\coker(f_{!!}) \simeq \coker(f)_!$.
    \item\label{lemma:properties-of-!!-4} The functor $(-)_{!!}\colon \bf{Alg}_R \to \bf{Alg}_R$ commutes with tensor products.
\end{enumerate}
\end{lemma}
\begin{proof}
\ref{lemma:properties-of-!!-1}: We recall that the morphism $A_! \to A$ is almost isomorphism. In particular, it is almost surjective. Thus, $A_{!!} \to A$ is also almost surjective. Now we check almost injectivity. Suppose $\eta(\ov{a})=0$ where $a=r\oplus \sum_{i=1}^k m_i\otimes n_i \otimes a_i\in R\oplus \widetilde{\m}\otimes A$ and $\ov{a}\in A_{!!}$ is the class of $a$ in $A_{!!}$. Then the condition $\eta(\ov{a})=0$ implies that there is an equality
\[
r + \sum_{i=1}^k m_in_ia_i = 0
\]
in $A$. In particular, for every $\e\in \m$, we have $\e r= \sum_{i=1}^k(-m_i)(\e n_ia_i)$ in $A$. Thus, we see that 
\begin{align*}
\e a = & \e r \oplus \sum_{i=1}^k m_i\otimes n_i\otimes \e a_i \\
& = \sum_{i=1}^k (-m_i)(\e n_i a_i)\oplus \sum_{i=1}^k m_i\otimes n_i \e a_i \otimes 1_A\\
& = \sum_{i=1}^k\left( \left(-m_i\right)\left(\e n_i a_i\right)\oplus m_i\otimes \e n_i a_i \otimes 1_A \right)\in I_A.
\end{align*}
Therefore, $\e\ov{a}=0$ for every $\e\in \m$. In particular, $\eta$ is almost injective. \smallskip 

\ref{lemma:properties-of-!!-2} and \ref{lemma:properties-of-!!-3}: Consider a commutative diagram
\[
\begin{tikzcd}
A_{!!} \arrow{r}{f_{!!}} \arrow{d}{\eta_{A}}& B_{!!}\arrow{d}{\eta_{B}}\\
A \arrow{r}{f} & B.
\end{tikzcd}
\]
Since $\eta_{A}$ and $\eta_{B}$ are almost isomorphism, we see that $f$ is almost injective if and only if $f_{!!}$ is almost injective. So we are left to show that $f_{!!}$ is injective if $f$ is almost injective, and $\coker(f_{!!})=\coker(f)_!$. For this, we consider a commutative diagram of short exact sequences
\[
\begin{tikzcd}
0 \arrow{r}& I_{A} \arrow{r}\arrow{d}{\a} & R\oplus A_{!} \arrow{r} \arrow{d}{\rm{Id}\oplus f_!} & A_{!!}\arrow{d}{f_{!!}}\arrow{r} & 0 \\
0 \arrow{r}& I_{B}\arrow{r} & R\oplus B_{!}\arrow{r} & B_{!!} \arrow{r} & 0.
\end{tikzcd}
\]
Clearly, $\a$ is surjective, $\ker (\rm{Id}\oplus f_!) = \ker(f_!)=\ker(f)_!$, and $\coker (\rm{Id}\oplus f_!) = \coker(f_!) =\coker(f)_!$. Thus, the Snake Lemma implies that
\[
\ker(f)_! \to \ker(f_{!!})
\]
is surjective and
\[
\coker(f_{!!}) \to \coker(f)_!
\]
is an isomorphism. Thus $f_{!!}$ is injective if $f$ is almost injective, and $\coker(f_{!!})=\coker(f)_!$. \smallskip

\ref{lemma:properties-of-!!-4}: This is an elementary but pretty tedious computation. We leave it to the interested reader. 
\end{proof}

\begin{cor}\label{cor:!!-mod-equi} For any $R$-algebra $A$, the forgetful functor $\bf{Mod}_{A^a}^* \to \bf{Mod}_{A_{!!}^a}^*$ is an equivalence for $*\in \{``\text{ ''}, \rm{aft}, \rm{afp}, \rm{acoh}\}$.
\end{cor}
\begin{proof}
    For $*=``\text{ ''}$,the claim follows from  Lemma~\ref{lemma:properties-of-!!}~\ref{lemma:properties-of-!!-1}, Corollary~\ref{cor:almost-iso-almost-equi}, and Lemma~\ref{lemma:almost-iso-almost-equi-coh}.
\end{proof}

\begin{cor}\label{cor:!!-faithfully-flat} Let $f\colon A\to B$ be an almost faithfully flat morphism of $R$-algebras. Then $f_{!!} \colon A_{!!} \to B_{!!}$ is faithfully flat.
\end{cor}
\begin{proof}
    Denote by $Q$ the cokernel $f$ as an $A$-module. Then Lemma~\ref{lemma:almost-ff-almost-injective} and Lemma~\ref{lemma:properties-of-!!}\ref{lemma:properties-of-!!-2},\ref{lemma:properties-of-!!-3} ensure that $f_{!!}\colon A_{!!} \to B_{!!}$ is injective and $\coker(f_{!!})=\coker(f)_!$. Now Corollary~\ref{cor:almost-ff-cokernel-flat} and Lemma~\ref{lemma:where-sends-2} applied to $A_{!!}^a\simeq A^a$ imply that $\coker(f_{!!})=\coker(f)_!$ is a flat $A_{!!}$-module. This already implies that $B$ is a flat $A_{!!}$-module as an extension of two flat $A_{!!}$-modules. To see that it is faithfully flat, we note that flatness of $\coker(f_{!!})$ implies that 
    \[
    M \to M\otimes_{A_{!!}} B_{!!}
    \]
    is injective for any $A_{!!}$-module $M$. So $M\otimes_{A_{!!}} B_{!!} \simeq 0$ if and only if $M\simeq 0$. In other words, $B_{!!}$ is a faithfully flat $A_{!!}$-module.
\end{proof}

\begin{warning} The functor $(-)_{!!}$ does not send flat $A$-algebras to flat $A_{!!}$-algebras. See \cite[Remark 3.1.3]{GR}.
\end{warning}

For future reference, we also show that the base change functor interacts especially well with the $\rm{Hom}$-functor in the almost flat situation.

\begin{lemma}\label{flat-base-change-hom} Let $R\to S$ be an almost flat morphism of rings, let $M$ be an almost finitely presented $R$-module, and let $N$ be an $R$-module. Then the natural map
\[
\rm{Hom}_R(M, N) \otimes_R S \to \rm{Hom}_S(M\otimes_R S , N\otimes_R S) 
\]
is an almost isomorphism.
\end{lemma}
\begin{proof}
This follows from the classical $\otimes$-$\rm{Hom}$ adjunction and Lemma~\ref{base-change-hom-2var}.
\end{proof}

\begin{lemma}\label{base-change-hom-derived} Let $R$ be an almost coherent ring, let $R\to S$ be an almost flat map of rings, and let $M\in \bf{D}^-_{acoh}(R)$, $N\in \bf{D}^+(R)$. Then the natural map
\[
\bf{R}\rm{Hom}_R(M, N) \otimes^L_R S \to \bf{R}\rm{Hom}_S(M\otimes^L_R S , N\otimes^L_R S) 
\]
is an almost isomorphism.
\end{lemma}
\begin{proof}
We recall that we always have a canonical isomorphism $\bf{R}\rm{Hom}_R(K, L) \simeq \bf{R}\rm{Hom}_S(K\otimes^L_R S, L)$ for any $K\in \bf{D}^-(R)$ and any $L\in \bf{D}^+(S)$. This implies that it suffices to show that the natural map
\[
\bf{R}\rm{Hom}_R(M, N) \otimes^L_R S \to \bf{R}\rm{Hom}_R(M, N\otimes^L_R S) 
\]
is an almost isomorphism. This follows from Proposition~\ref{base-change-hom-2var-derived}.
\end{proof}

\subsection{Almost faithfully flat descent}

The main goal of this section is to show almost faithfully flat descent for almost modules. \smallskip

For the rest of the section, we fix a ring $R$ with an ideal of almost mathematics $\m$. \smallskip

In this section, for any morphism $A\to B$ of $R$-algebras, we denote the tensor product functor $-\otimes_{A^a} B^a$ simply by 
\[
f^*\colon \bf{Mod}_A^a \to \bf{Mod}_B^a.
\]
In particular, if $A\to B$ is a morphism of $R$-algebras, the canonical ``co-projection'' morphisms $p_{i}\colon B \to B\otimes_A B$ induce morphisms 
\[
p_i^*\colon \bf{Mod}_B^a \to \bf{Mod}_{B\otimes_A B}^a 
\]
for $i\in \{1, 2\}$. The same applies to the ``co-projections'' 
\[
p_{i, j}^*\colon \bf{Mod}_{B\otimes_A B}^a \to \bf{Mod}_{B\otimes_A B\otimes_A B}^a
\]
for $i\neq j \in \{1, 2\}$.

\begin{defn} An {\it almost descent category $\bf{Desc}_{B/A}^a$} for a morphism of $R$-algebras $A \to B$ is a category whose objects are pairs $(M^a, \phi)$, where $M^a\in \bf{Mod}_B^a$ and 
\[
\phi\colon p_1^*(M^a) \to p_2^*(M^a)
\]
in an isomorphism of $(B\otimes_A B)^a$-modules such that $p_{1, 3}^*(\phi) = p_{2, 3}^*(\phi) \circ p_{1, 2}^*(\phi)$. Morphisms between $(M^a, \phi_M)$ and $(N^a, \phi_N)$ are defined to be $B^a$-linear homomorphisms $f\colon M^a\to N^a$ such that the diagram
\[
\begin{tikzcd}
p_1^*(M^a) \arrow{r}{\phi_M}\arrow{d}{p_1^*(f)} & p_2^*(M^a)\arrow{d}{p_2^*(f)}\\
p_1^*(N^a) \arrow{r}{\phi_N} & p_2^*(N^a) 
\end{tikzcd}
\]
commutes. 
\end{defn}

\begin{rmk} Explicitly, an object of the descent category $\bf{Desc}_{B/A}^a$ is a $B^a$-module $M^a$ with a $(B\otimes_{A} B)^a$-linear homomorphism $\phi\colon M^a\otimes_{A^a} B^a \to B^a\otimes_{A^a} M^a$ satisfying the ``cocycle condition''. 
\end{rmk}

There is a natural functor 
\[
\rm{Ind}\colon \bf{Mod}_A^a \to \bf{Desc}_{B/A}^a
\]
that sends $M^a$ to $f^*(M^a)=M^a\otimes_{A^a} B^a$ with the canonical identification $\phi\colon p_1^*f^*\left(M^a\right) \simeq p_2^* f^*\left(M^a\right)$ coming from the equality $f\circ p_1 = f\circ p_2$. \smallskip

To define the functor in the other direction, we note that we have the natural $B^a$-module morphisms $\iota_i\colon M^a \to p_i^*\left(M^a\right)$ for $i\in \{1, 2\}$. Explicitly, they are defined as morphisms induced by $\iota_1(m)=m\otimes 1$ and $\iota_2(m)=1\otimes m$. Therefore, given a descent data $(M^a, \phi) \in \bf{Desc}_{B/A}^a$, we can define an $A^a$-module
\[
\ker(M^a, \phi) \coloneqq \ker(M^a\xr{i_1-\phi^{-1}i_2} M^a \otimes_{A^a} B^a)
\]
that is functorial in $\bf{Desc}_{B/A}^a$. Therefore, this defines a functor
\[
\ker \colon \bf{Desc}_{B/A}^a \to \bf{Mod}_A^a.
\]

We show that $\ker$ and $\rm{Ind}$ are quasi-inverse to each other and induce an equivalence between $\bf{Desc}_{B/A}^a$ and $\bf{Mod}_A^a$ for an almost faithfully flat morphism $f\colon A \to B$.

\begin{thm}\label{thm:faithuflly-flat descent} Let $f\colon A \to B$ be an almost faithfully flat morphism. Then 
\[
\rm{Ind}\colon \bf{Mod}_A^a \to \bf{Desc}_{B/A}^a
\]
is an equivalence, and its quasi-inverse is given by the functor $\ker\colon \bf{Desc}_{B/A}^a\to \bf{Mod}_A^a$.
\end{thm}
\begin{proof}
    Corollary~\ref{cor:!!-mod-equi} and Corollary~\ref{cor:!!-faithfully-flat} imply that we may replace $f$ with $f_{!!}$ to assume that $f$ is faithfully flat. Then the claim follows from the classical faithfully flat descent (see \cite[Theorem 6.1/4]{Neron}) and the observation that the classical versions of $\rm{Ind}$ and $\ker$ carry almost isomorphisms to almost isomorphisms. 
\end{proof}

On a similar note, we show that the Amitsur complex for an almost faithfully flat morphism is acyclic. 

\begin{lemma}\label{lemma:amitsur-acyclic} Let $f\colon A \to B$ be an almost faithfully flat morphism of $R$-algebras, and $M\in \bf{Mod}_B^a$. Then the Amitsur complex 
\[
0 \to M^a \to M^a\otimes_{A^a} B^a \to M^a\otimes_{A^a} B^a \otimes_{A^a} B^a \to \dots
\]
is an exact complex of $\bf{Mod}_B^a$-modules (see the discussion around \cite[\href{https://stacks.math.columbia.edu/tag/023K}{Tag 023K}]{stacks-project} for the precise definition of differentials in this complex).
\end{lemma}
\begin{proof}
    Corollary~\ref{cor:!!-mod-equi} and Corollary~\ref{cor:!!-faithfully-flat} imply that we may replace $f$ with $f_{!!}$ to assume that $f$ is faithfully flat. Then the claim follows from \cite[\href{https://stacks.math.columbia.edu/tag/023M}{Tag 023M}]{stacks-project}. 
\end{proof}

Now we show that some properties of $A^a$-modules can be verified after a faithfully flat base change.

\begin{lemma}\label{almost-flat-descent} Let $f\colon A \to B$ be an almost faithfully flat morphism of $R$-algebras, and let $M^a$ be an $A^a$-module. Then $M^a$ is an almost finitely generated (resp. almost finitely presented) $A^a$-module if and only if $M^a\otimes_{A^a} B^a$ is an almost finitely generated (resp. almost finitely presented) $B^a$-module. 
\end{lemma}
\begin{proof} 

Corollary~\ref{cor:!!-mod-equi} and Corollary~\ref{cor:!!-faithfully-flat} imply that we may replace $f$ with $f_{!!}$ to assume that $f$ is a faithfully flat morphism. Then a standard argument reduces the questions to the case of an honest $A$-module $M$, i.e., we show that an $A$-module $M$ is almost finitely generated (resp. almost finitely presented) if so is the $B$-module $M\otimes_{A} B$. \smallskip

We start with the almost finitely generated case. So we assume that $M\otimes_A B$ is almost finitely generated over $B$ and wish to show that $M$ is almost finitely generated over $A$. Our assumption implies that, for any $\e \in \m$, we can choose a morphism $g \colon B^n \to M\otimes_A B$ such that $\e(\coker g)=0$. Let us consider the standard basis $e_1, \dots, e_n$ of $B^n$, and we write 
\[
g(e_i) =\sum_j m_{i,j}\otimes b_{i,j} \text{ with } m_{i,j}\in M,  b_{i,j}\in B.
\]
We define the $A$-module $F$ as the finite free $A$-module with the basis $e_{i,j}$. Then we define the morphism 
\[
h\colon F \to M
\]
as the unique $A$-linear homomorphism with $h(e_{i,j})=m_{i,j}$. It is easy to see that $\e(\coker (h\otimes_A B))=0$. Since $f$ is faithfully flat, this implies that $\e(\coker h)=0$. We conclude that $M$ is almost finitely generated as $\e$ was an arbitrary element of $\m$. \smallskip

Now we deal with the almost finitely presented case. We pick some finitely generated ideal $\m_0\subset \m$, and another finitely generated ideal $\m_1\subset \m$ such that $\m_0 \subset \m_1\m$. We try to find a three-term complex
\[
A^m \xr{g} A^n \xr{f} M
\]
such that $\m_0(\coker f)=0$ and $\m_0(\ker f)\subset \operatorname{Im} g$. \smallskip

The almost finitely generated case established above implies that $M$ is almost finitely generated. In particular, we have some morphism
\[
A^n \xr{f} M
\]
such that $\m_1(\coker f)=0$, thus $\m_1(\coker (f\otimes_A B))=0$ as well. Therefore, we can apply Lemma~\ref{claim} to find a homomorphism $g'\colon B^m \to B^n$ such that $\m_0(\ker (f\otimes_A B))\subset \operatorname{Im}(g')$ and $(f\otimes_A B)\circ g'=0$. This implies that $g'$ lands inside $\ker(f\otimes_A B)=\ker(f)\otimes_A B$ due to $A$-flatness of $B$. \smallskip

Now we do the same trick as above: we write 
\[
g(e_i)=\sum_j m_{i,j}\otimes b_{i,j} \text{ with } m_{i,j}\in \ker(f),  b_{i,j}\in B,
\]
we define an $R$-module $F$ as a finite free $A$-module with a basis $e_{i,j}$, and then we define the morphism
\[
g\colon F \to \ker(f)
\]
as the unique $A$-linear morphism such that $g(e_{i,j})=m_{i,j}$. Then we see that $\m_0(\ker(f\otimes_A B))\subset \operatorname{Im}(g\otimes_A B)$. Since $B$ is faithfully flat, we conclude that $\m_0(\ker f)\subset \operatorname{Im}(g)$ as well. This shows that a three-term complex
\[
F\xr{g} A^n \xr{f} M
\]
does the job. Therefore, $M$ is an almost finitely presented $A$-module.
\end{proof}

\begin{cor}\label{almost-flat-descent-coh} Let $f\colon A \to B$ be an almost faithfully flat morphism of $R$-algebras, let $M^a$ be an $A^a$-module. Suppose that $M^a\otimes_{A^a} B^a$ is an almost coherent $B^a$-module. Then so is $M^a$.
\end{cor}
\begin{proof}
This follows directly from Lemma~\ref{lemma:almost-iso-almost-equi-coh} and Lemma~\ref{almost-flat-descent}.
\end{proof}

\begin{lemma}\label{lemma:almost-flat-descent-flat} Let $f\colon A \to B$ be an almost faithfully flat morphism of $R$-algebras, and let $M^a$ be an $A^a$-module. Then $M^a$ is a flat (resp. faithfully flat) $A^a$-module if and only if $M^a\otimes_{A^a} B^a$ is a flat (resp. faithfully flat) $B^a$-module. 
\end{lemma}
\begin{proof} 
The classical proof works verbatim in the almost world. We leave the details to the reader. 
\end{proof}

\subsection{(Topologically) Finite type $K^+$-algebras}\label{section:Examples}

This section is devoted to the proof that (topologically) finite type algebras over a perfectoid valuation ring $K^+$ are almost noetherian. We refer to Appendix~\ref{Section:perfectoid-valuation} for the relevant background on perfectoid valuation rings. \smallskip

For the rest of the section, we fix a perfectoid valuation ring $K^+$ (see Definition~\ref{defn:valuation-perfectoid}) with perfectoid fraction field $K$, associated rank-$1$ valuation ring $\O_K=K^\circ$ (see Remark~\ref{rmk:perfectoid-microbial}), and ideal of topologically nilpotent elements $\m=K^{\circ\circ} \subset K^+$. Lemma~\ref{lemma:almost-setup} ensures that $\m$ is flat over $K^+$ and $\widetilde{\m}\simeq \m^2=\m$. Therefore, it makes sense to do almost mathematics with respect to the pair $(K^+, \m)$. In what follows, we always do almost mathematics on $K^+$-modules with respect to this ideal.

\begin{warning1}\label{warning:not-maximal-ideal-appendix} The ideal $\m \subset K^+$ is not the maximal ideal of $K^+$. Instead, it is the maximal ideal of the associated rank-$1$ valuation ring $\O_K$.
\end{warning1}

\begin{lemma}\label{lemma:perfectoid-rank-1} Let $K^+$ be a perfectoid valuation ring. Then the natural inclusion $\iota\colon K^+ \to \O_K$ is an almost isomorphism.
\end{lemma}
\begin{proof}
    Clearly, the map $\iota \colon K^+ \to \O_K$ is injective, so it suffices to show that its cokernel is almost zero, i.e. annihilated by any $\e\in \m$. Pick an element $x\in \O_K$, then $\e x\in \m \subset K^+$. Therefore, we conclude that $\e(\coker \iota)=0$ finishing the proof. 
\end{proof}

The first main result of this section is that any (topologically) finite type algebra over $K^+$ is almost noetherian. 

\begin{lemma}\label{lemma:tate-algebra-almost-noetherian} Let $K^+$ be a perfectoid valuation ring, and $n\geq 0$ an integer. Then the Tate algebra $K^+\langle T_1, \dots, T_n\rangle$ is almost noetherian. 
\end{lemma}
\begin{proof}
    First, we note that $\O_K \langle T_1, \dots, T_n\rangle \simeq K^+\langle T_1, \dots, T_n\rangle \otimes_{K^+} \O_K$. Therefore, Lemma~\ref{lemma:perfectoid-rank-1} implies that the natural morphism
    \[
    K^+ \langle T_1, \dots, T_n\rangle \to \O_K \langle T_1, \dots, T_n\rangle
    \]
    is an almost isomorphism. So Corollary~\ref{cor:almost-noetherian-iso} ensures that it suffices to show that $\O_K\langle T_1, \dots, T_n\rangle$ is almost noetherian. \smallskip
    
    Pick any ideal $I \subset \O_K \langle T_1, \dots, T_n\rangle = K\langle T_1, \dots, T_n\rangle^\circ$ and $0\neq \e \in \m$. Now \cite[Satz 5.1]{Kiehl-finiteness} (or \cite[Lemma 6.4/5]{B}) applied to $B=K\langle T_1, \dots, T_n\rangle$, $E=\O_K\langle T_1, \dots, T_n\rangle$, $E' =I$, and $\alpha=|\e|_K$ guarantees that there is a finite submodule $E'' \subset I$ such that $\e I \subset E''$. Since $\e$ was an arbitrary element of $\m$, we conclude that $I$ is indeed almost finitely generated. 
\end{proof}

\begin{cor}\label{cor:almost-noetherian-mod-p} Let $K^+$ be a perfectoid valuation ring, $\varpi\in \m$, and $n\geq 0$ an integer. Then the polynomial algebra $(K^+/\varpi^m)[T_1, \dots, T_n]$ is almost noetherian for any $m\geq 1$. 
\end{cor}
\begin{proof}
    It easily follows from Lemma~\ref{lemma:tate-algebra-almost-noetherian}, Corollary~\ref{cor:almost-noetherian-aft=afpr}, and Lemma~\ref{extension}. 
\end{proof}

\begin{thm}\label{thm:top-ft-almost-noetherian} Let $K^+$ be a perfectoid valuation ring, and $A$ a topologically finite type $K^+$-algebra. Then $A$ is almost noetherian.
\end{thm}
\begin{proof}
    Since $A$ is topologically finite type over $K^+$, there exists a surjection
    \[
    f\colon K^+\langle T_1, \dots, T_n\rangle \to A \to 0.
    \]
    Pick an ideal $I \subset A$ and consider its preimage $J=f^{-1}(I)$. Then $J$ is almost finitely generated over $K^+\langle T_1, \dots, T_n\rangle$ by Lemma~\ref{lemma:tate-algebra-almost-noetherian}. Therefore, Lemma~\ref{main-almost}\ref{main-almost-1} ensures that $I$ is almost finitely generated over $K^+\langle T_1, \dots, T_n\rangle$. Finally, Lemma~\ref{extension} ensures that $I$ is almost finitely generated over $A$. 
\end{proof}

Now we are going to show that any finite type $K^+$-algebra is almost noetherian. Before doing this, we need a couple of preliminary lemmas. 

\begin{lemma}\label{lemma:bounded-torsion} Let $R$ be a rank-$1$ valuation ring with a non-zero topologically nilpotent element $\varpi \in R$, and $M$ a finite $R[T_1, \dots, T_n]$-module. Then $M[\varpi^\infty]=M[\varpi^c]$ for some $c\geq 0$.
\end{lemma}
\begin{proof}
    The $R[T_1, \dots, T_n]$-module $M'\coloneqq M/M[\varpi^\infty]$ is finitely generated. Furthermore, $M'$ is $R$-flat because it is torsion-free (and $R$ is a valuation ring). Therefore, \cite[\href{https://stacks.math.columbia.edu/tag/053E}{Tag 053E}]{stacks-project} ensures that $M'$ is finitely presented over $R[T_1, \dots, T_n]$. Thus, we conclude that $M[\varpi^\infty]$ is finitely generated. In particular, $M[\varpi^\infty]=M[\varpi^c]$ for some $N$. 
\end{proof}

\begin{lemma}\label{lemma:similar-topologies} Let $R$ be a rank-$1$ valuation ring with a non-zero topologically nilpotent element $\varpi \in R$, $M$ a finite $R[T_1, \dots, T_n]$-module, and $N\subset M$ an $R[T_1, \dots, T_n]$-submodule. Then there is $c$ such that 
\[
N\cap \varpi^{m+c}M = \varpi^m(N\cap \varpi^c M)
\]
for every $m\geq 0$. 
\end{lemma}
\begin{proof}
    Lemma~\ref{lemma:bounded-torsion} ensures that there $c$ such that $(M/N)[\varpi^\infty]=(M/N)[\varpi^c]$. Therefore, 
    \cite[Lemma 0.8.2.14]{FujKato} guarantees that, indeed,
    \[
        N\cap \varpi^{m+c}M = \varpi^m(N\cap \varpi^c M)
    \]
    for every $m\geq 0$. 
\end{proof}

\begin{lemma}\label{lemma:polynomials-almost-noetherian} Let $K^+$ be a perfectoid valuation ring, and $n\geq 0$ an integer. Then the polynomial algebra $K^+[T_1, \dots, T_n]$ is almost noetherian. 
\end{lemma}
\begin{proof}
    Similarly to the proof of Lemma~\ref{lemma:tate-algebra-almost-noetherian}, it suffices to treat the case $K^+=\O_K$ a perfectoid valuation ring of rank-$1$ with a pseudo-uniformizer $\varpi$. \smallskip
    
    Now we fix an ideal $I \subset A\coloneqq \O_K[T_1, \dots, T_n]$ and wish to show that $I$ is almost finitely generated. Recall that the polynomial algebra $K[T_1, \dots, T_n]$ is noetherian by Hilbert's Nullstellensatz. Therefore, the ideal 
    \[
    I\left[\frac{1}{\varpi}\right] \subset K[T_1, \dots, T_n]
    \]
    is finitely generated. So we can choose a finitely generated sub-ideal $J\subset I$ such that any element of $I/J$ is annihilated by a power of $\varpi$, i.e. $(I/J)[\varpi^\infty]=I/J$. Clearly $I/J$ is a submodule of a finite $A$-module $A/J$, so Lemma~\ref{lemma:bounded-torsion} easily implies that 
    \[
    I/J=(I/J)[\varpi^\infty]=(I/J)[\varpi^c]
    \]
    for some $c\geq 0$. In other words, $\varpi^c I \subset J$. Now we use Lemma~\ref{lemma:similar-topologies} to get an integer $c'$ such that 
    \[
    I \cap \varpi^{c'} A \subset \varpi^c I \subset J.
    \]
    We note that $I/(I\cap \varpi^{c'} A)$ is an ideal in $A/\varpi^{c'}A$, and therefore it is almost finitely generated over $A/\varpi^{c'}A$ by Corollary~\ref{cor:almost-noetherian-mod-p}. Lemma~\ref{extension} guarantees that it is also almost finitely generated over $A$. \smallskip
    
    The inclusion $I\cap \varpi^{c'} A \subset J$ implies that $I/J$ is a quotient of an almost finitely generated $A$-module $I/(I\cap \varpi^{c'}A)$, and so is also almost finitely generated. Finally, the short exact sequence
    \[
    0 \to J \to I \to I/J \to 0
    \]
    and Lemma~\ref{main-almost}\ref{main-almost-2} imply that $I$ is almost finitely generated as well.
\end{proof}

\begin{thm}\label{thm:finite-type-almost-noetherian} Let $K^+$ be a perfectoid valuation ring, and $A$ a finite type $K^+$-algebra. Then $A$ is almost noetherian.
\end{thm}
\begin{proof}
    It follows from Lemma~\ref{lemma:polynomials-almost-noetherian} similar to how Theorem~\ref{thm:top-ft-almost-noetherian} follows from Lemma~\ref{lemma:tate-algebra-almost-noetherian}.
\end{proof}

\subsection{Almost finitely generated modules over adhesive rings}\label{adhesive-almost}

This section discusses some basic aspects of almost finitely generated modules over adhesive rings. The results of this Section would be crucial in defining and verifying certain good properties of adically quasi-coherent, almost coherent sheaves on ``good'' formal schemes in Section~\ref{acoh-sheaves-formal}. One of the essential ingredients that we will need later is the ``weak'' version of the Artin--Rees Lemma (Lemma~\ref{Artin-Rees}) and Lemma~\ref{completion-finitely-generated}. Recall that these properties are already known for finite modules over adhesive rings. This is explained in a beautiful paper \cite{FGK}. The main goal of this section is to extend these results to the case of almost finitely generated modules. \smallskip

That being said, let us introduce the setup for this section. We start with the definition of an adhesive ring:

\begin{defn}\label{defn:adhesive}\cite[Definition 7.1.1]{FGK} An adically topologized ring $R$ endowed with the adic topology defined by a finitely generated ideal $I \subset R$ is said to be {\it ($I$-adically) adhesive} if it is Noetherian outside\footnote{By definition, this means that the scheme $\Spec A \setminus V(I)$ is noetherian.} $I$ and satisfies the following condition: for any finitely generated $R$-module $M$, its $I^{\infty}$-torsion part $M[I^{\infty}]$ is finitely generated.
\end{defn}

\begin{rmk} Following the convention of \cite{FGK} we do not require a ring $R$ with adic topology to be either $I$-adically complete or separated.
\end{rmk}

\begin{setup}\label{set-up1} We fix an $I$-adically adhesive  ring $R$ with an ideal $\m$ such that $I\subset \m$, $\m^2=\m$ and $\widetilde{\m}\coloneqq \m\otimes_R \m$ is flat. We always do almost mathematics with respect to the ideal $\m$. 
\end{setup}

The main example of an adhesive ring is a (topologically) finitely presented algebra over a complete microbial valuation ring. This follows from \cite[Proposition 7.2.2]{FGK} and \cite[Theorem 7.3.2]{FGK}. For example, any topologically finitely presented algebra over a complete rank-$1$ valuation ring is adhesive.

\begin{lemma}\label{saturated} Let $R$ be as in Setup~\ref{set-up1}, and let $M$ be an $I$-torsionfree almost finitely generated module. Then $M$ is almost finitely presented. Similarly, any saturated submodule\footnote{A submodule $N\subset M$ is {\it saturated} if $M/N[I^{\infty}]=0$.} of an almost finitely generated $R$-module is almost finitely generated. 
\end{lemma}
\begin{proof}
As $M$ is almost finitely generated, we can find a finitely generated submodule $N\subset M$ that contains $\m_0M$ for a choice of a finitely generated ideal $\m_0 \subset \m$. Since $N$ is a submodule of $M$, it is itself $I$-torsionfree. Then \cite[Proposition 7.1.2]{FGK} shows that $N$ is finitely presented. Then Lemma~\ref{almost-finitely-presented}\ref{almost-finitely-presented-2} implies that $M$ is almost finitely presented. 

Now let $M$ be an almost finitely generated $R$-module, and let $M'\subset M$ be a saturated submodule. Then $M/M'$ is almost finitely generated by Lemma~\ref{main-almost}\ref{main-almost-1} and it is $I$-torsionfree. Therefore, it is almost finitely presented by the argument above. Then Lemma~\ref{main-almost}\ref{main-almost-3} guarantees that $M'$ is almost finitely generated. 
\end{proof}

\begin{lemma}\label{bounded-torsion} Let $R$ be as in Setup~\ref{set-up1}, and let $M$ be an almost finitely generated $R$-module. Then the $I^{\infty}$-torsion module $M[I^{\infty}]$ is bounded (i.e. there is an integer $n$ such that $M[I^n]=M[I^{\infty}]$).
\end{lemma}
\begin{proof}
Since $M$ is almost finitely generated and the ideal $I\subset \m$ is finitely generated, we conclude that there exists a finitely generated submodule $N \subset M$ such that $IM \subset N$. Then $I(M[I^{\infty}]) \subset N[I^{\infty}]$, and $N[I^{\infty}]$ is finitely generated by adhesiveness of the ring $R$. In particular, there is an integer $n$ such that $N[I^{\infty}]$ is annihilated by $I^n$. This implies that any element of $M[I^{\infty}]$ is annihilated by $I^{n+1}$.
\end{proof}

\begin{lemma}\label{Artin-Rees} Let $R$ be as in Setup~\ref{set-up1}, and let $M$ be an almost finitely generated $R$-module. Suppose that $N\subset M$ is a submodule of $M$. For any integer $n$, there is an integer $m$ such that $N \cap I^mM \subset I^nN$. In particular, the induced topology on the module $N$ coincides with the $I$-adic one.
\end{lemma}
\begin{proof}
If $M$ is finitely generated, then this is \cite[Theorem 4.2.2]{FGK}. In general, we use the definition of almost finitely generated module to find a submodule $M' \subset M$ such that $M'$ is finitely generated and $IM\subset M'$. We define $N'\coloneqq N\cap M'$ as the intersection of those modules. Then the established ``weak'' form of the Artin--Rees Lemma for finitely generated $R$-modules provides us with an integer $m$ such that $N'\cap I^mM' \subset I^nN'$. In particular, we have
\[
I^{m+1}M \cap N' \subset I^mM' \cap N' \subset I^nN' \subset I^nN.
\]
Then we conclude that 
\[
I^{m+2}M \cap N \subset I^{m+1}M \cap M' \cap N \subset I^{m+1}M \cap N' \subset I^nN.
\]
Since $n$ was arbitrary, we conclude the claim. 
\end{proof}

\begin{lemma}\label{completion-finitely-generated} Let $R$ be as in Setup~\ref{set-up1}, and let $M$ be an almost finitely generated $R$-module. Then the natural morphism $M\otimes_R \widehat{R} \to \widehat{M}$ is an isomorphism. In particular, any almost finitely generated module over a complete adhesive ring is complete.
\end{lemma}
\begin{proof}
We know that the claim holds for finitely generated modules by \cite[Proposition 4.3.4]{FGK}. Now we deal with the almost finitely generated case. We choose a finitely generated submodule $N\subset M$ such that $IM \subset N$. Lemma~\ref{Artin-Rees} implies that the induced topology on $N$ coincides with the $I$-adic topology on $N$. Thus the short exact sequence 
\[
0\to N \to M \to M/N \to 0
\]
remains exact after completion. Since $R \to \widehat{R}$ is flat by \cite[Proposition 4.3.4]{FGK}, we conclude that we have a morphism of short exact sequences
\[
\begin{tikzcd}
0 \arrow{r} & N\otimes_R \wdh{R} \arrow{r} \arrow{d}{\varphi_{N}} & M\otimes_R \wdh{R} \arrow{r} \arrow{d}{\varphi_M} & (M/N)\otimes_R \wdh{R} \arrow{r} \arrow{d}{\varphi_{M/N}} & 0 \\
0 \arrow{r} & \wdh{N} \arrow{r}  & \wdh{M} \arrow{r}& \wdh{M/N} \arrow{r} & 0
\end{tikzcd}
\]
Note that $\varphi_N$ is an isomorphism as $N$ is finitely generated, and $\varphi_{M/N}$ is isomorphism since it is an $I$-torsion module so $M/N\simeq (M/N)\otimes_R \wdh{R} \simeq \wdh{M/N}$. The five-lemma implies that $\varphi_M$ is an isomorphism as well.
\end{proof}

\begin{cor}\label{almost-coh-derived-complete} Let $R$ be as in Setup~\ref{set-up1}, and let $M\in \bf{D}_{acoh}(R)$. Suppose that $R$ is $I$-adically complete. Then $M$ is $I$-adically derived complete\footnote{Look at \cite[\href{https://stacks.math.columbia.edu/tag/091N}{Tag 091N}]{stacks-project} for the definition of derived completeness (or Definition~\ref{defn:derived-complete} in case of a principal ideal $I$).}.
\end{cor}
\begin{proof}
First of all, we note that \cite[\href{https://stacks.math.columbia.edu/tag/091P}{Tag 091P}]{stacks-project} implies that $M$ is derived complete if and only if so are $\rm{H}^i(M)$ for any integer $i$. So it suffices to show that any almost coherent $R$-module is derived complete. Lemma~\ref{completion-finitely-generated} gives that any such module is classically complete, and \cite[\href{https://stacks.math.columbia.edu/tag/091T}{Tag 091T}]{stacks-project} ensures that any classically complete module is derived complete. 
\end{proof}

\subsection{Modules over topologically finite type $K^+$-algebras}

The main goal of this section is to show that almost coherentness of derived complete modules over a topologically finite type $K^+$-algebras can be checked modulo the pseudo-uniformizer.\smallskip

For the rest of the section we fix a valuation perfectoid ring $K^+$ (see Definition~\ref{defn:valuation-perfectoid}) with perfectoid fraction field $K$, associated rank-$1$ valuation ring $\O_K=K^\circ$ (see Remark~\ref{rmk:perfectoid-microbial}), and ideal of topologically nilpotent elements $\m=K^{\circ\circ} \subset K^+$ with a pseudo-uniformizer $\varpi \in \m$ as in Lemma~\ref{lemma:roots-of-pseudounformizer} (in particular, $\m = \bigcup_n \varpi^{1/p^n}K^+$). Lemma~\ref{lemma:almost-setup} ensures that $\m$ is flat over $K^+$ and $\widetilde{\m}\simeq \m^2=\m$. Therefore, it makes sense to do almost mathematics with respect to the pair $(K^+, \m)$. In what follows, we always do almost mathematics on $K^+$-modules with respect to this ideal. 

\begin{lemma}\label{lemma:check-almost-coh-mod-p-discrete} Let $R$ be a topologically finite type $K^+$-algebra, and $M$ an $R$-module that is $\varpi$-adically derived complete. Suppose that $M/\varpi M$ is almost coherent, then $M$ is almost coherent as well. 
\end{lemma}
\begin{proof}
    Theorem~\ref{thm:top-ft-almost-noetherian} ensures that $R$ is almost noetherian, and so Corollary~\ref{cor:almost-noetherian-almost-fg=almost-coh} implies that it suffices to check that $M$ is almost finitely generated. Recall that $\m= \bigcup_{n}\varpi^{1/p^n}K^+$ for a pseudo-uniformizer $\varpi$ as in Lemma~\ref{lemma:roots-of-pseudounformizer}. \smallskip
    
    The assumption on $M$ says that $M/\varpi M$ is almost coherent. Therefore, there is a morphism
    \[
    \ov{g} \colon (R/\varpi R)^{c} \to M/\varpi M
    \]
    such that $\varpi^{1/p}(\coker \ov{g})=0$. We denote its cokernel by $\ov{Q}\coloneqq \coker(\ov{g})$. Now we lift $\ov{g}$ to a morphism
    \[
    g\colon R^c \to M
    \]
    and denote is cokernel by $Q\coloneqq \coker(g)$. \smallskip
    
    {\it Step $1$: $Q$ is annihilated by $\varpi^{1/p}$.} Suppose that $\varpi^{1/p}Q\neq 0$, so there is $x_0\in Q$ such that $\varpi^{1/p}x_0 \neq 0$. Firstly, we note that $Q/\varpi \simeq \ov{Q}$ is annihilated by $\varpi^{1/p}$, so 
    \[
    \varpi^{1/p}x_0 = \varpi x_1 = \varpi^{1-1/p}(\varpi^{1/p}x_1)
    \]
    Now we apply the same thing to $x_1$ to get 
    \[
    \varpi^{1/p}x_0 = \varpi^{1-1/p}(\varpi^{1/p}x_1)=(\varpi^{1-1/p})^2(\varpi^{1/p}x_2).
    \]
    Continue the process to get a sequence of elements $x_n \in Q$ such that 
    \[
    \varpi^{1-1/p}(\varpi^{1/p}x_n)=\varpi^{1/p}x_{n-1}.
    \]
    The sequence $\{\varpi^{1/p}x_i\}$ gives an element of 
    \[
    T^0(Q, \varpi^{1-1/p}) \coloneqq \lim_n ( \dots \xr{\varpi^{1-1/p}} Q \xr{\varpi^{1-1/p}} Q) 
    \]
    that is non-trivial because $\varpi^{1/p}x_0\neq 0$. Now we note that $R^c$ is derived $\varpi$-adically complete since $R$ is classically $\varpi$-adically complete by \cite[Corollary 7.3/9]{B} and any classically complete module is derived complete by \cite[\href{https://stacks.math.columbia.edu/tag/091T}{Tag 091T}]{stacks-project}. Therefore, $Q$ is $\varpi$-adically derived complete derived complete as a cokernel of derived complete modules (see \cite[\href{https://stacks.math.columbia.edu/tag/091U}{Tag 091U}]{stacks-project}). Now \cite[\href{https://stacks.math.columbia.edu/tag/091S}{Tag 091S}]{stacks-project}, Remark~\ref{rmk:our-defn-stacks-project}, and 
    \cite[\href{https://stacks.math.columbia.edu/tag/091Q}{Tag 091Q}]{stacks-project} imply that $T^0(Q, \varpi^{1-1/p})$ must be zero leading to the contradiction.\smallskip 

    {\it Step~$2$: $M$ is almost coherent.} Note that $\ov{Q}\simeq Q/\varpi Q$ and $Q$ is $\varpi^{1/p}$-torsion, so $\ov{Q}\simeq Q$. We know that $\ov{Q}$ is almost finitely generated over $R/\varpi R$ because it is a quotient of an almost finitely generated module $M/\varpi M$. Therefore, $Q\simeq \ov{Q}$ is almost finitely generated over $R$ by Lemma~\ref{extension}. Now $M$ is an extension of a finite $R$-module $\rm{Im}(g)$ by an almost finitely generated $R$-module $Q$, so it is also almost finitely generated by Lemma~\ref{main-almost}~\ref{main-almost-2}. In particular, it is almost coherent since $R$ is almost noetherian. 
\end{proof}

\begin{thm}\label{thm:check-almost-coh-mod-p} Let $R$ be a topologically finite type $K^+$-algebra, and $M\in \bf{D}(R)$ a $\varpi$-adically derived complete complex. Suppose that $[M/\varpi] \in \bf{D}^{[c,d]}_{acoh}(R/\varpi)$, then $M\in \bf{D}^{[c,d]}_{acoh}(R)$.
\end{thm}
\begin{proof}
    Lemma~\ref{lemma:derived-complete-ampl-mod-p} guarantees that $M\in \bf{D}^{[c,d]}(R)$, so we only need to show that cohomology groups of $M$ are almost coherent over $R$. \smallskip
    
    We argue by induction on $d-c$. If $c=d$, then $\rm{H}^d(M)/\varpi \simeq \rm{H}^d([M/\varpi])$ is almost coherent. Therefore, $M\simeq \rm{H}^d(M)[-d]$ is almost coherent by Lemma~\ref{lemma:check-almost-coh-mod-p-discrete}. \smallskip
    
    If $d>c$, we consider an exact triangle
    \[
    \tau^{\leq d-1} M \to M \to \rm{H}^d(M)[-d].
    \]
    We see that both $\tau^{\leq d-1} M$ and $\rm{H}^d(M)$ are derived complete by \cite[\href{https://stacks.math.columbia.edu/tag/091P}{Tag 091P}]{stacks-project} and \cite[\href{https://stacks.math.columbia.edu/tag/091S}{Tag 091S}]{stacks-project}. Moreover, we know that $\rm{H}^d(M)/\varpi \simeq \rm{H}^d([M/\varpi])$ is almost coherent. Therefore, $\rm{H}^d(M)$ is almost coherent by Lemma~\ref{lemma:check-almost-coh-mod-p-discrete}. Finally, 
    \[
    [\tau^{\leq d-1}M/\varpi] \simeq \rm{cone}\left([M/\varpi] \to [\rm{H}^d(M)/\varpi][-d]\right)[1]
    \]
    is a (shifted) cone of a morphism in $\bf{D}^b_{acoh}(R/\varpi)$, therefore, $[\tau^{\leq d-1}M/\varpi]$ also lies in $\bf{D}^b_{acoh}(R/\varpi)$. By the induction hypothesis, we conclude that $\tau^{\leq d-1} M\in \bf{D}^{[c,d-1]}_{acoh}(R)$. So $M\in \bf{D}^{[c, d]}_{acoh}(R)$.
\end{proof}

\begin{cor}\label{cor:check-almost-coh-mod-p-2} Let $R$ be a topologically finite type $K^+$-algebra, and $M\in \bf{D}(R)$ a $\varpi$-adically derived complete complex. Suppose that $[M^a/\varpi] \in \bf{D}^{[c,d]}_{acoh}(R/\varpi)^a$, then $M^a\in \bf{D}^{[c,d]}_{acoh}(R)^a$.
\end{cor}
\begin{proof}
    Note that $\m\otimes M$ is derived complete by Lemma~\ref{lemma:derived-complete-after-!}. So the claim follows from Theorem~\ref{thm:check-almost-coh-mod-p} applied to $\m\otimes M$.
\end{proof}

\section{Almost mathematics on ringed sites}\label{almost-sheaves}

The main goal of this Chapter is to ``globalize'' the results of Chapter \ref{almost-commutative-algebra}. The two main cases of interest are almost coherent sheaves on schemes and ``good'' formal schemes. In order to treat those cases uniformly, we define the notion of almost sheaves in the most general set-up of ringed sites and check their basic properties. This is the content of Section~\ref{almost-locally-ringed}. Sections~\ref{acoh-sheaves} and \ref{acoh-sheaves-formal} are devoted to establishing the foundations of almost coherent sheaves on schemes and formal schemes, respectively. In particular, we show that the notion of almost finitely generated (resp. presented, resp. coherent) module globalizes well on schemes and some ``good'' formal schemes. Then we discuss the derived category of almost sheaves and various functors on the derived categories of almost sheaves. 

\subsection{The category of $\O_X^a$-modules}\label{almost-locally-ringed}

We start this section by fixing a ring $R$ with an ideal $\m$ such that $\m=\m^2$ and $\widetilde{\m}=\m\otimes_R \m$ is $R$-flat. We always do almost mathematics with respect to this ideal. The main goal of this section is to globalize the notion of almost mathematics to the case of ringed $R$-sites. \smallskip

In this section, we fix an {\it ringed $R$-site} $(X, \O_X)$, i.e., a ringed site $(X, \O_X)$ where $\O_X$ is a sheaf of $R$-algebras on $X$. Note that any ringed site $(X, \O_X)$ is, in particular, a ringed $\O_X(X)$-site. The main goal of this section is to develop foundations of almost mathematics on ringed $R$-sites.\smallskip

We note that, on each open $U\in X$, it makes sense to speak of almost $\O_X(U)$-modules with respect to the ideal $\m\O_X(U)$; we refer to Lemma~\ref{base-change} for the details. In what follows, we extend the definition of almost modules to the category of $\O_X$-modules.

\begin{defn} Let $(X, \O_X)$ be a ringed $R$-site, and let $\F$ be any $\O_X$-module. Then the {\it sheaf of almost section} $\widetilde{\m}\otimes \F$ is the sheafification of the presheaf defined via the formula
\[
U \mapsto \widetilde{\m}\otimes_{R} \F(U).
\]
\end{defn}

\begin{rmk}\label{another-almost} We note that this definition coincides with the tensor product $\ud{\widetilde{\m}} \otimes_{R} \F$, where $\ud{\widetilde{\m}}$ is the constant sheaf associated with the $R$-module $\m$. Alternatively, we see that $\widetilde{\m}\otimes \F \simeq \ud{\widetilde{\m}}_X\otimes_{\O_X}\F$ where $\ud{\widetilde{\m}}_X= \ud{\widetilde{\m}}\otimes_R \O_X$.
\end{rmk}

We also note that flatness of the $R$-module $\widetilde{\m}$ implies that the functor $-\otimes \widetilde{\m}$ is exact and descends to a functor
\[
-\otimes \widetilde{\m}\colon \mathbf{D}(X) \to \mathbf{D}(X),
\]
where $\mathbf{D}(X)$ is the derived category of $\O_X$-modules.

\begin{defn} An $\O_X$-module $\F$ is {\it almost zero} if $\widetilde{\m}\otimes \F$ is zero. We denote the category of almost zero $\O_X$-modules by $\Sigma_X$. 
\end{defn}

\begin{rmk} Since $\widetilde{\m}$ is an $R$-flat module, we easily see that the category of almost zero $\O_X$-modules form a Serre subcategory of $\bf{Mod}_{\O_X}=\bf{Mod}_{X}$. 
\end{rmk}

\begin{lemma}\label{almost-zero-sheaf} Let $(X, \O_X)$ be a ringed $R$-site, and let $\F$ be an $\O_X$-module. Suppose that $\mathcal U$ is a base of topology on $X$. Then the following conditions are equivalent:
\begin{enumerate}[label=\textbf{(\arabic*)}]
\item\label{almost-zero-sheaf-1} $\F\otimes\widetilde{\m}$ is the zero sheaf.
\item\label{almost-zero-sheaf-2} For any $\e \in \m$, $\e \F=0$.
\item\label{almost-zero-sheaf-3} For any $U\in \mathcal U$, the module $\widetilde{\m} \otimes \F(U)$ is zero.
\item\label{almost-zero-sheaf-4} For any $U\in \mathcal U$, the module $\m\otimes \F(U)$ is zero.
\item\label{almost-zero-sheaf-5} For any $U\in \mathcal U$, the module $\m\left(\F(U)\right)$ is zero.
\end{enumerate}
\end{lemma}
\begin{proof}
We first show that \ref{almost-zero-sheaf-1} implies \ref{almost-zero-sheaf-2}. We pick an element $\e\in \m=\m^2$ and write it as $\e=\sum x_i\cdot y_i$ for some $x_i, y_i \in \m$. So the multiplication by $\e$ map can be decomposed as 
\[
\F \xr{s\mapsto s\otimes \sum x_i \otimes y_i} \F\otimes \widetilde{\m} \xr{m} \F
\] 
where the last map is induced by the multiplication by $\widetilde{\m} \to R$. Then if $\F\otimes \widetilde{\m}=0$, then the multiplication by $\e$ map is zero for any $\e \in \m$. Now \ref{almost-zero-sheaf-2} easily implies \ref{almost-zero-sheaf-5}. Lemma~\ref{almost-zero} ensures that  \ref{almost-zero-sheaf-3}, \ref{almost-zero-sheaf-4}, and \ref{almost-zero-sheaf-5} are equivalent. Finally, \ref{almost-zero-sheaf-3} clearly implies \ref{almost-zero-sheaf-1}.
\end{proof}

\begin{lemma}\label{global-sections-almost-nothing} Let $(X, \O_X)$ be a ringed $R$-site, and let $\F$ be an almost zero $\O_X$-module. Then $\rm{H}^i(U, \F)\cong^a 0$ for any open\footnote{An open $U\in X$ is by definition an object $U\in \rm{Ob}(X)$ of the category underlying the site $X$.} $U\in X$ and any $i\geq 0$.
\end{lemma}
\begin{proof}
If $\F$ is almost zero, then $\e \F=0$ for any $\e \in \m$ by Lemma~\ref{almost-zero-sheaf}. Since the functors $\rm{H}^i(X, -)$ are $R$-linear, we conclude that $\e \rm{H}^i(U, \F)=0$ for any open $U$ and any $\e \in \m, i\geq 0$. Thus Lemma~\ref{almost-zero} ensures that $\rm{H}^i(U, \F)\cong^a 0$. 
\end{proof}

\begin{defn} We say that a homomorphism $\varphi \colon \F \to \G$ of $\O_X$-modules is an {\it almost isomorphism} if $\ker(\varphi)$ and $\coker(\varphi)$ are almost zero.
\end{defn}

\begin{lemma}\label{almost-iso-sections} A homomorphism $\varphi\colon \F\to \G$ of $\O_X$-modules is an almost isomorphism if and only if $\varphi(U)\colon \F(U) \to \G(U)$ is an almost isomorphism of $\O_X(U)$-modules for any open $U\in X$. 
\end{lemma}
\begin{proof}
The $\Leftarrow$ implication is clear from the definitions. We give a proof of the $\Rightarrow$ implication. \smallskip

Suppose that $\varphi$ is an almost isomorphism. We define the auxilliary $\O_X$-modules:
$
\mathcal K\coloneqq \ker(\varphi), \F'\coloneqq \rm{Im}(\varphi), \mathcal Q\coloneqq \coker(\varphi)
$. Lemma~\ref{global-sections-almost-nothing} implies that the maps
\[
\F(U) \to \F'(U) \text{ and } \F'(U) \to \G(U)
\]
are almost isomorphisms. In particular, the composition $\F(U) \to \G(U)$ must also be an almost isomorphism. 
\end{proof}

Now we discuss the notion of almost $\O_X$-modules on a ringed $R$-site $(X, \O_X)$. This notion can be defined in two different ways: either as the quotient of the category of $\O_X$-modules by the Serre subcategory of almost zero modules or as modules over the almost structure sheaf $\O_X^a$. Now we need to explain these two notions in more detail.

\begin{defn}\label{defn:almost-sheaves-1} We define the {\it category of almost $\O_X$-modules} as the quotient category 
\[
\bf{Mod}_{\O_X}^a\coloneqq \bf{Mod}_{\O_X}/\Sigma_X.
\] 
\end{defn}

Now we define the category $\bf{Mod}_{\O_X^a}$ of $\O_X^a$-modules that we will show to be equivalent to $\bf{Mod}_{\O_X}^a$. We recall that the almostification functor $(-)^a$ is exact and commutes with arbitrary products. This allows us to define the almost structure sheaf:

\begin{defn}\label{defn:almost-sheaves-2} The {\it almost structure sheaf} $\O_X^a$ is the sheaf\footnote{It is a sheaf exactly because $(-)^a$ is exact and commutes with arbitrary products.} of $R^a$-modules 
$
\O_X^a\colon (\rm{Ob}(X))^{op} \to \bf{Mod}_{R}^a
$ defined via the formula
$
U \mapsto \O_X(U)^a 
$. 
\end{defn}

\begin{defn} We define the {\it category of $\O_X^a$-modules} $\bf{Mod}_{\O_X^a}$ as the category of modules over $\O_X^a\in \bf{Shv}(X, R^a)$ in the categorical sense. More precisely, the objects are sheaves of $R^a$-modules $\F$ with a map $\F\otimes_{R^a} \O_X^a \to \F$ over $R^a$ satisfying the usual axioms for a module. Morphisms are defined in the evident way. 
\end{defn}

We now define the functor 
\[
(-)^a\colon \bf{Mod}_{\O_X} \to \bf{Mod}_{\O_X^a}
\]
that sends a sheaf to its ``almostification'', i.e. it applies the functor $(-)^a\colon \bf{Mod}_R \to \bf{Mod}_R^a$ section-wise. Since the almostification functor $(-)^a$ is exact and commutes with arbitrary product, it is evident that $\F^a$ is actually a sheaf for any $\O_X$-module $\F$. Moreover, it is clear that $\F^a\simeq 0$ for any almost zero $\O_X$-module $\F$. Thus, it induces the functor 
\[
(-)^a\colon \bf{Mod}_{\O_X}^a \to \bf{Mod}_{\O_X^a}.
\]
The claim is that this functor induces the equivalence of categories. The first step towards the proof is to construct the right adjoint to $(-)^a\colon \bf{Mod}_{\O_X} \to \bf{Mod}_{\O_X^a}$. Our construction of the right adjoint functor will use the existence of the left adjoint functor. So we slightly postpone the proof of the equivalence mentioned above and first discuss adjoints to $(-)^a$.
\smallskip

We start with the definition of the left adjoint functor. The idea is to apply the functor $(-)_!\colon \bf{Mod}_{\O_X^a} \to \bf{Mod}_{\O_X}$ section-wise, though this strategy does not quite work as $(-)_!$ does not commute with infinite products.

\begin{defn} 
\begin{itemize}
\item We define the functor $(-)^p_!\colon \bf{Mod}_{\O_X^a} \to \bf{Mod}^p_{\O_X}$\footnote{$\bf{Mod}^p_{\O_X}$ stands for the category of  modules over $\O_X$ in the category of presheaves} as 
\[
\F \mapsto \left(U\mapsto \F\left(U\right)_!\right)
\] 
\item We define the functor $(-)_!\colon \bf{Mod}_{\O_X^a} \to \bf{Mod}_{\O_X}$ as the composition $(-)_!\coloneqq (-)^\#\circ (-)_!^p$, where $(-)^\#$ is the sheafification functor.  
\end{itemize}
\end{defn}

\begin{lemma}\label{adjoint-almost-sheaf-!} Let $(X, \O_X)$ be a ringed $R$-site. Then
 \begin{enumerate}[label=\textbf{(\arabic*)}]
	\item\label{adjoint-almost-sheaf-!-1} The functor 
	\[
	(-)_!\colon \bf{Mod}_{\O_X^a} \to \textbf{Mod}_{\O_X}
	\]
	is the left adjoint to the localization functor $(-)^a\colon\bf{Mod}_{\O_X} \to \bf{Mod}_{\O_X^a}$. In particular, we have a functorial isomorphism 
	\[
	\rm{Hom}_{\O_X^a}(\F, \G^a)\simeq \rm{Hom}_{\O_X}(\F_!, \G)
	\] 
	for any $\F\in \bf{Mod}_{\O_X^a}, \G\in \bf{Mod}_{\O_X}$.
        \item\label{adjoint-almost-sheaf-!-2} The functor $(-)_!\colon \bf{Mod}_{\O_X^a} \to \textbf{Mod}_{\O_X}$ is exact.
	\item\label{adjoint-almost-sheaf-!-3} The counit morphism $(\F^a)_! \to \F$ is an almost isomorphism for any $\F\in \bf{Mod}_{\O_X}$. The unit morphism $\G \to (\G_!)^a$ is an isomorphism for any $\G \in \textbf{Mod}_{\O_X^a}$. In particular, the functor $(-)^a$ is essentially surjective. 
\end{enumerate}
\end{lemma}
\begin{proof}
\ref{adjoint-almost-sheaf-!-1} follows from Lemma~\ref{lemma:adjoint-almost}\ref{adjoint-almost-3} and the adjunction between sheafication and the forgetful functor. More precisely, we have the following functorial isomorphisms
\[
\rm{Hom}_{\O_X^a}(\F, \G^a)\simeq \rm{Hom}_{\bf{Mod}^p_{\O_X}}(\F_!^p, \G)\simeq \rm{Hom}_{\O_X}(\F_!, \G).
\]

We show \ref{adjoint-almost-sheaf-!-2}. It is easy to see that $(-)_!$ is left exact from Lemma~\ref{lemma:adjoint-almost}\ref{adjoint-almost-4} and the exactness of the sheafification functor. It is also right exact since it is a left adjoint functor to $(-)^a$. \smallskip

Now we show \ref{adjoint-almost-sheaf-!-3}. Lemma~\ref{lemma:adjoint-almost}\ref{adjoint-almost-5} ensures that the kernel and cokernel of the counit map of presheaves $(\F^a)_!^p \to \F$ are annihilated by any $\e\in \m$. Then the same holds after sheafification, proving the $(\F^a)_!^p \to \F$ is an almost isomorphism by Lemma~\ref{almost-zero-sheaf}. \smallskip

We consider the unit map $\G \to (\G_!)^a$, we note that using the adjuction $((-)_!, (-)^a)$ section-wise, we can refine this map
\[
\G \to (\G_!^p)^a \to (\G_!)^a.
\]
It suffices to show that both maps are isomorphisms; the first map is an isomorphism by Lemma~\ref{lemma:adjoint-almost}\ref{adjoint-almost-5}. In particular, this implies that $(\G_!^p)^a$ is already a sheaf of almost $R^a$-modules, but then we see that the natural map $(\G_!^p)^a \to (\G_!)^a$ must also be an isomorphism as it coincides with the sheafification in the category of presheaves of $R^a$-modules. 
\end{proof}

\begin{rmk} In what follows, we denote the objects of $\rm{Mod}_{\O_X^a}$ by $\F^a$ to distinguish $\O_X$ and $\O_X^a$-modules. This notation does not cause any confusion as $(-)^a$ is essentially surjective.
\end{rmk}


Now we construct the right adjoint functor to $(-)^a$. The naive idea of applying $(-)_*$ section-wise works well in this case. 

\begin{defn} The functor of {\it almost sections} $(-)_*\colon \bf{Mod}_{\O_X^a} \to \bf{Mod}_{\O_X}$ is defined as 
\[
\F^a \mapsto \Big(U \mapsto \rm{Hom}_{R}\big(\widetilde{\m}, \F^a\left(U\right)_!\big)= \rm{Hom}_{R}\big(\widetilde{\m}, \F\left(U\right)\big) \Big),
\]
where the equality comes from Lemma~\ref{prop-almost}\ref{prop-almost-2}. 
\end{defn}

\begin{rmk} The functor $(-)_*$ is well-defined, i.e. it defines a {\it sheaf} of $\O_X$-modules. This follows from the fact that $\rm{Hom}_{R}\left(\widetilde{\m}, -\right)$ is left exact and commutes with arbitrary products. 
\end{rmk}

\begin{lemma}\label{adjoint-almost-sheaf-*} Let $(X, \O_X)$ be a ringed $R$-site. Then
 \begin{enumerate}[label=\textbf{(\arabic*)}]
	\item\label{adjoint-almost-*-1} The functor $(-)_*\colon\bf{Mod}_{\O_X^a} \to \bf{Mod}_{\O_X}$ is the right adjoint to the exact localization functor $(-)^a\colon \bf{Mod}_{\O_X} \to \bf{Mod}_{\O_X^a}$. In particular, it is left exact. 
	\item\label{adjoint-almost-*-2} The unit morphism $\F \to (\F^a)_*$ is an almost isomorphism for any $\F\in \bf{Mod}_{\O_X}$. The counit morphism $(\G^a_*)^a \to \G^a$ is an isomorphism for any $\G^a\in \bf{Mod}_{\O_X^a}$.
\end{enumerate}
\end{lemma}
\begin{proof}
It is sufficient to check both claims section-wise. This, in turn, follows from Lemma~\ref{lemma:adjoint-almost}\ref{adjoint-almost-1} and Lemma~\ref{lemma:adjoint-almost}\ref{adjoint-almost-2} respectively. 
\end{proof}

\begin{cor}\label{cor:limits-colimits-almost-sheafy} The functor $(-)^a\colon \bf{Mod}_{\O_X}\to \bf{Mod}_{\O_X^a}$ commutes with limits and colimits. In particular, $\bf{Mod}_{\O_X^a}$ is complete and cocomplete, and filtered colimits and (finite) products are exact in $\bf{Mod}_{\O_X^a}$.
\end{cor}
\begin{proof}
    The first claim follows from the fact that $(-)^a$ admits left and right adjoints. The second claim follows the first claim, the exactness of $(-)^a$, and analogous exactness properties in $\bf{Mod}_R$. 
\end{proof}

\begin{cor} Let $(X, \O_X)$ be a ringed $R$-site. Then the functor
\[
(-)^a\colon \bf{Mod}_{\O_X} \to \bf{Mod}_{\O_X^a}
\]
is exact.
\end{cor}
\begin{proof}
The functor $(-)^a$ is exact as it has both left and right adjoints. 
\end{proof}


\begin{thm}\label{almost-two-different-variants} Let $(X, \O_X)$ be a ringed $R$-site. Then the functor
\[
(-)^a\colon \bf{Mod}_{\O_X}^a \to \bf{Mod}_{\O_X^a}
\]
is an equivalence of categories.
\end{thm}
\begin{proof}
Lemma~\ref{adjoint-almost-sheaf-*} implies that the functor $(-)^a\colon \bf{Mod}_{\O_X} \to \bf{Mod}_{\O_X^a}$ has a right adjoint functor $(-)_*$ such that the counit morphism $(-)^a\circ (-)_* \to \rm{Id}$ is an isomorphism of functors. Moreover, the exactness of $(-)^a$ implies that a morphism $\varphi\colon \F\to \G$ is an almost isomorphism if and only if $\varphi^a\colon \F^a\to \G^a$ is an isomorphism. Thus, \cite[Proposition 1.3]{GZ} guarantees that the induced functor $(-)^a\colon \bf{Mod}_{\O_X}^a \to \bf{Mod}_{\O_X^a}$ is an equivalence. 
\end{proof}

\begin{rmk} In what follows, we do not distinguish $\bf{Mod}_{\O_X^a}$ and $\bf{Mod}_{\O_X}^a$. Moreover, we sometimes denote both categories by $\bf{Mod}_X^a$ or $\bf{Mod}_{X^a}$ to simplify the notation.
\end{rmk}

\subsection{Basic functors on categories Of $\O_X^a$-modules}\label{basic-functors-sheaves} We discuss how to define certain basic functors on $\bf{Mod}_{X}^a$. Our main functors of interest are $\rm{Hom}$, $\rm{alHom}$, $\otimes$, $f^*$, and $f_*$. We define their almost analogues and discuss the relation with their classical versions. As a by-product, we give a slightly more intrinsic definition of $(-)_*\colon \bf{Mod}_{X}^a \to \bf{Mod}_{X}$ along the lines of the definition of the $\bf{Mod}_R^a$-version of this functor. \smallskip

For the rest of the section, we fix a ring $R$ with an ideal $\m$ such that $\m=\m^2$ and $\widetilde{\m}=\m\otimes_R \m$ is $R$-flat. We also fix an $R$-ringed site $(X, \O_X)$ that we also consider as a ringed $\O_X(X)$-site. 

\begin{defn}  \begin{itemize}
\item The {\it global Hom functor} 
\[
\rm{Hom}_{\O_X^a}(-, -) \colon \bf{Mod}_{X^a}^{op} \times \bf{Mod}_{X^a} \to \bf{Mod}_{\O_X(X)}
\]
is defined as
$
(\F^a, \G^a) \mapsto \rm{Hom}_{\O_X^a}(\F^a, \G^a)
$.
\item The {\it local Hom functor} 
\[
\ud{\mathcal{H}om}_{\O_X^a}(-, -) \colon \bf{Mod}_{X^a}^{op} \times \bf{Mod}_{X^a} \to \bf{Mod}_{X}
\]
is defined as  
$
(\F^a, \G^a) \mapsto \left(U \mapsto \rm{Hom}_{\O_U^a}(\F^a|_U, \G^a|_U)\right)
$. The standard argument shows that this functor is well-defined, i.e. $\ud{\mathcal{H}om}_{\O_X^a}(\F, \G)$ is indeed a sheaf of $\O_X$-modules. 
\end{itemize}
\end{defn}

\begin{lemma}\label{global-local-hom} Let $U$ be an open in $X$, and let $\F^a, \G^a$ be $\O_X^a$-modules. Then the natural map 
\[
\Gamma\left(U, \ud{\mathcal{H}om}_{\O_X^a}\left(\F^a, \G^a\right)\right) \to \rm{Hom}_{\O_U^a}\left(\F^a|_U, \G^a|_U\right)
\]
is an isomorphism of $\O_X(U)$-modules.
\end{lemma}
\begin{proof}
This is evident from the definition. 
\end{proof}

\begin{lemma}\label{local-adjunction} Let $(X, \O_X)$ be a ringed $R$-site. Then there is a functorial isomorphism of $\O_X$-modules
\[
\ud{\cal{H}om}_{\O_X^a}(\F^a, \G^a) \xrightarrow{\sim} \ud{\cal{H}om}_{\O_X}((\F^a)_!, \G)
\]
for $\F^a\in \bf{Mod}_X^a$ and $\G\in \bf{Mod}_X$.
\end{lemma}
\begin{proof}
Lemma~\ref{global-local-hom} and Lemma~\ref{adjoint-almost-sheaf-!} ensure that the desired isomorphism exists section-wise. It glues to a global isomorphism of sheaves since these section-wise isomorphisms are functorial in $U$.
\end{proof}

Now we move on to show a promised more intrinsic definition of the functor $(-)_*$. As a warm-up, we need the following result:

\begin{lemma}\label{local-hom-free-prel} Suppose that the ringed $R$-site $(X, \O_X)$ has a final object that (by slight abuse of notation) we denote by $X$. Then the evaluation map 
\[
\rm{ev}_X\colon \rm{Hom}_{\O_X^a}\left(\O_X^a, \G^a\right) \to \rm{Hom}_{\O_X\left(X\right)^a}\left(\O_X^a\left(X\right), \G^a\left(X\right)\right)
\]
\[
\varphi \mapsto \varphi(X)
\]
is an isomorphism of $\O_X(X)$-modules for any $\G^a\in \bf{Mod}_{X}^a$.
\end{lemma}
\begin{proof}
As $(-)^a$ is essentially surjective by Lemma~\ref{adjoint-almost-sheaf-!}\ref{adjoint-almost-sheaf-!-3}, there exists some $\O_X$-module $\G$ with almostification being equal to $\G^a$. Now we recall that the data of an $\O_X^a$-linear homomorphism $\varphi\colon \O_X^a \to \G^a$ is equivalent to the data of $\O_X(U)^a$-linear homomorphisms $\varphi_U \in \rm{Hom}_{\O_X(U)^a}\left(\O_X^a(U), \G^a(U)\right)$ for each open $U$ in $X$ such that the diagram
\[
\begin{tikzcd}
\O_X(U)^a \arrow{r}{\varphi_U} \arrow{d}{r_{\O_X^a}|^U_V}& \G(U)^a \arrow{d}{r_{\G^a}|^U_V} \\
\O_X(V)^a \arrow{r}{\varphi_V} & \G(V)^a  \\
\end{tikzcd}
\]
commutes for any $V\subset U$. Now we note that an $\O_X(U)^a$-linear homomorphism $\varphi_U$ uniquely determines an $\O_X(V)^a$-linear homomorphism $\varphi_V$ in such a diagram. Indeed, this follows from the equality
\begin{align*}
\rm{Hom}_{\O_X\left(V\right)^a}\left(\O_X\left(V\right)^a, \G\left(V\right)^a\right)&=\rm{Hom}_{\O_X\left(V\right)}\left(\widetilde{\m}\otimes \O_X\left(V\right), \G\left(V\right)\right)\\
&=\rm{Hom}_{\O_X\left(V\right)}\left(\widetilde{\m}\otimes\O_X\left(U\right)\otimes_{\O_X\left(U\right)} \O_X\left(V\right), \G\left(V\right)\right)\\
&=\rm{Hom}_{\O_X\left(U\right)}\left(\widetilde{\m}\otimes\O_X\left(U\right), \G\left(V\right)\right)\\
&=\rm{Hom}_{\O_X\left(U\right)^a}\left(\O_X\left(U\right)^a, \G\left(V\right)^a\right).
\end{align*}
Now we use the assumption that $X$ is the final object to conclude that any homomorphism $\varphi\colon \O_X^a \to \G^a$ is uniquely defined by $\varphi(X)$.
\end{proof}

\begin{cor}\label{local-hom-free} Let $(X, \O_X)$ be an $R$-ringed site, and let $U\in X$ be an open. Then the evaluation map 
\[
\rm{ev}_U\colon \rm{Hom}_{\O_U^a}\left(\O_U^a, \G|_U^a\right) \to \rm{Hom}_{\O_U\left(U\right)^a}\left(\O_U^a\left(U\right), \G^a\left(U\right)\right)
\]
\[
\varphi \mapsto \varphi(U)
\]
is an isomorphism of $\O_X(U)$-modules for any $\G^a\in \bf{Mod}_{X}^a$.
\end{cor}
\begin{proof}
For the purpose of the proof, we can change the site $X$ by the slicing site $X/U$ of objects over $U$. Then $U$ automatically becomes the final object in $X/U$, so we can apply Lemma~\ref{local-hom-free-prel} to finish the proof.
\end{proof}

Now we are ready to prove a new description of the sheaf version of the functor $(-)_*$. 

\begin{lemma} Let $(X,\O_X)$ be a ringed $R$-site. Then there is a functorial isomorphism of $\O_X$-modules 
\[
\ud{\mathcal{H}om}_{\O_X^a}(\O_X^a, \F^a) \to \F^a_*
\]
for $\F^a\in \bf{Mod}_{X}^a$.
\end{lemma}
\begin{proof}
Lemma~\ref{global-local-hom} and Corollary~\ref{local-hom-free} imply that there is an isomorphism of $\O_X(U)$-modules
\[
\Gamma\left(U, \ud{\mathcal{H}om}_{\O_X^a}\left(\O_X^a, \F^a\right)\right) \xrightarrow{\sim} \rm{Hom}_{\O_U\left(U\right)^a}\left(\O_U^a\left(U\right), \F^a\left(U\right)\right)
\]
that is functorial in both $U$ and $\F^a$. Now we use the functorial isomorphism of $\O_X(U)$
\[
\rm{Hom}_{\O_U\left(U\right)^a}\left(\O_U\left(U\right)^a, \F^a\left(U\right)\right) \simeq \rm{Hom}_{R^a}\left(R^a, \F^a\left(U\right)\right) = (\F^a)_*(U)
\]
to construct a functorial isomorphism
\[
\Gamma\left(U, \ud{\mathcal{H}om}_{\O_X^a}\left(\O_X^a, \F^a\right)\right) \xrightarrow{\sim} (\F^a)_*(U).
\]
Functoriality in $U$ ensures that it glues to the global isomorphism of $\O_X$-modules
\[
\ud{\mathcal{H}om}_{\O_X^a}\left(\O_X^a, \F^a\right) \xrightarrow{\sim} \F^a_*. \qedhere
\]
\end{proof}

Now we discuss the functor of almost homomorphisms. 

\begin{defn}\label{defn:alhom-sheafy}\begin{itemize}
\item The {\it global alHom functor} 
\[
\rm{alHom}_{\O_X^a}(-, -) \colon \bf{Mod}_{X^a}^{op} \times \bf{Mod}_{X^a} \to \bf{Mod}_{R^a}
\]
is defined as
\[
(\F^a, \G^a) \mapsto \rm{Hom}_{\O_X^a}(\F^a, \G^a)^a\simeq \rm{Hom}_{\O_X}\left(\left(\F^a\right)_!, \G\right)^a.
\]
\item The {\it local alHom functor} 
\[
\ud{al\mathcal{H}om}_{\O_X^a}(-, -) \colon \bf{Mod}_{X^a}^{op} \times \bf{Mod}_{X^a} \to \bf{Mod}_{X^a}
\]
is defined as  
\[
(\F^a, \G^a) \mapsto \left(U \mapsto \rm{alHom}_{\O_U^a}(\F^a|_U, \G^a|_U)^a\right).
\]
\end{itemize}
\end{defn}

\begin{rmk} At this point we have not checked that $\ud{al\mathcal{H}om}_{\O_X^a}(\F^a, \G^a)$ is actually a sheaf. However, this follows from the following lemma.
\end{rmk}

\begin{lemma} The natural map 
\[
\ud{\mathcal{H}om}_{\O_X}(\widetilde{\m}\otimes \F, \G)^a \to \ud{al\mathcal{H}om}_{\O_X^a}(\F^a, \G^a)
\]
is an almost isomorphism of $\O_X^a$-modules for any $\F^a, \G^a\in \bf{Mod}_{X}^a$. In particular, $\ud{al\mathcal{H}om}_{\O_X^a}(\F^a, \G^a)$ is a sheaf of $\O_X^a$-modules.
\end{lemma}
\begin{proof}
This follows from the sequence of functorial in $U$ isomorphisms: 
\begin{align*}
\ud{\mathcal{H}om}_{\O_X}(\widetilde{\m}\otimes \F, \G)(U)^a & \simeq^a \rm{Hom}_{\O_U}(\widetilde{\m} \otimes \F|_U, \G|_U)^a \\
& \simeq^a \rm{alHom}_{\O_U^a}(\F^a|_U, \G^a|_U) \\
& \simeq^a \ud{al\mathcal{H}om}_{\O_X^a}(\F^a, \G^a)(U) \qedhere
\end{align*}
\end{proof}

In order to make Definition~\ref{defn:alhom-sheafy} computable, we need to show that these functors can be computed by using any representative for $\F^a$ and $\G^a$.

\begin{prop}\label{many-functors-sheaf-alhom} \begin{enumerate}[label=\textbf{(\arabic*)}] Let $(X, \O_X)$ be a ringed $R$-site. Then:
\item\label{many-functors-sheaf-alhom-1} 
There is a natural transformation of functors 
\[
\begin{tikzcd}[column sep = 5em]
\bf{Mod}_X^{op}\times \bf{Mod}_X \arrow{r}{\rm{Hom}_{\O_X}(-, -)}\arrow{d}{(-)^a\times (-)^a} &  \bf{Mod}_X  \arrow[dl, Rightarrow, "\rho"] \arrow{d}{(-)^a}\\
\bf{Mod}_{X^a}^{op}\times \bf{Mod}_{X^a}   \arrow{r}{\rm{alHom}_{\O_X^a}(-, -)}& \bf{Mod}_{X^a}
\end{tikzcd}
\]
that makes the diagram $(2, 1)$-commutative. In particular, $\rm{alHom}_{\O_X^a}(\F^a, \G^a) \simeq \rm{Hom}_{\O_X}(\F, \G)^a$ for any $\F, \G \in \bf{Mod}_{X}$.
\item\label{many-functors-sheaf-alhom-2} 
Then there is a natural transformation of functors 
\[
\begin{tikzcd}[column sep = 5em]
\bf{Mod}_X^{op}\times \bf{Mod}_X \arrow{r}{\ud{\mathcal{H}om}_{\O_X}(-, -)}\arrow{d}{(-)^a\times (-)^a} &  \bf{Mod}_X  \arrow[dl, Rightarrow, "\rho"] \arrow{d}{(-)^a}\\
\bf{Mod}_{X^a}^{op}\times \bf{Mod}_{X^a}   \arrow{r}{\ud{al\mathcal{H}om}_{\O_X^a}(-, -)}& \bf{Mod}_{X^a}
\end{tikzcd}
\]
that makes the diagram $(2, 1)$-commutative. In particular, $\ud{al\mathcal{H}om}_{\O_X^a}(\F^a, \G^a) \simeq \ud{\mathcal{H}om}_{\O_X}(\F, \G)^a$ for any $\F, \G \in \bf{Mod}_{X}$.
\end{enumerate}
\end{prop}
\begin{proof}
The proof is similar to the proof of Proposition~\ref{many-functors}\ref{many-functors-3}. The only new thing is that we need to prove an analogue of Corollary~\ref{cor:preserve-almost-iso}, i.e. that the functors $\rm{alHom}_{\O_X}(-, \G)$, $\ud{al\cal{H}om}_{\O_X}(-, \G)$ preserve almost isomorphisms. It essentially boils down to showing that $\rm{Ext}^i_{\O_X}(\K, \G)\cong^a 0$ and $\ud{\mathcal{E}xt}^i_{\O_X}(\K, \G)\cong^a 0$ for any $\K \in \Sigma_X, \G \in \bf{Mod}_X,$ and an integer $i\geq 0$. \smallskip

Now Lemma~\ref{almost-zero-sheaf} implies that $\e \K=0$ for any $\e \in \m$. Thus, we see that $\rm{Ext}^i_{\O_X}(\K, \G)$ and $\ud{\mathcal{E}xt}^i_{\O_X}(\K, \G)$ are also annihilated by any $\e \in \m$ since the functors $\rm{Ext}^i_{\O_X}(-, \G)$, $\ud{\mathcal{E}xt}^i_{\O_X}(-, \G)$ are $R$-linear. Thus, $\rm{Ext}^i_{\O_X}(\K, \G)$ and $\ud{\mathcal{E}xt}^i_{\O_X}(\K, \G)$ are almost zero by Lemma~\ref{almost-zero} and Lemma~\ref{almost-zero-sheaf} respectively.
\end{proof}

\begin{defn} The {\it tensor product functor} $- \otimes_{\O_X^a} -  \colon \bf{Mod}_{X}^{a} \times \bf{Mod}_X^a \to \bf{Mod}_X^a $
is defined as
\[
(\F^a, \G^a) \mapsto \F^a_! \otimes_{\O_X} \G^a_!.
\]
\end{defn}

\begin{prop}\label{many-functors-tensor} There is a natural transformation of functors 
\[
\begin{tikzcd}[column sep = 5em]
\bf{Mod}_{X}\times \bf{Mod}_X \arrow{r}{-\otimes_{\O_X} -}\arrow{d}{(-)^a\times (-)^a} &  \bf{Mod}_X  \arrow{d}{(-)^a}\\
\bf{Mod}_X^a\times \bf{Mod}_X^a  \arrow[ru, Rightarrow, "\rho"] \arrow{r}{-\otimes_{\O_X^a}-}& \bf{Mod}_X^a
\end{tikzcd}
\]
that makes the diagram $(2, 1)$-commutative. In particular, there is a functorial isomorphism
\[
(\F\otimes_{\O_X} \G)^a \simeq \F^a\otimes_{\O_X^a} \G^a
\]
for any $\F, \G\in \bf{Mod}_X$.
\end{prop}
\begin{proof}
The proof is analogous to that of Propisition~\ref{many-functors}\ref{many-functors-1}. 
\end{proof} 

The tensor product is adjoint to $\ud{\mathcal{H}om}$ as it happens in the case of $R^a$-modules. We give a proof of the local version of this statement.

\begin{lemma} Let $(X, \O_X)$ be a ringed $R$-site, and let $\F^a, \G^a, \cal{H}^a$ be $\O_X^a$-modules. Then there is a functorial isomorphism
\[
\ud{\cal{H}om}_{\O_X^a}(\F^a\otimes_{\O_X^a}\G^a, \cal{H}^a)\simeq \ud{\cal{H}om}_{\O_X^a}(\F^a, \ud{al\cal{H}om}_{\O_X^a}(\G^a, \cal{H}^a)).
\]
After passing to the global sections, this gives the isomorphism
\[
\rm{Hom}_{\O_X^a}(\F^a\otimes_{\O_X^a}\G^a, \cal{H}^a)\simeq \rm{Hom}_{\O_X^a}(\F^a, \ud{al\mathcal{H}om}_{\O_X^a}(\G^a, \cal{H}^a)). 
\]
And after passing to the almostifications, it gives an isomorphism
\[
\ud{al\mathcal{H}om}_{\O_X^a}(\F^a\otimes_{\O_X^a}\G^a, \cal{H}^a)\simeq \ud{al\mathcal{H}om}_{\O_X^a}(\F^a, \ud{al\mathcal{H}om}_{\O_X^a}(\G^a, \cal{H}^a)).
\]
\end{lemma}
\begin{proof}
We compute $\Gamma(U, \ud{\cal{H}om}_{\O_X^a}(\F^a\otimes_{\O_X^a}\G^a, \cal{H}^a))$ by using Lemma~\ref{global-local-hom} and the standard $\otimes$-$\ud{\cal{H}om}$ adjunction. Namely,
\begin{align*}
\Gamma\left(U, \ud{\cal{H}om}_{\O_X^a}\left(\F^a\otimes_{\O_X^a}\G^a, \cal{H}^a\right)\right) &\simeq \rm{Hom}_{\O_U^a}\left(\F^a|_U\otimes_{\O_U^a}\G^a|_U, \cal{H}^a|_U\right) & \rm{Lemma}~\ref{global-local-hom} \\
& \simeq \rm{Hom}_{\O_U^a}\left(\left(\F|_U\otimes_{\O_U}\G|_U\right)^a, \cal{H}^a|_U\right) & \rm{Proposition}~\ref{many-functors-tensor}\\
& \simeq \rm{Hom}_{\O_U}\left(\widetilde{\m}\otimes \left(\F|_U\otimes_{\O_U}\G|_U\right), \cal{H}|_U\right) & \rm{Lemma}~\ref{adjoint-almost-sheaf-!}  \\
&\simeq \rm{Hom}_{\O_U}\left(\left(\widetilde{\m}\otimes \F|_U\right)\otimes_{\O_U}\left(\widetilde{\m}\otimes \G|_U\right), \cal{H}|_U\right) & \widetilde{\m}^{\otimes 2}\simeq \widetilde{\m}\\
&\simeq \rm{Hom}_{\O_U}\left(\widetilde{\m}\otimes \F|_U, \ud{\cal{H}om}_{\O_U}\left(\widetilde{\m}\otimes \G|_U, \cal{H}|_U\right) \right) & \otimes-\cal{H}om\text{ adjunction} \\
&\simeq \rm{Hom}_{\O_U^a}\left(\F^a|_U, \ud{al\cal{H}om}_{\O_U}\left(\widetilde{\m}\otimes \G|_U, \cal{H}|_U\right) \right) & \rm{Lemma}~\ref{adjoint-almost-sheaf-!} \\
&\simeq \Gamma\left(U, \ud{\cal{H}om}_{\O_X^a}\left(\F^a, \ud{al\cal{H}om}_{\O_X^a}\left(\G^a, \cal{H}^a\right)\right)\right) & \rm{Lemma}~\ref{global-local-hom} 
\end{align*}
Since these identifications are functorial in $U$, we can glue them to a global isomorphism 
\[
\ud{\cal{H}om}_{\O_X^a}(\F^a\otimes_{\O_X^a}\G^a, \cal{H}^a)\simeq \ud{\cal{H}om}_{\O_X^a}(\F^a, \ud{al\cal{H}om}_{\O_X^a}(\G^a, \cal{H}^a)).
\]
This finishes the proof. 
\end{proof}

\begin{cor} Let $(X, \O_X)$ be a ringed $R$-site, and let $\F^a$ be an $\O_X^a$-module. Then the functor $-\otimes_{\O_X^a}\F^a$ is left adjoint to $\ud{al\cal{H}om}_{\O_X^a}(\F^a, -)$. 
\end{cor}

For what follows, we fix a map $f\colon (X, \O_X) \to (Y, \O_Y)$ of ringed $R$-sites. We are going to define the almost version of the pullback and pushforward functors. 

\begin{defn} The {\it pullback functor} $f^*_a \colon \bf{Mod}_{X}^{a} \to \bf{Mod}_Y^a $
is defined as
\[
\F^a \mapsto \left(f^*\left(\F^a_!\right)\right)^a.
\]
In what follows, we will often abuse notation and simply write $f^*$ instead of $f^*_a$. This is ``allowed'' by Proposition~\ref{many-functors-pullback}.
\end{defn}

As always, we want to show that this functor can be actually computed by applying $f^*$ to {\it any} representative of $\F^a$. The main ingredient is to show that $f^*$ sends almost isomorphisms to almost isomorphisms. The following lemma shows slightly more, and will be quite useful later on. 

\begin{lemma}\label{pullback-proj} Let $f\colon (X, \O_X) \to (Y, \O_Y)$ be a morphism of ringed $R$-sites. Then for any $\O_X$-module $\F$, there is a natural isomorphism $\varphi_{f}(\F)\colon f^{*}(\widetilde{\m}\otimes \F) \to \widetilde{\m}\otimes f^*\F$ functorial in $\F$.
\end{lemma}
\begin{proof}
We use Remark~\ref{another-almost} to say that $\widetilde{\m}\otimes \F$ is functorially isomorphic to $\widetilde{\m}_Y\otimes_{\O_Y} \F$, where $\ud{\widetilde{\m}}_Y\coloneqq \ud{\widetilde{\m}}\otimes_R \O_Y$. Now we note that $f^*(\ud{\widetilde{\m}}_Y) \simeq \ud{\widetilde{\m}}_X$ as can be easily seen (using the $\widetilde{\m}$ is $R$-flat) from the very definitions. Therefore, $\varphi_{f}(\F)$ comes from the fact that the pullback functor commutes with the tensor product. More precisely, we define it as the composition
\[
f^{*}(\widetilde{\m}\otimes \F) \xrightarrow{\sim} f^*(\ud{\widetilde{\m}}_Y\otimes_{\O_Y} \F) \xrightarrow{\sim} f^*(\ud{\widetilde{\m}}_Y) \otimes_{\O_X}f^*(\F) \xrightarrow{\sim} \ud{\widetilde{\m}}_X \otimes_{\O_X} f^*(\F).
\]
\end{proof}

We now also show a derived version of Lemma~\ref{pullback-proj} that will be used later in the text. 

\begin{lemma}\label{derived-pull-projection} Let $f\colon (X, \O_X) \to (Y, \O_Y)$ be a morphism of ringed $R$-sites. Then for any $\F\in \bf{D}(X)$, there is a natural isomorphism 
\[
\varphi_{f}(\F)\colon \bf{L}f^{*}(\widetilde{\m}\otimes \F) \to \widetilde{\m}\otimes \bf{L}f^*\F
\]
 functorial in $\F$. 
\end{lemma}
\begin{proof}
Similarly, we use Remark~\ref{another-almost} to say that $\widetilde{\m}\otimes \F$ is functorially isomorphic to $\widetilde{\m}_Y\otimes_{\O_Y} \F$, where $\ud{\widetilde{\m}}_Y\coloneqq \ud{\widetilde{\m}}\otimes_R \O_Y$. Now we note that $\bf{L}f^*(\ud{\widetilde{\m}}_Y) \simeq f^*(\ud{\widetilde{\m}}_Y)\simeq \ud{\widetilde{\m}}_X$ as $\widetilde{\m}$ is $R$-flat. The rest of the proof is the same using the $\bf{L}f^*$ functorially commutes with the derived tensor product.
\end{proof}

\begin{cor}\label{pull-proj-cor} Let $f\colon (X, \O_X) \to (Y, \O_Y)$ be a morphism of ringed $R$-sites, and let $\varphi \colon \F \to \G$ be an almost isomorphism of $\O_Y$-modules. Then the homomorphism $f^*(\varphi)\colon f^*(\F) \to f^*(\G)$ is an almost isomorphism. 
\end{cor}
\begin{proof}
The question boils down to showing that the homomorphism 
\[
\widetilde{\m}\otimes f^*(\F) \to \widetilde{\m}\otimes f^*(\G)
\]
is an isomorphism. Lemma~\ref{pullback-proj} ensures that it is sufficient to prove that the map
\[
f^*(\widetilde{\m}\otimes \F) \to f^*(\widetilde{\m}\otimes \G)
\]
is an isomorphism. But this is clear because the map $\widetilde{\m}\otimes \F \to \widetilde{\m}\otimes \G$ is already an isomorphism. 
\end{proof}

\begin{prop}\label{many-functors-pullback} Let $f\colon(X, \O_X) \to (Y, \O_Y)$ be a morphism of ringed $R$-sites. Then there is a natural transformation of functors 
\[
\begin{tikzcd}[column sep = 5em]
\bf{Mod}_{Y} \arrow{r}{f^*}\arrow{d}{(-)^a} &  \bf{Mod}_X  \arrow{d}{(-)^a}\\
\bf{Mod}_Y^a \arrow[ru, Rightarrow, "\rho"] \arrow{r}{f^*_a}& \bf{Mod}_X^a
\end{tikzcd}
\]
that makes the diagram $(2, 1)$-commutative. In particular, there is a functorial isomorphism $(f^*\F)^a \simeq f_a^*(\F^a)$ for any $\F \in \bf{Mod}_X$.
\end{prop}
\begin{proof}
The proof is similar to that of Proposition~\ref{many-functors}. We define $\rho_\F\colon f^*(\widetilde{\m}\otimes \F)^a \to f^*(\F)^a$ as the map induced by the natural homomorphism $\widetilde{\m}\otimes \F \to \F$. It is clearly functorial in $\F$, and it is an isomorphism by Corollary~\ref{pull-proj-cor}.
\end{proof} 

\begin{defn} The {\it pushforward functor} $f_*^a \colon \bf{Mod}_{X}^{a} \to \bf{Mod}_Y^a $
is defined as
\[
\F^a \mapsto \left(f_*\left(\F^a_!\right)\right)^a.
\]
In what follows, we will often abuse the notation and simply write $f_*$ instead of $f_*^a$. This is ``allowed'' by Proposition~\ref{many-functors-pushforward}.
\end{defn}

\begin{defn} The {\it global sections functor} $\Gamma^a(X, -) \colon \bf{Mod}_{X}^{a} \to \bf{Mod}_R^a $
is defined as
\[
\F^a \mapsto \Gamma(X, \F^a_!)^a.
\]
In what follows, we will often abuse the notation and simply write $\Gamma$ instead of $\Gamma^a$. This is also ``allowed'' by Proposition~\ref{many-functors-pushforward}.
\end{defn}

\begin{rmk} The global section functor can be realized as the pushforward along the map $(X, \O_X) \to (*, R)$. 
\end{rmk}

\begin{lemma}\label{no-projection} Let $f\colon (X, \O_X)\to (Y, \O_Y)$ be a morphism of ringed $R$-sites, and let $\varphi \colon \F \to \G$ be an almost isomorphism. Then the morphism $f_*(\varphi)\colon f_*(\F) \to f_*(\G)$ is an almost isomorphism.
\end{lemma}
\begin{proof}
The standard argument considering the kernel and cokernel of $\varphi$ shows that it is sufficient to prove that $f_*\K\cong^a 0$, $\rm{R}^1f_*\K \cong^a 0$ for any almost zero $\O_X$-module $\K$. This follows from $R$-linearity of $f_*$ and Lemma~\ref{almost-zero-sheaf}.
\end{proof}

\begin{prop}\label{many-functors-pushforward} Let $f\colon(X, \O_X) \to (Y, \O_Y)$ be a morphism of ringed $R$-spaces. Then there is a natural transformation of functors 
\[
\begin{tikzcd}[column sep = 5em]
\bf{Mod}_{X} \arrow{r}{f_*}\arrow{d}{(-)^a} &  \bf{Mod}_Y  \arrow{d}{(-)^a}\\
\bf{Mod}_X^a \arrow[ru, Rightarrow, "\rho"] \arrow{r}{f_*^a}& \bf{Mod}_Y^a
\end{tikzcd}
\]
that makes the diagram $(2, 1)$-commutative. In particular, there is a functorial isomorphism $(f_*\F)^a \simeq f^a_*(\F^a)$ for any $\F \in \bf{Mod}_X$. The same results hold for $\Gamma^a(X, -)$.
\end{prop}
\begin{proof}
We define $\rho_\F\colon f_*(\widetilde{\m}\otimes \F)^a \to f_*(\F)^a$ as the map induced by the natural homomorphism $\widetilde{\m}\otimes \F \to \F$. It is clearly functorial in $\F$, and it is an isomorphism by Lemma~\ref{no-projection}.
\end{proof}

\begin{lemma}\label{lemma:alhom-global-sections} Let $(X, \O_X)$ be a ringed $R$-site, and let $\F, \G$ be $\O_X^a$-modules. Then there is a natural morphism
\[
\Gamma\left(U, \ud{al\cal{H}om}_{\O_X^a}(\F^a, \G^a)\right) \to \rm{alHom}_{\O_U^a}(\F^a|_U, \G^a|_U)
\]
is an isomorphism of $R^a$-modules for any open $U\subset X$. 
\end{lemma}
\begin{proof}
The claim easily follows from Lemma~\ref{global-local-hom}, Proposition~\ref{many-functors-sheaf-alhom}\ref{many-functors-sheaf-alhom-2}, and Proposition~\ref{many-functors-pushforward}
\end{proof}

\begin{lemma}\label{push-pull-local-adj} Let $f\colon (X, \O_X)\to (Y, \O_Y)$ be a morphism of ringed $R$-site, and let $\F^a\in \bf{Mod}_Y^a$, and  $\G^a \in \bf{Mod}_X^a$. Then there is a functorial isomorphism of $\O_Y$-modules
\[
f_*\ud{\cal{H}om}_{\O_X^a}(f^*(\F^a), \G^a)\simeq \ud{\cal{H}om}_{\O_Y^a}(\F^a, f_*(\G^a)).
\]
After passing to the global sections, this gives the isomorphism of $\O_Y(Y)$-modules
\[
\rm{Hom}_{\O_X^a}(f^*(\F^a), \G^a)\simeq \rm{Hom}_{\O_Y^a}(\F^a, f_*(\G^a)). 
\]
And after passing to the almostifications, it gives the isomorphism of $\O_Y^a$-modules
\[
f_*\ud{al\mathcal{H}om}_{\O_X^a}(f^*(\F^a), \G^a)\cong^a \ud{al\mathcal{H}om}_{\O_Y^a}(\F^a, f_*(\G^a)).
\]
\end{lemma}
\begin{proof}
This is a combination of the classical $(f^*, f_*)$-adjunction, Lemma~\ref{adjoint-almost-sheaf-!}, Lemma~\ref{pullback-proj}, Proposition~\ref{many-functors-pullback}, and Proposition~\ref{many-functors-pushforward}. Indeed, we choose an open $U\subset Y$ and denote its preimage by $V\coloneqq f^{-1}(U)$. We also define $\F^a_U\coloneqq \F^a|_U$ and $\G^a_V\coloneqq \G^a|_V$. The claim follows from the sequence of functorial isomorphisms
\begin{align*}
\Gamma\left(U, \ud{\cal{H}om}_{\O_Y^a}\left (\F^a, f_*\left(\G^a\right)\right)\right)&\simeq \rm{Hom}_{\O_U^a}\left(\F_U^a, f_*\left(\G_V^a\right)\right) & \rm{Lemma}~\ref{global-local-hom} \\
& \simeq \rm{Hom}_{\O_U^a}\left(\F_U^a, f_*\left(\G_V\right)^a\right) & \rm{Proposition}~\ref{many-functors-pushforward}\\
& \simeq \rm{Hom}_{\O_U}\left(\widetilde{\m}\otimes \F_U, f_*\left(\G_V\right)\right) & \rm{Lemma}~\ref{adjoint-almost-sheaf-!} \\
& \simeq \rm{Hom}_{\O_V}\left(f^*\left(\widetilde{\m}\otimes \F_U\right), \G_V\right) & (f^*, f_*)\text{-adjunction} \\ 
& \simeq \rm{Hom}_{\O_V}\left(\widetilde{\m}\otimes f^*\left(\F_U\right), \G_V\right) & \rm{Lemma}~\ref{pullback-proj} \\
& \simeq \rm{Hom}_{\O_V^a}\left(f^*\left(\F_U\right)^a, \G_V^a\right) & \rm{Lemma}~\ref{adjoint-almost-sheaf-!} \\
& \simeq \rm{Hom}_{\O_V^a}\left(f^*\left(\F_U^a\right), \G_V^a\right) & \rm{Proposition}~\ref{many-functors-pullback} \\
& \simeq \Gamma\left(U, f_*\ud{\cal{H}om}_{\O_X^a}\left(f^*\left(\F^a\right), \G^a\right)\right). & \rm{Lemma}~\ref{global-local-hom}
\end{align*}
Since these identifications are functorial in $U$, we can glue them to a global isomorphism 
\[
f_*\ud{\cal{H}om}_{\O_X^a}(f^*(\F^a), \G^a)\simeq \ud{\cal{H}om}_{\O_Y^a}(\F^a, f_*(\G^a)).
\]
\end{proof}

\begin{cor} 
Let $f\colon (X, \O_X)\to (Y, \O_Y)$ be a morphism of ringed $R$-site. Then the functors     $
        \begin{tikzcd}
            \bf{Mod}_X^a\arrow[r, swap, shift right=.75ex, "f_*"] & \bf{Mod}_Y^a\arrow[l, swap, shift right=.75ex, "f^*"] 
        \end{tikzcd}
    $ are adjoint.
\end{cor}

\subsection{Digression: the projection formula}\label{section:projection-formula}

In this section, we show that the tensor product $\widetilde{\m}\otimes -$ behaves especially well on locally spectral spaces\footnote{We refer to \cite[\href{https://stacks.math.columbia.edu/tag/08YF}{Tag 08YF}]{stacks-project} and \cite[\textsection 3]{wedhorn} for a comprehensive discussion of (locally) spectral spaces}. For instance, we show that we can explicitly describe sections of $\widetilde{\m}\otimes \F$ on a basis of opens for such spaces, and verify a version of the projection formula for this tensor product. 



\begin{lemma}\label{tensor-product} Let $(X,\O_X)$ be a locally spectral, locally ringed $R$-space. Then for any spectral\footnote{We remind the reader that any quasi-compact quasi-separated open subset of a locally spectral space is spectral. This can be easily seen from the definitions.} open subset $U\subset X$ the natural morphism
\[
\widetilde{\m}\otimes_R \F(U) \to (\widetilde{\m}\otimes \F)(U)
\]
is an isomorphism of $\O_X(U)$-modules.
\end{lemma}
\begin{proof}
As spectral subspaces form a basis of topology on $X$, it suffices to show that the functor
\[
U \to \widetilde{\m}\otimes_R \F(U)
\]
the sheaf condition on spectral open subsets. In particular, we can assume that $X$ itself is spectral.  

As any open spectral $U$ is quasi-compact, we conclude that any open covering $U=\bigcup_{i\in I} U_i$ admits a refinement by a finite one. Thus, it is sufficient to check the sheaf condition for finite coverings of a spectral spaces by spectral open subspaces. Thus, we need to show that, for any finite covering $U=\bigcup_{i\in I} U_i$, the sequence
\[
0 \to \widetilde{\m}\otimes_R \F(U) \to \prod_{i=1}^n (\widetilde{\m}\otimes_R \F(U_i)) \to \prod_{i,j=1}^n (\widetilde{\m}\otimes_R \F(U_i \cap U_j)).
\]
is exact. This follows from flatness of $\widetilde{\m}$ and the fact that tensor product commutes with {\it finite} direct products. 
\end{proof}

Now we want to show a version of the projection formula for the functor $\widetilde{\m}\otimes -$, it will take some time to rigorously prove it. We recall that a map of locally spectral spaces is called {\it spectral}, if the pre-image of any spectral open subset is spectral. 

\begin{lemma}\label{acyclic-1} Let $(X, \O_X)$ be a spectral locally ringed $R$-space. Then for any injective $\O_X$-module $\I$ the $\O_X$-module $\widetilde{\m}\otimes \I$ is an $\rm{H}^0(X, -)$-acyclic.
\end{lemma}
\begin{proof}
We note that spectral open subspaces form a basis for the topology on $X$. Thus \cite[\href{https://stacks.math.columbia.edu/tag/01EV}{Tag 01EV}]{stacks-project} and \cite[\href{https://stacks.math.columbia.edu/tag/0A36}{Tag 0A36}]{stacks-project} imply that it suffices to show that 
\[
(\widetilde{\m}\otimes \I)(V) \xr{r_{\widetilde{\m}\otimes \I}|_U^V} (\widetilde{\m}\otimes \I)(U)
\]
is surjective for any inclusion of any {\it spectral} open subsets $U \hookrightarrow V$. Lemma \ref{tensor-product} says that this map $r_{\widetilde{\m}\otimes \I}|_U^V$ is identified with the map
\[
\widetilde{\m}\otimes_R \I(V) \xr{\widetilde{\m}\otimes_R r_{\I}|_U^V} \widetilde{\m}\otimes_R \I(U).
\]

But now we note that $r_{\I}|_U^V$ is surjective since any injective $\O_X$-module is flasque by \cite[\href{https://stacks.math.columbia.edu/tag/01EA}{Tag 01EA}]{stacks-project}, and therefore the map $\widetilde{\m}\otimes_R r_{\I}|_U^V$ is surjective as well.
\end{proof}

\begin{Cor}\label{acyclic-2}  Let $f\colon (X, \O_X)\to (Y, \O_Y)$ be a spectral morphism of locally spectral, locally ringed $R$-spaces, and let $\I$ be an injective $\O_X$-module. Then $\widetilde{\m}\otimes \I$ is an $f_*(-)$-acyclic
\end{Cor}
\begin{proof}
It suffices to show that for any open spectral $U\subset Y$ the higher cohomology groups 
\[
\rm{H}^{i}(X_U, (\widetilde{\m}\otimes \I)|_{X_U})
\]
vanish. This follows from Lemma~\ref{acyclic-1} since $X_U$ is spectral because both $f$ and $U$ are spectral.
\end{proof}

\begin{lemma}\label{deg-0} Let $f\colon (X, \O_X)\to (Y, \O_Y)$ be a spectral morphism of locally spectral, locally ringed $R$-spaces, and let $\F$ be an $\O_X$-module. Then there is an isomorphism
\[
\beta:\widetilde{\m}\otimes f_*\F \to f_*(\widetilde{\m}\otimes \F)
\] 
functorial in $\F$.
\end{lemma}
\begin{proof}
It suffices to define a morphism on a basis of spectral open subspaces $U\subset Y$. For any such $U\subset Y$, we define 
\[
\beta_U\colon (\widetilde{\m}\otimes f_*\F)(U) \to f_*(\widetilde{\m}\otimes \F)(U)
\]
as the composition of isomorphisms
\[
(\widetilde{\m}\otimes f_*\F)(U) \xr{\a^{-1}_{U}} \widetilde{\m}\otimes_R (f_*\F)(U) = \widetilde{\m}\otimes_R \F(X_U) \xr{\a_{X_U}} (\widetilde{\m}\otimes \F)(X_U) = f_*(\widetilde{\m}\otimes \F)(U)
\]
with $\a_U$ and $\a_{X_U}$ being isomorphisms from Lemma \ref{tensor-product}. Since the construction of $\a$ is functorial in $U$, we conclude that $\beta$ defines a morphism of sheaves. It is an isomorphism because $\beta_U$ is an isomorphism an a basis of $Y$.
\end{proof} 

\begin{lemma}\label{projection} Let $f\colon (X, \O_X) \to (Y, \O_Y)$ be a spectral morphism of locally spectral, locally ringed $R$-spaces. Then for any $\F\in \mathbf{D}(X)$, there is a morphism 
\[
\rho_f(\F)\colon \widetilde{\m} \otimes \mathbf{R} f_*\F \to \mathbf{R}f_*(\widetilde{\m}\otimes \F)
\]
functorial in $\F$. This map is an isomorphism in either of the following cases:
\begin{itemize}
\item The complex $\F$ is bounded below, i.e. $\F \in \bf{D}^+(X)$, or
\item The space $X$ is locally of uniformly bounded Krull dimension and $\F\in \bf{D}(X)$.
\end{itemize}
\end{lemma}
\begin{proof} 
We start the proof by constructing the map $\rho_f(\F)$. Note that by adjunction it suffices to construct a map
\[
\mathbf{L}f^*(\widetilde{\m} \otimes \mathbf{R} f_*\F) \to \widetilde{\m}\otimes \F
\]
We define this map as the composition
\[
\mathbf{L}f^*(\widetilde{\m} \otimes \mathbf{R} f_*\F) \xr{\varphi_f(\mathbf{R}f_*\F)} \widetilde{\m}\otimes \mathbf{L}f^*\mathbf{R}f_*\F \xr{\widetilde{\m}\otimes \eta_{\F}} \widetilde{\m}\otimes \F
\]
where the first map is the isomorphism coming from Lemma~\ref{derived-pull-projection} and the second map comes from the counit $\eta_{\F}$ of the $(\bf{L}f^*, \bf{R}f_*)$-adjunction. \smallskip


Now we show that $\rho_f(\F)$ is an isomorphism for $\F\in \mathbf{D}^+(X)$. We choose an injective resolution $\F \to \mathcal I^{\bullet}$. In this case, we use Corollary~\ref{acyclic-2} to note that $\beta$ is the natural map 
\[
\widetilde{\m}\otimes f_*(\mathcal I^{\bullet})\to f_*(\widetilde{\m}\otimes \mathcal I^{\bullet})
\]
that is an isomorphism by Lemma~\ref{deg-0}. \\

The last thing we need to show is that $\rho_f(\F)$ is an isomorphism for any $\F\in \bf{D}(X)$ if $X$ is locally of uniformly bounded Krull dimension. The claim is local, so we may and do assume that both $X$ and $Y$ are spectral spaces. As $X$ is quasi-compact (because it is spectral) and locally of finite Krull dimension, we conclude that $X$ has finite Krull dimension, say $N\coloneqq \dim X$. Then \cite[Corollary 4.6]{Scheiderer} (or \cite[\href{https://stacks.math.columbia.edu/tag/0A3G}{Tag 0A3G}]{stacks-project}) implies that $\rm{H}^i(U, \G)=0$ for any open spectral $U\subset X$, $\G\in \bf{Mod}_X$, and $i>N$. In particular, $\rm{R}^if_*\G=0$ for any $\G \in \bf{Mod}_X$, and $i>N$. Thus we see that the assumptions of \cite[\href{https://stacks.math.columbia.edu/tag/0D6U}{Tag 0D6U}]{stacks-project} are verified in this case (with $\cal A=\bf{Mod}_X$ and $\cal A' = \bf{Mod}_Y$), so the natural map
\[
\cal{H}^j\left(\bf{R}f_*\F\right) \to \cal{H}^j\left(\bf{R}f_*\left(\tau^{\geq -n}\F\right)\right)
\] 
is an isomorphism for any $\F\in \bf{D}(X)$, $j\geq N-n$. As $\widetilde{\m}$ is $R$-flat, we get the commutative diagram
\[
\begin{tikzcd}[column sep = 5em]
\cal{H}^j\left(\widetilde{\m}\otimes \bf{R}f_*\F\right) \arrow{r}{\cal{H}^j(\rho_\F)} \arrow{d}{\sim} & \cal{H}^j\left(\bf{R}f_*\left(\widetilde{\m}\otimes \F\right)\right)\arrow{d}{\sim} \\
\cal{H}^j\left(\widetilde{\m}\otimes \bf{R}f_*\left( \tau^{\geq -n}\F\right)\right) \arrow{r}{\cal{H}^j(\rho_{\tau^{\geq -n} \F})} & \cal{H}^j\left(\bf{R}f_*\left(\widetilde{\m}\otimes \tau^{\geq -n}\F\right)\right)
\end{tikzcd}
\]
with the vertical arrows being isomorphisms for $j\geq N-n$, and the bottom horizontal map is an isomorphism as $\tau^{\geq -n} \F \in \bf{D}^+(X)$. Thus, by choosing an appropriate $n\geq 0$, we see that $\cal{H}^j(\rho_\F)$ is an isomorphism for any $j$; so $\rho_\F$ is an isomorphism itself.
\end{proof}

\subsection{Derived category of $\O_X^a$-modules}\label{section:derived-sheaves}

This section is a global analogue of Section~\ref{der-category-modules}. We give two different definitions of the derived category of almost $\O_X$-modules and then show that they coincide. \smallskip

For the rest of the section, we fix a ring $R$ with an ideal $\m$ such that $\m=\m^2$ and $\widetilde{\m}=\m\otimes_R \m$ is $R$-flat. We also fix an $R$-ringed site $(X, \O_X)$.

\begin{defn} We define the {\it derived category of $\O_X^a$-modules} as $\bf{D}(X^a)\coloneqq \bf{D}(\bf{Mod}_{X}^a)$. \smallskip 
\end{defn}

We define the bounded version of the derived category of almost $R$-modules $\bf{D}^*(X^a)$ for $*\in \{+, -, b\}$ as the full subcategory of $\bf{D}(X^a)$ consisting of bounded below (resp. bounded above, resp. bounded) complexes. 

\begin{defn} We define the {\it almost derived category of $\O_X$-modules} as the Verdier quotient\footnote{We refer to \cite[\href{https://stacks.math.columbia.edu/tag/05RA}{Tag 05RA}]{stacks-project} for an extensive discussion of Verdier quotients of triangulated categories.} $\bf{D}(X)^a\coloneqq \bf{D}(\bf{Mod}_{X})/\bf{D}_{\Sigma_X}(\bf{Mod}_{X})$. 
\end{defn}

\begin{rmk} We recall that $\Sigma_X$ is the Serre subcategory of $\bf{Mod}_{X}$ that consists of almost zero $\O_X$-modules. 
\end{rmk}

We note that the functor $(-)^a\colon \bf{Mod}_{X} \to \bf{Mod}_X^a$ is exact and additive. Thus, it can be derived to the functor $(-)^a\colon \bf{D}(X) \to \bf{D}(X^a)$. Similarly, the functor $(-)_!\colon \bf{Mod}_{X}^a \to \bf{Mod}_X$ can be derived to the functor $(-)_!\colon \bf{D}(X^a) \to \bf{D}(X)$. The standard argument shows that $(-)_!$ is a left adjoint functor to the functor $(-)^a$ as this already happens on the level of abelian categories. \smallskip

We also want to establish a derived version of the functor $(-)_*$. But since the functor is only left exact, we do need to do some work to derive it. Namely, we need to ensure that $\O_X^a$-modules admit enough $K$-injective complexes.  

\begin{defn} We say that a complex of $\O_X^a$-module $I^{\bullet, a}$ is {\it K-injective} if $\rm{Hom}_{K(\O_X^a)}(C^{\bullet, a}, I^{\bullet, a})=0$ for any acyclic complex $C^{\bullet, a}$ of $R^a$-modules.
\end{defn}

\begin{rmk} We remind the reader that $K(\O_X^a)$ stands for the homotopy category of $\O_X^a$-modules.
\end{rmk}

\begin{lemma}\label{K-inj-sheaf} The functor $(-)^a\colon \bf{Comp}(\O_X) \to \bf{Comp}(\O_X^a)$ sends $K$-injective $\O_X^a$-complexes to $K$-injective $\O_X^a$-complexes. 
\end{lemma}
\begin{proof}
We note that $(-)^a$ admits an exact left adjoint $(-)_!$ thus \cite[\href{https://stacks.math.columbia.edu/tag/08BJ}{Tag 08BJ}]{stacks-project} ensures that $(-)^a$ preserves $K$-injective complexes. 
\end{proof}

\begin{cor}\label{enough-K-inj-sheaf} Let $(X, \O_X)$ be a ringed $R$-site. Then every object $\F^{\bullet, a}\in \bf{Comp}(\O_X^a)$ is quasi-isomorphic to a $K$-injective complex. 
\end{cor}
\begin{proof}
The proof of Corollary~\ref{enough-K-inj} works verbatim with the only exception that one needs to use \cite[\href{https://stacks.math.columbia.edu/tag/079P}{Tag 079P}]{stacks-project} instead of \cite[\href{https://stacks.math.columbia.edu/tag/090Y}{Tag 090Y}]{stacks-project}.
\end{proof}

Now, similarly to the case of $R^a$-modules, we define the functor $(-)_*\colon \bf{D}(X^a) \to \bf{D}(X)$ as the derived functor of $(-)_*\colon \bf{Mod}_{X}^a \to \bf{Mod}_X$. This functor exists by \cite[\href{https://stacks.math.columbia.edu/tag/070K}{Tag 070K}]{stacks-project}.
 \smallskip

\begin{lemma}\label{adjoint-derived-sheafy} Let $(X, \O_X)$ be a ringed $R$-site. Then 
\begin{enumerate}[label=\textbf{(\arabic*)}]
	\item\label{adjoint-derived-sheafy-1} The functors
    $
        \begin{tikzcd}
            \bf{D}(X)\arrow[r, swap, shift right=.75ex, "(-)^a"] & \bf{D}(X^a)\arrow[l, swap, shift right=.75ex, "(-)_!"] 
        \end{tikzcd}
    $ are adjoint. Moreover, the counit (resp. unit) morphism 
    \[
   (\F^a)_! \to \F \text{ (resp. } \G \to (\G_!)^a)
    \] is an almost isomorphism (resp. isomorphism) for any $\F\in \bf{D}(X), \G\in \bf{D}(X^a)$. In particular, the functor $(-)^a$ is essentially surjective.
    \item\label{adjoint-derived-sheafy-2} The functor $(-)^a\colon \bf{D}(X)\to \bf{D}(X^a)$ also admits a right adjoint functor $(-)_*\colon \bf{D}(X^a) \to \bf{D}(X)$. Moreover, the unit (resp. counit) morphism
    \[
    \F \to (\F^a)_* \text{ (resp. }(\G_*)^a \to \G)
    \] 
    is an almost isomorphism (resp. isomorphism) for any $\F\in \bf{D}(X), \G\in \bf{D}(X^a)$.
\end{enumerate}
\end{lemma}
\begin{proof}
The proof is similar to that of Lemma~\ref{adjoint-derived}.
\end{proof}

\begin{thm}\label{derived-the-same-sheafy} The functor $(-)^a\colon \bf{D}(X) \to \bf{D}(X^a)$ induces an equivalence of triangulated categories $(-)^a\colon \bf{D}(X)^a \to \bf{D}(X^a)$.
\end{thm}
\begin{proof}
The proof is similar to that of Theorem~\ref{derived-the-same}.
\end{proof}

\begin{rmk} Theorem~\ref{derived-the-same-sheafy} shows that the two notions of the derived category of almost modules are the same. In what follows, we do not distinguish $\bf{D}(X^a)$ and $\bf{D}(X)^a$ anymore. 
\end{rmk}

\subsection{Basic functors on derived categories of $\O_X^a$-modules}\label{section:derived-sheaves-functors}

Now we can ``derive'' certain functors constructed in Section~\ref{basic-functors-sheaves}. For the rest of the section, we fix a ringed $R$-site $(X, \O_X)$. The section follows the exposition of Section~\ref{basic-functors-derived-modules} very closely.

\begin{defn} We define the {\it derived Hom} functors 
\[
\bf{R}\ud{\cal{H}om}_{\O_X^a}(-, -) \colon \bf{D}(X^a)^{ op}\times \bf{D}(X^a) \to \bf{D}(X), \text{ and }
\]
\[
\bf{R}\rm{Hom}_{\O_X^a}(-, -) \colon \bf{D}(X^a)^{op}\times \bf{D}(X^a) \to \bf{D}(R)
\]
as it is done in \cite[\href{https://stacks.math.columbia.edu/tag/08DH}{Tag 08DH}]{stacks-project} and \cite[\href{https://stacks.math.columbia.edu/tag/0B6A}{Tag 0B6A}]{stacks-project}, respectively.
\end{defn}

\begin{defn} We define the {\it global Ext-modules} as the $R$-modules 
\[
\rm{Ext}^i_{\O_X^a}(\F^a, \G^a)\coloneqq \rm{H}^i(\bf{R}\rm{Hom}_{\O_X^a}(\F^a, \G^a))
\]
for $\F^a, \G^a\in\bf{Mod}_X^a$. \smallskip

We define the {\it local Ext-sheaves} as the $\O_X$-modules $\ud{\cal{E}xt}^i_{\O_X^a}(\F^a, \G^a)\coloneqq \cal{H}^i(\bf{R}\ud{\cal{H}om}_{\O_X^a}(\F^a, \G^a))$ for $\F^a, \G^a\in\bf{Mod}_X^a$. 
\end{defn}

\begin{rmk} We see that \cite[\href{https://stacks.math.columbia.edu/tag/0A64}{Tag 0A64}]{stacks-project} implies that there is a functorial isomorphism 
\[
\rm{H}^i\left(\bf{R}\rm{Hom}_{\O_X^a}\left(\F^a, \G^a\right)\right) \simeq \rm{Hom}_{\bf{D}(R)^a}\left(\F^a, \G^a[i]\right) 
\]
for $\F^a, \G^a\in \bf{D}(X)^a$. 
\end{rmk}

\begin{rmk}\label{local-global-derived} The standard argument shows that there is a functorial isomorphism 
\[
\bf{R}\Gamma(U, \bf{R}\ud{\cal{H}om}_{\O_X^a}(\F^a, \G^a)) \simeq \bf{R}\rm{Hom}_{\O_U^a}(\F^a|_U, \G^a|_U)
\]
for any open $U\in  X$, $\F^a, \G^a\in \bf{D}(X)^a$. 
\end{rmk}

Now we show a local version of the $((-)_!, (-)^a)$-adjunction, and relate $\bf{R}\ud{\cal{H}om}$ (resp. $\bf{R}\rm{Hom}$) to a certain derived functor. This goes in complete analogy with the situation in the usual (not almost) world. 

\begin{lemma}\label{derived-hom-alg-sheaf} Let $(X, \O_X)$ be a ringed $R$-site. Then
\begin{enumerate}[label=\textbf{(\arabic*)}]
\item\label{derived-hom-alg-sheaf-1} There is a functorial isomorphism 
\[
\bf{R}\ud{\cal{H}om}_{\O_X^a}(\F^a, \G^a) \simeq \bf{R}\ud{\cal{H}om}_{\O_X}(\F^a_!, \G)
\]
for any $\F^a\in \bf{D}(X)^a$ and $\G\in \bf{D}(X)$. In particular, this isomorphism induces functorial isomorphisms
\[
\bf{R}\rm{Hom}_{\O_X^a}(\F^a, \G^a) \simeq \bf{R}\rm{Hom}_{\O_X}(\F^a_!, \G) \text{ and } \rm{Hom}_{\bf{D}(X)^a}(\F^a, \G^a) \simeq \bf{R}\rm{Hom}_{\bf{D}(X)}(\F^a_!, \G).
\]

\item\label{derived-hom-alg-sheaf-2} For any chosen $\F^a\in \bf{Mod}_X^a$, the functor $\bf{R}\rm{Hom}_{\O_X^a}(\F^a, -)\colon \bf{D}(X)^a \to \bf{D}(R)$ is isomorphic to the (right) derived functor of $\rm{Hom}_{\O_X^a}(\F^a, -)$.

\item\label{derived-hom-alg-sheaf-3} For any chosen $\F^a\in \bf{Mod}_X^a$, the functor $\bf{R}\ud{\cal{H}om}_{\O_X^a}(\F^a, -)\colon \bf{D}(X)^a \to \bf{D}(X)$ is isomorphic to the (right) derived functor of $\rm{Hom}_{\O_X^a}(\F^a, -)$.

\end{enumerate}
\end{lemma}
\begin{proof}
\ref{derived-hom-alg-sheaf-1}: Lemma~\ref{K-inj-sheaf} and the construction of derived homs ensure that 
\begin{align*}
\bf{R}\ud{\cal{H}om}_{\O_X^a}(\F^a, \G^a) \simeq \ud{\cal{H}om}_{\O_X^a}^\bullet(\F^{\bullet, a}, \I^{\bullet, a}) \\
\bf{R}\ud{\cal{H}om}_{\O_X}(\F^a_!, \G) \simeq \ud{\cal{H}om}_{\O_X}^\bullet(\F^{\bullet, a}_!, \I^{\bullet}),
\end{align*}
where $\cal{G} \to \cal{I}^\bullet$ is a $K$-injective resolution. Now we recall the term-wise equalities
\[
\ud{\cal{H}om}_{\O_X^a}^n(\F^{\bullet, a}, \I^{\bullet, a})=\prod_{p+q=n} \ud{\cal{H}om}_{\O_X^a}(\F^{-q, a}, \I^{p, a})
\]
\[
\ud{\cal{H}om}_{\O_X}^n(\F^{\bullet, a}_!, \I^{\bullet})=\prod_{p+q=n} \ud{\cal{H}om}_{\O_X}(\F^{-q, a}_!, \I^{p})
\]
Thus Lemma~\ref{local-adjunction} produces term-wise isomorphisms
\[
\kappa_n\colon \ud{\cal{H}om}_{\O_X^a}^n(\F^{\bullet, a}, \I^{\bullet, a}) \to \ud{\cal{H}om}_{\O_X}^n(\F^{\bullet, a}_!, \I^{\bullet})
\]
that commute with the differentials by inspection, therefore defining the desired isomorphism of complexes.

Parts~\ref{derived-hom-alg-sheaf-2} and~\ref{derived-hom-alg-sheaf-3} are identical to Lemma~\ref{derived-hom-alg}\ref{derived-hom-alg-2}.
\end{proof}

\begin{defn}\label{defn-almost-hom-sheaf-derived} We define the {\it derived almost Hom} functors 
\[
\bf{R}\ud{al\cal{H}om}_{\O_X^a}(-, -) \colon \bf{D}(X^a)^{op}\times \bf{D}(X^a) \to \bf{D}(X^a)
\] 
\[
\bf{R}\rm{alHom}_{\O_X^a}(-, -) \colon \bf{D}(X^a)^{op}\times \bf{D}(X^a) \to \bf{D}(R^a)
\] 
as 
\[
\bf{R}\ud{al\cal{H}om}_{\O_X^a}(\F^a, \G^a)\coloneqq \bf{R}\ud{\cal{H}om}_{\O_X^a}(\F^a, \G^a)^a=\bf{R}\ud{\cal{H}om}_{\O_X}(\F^a_!, \G)^a 
\]
\[
\bf{R}\rm{alHom}_{\O_X^a}(\F^a, \G^a)\coloneqq \bf{R}\rm{Hom}_{\O_X^a}(\F^a, \G^a)^a=\bf{R}\rm{Hom}_{\O_X}(\F^a_!, \G)^a 
\]
\end{defn}

\begin{defn}
We define the {\it global almost Ext modules} as the $R^a$-modules  $\rm{alExt}^i_{\O_X^a}(\F^a, \G^a)\coloneqq \rm{H}^i(\bf{R}\rm{alHom}_{\O_X^a}(\F^a, \G^a))$ for $\F^a, \G^a\in\bf{Mod}_X^a$. \smallskip

We define the {\it local almost Ext sheaves} as the $\O_X^a$-modules  $\ud{al\cal{E}xt}^i_{\O_X^a}(\F^a, \G^a)\coloneqq \cal{H}^i(\bf{R}\ud{al\cal{H}om}_{\O_X^a}(\F^a, \G^a))$ for $\F^a, \G^a\in\bf{Mod}_X^a$. 
\end{defn}

\begin{prop}\label{derived-al-hom-sheaf} Let $(X, \O_X)$ be a ringed $R$-site. Then:
\begin{enumerate}[label=\textbf{(\arabic*)}]
\item\label{derived-al-hom-sheaf-1} There is a natural transformation of functors   
\[
\begin{tikzcd}[column sep = 7em, row sep = 5em]
\bf{D}(X)^{op}\times \bf{D}(X) \arrow{r}{\bf{R}\ud{\cal{H}om}_{\O_X}(-, -)} \arrow{d}{(-)^a\times (-)^a}& \bf{D}(X)  \arrow{d}{(-)^a}\arrow[dl, Rightarrow, "\rho"]\\
\bf{D}(X^a)^{op}\times \bf{D}(X^a)  \arrow{r}{\bf{R}\ud{al\cal{H}om}_{\O_X^a}(-, -)}& \bf{D}(X^a)
\end{tikzcd}
\]
that makes the diagram $(2, 1)$-commutative. In particular, $\bf{R}\ud{al\cal{H}om}_{\O_X^a}(\F^a, \G^a)\simeq \bf{R}\ud{\cal{H}om}_{\O_X}(\F, \G)^a$ for any $\F, \G\in \bf{D}(X)$. 
\item\label{derived-al-hom-sheaf-2} For any chosen $\F^a\in \bf{Mod}_R^a$, the functor $\bf{R}\ud{al\cal{H}om}_{\O_X^a}(\F^a, -)\colon \bf{D}(X)^a \to \bf{D}(X)^a$ is isomorphic to the (right) derived functor of $\ud{al\cal{H}om}_{\O_X^a}(\F^a, -)$.
\item\label{derived-al-hom-sheaf-3} The analogous results hold true for the functor $\bf{R}\rm{alHom}_{\O_X^a}(-, -)$.
\end{enumerate}
\end{prop}
\begin{proof}
The proof is identical to that of Proposition~\ref{derived-al-hom}. One only needs to use Proposition~\ref{many-functors-sheaf-alhom} in place of Proposition~\ref{many-functors}\ref{many-functors-3}. 
\end{proof}

Now we deal with the case of the derived tensor product functor. We will show that our definition of the derived tensor product functor makes $\bf{R}\ud{al\cal{H}om}_{\O_X^a}(-, -)$ into the inner Hom functor. 

\begin{defn} We say that a complex of $\O_X^a$-module $\F^{\bullet, a}$ is {\it almost K-flat} if the naive tensor product complex $\mathcal C^{\bullet, a}\otimes^{\bullet}_{\O_X^a} \F^{\bullet, a}$ is acyclic for any acyclic complex $\mathcal C^{\bullet, a}$ of $\O_X^a$-modules.
\end{defn}


\begin{lemma}\label{almost-K-flat-sheaf} The functor $(-)^a\colon \bf{Comp}(\O_X) \to \bf{Comp}(\O_X^a)$ sends $K$-flat $\O_X$-complexes to almost $K$-flat $\O_X^a$-complexes. 
\end{lemma}
\begin{proof}
The proof Lemma~\ref{almost-K-flat} applies verbatim. 
\end{proof}

\begin{lemma}\label{K-flat-preserved} Let $f\colon (X,\O_X) \to (Y, \O_Y)$ be a morphism of ringed $R$-sites, and let $\F^{\bullet, a}\in \bf{Comp}(\O_Y^a)$ be an almost $K$-flat complex. Then $f^*(\F^{\bullet, a})\in \bf{Comp}(\O_X^a)$ is almost $K$-flat.
\end{lemma}
\begin{proof}
The proof of \cite[\href{https://stacks.math.columbia.edu/tag/06YW}{Tag 06YW}]{stacks-project} works verbatim in this situation.
\end{proof}

\begin{cor}\label{enough-K-sheaf} Every object $\F^{\bullet, a}\in \bf{Comp}(\O_X^a)$ is quasi-isomorphic to an almost $K$-flat complex. 
\end{cor}
\begin{proof}
The proof of Corollary~\ref{enough-K} applies verbatim with the only difference that one needs to use \cite[\href{https://stacks.math.columbia.edu/tag/06YF}{Tag 06YF}]{stacks-project} in place of \cite[\href{https://stacks.math.columbia.edu/tag/06Y4}{Tag 06Y4}]{stacks-project}. 
\end{proof}

\begin{defn} We define the {\it derived tensor product functor} 
\[
-\otimes^L_{\O_X^a}-\colon \bf{D}(X)^a\times \bf{D}(X)^a \to \bf{D}(X)^a
\] by the rule 
$
(\F^a, \G^a) \mapsto (\G_!\otimes^L_{\O_X} \G_!)^a
$
for any $\F^a, \G^a \in \bf{D}(X)^a$.
\end{defn}

\begin{prop}\label{derived-tensor-product-sheaf}
\begin{enumerate}[label=\textbf{(\arabic*)}]
\item\label{derived-tensor-product-sheaf-1} There is a natural transformation of functors 
\[
\begin{tikzcd}
\bf{D}(X)\times \bf{D}(X) \arrow{r}{-\otimes_{\O_X}^L -} \arrow[d, swap, "(-)^a\times (-)^a"]& \bf{D}(X) \arrow{d}{(-)^a} \\
\bf{D}(X)^a\times \bf{D}(X)^a \arrow{r}{-\otimes_{\O_X^a}^L -} \arrow[ru, Rightarrow, "\rho"]& \bf{D}(X)^a \\
\end{tikzcd}
\]
that makes the diagram $(2, 1)$-commutative. In particular, there is a functorial isomorphism $(\F\otimes^L_{\O_X} \G)^a\simeq \F^a\otimes^L_{\O_X^a} \G^a$ for any $\F, \G\in \bf{D}(X)$. 
\item\label{derived-tensor-product-sheaf-2} For any chosen $\F^a\in \bf{Mod}_X^a$, the functor $\F^a\otimes^L_{R^a}-\colon \bf{D}(X)^a \to \bf{D}(X)^a$ is isomorphic to the (left) derived functor of $\F^a\otimes_{\O_X^a}-$.
\end{enumerate}
\end{prop}
\begin{proof}
Again, the proof is identical to that of Proposition~\ref{derived-tensor-product-sheaf}. The only non-trivial input that we need is the existence of sufficiently many $K$-flat complexes of $\O_X^a$-modules. But this is guaranteed by Lemma~\ref{enough-K-sheaf}.
\end{proof}

\begin{rmk}\label{rmk:first-map-in-duality} For any $\F^a, \G^a \in \bf{D}(X)^a$, there is a canonical morphism
\[
\bf{R}\ud{al\cal{H}om}_{\O_X^a}(\F^a, \G^a) \otimes^L_{\O_X} \F^a \to \G^a
\]
that, after the identifications from Proposition~\ref{derived-al-hom-sheaf} and Proposition~\ref{derived-tensor-product-sheaf}, is the almostification of the canonical morphism
\[
\bf{R}\ud{\cal{H}om}_{\O_X}(\F^a_!, \G^a_!) \otimes^L_{\O_X} \F^a_! \to \G^a_!
\]
from \cite[\href{https://stacks.math.columbia.edu/tag/0A8V}{Tag 0A8V}]{stacks-project}. 
\end{rmk}

\begin{lemma}\label{o-hom-adj-derived-sheaf} Let $(X,\O_X)$ be a ringed $R$-site, and let $\F^a,\G^a,\cal{H}^a \in \bf{D}(X)^a$. Then we have a functorial isomorphism
\[
\bf{R}\ud{\cal{H}om}_{\O_X^a}(\F^a\otimes^L_{\O_X^a} \G^a, \cal{H}^a) \simeq \bf{R}\ud{\cal{H}om}_{\O_X^a}(\F^a, \bf{R}\ud{al\cal{H}om}_{R^a}(\G^a, \cal{H}^a)).
\]
This induces functorial isomorphisms
\[
\bf{R}\rm{Hom}_{\O_X^a}(\F^a\otimes^L_{\O_X^a} \G^a, \cal{H}^a) \simeq \bf{R}\rm{Hom}_{\O_X^a}(\F^a, \bf{R}\ud{al\cal{H}om}_{R^a}(\G^a, \cal{H}^a)) \ ,
\]
\[
\bf{R}\ud{al\cal{H}om}_{\O_X^a}(\F^a\otimes^L_{\O_X^a} \G^a, \cal{H}^a) \simeq \bf{R}\ud{al\cal{H}om}_{\O_X^a}(\F^a, \bf{R}\ud{al\cal{H}om}_{R^a}(\G^a, \cal{H}^a)) \ ,
\]
\[
\bf{R}\rm{alHom}_{\O_X^a}(\F^a\otimes^L_{\O_X^a} \G^a, \cal{H}^a) \simeq \bf{R}\rm{alHom}_{\O_X^a}(\F^a, \bf{R}\ud{al\cal{H}om}_{R^a}(\G^a, \cal{H}^a)).
\]
\end{lemma}
\begin{proof}
The proof of the first isomorphism is very similar to that of Lemma~\ref{o-hom-adj}. We leave the details to the interested reader. The second isomorphism comes from the fist one by applying the functor $\mathbf{R}\Gamma(X, -)$. The third and the fourth isomorphisms are obtained by applying $(-)^a$ to the first and the second isomorphisms respectively. Here, we implicitly use Proposition~\ref{derived-al-hom-sheaf}. 
\end{proof}

\begin{cor} Let $(X, \O_X)$ be a ringed $R$-site, and let $\G^a\in \bf{D}(X)^a$. Then the functors \[
        \begin{tikzcd}
          \bf{R}\ud{al\cal{H}om}_{\O_X^a}(\G^a, -)\colon   \bf{D}(X)^a\arrow[r, swap, shift right=.75ex] & \bf{D}(X)^a\colon -\otimes^L_{\O_X^a} \G^a \arrow[l, swap, shift right=.75ex] 
        \end{tikzcd}
    \] are adjoint. 
\end{cor}

Now we discuss the almost analogues of derived pullbacks and derived pushforwards. We start with the derived pullbacks:

\begin{defn} Let $f\colon (X, \O_X) \to (Y, \O_Y)$ be a morphism of ringed $R$-sites. We define the {\it derived pullback functor} 
\[
\bf{L}f^* \colon \bf{D}(Y)^a \to \bf{D}(X)^a
\] as the derived functor of the right exact, additive functor $f^*\colon \bf{Mod}_Y^a \to \bf{Mod}_X^a$.
\end{defn}

\begin{rmk} We need to explain why the desired derived functor exists and how it can be computed. It turns out that it can be constructed by choosing $K$-flat resolutions, the argument for this is identical to \cite[\href{https://stacks.math.columbia.edu/tag/06YY}{Tag 06YY}]{stacks-project}. We only emphasize that three main inputs are Lemma~\ref{K-flat-preserved}, Lemma~\ref{almost-K-flat-sheaf}, and an almost analogue of \cite[\href{https://stacks.math.columbia.edu/tag/06YG}{Tag 06YG}]{stacks-project}. 
\end{rmk}

\begin{prop}\label{derived-pullback} Let $f\colon (X, \O_X) \to (Y,\O_Y)$ be a morphism of ringed $R$-sites. Then 
there is a natural transformation of functors 
\[
\begin{tikzcd}
\bf{D}(Y) \arrow{r}{\bf{L}f^*} \arrow[d, swap, "(-)^a"]& \bf{D}(X) \arrow{d}{(-)^a} \\
\bf{D}(Y)^a \arrow{r}{\bf{L}f^*} \arrow[ru, Rightarrow, "\rho"]& \bf{D}(X)^a \\
\end{tikzcd}
\]
that makes the diagram $(2, 1)$-commutative. In particular, there is a functorial isomorphism $(\bf{L}f^*\F)^a\simeq \bf{L}f^*(\F^a)$ for any $\F\in \bf{D}(Y)$. 
\end{prop}
\begin{proof} We construct the natural transformation $\rho\colon \bf{L}f^*\circ (-)^a \Rightarrow (-)^a\circ \bf{L}f^*$ as follows. Pick any object $\F\in \bf{D}(Y)$ and its $K$-flat representative $\K^{\bullet}$, then $\K^{\bullet}$ is adapted to compute the usual derived pullback $\bf{L}f^*$. Lemma~\ref{K-flat-preserved} ensures $\K^{\bullet, a}$ is also adapted to compute the almost version of the derived pullback $\bf{L}f^*$. So we define the morphism
\[
\rho_{\F} \colon (f^*(\widetilde{\m}\otimes \K^{\bullet}))^a \to f^*(\K^{\bullet})^a
\]
as the natural morphism induced by $\widetilde{\m}\otimes \K^{\bullet} \to \K^{\bullet}$. This map is clearly functorial, so it defines a transformation of functors $\rho$. To show that it is an isomorphism of functors, it suffices to show that the map
\[
f^*(\widetilde{\m}\otimes \K^{\bullet}) \to f^*(\K^{\bullet})
\] 
is an almost isomorphism of complexes for any $K$-flat complex $K^{\bullet}$. But this is clear as $\widetilde{\m}\otimes \K^{\bullet} \to \K^{\bullet}$ is an almost isomorphism, and Corollary~\ref{pull-proj-cor} ensures that $f^*$ preserves almost isomorphisms. 
\end{proof}

\begin{defn} Let $f\colon (X, \O_X) \to (Y, \O_Y)$ be a morphism of ringed $R$-sites. We define the {\it derived pushforward functor} 
\[
\bf{R}f_* \colon \bf{D}(X)^a \to \bf{D}(Y)^a
\] as the derived functor of the left exact, additive functor $f_*\colon \bf{Mod}_X^a \to \bf{Mod}_Y^a$. \smallskip

We define the {\it derived global sections functor} $\bf{R}\Gamma(U, -)\colon \bf{D}(X)^a \to \bf{D}(R)^a$ in a similar way for any open $U\subset X$. 
\end{defn}

\begin{rmk} This functor exists by abstract nonsense (i.e. \cite[\href{https://stacks.math.columbia.edu/tag/070K}{Tag 070K}]{stacks-project}) as the category $\bf{Mod}_X^a$ has enough $K$-injective complexes by Lemma~\ref{enough-K-inj-sheaf}. 
\end{rmk}

\begin{prop}\label{derived-pushforward} Let $f\colon (X, \O_X) \to (Y,\O_Y)$ be a morphism of ringed $R$-sites. Then 
there is a natural transformation of functors 
\[
\begin{tikzcd}
\bf{D}(X) \arrow{r}{\bf{R}f_*} \arrow[d, swap, "(-)^a"]& \bf{D}(Y) \arrow{d}{(-)^a} \\
\bf{D}(X)^a \arrow{r}{\bf{R}f_*} \arrow[ru, Rightarrow, "\rho"]& \bf{D}(Y)^a \\
\end{tikzcd}
\]
that makes the diagram $(2, 1)$-commutative. In particular, there is a functorial isomorphism $(\bf{R}f_*\F)^a\simeq \bf{R}f_*(\F^a)$ for any $\F\in \bf{D}(X)$. The analogous results hold for the functor $\bf{R}\Gamma(U, -)$. 
\end{prop}
\begin{proof}
The proof is very similar to that of Proposition~\ref{derived-pullback}. The main essential ingredients are: $(-)^a$ sends $K$-injective complexes to $K$-injective complexes, and $f_*$ preserves almost isomorphisms. These two results were shown in Lemma~\ref{K-inj-sheaf} and Lemma~\ref{no-projection}.
\end{proof}

\begin{lemma} Let $(X, \O_X)$ be a ringed $R$-site, let $\F$ be an $\O_X^a$-module, and let $U\in X$ be an open object. Then we have a canonical isomorphism 
\[
\bf{R}\Gamma(U, \bf{R}\ud{al\cal{H}om}_{\O_X^a}(\F^a, \G^a)) \simeq \bf{R}\rm{alHom}_{\O_U^a}(\F^a|_U, \G^a|_U)
\]
\end{lemma}
\begin{proof}
This easily follows from Remark~\ref{local-global-derived}, Proposition~\ref{derived-al-hom-sheaf}, and Proposition~\ref{derived-pushforward}. 
\end{proof}

\begin{lemma} Let $f\colon (X, \O_X) \to (Y, \O_Y)$ be a morphism of ringed $R$-sites. Then there is a functorial isomorphism
\[
\bf{R}f_*\bf{R}\ud{\cal{H}om}_{\O_X^a}(\bf{L}f^*\F^a, \G^a) \simeq \bf{R}\ud{\cal{H}om}_{\O_Y^a}(\F^a, \bf{R}f_*\G^a) 
\]
for $\F^a \in \bf{D}(Y)^a$, $\G^a \in \bf{D}(X)^a$. This isomorphism induces isomorphisms
\[
\bf{R}f_*\bf{R}\ud{al\cal{H}om}_{\O_X^a}(\bf{L}f^*\F^a, \G^a) \simeq \bf{R}\ud{al\cal{H}om}_{\O_Y^a}(\F^a, \bf{R}f_*\G^a) \ ,
\]
\[
\bf{R}\rm{Hom}_{\O_X^a}(\bf{L}f^*\F^a, \G^a) \simeq \bf{R}\rm{Hom}_{\O_Y^a}(\F^a, \bf{R}f_*\G^a)  \ ,
\]
\[
\bf{R}\rm{alHom}_{\O_X^a}(\bf{L}f^*\F^a, \G^a) \simeq \bf{R}\rm{alHom}_{\O_Y^a}(\F^a, \bf{R}f_*\G^a).
\]
\end{lemma}
\begin{proof}
It is a standard exercise to show that the first isomorphism implies all other isomorphisms by applying certain functors to it, so we deal only with the first one. The proof of the first one is also quite standard and similar to Lemma~\ref{push-pull-local-adj}, but we spell it out for the reader's convenience. The desired isomorphism comes from a sequence of canonical identifications:
\begin{align*}
\bf{R}f_*\bf{R}\ud{\cal{H}om}_{\O_X^a}(\bf{L}f^*(\F^a), \G^a) & \simeq \bf{R}f_*\bf{R}\ud{\cal{H}om}_{\O_X^a}(\bf{L}f^*(\F)^a, \G^a) &      \rm{Proposition}~\ref{derived-pullback}   \\
&\simeq \bf{R}f_*\bf{R}\ud{\cal{H}om}_{\O_X}(\widetilde{\m}\otimes \bf{L}f^*(\F), \G) & \rm{Lemma}~\ref{derived-hom-alg-sheaf}\ref{derived-hom-alg-sheaf-1}  \\
&\simeq \bf{R}f_*\bf{R}\ud{\cal{H}om}_{\O_X}(\bf{L}f^*(\widetilde{\m}\otimes \F), \G) & \rm{Lemma}~\ref{derived-pull-projection} \\
&\simeq \bf{R}\ud{\cal{H}om}_{\O_Y}(\widetilde{\m}\otimes \F, \bf{R}f_*(\G)) & \rm{Classical} \\
&\simeq \bf{R}\ud{\cal{H}om}_{\O_Y^a}(\F^a, \bf{R}f_*(\G)^a) & \rm{Lemma}~\ref{derived-hom-alg-sheaf}\ref{derived-hom-alg-sheaf-1} \\
& \simeq \bf{R}\ud{\cal{H}om}_{\O_Y^a}(\F^a, \bf{R}f_*(\G^a)) & \rm{Proposition}~\ref{derived-pushforward}. 
\end{align*}
\end{proof}

\begin{cor}\label{cor:*-adjoint-derived} Let $f\colon (X, \O_X) \to (Y, \O_Y)$ be a morphism of ringed $R$-sites. Then the functors $
        \begin{tikzcd}
          \bf{R}f_*(-)\colon   \bf{D}(X)^a\arrow[r, shift right=.75ex] & \bf{D}(Y)^a\colon \bf{L}f^*(-) \arrow[l,  shift right=.75ex] 
        \end{tikzcd}
    $ are adjoint. 
\end{cor}

Now we discuss the projection formula in the world of almost sheaves. Suppose $f\colon (X, \O_X) \to (Y, \O_Y)$ a morphism of ringed $R$-sites, $\F^a\in \bf{D}(X)^a$, and $\G^a\in \bf{D}(Y)^a$. We wish to construct the projection morphism
\[
\rho \colon \bf{R}f_*(\F^a)\otimes^L_{\O_Y^a} \G^a \to \bf{R}f_*(\F^a\otimes^L_{\O_X^a} \bf{L}f^*(\G^a)).
\]
By Corollary~\ref{cor:*-adjoint-derived}, it is equivalent to constructing a morphism
\[
\pi\colon \bf{L}f^*(\bf{R}f_*(\F^a)\otimes^L_{\O_Y^a} \G^a) \to \F^a\otimes^L_{\O_X^a} \bf{L}f^*(\G^a).
\]
We define $\pi$ as the composition of the natural isomorphism
\[
\bf{L}f^*(\bf{R}f_*(\F^a)\otimes^L_{\O_Y^a} \G^a) \simeq \bf{L}f^*(\bf{R}f_*(\F^a)) \otimes^L_{\O_X^a} \bf{L}f^*(\G^a)
\]
and the morphism
\[
\bf{L}f^*(\bf{R}f_*(\F^a)) \otimes^L_{\O_X^a} \bf{L}f^*(\G^a) \xr{\e_{\F^a}\otimes \rm{Id}} \F^a\otimes^L_{\O_X^a} \bf{L}f^*(\G^a)
\]
induced by the co-unit of the $(\bf{L}f^*, \bf{R}f_*)$-adjunction. 

\begin{prop}\label{prop:projection-formula-perfect} Let $f\colon (X, \O_X) \to (Y, \O_Y)$ be a be a morphism of ringed $R$-sites, $\F^a\in \bf{D}(X)^a$, and $\G\in \bf{D}(Y)$ a perfect complex. Then the projection morphism
\[
\rho\colon \bf{R}f_*(\F^a)\otimes^L_{\O_Y^a} \G^a \to \bf{R}f_*(\F^a\otimes^L_{\O_X^a} \bf{L}f^*(\G^a))
\]
is an isomorphism in $\bf{D}(Y)^a$.
\end{prop}
\begin{proof}
The claim is local on $Y$, so we may assume that $\G$ is isomorphic to a bounded complex of finite free $\O_Y$-modules. Then an easy argument with naive filtrations reduces the question to the case when $\G=\O_Y^n$. This case is essentially obvious. 
\end{proof}

\section{Almost coherent sheaves on schemes and formal schemes}

In this Section, we develop a theory of almost coherent sheaves on schemes, and on a ``nice'' class of formal schemes.

\subsection{Schemes. The category of almost coherent $\O_X^a$-modules}\label{acoh-sheaves}

In this Section we discuss the notion of almost quasi-coherent, almost finite type, almost finitely presented and almost coherent sheaves on an arbitrary scheme. One of the main goals is to show that almost coherent sheaves form a weak Serre subcategory in $\O_X$-modules. Another important result is the ``approximation'' Corollary~\ref{cor:approximate-afpr}; it will play a key role in reducing many global results about almost finitely presented $\O_X$-modules to the classical case of finitely presented $\O_X$-modules. In particular, we follow this approach in our proof of the Almost Proper Mapping Theorem in Section~\ref{APMT-formal}. \smallskip

As always, we fix a ring $R$ with an ideal $\m$ such that $\m=\m^2$ and $\widetilde{\m}=\m\otimes_R \m$ is $R$-flat. We always do almost mathematics with respect to this ideal. In what follows, $X$ will always denote an $R$-scheme. Note that this implies that $X$ is a locally spectral, ringed $R$-site, so the results of Section~\ref{almost-sheaves} and Section~\ref{section:projection-formula} apply. \smallskip

We begin with some definitions:

\begin{defn} We say that an $\O_X^a$-module $\F^a$ is {\it almost quasi-coherent} if $\F^a_!\simeq \widetilde{\m}\otimes \F$ is a quasi-coherent $\O_X$-module. \smallskip

We say that an $\O_X$-module $\F$ is {\it almost quasi-coherent} if $\F^a$ is an almost quasi-coherent $\O_X^a$-module.
\end{defn}

\begin{rmk} Any quasi-coherent $\O_X$-module is almost quasi-coherent. 
\end{rmk}

\begin{rmk}\label{remark-quasi} We denote by $\bf{Mod}_{X^a}^{\aqc} \subset \bf{Mod}_{X^a}$ the full subcategory consisting of almost quasi-coherent $\O_X^a$-modules. It is straightforward\footnote{The proof is completely similar to the proof of Theorem~\ref{almost-two-different-variants} or Theorem~\ref{derived-the-same-sheafy}.} to see that the ``almostification'' functor induces an equivalence 
\[
\bf{Mod}_{X^a}^{\aqc} \simeq \bf{Mod}_{X}^{\qc}/(\Sigma_X\cap \bf{Mod}_X^{\qc}),
\] 
i.e. $\bf{Mod}_{X^a}^{\aqc}$ is equivalent to the quotient category of quasi-coherent $\O_X$-modules by the full subcategory of almost zero, quasi-coherent $\O_X$-modules.
\end{rmk}

\begin{defn}\label{defn:almost-finitely-resented-schemes} We say that an $\O_X^a$-module $\F^a$ is {\it of almost finite type} (resp. {\it almost finitely presented}) if $\F^a$ is almost quasi-coherent, and there is a covering of $X$ by open affines $\{U_i\}_{i\in I}$ such that $\F^a(U_i)$ is an almost finitely generated (resp. almost finitely presented) $\O_X^a(U_i)$-module. \smallskip

We say that an $\O_X$-module $\F$ is {\it of almost finite type} (resp. {\it almost finitely presented}) if so is $\F^a$. 
\end{defn}

\begin{rmk}\label{remark-quasi-finite} We denote by $\bf{Mod}^{\qc, \text{aft}}_X$ (resp. $\bf{Mod}^{\qc, \text{afp}}_X$) the full subcategory of $\bf{Mod}_X$ consisting of almost finite type (resp. almost finitely presented) quasi-coherent $\O_X$-modules. Similarly, we denote by $\bf{Mod}^{\text{aft}}_{X^a}$ (resp. $\bf{Mod}^{\text{afp}}_{X^a}$) the full subcategory of $\bf{Mod}_{X^a}$ consisting of almost finite type (resp. almost finitely presented) $\O_X^a$-modules. Again, it is straightforward to see that the ``almostification'' functors induce equivalences 
\[
    \bf{Mod}^{\aft}_{X^a} \simeq \bf{Mod}^{\qc, \aft}_X/(\Sigma_X\cap \bf{Mod}^{\qc, \aft}_X)  , \ \bf{Mod}^{\afp}_{X^a} \simeq \bf{Mod}^{\qc, \afp}_X/(\Sigma_X\cap \bf{Mod}^{\qc, \afp}_X).
\]  
\end{rmk}

\begin{rmk} Recall that, in the usual theory of $\O_X$-modules, finite type $\O_X$-modules are usually not required to be quasi-coherent. However, it is more convenient for our purposes to put almost quasi-coherence in the definition of almost finite type modules.
\end{rmk}

The first thing that we need to check is that these notions do not depend on a choice of an affine covering. 

\begin{lemma}\label{localize-ft} Let $\F^a$ be an almost finite type (resp. almost finitely presented) $\O_X^a$-module on an $R$-scheme $X$. Then $\F^a(U)$ is an almost finitely generated (resp. almost finitely presented) $\O_X^a(U)$-module for any open affine $U\subset X$.
\end{lemma}
\begin{proof}
First of all, Corollary~\ref{almost-almost-finitely-presented} and Lemma~\ref{tensor-product} imply that for any open affine $U$, $\F^a(U)$ is almost finitely generated (resp. almost finitely presented) if and only if so is $(\widetilde{\m}\otimes \F^a)(U)$. Thus, we can replace $\F^a$ by $\F^a_!\simeq \widetilde{\m}\otimes \F$ to assume that $\F$ is an honest quasi-coherent $\O_X$-module. \smallskip

Now Lemma~\ref{trivial-base-change} guarantees that the problem is local on $X$. So we can assume that $X=U$ is an affine scheme and we need to show that $\F(X)$ is almost finitely generated (resp. almost finitely presented).\smallskip

We pick some covering $X=\cup_{i=1}^n U_i$ by open affines $U_i$ such that $\F(U_i)$ is almost finitely generated (resp. almost finitely presented) as an $\O_X(U_i)$-module. We note that since $\F$ is quasi-coherent we have an isomorphism 
\[
\F(U_i) \simeq \F(X) \otimes_{\O_X(X)} \O_X(U_{i}).
\]
Now we see that a map $\O_X(X) \to \prod_{i=1}^n\O_X(U_i)$ is faithfully flat, and the module 
\[
\F(X)\otimes_{\O_X(X)} \left( \prod_{i=1}^n \O_X\left(U_i\right) \right) \simeq \left( \prod_{i=1}^n \O_X\left(U_i\right)\right) \otimes_{\O_X(X)}\F(X)
\] 
is almost finitely generated (resp. almost finitely presented) over $\prod_{i=1}^n \O_X(U_i)$. Therefore, Lemma \ref{almost-flat-descent} guarantees that $\F(X)$ is almost finitely generated (resp. almost finitely presented) as an $\O_X(X)$-module.
\end{proof}

\begin{cor}\label{criterion-almost-finite-type} Let $X=\Spec A$ be an affine $R$-scheme, and let $\F^a$ be an almost quasi-coherent $\O_X^a$-module. Then $\F^a$ is almost finite type (resp. almost finitely presented) if and only if $\Gamma(X, \F^a)$ is almost finitely generated (resp. almost finitely presented) $A$-module.
\end{cor}

Now we check that almost finite type and almost finitely presented $\O_X^a$ behave nicely in short exact sequences.

\begin{lemma}\label{main-almost-sheaf} Let $0 \to \F'^a \xr{\varphi} \F^a \xr{\psi} \F''^a \to 0$ be an exact sequence of $\O_X^a$-modules. Then
\begin{enumerate}[label=\textbf{(\arabic*)}]
\itemsep0.5ex
\item If $\F^a$ is almost of finite type and $\F''^a$ is almost quasi-coherent, then $\F''^a$ is almost finite type.
\item If $\F'^a$ and $\F''^a$ are of almost finite type (resp. finitely presented), then so is $\F^a$.
\item If $\F^a$ is of almost finite type and $\F''^a$ is almost finitely presented, then $\F'^a$ is of almost finite type.
\item If $\F^a$ is almost finitely presented and $\F'^a$ is of almost finite type, then $\F''^a$ is almost finitely presented.
\end{enumerate}
\end{lemma}
\begin{proof}
First of all, we apply the exact functor $(-)_!$ to all $\O_X^a$-modules in question to assume the short sequence is an exact sequence of $\O_X$-modules and all $\O_X$-modules in this sequence are quasi-coherent. Note that we implicitly use here that quasi-coherent modules form a Serre subcategory of all $\O_X$-modules by \cite[\href{https://stacks.math.columbia.edu/tag/01IE}{Tag 01IE}]{stacks-project}. Then we check the statement on the level of global sections on all open affine subschemes $U\subset X$ using that quasi-coherent sheaves have vanishing higher cohomology on affine schemes. And that is done in Lemma~\ref{main-almost}.
\end{proof}

\begin{defn} We say that an $\O_X^a$-module $\F^a$ is {\it almost coherent} if $\F^a$ is almost finite type, and for any open set $U$ any almost finite type $\O_U^a$-submodule $\G^a \subset (\F^a|_U)$ is an almost finitely presented $\O_U^a$-module. 

We say that an $\O_X$-module $\F$ is {\it almost coherent} if $\F^a$ is an almost coherent $\O_X^a$-module. 
\end{defn}

\begin{lemma}\label{criterion-almost-coherent} Let $X$ be an $R$-scheme, and let $\F^a$ be an $\O_X^a$-module. Then the following are equivalent:
\begin{enumerate}[label=\textbf{(\arabic*)}]
\item\label{criterion-almost-coherent-1} $\F^a$ is almost coherent.
\item\label{criterion-almost-coherent-2} $\F^a$ is almost quasi-coherent, and the $\O_X^a(U)$-module $\F^a(U)$ is almost coherent for any open affine subscheme $U\subset X$. 
\item\label{criterion-almost-coherent-3} $\F^a$ is almost quasi-coherent, and there is a covering of $X$ by open affine subschemes $(U_i)_{i\in I}$ such that $\F^a(U_i)$ is almost coherent for each $i$.
\end{enumerate} 
In particular, if $X=\Spec A$ is an affine $R$-scheme and $\F^a$ is an almost quasi-coherent $\O_X^a$-module. Then $\F^a$ is almost coherent if and only if $\F^a(X)$ is almost coherent as an $A$-module.
\end{lemma}
\begin{proof}
We start the proof by noting that we can replace $\F^a$ by $\F^a_!$ to assume that $\F$ is a quasi-coherent $\O_X$-module. \smallskip

First, we check that \ref{criterion-almost-coherent-1} implies \ref{criterion-almost-coherent-2}. Given any affine open $U\subset X$ and any almost finitely generated almost submodule $M^a\subset \F(U)^a$, we define an almost subsheaf $\widetilde{(M^a)_!} \subset (\F|_U)^a$. We see that $\widetilde{(M^a)_!}$ must be an almost finitely presented $\O_U$-module, so Lemma \ref{localize-ft} guarantees that $M^a_!$ is almost finitely presented $\O_X(U)$-module. Therefore, any almost finitely generated submodule $M^a \subset \F(U)^a$ is almost finitely presented. This shows that $\F(U)$ is almost coherent. \smallskip

Now we show that \ref{criterion-almost-coherent-2} implies \ref{criterion-almost-coherent-1}. Suppose that $\F$ is almost quasi-coherent and $\F(U)$ is almost coherent for any open affine $U\subset X$. First of all, it shows that $\F$ is of almost finite type, since this notion is local by definition. Now suppose that we have an almost finite type almost $\O_X$-submodule $\G \subset (\F|_U)^a$ for some open $U$. It is represented by a homomorphism
\[
\widetilde{\m} \otimes \G \xr{g} \F
\]
with $\G$ being an $\O_X$-module of almost finite type, and $\widetilde{\m} \otimes \ker g\simeq 0$. We want to show that $\G$ is almost finitely presented as $\O_X$-module. This is a local question, so we can assume that $U$ is affine. Then Lemma \ref{tensor-product} implies that the natural morphism 
\[
g(U): \widetilde{\m}\otimes_R \G(U) \to \F(U)
\]
defines an almost submodule of $\F(U)$. We conclude that $\widetilde{\m}\otimes_R \G(U)$ is almost finitely presented by the assumption on $\F(U)$. Since the notion of almost finitely presented $\O_X$-module is local, we see that $\G$ is almost finitely presented. \smallskip

Clearly, \ref{criterion-almost-coherent-2} implies \ref{criterion-almost-coherent-3}. And it is easy to see that Lemma~\ref{almost-flat-descent-coh} guarantees that \ref{criterion-almost-coherent-3} implies \ref{criterion-almost-coherent-2}.
\end{proof}

\begin{Cor}\label{main-coh-glob}  Let $X$ be an $R$-scheme, then:
\begin{enumerate}[label=\textbf{(\arabic*)}]
\item\label{main-coh-glob-1} Any almost finite type $\O_X^a$-submodule of an almost coherent $\O_X^a$-module is almost coherent.
\item\label{main-coh-glob-2} Let $\varphi\colon  \F^a \to \G^a$ be a homomorphism from an almost finite type $\O_X^a$-module to an almost coherent $\O_X^a$-module, then $\ker (\varphi)$ is an almost finite type $\O_X^a$-module.
\item\label{main-coh-glob-4} Let $\varphi\colon \F^a \to \G^a$ be a homomorphism of almost coherent $\O_X^a$-modules, then $\ker(\varphi)$ and $\coker(\varphi)$ are almost coherent $\O_X^a$-modules.
\item\label{main-coh-glob-5} Given a short exact sequence of $\O_X^a$-modules
\[0 \to \F'^a \to \F^a \to \F''^a \to 0\]
 if two out of three are almost coherent, so is the third.
\end{enumerate}
\end{Cor}
\begin{proof}
The proofs \ref{main-coh-glob-1}, \ref{main-coh-glob-2} and \ref{main-coh-glob-4} are quite straightforward. As usually, we apply $(-)_!$ to assume that all sheaves in the question are quasi-coherent $\O_X$-modules. Then the question is local and it is sufficient to work on global sections over all affine open subschemes $U\subset X$. So the problem is reduced to Lemma \ref{main-coh}.

The proof of part \ref{main-coh-glob-5} is similar, but we need to invoke that given a short exact sequence of $\O_X$-modules
\[
0 \to \F'^a_! \to \F^a_! \to \F''^a_! \to 0
\]
if two of these sheaves are quasi-coherent, so is the third one. This is proven in the affine case in  \cite[\href{https://stacks.math.columbia.edu/tag/01IE}{Tag 01IE}]{stacks-project}, the general case reduces to the affine one. The rest of the argument is the same.
\end{proof}

\begin{defn}\label{remark-quasi-coherent} We define the categories $\bf{Mod}_X^{\acoh}$ (resp. $\bf{Mod}_X^{\qc, \acoh}$, resp. $\bf{Mod}_{X^a}^{\acoh}$) as the full subcategory of $\bf{Mod}_X$ (resp. $\bf{Mod}_X$, resp. $\bf{Mod}_{X^a}$) consisting of almost coherent $\O_X$-modules (resp. quasi-coherent almost coherent modules, resp. almost coherent almost $\O_X$-modules). As always, it is straightforward to see that the ``almostification'' functor induces the equivalence
\[
\bf{Mod}_{X^a}^{\acoh} \simeq \bf{Mod}_X^{\qc, \acoh}/(\Sigma_X \cap \bf{Mod}_X^{\qc, \acoh}).
\]
Moreover, Corollary~\ref{main-coh-glob} ensures that $\bf{Mod}_X^{\acoh} \subset \bf{Mod}_X$, $\bf{Mod}_X^{\qc, \acoh} \subset \bf{Mod}_X$, and $\bf{Mod}_{X^a}^{\acoh} \subset \bf{Mod}_{X^a}$ are weak Serre subcategories. 
\end{defn}

The last thing that we discuss here is the notion of almost coherent schemes.

\begin{defn} We say that an $R$-scheme $X$ is {\it almost coherent} if the sheaf $\O_X$ is an almost coherent $\O_X$-module.
\end{defn}

\begin{lemma}\label{sch-coh-acoh} Let $X$ be a coherent $R$-scheme. Then $X$ is also almost coherent.
\end{lemma}
\begin{proof}
The structure sheaf $\O_X$ is quasi-coherent by definition. Lemma~\ref{criterion-almost-coherent} says that it suffices to show that $\O_X(U)$ is an almost coherent $\O_X(U)$-module for any open affine $U\subset X$. Since $X$ is coherent, we conclude that $\O_X(U)$ is coherent as an $\O_X(U)$-module. Then Lemma~\ref{coh-acoh} implies that it is actually almost coherent.
\end{proof}

\begin{lemma}\label{sheaf-acoh-afp} Let $X$ be an almost coherent $R$-scheme. Then an $\O_X^a$-module $\F^a$ is almost coherent if and only if it is of almost finite presentation.
\end{lemma}
\begin{proof}
The ``only if'' part is easy since any almost coherent $\O_X^a$-module is of almost finite presentation by definition. The ``if'' part is a local question, so we can assume that $X$ is affine, then the claim follows from Lemma~\ref{acoh-fp-coh}.
\end{proof}

\subsection{Schemes. Basic functors on almost coherent $\O_X^a$-modules}\label{basic-functors-schemes} 
This section is devoted to study how certain functors defined in Section~\ref{basic-functors-sheaves} interact with the notions of almost (quasi-) coherent $\O_X^a$-modules defined in the previous section. \smallskip

As always, we fix a ring $R$ with an ideal $\m$ such that $\m=\m^2$ and $\widetilde{\m}=\m\otimes_R \m$ is $R$-flat. We always do almost mathematics with respect to this ideal.\smallskip

We start with the affine situation, i.e., $X=\Spec A$. In this case, we note that the functor $\widetilde{(-)}\colon \bf{Mod}_A \to \bf{Mod}^{\qc}_X$ sends almost zero $A$-modules to almost zero $\O_X$-modules. Thus, it induces the functor
\[
\widetilde{(-)}\colon \bf{Mod}_{A^a} \to \bf{Mod}^{\aqc}_{X^a}.
\]
\begin{lemma}\label{sheafify-modules} Let $X=\Spec A$ be an affine $R$-scheme. Then the functor $\widetilde{(-)}\colon \bf{Mod}_A \to \bf{Mod}^{\qc}_X$ induces equivalences $\widetilde{(-)}\colon \bf{Mod}_A^* \to \bf{Mod}^{\qc, *}_{X}$ for any $*\in \{``\text{ ''}, \text{aft}, \text{\text{afp}}, \text{acoh}\}$. The quasi-inverse functor is given by $\Gamma(X, -)$. 
\end{lemma}
\begin{proof}
We note that the functor $\widetilde{(-)}\colon \bf{Mod}_A \to \bf{Mod}^{\qc}_X$ is an equivalence with the quasi-inverse $\Gamma(X, -)$. Now we note that Lemma~\ref{criterion-almost-finite-type} and Lemma~\ref{criterion-almost-coherent} guarantee that a quasi-coherent $\O_X$-module $\F$ is almost finite type (resp. almost finitely presented, resp. almost coherent) if $\Gamma(X, \F)$ is almost finitely generated (resp. almost finitely presented, resp. almost coherent) as an $A$-module.  
\end{proof}

\begin{lemma}\label{sheafify-modules-almost} Let $X=\Spec A$ be an affine $R$-scheme. Then the functor $\widetilde{(-)}\colon \bf{Mod}_{A^a} \to \bf{Mod}^{\aqc}_{X^a}$ induces an equivalence $\widetilde{(-)}^a\colon \bf{Mod}_{A^a} \to \bf{Mod}_{X^a}^{\aqc}$ and restricts to equivalences $\widetilde{(-)}^a\colon \bf{Mod}_{A^a}^* \to \bf{Mod}_{X^a}^*$ for any $*\in \{\text{aft}, \text{afp}, \text{acoh}\}$. The quasi-inverse functor is given by $\Gamma(X, -)$. 
\end{lemma}
\begin{proof}
We note that $\widetilde{(-)}\colon \bf{Mod}_A \to \bf{Mod}^{\qc}_X$ induces an equivalence between almost zero $A$-modules and almost zero, quasi-coherent $\O_X$-modules. Thus, the claim follows from Lemma~\ref{sheafify-modules}, Remark~\ref{remark-quasi}, Remark~\ref{remark-quasi-finite}, Definition~\ref{remark-quasi-coherent} and the analogous presentations of $\bf{Mod}_{A^a}^{*}$ as quotients of $\bf{Mod}_{A^a}$ for any $*\in \{\text{aft}, \text{afp}, \text{acoh}\}$. 
\end{proof}

Now we show that the pullback functor preserves almost finite type and almost finitely presented $\O_X^a$-modules.

\begin{lemma}\label{pullback-almost-coh} Let $f\colon X \to Y$ be a morphism of $R$-schemes.  
\begin{enumerate}[label=\textbf{(\arabic*)}]
\item\label{pullback-almost-coh-1} Suppose that $X=\Spec B$, $Y=\Spec A$ are affine $R$-schemes. Then $f^*(\widetilde{M^a})$ is functorially isomorphic to $\widetilde{M^a\otimes_{A^a} B^a}$ for any $M^a\in \bf{Mod}_A^a$.
\item\label{pullback-almost-coh-2} The functor $f^*$ preserves almost quasi-coherence (resp. almost finite type, resp. almost finitely presented) for $\O$-modules.
\item\label{pullback-almost-coh-3} The functor $f^*$ preserves almost quasi-coherence (resp. almost finite type, resp. almost finitely presented) for $\O^a$-modules.
\end{enumerate}
\end{lemma}
\begin{proof}
\ref{pullback-almost-coh-1} follows from Proposition~\ref{many-functors-pullback} and the analogous result for quasi-coherent $\O_Y$-modules. More precisely, Proposition~\ref{many-functors-pullback} provides with the functorial isomorphism
\[
f^*\left(\widetilde{M^a}\right) \simeq \left(f^*(\widetilde{M})\right)^a\simeq \left(\widetilde{M\otimes_A B}\right)^a\simeq \widetilde{M^a\otimes_{A^a} B^a}.
\]
Now we note that \ref{pullback-almost-coh-2} and \ref{pullback-almost-coh-3} are local on $X$ and $Y$, so we may and do assume that $X=\Spec B$, $Y=\Spec A$ are affine $R$-schemes. Furthermore, it clearly suffices to prove \ref{pullback-almost-coh-2} as \ref{pullback-almost-coh-3} follows formally from it. \smallskip

Now Lemma~\ref{sheafify-modules-almost} guarantees that any almost quasi-coherent $\O_X^a$-module is isomorphic to $\widetilde{M^a}$ for some $A^a$-module $M^a$. Now \ref{pullback-almost-coh-1} ensures that $f^*(\widetilde{M^a})\simeq \widetilde{M^a\otimes_{A^a} B^a}$ is an almost quasi-coherent $\O_X^a$-module. The other claims from \ref{pullback-almost-coh-2} are proven similarly using Lemma~\ref{sheafify-modules-almost} and Lemma~\ref{trivial-base-change}.
\end{proof}

Now we discuss that tensor product preserves certain finiteness properties of almost sheaves. 

\begin{lemma}\label{tensor-almost-coh} Let $X$ be an $R$-scheme.  
\begin{enumerate}[label=\textbf{(\arabic*)}]
\item\label{tensor-almost-coh-1} Suppose that $X=\Spec A$ is an affine $R$-scheme. Then $\widetilde{M^a}\otimes_{\O_X^a}\widetilde{N^a}$ is functorially isomorphic to $\widetilde{M^a\otimes_{A^a} N^a}$ for any $M^a, N^a\in \bf{Mod}_A^a$. 
\item\label{tensor-almost-coh-2} Let $\F^a, \G^a$ be two almost finite type (resp. almost finitely presented) $\O_X^a$-modules. Then the $\O_X^a$-module $\F^a\otimes_{\O_X^a}\G^a$ is almost finite type (resp. almost finitely presented). The analogous result holds for $\O_X$-modules $\F, \G$.  
\item\label{tensor-almost-coh-3} Let $\F^a$ be an almost coherent $\O_X^a$-module, and let $\G^a$ be an almost finitely presented $\O_X^a$-module. Then $\F^a\otimes_{\O_X^a}\G^a$ is an almost coherent $\O_X^a$-module. The analogous result holds for $\O_X$-modules $\F, \G$. 
\end{enumerate}
\end{lemma}
\begin{proof}
The proof is similar to the proof of Lemma~\ref{pullback-almost-coh}. The only difference is that one needs to use Proposition~\ref{many-functors-tensor} in place of Proposition~\ref{many-functors-pullback} to prove Part~\ref{tensor-almost-coh-1}. Part~\ref{tensor-almost-coh-2} follows from Lemma~\ref{tensor-of-fp}, and Part~\ref{tensor-almost-coh-3} follows from Corollary~\ref{tensor-hom-coh}.
\end{proof}

We show $f_*$ preserves almost quasi-coherence of $\O^a$-modules for a quasi-compact and quasi-separated morphism $f$. Later on, we will be able to show that $f_*$ also preserves almost coherence of $\O^a$-modules for certain proper morphisms. 

\begin{lemma}\label{pushforward-almost-coh} Let $f\colon X\to Y$ be a quasi-compact and quasi-separated morphism of $R$-schemes. Then  
\begin{enumerate}[label=\textbf{(\arabic*)}]
\item\label{pushforward-almost-coh-1} The $\O_Y$-module $f_*(\F)$ is almost quasi-coherent for any almost quasi-coherent $\O_X$-module $\F$.
 \item\label{pushforward-almost-coh-2} The $\O_Y^a$-module $f_*(\F^a)$ is almost quasi-coherent for any almost quasi-coherent $\O_X^a$-module $\F^a$.
\end{enumerate}
\end{lemma}
\begin{proof}
Since $\F$ is almost quasi-coherent, we conclude that $\widetilde{\m}\otimes\F$ is a quasi-coherent $\O_X$-module. Thus $f_*(\widetilde{\m}\otimes\F)$ is a quasi-coherent $\O_Y$-module by \cite[\href{https://stacks.math.columbia.edu/tag/01LC}{Tag 01LC}]{stacks-project}. Recall that the projection formula (Lemma~\ref{deg-0}) ensures that 
\[
f_*(\widetilde{\m}\otimes\F) \simeq \widetilde{\m}\otimes f_*\F.
\]
Thus, we see that $\widetilde{\m}\otimes f_*\F \simeq f_*(\F^a)_!$ is a quasi-coherent $\O_Y$-module. This shows that both $f_*(\F)$ and $f_*(\F^a)$ are almost quasi-coherent over $\O_Y$ and $\O_Y^a$ respectively. This finishes the proof.  
\end{proof}

Finally, we deal with the $\ud{\cal{H}om}_{\O_X^a}(-, -)$ functor. We start with the following preparatory lemma:

\begin{lemma}\label{hom-alcoh} Let $X$ be an $R$-scheme.  
\begin{enumerate}[label=\textbf{(\arabic*)}]
\item\label{hom-alcoh-1} Suppose $X=\Spec A$ is an affine $R$-scheme. Then the canonical map 
\begin{equation}\label{1-1-1}
\widetilde{\rm{Hom}_A(M,N)} \to \ud{\cal{H}om}_{\O_X}(\widetilde{M}, \widetilde{N})
\end{equation}
is an almost isomorphism of $\O_X$-modules for any almost finitely presented $A$-module $M$ and any $A$-module $N$. 
\item\label{hom-alcoh-2} Suppose $X=\Spec A$ is an affine $R$-scheme. Then there is a functorial isomorphism 
\begin{equation}\label{1-1-3}
\widetilde{\rm{alHom}_{A^a}(M^a,N^a)} \simeq \ud{al\cal{H}om}_{\O_X^a}(\widetilde{M^a}, \widetilde{N^a})
\end{equation}
of $\O_X^a$-modules for any almost finitely presented $A^a$-module $M^a$, and any $A^a$-module $N^a$. We also get a functorial almost isomorphism
\begin{equation}\label{1-1-2-e}
\widetilde{\rm{Hom}_A(M,N)} \simeq^a \ud{\cal{H}om}_{\O_X^a}(\widetilde{M^a}, \widetilde{N^a})
\end{equation}
of $\O_X$-modules for any almost finitely presented $A$-module $M$, and any $A$-module $N$. 
\item\label{hom-alcoh-3} Suppose $\F$ is an almost finitely presented $\O_X$-module and $\G$ is an almost quasi-coherent $\O_X$-module, then $\ud{\cal{H}om}_{\O_X}(\F, \G)$ is an almost quasi-coherent $\O_X$-module.
\item\label{hom-alcoh-4} Suppose $\F^a$ is an almost finitely presented $\O_X^a$-module and $\G^a$ is an almost quasi-coherent $\O_X^a$-module, then $\ud{\cal{H}om}_{\O_X^a}(\F^a, \G^a)$ (resp. $\ud{al\cal{H}om}_{\O_X^a}(\F^a, \G^a)$) is an almost quasi-coherent $\O_X$-module (resp. $\O_X^a$-module).
\end{enumerate}
\end{lemma}
\begin{proof}
\ref{hom-alcoh-1}: We note that we have a canonical isomorphism $\rm{Hom}_A(M,N) \to \rm{Hom}_{\O_X}(\widetilde{M}, \widetilde{N})$ for any $A$-modules $M$, $N$. This induces a morphism 
\[
\widetilde{\rm{Hom}_A(M,N)} \to \ud{\cal{H}om}_{\O_X}(\widetilde{M}, \widetilde{N}).
\]
In order to show that it is an almost isomorphism for an almost finitely presented $M$, it suffices to show that the natural map
\[
\rm{Hom}_A(M,N) \otimes_A A_f \to \rm{Hom}_{A_f}(M\otimes_A A_f, N\otimes_A A_f)
\]
is an almost isomorphism for any $f\in A$. This follows from Lemma~\ref{flat-base-change-hom}. \medskip

\ref{hom-alcoh-2} follows easily from \ref{hom-alcoh-1}. Indeed, we apply the functorial isomorphism 
\[
\ud{\cal{H}om}_{\O_X}(\F,\G)^a \simeq \ud{al\cal{H}om}_{\O_X^a}(\F^a,\G^a)
\] 
from Proposition~\ref{many-functors-sheaf-alhom}\ref{many-functors-sheaf-alhom-2} to the almost isomorphism in Part~\ref{hom-alcoh-1} to get a functorial isomorphism
\[
\widetilde{\rm{Hom}_A(M,N)^a} \simeq \ud{al\cal{H}om}_{\O_X^a}(\widetilde{M^a}, \widetilde{N^a}).
\]
Now we use Proposition~\ref{many-functors}\ref{many-functors-3} to get a functorial isomorphism 
\[
\rm{alHom}_{A^a}(M^a,N^a) \simeq \rm{Hom}_A(M,N)^a.
\]
Applying the functor $\widetilde{(-)}$ to this isomorphism and composing it with the isomorphism above, we get a functorial isomorphism
\[
\widetilde{\rm{alHom}_{A^a}(M^a,N^a)} \simeq \ud{al\cal{H}om}_{\O_X^a}(\widetilde{M^a}, \widetilde{N^a}).
\]
The construction of the isomorphism~(\ref{1-1-2-e}) is similar and even easier. \medskip

\ref{hom-alcoh-3} is a local question, so we can assume that $X=\Spec A$. We note that 
\[
\ud{\cal{H}om}_{\O_X}(\F, \G) \simeq^a \ud{\cal{H}om}_{\O_X}(\widetilde{\m}\otimes \F, \widetilde{\m}\otimes \G)
\]
by Proposition~\ref{many-functors-sheaf-alhom}\ref{many-functors-sheaf-alhom-2}. Thus, we can assume that both $\F$ and $\G$ are quasi-coherent. Then the claim follows from \ref{hom-alcoh-1} and Lemma~\ref{sheafify-modules}. \medskip

\ref{hom-alcoh-4} is similarly just a consequence of \ref{hom-alcoh-2} and Lemma~\ref{sheafify-modules-almost}.
\end{proof}

\begin{cor}\label{alhom-acoh} Let $X$ be an $R$-scheme. 
\begin{enumerate}[label=\textbf{(\arabic*)}]
\item\label{alhom-acoh-1} Let $\F$ be an almost finitely presented $\O_X$-module, and let $\G$ be an almost coherent $\O_X$-module. Then $\ud{\cal{H}om}_{\O_X}(\F, \G)$ is an almost coherent $\O_X$-module. 
\item\label{alhom-acoh-2} Let $\F^a$ be an almost finitely presented $\O_X^a$-module, and let $\G^a$ be an almost coherent $\O_X^a$-module. Then $\ud{\cal{H}om}_{\O_X^a}(\F^a, \G^a)$ (resp. $\ud{al\cal{H}om}_{\O_X^a}(\F^a, \G^a)$) is an almost coherent $\O_X$-module (resp. $\O_X^a$-module). 
\end{enumerate}
\end{cor}
\begin{proof}
We start the proof by observing that $\ud{\cal{H}om}_{\O_X}(\F, \G) \simeq^a \ud{\cal{H}om}_{\O_X}(\widetilde{\m}\otimes \F, \widetilde{\m}\otimes \G)$ by Proposition~\ref{many-functors-sheaf-alhom}\ref{many-functors-sheaf-alhom-2}. Thus we can assume that both $\F$ and $\G$ are actually quasi-coherent. In that case we use Lemma~\ref{hom-alcoh}\ref{hom-alcoh-1} and Lemma~\ref{criterion-almost-coherent} to reduce the question to showing that $\ud{\cal{H}om}_A(M,N)$ is almost coherent for any almost finitely presented $M$ and almost coherent $N$. However, this has already been done in Corollary~\ref{tensor-hom-coh}. \medskip

Part~\ref{alhom-acoh-2} follows from Part~\ref{alhom-acoh-1} as $\ud{\cal{H}om}_{\O_X^a}(\F^a, \G^a) \simeq \ud{\cal{H}om}_{\O_X}(\F^a_!, \G)$ and $\ud{al\cal{H}om}_{\O_X^a}(\F^a, \G^a)\simeq \ud{\cal{H}om}_{\O_X}(\F, \G)^a$.
\end{proof}

\subsection{Schemes. Approximation of almost finitely presented $\O_X^a$-modules}

One of the defects of our definition of almost finitely presented $\O_X$-modules is that is (Zariski)-local on $X$; we require the existence of an approximation by finitely presented $\O_X$-modules only Zariski-locally on $X$. In particular, this definition is not well adapted to proving global statements such as the Almost Proper Mapping Theorem. We resolve this issue by showing that (on a quasi-compact quasi-separated scheme) any almost finitely presented $\O_X^a$-module can be {\it globally} approximated by finitely presented $\O_X$-modules. \smallskip

As always, we fix a ring $R$ with an ideal $\m$ such that $\m=\m^2$ and $\widetilde{\m}=\m\otimes_R \m$ is $R$-flat. We always do almost mathematics with respect to this ideal.\smallskip

\begin{lemma}\label{lemma:schemes-afp-compact-1} Let $X$ be an $R$-scheme, and $\{\G_i^a\}_{i\in I}$ a filtered diagram of almost quasi-coherent $\O_X^a$-modules. 
    \begin{enumerate}[label=\textbf{(\arabic*)}]
        \item\label{lemma:schemes-afp-compact-1-1} The natural morphism 
        \[
        \gamma^0_\F\colon \colim_I \ud{al\cal{H}om}_{\O_X}(\F^a, \G^a_i) \to \ud{al\cal{H}om}_{\O_X}(\F^a, \colim_I \G^a_i)
        \]
        is injective for an almost finitely generated $\O_X^a$-module $\F^a$;
        \item\label{lemma:schemes-afp-compact-1-2} The natural morphism 
        \[
        \gamma^0_\F\colon \colim_I \ud{al\cal{H}om}_{\O_X}(\F^a, \G^a_i) \to \ud{al\cal{H}om}_{\O_X}(\F^a, \colim_I \G^a_i)
        \]
        is an almost isomorphism for an almost finitely presented $\O_X^a$-module $\F^a$.
    \end{enumerate}
\end{lemma}
\begin{proof}
    The statement is local, so we can assume that $X=\Spec A$ is an affine scheme. Then Lemma~\ref{sheafify-modules-almost} implies that $\F^a \simeq M^a$ and $\G_i^a \simeq N_i^a$ for an almost finitely generated (resp. almost finitely presented) $A$-module $M$. Then \cite[\href{https://stacks.math.columbia.edu/tag/009F}{Tag 009F}]{stacks-project} and Lemma~\ref{hom-alcoh} imply that it suffices to show that 
    \[
    \gamma^0_M \colon \colim_i \rm{alHom}_{A^a}(M^a, N_i^a) \to \rm{alHom}_{A^a}(M^a, \colim N_i^a)
    \]
    is injective (resp. an isomorphism) in $\bf{Mod}_R^a$. But this is exactly Corollary~\ref{cor:compact-criterion}. 
\end{proof}

\begin{cor}\label{cor:schemes-afp-compact} Let $X$ be a quasi-compact and quasi-separated $R$-scheme, and $\{\G_i^a\}_I$ a filtered diagram of almost quasi-coherent $\O_X^a$-modules. 
    \begin{enumerate}[label=\textbf{(\arabic*)}]
        \item The natural morphism 
        \[
        \gamma^0_\F\colon \colim_I \rm{alHom}_{\O_X}(\F^a, \G^a_i) \to \rm{alHom}_{\O_X}(\F^a, \colim_I \G^a_i)
        \]
        is injective for an almost finitely generated $\O_X^a$-module $\F^a$;
        \item The natural morphism 
        \[
        \gamma^0_\F\colon \colim_I \rm{alHom}_{\O_X}(\F^a, \G^a_i) \to \rm{alHom}_{\O_X}(\F^a, \colim_I \G^a_i)
        \]
        is an almost isomorphism for an almost finitely presented $\O_X^a$-module $\F^a$.
    \end{enumerate}
\end{cor}
\begin{proof}
    It formally follows from Lemma~\ref{lemma:alhom-global-sections}, Lemma~\ref{lemma:schemes-afp-compact-1}, and \cite[\href{https://stacks.math.columbia.edu/tag/009F}{Tag 009F}]{stacks-project} (and Corollary~\ref{cor:limits-colimits-almost-sheafy}). 
\end{proof}

\begin{defn} An $\O_X$-module $\F$ is {\it globally almost finitely generated} (resp. {\it globally almost finitely presented}) if, for every finitely generated ideal $\m_0 \subset \m$, there is a quasi-coherent finitely generated (resp. finitely presented) $\O_X$-module $\G$ and a morphism $f\colon \G \to \F$ such that $\m_0(\ker f)=0$, $\m_0(\coker f)=0$.
\end{defn}

\begin{lemma}\label{lemma:schemes-afp-compact-2} Let $X$ be a quasi-compact and quasi-separated $R$-scheme, and $\F$ an almost adically quasi-coherent $\O_X$-module.
    \begin{enumerate}[label=\textbf{(\arabic*)}]
        \item\label{lemma:schemes-afp-compact-2-1} If, for any filtered diagram of adically quasi-coherent $\O_X$-modules $\{\G_i\}_{i\in I}$, the natural morphism 
        \[
         \colim_I \rm{Hom}_{\O_X}(\F, \G_i) \to \rm{Hom}_{\O_X}(\F, \colim_I \G_i)
        \]
        is almost injective, then $\F$ is globally almost finitely generated.
        \item\label{lemma:schemes-afp-compact-2-2} If, for any filtered system of adically quasi-coherent $\O_X$-modules $\{\G_i\}_{i\in I}$, the natural morphism
        \[
         \colim_I \rm{Hom}_{\O_X}(\F, \G_i) \to \rm{Hom}_{\O_X}(\F, \colim_I \G_i)
        \]
        is an almost isomorphism, then $\F$ is globally almost finitely presented.
    \end{enumerate}
\end{lemma}
\begin{proof}
    Lemma~\ref{lemma:alhom-global-sections} and Corollary~\ref{cor:limits-colimits-almost-sheafy} ensure that we can replace $\F$ with $\F^a_!$ without loss of generality. So we may and do assume that $\F$ is are quasi-coherent. Then the proof of Lemma~\ref{lemma:afp-compact-2} works essentially verbatim. We repeat it for the reader's convenience. 
    
    $\textbf{(1)}:$ Note that $\F \simeq \colim_I \F_i$ is a filtered colimit of its finitely generated {\it $\O_X$-submodules} (see \cite[\href{https://stacks.math.columbia.edu/tag/01PG}{Tag 01PG}]{stacks-project}). Therefore, we see that 
    \[
    \colim_I \rm{Hom}_{\O_X}(\F, \F/\F_i) \simeq^a \rm{Hom}_{\O_X}(\F, \colim_I (\F/\F_i)) \simeq 0.
    \]
    Consider an element $\alpha$ of the colimit that has a representative the quotient morphism $\F \to \F/\F_i$ (for some choice of $i$). Then, for every $\e\in \m$, $\e \alpha =0$ in $\colim_I \rm{Hom}_{\O_X}(\F, \F/\F_i)$. Explicitly this means that there is $j\geq i$ such that $\e \F \subset \F_j$. Now note that this property is preserved by replacing $j$ with any $j'>j$. Therefore, for any $\m_0=(\e_1, \dots, \e_n)$, we can find a finitely generated $\O_X$-submodule $\F_i \subset \F$ such that $\m_0\F \subset \F_i$. Therefore, $\F$ is almost finitely generated. \smallskip
    
    $\textbf{(2)}:$ Fix any finitely generated sub-ideal $\m_0=(\e_1, \dots, \e_n) \subset \m$. We use \cite[\href{https://stacks.math.columbia.edu/tag/01PJ}{Tag 01PJ}]{stacks-project} to write $\F\simeq \colim_{\Lambda} \F_\lambda$ as a filtered colimit of {\it finitely presented} $\O_X$-modules. By assumption, the natural morphism
    \[
    \colim_{\Lambda} \rm{Hom}_{\O_X}(\F, \F_\lambda) \to \rm{Hom}_{\O_X}(\F, \colim_{\Lambda} \F_\lambda)=\rm{Hom}_{\O_X}(\F, \F)
    \]
    is an almost isomorphism. In particular, $\e_i \rm{Id}_{\F}$ is in the image of this map for every $i=1,\dots, n$. This means that, for every $\e_i$, there is $\lambda_i\in \Lambda$ and a morphism $g_i \colon \F \to \F_{\lambda_i}$ such that the composition 
    \[
    f_{\lambda_i} \circ g_i=\e_i \rm{Id}_\F,
    \]
    where $f_{\lambda_i}\colon \F_{\lambda_i} \to \F$ is the natural morphism to the colimit. Note that existence of such $g_i$ is preserved by replacing $\lambda_i$ with any $\lambda'_i \geq \lambda_i$. Therefore, using that $\{\F_\lambda\}$ is a filtered diagram, we can find an index $\lambda$ with maps 
    \[
    g_i \colon \F \to \F_{\lambda}
    \]
    such that $f_\lambda \circ g_i = \e_i \rm{Id}_{\F}$. Now we consider the morphism 
    \[
    G_i\coloneqq g_i\circ f_\lambda - \e_i \rm{Id}_{\F_\lambda} \colon \F_\lambda \to \F_\lambda.
    \]
    Note that $\rm{Im}(G_i) \subset \ker(f_\lambda)$ because
    \[
    f_\lambda \circ g_i \circ f_\lambda - f_\lambda \e_i \rm{Id}_{\F_\lambda} = \e_i f_\lambda - \e_i f_\lambda = 0.
    \]
    We also have that $\e_i \ker(f_\lambda) \subset \rm{Im}(G_i)$ because $G_i|_{\ker(f_\lambda)} = \e_i \rm{Id}$. Therefore, $\sum_i \rm{Im}(G_i)$ is a quasi-coherent finitely generated $\O_X$-module such that 
    \[
    \m_0(\ker f_\lambda) \subset \sum_i \rm{Im}(G_i) \subset \ker(f_\lambda).
    \]
    Therefore, $f\colon \F'\coloneqq \F_\lambda/(\sum_i \rm{Im}(G_i)) \to \F$ is a morphism such that $\F'$ is finitely presented, $\m_0(\ker f)=0$, and $\m_0(\coker f)=0$. Since $\m_0 \subset \m$ was an arbitrary finitely generated sub-ideal, we conclude that $\F$ is globally almost finitely presented. 
\end{proof}

\begin{cor}\label{cor:approximate-afpr} Let $X$ be a quasi-compact and quasi-separated $R$-scheme, and let $\F$ be an almost quasi-coherent $\O_X$-module. Then $\F$ is almost finitely presented (resp. almost finitely generated) if and only if for any finitely generated ideal $\m_0\subset \m$ there is a morphism $f\colon \G \to \F$ such that $\G$ is a quasi-coherent finitely presented (resp. finitely generated) $\O_X$-module , $\m_0(\ker f)=0$ and $\m_0(\coker f)=0$.
\end{cor}
\begin{proof}
    Corollary~\ref{cor:schemes-afp-compact} ensures that $\F$ satisfies the conditions of Lemma~\ref{lemma:schemes-afp-compact-2}. Lemma~\ref{lemma:schemes-afp-compact-2} now gives the desired result.  
\end{proof}

\begin{cor}\label{cor:schemes-compact-criterion-2} Let $X$ be a quasi-compact and quasi-separated $R$-scheme, and $\F^a$ an almost quasi-coherent $\O^a_X$-module. Then
    \begin{enumerate}[label=\textbf{(\arabic*)}]
        \item $\F^a$ is almost finitely generated if and only if, for every filtered diagram $\{\G^a_i\}_{i\in I}$ of almost quasi-coherent $\O^a_X$-modules, the natural morphism
        \[
         \colim_I \rm{alHom}_{\O_X^a}(\F^a, \G^a_i) \to \rm{alHom}_{\O_X^a}(\F^a, \colim_I \G^a_i)
        \]
        is injective in $\bf{Mod}_R^a$;
        \item $\F^a$ is almost finitely presented if and only if, for every filtered diagram $\{\G^a_i\}$ of almost quasi-coherent $\O^a_X$-modules, the natural morphism
        \[
         \colim_I \rm{alHom}_{\O_X^a}(\F^a, \G^a_i) \to \rm{alHom}_{\O_X^a}(\F^a, \colim_I \G^a_i)
        \]
        is an isomorphism in $\bf{Mod}_R^a$;    
    \end{enumerate}
\end{cor}

\subsection{Schemes. Derived category of almost coherent $\O_X^a$-modules}\label{section:derived-category-schemes}

The goal of this section is to define different versions of the ``derived category of almost coherent sheaves''. Namely, we define the categories $\bf{D}_{acoh}(X)$, $\bf{D}_{qc, acoh}(X)$, and $\bf{D}_{acoh}(X)^a$. Then we show that many functors of interest preserve almost coherence in an appropriate sense. 

\begin{defn}\label{defn:derived-almost-category-1} We define $\bf{D}_{aqc}(X)$ (resp. $\bf{D}_{aqc}(X)^a$) to be the full triangulated subcategory of $\bf{D}(X)$ (resp. $\bf{D}(X)^a$) consisting of complexes with almost quasi-coherent cohomology sheaves. 
\end{defn} 

\begin{defn}\label{defn:derived-almost-category-2} We define $\bf{D}_{acoh}(X)$ (resp. $\bf{D}_{qc, acoh}(X)$, resp. $\bf{D}_{acoh}(X)^a$) to be the full triangulated subcategory of $\bf{D}(X)$ (resp. $\bf{D}(X)$, resp. $\bf{D}(X)^a$) consisting of complexes with almost coherent (resp. quasi-coherent and almost coherent, resp. almost coherent) cohomology sheaves. 
\end{defn} 

\begin{rmk} The definition above makes sense as the categories $\bf{Mod}^{\acoh}_X$, $\bf{Mod}_X^{\qc, \text{acoh}}$, and $\bf{Mod}_{X^a}^{\text{acoh}}$ are weak Serre subcategories of $\bf{Mod}_X$, $\bf{Mod}_X$, and $\bf{Mod}_X^a$ respectively. 
\end{rmk}

Now suppose that $X=\Spec A$ is an affine $R$-scheme. Then we note that the functor
\[
\widetilde{(-)}\colon \bf{Mod}_A \to \bf{Mod}_X
\]
is additive and exact, thus it can be easily derived to the functor 
\[
\widetilde{(-)}\colon \bf{D}(A) \to \bf{D}_{qc}(X).
\]
\begin{lemma}\label{derived-schemes} Let $X=\Spec A$ be an affine $R$-scheme. Then the functor 
\[
\widetilde{(-)} \colon \bf{D}(A) \to \bf{D}_{qc}(X)
\] is a $t$-exact equivalence of triangulated categories\footnote{with respect to the standard $t$-structures.} with the quasi-inverse given by $\bf{R}\Gamma(X, -)$. Moreover, these two functors induce the equivalence
\[
        \begin{tikzcd}
          \widetilde{(-)} \colon \bf{D}^*_{acoh}(A)\arrow[r, swap, shift right=.75ex] & \bf{D}^*_{qc, acoh}(X)\arrow[l, swap, shift right=.75ex] \colon \bf{R}\Gamma(X, -)
        \end{tikzcd}
\]
for any $*\in \{``\text{ ''}, +, -, b\}$.
\end{lemma}
\begin{proof}
The first part is just \cite[\href{https://stacks.math.columbia.edu/tag/06Z0}{Tag 06Z0}]{stacks-project}. In particular, it shows that $\rm{H}^i(\bf{R}\Gamma(X, \F))\simeq \rm{H}^0(X, \cal{H}^i(\F))$ for any $\F\in \bf{D}_{qc}(X)$. Now Lemma~\ref{criterion-almost-coherent} implies that $\cal{H}^i(\F)$ is almost coherent if and only if so is $\rm{H}^0(X, \cal{H}^i(\F))$. So the functor $\bf{R}\Gamma(X, -)$ sends $\bf{D}^*_{qc, acoh}(X)$ to $\bf{D}^*_{acoh}(A)$. \smallskip

We also observe that the functor $\widetilde{(-)}$ clearly sends $\bf{D}_{acoh}(A)$ to $\bf{D}_{qc, acoh}(X)$. Thus, we conclude that $\widetilde{(-)}$ and $\bf{R}\Gamma(X, -)$ induce an equivalence between $\bf{D}_{acoh}(A)$ and $\bf{D}_{qc, acoh}(X)$. The bounded versions follow from $t$-exactness of both functors. 
\end{proof}

\begin{lemma}\label{quotient-derived-acoh-je} Let $X=\Spec A$ be an affine $R$-scheme. Then the almostification functor 
\[
(-)^a\colon \bf{D}^*_{qc}(X) \to \bf{D}^*_{aqc}(X)^a
\] induces an equivalence
$
\bf{D}^*_{qc}(X)/\bf{D}^*_{qc, \Sigma_X}(X) \xrightarrow{\sim} \bf{D}^*_{aqc}(X)^a
$
for any $*\in \{``\text{ ''}, +, -, b\}$. Similarly, the induced functor 
\[
\bf{D}^*_{qc, acoh}(X)/\bf{D}^*_{qc, \Sigma_X}(X) \xrightarrow{\sim} \bf{D}^*_{acoh}(X)^a
\]
is an equivalence for any $*\in \{``\text{ ''}, +, -, b\}$.
\end{lemma}
\begin{proof}
The functor $(-)_!\colon \bf{D}^*_{aqc}(X)^a \to \bf{D}^*_{qc}(X)$ gives the left adjoint to $(-)^a$ such that $\rm{Id}\to (-)_!\circ(-)^a$ is an isomorphism and the kernel of $(-)^a$ consists exactly of the morphisms $f$ such that $\rm{cone}(f)\in \bf{D}_{qc, \Sigma_X}(X)$. Thus, the dual version of \cite[Proposition 1.3]{GZ} finishes the proof of the first claim. The proof of the second claim is similar once one notices that $\widetilde{M^a}$ is almost coherent for any almost coherent $A^a$-module $M^a$. The latter fact follows from Lemma~\ref{criterion-almost-coherent}.
\end{proof}

Lemma~\ref{criterion-almost-coherent} ensures that $\bf{D}(A)^a\simeq \bf{D}(A)/\bf{D}_{\Sigma_A}(A)$. As $\widetilde{(-)}$ clearly sends $\bf{D}_{\Sigma_A}(A)$ to $\bf{D}^*_{qc, \Sigma_X}(X)$, so we conclude that it induces a functor 
\[
\widetilde{(-)} \colon \bf{D}^*(A)^a \to \bf{D}_{aqc}^*(X)^a.
\]

\begin{thm}\label{derived-schemes-2} Let $X=\Spec A$ be an affine $R$-scheme. Then the functor 
\[
\widetilde{(-)} \colon \bf{D}(A)^a \to \bf{D}_{aqc}(X)^a
\] is a $t$-exact equivalence of triangulated categories with the quasi-inverse given by $\bf{R}\Gamma(X, -)$. Moreover, these two functors induce equivalences
\[
        \begin{tikzcd}
          \widetilde{(-)} \colon \bf{D}^*_{acoh}(A)^a\arrow[r, swap, shift right=.75ex] & \bf{D}^*_{acoh}(X)^a\arrow[l, swap, shift right=.75ex] \colon \bf{R}\Gamma(X, -)
        \end{tikzcd}
\]
for any $*\in \{``\text{ ''}, +, -, b\}$.
\end{thm}
\begin{proof}
We note that Lemma~\ref{derived-schemes} ensures that $\widetilde{(-)}\colon \bf{D}^*_{qc, acoh}(X) \to \bf{D}^*_{acoh}(X)^a$ is an equivalence with quasi-inverse $\bf{R}\Gamma(X, -)$. Moreover, $(-)^a$ induces an equivalence between $\bf{D}_{\Sigma_A}(A)$ and $\bf{D}_{qc, \Sigma_X}(X)$; we leave the verification to the interested reader. Thus, Lemma~\ref{quotient-derived-acoh-je} ensures that $\widetilde{(-)}$ gives an equivalence
\[
\bf{D}(A)^a \simeq \bf{D}(A)/\bf{D}_{\Sigma_A}(A) \xrightarrow{\sim} \bf{D}_{qc}(X)/\bf{D}_{qc, \Sigma_X}(X) \simeq \bf{D}_{aqc}(X)^a.
\]
Its quasi-inverse is given by the functor  $\bf{R}\Gamma(X, -) \colon \bf{D}_{aqc}(X)^a \to \bf{D}(A)^a$ by Proposition~\ref{derived-pushforward}. \smallskip

The version with almost coherent cohomology sheaves is similar to the analogous statement from Lemma~\ref{derived-schemes}.
\end{proof}

\begin{lemma}\label{pullback-almost-coh-derived} Let $f\colon X \to Y$ be a morphism of $R$-schemes.  
\begin{enumerate}[label=\textbf{(\arabic*)}]
\item\label{pullback-almost-coh-derived-1} Suppose that $X=\Spec B$, $Y=\Spec A$ are affine $R$-schemes. Then $\bf{L}f^*(\widetilde{M^a})$ is functorially isomorphic to $\widetilde{M^a\otimes^L_{A^a} B^a}$ for any $M^a\in \bf{D}(A)^a$.
\item\label{pullback-almost-coh-derived-2} The functor $\bf{L}f^*$ carries an object of $\bf{D}^*_{aqc}(Y)$ to an object of $\bf{D}^*_{aqc}(X)$ for $*\in \{``\text{ ''}, -\}$.

\item\label{pullback-almost-coh-derived-3} The functor $\bf{L}f^*$ carries an object of $\bf{D}^*_{aqc}(Y)^a$ to an object of $\bf{D}^*_{aqc}(X)^a$ for $*\in \{``\text{ ''}, -\}$.

\item\label{pullback-almost-coh-derived-4} Suppose that $X$ and $Y$ are almost coherent $R$-schemes. Then the functor $\bf{L}f^*$ carries an object of $\bf{D}^-_{qc, acoh}(Y)$ (resp. $\bf{D}^-_{acoh}(Y)$) to an object of $\bf{D}^-_{qc, acoh}(X)$ (resp. $\bf{D}^-_{acoh}(X)$).

\item\label{pullback-almost-coh-derived-5} Suppose that $X$ and $Y$ are almost coherent $R$-schemes. Then the functor $\bf{L}f^*$ carries an object of $\bf{D}^-_{acoh}(Y)^a$ to an object of $\bf{D}^-_{acoh}(X)^a$.
\end{enumerate}
\end{lemma}
\begin{proof}
We start with \ref{pullback-almost-coh-derived-1}. We use Proposition~\ref{derived-pullback} to see that $\bf{L}f^*(\widetilde{M^a}) \simeq (\bf{L}f^*(\widetilde{M}))^a$. Proposition~\ref{derived-base-change} ensures that $(\widetilde{M\otimes^L_{A} B})^a\simeq \widetilde{M^a\otimes^L_{A^a} B^a}$, so it suffices to show $\bf{L}f^*(\widetilde{M}) \simeq \widetilde{M\otimes^L_{A} B}$. But this is a classical fact about quasi-coherent sheaves. \smallskip

Now we show \ref{pullback-almost-coh-derived-2}. We note that Lemma~\ref{derived-pull-projection} implies that $\bf{L}f^*(\widetilde{\m}\otimes \F) \simeq \widetilde{\m}\otimes \bf{L}f^*(\F)$ for any $\F\in \bf{D}(Y)$. Thus, we can replace  $\F$ with $\widetilde{\m}\otimes \F$ to assume that it is quasi-coherent. Then it is a standard fact that $\bf{L}f^*$ sends $\bf{D}^*_{qc}(Y)$ to $\bf{D}^*_{qc}(X)$ for $*\in \{``\text{ ''}, -\}$.  \smallskip

\ref{pullback-almost-coh-derived-3} follows from \ref{pullback-almost-coh-derived-2} by noting that $\bf{L}f^*(\F^a)\simeq \left(\bf{L}f^*\left(\F^a_!\right)\right)^a$ according to Proposition~\ref{derived-pullback}. \smallskip

To prove \ref{pullback-almost-coh-derived-4}, we again use the isomorphism $\bf{L}f^*(\widetilde{\m}\otimes \F) \simeq \widetilde{\m}\otimes \bf{L}f^*(\F)$ to assume that $\F$ is in $\bf{D}^-_{qc, acoh}(X)$. Lemma~\ref{derived-schemes} guarantees that there exists $M\in \bf{D}_{coh}^-(A)$ such that $\widetilde{M}\simeq \F$. Thus Part~\ref{pullback-almost-coh-derived-1} and Lemma~\ref{criterion-almost-coherent} ensure that it is sufficient to show that $M^a\otimes^L_{A^a} B^a\simeq (M\otimes^L_A B)^a$ has almost finitely presented cohomology modules. This is exactly the content of Corollary~\ref{derived-trivial-base-change}. \medskip

\ref{pullback-almost-coh-derived-5} follows from \ref{pullback-almost-coh-derived-4} as $\bf{L}f^*(\F^a)\simeq \left(\bf{L}f^*\left(\F^a_!\right)\right)^a$.
\end{proof}

\begin{lemma}\label{tensor-almost-coh-derived} Let $X$ be an $R$-scheme.  
\begin{enumerate}[label=\textbf{(\arabic*)}]

\item\label{tensor-almost-coh-derived-1} Suppose that $X=\Spec A$ is an affine $R$-scheme. Then $\widetilde{M^a}\otimes^L_{\O_X^a}\widetilde{N^a}$ is functorially isomorphic to $\widetilde{M^a\otimes^L_{A^a} N^a}$ for any $M^a, N^a\in \bf{D}(A)^a$. 
\item\label{tensor-almost-coh-derived-2} Let $\F, \G\in \bf{D}^*_{aqc}(X)$, then $\F\otimes^L_{\O_X}\G \in \bf{D}_{aqc}(X)$ for $*\in \{`` \text{ ''}, -\}$. 
\item\label{tensor-almost-coh-derived-3} Let $\F^a, \G^a\in \bf{D}^*_{aqc}(X)^a$, then $\F^a\otimes^L_{\O_X^a}\G^a \in \bf{D}_{aqc}(X)^a$ for $*\in \{`` \text{ ''}, -\}$. 
\item\label{tensor-almost-coh-derived-4} Suppose that $X$ is an almost coherent $R$-scheme, and let $\F, \G\in \bf{D}^-_{qc, acoh}(X)$ (resp. $\bf{D}^-_{acoh}(X)$). Then $\F\otimes^L_{\O_X}\G \in \bf{D}^-_{qc, acoh}(X)$ (resp. $\bf{D}^-_{acoh}(X)$). 
\item\label{tensor-almost-coh-derived-5} Suppose that $X$ is an almost coherent $R$-scheme, and let $\F^a, \G^a\in \bf{D}^-_{acoh}(X)^a$. Then $\F^a\otimes^L_{\O_X^a}\G^a \in \bf{D}^-_{acoh}(X)^a$. 
\end{enumerate}
\end{lemma}
\begin{proof}
The proof is basically identical to that of Lemma~\ref{pullback-almost-coh-derived} and is left to the reader. We only mention that one has to use Proposition~\ref{tensor-product-coh} in place of Corollary~\ref{derived-trivial-base-change}. 
\end{proof}

\begin{lemma}\label{pushforward-almost-coh-derived} Let $f\colon X\to Y$ be a quasi-compact and quasi-separated morphism of $R$-schemes. Suppose that $Y$ is quasi-compact. Then  
\begin{enumerate}[label=\textbf{(\arabic*)}]
\item\label{pushforward-almost-coh-derived-1} The functor $\bf{R}f_*$ carries $\bf{D}^*_{aqc}(X)$ to $\bf{D}^*_{aqc}(Y)$ for any $*\in \{``\text{ ''}, -, +, b\}$.
\item\label{pushforward-almost-coh-derived-2} The functor $\bf{R}f_*$ carries $\bf{D}^*_{aqc}(X)^a$ to $\bf{D}^*_{aqc}(Y)^a$ for any $*\in \{``\text{ ''}, -, +, b\}$.
\end{enumerate}
\end{lemma}
\begin{proof}
    Proposition~\ref{derived-pushforward} guarantees that $(\bf{R}f_*\F)^a \simeq \bf{R}f_*\F^a$. Since $(\widetilde{\m}\otimes \F)^a \simeq \F^a$, we see that it suffices to show that the functor
    \[
    \bf{R}f_*(\widetilde{\m} \otimes -)
    \]
    carries $\bf{D}^*_{aqc}(X)$ to $\bf{D}^*_{aqc}(Y)$ for any $*\in \{``\text{ ''}, -, +, b\}$. Since $\widetilde{\m}\otimes \F$ is in $\bf{D}_{qc}(X)$, we conclude that it is enough to show that $\bf{R}f_*(-)$ carries $\bf{D}^*_{qc}(X)$ to $\bf{D}^*_{qc}(Y)$ for any $*\in \{``\text{ ''}, -, +, b\}$. This is proven in \cite[\href{https://stacks.math.columbia.edu/tag/08D5}{Tag 08D5}]{stacks-project}.
\end{proof}

Before going to the case of the derived Hom-functors, we recall the construction of the functorial map
\[
\psi\colon \widetilde{\bf{R}\rm{Hom}_A(M, N)} \to \bf{R}\ud{\cal{H}om}_{\O_X}(\widetilde{M}, \widetilde{N})
\]
for any $M\in \bf{D}^-(A)$, $N\in \bf{D}^+(A)$, and an affine scheme $X=\Spec A$. For this, we recall that the functor $\widetilde{(-)}$ is left adjoint to the global section functor $\Gamma(X, -)$ on the abelian level. Thus, after deriving these functors, we see that $\widetilde{-}$ is left adjoint to $\bf{R}\Gamma(X, -)$. Therefore, for any $\F \in \bf{D}(X)$, there is a canonical morphism $\widetilde{\bf{R}\Gamma(X, \F)} \to \F$. We apply it to $\F = \bf{R}\ud{\cal{H}om}_{\O_X}(\widetilde{M}, \widetilde{N})$ to obtain the desired morphism
\[
\psi\colon \widetilde{\bf{R}\rm{Hom}_A(M, N)} \to \bf{R}\ud{\cal{H}om}_{\O_X}(\widetilde{M}, \widetilde{N}).
\]

\begin{lemma}\label{hom-alcoh-derived} Let $X$ be an almost coherent $R$-scheme.  
\begin{enumerate}[label=\textbf{(\arabic*)}]
\item\label{hom-alcoh-derived-1} Suppose $X=\Spec A$ is an affine $R$-scheme. The canonical map 
\[
\psi \colon \widetilde{\bf{R}\rm{Hom}_A(M,N)} \to \bf{R}\ud{\cal{H}om}_{\O_X}(\widetilde{M}, \widetilde{N})
\]
is an almost isomorphism for $M\in \bf{D}^-_{acoh}(A)$, $N\in \bf{D}^+(A)$. 
\item\label{hom-alcoh-derived-2} Suppose $X=\Spec A$ is an affine $R$-scheme. There is a functorial isomorphism 
\[
\widetilde{\bf{R}\rm{alHom}_{A^a}(M^a,N^a)} \simeq \bf{R}\ud{al\cal{H}om}_{\O_X^a}(\widetilde{M^a}, \widetilde{N^a})
\]
for $M^a\in \bf{D}^-_{acoh}(A)^a$, $N^a\in \bf{D}^+(A)^a$. We also get a functorial almost isomorphism
\[
\widetilde{\bf{R}\rm{Hom}_{A^a}(M^a,N^a)} \simeq^a \bf{R}\ud{\cal{H}om}_{\O_X^a}(\widetilde{M^a}, \widetilde{N^a})
\]
for $M\in \bf{D}^-_{acoh}(A)$, $N\in \bf{D}^+(A)$. 
\item\label{hom-alcoh-derived-3} Suppose $\F\in \bf{D}^-_{acoh}(X)$ and $\G\in \bf{D}^+_{aqc}(X)$. Then $\bf{R}\ud{\cal{H}om}_{\O_X}(\F, \G) \in \bf{D}^+_{aqc}(X)$.
\item\label{hom-alcoh-derived-4} Suppose $\F^a\in \bf{D}^-_{acoh}(X)^a$ and $\G^a\in \bf{D}^+_{aqc}(X)^a$. Then $\bf{R}\ud{\cal{H}om}_{\O_X^a}(\F^a, \G^a)\in \bf{D}^+_{aqc}(X)$ and $\bf{R}\ud{al\cal{H}om}_{\O_X^a}(\F^a, \G^a) \in \bf{D}^+_{aqc}(X)^a$.
\end{enumerate}
\end{lemma}
\begin{proof}
We start with \ref{hom-alcoh-derived-1}. We use the convergent compatible spectral sequences
\[
\rm{E}^{p,q}_2=\widetilde{\rm{Ext}^p_A(\rm{H}^{-q}(M), N)} \Rightarrow \widetilde{\rm{Ext}^{p+q}_A(M, N)}
\]
\[
\rm{E'}^{p,q}_2=\cal{E}xt^p_{\O_X}\left(\widetilde{\rm{H}^{-q}(M)}, \widetilde{N}\right) \Rightarrow \cal{E}xt^{p+q}_{\O_X}\left(\widetilde{M}, \widetilde{N}\right)
\]
to reduce to the case when $M\in \bf{Mod}^{\acoh}_A$ is an $A$-module concentrated in degree $0$. Similarly, we use the compatible spectral sequences
\[
\rm{E}^{p,q}_2=\widetilde{\rm{Ext}^q_A(M, \rm{H}^{p}(N))} \Rightarrow \widetilde{\rm{Ext}^{p+q}_A(M, N)}
\]
\[
\rm{E'}^{p,q}_2=\cal{E}xt^q_{\O_X}(\widetilde{M}, \widetilde{\rm{H}^{p}(N)}) \Rightarrow \cal{E}xt^{p+q}_{\O_X}(\widetilde{M}, \widetilde{N})
\] 
to assume that $N\in \bf{Mod}_A$. Thus, the claim boils down to showing that the natural map 
\[
\widetilde{\rm{Ext}^{p}_A(M, N)} \to \cal{E}xt^{p}_{\O_X}(\widetilde{M}, \widetilde{N})
\]
is an almost isomorphism for any $M\in \bf{Mod}^{\acoh}_A$, $N\in \bf{Mod}_A$, and $p\geq 0$. Lemma~\ref{almost-zero-sheaf} says that it is sufficient to show that the kernel and cokernel are annihilated by any finitely generated sub-ideal $\m_0\subset \m$. 

Recall that, for any $\O_X$-modules $\F$, $\G$, the sheaf $\cal{E}xt^p_{\O_X}(\F, \G)$ is canonically isomorphic to the sheafification of the presheaf
\[
U\mapsto \rm{Ext}^p_{\O_U}(\F|_U, \G|_U).
\]
Thus, in order to show that the map $\widetilde{\rm{Ext}^{p}_A(M, N)} \to \cal{E}xt^{p}_{\O_X}(\widetilde{M}, \widetilde{N})$ is an almost isomorphism, it suffices to show that
\[
\rm{Ext}^{p}_A(M, N) \otimes_A A_f \to \rm{Ext}^p_{\O_{X_f}}(\widetilde{M_f}, \widetilde{N_f})
\]
is an almost isomorphism. Now we use canonical isomorphisms
\begin{align*}
    \rm{Ext}^p_{\O_{X_f}}(\widetilde{M_f}, \widetilde{N_f})&\simeq \rm{Hom}_{\bf{D}(X_f)}(\widetilde{M_f}, \widetilde{N_f}[p]) \\
    & \simeq \rm{Hom}_{\bf{D}(A_f)}(M_f, N_f[p]) \\
    & \simeq \rm{Ext}^p_{A_f}(M_f, N_f), 
\end{align*}
where the second isomorphism uses that $(\widetilde{-})$ induces a $t$-exact equivalence $(\widetilde{-}) \colon \bf{D}(A_f) \xr{~} \bf{D}_{qc}(\Spec A_f)$. Thus, the question boils down to showing that the natural map
\[
\rm{Ext}^{p}_A(M, N) \otimes_A A_f \to \rm{Ext}^p_{A_f}(M_f, N_f)
\]
is an almost isomorphism. This follows from Proposition~\ref{base-change-hom-derived}. \smallskip

\ref{hom-alcoh-derived-2} formally follows from \ref{hom-alcoh-derived-1} by using Proposition~\ref{derived-al-hom-sheaf}\ref{derived-al-hom-sheaf-1}.  \medskip

\ref{hom-alcoh-derived-3} is also a basic consequence of \ref{hom-alcoh-derived-2}. Indeed, the claim is local, so we can assume that $X=\Spec A$ is an affine $R$-scheme. In that case, we use Theorem~\ref{derived-schemes-2} to say that $\F\simeq \widetilde{M}$, $\G\simeq \widetilde{N}$ for some $M\in \bf{D}^-_{acoh}(A)$, $N\in \bf{D}^+(A)$. Then $\bf{R}\ud{\cal{H}om}_{\O_X}(\F, \G) \simeq \widetilde{\bf{R}\rm{Hom}_A(M,N)}$ by \ref{hom-alcoh-derived-2}, and the latter complex has almost quasi-coherent cohomology sheaves by design. \medskip 

\ref{hom-alcoh-derived-4} easily follows from \ref{hom-alcoh-derived-3} and the isomorphisms
\[
\bf{R}\ud{\cal{H}om}_{\O_X^a}(\F^a, \G^a) \simeq \bf{R}\ud{\cal{H}om}_{\O_X}(\F^a_!, \G)
\]
\[
\bf{R}\ud{al\cal{H}om}_{\O_X^a}(\F^a, \G^a) \simeq \bf{R}\ud{\cal{H}om}_{\O_X}(\F^a_!, \G)^a
\]
that come from Lemma~\ref{derived-hom-alg-sheaf-1} and Definition~\ref{defn-almost-hom-sheaf-derived}. 
\end{proof}

\begin{cor}\label{alhom-acoh-derived} Let $X$ be an almost coherent $R$-scheme. 
\begin{enumerate}[label=\textbf{(\arabic*)}]
\item\label{alhom-acoh-derived-1} Let $\F\in \bf{D}^-_{aqc, acoh}(X)$, $\G\in \bf{D}^+_{aqc, acoh}(X)$. Then $\bf{R}\ud{\cal{H}om}_{\O_X}(\F, \G) \in \bf{D}^+_{aqc, acoh}(X)$. 
\item\label{alhom-acoh-derived-2} Let $\F^a\in \bf{D}^-_{acoh}(X)^a$, $\G^a\in \bf{D}^+_{acoh}(X)^a$. Then $\bf{R}\ud{al\cal{H}om}_{\O_X^a}(\F^a, \G^a) \in \bf{D}^+_{acoh}(X)^a$. 
\end{enumerate}
\end{cor}
\begin{proof} The question is local on $X$, so we can assume that $X=\Spec A$ is affine. Then Lemma~\ref{hom-alcoh-derived}, Theorem~\ref{derived-schemes-2}, and Lemma~\ref{criterion-almost-coherent} reduce both questions to showing that $\bf{R}\rm{Hom}_A(M, N)\in \bf{D}^+_{acoh}(A)$ for $M\in \bf{D}^-_{acoh}(A)$ and $N\in \bf{D}^+_{acoh}(A)$. This is the content of Proposition~\ref{alHom-derived-coh}. 
\end{proof}

\begin{prop}\label{prop:projection-formula-quasicoh} Let $f\colon X\to Y$ be a quasi-compact quasi-separated morphism of $R$-schemes, $\F^a\in \bf{D}_{aqc}(X)^a$, and $\G\in \bf{D}_{aqc}(Y)^a$. Then the projection morphism (see the discussion before Proposition~\ref{prop:projection-formula-perfect})
\[
\rho\colon \bf{R}f_*(\F^a)\otimes^L_{\O_Y^a} \G^a \to \bf{R}f_*(\F^a\otimes^L_{\O_X^a} \bf{L}f^*(\G^a))
\]
is an isomorphism in $\bf{D}(Y)^a$.
\end{prop}
\begin{proof}
    Proposition~\ref{derived-tensor-product-sheaf}, Proposition~\ref{derived-pullback}, and Proposition~\ref{derived-pushforward} imply that we can replace $\F^a$ (resp. $\G^a$) with $\F^a_!\in \bf{D}_{qc}(X)^a$ (resp. $\G^a_!\in \bf{D}_{qc}(Y)^a$). So it suffices to show the analogous result for modules with {\it quasi-coherent} cohomology sheaves. This is proven in \cite[\href{https://stacks.math.columbia.edu/tag/08EU}{Tag 08EU}]{stacks-project}. 
\end{proof}

\subsection{Formal Schemes. The category of almost coherent $\O_\X^a$-modules}\label{acoh-sheaves-formal}

In this Section, we discuss the notion of almost coherent sheaves on ``good'' formal schemes. One of the main complications compared to the case of usual schemes is that there is no good notion of a ``quasi-coherent'' sheaf on a formal scheme. Namely, even though there is a notion of adically quasi-coherent sheaves on a large class of formal schemes due to \cite[\textsection I.3]{FujKato}, this notion does not behave well. In particular, this category is not a weak Serre subcategory of $\O_\X$-modules for ``nice'' a formal scheme $\X$. \smallskip

Another (related) difficulty comes from the lack of the Artin--Rees lemma for not finitely generated modules. More precisely, many operations with adically quasi-coherent sheaves require taking completions, but this operation is usually not exact without the presence of the Artin--Rees lemma. \smallskip 

We deal with this issue by using a version of the Artin--Rees lemma (Lemma~\ref{Artin-Rees}) for almost finitely generated modules over ``good'' rings. The presence of the Artin--Rees lemma suggests that it is reasonable to expect that we might have a good notion of adically quasi-coherent, almost coherent $\O_\X$-modules on some ``good'' class of formal schemes. \smallskip

We start by giving a setup in which we can develop the theory of almost coherent sheaves

\begin{setup}\label{set-up3} We fix a ring $R$ with a finitely generated ideal $I$ such that $R$ is $I$-adically complete, $I$-adically topologically universally adhesive\footnote{This means that the algebra $R\langle X_1, \dots, X_n\rangle[T_1, \dots, T_m]$ is $I$-adically adhesive for any $n$ and $m$}, and $I$-torsion free with an ideal $\m$ such that $I\subset \m$, $\m^2=\m$ and $\widetilde{\m}\coloneqq \m\otimes_R \m$ is $R$-flat. 
\end{setup}

The basic example of such a ring is a complete microbial\footnote{A valuation ring $R$ is microbial if there is a non-zero topologically nilpotent element $\varpi \in R$. Any such element is called a pseudo-uniformizer.} valuation ring $R$ with algebraically closed fraction field $K$. We pick a pseudo-uniformizer $\varpi$ and define $I\coloneqq (\varpi)$, $\m\coloneqq \cup_{i=1}^{\infty} (\varpi^{1/n})$ for some compatible choice of roots of $\varpi$. We note that $R$ is topologically universally adhesive by \cite[Theorem 7.3.2]{FGK}. \smallskip 

We note that the assumptions in Setup~\ref{set-up3} imply that any finitely presented algebra over a topologically finitely presented $R$-algebra is coherent and $I$-adically adhesive. Coherence follows from \cite[Proposition 7.2.2]{FGK} and adhesiveness follows from the definition. In what follows, we will use those facts without saying. \smallskip 



In what follows, $\X$ always means a topologically finitely presented formal $R$-scheme. We will denote by $\X_k\coloneqq \X\times_{\Spf R}\Spec R/I^{k+1}$ the ``reduction'' schemes. They come equipped with a closed immersion $i_k\colon \X_k \to \X$. Also, given any $\O_\X$-module $\F$, we will always denote its ``reduction'' $i_k^*\F$ by $\F_k$. \smallskip

\begin{defn}\label{defn:adically-quasi-coherent-sheaves} \cite[Definition I.3.1.3]{FujKato} An $\O_\X$-module $\F$ on a formal scheme $\X$ of finite ideal type is called {\it adically quasi-coherent} if $\F \to \lim_n \F_n$ is an isomorphism and, for any open formal subscheme $\sU \subset \F$ and any ideal of definition $\mathcal{I}$ of finite type, the sheaf $\F/\mathcal{I}\F$ is a quasi-coherent sheaf on the scheme $(\sU, \O_{\sU}/\mathcal{I})$. 

We denote by $\bf{Mod}_\X^{qc}$ the full subcategory of $\bf{Mod}_\X$ consisting of adically quasi-coherent $\O_\X$-modules. 
\end{defn}

\begin{defn} We say that an $\O_\X^a$-module $\F^a$ is {\it almost adically quasi-coherent} if $\F^a_!\simeq \widetilde{\m}\otimes \F$ is an adically quasi-coherent $\O_\X$-module. We denote by $\bf{Mod}^{\aqc}_{\X^a}$ the full subcategory of $\bf{Mod}_{\X^a}$ consisting of almost adically quasi-coherent $\O_\X$-modules. \smallskip

We say that an $\O_\X$-module $\F$ is {\it almost adically quasi-coherent} if $\F^a$ is an almost quasi-coherent $\O_\X^a$-module. We denote by $\bf{Mod}_\X^{\aqc}$ the full subcategory of $\bf{Mod}_\X$ consisting of adically quasi-coherent $\O_\X$-modules.
\end{defn}

\begin{rmk}\label{problem} In general, we can not say that an adically quasi-coherent $\O_\X$-module $\F$ is almost adically quasi-coherent. The problem is that the sheaf $\widetilde{\m}\otimes \F$ might not be complete, i.e. the map $\widetilde{\m}\otimes \F \to \lim_k \widetilde{\m}\otimes \F_k$ is a priori only an almost isomorphism. 
\end{rmk}

\begin{lemma}\label{reduction-almost-qcoh} Let $\X$ be a topologically finitely presented formal $R$-scheme for $R$ as in Setup~\ref{set-up3}, and let $\F^a$ be an almost adically quasi-coherent $\O_\X^a$-module. Then $\F_k^a$ is almost quasi-coherent for all $k$. Moreover, if an $\O_\X^a$-module $\G^a$ is annihilated by some $I^{n+1}$, then $\G^a$ is almost adically quasi-coherent if and only if so is $\G_n^a$.
\end{lemma}
\begin{proof}
To prove the first claim, it is sufficient to show that $\widetilde{\m}\otimes \F_k$ is quasi-coherent provided that $\widetilde{\m}\otimes \F$ is adically quasi-coherent. We use Corollary~\ref{pull-proj-cor} to say that $\widetilde{\m}\otimes \F_k \simeq (\widetilde{\m}\otimes \F)_k$ and the reduction of an adically quasi-coherent module is quasi-coherent. Therefore, each $\F_k^a$ is almost adically quasi-coherent.

Now if $\G$ is annihilated by $I^{n+1}$ then $\G=i_{n,*}\G_n$. We use the Projection Formula (Corollary~\ref{projection}) to see that $\widetilde{\m}\otimes \G \simeq i_{n,*}(\G_n\otimes \widetilde{\m})$. Clearly, $i_{n,*}$ sends quasi-coherent sheaves to adically quasi-coherent sheaves. So $\G^a$ is almost adically quasi-coherent if so is $\G_n^a$.
\end{proof}

\begin{defn} We say that an $\O_\X^a$-module $\F^a$ is {\it of almost finite type} (resp. {\it almost finitely presented}) if $\F^a$ is almost adically quasi-coherent, and there is a covering of $\X$ by open affines $\{\sU_i\}_{i\in I}$ such that $\F^a(\sU_i)$ is an almost finitely generated (resp. almost finitely presented) $\O_\X^a(\sU_i)$-module. We denote these categories by $\bf{Mod}_{\X^a}^{\text{aft}}$ and $\bf{Mod}_{\X^a}^{\text{afp}}$ respectively. \smallskip

We say that an $\O_\X$-module $\F$ is {\it of almost finite type} (resp. {\it almost finitely presented}) if so is $\F^a$. We denote these categories by $\bf{Mod}_{\X}^{\text{aft}}$ and $\bf{Mod}_{\X}^{\text{afp}}$ respectively. \smallskip
\end{defn}

\begin{defn}\label{defn:aqc-acoh} We say that an $\O_\X$-module $\F$ is {\it adically quasi-coherent of almost finite type} (resp. {\it adically quasi-coherent almost finitely presented}) if it is adically quasi-coherent and there is a covering of $\X$ by open affines $\{\sU_i\}_{i\in I}$ such that $\F(\sU_i)$ is an almost finitely generated (resp. almost finitely presented) $\O_\X(\sU_i)$-module. We denote these categories by $\bf{Mod}_\X^{\qc, \aft}$ and $\bf{Mod}_\X^{\qc, \afp}$ respectively.
\end{defn}

\begin{rmk} If $\F^a$ is a finite type (resp. finitely presented) $\O_\X^a$-module, then $(\F^a)_!$ is adically quasi-coherent of almost finite type (resp. almost finite presentation).
\end{rmk}

\begin{rmk} We note that, a priori, it is not clear if $\F^a$ is an almost finite type (resp. almost finitely presented) $\O_\X^a$-module for an adically quasi-coherent almost finite type (resp. almost finitely presented) $\O_\X$-module $\F$. The problem comes from the fact that we do not require $\widetilde{\m}\otimes \F$ to be adically quasi-coherent in Definition~\ref{defn:aqc-acoh}. However, we will show in Lemma~\ref{qc-aqc} that it is indeed automatic in this case.  
\end{rmk}

\begin{lemma}\label{qc-aqc} Let $\X$ be a topologically finitely presented formal $R$-scheme for $R$ as in Setup~\ref{set-up3}, and let $\F$ be an adically quasi-coherent $\O_\X$-module of almost finite type (resp. of almost finite presentation). Then $\widetilde{\m}\otimes \F$ is adically quasi-coherent. In particular, $\F$ is almost of finite type (resp. almost finitely presented). 
\end{lemma}
\begin{proof}
Corollary~\ref{almost-almost-finitely-presented} and Lemma~\ref{tensor-product} imply that the only condition we really need to check is that $\widetilde{\m}\otimes \F$ is adically quasi-coherent. Therefore, it suffices to prove the result for an adically quasi-coherent, almost finite type $\O_\X$-module $\F$.

Since the question is local on $\X$, we can assume that $\X=\Spf A$ is affine and $M\coloneqq \F(X)$ is almost finitely generated over $A$. Then we use \cite[Theorem I.3.2.8]{FujKato} to say that $\F$ is isomorphic to $M^{\Updelta}$. We claim that $\widetilde{\m}\otimes \F$ is isomorphic to $(\widetilde{\m}\otimes_A M)^{\Updelta}$ as that would imply that $\widetilde{\m}\otimes \F$ is adically quasi-coherent by \cite[Proposition I.3.2.2]{FujKato}. In order to show that $\widetilde{\m}\otimes \F$ is isomorphic to $(\widetilde{\m}\otimes_R M)^{\Updelta}$ we need to check two things: for any open affine $\Spf B=\sU \subset \X$ the $B$-module $(\widetilde{\m}\otimes \F)(\sU)$ is $I$-adically complete, and then the natural map $(\widetilde{\m}\otimes_R M)\widehat{\otimes}_A B \to (\widetilde{\m}\otimes \F)(\sU)$ is an isomorphism. \smallskip

We start with the first claim. Lemma~\ref{tensor-product} says that $(\widetilde{\m}\otimes \F)(\sU)$ is isomorphic to $\widetilde{\m}\otimes_R \F(\sU)$. Since $\F$ is adically quasi-coherent, $\F(\sU) \simeq M\wdh{\otimes}_A B$, so  $(\widetilde{\m}\otimes \F)(\sU) \simeq \widetilde{\m}\otimes_R (M\widehat{\otimes}_A B)$. Lemma~\ref{trivial-base-change} says that $M{\otimes}_A B$ is almost finitely generated over $B$, so it is already $I$-adically complete by Lemma~\ref{completion-finitely-generated}. Therefore, we see that $\widetilde{\m}\otimes_R \F(\sU) \simeq \widetilde{\m}\otimes_R (M{\otimes}_A B)$, and the latter is almost finitely generated over $B$ by Corollary~\ref{almost-almost-finitely-presented}. Thus, we use Lemma~\ref{completion-finitely-generated} once again to show its completeness. \smallskip

Now we show that the natural morphism $(\widetilde{\m}\otimes_R M)\widehat{\otimes}_A B \to (\widetilde{\m}\otimes \F)(\sU)$ is an isomorphism. Again, using the same results as above, we can get rid of any completions and identify this map with the ``identity'' map
\[
(\widetilde{\m}\otimes_R M){\otimes}_A B \to \widetilde{\m}\otimes_R (M{\otimes}_A B)
\]
This finishes the proof.
\end{proof}

\begin{lemma}\label{reduction-ft-formal} Let $\X$ be a topologically finitely presented formal $R$-scheme for $R$ as in Setup~\ref{set-up3}, and let $\F^a$ be an almost finite type (resp. almost finitely presented) $\O_\X^a$-module. Then the $\O_{\X_k}^a$-module $\F^a_k$ is almost finite type (resp. almost finitely presented) for any integer $k$.
\end{lemma}
\begin{proof}
Lemma~\ref{reduction-almost-qcoh} implies that each $\F^a_k$ is an almost quasi-coherent $\O_{X_k}$-module. So it is sufficient to find a covering of $\X_k$ by open affines $\sU_{i, k}$ such that $\F^a_k(\sU_{i,k})$ is almost finitely generated (resp. almost finitely presented) over $\O_{\X_k}^a(\sU_{i,k})$. We choose a covering of $\X$ by open affines $\sU_i$ such that $\F^a(\sU_{i})$ are almost finitely generated (resp. almost finitely presented) over $\O_\X^a(\sU_{i})$. Since the underlying topological spaces of $\X$ and $\X_k$ are the same, we conclude that $\sU_{i,k}$ form an affine open covering of $\X_k$. Then using the vanishing result for higher cohomology groups of adically quas-coherent sheaves on affine formal schemes of finite type \cite[Theorem I.7.1.1]{FujKato} and Lemma~\ref{tensor-product}, we deduce that 
\[
\F^a_k(\sU_{i,k}) \simeq (\widetilde{\m}\otimes \F^a_k)(\sU_{i,k})^a \simeq \left(\widetilde{\m}\otimes \F\left(\sU_{i}\right)/I^{k+1}\right)^a
\]
is almost finitely generated (resp. almost finitely presented) over $\O_{\X_k}(\sU_{i,k})$.
\end{proof}

\begin{lemma}\label{localize-ft-formal} Let $\X$ be a locally topologically finitely presented formal $R$-scheme for $R$ as in Setup~\ref{set-up3}, and let $\F^a$ be an almost finite type (resp. almost finitely presented) $\O_\X^a$-module. Then $\F^a(\sU)$ is an almost finitely generated (resp. almost finitely presented) $\O_\X^a(\sU)$-module for any open affine $\sU\subset \X$.
\end{lemma}
\begin{proof}
Corollary~\ref{almost-almost-finitely-presented} and Lemma~\ref{tensor-product} guarantee that we can replace $\F$ with $\widetilde{\m}\otimes \F$ for the purpose of the proof. Thus, we may and do assume that $\F$ is an adically quasi-coherent almost finitely generated (resp. almost finitely presented) $\O_\X$-module. Then, using Lemma~\ref{trivial-base-change} and Lemma~\ref{completion-finitely-generated}, we can use the argument as in the proof of Lemma~\ref{qc-aqc} to show that the restriction of $\F$ to any open formal subscheme is still adically quasi-coherent of almost finite type (resp. finitely presented). Thus, we may and do assume that $\X=\Spf A$ is an affine formal $R$-scheme. \smallskip

Now we have an affine topologically finitely presented formal $R$-scheme $\X=\Spf A$, a finite\footnote{We implicitly use here that $\X$ is quasi-compact.} covering of $\X$ by affines $\sU_i=\Spf A_i$, and an adically quasi-coherent $\O_\X$-module $\F$ such that $\F(\sU_i)$ is almost finitely generated (resp. almost finitely presented) $A_i$-module. We want to show that $\F(\X)$ is an almost finitely generated (resp. almost finitely presented) $A$-module. \smallskip

First, we deal with the {\it almost finitely generated case.} We note that Lemma~\ref{localize-ft}, Lemma~\ref{reduction-ft-formal}, and \cite[Theorem I.7.1.1]{FujKato} imply that $\F(\X)/I$ is almost finitely generated. We know that $\F$ is adically quasi-coherent so $\F(\X)$ must be an $I$-adically complete $A$-module. Therefore, Corollary~\ref{finite-mod-ideal} guarantees that $\F(X)$ is an almost finitely generated $A$-module. \smallskip

Now we move to the {\it almost finitely presented case}. We already know that $\F(\X)$ is almost finitely generated over $A$. Thus, the standard argument with Lemma~\ref{completion-finitely-generated} implies that $\F(\sU_i)=\F(\X)\otimes_{A} A_i$ for any $i$. Recall that \cite[Proposition I.4.8.1]{FujKato} implies\footnote{Topologically universally adhesive rings are by definition ``t. u. rigid-Noetherian''} that each $A \to A_i$ is flat. Since $\Spf A_i$ form a covering of $\Spf A$, we conclude that $A \to \prod_{i=1}^n A_i$ is faithfully flat. Now the result follows from faithfully flat descent for almost finitely presented modules, which is proven in Lemma~\ref{almost-flat-descent}.
\end{proof}

\begin{cor}\label{criterion-almost-finite-formal} Let $\X=\Spf A$ be a topologically finitely presented affine formal $R$-scheme for $R$ as in Setup~\ref{set-up3}, and let $\F^a$ be an almost adically quasi-coherent $\O_\X^a$-module. Then $\F^a$ is almost finite type (resp. almost finitely presented) if and only if $\F^a(\X)$ is almost finitely generated (resp. almost finitely presented) $A^a$-module. 

Similarly, an adically quasi-coherent $\O_\X$-module $\F$ is of almost finite type (resp. almost finitely presented) if and only if $\F(\X)$ is almost finitely generated (resp. almost finitely presented) $A$-module.
\end{cor}

\begin{lemma}\label{ker-coker-3} Let $\X=\Spf A$ be a topologically finitely presented affine formal $R$-scheme for $R$ as in Setup~\ref{set-up3}, let $\varphi\colon N \to M$ be a homomorphism of almost finitely generated $A$-modules. Then the following sequence
\[
0 \to (\ker \varphi)^\Updelta \to N^{\Updelta} \xr{\varphi^\Updelta} M^{\Updelta} \to (\coker \varphi)^{\Updelta} \to 0
\]
is exact. Moreover, $\Imm(\varphi)^\Updelta \simeq \Imm(\varphi^\Updelta)$. 
\end{lemma}
\begin{proof}
We denote the kernel $\ker \phi$ by $K$, the image $\rm{Im}(\varphi)$ by $M'$, and  the cokernel $\coker \phi$ by $Q$. \smallskip 

{\it We start with $\ker \varphi^\Updelta$}: We note that $(\ker \varphi^\Updelta) (\X)=K$, this induces a natural morphism $\a\colon K^{\Updelta} \to \ker \varphi^\Updelta$. In order to show that it is an isomorphism; it suffices to check that it induces an isomorphism on values over a basis of principal open subsets. Now recall that for any $A$-module $L$, we have an equality $L^{\Updelta}(\Spf A_{\{f\}})=\widehat{L_f}$ where the completion is taken with respect to the $I$-adic topology. Thus, in order to check that $\a$ is an isomorphism, it suffices to show that $\wdh{K_f}$ is naturally identified with $(\ker \varphi)(\Spf A_{\{f\}})=\ker (\wdh{N_f} \to \wdh{M_f})$. Using the Artin--Rees Lemma~\ref{Artin-Rees} over the adhesive ring $A_f$, we conclude that the induced topologies on $K_f$ and $M'_f$ coincide with the $I$-adic ones. This implies that
\[
\wdh{K_f}=\lim K_f/I^nK_f = \lim K_f/(I^nN_f \cap K_f) \text { and } \wdh{M'_f}=\lim M'_f/I^nM'_f = \lim M'_f/(I^nM_f \cap M'_f)
\]
This guarantees that we have two exact sequences:
\[
0 \to \wdh{K_f} \to \wdh{M_f} \to \wdh{M'_f} \to 0,
\]
\[
0 \to \wdh{M'_f} \to \wdh{N_f}
\]
In particular, we get that the natural map $\wdh{K_f}  \to \ker(\wdh{M_f} \to \wdh{N_f})$ is an isomorphism. This shows that $K^{\Updelta} \simeq \ker (\varphi^\Updelta)$. \smallskip

{\it We prove the claim for $\Imm \varphi^\Updelta$}: We note that since the category of $\O_\X$-modules is abelian, we can identify $\Imm \varphi^\Updelta \simeq \coker (K^{\Updelta} \to N^{\Updelta})$. We observe that \cite[Theorem I.7.1.1]{FujKato} and the established fact above that $\ker \varphi$ is adically quasi-coherent imply that the natural map $\F(\sU)/K^{\Updelta}(\sU) \to (\Imm \varphi^\Updelta)(\sU)$ is an isomorphism for any affine open formal subscheme $\sU$. In particular, we have $(\Imm \varphi^\Updelta)(\X)=M/K=M'$. Therefore, we have a natural map $(M')^{\Updelta} \to \Imm \varphi^\Updelta$ and we show that it is an isomorphism. It suffices to show that this map is an isomorphism on values over a basis of principal open subsets. Then we use the identification $\F(\sU)/K^{\Updelta}(\sU) \simeq (\Imm \varphi)(\sU)$ and the short exact sequence
\[
0 \to \wdh{K_f} \to \wdh{M_f} \to \wdh{M'_f} \to 0,
\]
to finish the proof. \smallskip

{\it We show the claim for $\coker \varphi^\Updelta$}: The argument is identical to the argument for $\Imm \varphi$ once we know that $\Imm \varphi=\ker (\G \to \coker \varphi)$ is adically quasi-coherent.
\end{proof}

\begin{cor}\label{ker-coker-4} Let $\X=\Spf A$ be a topologically finitely presented affine formal $R$-scheme for $R$ as in Setup~\ref{set-up3}, let $M$ an almost finitely generated $A$-module, and let $N$ be any $A$-submodule of $M$. Then the following sequence
\[
0 \to N^{\Updelta} \xr{\varphi^\Updelta} M^{\Updelta} \to (M/N)^{\Updelta} \to 0
\]
is exact. 
\end{cor}
\begin{proof}
We just apply Lemma~\ref{ker-coker-3} to the homomorphism $M \to M/N$ of almost finitely generated $A$-modules. 
\end{proof}

\begin{cor}\label{ker-coker} Let $\X$ be a topologically finitely presented formal $R$-scheme for $R$ as in Setup~\ref{set-up3}, and let $\varphi \colon \F\to \G$ be a morphism of adically quasi-coherent, almost finite type $\O_\X$-modules. Then $\ker \varphi$ is an adically quasi-coherent $\O_\X$-module, $\coker \varphi$ and $\Imm \varphi$ are adically quasi-coherent $\O_\X$-modules of almost finite type.
\end{cor}

\begin{cor}\label{ker-coker-2} Let $\X$ be a topologically finitely presented formal $R$-scheme for $R$ as in Setup~\ref{set-up3}, and let $\varphi \colon \F^a\to \G^a$ be a morphism of almost almost finite type $\O_\X^a$-modules. Then $\ker \varphi$ is an almost adically quasi-coherent $\O_\X^a$-module, $\coker \varphi$ and $\Imm \varphi$ are $\O_\X^a$-modules of almost finite type.
\end{cor}
\begin{proof}
We apply the exact functor $(-)_!$ to the map $\varphi$ and reduce the claim to Corollary~\ref{ker-coker}.
\end{proof}

Now we are ready to show that almost finite type and almost finitely presented $\O_\X$-modules share many good properties as we would expect. The only subtle thing is that we do not know whether an adically quasi-coherent quotient of an adically quasi-coherent, almost finite type $\O_\X$-module is of almost finite type. The main extra complication here is that we need to be very careful with the adically quasi-coherent condition in the definition of almost finite type (resp. almost finitely presented) modules since that condition does not behave well in general.

\begin{lemma}\label{main-almost-sheaf-formal} Let $0 \to \F' \xr{\varphi} \F \xr{\psi} \F'' \to 0$ be an exact sequence of $\O_\X$-modules, then
\begin{enumerate}[label=\textbf{(\arabic*)}]
\itemsep0.5ex
\item\label{main-almost-sheaf-formal-1} If $\F$ is adically quasi-coherent of almost finite type, and $\F'$ is adically quasi-coherent then $\F''$ is adically quasi-coherent of almost finite type.
\item\label{main-almost-sheaf-formal-2} If $\F'$ and $\F''$ are adically quasi-coherent of almost finite type (resp. almost finitely presented), then so is $\F$.
\item\label{main-almost-sheaf-formal-3} If $\F$ is adically quasi-coherent of almost finite type and $\F''$ is adically quasi-coherent almost finitely presented then $\F'$ is adically quasi-coherent of almost finite type.
\item\label{main-almost-sheaf-formal-4} If $\F$ is adically quasi-coherent of almost finitely presented and $\F'$ is adically quasi-coherent of almost finite type then $\F''$ is adically quasi-coherent, almost finitely presented.
\end{enumerate}
\end{lemma}
\begin{proof}
\ref{main-almost-sheaf-formal-1}: Without loss of generality, we can assume that $\X=\Spf A$ is an affine formal scheme. Then $\F\cong M^\Updelta$ for some almost finitely generated $A$-module $M$, and $\F' \cong N^\Updelta$ for some $A$-submodule $N \subset M$. Then Corollary~\ref{ker-coker-4} ensures that $\F'' \simeq (M/N)^\Updelta$. In particular, it is adically quasi-coherent. The claim is then an easy consequence of the vanishing theorem \cite[Theorem I.7.1.1]{FujKato} and Lemma~\ref{main-almost}\ref{main-almost-1}. \\

\ref{main-almost-sheaf-formal-2}: The difficult part is to show that $\F$ is adically quasi-coherent. In fact, once we know that $\F$ is adically quasi-coherent, it is automatically of almost finite type (resp. almost finitely presented) by \cite[Theorem I.7.1.1]{FujKato} and Lemma~\ref{main-almost}\ref{main-almost-2}. 

To show that $\F$ is adically quasi-coherent, we may and do assume that $\X=\Spf A$ is an affine formal $R$-scheme for some adhesive ring $A$. Then let us introduce $A$-modules $M'\coloneqq \F'(\X), \ M\coloneqq \F(\X)$, and $M''\coloneqq \F''(\X)$. We have the natural morphism $M^{\Updelta} \to \F$ and we show that it is an isomorphism. The vanishing theorem \cite[Theorem I.7.1.1]{FujKato} implies that we have a short exact sequence:
\[
0 \to M' \to M \to M'' \to 0
\]
Thus $M$ is almost finitely generated (resp. almost finitely presented) by Lemma~\ref{main-almost}\ref{main-almost-2}. Then Corollary~\ref{ker-coker-3} gives that we have a short exact sequence 
\[
0 \to M'^{\Updelta} \to M^{\Updelta} \to M''^{\Updelta} \to 0
\] 
Using the vanishing theorem \cite[Theorem I.7.1.1]{FujKato} once again we get a commutative diagram
\[
\begin{tikzcd}
0 \arrow{r} & M'^{\Updelta} \arrow{r} \arrow{d} &M^{\Updelta} \arrow{r} \arrow{d} & M''^{\Updelta} \arrow{r}\arrow{d} & 0 \\
0 \arrow{r} & \F' \arrow{r}  &\F \arrow{r}  & \F'' \arrow{r} & 0
\end{tikzcd}
\]
where the rows are exact, and the left and right vertical arrows are isomorphisms. That implies that the map $M^{\Updelta} \to \F$ is an isomorphism. \\

\ref{main-almost-sheaf-formal-3}: This easily follows from Lemma~\ref{main-almost}\ref{main-almost-3}, Lemma~\ref{ker-coker}, and \cite[Theorem I.7.1.1]{FujKato}. \\

\ref{main-almost-sheaf-formal-4}: This also easily follows from Lemma~\ref{main-almost}\ref{main-almost-4}, Lemma~\ref{ker-coker}, and \cite[Theorem I.7.1.1]{FujKato}.
\end{proof}

We also give the almost version of this lemma:

\begin{cor}\label{cor:main-almost-sheaf-formal-2} Let $0 \to \F'^a \xr{\varphi} \F^a \xr{\psi} \F''^a \to 0$ be an exact sequence of $\O_\X^a$-modules, then
\begin{enumerate}[label=\textbf{(\arabic*)}]
\itemsep0.5ex
\item If $\F^a$ is of almost finite type, and $\F'^a$ is almost adically quasi-coherent then $\F''^a$ is of almost finite type;
\item If $\F'^a$ and $\F''^a$ are of almost finite type (resp. almost finitely presented), then so is $\F^a$;
\item If $\F^a$ is  of almost finite type and $\F''^a$ is almost finitely presented, then $\F'^a$ is of almost finite type;
\item If $\F^a$ is of almost finitely presented and $\F'^a$ is  of almost finite type then $\F''^a$ is  almost finitely presented.
\end{enumerate}
\end{cor}

\begin{defn} We say that an $\O_\X^a$-module $\F^a$ is {\it almost coherent} if $\F^a$ is almost finite type and for any open set $\sU$ any finite type $\O_\X^a$-submodule $\G^a \subset (\F|_\sU)^a$ is an almost finitely presented $\O_\sU$-module. \\

We say that $\F$ is {\it (adically quasi-coherent) almost coherent} $\O_\X$-module if $\F^a$ almost coherent (and $\F$ is adically quasi-coherent).
\end{defn}

\begin{rmk} Lemma~\ref{qc-aqc} ensures that any adically quasi-coherent almost coherent $\O_X$-module $\F$ is almost coherent. 
\end{rmk}

\begin{lemma}\label{criterion-almost-coherent-formal} Let $\F^a$ be an $\O_\X^a$-module on a topologically finitely presented formal $R$-scheme $\X$. Then the following are equivalent:
\begin{enumerate}[label=\textbf{(\arabic*)}]
\item $\F^a$ is almost coherent;
\item $\F^a$ is almost quasi-coherent and the $\O_\X^a(\sU)$-module $\F^a(\sU)$ is almost coherent for any open affine formal subscheme $\sU\subset \X$; 
\item $\F^a$ is almost quasi-coherent and there is a covering of $\X$ by open affine subschemes $(\sU_i)_{i\in I}$ such that $\F^a(\sU_i)$ is almost coherent for each $i$.
\end{enumerate} 
In particular, an $\O_\X^a$-module $\F^a$ is almost coherent if and only if it is almost finitely presented.
\end{lemma}
\begin{proof}
The proof that these three notions are equivalent is identical to the proof of Lemma~\ref{criterion-almost-coherent-formal} modulo facts that we have already established in this chapter, especially Corollary~\ref{ker-coker-3}. 

As for the last claim, we recall that $\X$ is topologically finitely presented over a topologically universally adhesive ring, so $\O_\X(\sU)$ is coherent for any open affine $\sU$ \cite[Prop. 0.8.5.23, Lemma I.1.7.4, Prop. I.2.3.3]{FujKato}. Then Lemma~\ref{coh-acoh-2} and Lemma~\ref{acoh-fp-coh-2} prove the equivalence.
\end{proof}

Although Lemma~\ref{criterion-almost-coherent-formal} says that the notion of almost coherence coincides with the notion of almost finite presentation, it shows that almost coherence is morally ``the correct'' definition. In what follows, we prefer to use the terminology of almost coherent sheaves as it is shorter and gives a better intuition from our point of view. 

\begin{lemma}\label{main-coh-glob-formal}  
\begin{enumerate}[label=\textbf{(\arabic*)}]
\item Any almost finite type $\O_\X^a$-submodule of an almost coherent $\O_\X^a$-module is almost coherent;
\item Let $\varphi\colon \F^a \to \G^a$ be a homomorphism from an almost finite type $\O_\X^a$-module to an almost coherent $\O_\X^a$-module. Then $\ker \varphi$ is an almost finite type $\O_\X^a$-module;
\item\label{main-coh-formal-4} Let $\varphi\colon \F^a \to \G^a$ be a homomorphism of almost coherent $\O_\X^a$-modules. Then $\ker \varphi$ and $\coker \varphi$ are almost coherent $\O_\X$-modules;
\item\label{main-coh-formal-5} Given a short exact sequence of $\O_\X^a$-modules
\[0 \to \F'^a \to \F^a \to \F''^a \to 0,\]
 if two out of three are almost coherent, then so is the third one.
\end{enumerate}
\begin{rmk}
    There is also an evident version of this corollary for adically quasi-coherent almost coherent $\O_\X$-modules.
\end{rmk}
\end{lemma}
\begin{proof}
The proof is identical to the proof of Corollary~\ref{main-coh-glob} once we have Lemma~\ref{ker-coker} and the equivalence of almost coherent and almost finitely presented $\O_\X$-modules from Lemma~\ref{criterion-almost-coherent-formal}.
\end{proof}

\begin{cor} Let $\X$ be topologically finitely presented formal $R$-scheme for $R$ as in Setup~\ref{set-up3}. Then the category $\bf{Mod}_\X^{\text{acoh}}$ (resp. $\bf{Mod}_\X^{\qc, \text{acoh}}$, $\bf{Mod}_{\X^a}^{\text{acoh}}$) of almost coherent $\O_\X$-modules (resp. adically quasi-coherent, almost coherent $\O_\X$-modules, resp. almost coherent $\O_\X^a$-modules) is a weak Serre subcategory of $\bf{Mod}_\X$ (resp. $\bf{Mod}_\X$, resp. $\bf{Mod}_\X^a$). 
\end{cor}

\subsection{Formal Schemes. Basic functors on almost coherent $\O_\X^a$-modules}\label{basic-formal}

In this section, we study the interaction between the functors defined in Section~\ref{basic-functors-sheaves} and the notion of almost (quasi-)coherent $\O_\X^a$-modules. The exposition follows Section~\ref{basic-functors-schemes} very closely. \smallskip

We start with an affine situation, i.e. $\X=\Spf A$. In this case, we note that the functor $(-)^\Updelta\colon \bf{Mod}_A \to \bf{Mod}^{\qc}_\X$ sends almost zero $A$-modules to almost zero $\O_\X$-modules. Thus, it induces a functor
\[
(-)^\Updelta \colon \bf{Mod}_{A^a} \to \bf{Mod}_{\X^a}.
\]

\begin{lemma}\label{sheafify-modules-formal} Let $\X=\Spf A$ be an affine formal $R$-scheme for $R$ as in Setup~\ref{set-up3}. Then the functor $(-)^\Updelta \colon \bf{Mod}_A \to \bf{Mod}^{\qc}_\X$ induces an equivalence $(-)^\Updelta \colon \bf{Mod}_A^* \to \bf{Mod}_\X^{\qc, *}$ for any $*\in \{\text{aft}, \text{acoh}\}$. The quasi-inverse functor is given by $\Gamma(\X, -)$. 
\end{lemma}
\begin{proof}
We note that the functor $(-)^\Updelta \colon \bf{Mod}_A \to \bf{Mod}^{\qc}_\X$ induces an equivalence between the category of $I$-adically complete $A$-modules and adically quasi-coherent $\O_\X$-modules by \cite[Theorem I.3.2.8]{FujKato}. Recall that almost finite type modules are complete due to Lemma~\ref{completion-finitely-generated}. Thus, it suffices to show that an adically quasi-coherent $\O_\X$-module is almost finitely generated (resp. almost coherent) if and only if so is $\Gamma(X, \F)$. Now this follows from Lemma~\ref{criterion-almost-finite-formal} and Lemma~\ref{criterion-almost-coherent-formal}.  
\end{proof}

\begin{lemma}\label{sheafify-modules-almost-formal} Let $\X=\Spf A$ be an affine formal $R$-scheme for $R$ as in Setup~\ref{set-up3}. Then the functor $(-)^\Updelta \colon \bf{Mod}_A \to \bf{Mod}^{\qc}_\X$ induces equivalences $(-)^\Updelta \colon \bf{Mod}_{A^a}^{*} \to \bf{Mod}_{\X^a}^{*}$ for any $*\in \{\text{aft}, \text{acoh}\}$. The quasi-inverse functor is given by $\Gamma(X, -)$. 
\end{lemma}
\begin{proof}
The proof is analogous to Lemma~\ref{sheafify-modules-almost} once Lemma~\ref{sheafify-modules-formal} is verified. 
\end{proof}

Now we recall that, for any $R$-scheme $X$, we can define the $I$-adic completion of $X$ as a colimit $\colim (X_k, \O_{X_k})$ of the reductions $X_k\coloneqq X\times_R \Spec R/I^{k+1}$ in the category of formal schemes. We refer to \cite[\textsection 1.4(c)]{FujKato} for more details. This completion comes with a map of locally ringed spaces
\[
c\colon \wdh{X} \to X.
\]
In the affine case, we note that $\wdh{\Spec A}=\Spf \wdh{A}$ for any $R$-algebra $A$\footnote{We note that $\wdh{A}$ is $I$-adically complete due to \cite[\href{https://stacks.math.columbia.edu/tag/05GG}{Tag 05GG}]{stacks-project}.}. We study properties of the completion map for (topologically) finitely presented $R$-algebra $A$.

\begin{lemma}\label{completion-affine} Let $X=\Spec A$ be an affine $R$-scheme for $R$ as in Setup~\ref{set-up3}. Suppose that $A$ is either finitely presented or topologically finitely presented over $R$. Then the morphism $c\colon \wdh{X }\to X$ is flat, and there is a functorial isomorphism $M^\Updelta \cong c^*(\widetilde{M})$ for any almost finitely generated $A$-module $M$.
\end{lemma}
\begin{proof}
The flatness assertion is proven in \cite[Proposition I.1.4.7 (2)]{FujKato}. Now the natural map
\[
M\to \rm{H}^0(\X, c^*(\widetilde{M}))
\] 
induces the map $M^{\Updelta} \to c^*(\widetilde{M})$. To show that it is an isomorphism, it suffices to show that the map
\[
\wdh{M_f} \to M_f \otimes_{A_f} \wdh{A_{f}}
\]
is an isomorphism for any $f\in A$. This follows from Lemma~\ref{completion-finitely-generated}, as each such $A_f$ is $I$-adically adhesive. 
\end{proof}

\begin{cor}\label{completion-global} Let $X$ be a locally finitely presented $R$-scheme for $R$ as in Setup~\ref{set-up3}. Then the morphism $c\colon \wdh{X}\to X$ is flat and $c^*$ sends almost finite type $\O_X^a$-modules (resp. almost coherent $\O_X^a$-modules) to almost finite type $\O_\X^a$-modules (resp. almost coherent $\O_\X^a$-modules). \smallskip

Similarly, $c^*$ sends quasi-coherent almost finite type $\O_X$-modules (resp. quasi-coherent almost coherent $\O_X$-modules) to adically quasi-coherent almost finite type $\O_\X$-modules (resp. adically quasi-coherent almost coherent $\O_\X$-modules)
\end{cor}
\begin{proof}
The statement is local, so we can assume that $X=\Spec A$. Then the claim follows from Lemma~\ref{completion-affine}.
\end{proof}

Now we show that the pullback functor preserves almost finite type and almost coherent $\O_\X^a$-modules.

\begin{lemma}\label{pullback-almost-coh-formal} Let $\mf\colon \X \to \Y$ be a morphism of locally finitely presented formal $R$-scheme for $R$ as in Setup~\ref{set-up3}.  
\begin{enumerate}[label=\textbf{(\arabic*)}]
\item\label{pullback-almost-coh-formal-1} Suppose that $\X=\Spf B$, $\Y=\Spf A$ are affine formal $R$-schemes. Then $\mf^*(M^{\Updelta})$ is functorially isomorphic to $(M\otimes_{A} B)^\Updelta$ for any $M\in \bf{Mod}_A^{\text{aft}}$.
\item\label{pullback-almost-coh-formal-1.5} Suppose that $\X=\Spf B$, $\Y=\Spf A$ are affine formal $R$-schemes. Then $\mf^*(M^{a, \Updelta})$ is functorially isomorphic to $(M^a\otimes_{A^a} B^a)^\Updelta$ for any $M^a\in \bf{Mod}_A^{a, \text{aft}}$.
\item\label{pullback-almost-coh-formal-2} The functor $\mf^*$ sends $\bf{Mod}_\Y^{\qc, \text{aft}}$ (resp. $\bf{Mod}_\Y^{\qc, \text{acoh}}$) to $\bf{Mod}_\X^{\text{aft}}$ (resp. $\bf{Mod}_\X^{\qc, \text{acoh}}$).
\item\label{pullback-almost-coh-formal-3} The functor $\mf^*$ sends $\bf{Mod}_{\Y^a}^{\text{aft}}$ (resp. $\bf{Mod}_{\Y^a}^{\text{acoh}}$) to $\bf{Mod}_{\X^a}^{\text{aft}}$ (resp. $\bf{Mod}_{\X^a}^{\text{acoh}}$).
\end{enumerate}
\end{lemma}
\begin{proof}
We prove \ref{pullback-almost-coh-formal-1}, the other parts follow from this (as in the proof of Lemma~\ref{pullback-almost-coh}). \smallskip

We consider a commutative diagram
\[
\begin{tikzcd}
\Spf B \arrow{r}{c_B} \arrow{d}{\mf} & \Spec B \arrow{d}{f} \\
\Spf A \arrow{r}{c_A} & \Spec Am
\end{tikzcd}
\]
where the map $f\colon \Spec B \to \Spec A$ is the map induced by $\mf^\#\colon A \to B$. Then we have that $M^\Updelta \simeq c_A^*\widetilde{M}$ by Lemma~\ref{completion-affine}. Therefore,
\[
\mf^*(M^\Updelta) \simeq c_B^*(f^* \widetilde{M}) \simeq c_B^*(\widetilde{M \otimes_A B}) \simeq (M \otimes_A B)^\Updelta
\]
where the last isomorphism follows from Lemma~\ref{completion-affine}.
\end{proof}

The next thing we discuss is the intersection of tensor products and almost coherent sheaves.

\begin{lemma}\label{tensor-almost-coh-formal} Let $\X$ be a topologically finitely presented formal $R$-scheme for $R$ as in Setup~\ref{set-up3}.  
\begin{enumerate}[label=\textbf{(\arabic*)}]
\item\label{tensor-almost-coh-formal-1} Suppose that $\X=\Spf A$ is affine. Then $M^{\Updelta}\otimes_{\O_\X}N^\Updelta$ is functorially isomorphic to $(M\otimes_{A} N)^\Updelta$ for any $M, N\in \bf{Mod}_A^{\text{aft}}$. 
\item\label{tensor-almost-coh-formal-2} Suppose that $\X=\Spf A$ is affine. Then $M^{a, \Updelta}\otimes_{\O_\X^a}N^{a, \Updelta}$ is functorially isomorphic to $(M^a\otimes_{A^a} N^a)^{\Updelta}$ for any $M^a, N^a\in \bf{Mod}_{A^a}^{ \text{aft}}$. 
\item\label{tensor-almost-coh-formal-3} Let $\F, \G$ be two adically quasi-coherent almost finite type (resp. almost finitely presented) $\O_\X$-modules. Then the $\O_X$-module $\F\otimes_{\O_\X}\G$ is adically quasi-coherent of almost finite type (resp. almost finitely presented).   
\item\label{tensor-almost-coh-formal-4} Let $\F^a, \G^a$ be two almost finite type (resp. almost coherent) $\O_\X^a$-modules. Then the $\O_\X^a$-module $\F^a\otimes_{\O_\X^a}\G^a$ is of almost finite type (resp. almost coherent). The analogous result holds for $\O_X$-modules $\F, \G$.  
\end{enumerate}
\end{lemma}
\begin{proof}
Again, we only show \ref{tensor-almost-coh-formal-1} as the other parts follow from this similarly to the proof of Lemma~\ref{tensor-almost-coh}. \smallskip

The proof of \ref{tensor-almost-coh-formal-1} is, in turn, similar to that of Lemma~\ref{pullback-almost-coh-formal}~\ref{pullback-almost-coh-formal-1}. We consider the completion morphism $c\colon \Spf A \to \Spec A$. Then we have a sequence of isomorphisms
\[
M^{\Updelta}\otimes_{\O_\X}N^\Updelta \simeq c^*(\widetilde{M}) \otimes_{\O_\X} c^*(\widetilde{N}) \simeq c^*(\widetilde{M} \otimes_{\O_{\Spec A}} \widetilde{N}) \simeq c^*(\widetilde{M\otimes_A N}) \simeq (M\otimes_A N)^{\Updelta}. \qedhere
\]
\end{proof}

Finally, we deal with the functor $\ud{\cal{H}om}_{\O_X^a}(-, -)$. We start with the following preparatory lemma:

\begin{lemma}\label{hom-alcoh-formal} Let $\X$ be a locally topologically finitely presented formal $R$-scheme for $R$ as in Setup~\ref{set-up3}.  
\begin{enumerate}[label=\textbf{(\arabic*)}]
\item\label{hom-alcoh-formal-1} Suppose $\X=\Spf A$ is affine. Then the canonical map 
\begin{equation}\label{1-1-1-formal}
\rm{Hom}_A(M,N)^{\Updelta} \to \ud{\cal{H}om}_{\O_\X}(M^\Updelta, N^\Updelta)
\end{equation}
is an almost isomorphism for any almost coherent $A$-modules $M$ and $N$. 
\item\label{hom-alcoh-formal-2} Suppose $\X=\Spf A$ is affine. Then there is a functorial isomorphism 
\begin{equation}\label{1-1-3-formal}
\rm{alHom}_{A^a}(M^a,N^a)^{\Updelta} \simeq \ud{al\cal{H}om}_{\O_\X^a}(M^{a,\Updelta}, N^{a, \Updelta})
\end{equation}
for any almost coherent $A^a$-modules $M^a$ and $N^a$. We also get a functorial almost isomorphism
\begin{equation}\label{1-1-2}
\rm{Hom}_{A^a}(M^a,N^a)^{\Updelta} \simeq^a \ud{\cal{H}om}_{\O_\X^a}(M^{a,\Updelta}, N^{a, \Updelta})
\end{equation}
for any almost coherent $A^a$-modules $M^a$ and $N^a$.
\item\label{hom-alcoh-formal-3} Suppose $\F$ and $\G$ are almost coherent $\O_\X$-modules. Then $\ud{\cal{H}om}_{\O_\X}(\F, \G)$ is an almost coherent $\O_\X$-module.
\item\label{hom-alcoh-formal-4} Suppose $\F^a$ and $\G^a$ are almost coherent $\O_\X^a$-modules. Then 
\[
\ud{\cal{H}om}_{\O_\X^a}(\F^a, \G^a) \text{ (resp. } \ud{al\cal{H}om}_{\O_\X^a}(\F^a, \G^a)\text{)}
\]
is an almost coherent $\O_\X$-module (resp. $\O_X^a$-module).
\end{enumerate}
\end{lemma}
\begin{proof}
Again, the proof is analogous to that of Lemma~\ref{hom-alcoh} and Corollary~\ref{alhom-acoh} once \ref{hom-alcoh-formal-1} is proven. So we only give a proof of \ref{hom-alcoh-formal-1} here. \smallskip

We note that both $M$ and $N$ are $I$-adically complete by Lemma~\ref{completion-finitely-generated}. Now we use \cite{FujKato} to say that the natural map $\rm{Hom}_A(M,N) \to \rm{Hom}_{\O_\X}(M^\Updelta, N^\Updelta)$ is an isomorphism. This induces a morphism 
\[
\rm{Hom}_A(M,N)^\Updelta \to \ud{\cal{H}om}_{\O_\X}(M^\Updelta, N^\Updelta).
\]
In order to prove that it is an almost isomorphism, it suffices to show that the natural map
\[
\rm{Hom}_A(M,N) \wdh{\otimes}_A A_{\{f\}} \to \rm{Hom}_{A_{\{f\}}}(M\wdh{\otimes}_A A_{\{f\}}, N\wdh{\otimes}_A A_{\{f\}})
\]
is an almost isomorphism for any $f\in A$. Now we note that $\rm{Hom}_A(M, N)$ is almost coherent by Corollary~\ref{tensor-hom-coh}. Thus, $\rm{Hom}_A(M,N) \otimes_A A_{\{f\}}$ is already complete, so the completed tensor product coincides with the usual one. Similarly, $M\wdh{\otimes}_A A_{\{f\}}\simeq M \otimes_A A_{\{f\}}$ and $N\wdh{\otimes}_A A_{\{f\}}\simeq N \otimes_A A_{\{f\}}$. Therefore, the question boils down to showing that the natural map
\[
\rm{Hom}_A(M,N) \otimes_A A_{\{f\}} \to \rm{Hom}_{A_{\{f\}}}(M\otimes_A A_{\{f\}}, N\otimes_A A_{\{f\}})
\]
is an almost isomorphism. This, in turn, follows from Lemma~\ref{flat-base-change-hom}.
\end{proof}

\subsection{Formal Schemes. Approximation of almost coherent $\O_\X^a$-modules}

In this section, we fix a ring $R$ as in the Set-up~\ref{set-up3}, and $\X$ a topologically finitely presented formal $R$-scheme. \smallskip

The main goal of this section is to establish an analogue of Corollary~\ref{cor:approximate-afpr} in the context of formal schemes. More precisely, we show that, for any finitely generated ideal $\m_0\subset \m$, an almost coherent $\O_\X$-module $\F$ can be ``approximated'' by a coherent $\O_\X$-module $\G_{\m_0}$ up to $\m_0\subset \m$ torsion. It turns out that this result is more subtle than its algebraic counterpart because, in general, we do not know if we can present an adically quasi-coherent $\O_\X$-module as a filtered colimit of finitely presented $\O_\X$-modules. Also, colimits are much more subtle in the formal set-up due to the presence of topology. This seems unlikely that the method used in the proof Corollary~\ref{cor:approximate-afpr} can be used in the formal set-up. Instead, we take another route and, instead, we first approximate $\F$ up to bounded torsion and then reduce to the algebraic case. \smallskip

\begin{defn} A map of $\O_\X$-modules $\phi\colon \G \to \F$ is an {\it FP-approximation} if  $\G$ is a finitely presented $\O_\X$-module, and $I^n(\Ker \phi)=0$, $I^n(\coker \phi)=0$ for some $n>0$. 

If $\m_0 \subset \m$ is a finitely generated sub-ideal of $\m_0$, a map of $\O_\X$-modules $\phi\colon \G \to \F$ is an {\it FP-$\m_0$-approximation} if it is an FP-approximation and $\m_0(\coker \phi)=0$.
\end{defn}

\begin{lemma}\label{lemma:FP-approximation-affine} Let $\X=\Spf A$ be an affine topologically finitely presented formal $R$-scheme, and $\F$ an adically quasi-coherent $\O_\X$-module of almost finite type. Then, for any finitely generated ideal $\m_0\subset \m$, $\F$ admits an FP-$\m_0$-approximation. 
\end{lemma}
\begin{proof}
    Lemma~\ref{sheafify-modules-almost-formal} guarantees that $\F=M^\Updelta$ for some almost finitely generated $A$-module $M$. Then, by definition, there is a submodule $N\subset M$ such that $\m_0(M/N)$. By assumption, $U\coloneqq \Spec A \setminus \rm{V}(I)$ is noetherian, so $\widetilde{N}|_U$ is a finitely presented $\O_U$-module. Then \cite[Lemma 0.8.1.6(2)]{FujKato} guarantees that there is a finitely presented $A$-module $N'$ with a surjective map $N' \to N$ such that its kernel $K$ is $I^\infty$-torsion. In particular, $K \subset N'[I^\infty]$. But since $A$ is $I$-adically complete and noetherian outside $I$, \cite[Theorem 5.1.2 and Definition 4.3.1]{FGK} guarantee that $N'[I^\infty]=N'[I^n]$ for some $n\geq 0$. In particular, $K$ is an $I^n$-torsion module. \smallskip
    
    Therefore, we have an exact sequence
    \[
    0 \to K \to N' \to M \to Q \to 0,
    \]
    where $N'$ is finitely presented, $M$ is almost finitely generated,  $\m_0Q=0$, and $I^nK=0$ for some $n\geq 1$. Now Lemma~\ref{ker-coker-3} says that the following sequence is exact:
    \[
    0 \to K^\Updelta \to N'^\Updelta \to M^\Updelta \to Q^\Updelta \to 0.
    \]
    In particular, $N'^\Updelta$ is a finitely presented $\O_\X$-module, $\m_0(Q^\Updelta)=0$, and $I^n(K^\Updelta)$.    
\end{proof}

\begin{lemma}\label{lemma:approximate-finite-type}\cite[Exercise I.3.4]{FujKato} Let $\X$ be a finitely presented formal $R$-scheme, $\F$ an adically quasi-coherent $\O_\X$-module of finite type, and $\G \subset \F$ an adically quasi-coherent $\O_\X$-submodule. Then $\G$ is a filtered colimit $\G = \colim_{\lambda \in \Lambda} \G_{\lambda}$ of adically quasi-coherent $\O_\X$-submodules of finite type such that, for all $\lambda\in \Lambda$,  $\G/\G_\lambda$ is annihilated by $I^n$ for
a fixed $n>0$.
\end{lemma}

\begin{lemma}\label{lemma:FP-approximations-filtered} Let $\X$ be a finitely presented formal $R$-scheme, $\F$ an adically quasi-coherent, almost finitely generated $\O_\X$-module, and $\phi_i \colon \G_i \to \F$ for $i=1, 2$ two FP-$\m_0$-approximations of $\F$ for some finitely generated sub-ideal $\m_0\subset \m$. Then there is a commutative diagram
\[
\begin{tikzcd}
\G_1 \arrow[d, swap, "q_1"] \arrow{dr}{\phi_1} &  \\
\cal{H} \arrow{r}{\phi} & \F \\
\G_2 \arrow{u}{q_2} \arrow[ur, swap, "\phi_2"]& 
\end{tikzcd}
\]
where $\phi$ and $q_i$ are FP-$\m_0$-approximations for $i=1, 2$. 
\end{lemma}
\begin{proof}
    By assumption, there is an integer $c>0$ such that $\ker(\phi_i)$ and $\coker(\phi_i)$ are annihilated by $I^c$ for $i=0,1$. Therefore, we may replace $\m_0$ by $\m_0+I^c$ to assume that $\m_0$ contains $I^c$. \smallskip
    
    Now we define $\cal{K}$ to be the kernel of the natural morphism $\G_1\oplus \G_2 \to \F$. Note that it is an adically quasi-coherent $\O_\X$-submodule of $\G_1\oplus \G_2$ due to Lemma~\ref{ker-coker-3}. Therefore, Lemma~\ref{lemma:approximate-finite-type} applies to the inclusion $\cal{K} \subset \G_1\oplus \G_2$, so we can write $\cal{K}=\colim_{\lambda\in \Lambda} \cal{K}_\lambda$ as a filtered colimit of adically quasi-coherent, finite type $\O_\X$-submodules of $\G_1\oplus \G_2$ with $I^{m}(\cal{K}/\cal{K}_\lambda)=0$ for some fixed $m>0$ and every $\lambda\in \Lambda$. We define $\cal{H}_\lambda=(\G_1\oplus \G_2)/\cal{K}_\lambda$, it comes with the natural morphisms
    \[
    \phi_\lambda \colon \cal{H}_{\lambda} \to \F, 
    \]
    \[
    q_{i, \lambda} \colon \G_i \to \cal{H}_\lambda 
    \]
    for $i=1,2$. We claim that these morphisms satisfy the claim of the lemma for some $\lambda \in \Lambda$, i.e. $\phi_{\lambda}$, and $q_{i, \lambda}$ are FP-$\m_0$-approximations. \smallskip
    
    Since $\X$ is topologically finitely presented (in particular, it is quasi-compact and quasi-separated), these claims can be checked locally. So we may and do assume that $\X=\Spf A$ is affine. Then we use Lemma~\ref{sheafify-modules-almost-formal}, \cite[Theorem I.3.2.8, Proposition I.3.5.4]{FujKato} to reduce to the situation where $\X=\Spf A$, $\F=M^\Updelta$, $\G_1=N_1^\Updelta$, $\G_2=N_2^\Updelta$ for some almost finitely generated $A$-module $M$, and finitely presented $A$-modules $N_1$, $N_2$ with maps of sheaves induced by homomorphisms $N_1 \to M$ and $N_2 \to M$. Then Lemma~\ref{ker-coker-3} guarantees that $\cal{K}=K^\Updelta$ for $K=\ker (N_1\oplus N_2 \to M)$, and $K=\colim_{\lambda\in \Lambda} K_\lambda$ for finitely generated $A$-submodules\footnote{Here, $K_\lambda=\Gamma(\X, \cal{K}_\lambda)$, so the equality follows from \cite[\href{https://stacks.math.columbia.edu/tag/009F}{Tag 009F}]{stacks-project}.} $K_\lambda$ with $I^m(K/K_\lambda)=0$ for some fixed $m>0$ and all $\lambda\in \Lambda$. So one can use Lemma~\ref{ker-coker-3} once again to conclude that it suffices (due to the assumption that $I^c \subset \m_0$) to show that, for some $\lambda\in \Lambda$, the natural morphisms $(N_1\oplus N_2)/K_\lambda \to M$, $N_i \to (N_1\oplus N_2)/K_\lambda$ have kernels annihilated by some power of $I$, and cokernels annihilated by $\m_0$. \smallskip
    
    The kernels of $N_i \to (N_1\oplus N_2)/K_\lambda$ (for $i=1, 2$) embed into the respective kernels for the natural morphisms $N_i \to M$, so they are automatically annihilated by some power of $I$ for any $\lambda\in \Lambda$. Also, clearly, the morphism $(N_1\oplus N_2)/K_\lambda \to M$ has kernel $K/K_{\lambda}$ that is annihilated by $I^m$ by the choice of $K_{\lambda}$. \smallskip
    
    Therefore, it suffices to show that we can choose $\lambda\in \Lambda$ such that $q_{i, \lambda}\colon N_i\to (N_1\oplus N_2)/K_\lambda$ (for $i=1,2$) and $\phi_\lambda \colon (N_1\oplus N_2)/K_\lambda \to M$ have cokernels annihilated by $\m_0$. The latter case is automatic and actually holds for any $\lambda\in \Lambda$. So the only non-trivial thing we need to check is that $\m_0(\coker q_{i, \lambda})=0$ for some $\lambda \in \Lambda$. \smallskip
    
    Let $(m_1, \dots, m_d) \in \m_0$ be a finite set of generators, and $\{y_{i, j}\}_{j\in J_i}$ a finite set of generators of $N_i$ for $i=1, 2$. Denote by $\ov{y_{i,j}}$ the image of $y_{i, j}$ in $M$. Define $x_{i, j, k}\in N_{2-i}$ to be a lift of $m_k \ov{y_{i, j}}\in M$ in $N_{2-i}$ for $k=1, \dots, d$, $i=1, 2$ and $j\in J_i$. Note that elements $(m_k y_{1, j}, x_{1,j,k})\in N_1\oplus N_2$ and $(x_{2, j, k}, m_k y_{2, j}) \in N_1\oplus N_2$ lie in $K$. So for some $\lambda\in \Lambda$, $K_\lambda$ contains the elements $(m_k y_{1, j}, x_{1, j, k})$ and $(x_{2, j, k}, m_k y_{2, j})$. Then it is easy to see that the cokernels of  $N_i \to  (N_1\oplus N_2)/K_\lambda$ are annihilated by $\m_0$. This finishes the proof.
\end{proof}

\begin{lemma}\label{lemma:global-FP-approximation} Let $\X$ be a finitely presented formal $R$-scheme, $\F$ an adically quasi-coherent, almost finitely generated $\O_\X$-module. Then, for any finitely generated ideal $\m_0\subset \m$, $\F$ is FP-$\m_0$-approximated.
\end{lemma}
\begin{proof}
   First, we note that Lemma~\ref{lemma:FP-approximation-affine} guarantees that the claim holds if $\X$ is affine. Now choose a covering of $\X$ by open affines $\X=\bigcup_{i=1}^n \V_i$. We know that claim on each $\V_i$, so it suffices to show that, if $\X = \sU_1 \cup \sU_2$ is a union of two finitely presented open formal subschemes and $\F$ is FP-$\m_0$-approximated on $\sU_1$ and $\sU_2$, then $\F$ is FP-$\m_0$-approximated on $\X$. \smallskip
    
    Suppose that $\G_i \to \F|_{\sU_i}$ are FP-$\m_0$-approximations on $\sU_i$ for $i=1, 2$. Then the intersection $\sU_{1,2}\coloneqq \sU_1 \cap \sU_2$ is again a topologically finitely presented formal $R$-scheme because $\X$ is so. Therefore, Lemma~\ref{lemma:FP-approximations-filtered} guarantees that we can find another FP-$\m_0$-approximation $\cal{H} \to \F|_{\sU_{1,2}}$ that is dominated by both $\G_i|_{\sU_{1,2}} \to \F|_{\sU_{1,2}}$ for $i=1,2$. Consider the $\O_{\sU_{1,2}}$-modules 
    \[
    \cal{K}_i \coloneqq \ker(\G_i|_{\sU_{1,2}} \to \cal{H}) \text{ for $i=1, 2$}.
    \]
    Lemma~\ref{ker-coker-3} guarantees that both $\cal{K}_i$ are adically quasi-coherent $\O_\X$-modules of finite type\footnote{Since they are kernels of morphisms between coherent $\O_\X$-modules}. The fact that $\G_i|_{\sU_{1,2}} \to \cal{H}$ are FP-$\m_0$-approximations ensures that both $\cal{K}_i$ are killed by some $I^m$ for $m\geq 1$. In particular, we see that $\cal{K}_i \subset \G_i[I^m]|_{\sU_{1,2}}$, so they are naturally quasi-coherent sheaves on $\X_{m-1}=\X\times_{\Spf R} \Spec R/I^{m}$. Therefore, one can use \cite[\href{https://stacks.math.columbia.edu/tag/01PF}{Tag 01PF}]{stacks-project} (applied to $\X_{m-1}$) to extend $\cal{K}_i$ to \[
        \widetilde{\cal{K}}_i \subset \G_i[I^m] \subset \G_i
    \]
    where $\widetilde{\cal{K}}_i$ adically quasi-coherent $\O_\X$-modules of finite type. Then $\G_i/\widetilde{\cal{K}}_i \to \F|_{\sU_i}$ are FP-$\m_0$-approximations of $\F|_{\sU_i}$ that are isomorphic on the intersection. Therefore, they glue to a global FP-$\m_0$-approximation $\G \to \F$.  
\end{proof}

\begin{thm}\label{thm:approximate-formal} Let $\X$ be a finitely presented formal $R$-scheme, $\F$ an almost finitely generated (resp. almost finitely presented) $\O_\X$-module. Then, for any finitely generated ideal $\m_0\subset \m$, there is an adically quasi-coherent, finitely generated (resp. finitely presented) $\O_\X$-module $\G$ and a map $\phi\colon \G \to \F$ such that $\m_0(\coker\phi)=0$ and $\m_0(\ker\phi)=0$.
\end{thm}
\begin{proof}
Without loss of generality, we can replace $\F$ by $\widetilde{\m}\otimes \F$, so we may and do assume that $\F$ is adically quasi-coherent. \smallskip

The case of almost adically quasi-coherent, almost finite type $\O_\X$-module $\F$ follows from Lemma~\ref{lemma:global-FP-approximation}. Indeed, there is an FP-$\m_0$-approximation $\phi'\colon \G' \to \F$, so we define $\phi\colon \G \to \F$ to be the natural inclusion $\G\coloneqq \rm{Im}(\phi') \to \F$. This gives the desired morphism as $\G$ is an adically quasi-coherent $\O_\X$-module of almost finite type by Corollary~\ref{ker-coker}. \smallskip

Now suppose $\F$ is an adically quasi-coherent, almost finitely presented $\O_\X$-module. Then we use Lemma~\ref{lemma:global-FP-approximation} to find an FP-$\m_0$-approximation $\phi'\colon \G' \to \F$. Now we note that any almost finitely presented $\O_\X$-module is almost coherent by Lemma~\ref{criterion-almost-coherent-formal}. Therefore, $\ker \phi$ is again adically quasi-coherent, almost finitely presented. Therefore, we can find an FP-$\m_0$-approximation $\phi''\colon \G'' \to  \ker(\phi')$ by Lemma~\ref{lemma:global-FP-approximation}. Denote by $\phi'''\colon \G'' \to G'$ the composition of $\phi''$ with the natural inclusion $\ker(\phi') \to \G'$. Now it is easy to check that $\phi\colon \coker (\phi''') \to \F$ gives the desired ``approximation''.
\end{proof}

\subsection{Formal Schemes. Derived category of almost coherent $\O_\X^a$-modules}\label{section:derived-category-formal-schemes}

We discuss the notion of the derived category of almost coherent sheaves on a formal scheme $\X$. One major issue is that, in this situation, the derived category of $\O_\X$-modules with adically quasi-coherent cohomology sheaves is not well-defined, as adically quasi-coherent sheaves do not form a weak Serre subcategory of $\bf{Mod}_\X$. \smallskip

To overcome this issue, we follow the strategy used in \cite{Lurie-spectral} and define another category $``\bf{D}_{qc}(\X)\text{''}$ completely on the level of derived categories. For the rest of the section, we fix a base ring $R$ as in Setup~\ref{set-up3}.

\begin{defn}\label{defn:derived-quasi-coherent-sheaves} Let $\X$ be a locally topologically finitely presented $R$-scheme. Then we define the {\it derived category of adically quasi-coherent sheave}s $``\bf{D}_{qc}(\X)\text{''}$ as a full subcategory of $\bf{D}(\X)$ consisting of objects $\F$ such that
\begin{itemize}
\item For every open affine $\sU \subset \X$, $\bf{R}\Gamma(\sU, \F)\in \bf{D}(\O_\X(\sU))$ is derived $I$-adically complete.
\item For every inclusion $\sU \subset \V$ of affine formal subschemes of $\X$, the natural morphism
\[
\bf{R}\Gamma(\V, \F) \wdh{\otimes}^L_{\O_\X(\V)} \O_{\X(\sU)} \to \bf{R}\Gamma(\sU, \F)
\]
is an isomorphism, where the completion is understood in the derived sense.
\end{itemize}
\end{defn}

\begin{rmk} We refer to \cite[\href{https://stacks.math.columbia.edu/tag/091N}{Tag 091N}]{stacks-project} and \cite[\href{https://stacks.math.columbia.edu/tag/0995}{Tag 0995}]{stacks-project} for a self-contained discussion on derived completions of modules and sheaves of modules respectively.
\end{rmk}

We now want to give an interpretation of $``\bf{D}_{qc}(\X)\text{''}$ in terms of $A$-modules for an affine formal scheme $\X=\Spf A$. We recall that in the case of schemes, we have a natural equivalence $\bf{D}_{qc}(\Spec A) \simeq \bf{D}(A)$ and the map is induced by $\bf{R}\Gamma(\Spec A, -)$. In the case of formal schemes, it is not literally true. We need to impose certain completeness conditions. 

\begin{defn} Let $A$ be a ring with a finitely generated ideal $I$. We define the {\it complete derived category} $\bf{D}_{comp}(A, I)\subset \bf{D}(A)$ as a full triangulated subcategory consisting of $I$-adically derived complete objects.
\end{defn}

Suppose that $\X=\Spf A$ is a topologically finitely presented affine formal $R$-scheme. We note that the natural functor $\bf{R}\Gamma(\X, -) \colon \bf{D}(\X) \to \bf{D}(A)$ induces a functor 
\[
\bf{R}\Gamma(\X, -) \colon ``\bf{D}_{qc}(\X)\text{''} \to \bf{D}_{comp}(A, I).
\]
We wish to show that this functor is an equivalence. For this, we need some preliminary lemmas:



\begin{lemma}\label{derived-Artin-Rees} Let $A$ be a topologically finitely presented $R$-algebra for $R$ as in Setup~\ref{set-up3}, let $f\in A$ be any element, and let $(x_1, \dots, x_d)=I$ be a choice of generators for the ideal of definition of $R$. Denote by $K(A_f; x_1^n, \dots, x_d^n)$ the Koszul complexes for the sequence $(x_1^n, \dots, x_d^n)$. Then the pro-systems $\{K(A_f; x_1^n, \dots, x_d^n)\}$ and $\{A_f/I^n\}$ are isomorphic in $\rm{Pro}(\bf{D}(A_f))$. 
\end{lemma}
\begin{proof}
The proof is the same \cite[\href{https://stacks.math.columbia.edu/tag/0921}{Tag 0921}]{stacks-project}. The only difference is that one needs to use \cite[Theorem 4.2.2(2)(b)]{FGK} in place of the usual Artin--Rees lemma. 
\end{proof}

\begin{lemma}\label{good-completion} Let $A$ be a topologically finitely presented $R$-algebra for $R$ as in Setup~\ref{set-up3}, let $f\in A$ be any element. Then the completed localization $A_{\{f\}}$ coincides with the $I$-adic derived completion of $A_f$.
\end{lemma}
\begin{proof}
Choose some generators $I=(x_1, \dots, x_d)$. Then we know that the derived completion completion of $A_f$ is given by $\bf{R}\lim_n K(A_f; x_1^n, \dots, x_d^n)$ where $K(A_f; x_1^n, \dots, x_d^n)$ is the Koszul complex for the sequence $(x_1^n, \dots, x_d^n)$. Lemma~\ref{derived-Artin-Rees} implies that the pro-systems $\{K(A_f; x_1^n, \dots, x_d^n)\}$ and $\{A_f/I^n\}$ are naturally pro-isomorphic. Thus we have an isomorphism
\[
\bf{R}\lim_n K(A_f; x_1^n, \dots, x_d^n) \cong \bf{R}\lim_n  A_f/I^n \simeq A_{\{f\}}.
\]
The last isomorphism uses the Mittag-Leffler criterion to ensure vanishing of $\lim^1$. 
\end{proof}

\begin{thm}\cite[Corollary 8.2.4.15]{Lurie-spectral}\label{updelta-equivalence} Let $\X=\Spf A$ be an affine, finitely presented formal scheme over $R$ as in Setup~\ref{set-up3}. Then the functor $\bf{R}\Gamma(\X, -) \colon ``\bf{D}_{qc}(\X)\text{''} \to \bf{D}_{comp}(A, I)$ is an equivalence of categories. 
\end{thm}
\begin{proof}
Lemma~\ref{good-completion} implies that the definition of $\Spf A$ in \cite{Lurie-spectral} is compatible with the classical one. Now \cite[Proposition 8.2.4.18]{Lurie-spectral} ensures that our definition of $``\bf{D}_{qc}(\X)\text{''}$ is equivalent to the homotopy category of $\rm{Qcoh}(\X)$ in the sense of \cite{Lurie-spectral}. Therefore, the result follows from \cite[Corollary 8.2.4.15]{Lurie-spectral} by passing to the homotopy categories.\smallskip

The proof of \cite[Corollary 8.2.4.15]{Lurie-spectral} can also be rephrased in our situation without using any derived geometry. However, it would require quite a long digression. 
\end{proof}

\begin{defn}\label{defn:updelta-derived} We denote by 
\[
(-)^{L\Updelta}\colon \bf{D}_{comp}(A, I) \to ``\bf{D}_{qc}(\X)\text{''}
\] 
the pseudo-inverse to $\bf{R}\Gamma(\X, -)\colon ``\bf{D}_{qc}(\X)\text{''} \to \bf{D}_{comp}(A, I)$. We note that 
\[
\bf{R}\Gamma(\Spf A_{\{f\}}, M^{L\Updelta}) \simeq M\wdh{\otimes}_A A_{\{f\}}
\]
for any $M\in \bf{D}_{comp}(A, I)$.
\end{defn}

\begin{rmk} The functor $(-)^{L\Updelta}$ {\it is not} compatible with the ``abelian'' functor $(-)^\Updelta$ used the previous sections. 
\end{rmk}

Now we define a category $\bf{D}_{qc, acoh}(\X)$ and show that it is equivalent to $\bf{D}_{acoh}(A)$. Theorem~\ref{updelta-equivalence} will be an important technical tool for establishing this equivalence.  

\begin{defn}\label{defn:almost-category-formal} We define $\bf{D}_{qc, acoh}(\X)$ (resp. $\bf{D}_{acoh}(\X)^a$) to be the full triangulated subcategory of $\bf{D}(\X)$ (resp. $\bf{D}(\X)^a$) consisting of complexes with adically quasi-coherent, almost coherent (resp. almost coherent) cohomology sheaves (resp. almost sheaves). 
\end{defn} 

\begin{rmk}\label{quotient-qc-acoh-formal} An argument similar to the one in the proof of Lemma~\ref{quotient-derived-acoh-je} shows that $\bf{D}_{acoh}(\X)^a$ is equivalent to the Verdier quotient $\bf{D}_{qc, acoh}(\X)/\bf{D}_{qc, \Sigma_\X}(\X)$.
\end{rmk}

In order to show an equivalence $\bf{D}_{qc, acoh}(\X)\simeq \bf{D}_{acoh}(A)$, our first goal is to show that $\bf{D}_{qc, acoh}$ lies inside $``\bf{D}_{qc}(\X)\text{''}$. This is not entirely obvious because the definition of $\bf{D}_{qc, acoh}(\X)$ imposes some restrictions on individual cohomology sheaves while the definition of $``\bf{D}_{qc}(\X)\text{''}$ on the whole complex itself. 

\begin{lemma}\label{derived-formal-schemes} Let $\X=\Spf A$ be an affine topologically finitely presented formal $R$-scheme for $R$ as in Setup~\ref{set-up3}. Then the functor $\mathbf{R}\Gamma(\X, -)\colon \mathbf{D}_{qc, acoh}(\X) \to \mathbf{D}(A)$ is $t$-exact (with respect to the evident $t$-structures on both sides) and factors through $\mathbf{D}_{acoh}(A)$. More precisely, there is an isomorphism 
\[
\rm{H}^i\left(\mathbf{R}\Gamma\left(\X, \F\right)\right) \simeq \mathrm{H}^0\left(\X, \mathcal{H}^i\left(\F\right)\right) \in \bf{Mod}_A^{\text{acoh}}
\]
for any object $\F\in \mathbf{D}_{qc, acoh}(\X)$. 
\end{lemma}
\begin{proof}
We note that the vanishing theorem \cite[Theorem I.7.1.1]{FujKato} implies that we can use \cite[\href{https://stacks.math.columbia.edu/tag/0D6U}{Tag 0D6U}]{stacks-project} with $N=0$. Thus, we see that the map $\rm{H}^i(\mathbf{R}\Gamma(\X, \F)) \to \mathrm{H}^i(\mathbf{R}\Gamma(\X, \tau^{\geq i}\F))$ is an isomorphism for any integer $i$, and that $\mathbf{R}\Gamma(\X, \F) \in \mathbf{D}_{acoh}(A)$ for any $\F \in \mathbf{D}_{qc, acoh}(\X)$. Combining it with the canonical isomorphism $\mathrm{H}^i(\mathbf{R}\Gamma(\X, \tau^{\geq i}\F)) \simeq \rm{H}^0(\X, \mathcal{H}^i(\F))$, we get the desired result.
\end{proof}

\begin{lemma}\label{acoh-qcoh} Let $\X$ be a locally topologically finitely presented formal $R$-scheme for $R$ as in Setup~\ref{set-up3}. Then $\bf{D}_{qc, acoh}(\X)$ is naturally a full triangulated subcategory of $``\bf{D}_{qc}(\X)\text{''}$.
\end{lemma}
\begin{proof}
Both $\bf{D}_{qc, acoh}(\X)$ and $``\bf{D}_{qc}(\X)\text{''}$ are full triangulated subcategories of $\bf{D}(\X)$. Thus, it suffices to show that any $\F\in \bf{D}_{qc, acoh}(\X)$ lies in $``\bf{D}_{qc}(\X)\text{''}$. \medskip 

Lemma~\ref{derived-formal-schemes} and Corollary~\ref{almost-coh-derived-complete} imply that $\bf{R}\Gamma(\sU, \F)\in \bf{D}_{comp}(A, I)$ for any open affine $\sU\subset \X$. Now suppose $\sU \subset \V$ is an inclusion of open affine formal subschemes in $\X$. We consider the natural morphism
\[
\bf{R}\Gamma(\V, \F) \wdh{\otimes}^L_{\O_\X(\V)} \O_{\X}(\sU) \to \bf{R}\Gamma(\sU, \F)
\]
We note that $\O_{\X}(\sU)$ is flat over $\O_\X(\V)$ by \cite[Proposition I.4.8.1]{FujKato}. Thus, the complex 
\[
    \bf{R}\Gamma(\V, \F) \otimes^L_{\O_\X(\V)} \O_{\X(\sU)}
\]
lies in $\bf{D}_{acoh}(\O_{\X}(\sU))$ by Lemma~\ref{trivial-base-change}. Therefore, it also lies in $\bf{D}_{comp}(A, I)$ by Corollary~\ref{almost-coh-derived-complete}. So we conclude that 
\[
\bf{R}\Gamma(\V, \F) \wdh{\otimes}^L_{\O_\X(\V)} \O_{\X}(\sU) \simeq \bf{R}\Gamma(\V, \F) \otimes^L_{\O_\X(\V)} \O_{\X}(\sU).
\]
Using $\O_\X(\V)$-flatness of $\O_\X(\sU)$, we conclude that the question boils down to showing that 
\[
\rm{H}^i(\V, \F) \otimes_{\O_\X(\V)} \O_{\X(\sU)} \to \rm{H}^i(\sU, \F)
\]
is an isomorphism for all $i$. Now Lemma~\ref{derived-formal-schemes} implies that this, in turn, reduces to showing that the natural map
\[
\Gamma(\V, \mathcal{H}^i(\F)) \otimes_{\O_\X(\V)} \O_{\X(\sU)} \to \Gamma(\sU, \mathcal{H}^i(\F))
\]
is an isomorphism. Without loss of generality, we may assume that $\X=\V=\Spf A$. Then $\mathcal{H}^i(\F)$ is an adically quasi-coherent, almost coherent $\O_\X$-module, so Lemma~\ref{sheafify-modules-formal} ensures that it is isomorphic to $M^\Updelta$ for some $M\in \bf{Mod}_A^{\text{acoh}}$. Thus, the desired claim follows from \cite[Lemma 3.6.4]{FujKato} and the observation that $M \otimes_{\O_\X(\V)} \O_{\X}(\sU)$ is already $I$-adically complete due to Lemma~\ref{completion-finitely-generated}.
\end{proof}

Now we show that the $(-)^{L\Updelta}$ functor sends $\bf{D}_{acoh}(A)$ to $\bf{D}_{qc, acoh}(\Spf A)$. This is also not entirely obvious as $(-)^{L\Updelta}$ is a priori different from the classical version of the $(-)^\Updelta$-functor studied in previous sections. However, we show that these functors coincide on $\bf{Mod}_A^{\text{acoh}}$. 

\begin{lemma}\label{updelta-derived} Let $\X=\Spf A$ be an affine topologically finitely presented formal $R$-scheme for $R$ as in Setup~\ref{set-up3}. Then the functor $(-)^{L\Updelta} \colon \bf{D}_{acoh}(A) \to ``\bf{D}_{qc}(\X)\text{''}$ factors through $\bf{D}_{qc, acoh}(\X)$. Moreover, for any $M\in \bf{D}_{acoh}(A)$ and an integer $i$, there is a functorial isomorphism
\[
\rm{H}^i(M)^{\Updelta} \simeq \mathcal{H}^i(M^{L\Updelta}).
\]
\end{lemma}
\begin{proof}
We note that $\rm{H}^i(\X, M^{L\Updelta})\simeq \rm{H}^i(M)$ due to its construction. Since $\mathcal{H}^i(M^{L\Updelta})$ is canonically isomorphic to the sheafification of the presheaf
\[
\sU \mapsto \rm{H}^i(\sU, M^{L\Updelta}),
\]
we get that there is a canonical map $\rm{H}^i(M) \to \Gamma(\X, \mathcal{H}^i(M^\Updelta))$. By the universal property of the classical $(-)^\Updelta$ functor, we get a functorial morphism
\[
\rm{H}^i(M)^\Updelta \to \mathcal{H}^i(M^{L\Updelta}).
\]
Since $\rm{H}^i(M)$ is almost coherent, we only need to show that this map is an isomorphism for any $i$. This boils down (using almost coherence of $\rm{H}^i(M)$) to showing that the natural morphism
\[
\rm{H}^i(M) \otimes_{A} A_{\{f\}} \to \mathrm{H}^i(\Spf A_{\{f\}}, M^{L\Updelta}).
\]
is an isomorphism for all $f\in A$. Now recall that $\bf{R}\Gamma(\Spf A_{\{f\}}, M^{L\Updelta})\simeq M\wdh{\otimes}^L_A A_{\{f\}}$ for any $f\in A$. Using that $M\in \bf{D}_{acoh}(A)$, $A_{\{f\}}$ is flat over $A$, and that almost coherent complexes are derived complete by Lemma~\ref{almost-coh-derived-complete}, we conclude that the natural map
\[
\rm{H}^i(M) \otimes_A A_{\{f\}} \to \rm{H}^i(\Spf A_{\{f\}}, M^{L\Updelta}) 
\]
is an isomorphism finishing the proof. 
\end{proof}

\begin{cor}\label{updelta-pr-coh} Let $\X=\Spf A$ be an affine topologically finitely presented formal $R$-scheme for $R$ as in Setup~\ref{set-up3}. Suppose that $M\in \bf{D}(A)$ has almost zero cohomology modules. Then $\cal{H}^i(M^{L\Updelta})$ is an almost zero, adically quasi-coherent $\O_\X$-module for each integer $i$. In particular, $(-)^{L\Updelta}$ induces a functor $(-)^{L\Updelta}\colon \bf{D}_{acoh}(A)^a \to \bf{D}_{acoh}(\X)^a$.
\end{cor}
\begin{proof}
We note that any almost zero $A$-modules is almost coherent, thus the result follows directly from the formula $\rm{H}^i(M)^\Updelta \simeq \mathcal{H}^i(M^{L\Updelta})$ established in Lemma~\ref{updelta-derived}. 
\end{proof}

\begin{thm}\label{equiv-derived-formal} Let $\X=\Spf A$ be an affine topologically finitely presented formal $R$-scheme for $R$ as in Setup~\ref{set-up3}. Then the functor $\mathbf{R}\Gamma(\X, -)\colon \mathbf{D}_{qc, acoh}(\X) \to \mathbf{D}_{acoh}(A)$ is a $t$-exact equivalence of triangulated categories with the pseudo-inverse $(-)^{L\Updelta}$.
\end{thm}
\begin{proof}
Lemma~\ref{derived-formal-schemes} implies that $\bf{R}\Gamma(\X, -)$ induces the functor $\mathbf{D}_{qc, acoh}(\X) \to \mathbf{D}_{acoh}(A)$ and that this functor is $t$-exact. Lemma~\ref{acoh-qcoh} and Theorem~\ref{updelta-equivalence} ensure that it is sufficient to show that $(-)^{L\Updelta}$ sends $\bf{D}_{acoh}(A)$ to $\bf{D}_{qc, acoh}(\X)$, this follows from Lemma~\ref{updelta-derived}.
\end{proof}

Now we can pass to the almost categories using Remark~\ref{quotient-qc-acoh-formal} to get the almost version of Theorem~\ref{equiv-derived-formal}.

\begin{cor}\label{cor:derived-almost-coherent-sheaves-formal-affine-schemes} Let $\X=\Spf A$ be an affine topologically finitely presented formal $R$-scheme for $R$ as in Setup~\ref{set-up3}. Then the functor $\mathbf{R}\Gamma(\X, -)\colon \mathbf{D}_{acoh}(\X)^a \to \mathbf{D}_{acoh}(A)^a$ is a $t$-exact equivalence of triangulated categories with the pseudo-inverse $(-)^{L\Updelta}$.
\end{cor}

\subsection{Formal Schemes. Basic functors on derived categories of $\O_\X^a$-modules}\label{section:derived-functors-formal-schemes}

We discuss the derived analogue of the main results of Section~\ref{basic-formal}. We show that the derived completion, derived tensor product, derived pullback, and derived almost Hom functors preserve complexes with almost coherent cohomology sheaves under certain conditions. For the rest of the section, we fix a ring $R$ as in Setup~\ref{set-up3}. \smallskip

We start with the completion functor. We recall that we have defined the morphism of locally ringed spaces $c\colon \wdh{X} \to X$ for any $R$-scheme $X$. If $X$ is locally finitely presented over $R$ or $X=\Spec A$ for a topologically finitely presented $R$-algebra $A$, then $c$ is a flat morphism as shown in Lemma~\ref{completion-affine} and Corollary~\ref{completion-global}. 

\begin{lemma}\label{affine-completion-derived} Let $X=\Spec A$ be an affine $R$-scheme for $R$ as in Setup~\ref{set-up3}. Suppose that $A$ is either finitely presented or topologically finitely presented over $R$. Suppose $M\in \bf{D}_{acoh}(A)$. Then $M^{L\Updelta} \simeq \bf{L}c^*(\widetilde{M})$.
\end{lemma}
\begin{proof}
First of all, we show that $\bf{L}c^*(\widetilde{M})\in \bf{D}_{qc, acoh}(\wdh{X})$. Indeed, the functor $c^*$ is exact as $c$ is flat. Thus, Lemma~\ref{completion-affine} guarantees that we have a sequence of isomorphisms 
\[
\mathcal{H}^i\left(\bf{L}c^*\left(\widetilde{M}\right)\right) \simeq c^*\left(\widetilde{\rm{H}^i\left(M\right)}\right)\simeq \left(\rm{H}^i\left(M\right)\right)^\Updelta.
\]
In particular, Theorem~\ref{updelta-equivalence} ensures that the natural morphism 
\[
M \simeq \bf{R}\Gamma(X, \widetilde{M}) \to \bf{R}\Gamma(\wdh{X}, \bf{L}c^*(\widetilde{M}))
\]
induces the morphism $M^{L\Updelta} \to \bf{L}c^*(\widetilde{M})$. As $c^*$ is exact, Lemma~\ref{updelta-derived} implies that it is sufficient to show that the natural map
\[
\rm{H}^i(M)^\Updelta \to c^*(\widetilde{\rm{H}^i(M)})
\]
is an isomorphism for all $i$. This follows from Lemma~\ref{completion-affine}.
\end{proof}

\begin{cor}\label{completion-global-derived} Let $X$ be a locally finitely presented $R$-scheme for $R$ as in Setup~\ref{set-up3}. Then $\bf{L}c^*$ induces functors $\bf{L}c^*\colon \bf{D}_{qc, acoh}^*(X) \to \bf{D}^*_{qc, acoh}(\wdh{X})$ (resp. $\bf{L}c^*\colon \bf{D}_{acoh}^*(X)^a \to \bf{D}^*_{acoh}(\wdh{X})^a$) for any $*\in \{``\text{ ''}, - , b , +\}$.
\end{cor}
\begin{proof}
The claim is local, so it suffices to assume that $X=\Spec A$. Then it follows from the exactness of $c^*$ and Lemma~\ref{affine-completion-derived}.
\end{proof}

\begin{lemma}\label{pullback-almost-coh-formal-derived} Let $\mf\colon \X \to \Y$ be a morphism of locally finitely presented formal $R$-schemes for $R$ as in Setup~\ref{set-up3}.  
\begin{enumerate}[label=\textbf{(\arabic*)}]
\item\label{pullback-almost-coh-formal-derived-1} Suppose that $\X=\Spf B$, $\Y=\Spf A$ are affine formal $R$-schemes. Then there is a functorial isomorphism
\[
\bf{L}\mf^*\left(M^{L\Updelta}\right) \simeq \left(M\otimes_{A} B\right)^{L\Updelta}
\]
for any $M\in \bf{D}_{acoh}(A)$.
\item\label{pullback-almost-coh-formal-derived-1.5} Suppose that $\X=\Spf B$, $\Y=\Spf A$ are affine formal $R$-schemes. Then there is a functorial isomorphism
\[
\bf{L}\mf^*\left(M^{a, L\Updelta}\right)\simeq \left(M^a\otimes_{A^a} B^a\right)^{L\Updelta}
\]
for any $M^a\in \bf{D}_{acoh}(A)$.
\item\label{pullback-almost-coh-formal-derived-2} The functor $\bf{L}\mf^*$ carries $\bf{D}^-_{qc, acoh}(\Y)$ to $\bf{D}^-_{qc, acoh}(\X)$.
\item\label{pullback-almost-coh-formal-derived-3} The functor $\bf{L}\mf^*$ carries $\bf{D}^-_{acoh}(\Y)^a$ to $\bf{D}^-_{acoh}(\X)^a$.
\end{enumerate}
\end{lemma}
\begin{proof}
The proof is similar to that of Lemma~\ref{pullback-almost-coh-formal}. We use Lemma~\ref{affine-completion-derived} and Lemma~\ref{updelta-derived} to reduce to the analogous algebraic facts that were already proven in Lemma~\ref{pullback-almost-coh}.
\end{proof}

\begin{lemma}\label{tensor-almost-coh-formal-derived} Let $\X$ be a locally topologically finitely presented formal $R$-scheme for $R$ as in Setup~\ref{set-up3}.  
\begin{enumerate}[label=\textbf{(\arabic*)}]
\item\label{tensor-almost-coh-formal-derived-1} Suppose that $\X=\Spf A$ is affine. Then there is a functorial isomorphism
\[
M^{L\Updelta}\otimes^L_{\O_\X}N^{L\Updelta}\simeq (M\otimes^L_{A} N)^{L\Updelta}
\]
for any $M$, $N\in \bf{D}_{acoh}(A)$. 
\item\label{tensor-almost-coh-formal-derived-2} Suppose that $\X=\Spf A$ is affine. Then there is a functorial isomorphism 
\[
M^{a, L\Updelta}\otimes^L_{\O_\X^a}N^{a, L\Updelta} \simeq (M^a\otimes^L_{A^a} N^a)^{L\Updelta}
\]
for any $M^a$, $N^a\in \bf{D}_{acoh}(A)^{a}$. 
\item\label{tensor-almost-coh-formal-derived-3} Let $\F$, $\G\in \bf{D}^-_{qc, acoh}(\X)$. Then $\F\otimes^L_{\O_\X}\G \in \bf{D}^-_{qc, acoh}(\X)$.   
\item\label{tensor-almost-coh-formal-derived-4} Let $\F^a$, $\G^a\in \bf{D}^-_{acoh}(\X)^a$. Then $\F^a\otimes^L_{\O_\X^a}\G^a \in \bf{D}^-_{acoh}(\X)^a$.  
\end{enumerate}
\end{lemma}
\begin{proof}
Similarly to Lemma~\ref{pullback-almost-coh-formal-derived}, we use Lemma~\ref{affine-completion-derived} and Lemma~\ref{updelta-derived}  to reduce to the analogous algebraic facts that were already proven in Lemma~\ref{tensor-almost-coh}.
\end{proof}

Now we discuss the $\bf{R}\ud{al\cal{H}om}_{\O_\X}(-, -)$-functor. Our strategy of showing that $\bf{R}\ud{al\cal{H}om}(-, -)$ preserves almost coherent complexes will be slightly different from the schematic case. The main technical problem is to define the map $\bf{R}\rm{alHom}_{A^a}(M^a, N^a)^{L\Updelta} \to \bf{R}\ud{al\cal{H}om}_{\O_\X^a}(M^{a, L\Updelta}, N^{a, L\Updelta})$ in the affine case. 

The main issue is that we do not know if $(-)^{L\Updelta}$ is a left adjoint to the functor of global section on the whole category $\bf{D}(\X)$; we only know that it becomes a pseudo-inverse to $\bf{R}\Gamma(\X, -)$ after restriction to $``\bf{D}_{qc}(\X)\text{''}$. However, the complex $\bf{R}\ud{\cal{H}om}_{\O_\X}(M^{L\Updelta}, N^{L\Updelta})$ itself usually does not lie inside $``\bf{D}_{qc}(\X)\text{''}$. To overcome this issue, we will show that 
\[
\widetilde{\m}\otimes \bf{R}\ud{\cal{H}om}_{\O_\X}(M^{L\Updelta}, N^{L\Updelta})
\]
does lie in $``\bf{D}_{qc}(\X)\text{''}$ for $M\in \bf{D}^-_{acoh}(A)$ and $N\in \bf{D}^+_{acoh}(A)$. \smallskip

Since $``\bf{D}_{qc}(\X)\text{''}$ was defined in a bit abstract way, it is probably the easiest way to show that $\widetilde{\m}\otimes \bf{R}\ud{\cal{H}om}_{\O_\X}(M^{L\Updelta}, N^{L\Updelta})$ actually lies in $\bf{D}_{qc, acoh}(\X)$. That is sufficient by Lemma~\ref{acoh-qcoh}. 

\begin{lemma}\label{prelim-rhom-formal} Let $\X=\Spf A$ be a topologically finitely presented formal $R$-scheme for $R$ as in Setup~\ref{set-up3}. Let $M, N \in \bf{Mod}_A^{\text{acoh}}$ there is a natural almost isomorphism
\[
\rm{Ext}^p_A(M, N)^\Updelta \xr{\sim} \cal{E}xt^p_{\O_\X}(M^\Updelta, N^\Updelta)
\]
for every integer $p$. 
\end{lemma}
\begin{proof}
We recall that $\cal{E}xt^p_{\O_\X}(M^\Updelta, N^\Updelta)$ is canonically isomorphic to the sheafification of the presheaf
\[
\sU \mapsto \rm{Ext}^p_{\O_{\sU}} (M^\Updelta|_{\sU}, N^\Updelta|_{\sU}).
\]
In particular, there is a canonical map $\rm{Ext}^p_{\O_\X}(M^{\Updelta}, N^\Updelta) \to \Gamma(\X, \cal{E}xt^p_{\O_\X}(M^\Updelta, N^\Updelta))$. It induces a morphism
\begin{equation}\label{Hom-map}
\rm{Ext}^p_{\O_\X}(M^{\Updelta}, N^\Updelta)^\Updelta \to \cal{E}xt^p_{\O_\X}(M^\Updelta, N^\Updelta).
\end{equation}
Now we note that the classical $(-)^\Updelta$ functor and the derived version coincide on almost coherent modules by Lemma~\ref{updelta-derived}. Hence, the equivalence $``\bf{D}_{qc}(\X)\text{''}\simeq \bf{D}_{comp}(A, I)$ coming from Theorem~\ref{updelta-equivalence} and Lemma~\ref{updelta-derived} ensures that $\rm{Ext}^p_{\O_\X}(M^{\Updelta}, N^\Updelta) \simeq \rm{Ext}^p_A(M, N)$. So the map~(\ref{Hom-map}) becomes a map
\[
\rm{Ext}^p_A(M, N)^\Updelta \to \cal{E}xt^p_{\O_\X}(M^\Updelta, N^\Updelta).
\]
We note that $\rm{Ext}^p_A(M, N)$ is an almost coherent $A$-module by Proposition~\ref{alHom-derived-coh}. Using that almost coherent modules are complete, we conclude that it suffices to show that
\[
\rm{Ext}^p_A(M, N) \otimes_A A_{\{f\}} \to \rm{Ext}^p_{\Spf A_{\{f\}}}(M^\Updelta|_{\Spf A_{\{f\}}}, N^\Updelta|_{\Spf A_{\{f\}}}) 
\]
is an almost isomorphism. Using Lemma~\ref{pullback-almost-coh-formal} and the equivalence $``\bf{D}_{qc}(\X)\text{''}\simeq \bf{D}_{comp}(A, I)$ as above, we see that the map above becomes the canonical map
\[
\rm{Ext}^p_A(M, N) \otimes_A A_{\{f\}}  \to \rm{Ext}^p_{A_{\{f\}}}(M\otimes_A A_{\{f\}}, N\otimes_{A} A_{\{f\}}).
\]
Finally, this map is an almost isomorphism by Proposition~\ref{base-change-hom-derived}. 
\end{proof}

\begin{cor}\label{rhom-qcoh} Let $\X$ be a locally topologically finitely presented formal $R$-scheme for $R$ as in Setup~\ref{set-up3}. Then
\[
\widetilde{\m}\otimes \bf{R}\ud{\cal{H}om}_{\O_\X}(\F, \G)\in \bf{D}^+_{qc, acoh}(\X)
\]
for $\F \in \bf{D}^-_{qc, acoh}(\X)$, and $\G \in \bf{D}^+_{qc, acoh}(\X)$.
\end{cor}
\begin{proof}
The claim is local, so we can assume that $\X=\Spf A$. Then we use the Ext-spectral sequence and Lemma~\ref{main-almost-sheaf-formal} to reduce to the case when $\F$ and $\G$ are in $\bf{Mod}_\X^{\rm{qc}, \rm{acoh}}$. Thus, Theorem~\ref{sheafify-modules-almost-formal} ensures that $\F=M^\Updelta$ and $\G=N^\Updelta$ for some $M$, $N\in \bf{Mod}_A^{\rm{acoh}}$. So Lemma~\ref{prelim-rhom-formal} guarantees that 
\[
\cal{H}^p\left(\bf{R}\ud{\cal{H}om}_{\O_\X}(\F, \G) \right) \simeq^a \rm{Ext}^p_A(M, N)^\Updelta.
\]
In other words, 
\[
\widetilde{\m}\otimes \cal{H}^p\left(\bf{R}\ud{\cal{H}om}_{\O_\X}(\F, \G) \right) \simeq \widetilde{\m}\otimes \rm{Ext}^p_A(M, N)^\Updelta.
\]
Now $\rm{Ext}^p_A(M, N)^\Updelta$ is an adically quasi-coherent, almost coherent $\O_\X$-module by Proposition~\ref{alHom-derived-coh} and Lemma~\ref{sheafify-modules-formal}. So Lemma~\ref{qc-aqc} guarantees that $\widetilde{\m}\otimes \rm{Ext}^p_A(M, N)^\Updelta$ is also adically quasi-coherent and almost coherent. Therefore, $\widetilde{\m}\otimes \bf{R}\ud{\cal{H}om}_{\O_\X}(\F, \G)\in \bf{D}^+_{qc, acoh}(\X)$.
\end{proof}

\begin{lemma}\label{hom-alcoh-formal-derived} Let $\X$ be a locally topologically finitely presented formal $R$-scheme for $R$ as in Setup~\ref{set-up3}. 
\begin{enumerate}[label=\textbf{(\arabic*)}]
\item\label{hom-alcoh-formal-derived-1} Suppose $\X=\Spf A$ is affine. Then there is a functorial isomorphism 
\begin{equation*}
\bf{R}\rm{alHom}_{A^a}(M^a, N^a)^{L\Updelta} \to \bf{R}\ud{al\cal{H}om}_{\O_\X^a}(M^{a, L\Updelta}, N^{a, L\Updelta})
\end{equation*}
for $M\in \bf{D}^-_{acoh}(A)^a$ and $N\in \bf{D}^+_{acoh}(A)^a$. 
\item\label{hom-alcoh-formal-derived-2} Suppose $\F^a\in \bf{D}^+_{acoh}(\X)^a$ and $\G^a\in \bf{D}^-_{acoh}(\X)$ are almost coherent $\O_\X^a$-modules. Then $\bf{R}\ud{al\cal{H}om}_{\O_\X^a}(\F^a, \G^a) \in \bf{D}^+_{acoh}(\X)^a$.
\end{enumerate}
\end{lemma}
\begin{proof}
We start with \ref{hom-alcoh-formal-derived-1}. Proposition~\ref{derived-al-hom-sheaf} implies that the map 
\[
(\widetilde{\m}\otimes \bf{R}\ud{\cal{H}om}_{\O_\X}(M^\Updelta, N^\Updelta))^a \to  \bf{R}\ud{al\cal{H}om}_{\O_\X^a}(M^{a, \Updelta}, N^{a, \Updelta})
\]
is an isomorphism in $\bf{D}(\X)^a$. Similarly, the map
\[
(\widetilde{\m} \otimes \bf{R}\rm{Hom}_{A}(M, N)^{\Updelta})^a \to \bf{R}\rm{alHom}_{A^a}(M^a, N^a)^{\Updelta}
\]
is an isomorphism by Lemma~\ref{prelim-rhom-formal}. Thus, it suffices to construct a functorial isomorphism
\[
\widetilde{\m} \otimes \bf{R}\rm{Hom}_{A}(M, N)^{L\Updelta} \to \widetilde{\m} \otimes \bf{R}\ud{\cal{H}om}_{\O_\X}(M^{L\Updelta}, N^{L\Updelta}).
\]

Now Lemma~\ref{updelta-derived} and Corollary~\ref{rhom-qcoh} guarantee that 
\[
\widetilde{\m} \otimes \bf{R}\ud{\cal{H}om}_{\O_\X}(M^{L\Updelta}, N^{L\Updelta}) \in \bf{D}_{qc, acoh}(\X).
\]
Proposition~\ref{alHom-derived-coh}, Lemma~\ref{sheafify-modules-formal}, and Lemma~\ref{qc-aqc} also guarantee that 
\[
\widetilde{\m} \otimes \bf{R}\rm{Hom}_{A}(M, N)^{\Updelta} \in \bf{D}_{qc, acoh}(\X).
\]
Thus, Theorem~\ref{updelta-equivalence} ensures that, in order to construct the desired isomorphism, it suffices to do it after applying $\bf{R}\Gamma(\X, -)$. The Projection Formula (see Theorem~\ref{projection}) and the definition of the $(-)^{L\Updelta}$-functor provide us with functorial isomorphisms
\[
\bf{R}\Gamma\left(\X, \widetilde{\m} \otimes \bf{R}\rm{Hom}_{A}(M, N)^{L\Updelta}\right)\simeq \widetilde{\m} \otimes \bf{R}\Hom_A(M, N)
\]
\begin{align*}
\bf{R}\Gamma\left(\X, \widetilde{\m} \otimes \bf{R}\ud{\cal{H}om}_{\O_\X}(M^{L\Updelta}, N^{L\Updelta})\right) & \simeq \widetilde{\m} \otimes \bf{R}\Gamma\left(\X, \bf{R}\ud{\cal{H}om}_{\O_\X}(M^{L\Updelta}, N^{L\Updelta})\right) \\
& \simeq \widetilde{\m} \otimes \bf{R}\rm{Hom}_{\O_X}(M^{L\Updelta}, N^{L\Updelta}) \\
&\simeq \widetilde{\m}\otimes \bf{R}\Hom_A\left(M, N\right)
\end{align*}
where the last isomorphism uses equivalence from Theorem~\ref{updelta-equivalence}. Thus, we see 
\[
\bf{R}\Gamma\left(\X, \widetilde{\m} \otimes \bf{R}\rm{Hom}_{A}(M, N)^{L\Updelta}\right) \simeq \bf{R}\Gamma\left(\X, \widetilde{\m} \otimes \bf{R}\ud{\cal{H}om}_{\O_\X}(M^{L\Updelta}, N^{L\Updelta})\right).
\]
As a consequence, we have a functorial isomorhism
\[
\widetilde{\m} \otimes \bf{R}\rm{Hom}_{A}\left(M, N\right)^{L\Updelta} \xr{\sim} \widetilde{\m} \otimes \bf{R}\ud{\cal{H}om}_{\O_\X}\left(M^{L\Updelta}, N^{L\Updelta}\right).
\]
This induces the desired isomorphism 
\[
\bf{R}\rm{alHom}_{A^a}\left(M^a, N^a\right)^{L\Updelta} \xr{\sim} \bf{R}\ud{al\cal{H}om}_{\O_\X^a}\left(M^{a, L\Updelta}, N^{a, L\Updelta}\right).
\]

\ref{hom-alcoh-formal-derived-2} is an easy consequence of \ref{hom-alcoh-formal-derived-1}, Proposition~\ref{alHom-derived-coh}, and Corollary~\ref{updelta-pr-coh}. 
\end{proof}

\section{Cohomological properties of almost coherent sheaves}\label{section:cohomological-properties}

The main goal of this section is to establish that almost coherent sheaves share similar cohomological properties as classical coherent sheaves. In particular, we prove almost versions of the Proper Mapping Theorem (both for schemes and nice formal schemes), of the formal GAGA Theorem, of the Formal Function Theorem, and of the Grothendieck Duality. The Formal GAGA Theorem is arguably quite surprising in the almost coherent context because almost coherent sheaves are rarely of finite type, so none of the classical proofs of the Formal GAGA Theorem applies in this situation. We resolve this issue by adapting a new approach to GAGA Theorems due to J.~Hall (see \cite{JH}).

\subsection{Almost Proper Mapping Theorem}\label{APMT-formal}

The main goal of this section is to prove the Almost Proper Mapping Theorem which says that derived pushforward along proper (topologically) finitely presented morphism of nice (formal) schemes preserves almost coherent sheaves. 

The idea of the proof is relatively easy: we approximate an almost finitely presented $\O_X$-module by finitely presented using Corollary~\ref{cor:approximate-afpr} or Theorem~\ref{thm:approximate-formal} and then use the usual Proper Mapping Theorem. For this, we will need a version of the Proper Mapping Theorem for a class of non-noetherian rings, which we review below. 

\begin{defn} We say that a scheme $Y$ is {\it universally coherent} if any scheme $X$ that is locally of finite presentation over $Y$ is coherent (i.e.~the structure sheaf $\O_X$ is coherent).
\end{defn} 

\begin{thm}[Proper Mapping Theorem]\label{proper-mapping}\cite[Theorem I.8.1.3]{FujKato} Let $Y$ be a universally coherent quasi-compact scheme, and let $f\colon X \to Y$ be a proper morphism of finite presentation. Then the functor $\mathbf{R}f_*$ sends $\mathbf{D}^*_{coh}(X)$ to  $\mathbf{D}^*_{coh}(Y)$ for any $*\in \{\text{`` ''}, +, -, b\}$.
\end{thm}

We want to generalize this theorem to the ``almost world''. So we pick a ring $R$ and a fixed ideal $\m\subset R$ such that $\m^2=\m$ and $\widetilde{\m}=\m\otimes_R \m $ is $R$-flat. In this section, we always consider almost mathematics with respect to this ideal. 

\begin{thm}[Almost Proper Mapping Theorem]\label{almost-proper-mapping} Let $Y$ be a universally coherent quasi-compact $R$-scheme, and let $f\colon X \to Y$ be a proper, finitely presented morphism. Then 
\begin{itemize}\itemsep0.5em
\item The functor $\mathbf{R}f_*$ sends $\mathbf{D}^*_{qc, acoh}(X)$ to $\mathbf{D}^*_{qc, acoh}(Y)$ for any $*\in \{\text{`` ''}, +, -, b\}$. 
\item The functor $\mathbf{R}f_*$ sends $\mathbf{D}^*_{acoh}(X)^a$ to $\mathbf{D}^*_{acoh}(Y)^a$ for any $*\in \{\text{`` ''}, +, -, b\}$. 
\item The functor $\mathbf{R}f_*$ sends $\mathbf{D}^+_{acoh}(X)$ to $\mathbf{D}^+_{acoh}(Y)$. 
\item If $Y$ has finite Krull dimension, then $\mathbf{R}f_*$ sends $\mathbf{D}^*_{acoh}(X)$  to $\mathbf{D}^*_{acoh}(Y)$ for any $*\in \{\text{`` ''}, +, -, b\}$.
\end{itemize}
\end{thm}


\begin{lemma}\label{finite-coh-dimension-sch} Let $Y$ be a quasi-compact scheme of finite Krull dimension, and let $f\colon X \to Y$ be a finite type, quasi-separated morphism. Then $X$ has finite Krull dimension, and $f_*$ has finite cohomological dimension on $\mathbf{Mod}_X$.
\end{lemma}
\begin{proof}

First of all, we show that $X$ has finite Krull dimension. Indeed, the morphism $f\colon X \to Y$ is quasi-compact, therefore $X$ is quasi-compact. So it suffices to show that locally $X$ has finite Krull dimension. Thus, we can assume that $X=\Spec B$ and $Y=\Spec A$ are affine, and the map is given by a finite type morphism $A\to B$. In this case, we have $\dim Y =\dim A$ and $\dim X = \dim B$. Thus, it is enough to show that the Krull dimension of a finite type $A$-algebra is finite. This readily reduces the question to the case of a polynomial algebra $\dim A[X_1, \dots, X_n]$. Now \cite[Chapter 11, Exercise 6]{AM} implies that $\dim A[X_1, \dots, X_n] \leq \dim A +2n$. 

Now we prove that $f_*$ has finite cohomological dimension. We note that it suffices to show that there is an integer $N$ such that, for any open affine $U\subset Y$, cohomology groups $\rm{H}^i(X_U, \F)$ vanish for $i\geq N$ and any $\O_{X_U}$-module $\F$. We recall that $f$ is quasi-separated, so $X_U$ is quasi-compact, quasi-separated and $\dim X_U \leq \dim X$ for any open $U\subset X$. Therefore, it suffices to show that on any spectral space $X$, we have $\rm{H}^i(X, \F)=0$ for $i > \dim X$ and $\F \in \mathcal{A}b(X)$. This is proven in \cite[Corollary 4.6]{Scheiderer} (another reference is \cite[\href{https://stacks.math.columbia.edu/tag/0A3G}{Tag 0A3G}]{stacks-project}). Thus, we conclude that $N=\dim X$ does the job. 
\end{proof}

\begin{proof}[Proof of Theorem~\ref{almost-proper-mapping}]

{\it Step 0. Reduction to the case of bounded below derived categories:} We note that $f_*$ has a bounded cohomological dimension on $\mathbf{Mod}_X^{\qc}$. Indeed, for any quasi-compact separated scheme $X$ and $\F \in \mathbf{Mod}_X^{\qc}$, we can compute $\rm{H}^i(X, \F)$ via the {\it alternating} \v{C}ech complex for some finite affine covering of $X$. Therefore, if $X$ can be covered by $N$ affines, the functor $f_*$ restricted to $\mathbf{Mod}_X^{\qc}$ has cohomological dimension at most $N$. \smallskip

Now we use \cite[\href{https://stacks.math.columbia.edu/tag/0D6U}{Tag 0D6U}]{stacks-project} (alternatively, one can use \cite[Lemma 3.4]{benlim}) to reduce the question of proving the claim for any $\F\in \mathbf{D}_{qc, acoh}(X)$ to the question of proving the claim for all its truncations $\tau^{\geq a}\F$. In particular, we reduce the case of $\F\in \mathbf{D}_{qc, acoh}(X)$ to the case when $\F\in \bf{D}^+_{qc, acoh}(X)$. Similarly (using Proposition~\ref{derived-pushforward}), we reduce the case of $\F^a\in \bf{D}_{acoh}(X)^a$ to the case of $\F^a\in \bf{D}^+_{acoh}(X)^a$. \smallskip

Using Lemma~\ref{finite-coh-dimension-sch}, a similar argument also allows us to reduce the case of $\F\in \bf{D}_{acoh}(X)$ to the case of $\F\in \bf{D}^+_{acoh}(Y)$ when $Y$ has finite Krull dimension. \smallskip



{\it Step 1. Reduction to the case of quasi-coherent almost coherent sheaves:} Using the Projection Formula (Lemma~\ref{projection}) (resp. Proposition~\ref{derived-pushforward}), we see that, in order to show $\mathbf{R}f_*$ sends $\mathbf{D}^+_{acoh}(X)$ to $\mathbf{D}^+_{acoh}(Y)$ (resp. $\mathbf{D}^+_{acoh}(X)^a$ to $\mathbf{D}^+_{acoh}(Y)^a$), it suffices to show the analogous result for $\mathbf{D}^+_{qc, acoh}(X)$. Moreover, we can use the spectral sequence
\[
\rm{E}^{p,q}_2= \mathrm{R}^pf_*\mathcal{H}^q(\F) \Rightarrow \mathbf{R}^{p+q}f_*(\F)
\]
to reduce the claim to the fact that higher derived pushforwards of a quasi-coherent, almost coherent sheaf are quasi-coherent and almost coherent. \smallskip

 {\it Step 2. The case of a quasi-coherent, almost coherent $\O_X$-module $\F$}: We show that $\mathrm{R}^if_*\F$ is a quasi-coherent, almost coherent $\O_Y$-module for any quasi-coherent, almost coherent $\O_X$-module $\F$ and any $i$. First, we note that $\mathrm{R}^if_*\F$ is quasi-coherent as higher pushforwards along quasi-compact, quasi-separated morphisms preserve quasi-coherence. \smallskip

Now we show that $\rm{R}^i f_*\F$ is almost coherent. Note that it is sufficient to show that $\rm{R}^i f_*\F$ is almost finitely presented as $Y$ is a coherent scheme (this follows from Lemma \ref{sch-coh-acoh} and Lemma \ref{sheaf-acoh-afp}). We choose some finitely generated ideal $\m_0\subset \m$ and another finitely generated ideal $\m_1\subset \m$ such that $\m_0 \subset \m_1^2$. Then we use Corollary~\ref{cor:approximate-afpr} to find a finitely presented $\O_X$-module $\G$ and a morphism 
\[
\varphi\colon \G \to \F
\]
such that $\ker(\varphi)$ and $\coker(\varphi)$ are annihilated by $\m_1$. We define $\O_X$-modules \[
\K\coloneqq \ker \varphi, \ \M\coloneqq \Imm \varphi \text{ and } \sQ\coloneqq \coker \varphi,
\]
so we have two short exact sequences
\begin{align*}
0 \to \K \to \G \to \M \to 0, \\
0 \to \M \to \F \to \sQ \to 0
\end{align*}
with sheaves $\K$ and $\sQ$ killed by $\m_1$. This easily shows that the natural homomorphisms
\[
\rm{R}^if_*(\varphi)\colon \rm{R}^if_*\G \to \rm{R}^if_*\F
\]
have kernels and cokernels annihilated by $\m_1^2$. Since $\m_0 \subset \m_1^2$ we see that $\m_0(\ker \rm{R}^if_*(\varphi))=0$ and $\m_0(\coker \rm{R}^if_*(\varphi))=0$.  Moreover, we know that $\rm{R}^if_*\G$ is a finitely presented $\O_Y$-module by Theorem \ref{proper-mapping} ($\G$ is a coherent $\O_X$-module since $X$ is a coherent scheme). Therefore, we use Corollary~\ref{cor:approximate-afpr} to conclude that $\rm{R}^i f_*\F$ is an almost finitely presented $\O_Y$-module for any $i\geq 0$. .
\end{proof}

The next goal is to prove a version of the Almost Proper Mapping Theorem for nice formal schemes. But before doing this, we need to establish a slightly more precise version of the usual Proper Mapping Theorem for formal schemes than the one in \cite{FujKato}. 

\begin{thm}[Proper Mapping Theorem]\label{thm:proper-mapping-formal} Let $R$ be as in Set-up~\ref{set-up3}, $A$ a topologically finitely presented $R$-algebra, $\mf\colon \X \to \Spf A$ a topologically finitely presented, proper morphism, and $\F$ a coherent $\O_\X$-module. Then $\rm{H}^i(\X, \F)$ is a coherent $A$-module and the natural morphism
\[
\rm{H}^i(\X, \F)^{\Updelta} \to \rm{R}^i\mf_*\left(\F\right)
\]
is an isomorphism for any $i\geq 0$.
\end{thm}
\begin{proof}
    First, we use \cite[Theorem I.11.1.2]{FujKato} to conclude that $\bf{R}\mf_*\F\in \bf{D}^+_{coh}(\Spf A)$. Therefore, Theorem~\ref{equiv-derived-formal} implies that  $M\coloneqq \bf{R}\Gamma(\Spf A, \bf{R}\mf_*\F)$ lies in $\bf{D}^+_{acoh}(A)$, and 
    \[
    M^{L\Updelta} \simeq \bf{R}\mf_*\F.
    \]
    Moreover, Lemma~\ref{updelta-derived} implies that the natural map
    \[
    \rm{H}^i(\X, \F)^\Updelta \simeq \rm{H}^i(M)^{\Updelta}\to \rm{R}^i\mf_*\F
    \]
    is an isomorphism. Finally, we conclude that 
    \[
    \rm{H}^i(\X, \F)\simeq \rm{H}^0\left(\X, \rm{H}^i(\X, \F)^\Updelta\right) \simeq \rm{H}^0(\X, \rm{R}^i\mf_*\F)
    \]
    must be coherent because $\rm{R}^i\mf_*\F$ is coherent. 
\end{proof}

\begin{thm}[Almost Proper Mapping Theorem]\label{almost-proper-mapping-formal} Let $\Y$ be a topologically finitely presented formal $R$-scheme for $R$ as in Setup~\ref{set-up3}. And let $\mf\colon \X \to \Y$ be a proper, topologically finitely presented morphism. Then 
\begin{itemize}\itemsep0.5em
\item The functor $\mathbf{R}\mf_*$ sends $\mathbf{D}^*_{qc, acoh}(\X)$ to $\mathbf{D}^*_{qc, acoh}(\Y)$ for any $*\in \{\text{`` ''}, +, -, b\}$. 
\item The functor $\mathbf{R}\mf_*$ sends $\mathbf{D}^*_{acoh}(\X)^a$ to $\mathbf{D}^*_{acoh}(\Y)^a$ for any $*\in \{\text{`` ''}, +, -, b\}$. 
\item The functor $\mathbf{R}\mf_*$ sends $\mathbf{D}^+_{acoh}(\X)$ to $\mathbf{D}^+_{acoh}(\Y)$. 
\item If $Y_0\coloneqq \Y\times_{\Spf R} (\Spec R/I)$ has finite Krull dimension, then $\mathbf{R}\mf_*$ sends $\mathbf{D}^*_{acoh}(\X)$  to $\mathbf{D}^*_{acoh}(\Y)$ for any $*\in \{\text{`` ''}, +, -, b\}$.
\end{itemize}
Moreover, if $\Y=\Spf A$ is an affine scheme and $\F$ is an adically quasi-coherent, almost coherent $\O_\X$-module, then $\rm{H}^n(\X, \F)$ is almost coherent over $A$, and the natural map $\rm{H}^n(\X, \F)^{\Updelta} \to \rm{R}^n\mf_*\F$ is an isomorphism of $\O_\Y$-modules for $n\geq 0$.
\end{thm}

\begin{lemma}\label{finite-coh-dimension-2} Let $\Y$ be a quasi-compact adic formal $R$-scheme, and let $\mf\colon \X \to \Y$ be a topologically finite type, quasi-separated morphism. Suppose that the reduction $Y_0=\Y\times_{\Spf R} (\Spec R/I)$ (or equivalently the ``special fiber'' $\ov{\Y}=\Y\times_{\Spf R} \Spec R/\rm{Rad}(I)$) is of finite Krull dimension. Then $\X$ has finite Krull dimension, and $\mf_*$ is of finite cohomological dimension on $\mathbf{Mod}_\X$.
\end{lemma}
\begin{proof} The proof is identical to that of Lemma~\ref{finite-coh-dimension-sch} once we notice that the underlying topological spaces of $\Y$, $Y_0$, and $\ov{\Y}$ are canonically identified. 
\end{proof}

Also, before starting the proof of Theorem~\ref{almost-proper-mapping-formal}, we need to establish the following preliminary lemma: 

\begin{lemma}\label{useful} Let $\mf\colon \X \to \Y=\Spf A$ be a morphism as in Theorem~\ref{almost-proper-mapping-formal} with affine $\Y$, and let $\F\in \mathbf{Mod}_\X$ be an adically quasi-coherent, almost coherent sheaf. Then  $\rm{R}^q\mf_*\F$ is an adically quasi-coherent, almost coherent $\O_\Y$-module if
\begin{enumerate}[label=\textbf{(\arabic*)}]
\item\label{hyp-1} the $A$-module $\rm{H}^q(\X, \F)$ is almost coherent for any $q\geq 0$,
\item for any $g\in A$ with $\sU =\Spf A_{\{g\}}$, the canonical map 
\[
\rm{H}^q(\X, \F) \otimes_A A_{\{g\}} \to \rm{H}^q(\X_{\sU}, \F),
\]
is an isomorphism for any $q\geq 0$.
\end{enumerate}
\end{lemma}
\begin{proof}
Consider an $A$-module $M \coloneqq \rm{H}^q(\X,\F)$ that is almost coherent by our assumption. So Lemma~\ref{completion-finitely-generated} guarantees that $M$ is $I$-adically complete, and so $M^{\Updelta}$ is an adically quasi-coherent, almost coherent $\O_\X$-module. Now note that $\rm{R}^q\mf_*\F$ is the sheafification of the presheaf 
\[
\sU \mapsto \rm{H}^q(\X_{\sU}, \F)
\]
Thus, there is a canonical map $M \to \rm{H}^0(\Y, \rm{R}^q\mf_*\F)$ that induces a morphism
\[
M^{\Updelta} \to \rm{R}^q\mf_*\F.
\]
The second assumption together with Lemma~\ref{trivial-base-change} and Lemma~\ref{completion-finitely-generated} ensure that this map is an isomorphism on stalks (as the sheafification process preserves stalks). Therefore, $M^{\Updelta} \to \rm{R}^q\mf_*\F$ is an isomorphism of $\O_\X$-modules. In particular, $\rm{R}^q\mf_*\F$ is adically quasi-coherent and almost coherent.
\end{proof}

\begin{proof}[Proof of Theorem~\ref{almost-proper-mapping-formal}] We use the same reductions as in the proof of Theorem~\ref{almost-proper-mapping} to reduce to the situation of an adically quasi-coherent, almost coherent $\O_\X$-module $\F$. Moreover, the statement is local on $\Y$, so we can assume that $\Y=\Spf A$ is affine. \smallskip

Now we show that both conditions in Lemma~\ref{useful} are satisfied in our situation. \smallskip

{\it Step $1$: $\rm{H}^q(\X,\F)$ is almost coherent for every $q\geq 0$.} Fix a finitely generated ideal $\m_0 \subset \m$ and another finitely generated ideal $\m_1 \subset \m$ such that $\m_0 \subset \m_1^2$. \smallskip

Theorem~\ref{thm:approximate-formal} guarantees that there is a coherent $\O_\X$-module $\G_{\m_1}$ and a morphism $\phi_{\m_1} \colon \G_{\m_1} \to \F$ such that its kernel and cokernel are annihilated by $\m_1$. Then it is easy to see that the natural morphism
\[
\rm{H}^q(\X, \G_{\m_1})\to \rm{H}^q(\X, \F)
\]
has kernel annihilated by $\m_1^2$ and cokernel annihilated by $\m_1$. In particular, both the kernel and co-kernel are annihilated by $\m_0$. Since $\m_0$ was an arbitrary finitely generated sub-ideal of $\m$, it suffices to show that $\rm{H}^q(\X, \G_{\m_1})$ are coherent $A$-modules for any $q\geq 0$. This follows from Theorem~\ref{thm:proper-mapping-formal}. \smallskip

{\it Step $2$: canonical maps $\rm{H}^q(\X, \F) \otimes_A A_{\{g\}} \to \rm{H}^q(\X_{\sU}, \F)$ are isomorphisms for any $g\in A$, $q\geq 0$, and $\sU=\Spf A_{\{g\}}$.} Lemma~\ref{lemma:global-FP-approximation} guarantees that $\F$ admits an FP-approximation $\phi \colon \G \to \F$. Using Lemma~\ref{ker-coker-3}, we get the short exact sequences of adically quasi-coherent sheaves
\[
0 \to \cal{K} \to \G \to \cal{M} \to 0,
\]
\[
0 \to \cal{M} \to \F \to \cal{Q} \to 0,
\]
where $\cal{K}$ and $\cal{Q}$ are annihilated by $I^{n+1}$ for some $n\geq 0$. So $\cal{K}$ and $\cal{Q}$ can be identified with quasi-coherent sheaves on $\X_n \coloneqq \X \times_{\Spf A} \Spec A/I^{n+1}$. Therefore, the natural morphisms
\[
\rm{H}^q(\X, \cal{K}) \otimes_A A_{\{g\}} \simeq \rm{H}^q(\X_n, \cal{K}) \otimes_{A/I^{n+1}} (A/I^{n+1})_{g} \to \rm{H}^q(\X_{\sU, n}, \cal{K}),
\]
\[
\rm{H}^q(\X, \cal{Q}) \otimes_A A_{\{g\}} \simeq \rm{H}^q(\X_n, \cal{Q}) \otimes_{A/I^{n+1}} (A/I^{n+1})_{g} \to \rm{H}^q(\X_{\sU, n}, \cal{Q}) 
\]
are isomorphisms for $q\geq 0$. The morphism
\begin{equation}\label{eqn:adically-quasi-coh}
\rm{H}^q(\X, \cal{G})\otimes_A A_{\{g\}} \to \rm{H}^q(\X_{\sU}, \cal{G})
\end{equation}
is an isomorphism by Theorem~\ref{thm:proper-mapping-formal}. Therefore, the five-lemma and $A$-flatness of $A_{\{g\}}$ imply that the morphism 
\[
\rm{H}^q(\X, \cal{M})\otimes_A A_{\{g\}} \to \rm{H}^q(\X_{\sU}, \cal{M})
\]
is an isomorphism for any $q\geq 0$ as well. Applying the five-lemma again (and $A$-flatness of $A_{\{g\}}$), we conclude that the morphism 
\[
\rm{H}^q(\X, \cal{F})\otimes_A A_{\{g\}} \to \rm{H}^q(\X_{\sU}, \cal{F})
\]
must be an isomorphism for any $q\geq 0$ as well.
\end{proof} 

\subsection{Characterization of quasi-coherent, almost coherent complexes}\label{derived}

The main goal of this Section is to show an almost analogue of \cite[\href{https://stacks.math.columbia.edu/tag/0CSI}{Tag 0CSI}]{stacks-project}. This gives a useful characterization of objects in $\bf{D}_{qc, acoh}^b(X)$ on a separated, finitely presented $R$-scheme for a universally coherent $R$. This will be crucially used in our proof of the almost version of the Formal GAGA Theorem (see Theorem~\ref{GAGA}).\smallskip

Our proof follows the proof of \cite[\href{https://stacks.math.columbia.edu/tag/0CSI}{Tag 0CSI}]{stacks-project} quite closely, but we need to make certain adjustments to make the arguments work in the almost coherent setting.

\begin{thm}\label{check-perfect-pn} Let $R$ be an universally coherent ring with an ideal $\m$ such that $\m^2=\m$ and $\widetilde{\m}\coloneqq \m\otimes_R \m$ is flat. Suppose that $\F\in \mathbf{D}_{qc}(\P^n_R)$ an element such that $\mathbf{R}\Hom_{\P^n}(\mathcal P, \F) \in \mathbf{D}^-_{acoh}(R)$ for $\mathcal P=\oplus_{i=0}^n \O(i)$. Then $\F\in \mathbf{D}^-_{qc,acoh}(\P^n_R)$.
\end{thm}
\begin{proof}
We follow the ideas of \cite[\href{https://stacks.math.columbia.edu/tag/0CSG}{Tag 0CSG}]{stacks-project}. Denote the dg algebra $\mathbf{R}\Hom_X(\mathcal P, \mathcal P)$ by $S$. A computation of cohomology groups of line bundles on $\P_R^n$ implies that $S$ is a ``discrete'' non-commutative algebra that is finite and flat over $R$. Now \cite[\href{https://stacks.math.columbia.edu/tag/0BQU}{Tag 0BQU}]{stacks-project}\footnote{Note that they have slightly different notations for $R$ and $S$} guarantees that the functor
\[
-\otimes_S^{\mathbf{L}} \mathcal P \colon \mathbf{D}(S) \to \mathbf{D}_{qc}(\P_R^n)
\]
is an equivalence of categories, and its quasi-inverse is given by
\[
\mathbf{R}\Hom(\mathcal P,-)\colon \mathbf{D}_{qc}(\P_R^n) \to \mathbf{D}(S)
\]
So if we define $M\coloneqq \mathbf{R}\Hom(\mathcal P,\F) \in \mathbf{D}(S)$, our assumption implies that the image of $M$ in $\mathbf{D}(R)$ lies $\mathbf{D}^-_{acoh}(R)$. \smallskip

Therefore, it suffices to show that, for any $N\in \bf{D}(S)$ such that its image in $\bf{D}(R)$ lies in $\bf{D}^-_{acoh}(R)$, we have $N\otimes^{\bf{L}}_S \mathcal{P}$ lies in $\mathbf{D}^-_{qc, acoh}(\P_R^n)$.  \smallskip

We use the convergence spectral sequence
\[
\rm{E}^{p,q}_2=\mathcal H^p(\rm{H}^q(N)\otimes^{\mathbf{L}}_S \mathcal P) \Rightarrow \mathcal H^{p+q}(N\otimes^{\mathbf{L}}_S \mathcal P)
\]
to conclude that it suffices to assume that $N$ is just an $S$-module. Now we fix a finitely generated ideal $\m_1\subset \m$ and a finitely generated ideal $\m_0\subset \m$ such that $\m_1\subset \m_0^2$. Then Lemma~\ref{extension-noncommutative} implies that there is a finitely presented right $S$-module $N'$ with a morphism $f\colon N'\to N$ such that $\ker f$ and $\coker f$ are annihilated by $\m_0$. Then the universal coherence of $R$ and \cite[\href{https://stacks.math.columbia.edu/tag/0CSF}{Tag 0CSF}]{stacks-project} imply that $N'\otimes^{\mathbf{L}}_S \mathcal P \in \mathbf{D}^{-}_{qc, coh}(\P_R^n)$. Now we note that the functor 
\[
-\otimes_S^{\mathbf{L}} \mathcal P \colon \mathbf{D}(S) \to \mathbf{D}_{qc}(\P_R^n)
\]
is $R$-linear, so the standard argument shows that the cone of the morphism
\[
f\otimes^{\mathbf{L}}_S \mathcal P\colon N'\otimes^{\mathbf{L}}_S \mathcal P \to N\otimes^{\mathbf{L}}_S \mathcal P
\]
has cohomology sheaves anihilated by $\m_1\subset \m_0^2$. Since $\m_1\subset \m$ was an arbitrary finitely generated ideal, Lemma~\ref{almost-finitely-presented} implies that $N\otimes^{\mathbf{L}}_S \mathcal P$ is in $\mathbf{D}^-_{qc, acoh}(\P_R^n)$ and finishes the proof.
\end{proof}

\begin{lemma}\label{check-perfect-prel} Let $R$ be a universally coherent ring, let $X$ be a finitely presented separated $R$-scheme, and let $K \in \mathbf{D}_{qc}(X)$. If $\mathbf{R}\Gamma(X,E\otimes^{\mathbf{L}}_{\O_X} K)$ is in $\mathbf{D}^-_{acoh}(R)$ for every $E \in \mathbf{D}^-_{coh}(X)$, then $K\in \mathbf{D}^-_{qc,acoh}(X)$.
\end{lemma}
\begin{proof}
We follow the proof of \cite[\href{https://stacks.math.columbia.edu/tag/0CSL}{Tag 0CSL}]{stacks-project}. First, we note that the condition that $K\in \mathbf{D}^-_{qc, acoh}(X)$ is local on $X$ because $X$ is quasi-compact. Therefore, we can prove it locally around each point $x$. We use \cite[\href{https://stacks.math.columbia.edu/tag/0CSJ}{Tag 0CSJ}]{stacks-project} to find 
\begin{itemize}\itemsep0.5em
\item An open subset $U\subset X$ containing $x$.
\item An open subset $V\subset \P^n_R$.
\item A closed subset $Z\subset X\times_R \P^n_R$ with a point $z\in Z$ lying over $x$
\item An object $E\in \mathbf{D}^-_{coh}(X\times_R \mathbf{P}^n_R)$.
\end{itemize}
with a lot of properties listed in the cited lemma. Even though the notation is pretty heavy, the only properties of these objects that we will use are that $x\in U$ and that
\[
\mathbf{R}q_*(\mathbf{L}p^*K\otimes^{\mathbf{L}} E)|_{V}=\mathbf{R}(U\to V)_*(K|_{U}).
\]
The last formula is proven in \cite[\href{https://stacks.math.columbia.edu/tag/0CSK}{Tag 0CSK}]{stacks-project} and we refer to this lemma for a discussion of the morphism $U\to V$ that turns out to be a finitely presented closed immersion. \smallskip

That being said, we note that the above argument shows that it is sufficient to show that $K|_{U}$ is almost coherent for each such $U$. Moreover, the formula $ \mathbf{R}q_*(\mathbf{L}p^*K\otimes^{\mathbf{L}} E)|_{V}=\mathbf{R}(U\to V)_*(K|_{U})$, the fact that $U\to V$ is a finitely presented closed immersion, and Lemma~\ref{extension-noncommutative} imply that it suffices to show that $\mathbf{R}(U\to V)_*(K|_{U})=\mathbf{R}q_*(\mathbf{L}p^*K\otimes^{\mathbf{L}} E)|_{V}$ lies in $\mathbf{D}^-_{qc, acoh}(V)$. In particular, it is enough to show that $\mathbf{R}q_*(\mathbf{L}p^*K\otimes^{\mathbf{L}} E) \in \mathbf{D}^-_{qc, acoh}(\P^n_R)$. \smallskip

Now we check this using Theorem~\ref{check-perfect-pn}. Namely, we define a sheaf $\mathcal{P}\coloneqq \bigoplus_{i=0}^n \O_{\mathbf P^n}(i)$ and observe that
\begin{align*}
\mathbf{R}\Hom_{\P^n}(\mathcal P, \mathbf{R}q_*(\mathbf{L}p^*K\otimes^{\mathbf{L}} E))&=\mathbf{R}\Gamma(\mathbf P^n,  \mathbf{R}q_*(\mathbf{L}p^*K\otimes^{\mathbf{L}} E) \otimes^{\mathbf{L}}_{\O_{\P^n}} \mathcal P^{\vee})\\
&=\mathbf{R}\Gamma(\mathbf P^n,  \mathbf{R}q_*(\mathbf{L}p^*K\otimes^{\mathbf{L}} E \otimes^{\mathbf{L}} \mathbf{L}q^*\mathcal P^{\vee}))\\
&=\mathbf{R}\Gamma(X\times_R \mathbf P^n_R,  \mathbf{L}p^*K\otimes^{\mathbf{L}} E \otimes^{\mathbf{L}} \mathbf{L}q^*\mathcal P^{\vee})  \\
&=\mathbf{R}\Gamma(X,  \mathbf{R}p_*(\mathbf{L}p^*K\otimes^{\mathbf{L}} E \otimes^{\mathbf{L}} \mathbf{L}q^*\mathcal P^{\vee})) \\
&=\mathbf{R}\Gamma(X,  K\otimes^{\mathbf{L}}_{\O_X}\mathbf{R}p_*( E \otimes^{\mathbf{L}} \mathbf{L}q^*\mathcal P^{\vee}))
\end{align*}
where the second and fifth equality come from the projection formula \cite[\href{https://stacks.math.columbia.edu/tag/08EU}{Tag 08EU}]{stacks-project}. Now we note that the Proper Mapping Theorem (see Theorem~\ref{proper-mapping}) implies that $\mathbf{R}p_*( E \otimes^{\mathbf{L}} \mathbf{L}q^*\mathcal P^{\vee})\in \mathbf{D}^-_{coh}(X)$. So our assumption on $K$ implies that 
\[
\mathbf{R}\Hom_{\P^n}(\mathcal P, \mathbf{R}q_*(\mathbf{L}p^*K\otimes^{\mathbf{L}} E))=\mathbf{R}\Gamma(X,  K\otimes^{\mathbf{L}}_{\O_X}\mathbf{R}p_*( E \otimes^{\mathbf{L}} \mathbf{L}q^*\mathcal P^{\vee})) \in \mathbf{D}^-_{acoh}(R)
\]
Now Theorem~\ref{check-perfect-pn} finishes the proof.
\end{proof}

\begin{thm}\label{check-perfect} Let $R$ be a universally coherent ring, let $X$ be a separated, finitely presented $R$-scheme, and let $\F \in \mathbf{D}^-_{qc}(X)$ be an object such that $\mathbf{R}\Hom_X(\mathcal P, \F) \in \mathbf{D}_{acoh}^-(R)$ for any $\mathcal P \in \rm{Perf}(X)$. Then $\F\in \mathbf{D}^-_{qc, acoh}(X)$. Analogously, if $\mathbf{R}\Hom_X(\mathcal P, \F) \in \mathbf{D}^b_{acoh}(R)$ for any $\mathcal P \in \rm{Perf}(X)$, then $\F\in \mathbf{D}^b_{qc, acoh}(X)$.
\end{thm}
\begin{proof}
Once we have have Lemma~\ref{check-perfect-prel} and the equality $\mathbf{R}\Hom_X(\mathcal P, \F)=\mathbf{R}\Gamma(X, \mathcal P^{\vee} \otimes_{\O_X}^{\mathbf L} \F)$, the first part of the Theorem is absolutely analogous to \cite[\href{https://stacks.math.columbia.edu/tag/0CSH}{Tag 0CSH}]{stacks-project}. The second part now follows directly from \cite[\href{https://stacks.math.columbia.edu/tag/09IS}{Tag 09IS}]{stacks-project} and \cite[Lemma 3.0.12]{BZNP}.
\end{proof}

\subsection{The GAGA Theorem}\label{GAGA-section}

The main goal of this section is to prove the formal GAGA Theorem for almost coherent sheaves. It roughly says that any adically quasi-coherent, almost coherent sheaf on a completion of a proper, finitely presented scheme admits an essentially unique algebraization, and the same holds for morphisms of those sheaves. \smallskip

We start by recalling the statement of the classical formal GAGA Theorem. We fix an $I$-adically complete noetherian ring $A$ and a proper $A$-scheme $X$. Then we consider the $I$-adic completion $\X$ as a formal scheme over $\Spf A$. It comes equipped with the natural morphism $c\colon \X \to X$ of locally ringed spaces that induces a functor 
\[
c^*\colon \mathbf{Coh}_{X} \to \mathbf{Coh}_{\X}
\]
The GAGA Theorem says that it is an equivalence of categories. Let us say a few words about the classical proof of this theorem. It consists of three essentially independent steps: the first is to show that the morphism $c$ is flat; the second is to show that the functor $c^*$ induces an isomorphism
\[
c^*\colon \rm{H}^i(X, \rm{F}) \to \rm{H}^i(\X, c^*\rm{F})
\]
for any $\rm{F} \in \mathbf{Coh}_X$ and any integer $i$. The last is to prove that any coherent sheaf $\G\in \mathbf{Coh}_{\bf{P}^N}$ admits a surjection of the form $\bigoplus_i \O(n_i)^{m_i} \to \G$. Though the first two steps generalize to our Setup, there is no chance of having an analogue of the last statement. The reason is easy: existence of such a surjection would automatically imply that the sheaf $\G$ is of finite type, however, almost coherent sheaves are usually not of finite type. \smallskip

This issue suggests that we should take another approach to GAGA Theorems recently developed by J.\,Hall in his paper \cite{JH}. The main advantage of this approach is that it first {\it constructs} a candidate for algebraization, and only then proves that this candidate indeed provides an algebraization. We adapt this strategy to our almost context. \smallskip

We start with the discussion of the GAGA functor in the almost world. In what follows, we assume that $R$ is a ring from Setup~\ref{set-up3}. We fix a finitely presented $R$-scheme $X$, and we consider its $I$-adic completion $\X$ that is a topologically finitely presented formal $R$-scheme. The formal scheme $\X$ comes equipped with the canonical morphism of locally ringed spaces
\[
c\colon \left(\X, \O_\X \right) \to \left(X, \O_X\right)
\]
that induces the pullback functor
\[
\mathbf{L}c^*\colon \mathbf{D}(X) \to \mathbf{D}(\X).
\]

We now want to check that this functor preserves quasi-coherent, almost coherent objects. This verification will be necessary even to formulate the GAGA statement. 

\begin{lemma}\label{preserve} Let $R$ be a ring as in Setup~\ref{set-up3}, $A$ a topologically finitely presented $R$-algebra, and $X$ a finitely presented $A$-scheme. Then the morphism $c$ is flat, and the funtor $c^*\colon \mathbf{Mod}_{X} \to \mathbf{Mod}_{\X}$ sends (quasi-coherent and) almost coherent sheaves to (adically quasi-coherent and) almost coherent sheaves. In particular, it induces functors
\[
\mathbf{L}c^*\colon \mathbf{D}_{qc, acoh}^*(X) \to \mathbf{D}_{qc, acoh}^*(\X)
\]
for any $*\in \{\text{`` ''}, +, -, b\}$.
\end{lemma}
\begin{proof}
The flatness assertion follows from \cite[Proposition I.1.4.7 (2)]{FujKato}. Flatness of $c$ implies that it suffices to show that $c^*(\rm{G})$ is adically quasi-coherent, almost coherent $\O_\X$-module for a quasi-coherent, almost coherent $\O_X$-module $\rm{G}$. This claim is Zariski-local on $X$. Thus we can assume that $X=\Spec A$ is affine, so $\rm{G}\simeq \widetilde{M}$ for some almost finitely presented $A$-module $M$. This case is done in Lemma~\ref{completion-affine}.
\end{proof}

\begin{thm}\label{GAGA} Let $R$ be a ring as in Setup~\ref{set-up3}, $A$ a topologically finitely presented $R$-algebra, and $X$ a finitely presented, proper $A$-scheme. Then the functor 
\[
\mathbf{L}c^*\colon \mathbf{D}_{qc, acoh}^*(X) \to \mathbf{D}_{qc, acoh}^*(\X)
\]
induces an equivalence of categories for $*\in  \{\text{`` ''}, +, -, b\}$. 
\end{thm}

\begin{cor}\label{cor:almost-GAGA} Let $R$, $A$ and $X$ be as in Theorem~\ref{GAGA}. Then the functor 
\[
\mathbf{L}c^*\colon \mathbf{D}_{acoh}^*(X)^a \to \mathbf{D}_{acoh}^*(\X)^a
\]
induces an equivalence of categories for $*\in  \{\text{`` ''}, +, -, b\}$. 
\end{cor}

\begin{cor}\label{cor:formal-functions-coh} Let $R$, $A$, and $X$ be as in Theorem~\ref{GAGA}, and let $K\in \bf{D}_{qc, acoh}(X)$. Then the natural map 
\[
\beta_{K}\colon \bf{R}\Gamma(X, K) \to \bf{R}\Gamma(\X, \bf{L}c^*K)
\]
is an isomorphism. Moreover, the map $\beta_K$ is an almost isomorphism for $K\in \bf{D}_{acoh}(X)$.
\end{cor}
\begin{proof}
Note that the case of $K\in \bf{D}_{acoh}(X)$ follows from the case of $K\in \bf{D}_{qc, acoh}(X)$ due to Lemma~\ref{derived-pull-projection} and Proposition~\ref{derived-pushforward}. So it suffices to prove for $K\in \bf{D}_{qc, acoh}(X)$. \smallskip

Now since we are allowed to replace $K$ with $K[i]$ for any integer $i$, it suffices to show that the map
\[
\rm{H}^0(\bf{R}\Gamma(X, K))\simeq\rm{Hom}_X(\O_X, K) \to \rm{Hom}_\X(\O_\X, \bf{L}c^*K)\simeq \rm{H}^0(\bf{R}\Gamma(\X, \bf{L}c^*K))
\]
is an isomorphism. This follows from Theorem~\ref{GAGA} and the observation that $\O_\X\simeq \bf{L}c^*\O_X$. 
\end{proof}

Our proof of Theorem~\ref{GAGA} will follow Jack Hall's proof of the GAGA Theorem very closely with some simplifications due to the flatness of the functor $c^*$. As he works entirely in the setting of the pseudo-coherent objects, and almost coherent sheaves may not be pseudo-coherent, we have repeat some arguments in our setting. \smallskip

Before we embark on the proof proof, we need to define the functor in the other direction. Recall that the morphism of locally ringed spaces $c$ defines the derived pushforward functor
\[
\mathbf{R}c_*\colon \mathbf{D}(\X) \to \mathbf{D}(X).
\]
This functor is $t$-exact as $c\colon \X \to X$ is topologically just a closed immersion. In particular, it preserves boundedness of complexes (in any direction). However, that functor usually does not preserve (almost) coherent objects as can be seen in the example of $\mathbf{R}c_*\O_\X=c_*\O_\X$. A way to fix it is to use the quasi-coherator functor
\[
\mathbf{R}Q_X\colon \mathbf{D}(X) \to \mathbf{D}_{qc}(X)
\]
that is defined as the right adjoint to the inclusion  $\iota\colon \mathbf{D}_{qc}(X) \to \mathbf{D}(X)$. It exists by \cite[\href{https://stacks.math.columbia.edu/tag/0CR0}{Tag 0CR0}]{stacks-project}. 
We define the functor
\[
\mathbf{R}c_{qc}\colon \mathbf{D}(\X) \to \mathbf{D}_{qc}(X)
\] 
as the composition $\mathbf{R}c_{qc}\coloneqq \mathbf{R}Q_X\circ \mathbf{R}c_*$. \smallskip

Combining the adjunctions $(\mathbf{L}c^*, \mathbf{R}c_*)$ and $(\iota, \mathbf{R}Q_X)$, we conclude that we have a pair of the adjoint functors:
\[
\adj{\mathbf{L}c^*}{\mathbf{D}_{qc}(X)}{\mathbf{D(\X)}}{\mathbf{R}c_{qc}}
\]
That gives us the unit and counit morphisms
\[
\eta\colon \rm{Id} \to \mathbf{R}c_{qc}\mathbf{L}c^* \text{  and  } \e\colon \mathbf{L}c^*\mathbf{R}c_{qc} \to \rm{Id}
\]  
For future reference, we also note that the above adjunction and the monoidal property of the functor $\mathbf{L}c^*$ define a projection morphism
\[
\pi_{\rm{G}, \F}\colon \rm{G}\otimes_{\O_X}^{\mathbf{L}} (\mathbf{R}c_{qc}\F) \to \mathbf{R}c_{qc}(\mathbf{L}c^*\rm{G}\otimes^{\mathbf{L}}_{\O_\X}\F)
\]
for any $\rm{G}\in \mathbf{D}_{qc}(X)$ and any $\F\in \mathbf{D}(\X)$. Before discussing the proof of Theorem~\ref{GAGA}, we need to establish some properties of these functors.

\begin{lemma}\label{finite-coh-dimension} Let $R$ be a ring as in Setup~\ref{set-up3}, $A$ a topologically finitely presented $R$-algebra, and $X$ a finitely presented $A$-scheme. Then there is an integer $N=N(X)$ such that $\mathbf{R}c_{qc}$ carries $\mathbf{D}^{\leq n}_{qc, acoh}(\X)$ to $\mathbf{D}^{\leq n+N}_{qc}(X)$ (resp. $\mathbf{D}^{[a, n]}_{qc, acoh}(\X)$ to $\mathbf{D}^{[a, n+N]}_{qc}(X)$) for any integer $n$. In particular, the natural map 
\[
\tau^{\geq a}\mathbf{R}c_{qc} \F \to \tau^{\geq a}(\mathbf{R}c_{qc} \tau^{\geq a-N} \F)
\]
is an isomorphism for any $\F \in \mathbf{D}_{qc, acoh}(\X)$ and any integer $a$.
\end{lemma}
\begin{proof}

We explain the proof that $\mathbf{R}c_{qc}$ carries $\mathbf{D}^{\leq n}_{qc, acoh}(\X)$ to $\mathbf{D}^{\leq n+N}_{qc}(X)$; the case of $\mathbf{D}^{[a, n]}_{qc, acoh}(\X)$ is similar. We fix an object $\F \in \mathbf{D}^{\leq n}_{qc, acoh}(\X)$. Then we note that $\bf{R}c_*\F = c_*\F$ since $c$ is topologically a closed immersion. Thus, \cite[\href{https://stacks.math.columbia.edu/tag/0CSA}{Tag 0CSA}]{stacks-project} implies\footnote{We note that the proof of \cite[\href{https://stacks.math.columbia.edu/tag/0CSA}{Tag 0CSA}]{stacks-project} works well with $a=-\infty$ as well.} that it suffices to show that 
\[
\mathrm{H}^i\big(\mathbf{R}\Gamma(U, \bf{R} c_*\F)\big) = \mathrm{H}^i\big(\mathbf{R}\Gamma(U, c_*\F)\big)=0
\]
for any open affine $U\subset X$ and any $i\geq n$. Therefore, we see that 
\[
\mathrm{H}^i\big(\mathbf{R}\Gamma(U, \bf{R} c_*\F)\big) = \mathrm{H}^i\big(\mathbf{R}\Gamma(\widehat{U}, \F)\big)=\mathrm{H}^i(\widehat{U}, \F|_{\widehat{U}}),
\]
and thus Lemma~\ref{derived-formal-schemes} implies that $\mathrm{H}^i(\widehat{U}, \F|_{\widehat{U}})=0$ for any $i\geq n$. This finishes the proof of the first claim in the lemma. \smallskip

The second claim of the lemma follows from the first claim and the distinguished triangle
\[
\tau^{\leq a-N-1} \F \to \F \to \tau^{\geq a-N} \F \to \tau^{\leq a-N-1} \F [1]
\]
Namely, we apply the exact functor $\mathbf{R}c_{qc}$ to this distinguished triangle to get that
\[
\mathbf{R}c_{qc}\left(\tau^{\leq a-N-1} \F\right) \to \mathbf{R}c_{qc}\F \to \mathbf{R}c_{qc}\left(\tau^{\geq a-N} \F\right) \to \mathbf{R}c_{qc} \left(\tau^{\leq a-N-1} \F [1]\right)
\]
is a distinguished triangle in $\mathbf{D}_{qc}(X)$ and that $\mathbf{R}c_{qc}(\tau^{\leq a-N-1} \F)\in \mathbf{D}^{\leq a-1}_{qc}(X)$. This implies that the map
\[
\tau^{\geq a}\mathbf{R}c_{qc}\F \to \tau^{\geq a}\mathbf{R}c_{qc}\left(\tau^{\geq a-N} \F\right)
\]
is an isomorphism.
\end{proof}

\begin{lemma}\label{trivial-cases} Let $X$ be as in Theorem~\ref{GAGA}, $\F\in \mathbf{D}^{-}_{qc, acoh}(\X)$ and $\rm{G}\in \mathbf{D}^-_{qc}(X)$. Suppose that for each $i$ there is $n_i$ such that $I^{n_i}\mathcal{H}^i(\F)=0$ and $I^{n_i}\mathcal{H}^i(\rm{G})=0$. Then the natural morphisms $\eta_{\rm{G}}$ and $\e_{\F}$ are isomorphisms.
\end{lemma}
\begin{proof}
We prove the claim only for $\F$ as the other claim is similar. \smallskip

{\it Reduction to the case when $\F\in \mathbf{D}^b_{qc,acoh}(\X)$}: First, we note that it suffices to show that the natural map
\[
\tau^{\geq a} \F \to \tau^{\geq a} \mathbf{L}c^*\mathbf{R}c_{qc}\F
\]
is an isomorphism for any integer $a$. Moreover, we note that  $t$-exactness of $\mathbf{L}c^*$ and Lemma~\ref{finite-coh-dimension} imply that there is an integer $N$ such that the natural map $\tau^{\geq a} \mathbf{L}c^*\mathbf{R}c_{qc}\F \to \tau^{\geq a} \mathbf{L}c^*\mathbf{R}c_{qc}\tau^{\geq a-N}\F$ is an isomorphism for any integer $a$. In particular, we have a commutative diagram
\[
\begin{tikzcd}
\tau^{\geq a-N} \F \arrow{d} \arrow{r} & \mathbf{L}c^*\mathbf{R}c_{qc}(\tau^{\geq a-N}\F) \arrow{d} \\
\tau^{\geq a} \F \arrow{r} &\tau^{\geq a} \mathbf{L}c^*\mathbf{R}c_{qc}\F \simeq \tau^{\geq a} \mathbf{L}c^*\mathbf{R}c_{qc}\tau^{\geq a-N}\F  
\end{tikzcd}
\]
where the vertical maps induce isomorphisms in degree $\geq a$. Therefore, it suffices to prove the claim for $\tau^{\geq a-N}\F$. So we may and do assume that $\F$ is bounded. \smallskip

{\it Proof for a bounded $\F$}: The case of a bounded $\F$ easily reduces to the case of an adically quasi-coherent, almost coherent $\O_\X$-module $\F$ concentrated in degree $0$. In that situation, we have $I^{k+1}\F=0$ for some $k$. This implies that $\F=i_{k,*}\F_k=\mathbf{R}i_{k,*}\F_k$ for the closed immersion $i_k\colon X_k \to \X$. Now it is straightforward to see that the canonical map
\[
\mathbf{R}i_{k,*}\F_k \to \mathbf{L}c^*\mathbf{R}c_{qc} (\mathbf{R}i_{k,*}\F_k)
\]
is an isomorphism. The key is flatness of $c$ and the observation that $\mathbf{R}c_*(\mathbf{R}i_{k,*}\F_k)$ is already quasi-coherent, so the quasi-coherator does nothing in this case.
\end{proof}


\begin{lemma}\label{perfect-ok} If $\rm{G} \in \mathbf{D}_{qc}(X)$ and $\F \in \mathbf{D}(\X)$, then the natural projection morphism
\[
\pi_{\rm{G},\F} \colon \rm{G}\otimes^{\mathbf{L}}_{\O_X} \mathbf{R}c_{qc}\F \to \mathbf{R}c_{qc}(\mathbf{L}c^* \rm{G} \otimes^{\mathbf L}_{\O_\X} \F)
\]
is an isomorphism if $G$ is perfect.
\end{lemma}
\begin{proof}
\cite[Lemma 4.3]{JH}.
\end{proof}

Now we come to the key input ingredient. Although $\mathbf{R}c_{qc}$ is quite abstract and difficult to compute in practice, it turns out that the Almost Proper Mapping Theorem allows us to check that this functor sends $\mathbf{D}^-_{qc, acoh}(\X)$ to $\mathbf{D}^-_{qc, acoh}(X)$. This would give us a candidate for an algebraization. 

\begin{lemma} Let $R$ be a ring as in Setup~\ref{set-up3}, $A$ a topologically finitely presented $R$-algebra, and $X$ a finitely presented, proper $A$-scheme. Then $\mathbf{R}c_{qc}$ sends $\mathbf{D}^*_{qc, acoh}(\X)$ to $\mathbf{D}^*_{qc, acoh}(X)$ for $*\in \{-, \ b\}$.
\end{lemma}
\begin{proof}
We prove only the bounded above case as the other one follows from this using Lemma~\ref{finite-coh-dimension}. We pick any $\F\in \mathbf{D}^-_{qc, acoh}(\X)$ and use Theorem~\ref{check-perfect} to say that it is sufficient to show that $\mathbf{R}\Hom_X(\rm{P}, \mathbf{R}c_* \F)\in \mathbf{D}^-_{acoh}(R)$ for any perfect complex $\rm{P}\in \rm{Perf}(X)$. For this, we consider the following sequence of isomorphisms
\begin{align*}
\mathbf{R}\Hom_X(\rm{P}, \mathbf{R}c_{qc}\F)&=\mathbf{R}\Hom_\X(\mathbf{L}c^*\rm{P}, \F) \\
&=\mathbf{R}\Hom_\X(\O_\X, (\mathbf{L}c^*\rm{P})^{\vee} \otimes^{\mathbf{L}}_{\O_{\X}} \F) \\&=\mathbf{R}\Gamma(\X, (\mathbf{L}c^*\rm{P})^{\vee} \otimes^{\mathbf{L}}_{\O_{\X}} \F).
\end{align*}
Then we note that $\cal{P}\coloneqq (\mathbf{L}c^*\rm{P})^{\vee}$ is a perfect complex of $\O_\X$-modules, and so $\cal{P}\otimes^{\bf{L}}_{\O_\X} \F$ lies in $\mathbf{D}^-_{qc, acoh}(\X)$. Thus, $\mathbf{R}\Gamma(\X, \cal{P} \otimes^{\mathbf{L}}_{\O_{\X}} \F)$ lies in $\mathbf{D}^-_{acoh}(R)$ due to the Almost Propper Mapping Theorem (see Theorem~\ref{almost-proper-mapping-formal}).
\end{proof}

Finally, we are ready to give a proof of the GAGA Theorem.

\begin{proof}[Proof of Theorem~\ref{GAGA}]

{\it Claim 0:} It suffices to show the theorem for $*=-$, i.e. for bounded above derived categories. Indeed, flatness of $c^*$ implies that $\mathbf{L}c^*$ preserve boundedness (resp. boundedness above, resp. boundedness below), so it suffices to show that the natural morphisms
\[
\eta_{\rm{G}} \colon \rm{G} \to \bf{R}c_{qc} \bf{L}c^* \rm{G}
\]
\[
\e_{\F}\colon \bf{L}c^* \bf{R}c_{qc} \F \to \F
\]
are isomorphisms for any $\rm{G}\in \bf{D}_{qc, acoh}(X)$ and $\F \in \bf{D}_{qc, acoh}(\X)$. \smallskip{}

We fix $N$ as in Lemma~\ref{finite-coh-dimension}. Then flatness of $c^*$ and Lemma~\ref{finite-coh-dimension} guarantee that 
\[
\bf{R}c_{qc} \bf{L}c^* \tau^{\geq a} \rm{G} \in \bf{D}^{[a, \infty]}(X)
\]
\[
\bf{L}c^* \bf{R}c_{qc} \tau^{\geq a} \F \in \bf{D}^{[a, \infty]}(\X).
\]
Therefore, we see that $\eta_{\rm{G}}$ is an isomorphism on $\mathcal{H}^{i}$ for $i<a$ if and only if the same holds for $\eta_{\tau^{\leq a-1}\rm{G}}$. Since $a$ was arbitrary, we conclude that it suffices to show that $\eta_{\rm{G}}$ is an isomorphism for $\rm{G} \in \bf{D}_{qc, acoh}^-(X)$. Similar argument shows that it suffices to show that $\e_{\F}$ is an isomorphism for $\F\in \bf{D}_{qc, acoh}^-(\X)$. So it suffices to prove the theorem for $*=-$. \smallskip

Before we formulate the next claim, we need to use the so-called ``approximation by perfect complexes'' \cite[\href{https://stacks.math.columbia.edu/tag/08EL}{Tag 08EL}]{stacks-project} to find some $\rm{P}\in \rm{Perf}(X)$ such that $\t^{\geq 0}\rm{P}\simeq \O_X/I\simeq \O_{X_0}$ and whose support is equal to $X_0$. We note that it implies that all cohomology sheaves $\mathcal{H}^i(\rm{P})$ are killed by some power of $I$. We also denote its (derived) pullback by $\mathcal P\coloneqq \mathbf{L}c^*\rm{P}$.\\

{\it Claim 1:} If $\rm{G} \in \mathbf{D}^-_{qc,acoh}(X)$ such that $\rm{G}\otimes^{\mathbf{L}}_{\O_X}\rm{P}\simeq 0$, then $\rm{G}\simeq 0$. Similarly, if $\F\in \mathbf{D}^{-}_{qc, acoh}(\X)$ such that $\F\otimes^{\mathbf{L}}_{\O_\X} \mathcal P\simeq 0$, then $\F\simeq 0$. \\

We choose a maximal $m$ (assuming that $\rm{G}\not\simeq
 0$) such that $\mathcal H^m(\rm{G})\neq 0$. Then we see that $\mathcal H^m(\rm{G}\otimes^{\mathbf{L}}_{\O_X}\rm{P})\simeq \mathcal H^m(\rm{G})\otimes_{\O_X}\O_{X_0}=\mathcal H^m(\rm{G})/I$. We have  $(\mathcal{H}^m(\rm{G})/I)(U)=\mathcal H^m(\rm{G})(U)/I \simeq 0$ on any open affine $U$. So Nakayama's Lemma~\ref{Nakayama} implies that $\mathcal H^m(\rm{G})(U) \simeq 0$ for any such $U$. This contradicts the choice of $m$. The proof in the formal setup is the same once we notice that $\mathcal{H}^0(\mathcal {P})=\O_{\X}/I$. \\

{\it Claim 2:} The map $\eta_{\rm{G}} \colon \rm{G} \to \mathbf{R}c_{qc}\mathbf{L}c^*\rm{G}$ is an isomorphism for any $\rm{G}\in \mathbf{D}^-_{qc, acoh}(X)$. \\

Claim $1$ implies that it is sufficient to show that the map
\begin{equation}\label{10}
\e_{\rm{G}}\otimes_{\O_X}^{\mathbf{L}} \rm{P} \colon \rm{G}\otimes_{\O_X}^{\mathbf{L}} \rm{P} \to 
\mathbf{R}c_{qc}\mathbf{L}c^*\rm{G} \otimes_{\O_X}^{\mathbf{L}} \rm{P}
\end{equation}
is an isomorphism. Recall that the cohomology sheaves of $\rm{P}$ are killed by some power of $I$. This property passes to $\rm{G}\otimes^{\mathbf{L}}_{\O_X} \rm{P}$, so we can use Lemma~\ref{trivial-cases} to get that the map 
\[
\e_{\rm{G}\otimes_{\O_X}^{\mathbf{L}} \rm{P}} \colon\rm{G}\otimes^{\mathbf{L}}_{\O_X} \rm{P}\to \mathbf{R}c_{qc}\left(\mathbf{L}c^*\left(\rm{G} \otimes_{\O_X}^{\mathbf{L}}P\right)\right)
\] 
is an isomorphism. Now comes the key: we fit  the morphism $\e_{\rm{G}\otimes_{\O_X}^{\mathbf{L}} \rm{P}}$ into the following commutative triangle:
\[
\begin{tikzcd}[column sep=10ex]
\rm{G}\otimes_{\O_X}^{\mathbf{L}}\rm{P}\arrow{r}{\e_{\rm{G}}\otimes_{\O_X}^{\mathbf{L}} \rm{P}} \arrow[d, swap, "\e_{\rm{G}\otimes_{\O_X}^{\mathbf{L}} \rm{P}}"] & \mathbf{R}c_{qc}\mathbf{L}c^*\rm{G} \otimes_{\O_X}^{\mathbf{L}} \rm{P}  \arrow{d}{\pi_{\rm{P}, \mathbf{L}c^*\rm{G}}} \\
\mathbf{R}c_{qc}(\mathbf{L}c^*(\rm{G} \otimes_{\O_X}^{\mathbf{L}}P)) \arrow{r}{\sim}& \mathbf{R}c_{qc}(\mathbf{L}c^*\rm{G}\otimes_{\O_\X}^{\mathbf{L}}\mathbf{L}c^*\rm{P})
\end{tikzcd}
\]
where the bottom horizontal arrow is the isomorphism map induced by the monoidal structure on $\mathbf{L}c^*$. Moreover, we have already established that the left vertical arrow is an isomorphism, and right vertical arrow is an isomorphism due to Lemma~\ref{perfect-ok}. That shows that the top horizontal must be also an isomorphism. \\

{\it Claim 3:} The map $\e_{\F} \colon \mathbf{L}c^*\mathbf{R}c_{qc}\F \to \F$ is an isomorphism for any $\F\in \mathbf{D}^-_{qc, acoh}(\X)$. \\

We use Claim $1$ again to say that it is sufficient to show that the map
\[
\e_{\F}\otimes_{\O_\X}^\mathbf{L} \mathbf{L}c^*\rm{P} \colon \mathbf{L}c^*\mathbf{R}c_{qc}\F \otimes_{\O_\X}^\mathbf{L} \mathbf{L}c^*\rm{P} \to \F \otimes_{\O_\X}^\mathbf{L} \mathbf{L}c^*\rm{P}
\]
is an isomorphism. But that map fits into the commutative diagram:

\[
\begin{tikzcd}[column sep=10ex, row sep=7ex]
 \mathbf{L}c^*\mathbf{R}c_{qc}\F \otimes_{\O_\X}^\mathbf{L} \mathbf{L}c^*\rm{P} \arrow{r}{\e_{\F}\otimes_{\O_\X}^\mathbf{L} \mathbf{L}c^*\rm{P}}\arrow{d}{\wr}  & \F \otimes_{\O_\X}^\mathbf{L} \mathbf{L}c^*\rm{P} \\
 \mathbf{L}c^*(\mathbf{R}c_{qc}\F \otimes_{\O_\X}^\mathbf{L}  \rm{P})  \arrow{r}{\mathbf{L}c^*(\pi_{\rm{P},\F})} & \mathbf{L}c^*\mathbf{R}c_{qc}(\F \otimes_{\O_\X}^\mathbf{L}  \mathbf{L}c^*\rm{P}) \arrow[u, swap, "\e_{\F\otimes_{\O_\X}^\mathbf{L} \mathbf{L}c^*\rm{P}}"]
\end{tikzcd}
\]
where the left vertical morphism is the canonical isomorphism induced by the monoidal structure on $\mathbf{L}c^*$, the bottom morphism is an isomorphism by Lemma~\ref{perfect-ok}, and the right vertical morphism is an isomorphism by Lemma~\ref{trivial-cases}. This implies that the top horizontal morphism is an isomorphism as well. This finishes the proof.
\end{proof}


\subsection{The formal function Theorem}

In this Section, we prove the Formal Function Theorem for almost coherent sheaves as a consequence of the Formal GAGA Theorem established in the previous section.  \smallskip

For the rest of the section, we fix a ring $R$ as in Setup~\ref{set-up3} and a finitely presented or a topologically finitely presented $R$-algebra $A$.

\begin{rmk}\label{rmk:A-good} Both $A$ and $\wdh{A}$ are topologically universally adhesive by \cite[Proposition 0.8.5.19]{FujKato}, and they are (topologically universally) coherent by \cite[Proposition 0.8.5.23]{FujKato}. 
\end{rmk}

For the next definition, we fix a finitely presented $A$-scheme $X$ and an $\O_X$-module $\F$.

\begin{defn}\label{defn:FP-approx-natural-topology} 
The {\it natural $I$-filtration} $\rm{F}^\bullet \rm{H}^i(X, \F)$ on $\rm{H}^i(X, \F)$ is defined via the formula 
\[
\rm{F}^n \rm{H}^i(X, \F) \coloneqq \rm{Im}\left(\rm{H}^i(X, I^n\F) \to \rm{H}^i(X, \F) \right)
\]

The {\it natural $I$-topology} on $\rm{H}^i(X, \F)$ is the topology induced by the natural $I$-filtration. 
\end{defn}

\begin{lemma}\label{lemma:FP-approx-subsheaf-top} Let $X$ be a finitely presented $A$-scheme, $\F$ a quasi-coherent almost finitely generated $\O_X$-module, and $\G \subset \F$ be a quasi-coherent $\O_X$-submodule of $\F$. Then, for any $n$, there is an $m$ such that $I^m\F\cap \G\subset I^n\G$. 
\end{lemma}
\begin{proof}
It suffices to assume that $X$ is affine, in which case it follows from Lemma~\ref{Artin-Rees}. 
\end{proof}

\begin{lemma}\label{lemma:different-top} Let $X$ be a finitely presented $A$-scheme, $\F$ and $\G$ quasi-coherent almost finitely generated $\O_X$-modules, and $\varphi\colon \G \to \F$ an $\O_X$-linear homomorphism such that $\ker(\varphi)$ and $\coker(\varphi)$ are annihilated by $I^c$ for some integer $c$. Then, for every $i\geq 0$, the natural $I$-topology on $\rm{H}^i(X, \F)$ coincides with the topology induced by the filtration 
\[
    \rm{Fil}_\G^n\rm{H}^i(X, \F)=\rm{Im}(\rm{H}^i(X, I^n\G) \to \rm{H}^i(X, \F)).
\]
\end{lemma}
\begin{proof}
Consider the short exact sequences
\[
0 \to \cal{K} \to \G \to \cal{H} \to 0,
\]
\[
0 \to \cal{H} \to \F \to \cal{Q} \to 0,
\]
where $\cal{K}$ and $\cal{Q}$ are annihilated by $I^c$. The first short exact sequence inducea the following short exact sequence
\[
0 \to \cal{K}\cap I^m\G \to I^m\G \to I^m\cal{H} \to 0
\]
for any $m\geq 0$. Lemma~\ref{lemma:FP-approx-subsheaf-top} implies that $\cal{K}\cap I^m\G \subset I^c\cal{K}=0$ for $m\gg 0$. Therefore, the natural map $I^m\G \to I^m\cal{H}$ is an isomorphism for $m\gg 0$. Note that $\cal{H}$ is almost finitely generated and quasi-coherent, so we can replace $\G$ with $\cal{H}$ to assume that $\varphi$ is injective. \smallskip

Clearly, $\rm{Fil}_\G^k\rm{H}^i(X, \F) \subset \rm{F}^k\rm{H}^i(X, \F)$ for every $k$. So it suffices to show that, for any $k$, there is $m$ such that $\rm{F}^m\rm{H}^i(X, \F) \subset \rm{Fil}_\G^k\rm{H}^i(X, \F)$. We consider the short exact sequence
\[
0 \to \G \cap I^m\F \to I^m\F \to I^m\cal{Q} \to 0.
\]
If $m\geq c$ we get that $\G \cap I^m \F = I^m\F$ because $I^c\cal{Q}\simeq 0$. Now we use Lemma~\ref{lemma:FP-approx-subsheaf-top} to conclude there is $m\geq c$ such that 
\[
I^m\F = \G\cap I^m\F \subset I^k\G
\]
Therefore, $\rm{F}^m\rm{H}^i(X, \F) \subset \rm{Fil}_\G^k\rm{H}^i(X, \F)$. 
\end{proof}

\begin{lemma}\label{lemma:FP-approx-enough} Let $X$ be a finitely presented $A$-scheme, $\F$ and $\G$ quasi-coherent almost finitely generated $\O_X$-modules, and $\varphi\colon \G \to \F$ an $\O_X$-linear homomorphism such that $\ker(\varphi)$ and $\coker(\varphi)$ are annihilated by $I^c$ for some integer $c$. Suppose that the natural $I$-topology on $\rm{H}^i(X, \G)$ is the $I$-adic topology. Then the same holds for $\rm{H}^i(X, \F)$. 
\end{lemma}
\begin{proof}
Clearly, $I^n\rm{H}^i(X, \F) \subset \rm{F}^n\rm{H}^i(X, \F)$. So it suffices to show that, for every $n$, there is an $m$ such that $\rm{F}^m \rm{H}^i(X, \F) \subset I^n\rm{H}^i(X, \F)$. \smallskip

The assumption that the natural $I$-topology on $\rm{H}^i(X, \G)$ coincides with the $I$-adic topology guarantees that $\rm{F}^k\rm{H}^i(X, \G) \subset I^n\rm{H}^i(X, \G)$ for large enough $k$. Pick such $k$. Lemma~\ref{lemma:different-top} implies that 
\[
\rm{F}^m\rm{H}^i(X, \F) \subset \rm{Im}(\rm{H}^i(X, I^k\G) \to \rm{H}^i(X, \F))
\]
for large enough $m$. So we get that 
\[
\rm{F}^m\rm{H}^i(X, \F) \subset \rm{Im}\left(\rm{H}^i(X, I^k\G) \to \rm{H}^i\left(X, \F\right)\right) \subset \rm{Im}\left(I^n\rm{H}^i\left(X, \G\right) \to \rm{H}^i\left(X, \F\right)\right) \subset I^n\rm{H}^i\left(X, \F\right)
\]
for a large enough $m$.
\end{proof}

\begin{thm}\label{thm:FP-approx-topology-cohomology} Let $X$ be a proper, finitely presented $A$-scheme, and $\F$ a quasi-coherent, almost coherent $\O_X$-module. Then the natural $I$-topology on $\rm{H}^i(X, \F)$ coincides with the $I$-adic topology for any $i \geq 0$.
\end{thm}
\begin{proof}
   Lemma~\ref{lemma:approximate-finite-type} guarantees that there is a finitely presented $\O_X$-module $\G$ and a morphism $\varphi\colon \G \to \F$ such that $I^c(\ker \varphi)=0$ and $I^c(\coker \varphi)=0$. Lemma~\ref{lemma:FP-approx-enough} then ensures that it suffices to prove the claim for $\G$. In this case, the claim follows \cite[Proposition I.8.5.2 and Lemma 0.7.4.3]{FujKato} and Remark~\ref{rmk:A-good}.
\end{proof}

Now we consider a proper, finitely presented $A$-scheme $X$, and an almost coherent $\O_X$-module $\F$. We denote the $I$-adic completion of $X$ by $\X$, so we have the following commutative diagram:
\[
\begin{tikzcd}\label{eqn:a-b-g-f}
\X \arrow{r}{c}\arrow{d}{\wdh{f}} & X \arrow{d}{f} \\
\Spf(\wdh{A}) \arrow{r} & \Spec A.
\end{tikzcd}
\]
Then we consider $4$ different cohomology groups
\[
\rm{H}^i(\X, c^*\F), \ \wdh{\rm{H}^i(X, \F)}, \ \rm{H}^i(X, \F)\otimes_A \wdh{A}, \text{ and }\lim_{n} \rm{H}^i(X_n, \F_n),
\]
and note that the are related via the following $A$-linear homomorphisms
\begin{equation}\label{diagram}
\begin{tikzcd}
\rm{H}^i(X, \F)\otimes_A \wdh{A} \arrow{r}{\a^i_{\F}} \arrow{d}{\beta^i_{\F}} & \wdh{\rm{H}^i(X, \F)} \arrow{d}{\phi^i_{\F}} \\
\rm{H}^i(\X, c^*\F) \arrow{r}{\gamma^i_{\F}} & \lim_{n} \rm{H}^i(X_n, \F_n).
\end{tikzcd}
\end{equation}

We show that all these morphisms are (almost) isomorphism:

\begin{thm}\label{formal-functions-thm} In the notation as above, all the maps $\a^i_{\F}, \ \beta^i_{\F}, \ \gamma^i_{\F}, \ \phi^i_{\F}$ are almost isomorphisms for any almost coherent $\O_X$-module $\F$. If $\F$ is quasi-coherent, almost coherent, then these maps are isomorphisms.
\end{thm}
\begin{proof}
{\it Step 0. Reduction to the case of a quasi-coherent, almost coherent sheaf $\F$}: We observe that Lemma~\ref{tensor-product}, Lemma~\ref{derived-pull-projection} and the fact that limits of two almost isomorphic direct systems are almost the same allow us to replace $\F$ with $\widetilde{\m}\otimes \F$ to assume that $\F$ is quasi-coherent and almost coherent. \smallskip

{\it Step 1. $\a_{\F}^i$ is an isomorphism}: This is just a consequence of Lemma~\ref{completion-finitely-generated} as we established in Theorem~\ref{almost-proper-mapping} that $\rm{H}^i(X, \F)$ is an almost coherent $A$-module. \smallskip

{\it Step 2. $\beta_{\F}^i$ is an isomorphism}: We note that the assumptions on $A$ imply that the map $A\to \wdh{A}$ is flat by \cite[Proposition 0.8.218]{FujKato}. Thus, flat base change for quasi-coherent cohomology groups implies that $\rm{H}^i(X, \F)\otimes_A \wdh{A}  \simeq \rm{H}^i(X_{\wdh{A}}, \F_{\wdh{A}})$. Therefore, we may and do assume that $A$ is $I$-adically complete. Then the map $\rm{H}^i(X, \F) \to \rm{H}^i(\X, c^*\F)$ is an isomorphism by Theorem~\ref{GAGA}.  \smallskip

{\it Step 3. $\alpha_{\F}^i$ is injective}: Theorem~\ref{thm:FP-approx-topology-cohomology} and Corollary~\ref{cor:formal-functions-coh} imply that the $I$-adic topology of $\rm{H}^i(X,\F)$ coincides with the natural $I$-topology. Therefore, 
\[
\widehat{\rm{H}^i(X, \F)} \simeq \lim_n \frac{\rm{H}^i(X, \F)}{\rm{Im}\left(\rm{H}^i(X, I^{n+1} \F)\to \rm{H}^i(X, \F)\right)}.
\]
Clearly, we have an inclusion 
\[
\frac{\rm{H}^i(X, \F)}{\rm{Im}\left(\rm{H}^i(X, I^{n+1} \F)\to \rm{H}^i(X, \F)\right)} \hookrightarrow \rm{H}^i(X_n, \F_n).
\]
Therefore, we conclude that $\alpha_\F^i$ is injective by left exactness of the limit functor. \smallskip

{\it Step 4. $\gamma_\F^i$ is surjective}: Recall that $\F\simeq \lim_k\F_k$ because $\F$ is adically quasi-coherent. Therefore, \cite[Corollary 0.3.2.16]{FujKato} implies that it is sufficient to show that there is a basis of opens $\cal{B}$ such that, for every $\sU \in \cal{B}$, 
\[
\rm{H}^i(\sU, \F)=0 \text{ for } i\geq 1, \text{ and}
\]
\[
\rm{H}^0(\sU,\F_{k+1}) \to \rm{H}^0(\sU, \F_k) \text{ is surjective for any } k\geq 0.
\]
Vanishing of the higher cohomology groups of adically quasi-coherent sheaves on affine formal schemes (see \cite[Theorem I.7.1.1]{FujKato}) implies that one can take $\cal{B}$ to be the basis consisting of open affine formal subschemes of $\X$. Therefore, we get that $\gamma_\F^i$ is indeed surjective for any $i\geq 0$. \smallskip

{\it Step 5.  $\alpha_{\F}^i$ and $\gamma_\F^i$ are isomorphisms}: This follows formally from commutativity of Diagram~\ref{eqn:a-b-g-f} and the previous steps.
\end{proof}

\subsection{Almost version of Grothendieck Duality}

For this section, we fix a universally coherent ring $R$ with an ideal $\m$ such that $\widetilde{\m}\coloneqq \m\otimes_R \m$ is $R$-flat and $\m^2=\m$. Since $R$ is universally coherent, there is a good theory of $f^!$ functor for morphisms between finitely presented, separated $R$-schemes\footnote{This theory does not seem to be addressed in the literature in this generality, however we all the arguments from \cite[\href{https://stacks.math.columbia.edu/tag/0DWE}{Tag 0DWE}]{stacks-project} can be adapted to this level generality with little or no extra work. See \cite[\textsection 2.1-2.2]{Zav4} for more detail.}.

\begin{prop}\label{duality-preserves-coh} Let $f\colon X \to Y$ be a morphism between separated, finitely presented $R$-schemes. Then $f^!$ sends $\bf{D}^+_{qc, acoh}(Y)$ to $\bf{D}^+_{qc, acoh}(X)$. 
\end{prop}
\begin{proof}
The only thing that we need to check here is that $f^!$ preserves almost coherence of cohomology sheaves. This statement is local, so we can assume that both $X$ and $Y$ are affine. Then we can choose a closed embedding $X \to \bf{A}^n_Y \to Y$. So, it suffices to prove the claim for a finitely presented closed immersion and for the morphism $\bf{A}^n_Y \to Y$. \smallskip

In the case $f\colon X \to Y$ a finitely presented closed immersion, we know that 
\[
f^!\F \simeq \bf{R}\ud{\cal{H}om}_{Y}(f_*\O_X, \F)
\]
for any $\F\in \bf{D}^+_{qc}(Y)$. Since $Y$ is a coherent scheme and $f$ is finitely presented, we conclude that $f_*\O_X$ is an almost coherent $\O_Y$-module. Therefore, $f^!\F= \bf{R}\ud{\cal{H}om}_{Y}(f_*\O_X, \F) \in \bf{D}_{qc, acoh}(X)$ by Corollary~\ref{alhom-acoh-derived}.

Now we consider the case of a relative affine space $f\colon X=\bf{A}^n_Y \to Y$. In this case, we have $f^!\F \simeq \mathbf{L}f^*\F \otimes^L_{\O_X} \Omega^n_{X/Y}[n]$. Then $\mathbf{L}f^*(\F) \in \bf{D}^+_{qc, acoh}(X)$ by Lemma~\ref{pullback-almost-coh-derived}\ref{pullback-almost-coh-derived-4}, and so $\mathbf{L}f^*\F \otimes^L_{\O_X} \Omega^n_{X/Y}[n] \in \bf{D}^+_{qc, acoh}(X)$ because $\Omega^n_{X/Y}$ is (non-canonically) isomorphic to $\O_X$. 
\end{proof}
 
Now we use Proposition~\ref{duality-preserves-coh} to define the almost version of the upper shriek functor: 
 
\begin{defn} Let $f\colon X \to Y$ be a morphism of separated, finitely presented $R$-schemes. We define the {\it almost upper shriek functor} $f^!_a\colon \bf{D}^+_{aqc}(Y)^a \to \bf{D}^+_{aqc}(X)^a$ as $f^!_a(\F) \coloneqq (f^!(\F_!))^a$.
\end{defn}
\begin{rmk} In what follows, we will usually denote the functor $f^!_a$ simply by $f^!$ as it will not cause any confusion. 
\end{rmk}

\begin{lemma} Let $f\colon X \to Y$ be a morphism between separated, finitely presented $R$-schemes. Then $f^!$ carries $\bf{D}^+_{acoh}(Y)^a$ to $\bf{D}^+_{acoh}(X)^a$.
\end{lemma}
\begin{proof}
This follows from Proposition~\ref{duality-preserves-coh}.
\end{proof}

\begin{thm}\label{thm:!-adjoint} Let $f\colon X \to Y$ be as above. Suppose that $f$ is proper. Then $f^!\colon \bf{D}^+_{aqc}(Y)^a \to \bf{D}^+_{aqc}(X)^a$ is a right adjoint to the functor $\bf{R}f_*\colon \bf{D}^+_{aqc}(Y)^a \to \bf{D}^+_{aqc}(X)^a$. 
\end{thm}
We note that the theorem makes sense as $\bf{R}f_*$ carries $\bf{D}^+_{aqc}(X)^a$ into $\bf{D}^+_{aqc}(Y)$ due to Lemma~\ref{pushforward-almost-coh-derived}.

\begin{proof}
This follows from a sequence of canonical isomorphisms:
\begin{align*}
\rm{Hom}_{\bf{D}(Y)^a}(\bf{R}f_*\F^a, \G^a) & \simeq \rm{Hom}_{\bf{D}(Y)}(\widetilde{\m} \otimes \bf{R}f_*\F, \G) & \text{Lemma}~\ref{adjoint-almost-sheaf-!} \\
& \simeq \rm{Hom}_{\bf{D}(Y)}(\bf{R}f_*(\widetilde{\m} \otimes  \F) , \G) & \text{Lemma}~\ref{projection}\\
& \simeq \rm{Hom}_{\bf{D}(X)}(\widetilde{\m} \otimes  \F , f^!(\G)) & \text{Grothendieck Duality}\\
& \simeq \rm{Hom}_{\bf{D}(X)^a}(\F^a , f^!(\G)^a) & \text{Lemma}~\ref{adjoint-almost-sheaf-!}. 
\end{align*}
\end{proof}

Now suppose that $f\colon X \to Y$ is a proper morphism of separated, finitely presented $R$-schemes, $\F^a\in \bf{D}^+_{aqc}(X)^a$, and $\G^a\in \bf{D}^+_{aqc}(Y)^a$. Then we want to construct a canonical morphism
\[
\bf{R}f_*\bf{R}\ud{al\cal{H}om}_{X}(\F^a, f^!(\G^a))\to \bf{R}\ud{al\cal{H}om}_Y(\bf{R}f_*(\F^a), \G^a). 
\]
Lemma~\ref{o-hom-adj-derived-sheaf} says that such a map is equivalent to a map
\[
\bf{R}f_*\bf{R}\ud{al\cal{H}om}_{X}(\F^a, f^!(\G^a)) \otimes^L_{\O_X} \bf{R}f_*(\F^a) \to \G^a. 
\]
We construct the latter map as the composition
\[
\bf{R}f_*\bf{R}\ud{al\cal{H}om}_{X}(\F^a, f^!(\G^a)) \otimes^L_{\O_X} \bf{R}f_*(\F^a) \to \bf{R}f_*\left(\bf{R}\ud{al\cal{H}om}_{X}(\F^a, f^!(\G^a)) \otimes^L_{\O_X} \F^a \right) \to \bf{R}f_*f^!\G^a \to \G^a
\]
where the first map is induced by the relative cup product (see \cite[\href{https://stacks.math.columbia.edu/tag/0B68}{Tag 0B68}]{stacks-project}), the second map comes from Remark~\ref{rmk:first-map-in-duality}, and the last map is the counit of the $(\bf{R}f_*, f^!)$-adjunction. 

\begin{lemma}\label{lemma:sheafy-duality-almost} Let $f\colon X \to Y$ be a proper morphism of separated, finitely presented $R$-schemes, $\F^a\in \bf{D}^-_{acoh}(X)^a$, and $\G^a\in \bf{D}^+_{aqc}(Y)^a$. Then the map
\[
\bf{R}f_*\bf{R}\ud{al\cal{H}om}_{X}(\F^a, f^!(\G^a))\to \bf{R}\ud{al\cal{H}om}_Y(\bf{R}f_*(\F^a), \G^a). 
\]
is an (almost) isomorphism in $\bf{D}_{aqc}^+(X)^a$. 
\end{lemma}
\begin{proof}
We note that $\bf{R}f_*\bf{R}\ud{al\cal{H}om}_{X}(\F^a, f^!(\G^a))$ lies in $\bf{D}^+_{aqc}(Y)^a$ due to Lemma~\ref{hom-alcoh-derived}~\ref{hom-alcoh-derived-4} and Lemma~\ref{pushforward-almost-coh-derived}. Likewise, $\bf{R}\ud{al\cal{H}om}_Y(\bf{R}f_*(\F^a), \G^a)$ lies in $\bf{D}^+_{aqc}(Y)^a$ by Theorem~\ref{almost-proper-mapping} and Lemma~\ref{hom-alcoh-derived}~\ref{hom-alcoh-derived-4}. Therefore, it suffices to show 
\[
\bf{R}\rm{Hom}_Y\left(\cal{H}^a, \bf{R}f_*\bf{R}\ud{al\cal{H}om}_{X}\left(\F^a, f^!\left(\G^a\right)\right)\right) \to \bf{R}\rm{Hom}_Y\left(\cal{H}^a, \bf{R}\ud{al\cal{H}om}_Y\left(\bf{R}f_*\left(\F^a\right), \G^a\right)\right)
\]
is an isomorphism for any $\cal{H}^a\in \bf{D}^+_{aqc}(Y)^a$. This follows from the following sequence of isomorphisms:
\begin{align*}
    \bf{R}\rm{Hom}_Y\left(\cal{H}^a, \bf{R}f_*\bf{R}\ud{al\cal{H}om}_{X}\left(\F^a, f^!\left(\G^a\right)\right)\right) & \simeq \bf{R}\rm{Hom}_X\left(\bf{L}f^*\cal{H}^a, \bf{R}\ud{al\cal{H}om}_{X}\left(\F^a, f^!\left(\G^a\right)\right)\right)  \\
    & \simeq \bf{R}\rm{Hom}_X\left(\bf{L}f^*\cal{H}^a\otimes^L_{\O_X} \F^a, f^!\left(\G^a\right)\right)\\
    & \simeq \bf{R}\rm{Hom}_Y\left(\bf{R}f_*\left(\bf{L}f^*\cal{H}^a\otimes^L_{\O_X} \F^a\right), \G^a\right) \\
    & \simeq \bf{R}\rm{Hom}_Y\left(\cal{H}^a \otimes \bf{R}f_*\left(\F^a \right), \G^a\right) \\
    & \simeq \bf{R}\rm{Hom}_Y\left(\cal{H}^a, \bf{R}\ud{al\cal{H}om}_Y\left(\bf{R}f_*\left(\F^a\right), \G^a\right)\right).
\end{align*}
The first isomorphism holds by Corollary~\ref{cor:*-adjoint-derived}. The second isomorphism holds by Corollary~\ref{o-hom-adj-derived-sheaf}. The third isomorphism holds by Theorem~\ref{thm:!-adjoint}. The fourth isomorphism holds by Proposition~\ref{prop:projection-formula-quasicoh}. The fifth equality holds by Corollary~\ref{o-hom-adj-derived-sheaf}.
\end{proof}

\begin{thm} Let $f\colon X \to Y$ be as above. Suppose that $f$ is smooth of pure dimension $d$. Then $f^!(-)\simeq \bf{L}f^*(-) \otimes^L_{\O_X} \Omega^d_{X/Y}[d]$
\end{thm}
\begin{proof}
It follows from the corresponding statement in the classical Grothendieck Duality. 
\end{proof}

We summarize all the results of this section in the following theorem:

\begin{thm}\label{thm:almost-grothendieck-duality} Let $R$ be a universally coherent ring with an ideal $\m$ such that $\widetilde{\m}\coloneqq \m\otimes_R \m$ is $R$-flat and $\m^2=\m$, and $\rm{FPS}_R$ be the category of finitely presented, separated $R$-schemes. Then there is a well-defined pseudo-functor $(-)^!$ from $\rm{FPS}_R$ into the $2$-category of categories such that
\begin{enumerate}[label=\textbf{(\arabic*)}]
    \item $(X)^!=\bf{D}^+_{aqc}(X)^a$,
    \item for a smooth morphism $f\colon X \to Y$ of pure relative dimension $d$, $f^! \simeq \bf{L}f^*(-)\otimes^L_{\O_X^a} \Omega^d_{X/Y}[d]$.
    \item for a proper morphism $f\colon X \to Y$, $f^!$ is right adjoint to $\bf{R}f_*\colon \bf{D}^+_{aqc}(X)^a\to \bf{D}^+_{aqc}(Y)^a$.
\end{enumerate}
\end{thm}

\newpage

\section{$\O^+/p$-modules}\label{section:proetale-v-sites}

The main goal of this Section is to discuss the comparison results between $\O^+/p$-modules in the \'etale, quasi-pro\'etale, and $v$-topologies. In particular, we show that the categories of $O^+/p$-vector bundles in all these topologies are canonically equivalent. Furthermore, one can compute cohomology groups with respect to any of these topologies (without passing to almost mathematics). A good understanding of $\O^+/p$-vector bundles in the $v$-topology will be crucial for our proof of almost coherence of nearby cycles for general $\O^+/p$-vector bundles (see Theorem~\ref{thm:main-thm-small}). We also discuss more general $\O^+/p$-modules in Section~\ref{section:etale-coefficients}. \smallskip

In this Section, we will freely use the notions of perfectoid spaces and their tilts as developed in \cite{Sch0} and \cite{Sch2}. \smallskip

\subsection{Recollection: the $v$-topology}\label{section:v-topology}

In this section, we discuss the $v$-topology on adic spaces and show some of its basic properties that seem difficult to find explicitly stated in the literature. \smallskip

Before we start this discussion, we recall the notion of a diamond and its relation to the notion of an adic space. To motivate this discussion, we remind the reader the two major problems with the category of adic spaces: the existence of non-sheafy (pre-)adic spaces, and the lack of (finite) limits in the category of adic spaces. It turns out that both of these problems go away if we consider a (pre-)adic space over $\Q_p$ as some kind of sheaf $X^\diamond$ on the category of perfectoid spaces of characteristic $p>0$. It could sound somewhat counter-intuitive to consider a $p$-adic rigid-analytic variety as a sheaf on characteristic $p$ objects, but it turns out to be quite useful in practice. The main idea is that an $S=\Spa (R, R^+)$-point of $X^\diamondsuit$ should be a choice of an untilt $S^\#$ of $S$ (this is a mixed characteristic object) and a morphism $S^\#\to X$. This procedure turns out to remember a lot of information about $X$ (e.g., \'etale cohomology), but not all information on $X$ (see Warning~\ref{warning:not-fully-faithful})

\begin{defn}\label{defn:v-topology}\cite[Definitions 8.1, 12.1, and 14.1]{Sch2}  The category $\rm{Perf}$ is the category of characteristic $p$ perfectoid spaces.\smallskip

The {\it $v$-topology} on $\rm{Perf}$ is defined by saying that a family $\{f_i\colon X_i \to X\}_{i\in I}$ of morphisms in $\rm{Perf}$ is a covering if, for any quasi-compact open $U\subset X$, there is a finite subset $I_0\subset I$ and quasi-compact opens $\{U_i\subset X_i\}_{i\in I_0}$ such that $U\subset \cup_{i\in I_0} f_i(U_i)$. \smallskip

A {\it small $v$-sheaf} is a $v$-sheaf Y on $\text{Perf}$ such that there is an epimorphism of
$v$-sheaves $Y' \to Y$ for some perfectoid space $Y'$. \smallskip

The {\it v-site} $Y_v$ of a small $v$-sheaf $Y$ is the site whose objects are all maps $Y' \to Y$ from small $v$-sheaves $Y'$, with coverings given by families $\{Y_i \to Y\}_{i\in I}$ such that $\sqcup_{i\in I} Y_i \to Y$ is an epimorphism of $v$-sheaves.
\end{defn}

\begin{rmk} The $v$-site of a small $v$-sheaf $Y$ has all finite limits by \cite[Proposition 12.10]{Sch2} and \cite[\href{https://stacks.math.columbia.edu/tag/002O}{Tag 002O}]{stacks-project}.
\end{rmk}

In what follows, we denote by $\bf{Ad}_{\Q_p}$ the category of adic spaces over $\Spa(\Q_p,\Z_p)$ and by $\bf{pAd}_{\Q_p}$ the category of pre-adic spaces over $\Spa(\Q_p, \Z_p)$ as defined in \cite[Definition 2.1.5]{Scholze-W} and \cite[Definition 8.2.3]{KedLiu1}\footnote{These spaces are called adic in \cite{Scholze-W}, we prefer to call them pre-adic to distinguish them with the usual adic spaces in the sense of Huber.}. The category of pre-adic spaces satisfies the following list of properties (see \cite[Proposition 2.1.6]{Scholze-W} or \cite[\textsection 8.2.3]{KedLiu1}):
\begin{enumerate}[label=\textbf{(\arabic*)}]
    \item The natural functor $\bf{Ad}_{\Q_p} \to \bf{pAd}_{\Q_p}$ is fully faithful,
    \item there is a functor $(\rm{Tate-Huber}_{(\Q_p, \Z_p)}^{\rm{comp}})^{\rm{op}} \to \bf{pAd}_{\Q_p}$ from the opposite category of complete Tate-Huber pairs over $(\Q_p, \Z_p)$ to pre-adic spaces over $\Spa(\Q_p, \Z_p)$. For each such $(A, A^+)$ it associates the pre-adic affinoid space\footnote{We follow \cite{KedLiu1} and use the notation $\widetilde{\Spa(A, A^+)}$ for affinoid {\it pre}-adic spaces. If $A$ is sheafy, we freely identify it with $\Spa(A, A^+)$.} $\widetilde{\Spa(A, A^+)}$;
    \item for an adic space $S$ and a pre-adic affinoid space $\widetilde{\Spa(A, A^+)}$, the set of morphisms is given by
    \[
    \rm{Hom}_{\bf{pAd}_{\Q_p}}(S,\widetilde{\Spa(A, A^+)}) = \rm{Hom}_{\rm{cont}}\big((A, A^+), (\O_S(S), \O_S^+(S))\big),
    \]
    \item $\bf{pAd}_{\Q_p}$ has all finite limits,
    \item for a pseudo-adic space $X$, one can functorially associate an underlying topological space $|X|$ such that it coincides with the usual underlying topological space $|X|$ when $X=\widetilde{\Spa(A, A^+)}$ is a pre-adic affinoid space or $X=(|X|, \O_X, \O^+_X)$ is an adic space,
    \item for every pre-adic space $X\in \bf{pAd}_{\Q_p}$, one can functorially associate an \'etale site $X_\et$ such that $X_\et$ coincides with the classical \'etale site when $X$ is a locally strongly noetherian space or a perfectoid space (see \cite[\textsection 2.1]{H3} and \cite[\textsection 7]{Sch0}).
\end{enumerate}

\begin{warning} In general, it is not true that $\rm{Hom}_{\bf{pAd}_{\Q_p}}\big(\Spa(B, B^+), \Spa(A, A^+)\big)$ is equal to $\rm{Hom}_{\rm{cont}}\big((A, A^+), (B, B^+)\big)$ unless $\Spa(B, B^+)$ is sheafy. In particular, the functor 
\[
    (\rm{Tate-Huber}_{(\Q_p, \Z_p)}^{\rm{comp}})^{\rm{op}} \to \bf{pAd}_{\Q_p}
\]
is {\it not} fully faithful. 
\end{warning}

\begin{defn}\label{defn:tilde-limit}\cite[Definition 2.4.1]{Scholze-W} Let $X_i$ be a cofiltered inverse system of pre-adic spaces with quasi-compact and quasi-separated transition maps, $X$ a pre-adic space, and $f_i\colon X \to X_i$ a compatible family of morphisms.

We say that $X$ is a {\it tilde-limit} of $X_i$, $X \sim \lim_I X_i$ if the map of underlying topological spaces $|X| \to \lim_I |X_i|$ is a homeomorphism and there is an open covering of $X$ by affinoids $\widetilde{\Spa(A, A+)} \subset X$, such that the map
\[
\colim_{\widetilde{\Spa(A_i,A_i^+)}\subset X_i} A_i \to A
\]
has dense image, where the filtered colimit runs over all open affinoid
\[
\widetilde{\Spa(A_i, A_i^+)} \subset X_i
\]
over which $\widetilde{\Spa(A, A+)} \subset X \to X_i$ factors.
\end{defn}


\begin{defn}\label{defn:diamondification}\cite[Definition 15.5]{Sch2} The {\it diamond associated} to $X\in \bf{pAd}_{\Q_p}$ is a presheaf 
\[
X^{\diamondsuit}\colon \rm{Perf}^{\rm{op}} \to \rm{Sets}
\]
such that, for any perfectoid space $S$ of characteristic $p$, we have
\[
X^{\diamondsuit}(S)=\left\{ \left(\big(S^{\sharp}, \iota\big), f\colon S^{\sharp} \to X \right)\right\}/\text{isom}
\]
where $S^{\sharp}$ is a perfectoid space, $\iota\colon (S^{\sharp})^{\flat} \to S$ is an isomorphism of the tilt of $S^\sharp$ with $S$, and $f\colon S^\sharp \to X$ is a morphism of pre-adic spaces. 

The {\it diamantine spectrum} $\Spd(A, A^+)$ of a Huber pair $(A, A^+)$ is a presheaf $\widetilde{\Spa(A, A^+)}^\diam$. 
\end{defn}

We list the main properties of this functor:

\begin{prop}\label{prop:properties-of-diamond} The diamondification functor factors through the category of $v$-sheaves. And the functor $(-)^\diam\colon \bf{pAd}_{\Q_p} \to \bf{Shv}(\rm{Perf}_v)$ satisfies the following list of properties:
\begin{enumerate}[label=\textbf{(\arabic*)}]
    \item\label{prop:properties-of-diamond-1} if $X$ is a perfectoid space, then $X^\diam \simeq X^{\flat}$,
    \item $X^\diam$ is a small $v$-sheaf for any $X\in \bf{pAd}_{\Q_p}$,\footnote{It is even a locally spatial diamond in the sense of \cite[Definition 11.17]{Sch2}.}
    \item\label{prop:properties-of-diamond-3} if $\{X_i \to X\}_{i\in I}$ is an open (resp. \'etale) covering in $\bf{pAd}_{\Q_p}$, then the family $\{X_i^\diam \to X^\diam\}_{i\in I}$ is an open (resp. \'etale) covering of $X^\diam$,
    \item there is a functorial homeomorphism $|X| \simeq |X^\diam|$ for any $X\in \bf{pAd}_{\Q_p}$,
    \item\label{prop:properties-of-diamond-5} if $X$ is a perfectoid space such that $X \sim \lim_I X_i$ in $\bf{pAd}_{\Q_p}$ with quasi-compact quasi-separated transition maps, then the natural functor $X^\diam \to \lim_{I} X^\diam_i$ is an isomorphism,
    \item\label{prop:properties-of-diamond-6} the functor $(-)^\diam \colon \colon \bf{pAd}_{\Q_p} \to \bf{Shv}(\rm{Perf}_v)$ commutes with fiber products. 
\end{enumerate}
\end{prop}
\begin{proof}
    The first claim follows from \cite[Corollary 3.20]{Sch2} and the definition of the diamondification functor. As for the second claim, \cite[Proposition 15.6]{Sch2} implies that $X^\diamondsuit$ is a diamond, and so it is a small $v$-sheaf due to \cite[Proposition 11.9]{Sch2} and the definition of a diamond (see \cite[Definition 11.1]{Sch2}). The third and fourth claims follow from \cite[Lemma 15.6]{Sch2}. The proof of the fifth claim is identical to that of \cite[Proposition 2.4.5]{Scholze-W} (the statement makes the assumption that $X$ and $X_i$ are defined over a perfectoid field, but it is not used in the proof). \smallskip
    
    Now we give a proof of the sixth claim. Let $U \to V$, $W\to V$ be morphisms in $\bf{pAd}_{\Q_p}$ with the fiber product $U\times_V W$. We fix a perfectoid space $S$ of characteristic $p$. Then we have a sequence of identifications 
    \begin{align*}
        (U\times_V W)^\diam(S) &=\left\{ \left(\big(S^{\sharp}, \iota\big), S^{\sharp} \to U\times_V W \right)\right\}/\text{isom} \\
        &=\left\{ \left(\big(S^{\sharp}, \iota\big), S^{\sharp} \to U \right)\right\}/\text{isom}\times_{\left\{ \left(\big(S^{\sharp}, \iota\big), S^{\sharp} \to V \right)\right\}/\text{isom}} \left\{ \left(\big(S^{\sharp}, \iota\big), S^{\sharp} \to W \right)\right\}/\text{isom} \\
        &=U^\diam (S) \times_{V^\diam (S)} W^\diam(S)
    \end{align*}
which is functorial in $S$. Therefore, this defines an isomorphism
\[
(U\times_V W)^{\diam}\xr{\sim} U^\diam \times_{V\diam} W^\diam. \qedhere
\]
\end{proof}

\begin{warning} The functor $(-)^\diam$ does not send the final object to the final object. In particular, it does not commute with all finite limits. 
\end{warning}

\begin{warning}\label{warning:not-fully-faithful} The functor $(-)^\diam \colon \bf{pAd}_{\Q_p} \to \bf{Shv}(\rm{Perf}_v)$ is not fully faithful. This observation is quite crucial for our proof of Theorem~\ref{coh-dimension}. In that proof, we exploit Theorem~\ref{perfectoidization} which guarantees that some non-perfectoid affinoid (pre-)adic spaces become perfectoid after diamondification.  
\end{warning}

The next goal is to discuss some examples of $v$-covers of $X^\diam$.

\begin{defn} A family of morphisms $\{f_i\colon X_i \to X\}_{i\in I}$ in $\bf{pAd}_{\Q_p}$ is a {\it naive $v$-covering} if, for any quasi-compact open $U\subset X$, there is a finite subset $I_0\subset I$ and quasi-compact opens $\{U_i\subset X_i\}_{i\in I_0}$ such that $|U|\subset \cup_{i\in I_0} |f_i|(|U_i|)$.
\end{defn}

\begin{rmk} Using that the natural morphism $|X\times_Y Z| \to |X|\times_{|Y|} |Z|$ is surjective, it is easy to see that a pullback of a naive $v$-covering is a naive $v$-covering.
\end{rmk}

\begin{lemma}\label{lemma:etale-open} Let $f\colon X \to Y$ be an \'etale morphism of pre-adic spaces in $\bf{pAd}_{\Q_p}$ (in the sense of \cite[Definition 8.2.19]{KedLiu1}). Then $f$ is an open map.
\end{lemma}
\begin{proof}
    By definition, open immersions induce open maps of underlying topological spaces. Therefore, after unravelling the definition of \'etale morphisms, it suffices to show that a map of pre-adic spaces $|\Spa(\varphi)|\colon |\Spa(B, B^+)| \to |\Spa(A, A^+)|$ is open when $\varphi \colon (A, A^+) \to (B, B^+)$ is a finite \'etale morphism of Tate-Huber pairs. In this case, Lemma~\ref{lemma:uniform-noetherian-approximation} and Corollary~\ref{cor:approximation-etale-general} allow us to assume that $(A, A^+)$ is strongly noetherian. In this case, the result follows from  \cite[Lemma 1.7.9]{H3} (alternatively, one can directly adapt the proof of \cite[Lemma 1.7.9]{H3} to work in the non-noetherian case).
\end{proof}

\begin{example}\label{exmpl:naive-v-covering} 
\begin{enumerate}
    \item A quasi-compact surjective morphism $X\to Y$ of pre-adic spaces over $\Spa(\Q_p, \Z_p)$ is a naive $v$-cover; 
    \item Lemma~\ref{lemma:etale-open} implies that a family of jointly surjective \'etale morphisms $\{X_i \to X\}$ of pre-adic spaces over $\Spa(\Q_p, \Z_p)$ is a naive $v$-cover.
\end{enumerate}
\end{example}

Our next goal is to show that the diamondification functor $(-)^\diamond$ sends naive $v$-covers to surjections of small $v$-sheaves.

\begin{lemma}\label{lemma:qcqs} Let $f\colon X \to Y$ be a quasi-compact (resp.~quasi-separated) morphism in $\bf{pAd}_{\Q_p}$. Then $f^\diam \colon X^\diam \to Y^\diam$ is quasi-compact (resp.~quasi-separated) in the sense of \cite[p.40]{Sch2}.
\end{lemma}
\begin{proof}
    We first deal with a quasi-compact $f$. To check that $f^\diam$ is quasi-compact, it suffices to show that $S\times_{Y^\diam} X^\diam$ is quasi-compact for any morphism $S \to Y$ with an affinoid perfectoid $S$. By definition, this morphism corresponds to a morphism $S^\sharp \to Y$ with an affinoid perfectoid source $S^\sharp$. By Proposition~\ref{prop:properties-of-diamond}, we have $S\times_{Y^\diam} X^\diam \simeq (S^\sharp \times_{Y} X)^\diam$, so \cite[Lemma 15.6]{Sch2} implies that 
   \[
    |S\times_{Y^\diam} X^\diam| \simeq |S^\sharp\times_Y X|
    \]
    quasi-compact by our assumption on $f$. Now $S\times_{Y\diam} X^\diam$ is quasi-compact due to the combination of \cite[Proposition 12.14(iii) and Lemma 15.6]{Sch2}. \smallskip
    
    The case of a quasi-separated $f$ follows from Proposition~\ref{prop:properties-of-diamond} and the quasi-compact case by considering the diagonal morphism $\Delta_f\colon X \to X\times_Y X$. 
\end{proof}

\begin{lemma}\label{lemma:naive-normal-v-covering} Let $\{f_i\colon X_i \to X\}_{i\in I}$ be a naive $v$-covering in $\bf{pAd}_{\Q_p}$. Then $\{f_i^\diam \colon X_i^\diam \to X^\diam\}_{i\in I}$ is a $v$-covering.
\end{lemma}
\begin{proof}
    We can find a covering $\{U_j \to X\}_{j\in J}$ by open affinoids. Since $\{U_j^\diam \to X^\diam\}$ is a $v$-covering by Proposition~\ref{prop:properties-of-diamond}, it suffices to show that $\{f_{i,j}\colon X_{i,j}\coloneqq X_i\times_X U_j \to U_j\}_{i\in I}$ is a $v$-covering for every $j\in J$. Since naive $v$-covers are preserved by open base change, we reduce to the case when $X$ is an affinoid. \smallskip
    
    Moreover, the proof of \cite[Proposition 15.4]{Sch2} ensures that there is a $v$-surjection $f\colon S \to X^\diam$ where $S$ is an affinoid perfectoid space. By definition, the map $f$ corresponds to a map $g\colon S^\sharp \to X$. Proposition~\ref{prop:properties-of-diamond} ensures that diamondization commutes with fiber products, so it suffices to show that $\{(X_i\times_X S^\sharp)^\diam \to (S^\sharp)^\diam\}_{i\in I}$ is a $v$-covering. In other words, we can assume that $X=S^\sharp$ is an affinoid perfectoid space. \smallskip

    Now we can find a covering $\{U_{i,j} \to X_i\}_{j\in J_i}$ by open affinoids for each $i\in I$. Then the family $\{U_{i, j}\to X\}_{i\in I, j\in J_i}$ is also a naive $v$-covering, and so it suffices to show that $\{U_{i,j}^\diam \to X^\diam\}_{i\in I, j\in J_i}$ is a $v$-covering. In other words, we can assume that $X$ is an affinoid perfectoid space and that all $X_i$ are affinoids. A similar argument allows us to assume that each $X_i$ is an affinoid perfectoid space. \smallskip
    
    Finally, we note that under our assumption that $X$ and $X_i$ are (affinoid) perfectoids, $\{X_i \to X\}_{i\in I}$ is a naive $v$-covering if and only if $\{X_i^\diam \to X^\diam\}_{i\in I}$ is a $v$-covering since $|X_i^\diam|\simeq |X_i|$ and $|X^\diam|=|X|$ by \cite[Lemma 15.6]{Sch2}. 
\end{proof}

\subsection{Recollection: the quasi-pro\'etale topology}

The main goal of this section is to recall the notions of a quasi-pro\'etale morphism and the quasi-pro\'etale topology. This topology will be a crucial intermediate tool to relate the $v$-topology to the \'etale topology.\smallskip

In this section, we will work only with strongly sheafy spaces in the sense of Definition~\ref{defn:strongly-sheafy}. We advise the reader to look at Section~\ref{section:appendix-etale} for basic definitions involving such spaces. Most likely, this discussion can be generalized to arbitrary affinoid pre-adic spaces, but we do not do this since we will never need this level of generality. \smallskip

For the next definition, we fix a morphism $f\colon X=\Spa(S, S^+) \to Y=\Spa(R, R^+)$ of strongly sheafy Tate-affinoid adic spaces. 

\begin{defn}\label{defn:proetale} A morphism $f\colon \Spa(S, S^+) \to \Spa(R, R^+)$ is an {\it affinoid strongly pro-\'etale morphism} if there is a cofiltered system of strongly \'etale morphisms (see Definition~\ref{defn:strongly-etale-spaces}) 
\[
    \Spa(R_i, R_i^+) \to \Spa(R, R^+)
\]
of strongly sheafy affinoid adic spaces such that  $(S, S^+) = \big(\wdh{(\colim_I R_i)}_{u}, \wdh{(\colim_I R_i)}_{u}^+\big)$ is the complete uniform filtered colimit of $(R_i, R^+_i)$ (see Definition~\ref{defn:uniform-colimit}). 

We will usually write $\Spa(S, S^+) \approx \lim_I \Spa(R_i, R_i^+) \to \Spa(R, R^+)$ for an affinoid strongly pro-\'etale presentation of $\Spa(S, S^+) \to \Spa(R, R^+)$.


\end{defn}

\begin{rmk} Explicitly, Remark~\ref{rmk:explicit-uniform-completion} implies that $S^+ = (\colim_I R_i^+)^{\wedge}_{\varpi}$ is equal to the $\varpi$-adic completion of $\colim_I R_i^+$ and $S=S^+\big[\frac{1}{\varpi}\big]$ for any choice of a pseudo-uniformizer $\varpi\in R^+$.
\end{rmk}

\begin{rmk} We note that Theorem~\ref{thm:approximation-etale}\ref{thm:approximation-etale-11} (see also \cite[Proposition 2.4.2]{Scholze-W}) implies that $\Spa(S, S^+) \sim \lim_I \Spa(R_i, R_i^+)$ for an affinoid strongly pro-\'etale morphism $\Spa(S, S^+) \approx \lim_I \Spa(R_i, R_i^+) \to \Spa(R, R^+)$.
\end{rmk}

\begin{warning}\label{warning:scholze-affinoid-proetale} Definition~\ref{defn:proetale} is more restrictive than \cite[Definition 7.8]{Sch2} when $\Spa(R, R^+)$ is an affinoid perfectoid space. 
\end{warning}



\begin{defn}\label{defn:totally-disconnected} 
A perfectoid space $X$ is {\it strictly totally disconnected} if $X$ is quasi-compact, quasi-separated, and every \'etale cover of $X$ splits.
\end{defn}



\begin{lemma}\label{lemma:good-properties-affinoid-proetale} Let each $X$, $Y$, $Y'$, and $Z$ be affinoid spaces over $\Spa(\Q_p, \Z_p)$. We assume that each of them is strongly sheafy. 
\begin{enumerate}[label=\textbf{(\arabic*)}]
\item\label{lemma:good-properties-affinoid-proetale-1} Let $f\colon X \to Y$ and $g\colon Y \to Z$ be affinoid strongly pro-\'etale morphisms. Then the composition $g\circ f \colon X \to Z$ is also an affinoid strongly pro-\'etale morphism;
\item\label{lemma:good-properties-affinoid-proetale-2} Let $f\colon X \to Y$ be an affinoid strongly pro-\'etale morphism, let $g\colon Y'\to Y$ be a morphism of adic spaces with $Y'$ being an affinoid perfectoid space (resp. strictly totally disconnected perfectoid space), and let $X_{Y'} \coloneqq X\times_Y Y'$ be the fiber product (in pre-adic spaces). Then $X_{Y'}^\diam$ is an affinoid perfectoid space (resp. strictly totally disconnected perfectoid space) and the morphism $f_{Y'}^\diam\colon X_{Y'}^\diam \to Y'^{\diam}$ is an affinoid pro-\'etale morphism in the sense of \cite[Definition 7.8]{Sch2}.
\end{enumerate}
\end{lemma}
\begin{proof}
    \ref{lemma:good-properties-affinoid-proetale-1} The proof of \cite[Lemma 2.5(1)]{Lucas-Annette} goes through if we use Theorem~\ref{thm:approximation-etale} in place of \cite[Proposition 6.4]{Sch2} (and \cite[Proposition 1.7.1]{H3}). \smallskip

    Now we show \ref{lemma:good-properties-affinoid-proetale-2}. We set $X=\Spa(S, S^+)$, $Y=\Spa(R, R^+)$, $Y'=\Spa(R', R'^+)$, and let $X\approx \lim_I \big(X_i=\Spa(R_i, R_i^+) \big) \to Y=\Spa(R, R^+)$ be an affinoid strongly pro-\'etale presentation of $X \to Y$. Then Proposition~\ref{prop:properties-of-diamond}\ref{prop:properties-of-diamond-5} implies that
    \[
    X^\diam = \lim_I X_i^\diam.
    \]
    Thus, $X^\diam_{Y'} = \lim_I (X_i \times_Y Y')^\diam \to Y'^\diam$. So it suffices to show that each $(X_i \times_Y Y')^\diam$ is represented by an affinoid perfectoid space and the morphism $f_i^\diam\colon (X_i \times_Y Y')^\diam \to Y'^\diam$ is \'etale. By construction, $f_i^\diam$ is \'etale. In particular, $(X_i \times_Y Y')^\diam$ is represented by a perfectoid space. Furthemore, $f_i^\diam$ is a composition of finite \'etale maps and finite disjoint unions of rational subdomainds. Therefore, $(X_i \times_Y Y')^\diam$ is an affinoid perfectoid space due to the combination of \cite[Theorem 6.3 and Theorem 7.9]{Sch0}. \smallskip

    If $Y'$ is strictly totally disconnected, then \cite[Lemma 7.19]{Sch2} implies that $X_{Y'}^\diam$ is also represented by a strictly totally disconnected perfectoid space. 
\end{proof}

\begin{warning} \cite[Lemma 2.5(1)]{Lucas-Annette} claims a stronger version of Lemma~\ref{lemma:good-properties-affinoid-proetale}\ref{lemma:good-properties-affinoid-proetale-2}. However, it seems to be false (see Warning~\ref{warning:perfectoid-only-after-diamondification}).
\end{warning}

Now we are ready to show that the issue raised in Warning~\ref{warning:scholze-affinoid-proetale} disappears when the target is a strictly totally disconnected perfectoid space.

\begin{lemma}\label{lemma:no-difference} Let $X=\Spa(R, R^+)$ be a strictly totally disconnected perfectoid space, and let $f\colon Y=\Spa(S, S^+) \to X=\Spa(R, R^+)$ be an affinoid pro-\'etale morphism (in the sense of \cite[Definition 7.8]{Sch2}). Then $f$ is an affinoid strongly pro-\'etale morphism. 
\end{lemma}
\begin{proof}
    The proof of \cite[Lemma 7.19]{Sch2} ensures that $f$ can be realized as a pro-(rational subdomain) inside the pro-(finite \'etale) morphism $X\times_{\pi_0(X)} \pi_0(Y)$. Each of these morphisms is an affinoid strongly pro-\'etale morphism. Therefore, Lemma~\ref{lemma:good-properties-affinoid-proetale}\ref{lemma:good-properties-affinoid-proetale-1} ensures that $f$ is an affinoid strongly pro-\'etale morphism as well. 
\end{proof}

For the next definition, we fix a morphism $f\colon X \to Y$ of adic spaces such that $X$ and $Y$ are strongly sheafy adic spaces over $\Spa(\Q_p, \Z_p)$. 

\begin{defn} A morphism $f\colon X \to Y$ is {\it strongly pro-\'etale} if, for every point $x\in X$, there is an open affinoid $x\in U\subset X$ and an open affinoid $f(x)\in V \subset Y$ such that $f|_{U} \colon U \to V$ is affinoid strongly pro-\'etale.
\end{defn}


Now we are ready to define quasi-pro\'etale morphisms.

\begin{defn}\cite[Definition 10.1 and 14.1]{Sch2} A morphism of small $v$-sheaves $f\colon X \to Y$ is {\it quasi-pro\'etale} if it is locally separated, and for every morphism $S \to Y$ with a strictly totally disconnected perfectoid $S$, the fiber product $X_S\coloneqq X\times_Y S$ is represented by a perfectoid space and $X_S\to S$ is pro-\'etale.

The {\it quasi-pro\'etale site} $X_\qproet$ of a small $v$-sheaf is the site whose objects are quasi-pro\'etale morphisms $Y \to X$, with coverings given by families $\{Y_i \to Y\}_{i\in I}$ such that $\sqcup_{i\in I} Y_i \to Y$ is a surjection of $v$-sheaves. 
\end{defn}

\begin{lemma}\label{lemma:proetale-quasiproetale} Let $f\colon X \to Y$ be a strongly pro-\'etale morphism such that both $X$ and $Y$ are strongly sheafy adic spaces over $\Spa(\Q_p, \Z_p)$. Then $f^\diam \colon X^\diam \to Y^\diam$ is quasi-pro\'etale. Furthermore, if $f$ is also a naive $v$-covering, then $f^\diam$ is a quasi-pro\'etale covering. 
\end{lemma}
\begin{proof}
    The question is local on the source and on the target, so we can assume that $f$ is an affinoid strongly pro-\'etale morphism. Then it is easy to see that $f^\diam \colon X^\diam \to Y^\diam$ is a separated morphism (for example, it is quasi-separated due to Lemma~\ref{lemma:qcqs} and then the valuative criterion \cite[Proposition 10.9]{Sch2} implies that it is separated). Therefore, it suffices to show that, for any strictly totally disconnected perfectoid $S$ and a morphism $S \to Y^\diam$, the fiber product $S\times_{Y^\diam}  X^\diam \to S$ is a pro-\'etale morphism of perfectoid spaces. \smallskip

    Now we recall that a morphism $f \colon S \to Y^\diam$ uniquely corresponds to a morphism $g\colon S^\sharp \to Y$. Proposition~\ref{prop:properties-of-diamond}~\ref{prop:properties-of-diamond-6} implies that 
    \[
    S\times_{Y^\diam} X^\diam\simeq (S^\sharp \times_Y X)^\diam.
    \]
    Therefore, Lemma~\ref{lemma:good-properties-affinoid-proetale}\ref{lemma:good-properties-affinoid-proetale-2} implies that $S \times_{Y^\diam} X^\diam \to S$ is affinoid pro-\'etale in the sense of \cite[Definition 7.8]{Sch2}. This finishes the proof that $f^\diam$ is quasi-pro\'etale. If we also assume that $f$ is a naive $v$-covering, then Lemma~\ref{lemma:naive-normal-v-covering} ensures that $f^\diam$ is a surjection of $v$-sheaves. Thus, $f^\diam$ is a quasi-pro\'etale covering in this case. 
\end{proof}

Finally, we wish to show that strongly sheafy Tate-affinoids $\Spa(A, A^+)$ admit affinoid strongly pro-\'etale covers by strictly totally disconnected perfectoid spaces. For this, we will need some preliminary lemmas:

\begin{lemma}\label{lemma:strictly-totally-disconnected-criterion} Let $(A, A^+)$ be an affinoid perfectoid pair. Suppose that every surjective (affinoid) strongly \'etale morphism $\Spa(B, B^+) \to \Spa(A, A^+)$ admits a section (see Definition~\ref{defn:strongly-etale-spaces}). Then $\Spa(A, A^+)$ is a strictly totally disconnected perfectoid space. 
\end{lemma}
\begin{proof}
    It suffices to show that every \'etale surjective morphism $X \to \Spa(A, A^+)$ admits a section. Any such morphism can be dominated by a surjective morphism of the form $\sqcup_{i\in I} X_i \to \Spa(A, A^+)$ where $X_i=\Spa(B_i, B_i^+) \to \Spa(A, A^+)$ is affinoid strongly \'etale and $I$ is a finite set. Then Remark~\ref{rmk:union-strongly-etale} implies that $\sqcup_{i\in I} X_i \to \Spa(A_i, A_i^+)$ is itself strongly \'etale (and affinoid), so it admits a section due to the assumption on $X$. Therefore, $X\to \Spa(A, A^+)$ also admits a section.
\end{proof}

\begin{lemma}\label{lemma:str-tot-disc-covering} Let $\Spa(A, A^+)$ be a strongly sheafy Tate-affinoid space over $(\Q_p, \Z_p)$. Then there is an affinoid strongly pro-\'etale covering $\Spa(A_\infty, A_\infty^+) \to \Spa(A, A^+)$ such that the fiber products $\Spd(A_\infty, A_\infty^+)^{j/\Spd(A, A^+)}$ are represented by strictly totally disconnected (affinoid) perfectoid spaces for $j\geq 1$. In particular, $\Spd(A_\infty, A_\infty^+)\to \Spd(A, A^+)$ is a quasi-pro\'etale covering by a strictly totally disconnected perfectoid space. 
\end{lemma}
\begin{proof}
    For the purpose of this proof, we say that a strongly \'etale morphism of complete Tate-Huber pairs $f\colon (R, R^+) \to (S, S^+)$ is a covering if $|\Spa(f)| \colon |\Spa(S, S^+)| \to |\Spa(R, R^+)|$ is surjective. \smallskip

    We fix a set of representatives of all strongly \'etale coverings $\{(A, A^+) \to (A_i, A_i^+)\}_{i\in I}$. Then, for each finite subset $S\subset I$, we define 
    \[
    (A_S, A_S^+) = \wdh{\otimes}_{s\in S} (A_s, A_s^+).
    \]
    Each $(A_S, A_S^+)$ is strongly \'etale covering of $(A, A^+)$. For each $S\subset S'$, we put $f_{S, S'}\colon (A_S, A_S^+) \to (A_{S'}, A_{S'}^+)$ to be the natural morphism induced by $S \hookrightarrow S'$. Then we see that $\{(A_S, A_S^+), f_{S, S'}\}_{S\subset I \text{ finite}}$ is a filtered system of strongly \'etale $(A, A^+)$-algebras. We put 
    \[
    (A(1), A^+(1)) =  \Big((\wdh{\colim_S A^+_S})\big[\frac{1}{p}\big], \wdh{\colim_S A^+_S}\Big)
    \]
    to be the completed uniform filtered colimit of $(A_S, A_S^+)$ (see Definition~\ref{defn:uniform-colimit}). Theorem~\ref{thm:approximation-etale} implies that every strongly \'etale covering $(A, A^+) \to (B, B^+)$ admits a splitting over $(A(1), A^+(1))$. We repeat the same construction to inductively define
    \[
    \Big(A(2), A^+(2)\Big) \coloneqq \Big(A(1)(1), A(1)^+(1)\Big), \Big(A(3), A^+(3)\Big) \coloneqq \Big(A(2)(1), A(2)^+(1)\Big), \dots,  
    \]
    \[
    \Big(A(n), A^+(n)\Big) \coloneqq \Big(A(n-1)(1), A(n-1)^+(1)\Big), \dots .
    \]
    Finally, we define $(A_\infty, A^+_\infty)$ to be the completed uniform filtered colimit of $(A(n), A(n)^+)$. Then Theorem~\ref{thm:approximation-etale} implies that any strongly \'etale covering of $(A_\infty, A_\infty^+)$ comes from a covering of some $(A(n), A^+(n))$, thus it admits a splitting over $(A(n+1), A^+(n+1))$. In particular, every strongly \'etale covering of $(A_\infty, A_\infty^+)$ admits a splitting. The proof of \cite[Lemma 15.3]{Sch2} implies that $(A_\infty, A_\infty^+)$ is a perfectoid pair. In particular, it is strongly sheafy. Furthermore, Lemma~\ref{lemma:strictly-totally-disconnected-criterion} ensures that it is strictly totally disconnected. We notice that the morphism $\Spa(A_\infty, A_\infty^+) \to \Spa(A, A^+)$ is an affinoid strongly pro-\'etale covering. Finally, we see that all fiber products $\Spd(A_\infty, A_\infty^+)^{j/\Spd(A, A^+)}$ are represented by strictly totally disconnected perfectoid spaces due to Lemma~\ref{lemma:good-properties-affinoid-proetale}\ref{lemma:good-properties-affinoid-proetale-2}.
\end{proof}

\begin{lemma}\label{lemma:basis-proetale} Let $X=\Spa(A, A^+)$ be a strongly sheafy Tate-affinoid over $\Spa(\Q_p, \Z_p)$. Then the set of all $f^\diam \colon Y^\diam \to X^\diam$ for affinoid perfectoid $Y$ with an affinoid strongly pro-\'etale morphism $f\colon Y \to X$ forms a basis of $X^\diam_\qproet$.
\end{lemma}
\begin{proof}
    Let $Z \to X^\diam$ be a quasi-pro\'etale morphism. We wish to show that it can be covered (in the quasi-pro\'etale topology) by elements of the form $Y^\diam \to X^\diam$ for an affinoid perfectoid $Y$ and an affinoid strongly pro-\'etale morphism $Y\to X$. \smallskip
    
    Then Lemma~\ref{lemma:str-tot-disc-covering} implies that we can find an affinoid strongly pro-\'etale covering $X' \to X$ such that $X'$ is a strictly totally disconnected perfectoid space. Since $Z\to X^\diam$ is quasi-pro\'etale, we conclude that $Z\times_{X^\diam} X'^\diam$ is a perfectoid space and $Z\times_{X^\diam} X'^\diam$ is pro-\'etale. Therefore, we can cover it (in the analytic topology) by affinoid perfectoid spaces $Z_i$ such that each $Z_i$ is an affinoid perfectoid space and $Z_i \to X'^\diam$ is affinoid pro-\'etale. By construction $\{Z_i \to Z\}_{i\in I}$ is a covering in the quasi-pro\'etale topology. \smallskip
    
    Now Lemma~\ref{lemma:no-difference} implies that each $Z_i \to X'^\diam$ is affinoid strongly pro-\'etale. Therefore, when we pass to the corresponding untilts, we get morphisms $Z^\sharp_i \to X'$ that are affinoid strongly pro-\'etale as well (we use \cite[Theorem 3.12 and Theorem 6.1]{Sch2}). Therefore, Lemma~\ref{lemma:good-properties-affinoid-proetale}\ref{lemma:good-properties-affinoid-proetale-1} implies that each $Z_i^\sharp \to X$ is an affinoid strongly pro-\'etale morphism (with an affinoid perfectoid $Z_i^\sharp$). By construction $(Z_i^\sharp)^\diam = Z_i \to X^\diam$ cover the morphism $Y \to X^\diam$.
\end{proof}


\subsection{Integral structure sheaves}

In this section, we define various structure sheaves associated with a (pre-)adic space over $\Q_p$. Then we discuss the relationship between some of these sheaves. We will continue the discussion between these sheaves (and their cohomology) in the next section. \smallskip 

First, we note that the \'etale, quasi-pro\'etale, and $v$-sites of a pre-adic space $X$ over $\Spa(\Q_p, \Z_p)$ are related via the following sequence of morphisms of sites
\begin{equation}\label{eqn:many-morphisms}
\begin{tikzcd}
X^\diam_v \arrow{r}{\lambda} & X^\diam_\qproet \arrow{r}{\mu} & X_\et,
\end{tikzcd}
\end{equation}
which essentially come from the fact that any \'etale covering is a quasi-pro\'etale covering, and any quasi-pro\'etale covering is a $v$-covering.\footnote{To show that the natural continuous functors $X_\et \to X^\diam_\qproet$ and $X^\diam_\qproet \to X^\diam_v$ induce morphisms of sites (in the other direction), one needs to verify that all these sites admit finite limits and these functors commute with all finite limits. We leave this as an exercise to the interested reader.} Now we define various structure sheaves on each of these sites:

\begin{defn}\label{defn:v-structure-sheaf} Let $X$ be a pre-adic space over $\Spa(\Q_p, \Z_p)$. 

The {\it integral ``untilted'' structure sheaf} $\O_{X^\diam}^+$ is a sheaf of rings on $X^\diam_v$ obtained as the sheafification of a pre-sheaf defined by the assignment
\[
\{S \to X^\diam\} \mapsto \O^+_{S^\sharp}(S^\sharp)
\]
for any perfectoid space $S \to X^\diam$ over $X^\diam$ (the transition maps are defined in the evident way\footnote{Recall that a morphism $S\to X^\diam$ is, by definition, a data of an untilt $S^\sharp$ with a morphism $S^\sharp \to X$ and an isomorphism $(S^\sharp)^\flat \simeq S$. Thus, a pair of morphisms $T \to S \to X^\diam$ defines a pair $T^\sharp \to S^\sharp \to X$}).

The {\it rational ``untilted'' structure sheaf} $\O_{X^\diam}$ is a sheaf of rings on $X^\diam_v$ given by the formula  $\O_{X^\diam} = \O_{X^\diam}^+[\frac{1}{p}]$.

The {\it mod-$p$ structure sheaf} $\O_{X^\diam}^+/p$ is the quotient of $\O_{X^\diam}^+$ by $p$ in the category of sheaves of rings on $X^\diam_v$.

The {\it quasi-pro\'etale integral ``untilted'' structure sheaf} $\O_{X^\diam_\qp}^+$ is the restriction of $\O_{X^\diam}^+$ to the quasi-pro\'etale site of $X^\diam$, i.e. $\O_{X^\diam_\qp}^+=\lambda_*\O_{X^\diam}^+$. 

The {\it quasi-pro\'etale mod-$p$ structure sheaf} $\O_{X^\diam_\qp}^+/p$ is the quotient of $\O_{X^\diam_\qp}^+$ by $p$ in the quasi-pro\'etale site $X^\diam_\qproet$.

If $X$ is a strongly sheafy space over $\Spa(\Q_p, \Z_p)$, the {\it \'etale mod-$p$ structure sheaf} $\O_{X_\et}^+/p$ is the quotient of $\O_{X_\et}^+$ by $p$ in the \'etale site $X_\et$ (see Definition~\ref{defn:etale-structure-pre-sheaves} and Lemma~\ref{lemma:etale-structure-sheaves}).
\end{defn}

\begin{rmk} We note that it is not, a priori, clear whether $\O_{X^\diam_\qp}^+/p\simeq \lambda_*\left(\O_{X^\diam}^+/p\right)$. The problem comes from the fact that $\lambda$ is not an exact functor, so it is not clear whether it commutes with quotiening by $p$.
\end{rmk}

\begin{rmk} The relation between $\O_{X^\diam_\qp}^+/p$ and $\O_{X_\et}^+/p$ is even more mysterious. The first sheaf is defined via descent from perfectoid spaces, so it seems subtle to control values of this sheaf on locally noetherian adic spaces. On the opposite, the second sheaf is defined using the \'etale topology on $X_\et$, so its definition has no direct relation to perfectoid spaces when $X$ is a locally noetherian  adic space. 
\end{rmk}

By definition, for a strongly sheafy adic space $X$ over $\Spa(\Q_p, \Z_p)$, we can promote Diagram~(\ref{eqn:many-morphisms}) to a diagram of morphisms of ringed sites:

\begin{equation}\label{eqn:many-morphisms-ringed}
\begin{tikzcd}
\left(X^\diam_v,\O_{X^\diam}^+/p\right) \arrow{r}{\lambda} & \left(X^\diam_\qproet, \O_{X^\diam_\qp}^+/p\right) \arrow{r}{\mu} & \left(X_\et, \O_{X_\et}^+/p\right).
\end{tikzcd}
\end{equation}

We also have ``tilted'' versions of the structure sheaves:

\begin{defn}\label{defn:v-structure-sheaf-tilted} Let $X$ be a pre-adic space over $\Spa(\Q_p, \Z_p)$. 

The {\it integral ``tilted'' structure sheaf} $\O_{X^\diam}^{\flat, +}$ is the sheaf of rings on $X^\diam_v$ obtained as the sheafification of a pre-sheaf defined by the assignment
\[
\{S \to X^\diam\} \mapsto \O^+_{S}(S)
\]
for any perfectoid space $S \to X^\diam$ over $X^\diam$.

If $X$ is a pre-adic space over a $p$-adic perfectoid pair $(R, R^+)$ with a good pseudo-uniformizer $\varpi\in R^+$ (see Definition~\ref{defn:good-unifor}), a {\it rational ``tilted'' structure sheaf} $\O_{X^\diam}^{\flat}$ is $\O_{X^\diam}^{\flat, +}[\frac{1}{\varpi^\flat}]$. 
\end{defn}

We start with some easy properties of these structure sheaves:

\begin{lemma}\label{lemma:first-properties-structure-sheaves} Let $X\in \bf{pAd}_{\Q_p}$ be a pre-adic space over $\Spa(\Q_p, \Z_p)$. Then 
\begin{enumerate}[label=\textbf{(\arabic*)}]
    \item\label{lemma:first-properties-structure-sheaves-1} for any affinoid perfectoid space $Y=\Spa(S, S^+)\to  X^\diam$, we have $\rm{H}^0(Y, \O_{X^\diam}^+)=S^{\sharp, +}$ and $\rm{H}^i(Y, \O_{X^\diam}^+)\simeq^a 0$ for $i\geq 1$;
    \item\label{lemma:first-properties-structure-sheaves-2} for any affinoid perfectoid space $Y=\Spa(S, S^+) \to X^\diam$, we have $\rm{H}^0(Y, \O_{X^\diam}^{+, \flat})=S^+$ and $\rm{H}^i(Y, \O_{X^\diam}^{+, \flat})\simeq^a 0$ for $i\geq 1$;
    \item\label{lemma:first-properties-structure-sheaves-3} the sheaf $\O_{X^\diam}^+$ is derived $p$-adically complete and $p$-torsionfree;
    \item\label{lemma:first-properties-structure-sheaves-4} if $X$ is pre-adic space over a perfectoid pair $(R, R^+)$ with a good pseudo-uniformizer $\varpi\in R^+$, the sheaf $\O_{X^\diam}^{+, \flat}$ is derived $\varpi^\flat$-adically complete and $\varpi^\flat$-torsionfree;
    \item\label{lemma:first-properties-structure-sheaves-5} if $X$ is pre-adic space over a perfectoid pair $(R, R^+)$ with a good pseudo-uniformizer $\varpi\in R^+$, there is a canonical isomorphism $\O_{X^\diam}^+/p \simeq \O_{X^\diam}^{+, \flat}/\varpi^\flat$.
\end{enumerate}
\end{lemma}
\begin{proof}
    \ref{lemma:first-properties-structure-sheaves-1} and \ref{lemma:first-properties-structure-sheaves-2} follow directly from \cite[Theorem 8.7 and Proposition 8.8]{Sch2}. \smallskip
    
    \ref{lemma:first-properties-structure-sheaves-3}: To show that $\O_{X^\diam}^+$ is $p$-torsionfree, it suffices to show that $\O_{X^\diam}^+(U)$ is $p$-torsionfree on a basis of $X_v^\diam$. Therefore, it is enough to show that
    \[
    \O_{X^\diam}^+\left(Y\right)
    \]
    is $p$-torsionfree for any affinoid perfectoid space $Y\to X^\diam$. This follows from \ref{lemma:first-properties-structure-sheaves-1}. \smallskip
    
    Lemma~\ref{lemma:derived-complete-global-sections} ensures that, for the purpose of proving that $\O_{X^\diam}^+$ is $p$-adically derived complete, it suffices to show that
    \[
    \bf{R}\Gamma(S, \O_{X^\diam}^+)
    \]
    is derived $p$-adically complete for any affinoid perfectoid space $Y=\Spa(S, S^+)\to X$. Then it suffices to show that each cohomology group $\rm{H}^i(Y, \O_{X^\diam}^+)$ is derived $p$-adically complete. Now \ref{lemma:first-properties-structure-sheaves-1} implies that
    \[
    \rm{H}^0(Y, \O_{X^\diam}^+)=S^{\sharp, +}
    \]
    is $p$-adically complete, and thus it is derived $p$-adically complete (see \cite[\href{https://stacks.math.columbia.edu/tag/091R}{Tag 091R}]{stacks-project}). Moreover, \ref{lemma:first-properties-structure-sheaves-1} implies that all higher cohomology groups
    \[
    \rm{H}^i(Y, \O_{X^\diam}^+)\simeq^a 0
    \]
    are almost zero. In particular, they are $p$-torsion, and so derived $p$-adically complete. Thus, $\bf{R}\Gamma(S, \O_{X^\diam}^+)$ is derived $p$-adically complete finishing the proof. \smallskip
    
    \ref{lemma:first-properties-structure-sheaves-4}: This is completely analogous to the proof of \ref{lemma:first-properties-structure-sheaves-3} using \ref{lemma:first-properties-structure-sheaves-2} in place of \ref{lemma:first-properties-structure-sheaves-1}. \smallskip
    
    \ref{lemma:first-properties-structure-sheaves-5}: Denote by $\F$ the {\it presheaf} quotient of $\O_{X^\diam}^+$ by $p$, and by $\G$ the presheaf quotient of $\O_{X^\diam}^{\flat, +}$ by $\varpi^\flat$. It suffices to construct a functorial isomorphism
    \[
    \F(U) \simeq \G(U)
    \]
    on a basis of $X^\diam_v$. Therefore, it suffices to construct such an isomorphism for any affinoid perfectoid space $U \to X^\diam$. Then \ref{lemma:first-properties-structure-sheaves-1} and \ref{lemma:first-properties-structure-sheaves-2} ensure that, for an affinoid perfectoid space $U=\Spa(S, S^+) \to X^\diam$, 
    \[
    \F(U) \simeq S^{\sharp, +}/pS^{\sharp, +} , \text{ and } \G(U) \simeq S^+/\varpi^\flat S^+.
    \]
    Essentially by the definition of a tilt, we have a canonical isomorphism 
    \[
    S^{\sharp, +}/pS^{\sharp, +} = S^{\sharp, +}/\varpi S^{\sharp, +} \simeq S^+/\varpi^\flat S^+
    \]
    finishing the proof. 
\end{proof}

\begin{rmk}\label{rmk:etale-derived-complete} The conclusion of Lemma~\ref{lemma:first-properties-structure-sheaves}\ref{lemma:first-properties-structure-sheaves-1},\ref{lemma:first-properties-structure-sheaves-3} holds for the sheaf $\O_{X^\diam_\qp}^+$ with a similar proof (using \cite[Theorem 8.5]{Sch2} in place of \cite[Theorem 8.7 and Proposition 8.8]{Sch2}). If $X$ is a perfectoid spaces, the same conclusions hold for $\O_{X_\et}^+$ with a similar proof (using \cite[Theorem 6.3]{Sch2} in palce of \cite[Theorem 8.7 and Proposition 8.8]{Sch2}). 
\end{rmk}

Our next goal is to discuss the precise relation between $\O_{X^\diam}^+/p$, $\O_{X^\diam_\qp}^+/p$, and $\O_{X_\et}^+/p$. If one is willing to work in the almost world, then one can quite easily see that each of this sheaves is obtained as the (derived) restriction of the previous one to the smaller site (this essentially boils down to Lemma~\ref{lemma:first-properties-structure-sheaves}). However, to understand the relation between the categories of $\O^+/p$-vector bundles in different topologies, it is essential to understand the relation between these sheaves on the nose. This turns out to be quite subtle and will be discussed in the rest of this and the next sections. 

\begin{lemma}\label{lemma:v-qproet-et} Let $X\in \bf{pAd}_{\Q_p}$ be a pre-adic space over $\Spa(\Q_p, \Z_p)$. Then the natural morphism
\[
\O_{X^\diam_\qp}^+/p \to \lambda_*\left(\O_{X^\diam}^+/p\right),
\]
is an isomorphism. If $X$ is a strongly sheafy adic space over $\Spa(\Q_p, \Z_p)$, then the natural morphisms\footnote{The functor $\mu^{-1} \colon \cal{A}b(X_\et) \to \cal{A}b(X^\diam_\qproet)$ denotes the pullback of sheaves of abelian groups.}
\[
\mu^{-1}\left(\O_{X_\et}^+/p\right)\to \O_{X^\diam_\qp}^+/p,
\]
\[
\O_{X_\et}^+/p \to \bf{R}\mu_*\left(\O_{X^\diam_\qp}^+/p\right)
\]
are isomorphisms as well. 
\end{lemma}
\begin{proof}
    The first result is \cite[Proposition 2.13]{Lucas-Annette}. For the second result, we note that \cite[Lemma 2.7]{Lucas-Annette} ensures\footnote{Strictly speaking, the proof of \cite[Lemma 2.7]{Lucas-Annette} assumes that $X$ is either locally noetherian or perfectoid. However, a similar proof works for any strongly sheafy $X$ if one uses Lemma~\ref{lemma:basis-proetale} in place of \cite[Lemma 2.6]{Lucas-Annette}.} that, for a strongly sheafy adic space $X$, the sheaf $\O_{X^\diam_\qp}^+$ is isomorphic to 
    \[
    \wdh{\O}^+_{X_\qp}\coloneqq \lim_n \mu^{-1}\left(\O_{X_\et}^+/p^n\right).\footnote{The sheaf $\O_{X^\diam_\qp}^+$ is denoted by $\wdh{\O}^+_{X^\diam}$ in \cite{Lucas-Annette}.}
    \]
    Now the quasi-pro\'etale site of a diamond is replete (in the sense of \cite[Definition 3.1.1]{Bhatt-Scholze}) due to \cite[Lemma 1.2]{Lucas-Annette}. Therefore, the fact that $\O_{X_\et}^+$ is $p$-torsionfree and \cite[Proposition 3.1.10]{Bhatt-Scholze} imply that
    \[
    \wdh{\O}^+_{X_\qp}\simeq \bf{R}\lim \mu^{-1}\left(\O_{X_\et}^+/p^n\right) \simeq \wdh{\mu^{-1}(\O_{X_\et}^+)}
    \]
    is the derived $p$-adic completion of $\mu^{-1}\left(\O_{X_\et}^+\right)$. Since $\O_{X^\diam_\qp}^+$ is also $p$-torsionfree by Lemma~\ref{lemma:first-properties-structure-sheaves}, the universal property of derived completion implies that 
    \begin{align*}
    \O_{X^\diam_\qp}^+/p & \simeq \left[\O_{X^\diam_\qp}^+/p\right]\\
    & \simeq [\wdh{\mu^{-1}(\O_{X_\et}^+)}/p]\\
    &\simeq \mu^{-1}\left(\O_{X_\et}^+/p\right).
    \end{align*}
    Finally, \cite[Proposition 14.8 and Lemma 15.6]{Sch2} imply that
    \[
    \O_{X_\et}^+/p \simeq \bf{R}\mu_*\mu^{-1}\left(\O_{X_\et}^+/p\right)\simeq \bf{R}\mu_*\left(\O_{X^\diam_\qp}^+/p\right). \qedhere
    \]
\end{proof}

Our next goal is to compare $\bf{R}\lambda_*\left(\O_{X^\diam}^+/p\right)$ with $\O_{X^\diam_\qp}^+/p$. To do this, we need a number of preliminary results. This will be done in the next section. 

\subsection{$v$-descent for \'etale cohomology of $\O^+/p$}

The main goal of this section is to show that the natural morphism
\[
\O_{X^\diam_\qp}^+/p \to \bf{R}\lambda_*\left(\O_{X^\diam}^+/p\right)
\]
is an isomorphism. However, our argument is a bit roundabout, and we first show that \'etale cohomology complex $\bf{R}\Gamma(X_\et, \O^+_{X_\et}/p)$ satisfies $v$-descent on affinoid perfectoid spaces. Even to formulate this precisely, we will need to use $\infty$-categories as developed in \cite{HTT}. In what follows, we denote by $\cal{D}(\Z)$ the $\infty$-enhancement of the triangulated derived category of abelian groups $\bf{D}(\Z)$. We are also going to slightly abuse the notation and identify a (usual) category $\cal{C}$ with its nerve $\rm{N}(\cal{C})$ (see \cite[\href{https://kerodon.net/tag/002M}{Tag 002M}]{kerodon}) considered as an $\infty$-category. \smallskip

We fix a category $\rm{PerfAff}_{\Q_p}$ of affinoid perfectoid spaces over $\Spa(\Q_p, \Z_p)$. For any morphism $Z\to Y$, we can consider its \v{C}ech nerve $\text{\v{C}}(Z/Y)$ as a simplicial object in $\rm{PerfAff}_{\Q_p}$, i.e. a functor
\[
\text{\v{C}}(Z/Y)\colon \Delta^{\rm{op}} \to \rm{PerfAff}_{\Q_p}.
\] 
More explicitly, the $n$-th term
\[
\text{\v{C}}(Z/Y)_n= Z^{n/Y}
\]
is the $n$-th fiber product of $Z$ over $Y$. Face and degeneracy maps are defined in an evident way. \smallskip

For any functor (in the $\infty$-categorical sense) $\cal{F}\colon \rm{PerfAff}_{\Q_p}^{\rm{op}} \to \cal{D}(\Z)$, we can compose $\cal{F}$ with $\text{\v{C}}(Z/Y)^{\rm{op}}$ to get a cosimplicial object $\text{\v{C}}(Z/Y, \F)$ in $\cal{D}(\Z)$, whose $n$-th term is given by
\[
\text{\v{C}}(Z/Y, \F)_n = \F(Z^{n/Y}).
\]
Now it makes sense to talk about (derived) limits over this cosimplicial object (see \cite[\href{https://kerodon.net/tag/02VY}{Tag 02VY}]{kerodon}
 for more detail).  

\begin{defn}\label{defn:v-sheaves} Let $\cal{F}\colon \rm{PerfAff}_{\Q_p}^{\rm{op}} \to \cal{D}(\Z)$ be a functor (in the $\infty$-categorical sense). 
\begin{itemize}
    \item A morphism $Z\to Y$ is {\it of $\F$-descent} if the natural morphism
    \[
    \F(Y) \to \rm{R}\lim_{n\in \Delta} \text{\v{C}}(Z/Y, \F)_n
    \]
    is an equivalence;
    \item a morphism $Z\to Y$ is {\it of universal $\F$-descent} if, for every morphism $Y'\to Y$, the base change $Z\times_Y Y' \to Y'$ is of $\F$-descent;
    \item $\F$ {\it satisfies $v$-descent} (resp. {\it quasi-pro\'etale descent}) if every $v$-covering\footnote{A morphism $f\colon Z \to Y$ in $\rm{PerfAff}_{\Q_p}$ is a $v$-covering (resp. a quasi-pro\'etale covering) if $f^\diam\colon X^\diam \to Y^\diam$ is so.} (resp. quasi-pro\'etale covering) $X\to Y$ is of (universal\footnote{We note that if every $v$-covering (resp.~quasi-pro\'etale covering) is of $\F$-descent, then they are automatically of universal $\F$-descent because $v$-coverings (resp.~quasi-pro\'etale coverings) are closed under pullbacks in $\rm{PerfAff}_{\Q_p}$.}) $\F$-descent. 
    \item $\F$ is a {\it (derived) $v$-sheaf} if $\F$ satisfies $v$-descent and the natural morphism $\F(Y_1\sqcup Y_2)\to \F(Y_1)\times \F(Y_2)$ is an equivalence for any $Y_1, Y_2\in \rm{PerfAff}_{\Q_p}$. 
\end{itemize}
\end{defn}

\begin{rmk} A functor $\F\colon \rm{PerfAff}_{\Q_p}^{\rm{op}} \to \cal{D}(\Z)$ is a (derived) $v$-sheaf in the sense of Definition~\ref{defn:v-sheaves} if and only if it is a $\cal{D}(\Z)$-valued sheaf on the (big) $v$-site $\rm{PerfAff}_{\Q_p}$ (see \cite[\textsection A.3.3]{Lurie-spectral} for the precise definition). See \cite[Proposition A.3.3.1]{Lurie-spectral} for a detailed proof of this fact. 
\end{rmk}

Our current goal is to give an explicit condition that ensures that a functor $\F$ satisfies $v$-descent. Later on, we will show that \'etale cohomology of the $\O^+/p$-sheaf satisfy this condition. This will be the crucial input to relate $\bf{R}\lambda_*\left(\O_{X^\diam}^+/p\right)$ to $\O_{X^\diam_\qp}^+/p$.

\begin{lemma}\label{lemma:basic-properties}\cite[Lemma 3.1.2]{stacks-sheaves} Let $\cal{F}\colon \rm{PerfAff}_{\Q_p}^{\rm{op}} \to \cal{D}(\Z)$ be a functor (in the $\infty$-categorical sense), and $f\colon Z\to Y$, $g\colon Z' \to Z$ be morphisms in $\rm{PerfAff}_{\Q_p}$. Then

\begin{enumerate}[label=\textbf{(\arabic*)}]
    \item if $f$ has a section, then it is of universal $\F$-descent;
    \item if $f$ is of universal $\F$-descent and $g$ is of universal $\F$-descent, then $f\circ g\colon Z' \to Y$ is of universal $\F$-descent;
    \item If $f\circ g$ of universal $\F$-descent, then $f$ is so. 
\end{enumerate}
\end{lemma}

\begin{lemma}\label{lemma:approximate-section} Let $Y$ be a strictly totally disconnected perfectoid space, and $Z \to Y$ a $v$-cover by an affinoid perfectoid space. Then there is a presentation $Z=\lim_I Z_i\to Y$ as cofiltered limit of affinoid perfectoid spaces over $Y$ such that each $Z_i \to Y$ admits a section.
\end{lemma}
\begin{proof}
    The proof of \cite[Lemma 2.11]{Lucas-Annette} carries over to this case if one uses \cite[Lemma 2.23]{Ben-line-bundles} in place of \cite[Lemma 9.5]{Sch2}.
\end{proof}

\begin{defn} A $v$-covering $Z\to Y$ of affinoid perfectoid spaces is {\it nice} if it can be written $Z=\lim_I Z_i\to Y$ as a cofiltered limit of affinoid perfectoid spaces over $Y$ such that each $Z_i\to Y$ admits a section.
\end{defn}

\begin{rmk}(\cite[Proposition 6.5]{Sch2}) We recall that the category of affinoid perfectoid spaces $\rm{PerfAff}$ admits cofiltered limits. Namely, the limit of the cofiltered system $\{\Spa(S_i, S^+_i)\}$ is given by $\Spa(S, S^+)$ where $S^+$ is the $\varpi$-adic completion of $\colim_I S_i^+$ (for some compatible choice of pseudo-uniformizers $\varpi$) and $S=S^+[\frac{1}{\varpi}]$. In particular, $\rm{PerfAff}_{\Q_p}$ also admits all cofiltered limits. Moreover, one can choose $\varpi=p$ in this case. 
\end{rmk}

\begin{lemma}\label{lemma:descent-criteria} Let $\cal{F}\colon \rm{PerfAff}_{\Q_p}^{\rm{op}} \to \cal{D}(\Z)$ be a functor (in the $\infty$-categorical sense) such that 
\begin{enumerate}
    \item $\F$ is universally bounded below, i.e. there is an integer $N$ such that $\F(Y)\in \cal{D}^{\geq -N}(\Z)$  for any $Y\in \rm{PerfAff}_{\Q_p}$;
    \item $\F$ satisfies quasi-pro\'etale descent;
    \item for an affinoid perfectoid space $Z=\lim_I Z_i$ that is a cofiltered limit of affinoid perfectoid spaces $Z_i$ over $\Spa(\Q_p, \Z_p)$, the natural morphism    
    \[
    \hocolim_I \F(Z_i) \to \F(Z)
    \]
    is an equivalence.
\end{enumerate}
Then $\F$ satisfies $v$-descent.
\end{lemma}
\begin{proof}
    By shifting, we can assume that $\F(Y)\in \cal{D}^{\geq 0}(\Z)$ for any $Y\in \rm{PerfAff}_{\Q_p}$. We pick a $v$-covering $f\colon Z \to Y$ in $\rm{PerfAff}_{\Q_p}$, and wish to show that it is of universal $\F$-descent. We use \cite[Lemma 7.18]{Sch2} to find a quasi-pro\'etale covering $g\colon Y' \to X$ such that $Y'$ is strictly totally disconnected. Then we consider the fiber product
    \[
    \begin{tikzcd}
    Z\times_Y Y'\arrow{r}{g'} \arrow{d}{f'} & Z\arrow{d}{f}\\
    Y' \arrow{r}{g} & Y.
    \end{tikzcd}
    \]
    Lemma~\ref{lemma:basic-properties} implies that $f$ is of universal $\F$-descent if $g$ and $f'$ are so. By assumption, $\F$ satisfies quasi-pro\'etale descent, so $g$ is of universal $\F$-descent. Therefore, it suffices to show that $f'$ is of universal $\F$-descent. \smallskip
    
    We rename $f'$ by $f$ to reduce the question to showing that any $v$-cover $f\colon Z\to Y$ with a strictly totally disconnected $Y$ is of universal $\F$-descent. Furthermore, Lemma~\ref{lemma:approximate-section} implies that  $f$ is nice, so it suffices to show that any nice $v$-cover (with an arbitrary affinoid perfectoid target space) is of universal $\F$-descent. The property of being nice is preserved by arbitrary pullbacks, so it suffices to show that a nice $v$-cover is of $\F$-descent. \smallskip
    
    After all these reductions, we are in the situation of a $v$-cover $f\colon Z\to Y$  that can be written as $Z=\lim_I Z_i \to Y$ a cofiltered limit of affinoid perfectoid spaces over $Y$ admitting a $Y$-section. Lemma~\ref{lemma:basic-properties} ensures that each $f_i\colon Z_i\to Y$ is of $\F$-descent since it has a section. We wish to show that 
    \[
    \F(Y) \to \rm{R}\lim_{n\in \Delta} \text{\v{C}}(Z/Y, \F)_n
    \]
    is an equivalence. By assumption, we know that the natural morphism
    \[
    \hocolim_I \text{\v{C}}(Z_i/Y, \F)_n \to \text{\v{C}}(Z/Y, \F)_n,
    \]
    is an equivalence for any $n\geq 0$. Now the claim follows from the fact that totalization of a coconnective cosimplisial object commutes with filtered (homotopy) colimits (for example, this follows from \cite[Corollary 3.1.13.]{Kubrak-Prikhodko} applied to $\cal{C}=\rm{Fun}(\Delta, \cal{D}(\Z))$, $\cal{D}=\cal{D}(\Z)$, and $F=\rm{hocolim}$). 
\end{proof}

The next goal is to show that the functor (in the $\infty$-categorical sense)
\[
\bf{R}\Gamma_\et(-, \O^+/p)\colon \rm{PerfAff}_{\Q_p}^{\rm{op}} \to \cal{D}(\Z)
\]
\[
Y\in \rm{PerAff}_{\Q_p} \mapsto \bf{R}\Gamma(Y_\et, \O_{Y_\et}^+/p)
\]
is a (derived) $v$-sheaf.

\begin{lemma}\label{lemma:etale-cohomology-quasi-proetale-descent} The functor $\bf{R}\Gamma_\et(-, \O^+/p)\colon \rm{PerfAff}_{\Q_p}^{\rm{op}} \to \cal{D}(\Z)$ satisfies quasi-pro\'etale descent. 
\end{lemma}
\begin{proof}
    By Lemma~\ref{lemma:v-qproet-et}, we have a functorial isomorphism
    \[
    \bf{R}\Gamma(Y_\et, \O_{Y_\et}^+/p) \simeq \bf{R}\Gamma(Y^\diam_{\qp}, \O_{Y^\diam_{\qp}}^+/p).
    \]
    Now quasi-pro\'etale cohomology satisfy quasi-pro\'etale descent by definition. 
\end{proof}

\begin{lemma}\label{lemma:approximation} Let $\{Z_i=\Spa(S_i, S_i^+)\}_{i\in I}$ be a cofiltered system of affinoid perfectoid spaces over $(\Q_p, \Z_p)$, and let $Z_\infty=\lim Z_i$ with morphisms $f_i\colon Z_{\infty}\to Z_i$. Then the natural morphism
\[
\colim_I f_i^{-1}\O_{Z_{i, \et}}^+/p\to \O_{Z_{\infty, \et}}^+/p
\]
is an isomorphism, where $f_i\colon X_\infty \to X_i$ are the natural projection morphisms and $f_i^{-1}$ is the pullback functor on small \'etale topoi. 
\end{lemma}
\begin{proof}
    Note that \cite[Proposition 6.5]{Sch2} implies that $Z_\infty=\Spa(S_\infty, S_\infty^+)$, where $S_\infty^+$ is the $p$-adic completion of $\colim_I S_i^+$ and $S_\infty=S_\infty^+[\frac{1}{p}]$. \smallskip
    
    Now we put $\F$ to be the sheaf $\colim_I f^{-1}_i \O_{Z_i, \et}^+$. Since filtered colimits are exact, we conclude that $\F/p = \colim_I f^{-1}_i \O_{Z_i, \et}^+/p$. Since affinoid perfectoid spaces $U_\infty \to Z_\infty$ \'etale over $Z_\infty$ form a basis of the \'etale site $Z_{\infty, \et}$, it thus suffices to show that the natural morphism
    \[
        \F(U_\infty)/p \to \O_{Z_{\infty, \et}}^+(U_\infty)/p
    \]
    is an isomorpism for any such $U_\infty \to Z_\infty$. Then \cite[Proposition 6.4(iv)]{Sch2} implies that, for some $i_0\in I$, there is an affinoid perfectoid space $U_{i_0}$ with an \'etale morphism $U_{i_0} \to Z_{i_0}$ such that
    \[
    U_{i_0} \times_{Z_{i_0}} Z_\infty \simeq U_{\infty}.
    \]
    For any $j\geq i_0$, we put $U_j\coloneqq U_{i_0}\times_{Z_{i_0}} Z_j$. Since fiber products commute with limits, we see that 
    \[
    U_\infty \simeq \lim_I U_i
    \]
    in the category of affinoid perfectoid spaces. Therefore, $\O_{Z_{\infty, \et}}^+(U_\infty) = \big(\colim_{i\geq i_0} \O_{Z_{i, \et}}^+(U_i)\big)^{\wedge}_{p}$. Now arguing as in  \cite[Proposition 14.9]{Sch2} (or as in \cite[Proposition 5.9.2]{Lei-Fu}), we conclude that $\F(U_\infty) = \colim_{i\geq i_0} \O_{Z_i, \et}^+(U_i)$. Therefore, \cite[\href{https://stacks.math.columbia.edu/tag/05GG}{Tag 05GG}]{stacks-project} ensures that the natural morphism
    \[
    \F(U_\infty)/p \to \O_{Z_{\infty, \et}}^+(U_\infty)
    \]
    is an isomorphism. 
\end{proof}

\begin{cor}\label{cor:cohomology-commute-colimit} Let $Z$ be an affinoid perfectoid space over $\Spa(\Q_p, \Z_p)$, and let $Z=\lim_I Z_i$ be a cofiltered limit of affinoid perfectoid spaces $Z_i$ over $\Spa(\Q_p, \Z_p)$. Then the natural morphism
\[
\hocolim_I \bf{R}\Gamma\left(Z_{i, \et}, \O_{Z_{i, \et}}^+/p\right) \to \bf{R}\Gamma\left(Z_\et, \O_{Z_\et}^+/p\right)
\]
is an equivalence.
\end{cor}
\begin{proof}
    The result is a formal consequence of Lemma~\ref{lemma:approximation} and \cite[Proposition 6.4]{Sch2} (for example, argue as in \cite[Proposition 5.9.2]{Lei-Fu}).
\end{proof}

\begin{cor}\label{cor:v-descent} The functor $\bf{R}\Gamma_\et(-, \O^+/p)\colon \rm{PerfAff}_{\Q_p}^{\rm{op}} \to \cal{D}(\Z)$ is a (derived) $v$-sheaf.
\end{cor}
\begin{proof}
    Clearly, $\bf{R}\Gamma_\et(-, \O^+/p)$ transforms disjoint unions into direct products, so it suffices to show that $\bf{R}\Gamma_\et(-, \O^+/p)$ satisfies $v$-descent. Then it suffices to show that $\bf{R}\Gamma_\et(-, \O^+/p)$ satisfies the conditions of Lemma~\ref{lemma:descent-criteria}. \smallskip
    
    By definition, $\bf{R}\Gamma(Y_\et,\O_{Y_\et}^+/p)\in \cal{D}^{\geq 0}(\Z)$ for any $Y\in \rm{AffPerf}_{\Q_p}$. Lemma~\ref{lemma:etale-cohomology-quasi-proetale-descent} implies that $\bf{R}\Gamma_\et(-, \O^+/p)$ satisfies quasi-pro\'etale descent, and Corollary~\ref{cor:cohomology-commute-colimit} ensures that it satisfies the third condition of Lemma~\ref{lemma:descent-criteria}. Thus, Lemma~\ref{lemma:descent-criteria} guarantees that $\bf{R}\Gamma_\et(-, \O^+/p)$ satisfies $v$-descent. 
\end{proof}


\begin{lemma}\label{lemma:strictly-tot-disconnected-hypercovers} Let $Y\in \rm{PerfAff}_{\Q_p}$, and let $K\to Y^\diam$ be a $v$-hypercover in $Y^\diam_v$ (in the sense of \cite[\href{https://stacks.math.columbia.edu/tag/01G5}{Tag 01G5}]{stacks-project}). Then there is a split (in the sense of \cite[\href{https://stacks.math.columbia.edu/tag/017P}{Tag 017P}]{stacks-project}) $v$-hypercover $K'\to Y$ such that each term $K'_n$ is a strictly totally disconnected perfectoid space, and there is a morphism $K'^\diam \to K$ of augmented (over $Y$) simplicial objects.
\end{lemma}
\begin{proof}
    This is a standard consequence of the fact that any $v$-small sheaf $X$ admits a $v$-covering $f\colon X' \to X$ with a strictly totally disconnected affinoid perfectoid space $X'$. Since this reduction is standard, we only indicate that one should argue as in \cite[\href{https://stacks.math.columbia.edu/tag/0DAV}{Tag 0DAV}]{stacks-project} or \cite[Theorem 4.16]{conrad-hyper} by inductively constructing a split $n$-truncated hypercover $K'$ with a morphism $K' \to K_{\leq n}$. For this inductive step, the crucial input is \cite[Theorem 4.12]{conrad-hyper} that allows us to construct morphisms from a split (truncated) hypercovering. 
\end{proof}

\begin{lemma}\label{lemma:vanishing-v-strictly-totally-disconnected} Let $Y=\Spa(S, S^+)$ be an affinoid perfectoid space over $\Spa(\Q_p, \Z_p)$. Then the natural morphism
\[
\bf{R}\Gamma(Y_\et, \O_{Y_\et}^+/p)\to \bf{R}\Gamma(Y^\diam_v, \O_{Y^\diam}^+/p)
\]
is an isomorphism. 
\end{lemma}
\begin{proof}
  {\it Step~$1$. Compute $\bf{R}\Gamma(Y^\diam_v, \O_{Y^\diam}^+/p)$ ``explicitly'' in terms of hypercovers (see \cite[\href{https://stacks.math.columbia.edu/tag/01G5}{Tag 01G5}]{stacks-project} for a definition of a hypercovering).} Let us denote by $\rm{HC}(Y^\diam)$ the category of all $v$-hypercovers of $Y^\diam$ up to homotopy\footnote{See \cite[\href{https://stacks.math.columbia.edu/tag/01GZ}{Tag 01GZ}]{stacks-project} for the precise definition.}. Likewise, we denote by $\rm{HC}(Y)$ the category of all $v$-hypercovers $Y$ in $\rm{PerfAff}_{\Q_p}$ up to homotopy, and by $\rm{HC}_{\rm{std}}(Y)$ the full subcategory of hypercovers $K \to Y$ such that each $K_n$ is strictly totally disconnected. \smallskip
  
  Then the diamondification functor naturally extends to a fully faithful functor $(-)^\diam\colon \rm{HC}_{\rm{std}}(Y) \to \rm{HC}(Y^\diam)$. Lemma~\ref{lemma:strictly-tot-disconnected-hypercovers} ensures that this functor is cofinal, and so  \cite[\href{https://stacks.math.columbia.edu/tag/01H0}{Tag 01H0}]{stacks-project} implies that, for every integer $i\geq 0$, we have a canonical isomorphism
  \begin{equation}\label{eqn:cech-v-cohomology}
  \rm{H}^i(Y^\diam_v, \O_{Y^\diam}^+/p)\simeq \colim_{K\in \rm{HC}_{\rm{std}}(Y)} \check{\rm{H}}^i(K^\diam_v, \O_{Y^\diam}^+/p),
  \end{equation}
  where $\check{\rm{H}}^i(K^\diam_v, \O_{Y^\diam}^+/p)$ is the \v{C}ech cohomology groups associated to a hypercover $K^\diam \to Y^\diam$ (see \cite[\href{https://stacks.math.columbia.edu/tag/01GU}{Tag 01GU}]{stacks-project}).\smallskip
    
  Moreover, for any affinoid perfectoid space $Z$ with a map $Z\to Y$, we have a natural isomorphism $\O_{Y^\diam}^+/p|_{Z^\diam} \simeq \O_{Z^\diam}^+/p$. Furthermore, Lemma~\ref{lemma:v-qproet-et} ensures that
  \[
  \rm{H}^0(Z^\diam_v, \O^+_{Y^\diam}/p)\simeq \rm{H}^0(Z^\diam_v, \O_{Z^\diam}^+/p)\simeq \rm{H}^0(Z_\et, \O_{Z_\et}^+/p)
  \]
  If $Z=\Spa(S, S^+)$ is strictly totally disconnected, we can simplify it even further by noting that all \'etale sheaves on $Z_\et$ have trivial higher cohomology groups, so 
  \[
  \rm{H}^0(Z_\et, \O_{Z_\et}^+/p) \simeq S^+/pS^+.
  \]

    Combining all these observations, we see that Equation~(\ref{eqn:cech-v-cohomology}) can be simplified to the following form:
    \begin{equation}\label{eqn:cech-v-simplified}
    \rm{H}^i(Y^\diam, \O_{Y^\diam}^+/p)\simeq \colim_{K\in \rm{HC}_{\rm{std}}(Y)} \rm{H}^i\left(S_{0, K}^+/p \to S_{1, K}^+/p \to \dots S^+_{n, K}/p \to \dots \right),
    \end{equation}
    where $K_n=\Spa(S_{n, K}, S_{n, K}^+)$ is a strictly totally disconnected perfectoid space, and the differentials are given by the usual \v{C}ech-type differentials. \smallskip
  
    {\it Step~$2$. $\bf{R}\Gamma_\et(-, \O^+/p)$ satisfies $v$-{\it hyperdescent}}. Firstly, we note that Corollary~\ref{cor:v-descent} and \cite[Proposition A.3.3.1]{Lurie-spectral} ensure that $\bf{R}\Gamma_\et(-, \O_{Y_\et}^+/p)$ is a $\cal{D}(\Z)$-valued $v$-sheaf on $\rm{PerfAff}_{\Q_p}$. Moreover, for any $Y\in \rm{PerfAff}_{\Q_p}$, we know that 
  \[
  \bf{R}\Gamma(Y_\et, \O_{Y_\et}^+/p) \in \cal{D}^{\geq 0}(\Z).
  \]
  Therefore, \cite[Lemma 6.5.2.9]{HTT} implies that $\bf{R}\Gamma_\et(-, \O^+/p)$ is a hypercomplete (derived) $v$-sheaf. Now \cite[Corollary 6.5.3.13]{HTT} implies that any hypercomplete (derived) $v$-sheaf $\F$ (in particular, $\bf{R}\Gamma_\et(-, \O^+/p)$) satisfies hyperdescent, i.e., for any $v$-hypercovering $K\to X$, the natural morphism
  \[
  \F(X) \to \bf{R}\lim_{n\in \Delta} \F(K_n)
  \]
  is an equivalence. \smallskip
  
  {\it Step~$3$. Compute $\bf{R}\Gamma(Y_\et, \O_{Y_\et}^+/p)$ ``explicitly'' in terms of hypercovers.} By Step~$2$, we know that, for any $v$-hypercovering $K\to Y$ in $\rm{PerfAff}_{\Q_p}$, the natural morphism
  \[
  \bf{R}\Gamma(Y_\et, \O_{Y_\et}^+/p) \to \bf{R}\lim_{n\in \Delta} \bf{R}\Gamma(K_{n, \et}, \O_{K_{n, \et}}^+/p)
  \]
  is an isomorphism. Now we assume that each term $K_n=\Spa(S_{n, K}, S_{n, K}^+)$ is strictly totally disconnected, so higher \'etale cohomology of any \'etale sheaf on $K_n$ vanish. Thus, we have
  \[
  \bf{R}\Gamma(K_{n, \et}, \O_{K_{n, \et}}^+/p) \simeq \rm{H}^0(K_{n, \et}, \O_{K_{n, \et}}^+/p)\simeq S_{n, K}^+/pS_{n, K}^+.  
  \]
    Therefore, in this case, the totalization $\bf{R}\lim_{n\in \Delta} \bf{R}\Gamma(K_{n, \et}, \O_{K_{n, \et}}^+/p)$ can be explicitly computed as the \v{C}ech cohomology associated to the hypercovering $K\to Y$. More explicitly, we see that, for every integer $i\geq 0$, we have 
    \begin{equation*}
        \rm{H}^i(Y_\et, \O_{Y_\et}^+/p)\simeq \rm{H}^i\left(S_{0, K}^+/p \to S_{1, K}^+/p \to \dots S^+_{n, K}/p \to \dots \right)
    \end{equation*}
    with standard \v{C}ech-type differentials. Since this formula holds for any $v$-hypercover $K \to Y$ with strictly totally disconnected terms $K_n$, we can pass to the filtered colimit\footnote{The category $\rm{HC}_{\rm{std}}(Y)$ is cofiltered because it is a cofinal category in the filtered category $\rm{HC}(Y^\diam)$. See Step~$1$ and \cite[\href{https://stacks.math.columbia.edu/tag/01GZ}{Tag 01GZ}]{stacks-project} for more detail.} over $\rm{HC}_{\rm{std}}(Y)$ to see that, for every integer $i\geq 0$, 
    \begin{equation}\label{eqn:cech-etale-simplified}
        \rm{H}^i(Y_\et, \O_{Y_\et}^+/p)\simeq \colim_{K\in \rm{HC}_{\rm{std}}(Y)}\rm{H}^i\left(S_{0, K}^+/p \to S_{1, K}^+/p \to \dots S^+_{n, K}/p \to \dots \right),
    \end{equation}
    where $K_n=\Spa(S_{n, K}, S_{n, K}^+)$ is a strictly totally disconnected perfectoid space, and the differentials are given by the usual \v{C}ech-type differentials. \smallskip
    
    {\it Step~$4$. Finish the proof.} Now Equations~(\ref{eqn:cech-v-simplified})~and~(\ref{eqn:cech-etale-simplified}) imply that the natural morphism
    \[
    \rm{H}^i(Y_\et, \O_{Y_\et}^+/p)\to \rm{H}^i(Y^\diam_v, \O_{Y^\diam}^+/p)
    \]
    is an isomorphism for every $i\geq 0$. In other words, the morphism
    \[
    \bf{R}\Gamma(Y_\et, \O_{Y_\et}^+/p) \to \bf{R}\Gamma(Y^\diam_v, \O_{Y^\diam}^+/p)
    \]
    is an isomorphism.
\end{proof}

\begin{cor}\label{cor:v-coh-qpro-coh} Let $X\in \bf{pAd}_{\Q_p}$ be a pre-adic space over $\Spa(\Q_p, \Z_p)$. Then the natural morphism
\[
\O_{X^\diam_\qp}^+/p \to \bf{R}\lambda_*\left(\O_{X^\diam}^+/p\right),
\]
is an isomorphism. 
\end{cor}
\begin{proof}
    Lemma~\ref{lemma:v-qproet-et} ensures that $\O_{X^\diam_\qp}^+/p\to \lambda_*\left(\O_{X^\diam}^+/p\right)$ is an isomorphism. Thus, it suffices to show that 
    \[
    \rm{R}^j\lambda_*\left(\O_{X^\diam}^+/p\right)\simeq 0
    \]
    for $j\geq 1$. Since strictly totally disconnected spaces form a basis for the quasi-pro\'etale topology of any diamond, it suffices to show that 
    \[
    \rm{H}^j(Y^\diam_v, \O_{Y^\diam}^+/p)=0
    \]
    for a totally strictly disconnected perfectoid $Y\to X$ and $j\geq 1$. Lemma~\ref{lemma:vanishing-v-strictly-totally-disconnected} implies that
    \[
    \rm{H}^j(Y^\diam_v, \O_{Y^\diam}^+/p) \simeq \rm{H}^j(Y_\et, \O_{Y_\et}^+/p).
    \]
    Now the latter group vanishes because any \'etale sheaf on a strictly totally disconnected perfectoid space has trivial higher cohomology groups. 
\end{proof}

\begin{cor}\label{cor:O+/p-cohomology} Let $X$ be a strongly sheafy adic space over $\Spa(\Q_p, \Z_p)$. Then the natural morphisms
\[
\bf{R}\Gamma(X, \O_{X_\et}^+/p) \to \bf{R}\Gamma(X_\qp^\diam, \O_{X_\qp^\diam}^+/p) \to \bf{R}\Gamma(X_v^\diam, \O_{X^\diam}^+/p)
\]
are isomorphisms.
\end{cor}
\begin{proof}
    It follows directly from Lemma~\ref{lemma:v-qproet-et} and Corollary~\ref{cor:v-coh-qpro-coh}.
\end{proof}

\begin{cor}\label{cor:no-v-cohomology} Let $X=\Spa(R, R^+)$ be a strictly totally disconnected perfectoid space over $\Spa(\Q_p, \Z_p)$. Then $\rm{H}^i(X^\diam_v, \O_{X^\diam}^+/p) \simeq 0$ for $i\geq 1$, and $\rm{H}^0(X^\diam_v, \O_{X^\diam}^+/p)\simeq R^+/p R^+$.
\end{cor}
\begin{rmk} We emphasize that Corollary~\ref{cor:no-v-cohomology} guarantees the actual vanishing of higher $v$-cohomology groups of $\O^+_{X^\diam}/p$ on a strictly totally disconnected perfectoid space $X$. This is quite surprising for two reasons: this vanishing holds on the nose (without passing to the almost category), the definition of strictly totally disconnected perfectoid spaces, a priori, guarantees vanishing only of \'etale cohomology groups (as opposed to the $v$-cohomology groups). 
\end{rmk}
\begin{proof}
    Corollary~\ref{cor:O+/p-cohomology} implies that
    \[
    \bf{R}\Gamma(X_v^\diam, \O_{X^\diam}^+/p) \simeq \bf{R}\Gamma(X, \O_{X_\et}^+/p).
    \]
    Since $X$ is a strictly totally disconnected space, so any \'etale sheaf has no higher cohomology groups. This implies that $\rm{H}^i(X_v^\diam, \O_{X^\diam}^+/p) \simeq 0$ for $i\geq 1$, and 
    \[
    \rm{H}^0(X_v^\diam, \O_{X^\diam}^+/p) \simeq \rm{H}^0(X, \O_{X_\et}^+)/p \simeq R^+/pR^+. \qedhere
    \]
\end{proof}

As an application, we get the following result: 

\begin{cor}\label{cor:base-change-result} Let $K$ be a $p$-adic non-archimedean field, let $K^+\subset K$ be an open and bounded valuation subring, and let $X$ be a locally noetherian adic space over $\Spa(K, K^+)$. Put $X^\circ\coloneqq X\times_{\Spa(K, K^+)} \Spa(K, \O_K)$. Then the natural morphism
\[
\bf{R}\Gamma(X_v^\diam, \O_{X^\diam}^+/p) \otimes_{K^+/p} \O_K/p \to \bf{R}\Gamma(X_v^{\circ, \diam}, \O_{X^{\circ, \diam}}^+/p)
\]
is an isomorphism. In particular, if $(K, K^+)$ is a perfectoid field pair, then the natural morphism
\[
\bf{R}\Gamma(X_v^\diam, \O_{X^\diam}^+/p) \to \bf{R}\Gamma(X_v^{\circ, \diam}, \O_{X^{\circ, \diam}}^+/p)
\]
is an almost isomorphism. 
\end{cor}
\begin{proof}
    Using the Mayer–Vietoris spectral sequence, we can localize the problem on $X$. Thus, we can assume that $X=\Spa(A, A^+)$ is affinoid. Then we can find a quasi-pro\'etale covering $\Spd(A_\infty, A_\infty^+) \to \Spd(A, A^+)$ such that all fiber products
    \[
    \Spd(A_\infty, A_\infty^+)^{j/\Spd(A, A^+)} = \Spd(B_j, B_j^+)
    \]
    are strictly totally disconnected (affinoid) perfectoid spaces for $j\geq 1$. Thus, Corollary~\ref{cor:no-v-cohomology} implies that
    \[
    \rm{H}^i\left(\Spd(B_j, B_j^+)_v, \O_{X^\diam}^+/p\right) \simeq 0
    \]
    for $i,j \geq 1$, and 
    \[
    \rm{H}^0\left(\Spd(B_j, B_j^+)_v, \O_{X^\diam}^+/p\right) \simeq B_j^+/pB_j^+
    \]
    for $j\geq 1$. Therefore, we can compute $\rm{H}^j(X_v^\diam, \O_{X^\diam}^+/p)$ via the \v{C}ech cohomology groups of the covering $\Spd(A_\infty, A_\infty^+) \to \Spa(A, A^+)$. So we get an isomorphism
    \[
    \rm{H}^i(X_v^\diam, \O_{X^\diam}^+/p) \simeq \rm{H}^i(B_{1}^+/p \to B_2^+/p\to \dots).
    \]
    Now the morphism $\Spa(K, \O_K) \to \Spa(K, K^+)$ is a pro-open immersion, so the fiber products
    \[
    \Spa(B_j, B_j^+) \times_{\Spa(K, K^+)} \Spa(K, \O_K)
    \]
    are strictly totally disconnected affinoid perfectoid spaces represented by\footnote{For example, the proof of Lemma~\ref{lemma:no-integral-closure} goes through without any changes as $\O_K$ is an algebraic localization of $K^+$.}
    \[
    \Spa(B_j, B_j\wdh{\otimes}_{K^+} \O_K).
    \]
    In particular, the same argument as above implies that the $\O^+/p$ cohomology of $X^{\circ, \diam}$ can be computed as follows:
    \[
    \rm{H}^i(X_v^{\circ, \diam}, \O_{X^\diam}^+/p) \simeq \rm{H}^i(B_{1}^+/p\otimes_{K^+/p} \O_K/p \to B_2^+/p \otimes_{K^+/p} \O_K/p \to \dots).
    \]
    Now \cite[Theorem 10.1]{M} implies that $\O_K$ is an algebraic localization of $K^+$, so $\O_K$ is $K^+$-flat. Thus, we get the desired isomorphism
    \[
    \bf{R}\Gamma(X_v^\diam, \O_{X^\diam}^+/p) \otimes_{K^+/p} \O_K/p \to \bf{R}\Gamma(X_v^{\circ, \diam}, \O_{X^{\circ, \diam}}^+/p).
    \]
    If $K$ is perfectoid, the almost isomorphism
    \[
    \bf{R}\Gamma(X_v^\diam, \O_{X^\diam}^+/p) \to \bf{R}\Gamma(X_v^{\circ, \diam}, \O_{X^{\circ, \diam}}^+/p)
    \]
    now follows from Lemma~\ref{lemma:perfectoid-almost}.
\end{proof}

\subsection{$\O^+/p$-vector bundles in different topologies}\label{section:different-topologies}

The main goal of this section is to show that the categories of $v$-, quasi-pro\'etale, and \'etale $\O^+/p$-vector bundles are all equivalent.\smallskip

The results of this section are mostly due to B.\,Heuer. A version of these results has also appeared in \cite{Heuer-G-torsor}. We present a slightly different argument that avoids non-abelian cohomology. We heartfully thank B.\,Heuer for various discussions around these questions and for allowing the author to present a variation of his ideas in this section. \smallskip

For the next definition, we fix an pre-adic space $X$ over $\Spa(\Q_p, \Z_p)$.

\begin{defn}\label{defn:vect-bundles} An $\O_{X^\diam}^+/p$-module $\cal{E}$ (in the $v$-topology on $X^\diam$) is an {\it $\O_{X^\diam}^+/p$-vector bundle} if, there is a $v$-covering $\{X_i \to X^\diam\}_{i\in I}$ such that $\cal{E}|_{X_i} \simeq (\O_{X^\diam}^+/p)|_{X_i}^{r_i}$ for some integers $r_i$. We denote {\it the category of $\O_{X^\diam}^+/p$-vector bundles} by $\rm{Vect}(X^\diam_v, \O_{X^\diam}^+/p)$. 

An $\O_{X^\diam_\qp}^+/p$-module $\cal{E}$ (in the quasi-pro\'etale topology on $X^\diam$) is an {\it $\O_{X^\diam_\qp}^+/p$-vector bundle} if, there is a quasi-pro\'etale covering $\{X_i \to X^\diam\}_{i\in I}$ such that $\cal{E}|_{X_i} \simeq (\O_{X^\diam_\qp}^+/p)|_{X_i}^{r_i}$ for some integers $r_i$. We denote {\it the category of $\O_{X^\diam_\qp}^+/p$-vector bundles} by $\rm{Vect}(X^\diam_\qp, \O_{X^\diam_\qp}^+/p)$.

Let $X$ be a strongly noetherian adic space over $\Spa(\Q_p, \Z_p)$. An $\O_{X_\et}^+/p$-module $\cal{E}$ (in the \'etale topology on $X$) is an {\it $\O_{X_\et}^+/p$-vector bundle} if, there is an \'etale covering $\{X_i \to X\}_{i\in I}$ such that $\cal{E}|_{X_i} \simeq (\O_{X_\et}^+/p)|_{X_i}^{r_i}$ for some integers $r_i$. We denote {\it the category of $\O_{X_\et}^+/p$-vector bundles} by $\rm{Vect}(X_\et, \O_{X_\et}^+/p)$.
\end{defn}

\begin{rmk}\label{rmk:big-small} Note that $\O_{X^\diam}^+/p$-vector bundles are ``big sheaves'', i.e. they are defined on the big $v$-site $X_v^\diam$. In contrast, $\O_{X^\diam_\qp}^+/p$ and $\O_{X_\et}^+/p$-vector bundles are ``small sheaves''; they are defined on the small quasi-pro\'etale $X^\diam_\qproet$ or the small \'etale site $X_\et$ respectively. 
\end{rmk}

The main goal of this section is to show that all these notions of $\O^+/p$-vector bundles are equivalent. 

First, we define the functors between these categories of $\O^+/p$-vector bundles which we later show to be equivalences. For this, we note that Lemma~\ref{lemma:v-qproet-et} implies that $\mu^{-1}\left(\O_{X_\et}^+/p\right) \simeq \O_{X^\diam_\qp}^+/p$. Therefore, $\mu^{-1}$ carries $\O_{X_\et}^+/p$-vector bundles to $\O_{X^\diam_\qp}^+/p$-vector bundles. In particular, it defines the functor
\[
\mu^*\coloneqq \mu^{-1}\colon \rm{Vect}(X_\et, \O_{X_\et}^+/p) \to \rm{Vect}(X^\diam_\qp, \O_{X^\diam_\qp}^+/p).
\]
Unfortunately, the natural morphism $\lambda^{-1}\left(\O_{X_\qp^\diam}^+/p\right) \to \O_{X^\diam_v}^+/p$ is not an isomorphism (see Remark~\ref{rmk:big-small}). Therefore, we define $\lambda^*$ to be the ``$\O^+/p$-module pullback'' functor
\[
\lambda^*\colon \rm{Vect}(X^\diam_\qp, \O_{X^\diam_\qp}^+/p) \to \rm{Vect}(X^\diam_v, \O_{X^\diam}^+/p).
\]
defined by the formula
\[
\lambda^* \cal{E}\coloneqq \lambda^{-1}\cal{E}\otimes_{\lambda^{-1}\O_{X_\qp^\diam}^+/p} \O_{X^\diam}^+/p.
\]

Our goal is to show that both $\lambda^*$ and $\mu^*$ are equivalences. Before we do this, we need some preliminary lemmas: 

\begin{lemma}\label{lemma:cohomology-of-bundles-commute-limit} Let $X$ be a pre-adic space over $\Spa(\Q_p, \Z_p)$, let $\cal{E}$ be an $\O_{X^\diam}^+/p$-vector bundle, and let $Z=\lim Z_i$ be a cofiltered limit of affinoid perfectoid spaces over $X$. Then the natural morphism
\[
\colim_I \rm{H}^0(Z^\diam_{i, v}, \cal{E}) \to \rm{H}^0(Z^\diam_v, \cal{E})
\]
is an isomorphism.
\end{lemma}
\begin{proof}
    Without loss of generality, we can assume that $I$ has a final object $0$. Then, by the sheaf condition and exactness of filtered colimits, it suffices to show the claim $v$-locally on $Z_0$. Therefore, we may assume that $\cal{E}|_{Z^\diam} \simeq (\O_{Z^\diam}^+/p)^d$ is a free vector bundle. The claim then follows from Corollary~\ref{cor:cohomology-commute-colimit}. 
\end{proof}

\begin{cor}\label{cor:cohomology-of-bundles-commute-tilde-limit} Let $X$ be a pre-adic space over $\Spa(\Q_p, \Z_p)$, let $\cal{E}$ be an $\O_{X^\diam}^+/p$-vector bundle, and let $Z \approx \lim_I Z_i \to Z_0$ be an affinoid strongly pro-\'etale morphism of strongly sheafy Tate-affinoid spaces over $X$. Then the natural morphism
\[
\colim_I \rm{H}^0(Z^\diam_{i, v}, \cal{E}) \to \rm{H}^0(Z^\diam_v, \cal{E})
\]
is an isomorphism.
\end{cor}
\begin{proof}
    We can prove the claim $v$-locally on $Z_0^\diam$. Therefore, we can choose a $v$-covering $\widetilde{Z}_0 \to Z_0$ with a strictly totally disconnected perfectoid space $\widetilde{Z}_0$. The proof of Lemma~\ref{lemma:good-properties-affinoid-proetale}\ref{lemma:good-properties-affinoid-proetale-2} ensures that each $\widetilde{Z}_i\coloneqq Z_i\times_{Z_0} \widetilde{Z}_0$ is a strictly totally disconnected affinoid space, and the diamond $(Z\times_{Z_0} \widetilde{Z}_0)^\diam$ is a strictly totally disconnected perfectoid space (of characteristic $p$). Therefore, we see that the natural morphism 
    \[
    \widetilde{Z}\coloneqq \Big((Z\times_{Z_0} \widetilde{Z}_0)^\diam\Big)^\sharp \to Z\times_{Z_0} \widetilde{Z}_0
    \]
    becomes an isomorphism after applying the diamondification functor, and 
    \[
    \widetilde{Z} \simeq \lim_I \widetilde{Z}_i
    \]
    in the category of perfectoid spaces over $X$. Since the question is $v$-local on $Z_0^\diam$ and depends only on the associated diamonds of $Z_i$ and $Z$, we can replace $Z_i$ and $Z$ with $\widetilde{Z}_i$ and $\widetilde{Z}$, respectively, to achieve that each $Z_i$ and $Z$ is an affinoid perfectoid. In this case, the result follows from Lemma~\ref{lemma:cohomology-of-bundles-commute-limit}.
\end{proof}

\begin{lemma}\label{lemma:spread-trivialization} Let $X$ be a pre-adic space over $\Spa(\Q_p, \Z_p)$, let $\cal{E}$ be an $\O_{X^\diam}^+/p$-vector bundle, and let $Z \approx \lim_I Z_i \to Z_0$ be an affinoid strongly pro-\'etale morphism of strongly sheafy Tate-affinoid adic spaces over $X$. If $\cal{E}|_{Z^\diam} \simeq (\O_{Z^\diam}^+/p)^d$ for some integer $d$, then there is $i\in I$ such that $\cal{E}|_{Z_i^\diam} \simeq (\O_{Z_i^\diam}^+/p)^d$.
\end{lemma}
\begin{proof}
    We choose an isomorphism 
    \[
    f\colon (\O_{Z^\diam}^+/p) \xr{\sim} \cal{E}|_{Z^\diam}
    \]
    and wish to descend it to a finite level. \smallskip
    
    {\it Step~$1$. We approximate $f$.} Corollary~\ref{cor:cohomology-of-bundles-commute-tilde-limit} ensures that we can find $i\in I$ and a morphism
    \[
    f_i\colon (\O_{Z_i^\diam}^+/p)^d \to \cal{E}|_{Z^\diam_i} 
    \]
    such that $f_i|_{Z^\diam}=f$. \smallskip
    
    {\it Step~$2$. Approximate $f^{-1}\colon \cal{E}|_{Z^\diam} \to (\O_{Z^\diam}^+/p)^d$.} We note that the dual sheaf
    \[
    \cal{E}^\vee = \cal{H}om_{\O_{X^\diam}^+/p}\left(\cal{E}, \O_{X^\diam}^+/p\right)
    \]
    is also an $\O_{X^\diam}^+/p$-vector bundle. So we can apply the same argument as in Step~$1$ to 
    \[
    (f^{-1})^{\vee}\colon (\O_{Z^\diam}^+/p)^d \to \cal{E}^{\vee}|_{Z^\diam} = \cal{H}om_{\O_{X^\diam}^+/p}\left(\cal{E}, \O_{X^\diam}^+/p\right)|_{Z^\diam}
    \]
    to find (after possible enlarging $i\in I$) a morphism
    \[
    g'_i\colon (\O_{Z_i^\diam}^+/p)^d \to \cal{E}^{\vee}|_{Z_i^\diam}
    \]
    such that ${g'_i}|_{Z^\diam}=(f^{-1})^{\vee}$. By dualizing, we get a morphism
    \[
    g_i\colon \cal{E}|_{Z^\diam_i}\to (\O_{Z_i^\diam}^+/p)^d  
    \]
    such that $g_i|_{Z^\diam}=f^{-1}$.\smallskip
    
    {\it Step~$3$. Show that $f_i\circ g_i=\rm{id}$ and $g_i\circ f_i=\rm{id}$ after possibly enlarging $i\in I$.} We show the first claim, the second is proven in the same way (and even easier). We consider $\rm{id}_{\cal{E}|_{Z^\diam_i}}$ and $f_i\circ g_i$ as sections of the internal Hom sheaf, i.e.
    \[
    \rm{id}_{\cal{E}|_{Z^\diam_i}}, f_i\circ g_i\in \big(\cal{E}nd_{\O_{X^\diam}^+/p}\left(\cal{E}\right)\big)\left(Z^\diam_i\right).
    \]
    
    For brevity, we denote $\cal{E}nd_{\O_{X^\diam}^+/p}\left(\cal{E}\right)$ by $\cal{E}nd$. Note that $\cal{E}nd$ is again an $\O_{X^\diam}^+/p$-vector bundle, and so Lemma~\ref{lemma:cohomology-of-bundles-commute-limit} ensures that
    \[
    \colim_I \cal{E}nd(Z^\diam_i) = \cal{E}(Z^\diam).
    \]
    Thus if $f_i\circ g_i$ and $\rm{id}$ are equal in the colimit, they are equal on $Z^\diam_i$ for some large index $i$. Similarly, $g_i \circ f_i =\rm{Id}$ for some $i\in I$. Therefore, $f_i \colon (\O_{Z_i^\diam}^+/p)^d \xr{\sim} \cal{E}|_{Z^\diam_i}$ is an isomorphism for $i\gg 0$. 
\end{proof}

\begin{lemma}\label{lemma:trivial-on-strictly-totally-disconnected} Let $Y$ be a strictly totally disconnected perfectoid space over $\Spa(\Q_p, \Z_p)$, and let $\cal{E}$ be an $\O_{Y^\diam}^+/p$-vector bundle. Then there is a finite clopen decomposition $Y = \sqcup_{i\in I} Y_i$ such that $\cal{E}|_{Y^\diam_i} \simeq (\O_{Y^\diam_i}^+/p)^{r_i}$ for some integers $r_i$.
\end{lemma}
\begin{proof}
    By assumption, there is a $v$-covering $\{f_j\colon Z_j \to Y^\diam=Y^\flat\}_{j\in J}$ by affinoid perfectoid spaces. Since $Y$ is quasi-compact, we can assume that $J$ is a finite set. \smallskip

    We put $Y'''_j\coloneqq f_j(Z_j) \subset Y^\flat$. This subset is pro-constructible due to \cite[\href{https://stacks.math.columbia.edu/tag/0A2S}{Tag 0A2S}]{stacks-project} and it is generalizing due to \cite[Lemma 1.1.10]{H3}. Therefore, \cite[Lemma 7.6]{Sch2} implies that there is a canonical structure of an affinoid perfectoid space on $Y'''_j$ such $\iota_j\colon Y'''_j \to Y^\flat$ is a pro-(rational subdomain). In particular, $Y'''_j$ is strictly totally disconnected for every $j\in J$ (for example, due to \cite[Lemma 7.19]{Sch2}). \smallskip

    Now Lemma~\ref{lemma:approximate-section} implies that, for each $j\in J$, we can write $Z_j=\lim_{\Lambda_j} Z_{j, \lambda} \to Y'''_j$ as a cofiltered limit of affinoid perfectoid spaces such that $Z_{j, \lambda} \to Y'''_j$ admits a section for each $\lambda\in \Lambda_j$. Therefore, Lemma~\ref{lemma:spread-trivialization} ensures that, for each $j\in J$, there is $\lambda_j\in \Lambda_j$ such that $\cal{E}|_{Z_{j, \lambda_j}}$ is a free $\O^+_{Y^\diam}/p$-vector bundle. Since each $Z_{j, \lambda_j} \to Y'''_j$ admits a section, we can pullback this trivialization along the section to conclude that $\cal{E}|_{Y'''_j}$ is a free $\O^+_{Y^\diam}/p$-vector bundle. \smallskip

    Now we use Lemma~\ref{lemma:spread-trivialization} and the fact that $\iota_j\colon Y'''_j \to Y^\flat$ is a pro-(rational subdomain) to find a rational open subdomain $Y''_j \subset Y^\flat$ such that $Y'''_j\subset Y''_j$ and $\cal{E}|_{Y''_j}$ is a free $\O^+_{Y^\diam}/p$-vector bundle of rank $r(j)$. Finally, for each integer $i$, we put $Y'_i$ to be the union of all $Y''_j$ such that $r(j)=i$ (in other words, it is the union of $Y''_j$ such that $\cal{E}|_{Y''_j}$ is free of rank $i$). Then all $Y'_i$ are disjoint and only finitely many of them are non-empty. Finally, we define $I$ to be the (finite) set of integers such that $Y'_i \neq \varnothing$. Then $Y^\diam = Y^\flat = \sqcup_{i\in I} Y'_i$ is a finite clopen decomposition such that $\cal{E}|_{Y'}$ is finite free. Then the set of untilts $Y_i \coloneqq Y'^\sharp_i \subset (Y^\flat)^\sharp=Y$ does the job. 
\end{proof}

\begin{thm}\label{thm:v-quasi-proetale} (see also \cite{Heuer-G-torsor}) Let $X$ be a pre-adic space over $\Q_p$. Then the functor
\[
\lambda^*\colon \rm{Vect}(X^\diam_\qp, \O_{X^\diam_\qp}^+/p) \to \rm{Vect}(X^\diam_v, \O_{X^\diam}^+/p)
\]
is an equivalence of categories. Furthermore, for any $\O_{X^\diam_\qp}^+/p$-vector bundle $\cal{E}$, the natural morphism
\[
\cal{E}\to \bf{R}\lambda_*\lambda^*\cal{E}
\]
is an isomorphism.
\end{thm}
\begin{proof}
    We start the proof by showing that the natural morphism
    \[
    \cal{E}\to \bf{R}\lambda_*\lambda^*\cal{E}
    \]
    is an isomorphism. The claim is quasi-pro\'etale local, so we can assume that $\cal{E}$ is a trivial $\O_{X^\diam_\qp}^+/p$-vector bundle. In this case, the claim follows from Corollary~\ref{cor:v-coh-qpro-coh}. \smallskip
    
    This already implies full faithfulness of $\lambda^*$. Indeed, it follows from a sequence of isomorphisms:
    \begin{align*}
        \rm{Hom}_{\O_{X^\diam}^+/p}\left(\lambda^*\cal{E}_1, \lambda^*\cal{E}_2\right) & \simeq \rm{Hom}_{\O_{X^\diam_\qp}^+/p}\left(\cal{E}_1, \lambda_*\lambda^*\cal{E}_2\right) \\
        \simeq \rm{Hom}_{\O_{X^\diam_\qp}^+/p}\left(\cal{E}_1, \cal{E}_2\right).
    \end{align*}
    
    To show that $\lambda^*$ is essentially surjective, it is enough to show that, for an $\O_{X^\diam}^+/p$-vector bundle $\cal{E}$, $\lambda_*\cal{E}$ is an $\O_{X^\diam_\qp}^+/p$-vector bundle and the natural morphism
    \[
    \cal{E} \to \lambda^*\lambda_*\cal{E}
    \]
    is an isomorphism. Both claims are quasi-pro\'etale local on $X^\diam$, so we can assume that $X$ is a strictly totally disconnected perfectoid space. Then we can assume that $\cal{E}$ is a free vector bundle due to Lemma~\ref{lemma:trivial-on-strictly-totally-disconnected}. Then $\lambda_*\cal{E}$ is a free $\O_{X^\diam_\qp}^+/p$-vector bundle by Lemma~\ref{lemma:v-qproet-et}. Thus, the natural morphism
    \[
    \cal{E} \to \lambda^*\lambda_*\cal{E}
    \]
    is evidently an isomorphism.
\end{proof}

\begin{lemma}\label{lemma:local-trivialization} Let $X$ be a strongly sheafy adic space over $\Spa(\Q_p, \Z_p)$, and let $\cal{E}$ be an $\O_{X^\diam_\qp}^+/p$-vector bundle (equivalently, an $\O_{X^\diam}^+/p$-vector bundle). Then there is an \'etale covering $\{X'_i \to X\}_{i\in I}$ such that $\cal{E}|_{X'^\diam_i} \simeq (\O^+_{X'^\diam_i}/p)^{r_i}$ for some integers $r_i$.
\begin{proof}
    The question is local on $X$. So we can assume that $X=\Spa(A, A^+)$ for a complete strongly sheafy Tate-Huber pair $(A, A^+)$. Then the result  follows directly from Lemma~\ref{lemma:str-tot-disc-covering}, Theorem~\ref{thm:approximation-etale}, Lemma~\ref{lemma:trivial-on-strictly-totally-disconnected}, and Lemma~\ref{lemma:cohomology-of-bundles-commute-limit}.
\end{proof}
\end{lemma}

\begin{thm}\label{thm:quasi-proetale-etale} (see also \cite{Heuer-G-torsor}) Let $X$ be a strongly sheafy adic space over $\Spa(\Q_p, \Z_p)$. Then the functor
\[
\mu^*\colon \rm{Vect}(X_\et, \O_{X_\et}^+/p) \to \rm{Vect}(X^\diam_\qp, \O_{X^\diam_\qp}^+/p)
\]
is an equivalence of categories. Furthermore, for any $\O_{X^\diam_\et}^+/p$-vector bundle $\cal{E}$, the natural morphism
\[
\cal{E}\to \bf{R}\mu_*\mu^*\cal{E}
\]
is an isomorphism.
\end{thm}
\begin{proof}
    The proof is completely analogous to the proof of Theorem~\ref{thm:v-quasi-proetale} using Lemma~\ref{lemma:local-trivialization} in place of Lemma~\ref{lemma:trivial-on-strictly-totally-disconnected}. 
\end{proof}

\subsection{Trivializing $\O^+/p$-vector bundles}\label{section:trivialize}

We recall that Theorem~\ref{thm:v-quasi-proetale} and Theorem~\ref{thm:quasi-proetale-etale} ensure that the categories of $\O^+/p$-vector bundles in the $v$, quasi-pro\'etale, and \'etale topologies are equivalent. In particular, any $\O^+/p$-vector bundle in the $v$-topology can be trivialized \'etale locally. The main goal of this section is to show that it suffices to consider some very specific \'etale covers. \smallskip

To do this, we need to start with the discussion of $\O^+/p$-vector bundles on some very specific adic spaces. \smallskip

\begin{lemma}\label{lemma:prufer-bijection} Let $X=\Spa(A, A^+)$ be a Tate affinoid pre-adic space such that $A^+$ is a Pr\"ufer domain (in the sense of \cite[Theorem 22.1 and the discussion before it]{Gilmer}). Then the specialization map $\rm{sp}_X \colon |X| \to |\Spf A^+| = |\Spec A^+/A^{\circ\circ}|$ is a homeomorphism.
\end{lemma}
\begin{proof}
    First, (the proof of) \cite[Theorem 8.1.2]{Bhatt-notes} implies that it suffices to show that $\Spec A^+$ does not admit any non-trivial admissible blow-ups. For this, it suffices to show that any finitely generated ideal $I\subset A^+$ is invertible. This is, in turn, one of the defining properties of Pr\"ufer domains (see \cite[Theorem 22.1]{Gilmer}).
\end{proof}

\begin{lemma}\label{lemma:formal-models-of-prufer-domains} Let $X=\Spa(K, K^+)$ be a Tate affinoid adic space such that $K$ is a non-archimedean field and $K^+$ is a Pr\"ufer domain\footnote{We do not assume that $K^+$ is a valuation ring.}. Then the morphism of locally ringed spaces 
\[
\rm{sp}_X\colon (X_{\rm{an}}, \O_{X}^+) \to (\Spf K^+, \O_{\Spf K^+})
\]
is an isomorphism. 
\end{lemma}
\begin{proof}
    Lemma~\ref{lemma:prufer-bijection} implies that $\rm{sp}_X$ is a homeomorphism. Therefore, it suffices to show that $\rm{sp}_X^\# \colon \O_{\Spf K^+} \to \rm{sp}_{X, *} \big(\O_{X}^+\big)$ is an isomorphism. It suffices to show that $\rm{sp}_X^\#(\rm{D}_f)$ is an isomorphism for any $f\in K^+$. \smallskip
    
    Since $K$ is a non-archimedean field, we conclude that $K^\circ = \O_K$ is a rank-$1$ valuation ring. Then we consider the inclusions $K^{\circ \circ} \subset K^+ \subset \O_K$ and fix a pseudo-uniformizer $\varpi\in K^+$. Since $\O_K$ is rank-$1$ valuation ring, we conclude that the induced topologies on $\O_K$ and $K^+$ coincide with the $\varpi$-adic topologies.\smallskip
    
    Now pick $f\in K^+$. If $f\in K^{\circ\circ}$, then $K^+\big[\frac{1}{f}\big]=K$ and so the principal open $\rm{D}(f)$ is empty. In particular,  $\rm{sp}_X^\#(\rm{D}_f)$ is clearly an isomorphism. Therefore, we can assume that $f\in K^+ \smallsetminus K^{\circ\circ}$. Then $K^+\big[\frac{1}{f}\big] \subset \O_K$ is an open subring, so $K^+\big[\frac{1}{f}\big]$ is already complete in the $\varpi$-adic topology. In particular, we conclude that $\O_{\Spf K^+}(\rm{D}(f))=K^+\big[\frac{1}{f}\big]$. Likewise, since $K^+\big[\frac{1}{f}\big] \subset K$ is already complete and integrally closed, we conclude that
    \[
    \big(\rm{sp}_{X, *} \O_{X}^+\big)\big(\rm{D}(f)\big) = \O_X^+\Big(X\big(\frac{1}{f}\big)\Big) = K^+\big[\frac{1}{f}\big].
    \]
    In particular, we conclude that $\rm{sp}_X^\#(\rm{D}_f)$ is an isomorphism. This finishes the proof.
\end{proof}

Now we recall that any locally noetherian analytic adic space $X$ comes with the natural morphism of ringed sites $i_X \colon (X_{\et}, \O_{X_{\et}}^+) \to (X_{\an}, \O_X^+)$. We show that this is an equivalence for some special $X$.

\begin{lemma}\label{lemma:rationaly-field} Let $X=\Spa(K, K^+)$ be a Tate affinoid adic space such that $K$ is a non-archimedean field\footnote{We do not assume that $K^+$ is a valuation ring.}, and let $U\subset X$ be a non-empty rational subdomain. Then $U=\Spa(K, K'^+)$ for some Tate-Huber pair $(K, K'^+)$.
\end{lemma}
\begin{proof}
    Since $U$ is an affinoid space, we only need to show that $\O_X(U)=K$. First, we choose a pseudo-uniformizer $\varpi\in K^+$. Then we note that $K^\circ= \O_K$ is a rank-$1$ valuation ring since $K$ is a non-archimedean field. In particular, we conclude that the induced topologies on both $K^\circ$ and $K^+$ coincide with the $\varpi$-adic topology (and both are complete with respect to this topology). \smallskip
    
    Now we consider the case $U = X\big(\frac{f}{g}\big)$ for some $f, g\in K^\times$. Since $K$ is a field, we can assume that $U = X\big(\frac{1}{f}\big)$ for some $f\in K^\times$. If $f\in K^{\circ\circ}$, then $X\big(\frac{1}{f}\big) = \varnothing$, so we can assume that $f\notin K^{\circ\circ}$. Then we recall that the induced topology on $K^+$ is equal to the $\varpi$-adic topology to conclude that (see \cite[\textsection 1]{H1})
    \[
    \O_X(U) = \Big(K^+\big[\frac{1}{f}\big]^{\wedge}_{(\varpi)}\Big)\big[\frac{1}{\varpi}\big],
    \]
    where $K^+\big[\frac{1}{f}\big]$ is the $K^+$-subalgebra of $K\big[\frac{1}{f}\big]=K$ generated by $\frac{1}{f}$. Since $f\notin K^{\circ\circ}$, we conclude that $K^+\big[\frac{1}{f}\big]\subset \O_K$ is an open subring of $\O_K$. Thus, it is already complete in the $\varpi$-adic topology. So we conclude that 
    \[
    \O_X(U) = \Big(K^+\big[\frac{1}{f}\big]\Big)\big[\frac{1}{\varpi}\big] = K.
    \]
    In general, a rational subdomain $U$ is equal to $X\big(\frac{f_1, \dots, f_n}{g}\big)$ for some $f_1, \dots, f_n, g\in K^\times$. Denote by $U_i$ the rational subdomain $X\big(\frac{f_i}{g}\big)$. Then
    \[
    U = U_1 \cap U_2 \cap \dots \cap U_n.
    \]
    Therefore, we see that  
    \[
    \O_X(U) \simeq  \O_X(U_1) \wdh{\otimes}_K \O_X(U_2) \wdh{\otimes}_K \dots \wdh{\otimes}_K  \O_X(U_n) \simeq K \wdh{\otimes}_K K \wdh{\otimes}_K \dots \wdh{\otimes}_K  K \simeq K. \qedhere
    \]    
\end{proof}

\begin{lemma}\label{lemma:etale-analytic-prufer} Let $X=\Spa(C, C^+)$ be a Tate affinoid adic space such that $C$ is an algebraically closed non-archimedean field\footnote{We do not assume that $C^+$ is a valuation ring.}. Let $\varpi\in C^+$ be a pseudo-uniformizer. Then the morphism of ringed topoi 
\[
i_X\colon (X_{\et}, \O_{X_\et}^+) \to (X_{\an}, \O_X^+)
\]
is an equivalence. In particular, the functor $i_X^{-1}$ induces an equivalence of categories
\[
i_X^{-1} \colon \rm{Vect}(X_{\an}, \O_X^+/\varpi) \xr{\sim} \rm{Vect}(X_{\et}, \O_{X_{\et}}^+/\varpi).
\]
\end{lemma}
\begin{proof}
    We verify conditions (a)-(d) of \cite[Corollary A.5]{H3}. Conditions (a) and (c) are clear. Condition (d) follows from the fact that \'etale maps are open. Indeed, in the notation of \cite[Corollary A.5]{H3}, we can take $I=\{0\}$, $X_0 =X$, $Y_0=\varphi(X)$, and $X_0 \to Y_0$ the map induced by $\varphi$. \smallskip

    Therefore, we are only left to check condition (c) of {\it loc.~cit.} In other words, we need to show that any \'etale morphism $f\colon Y \to X$ admits an \'etale covering $\{Y_i \to X\}_{i\in I}$ such that $Y_i \to X$ is an open immersion. \smallskip

    Without loss of generality, we can assume that $Y=\Spa(A, A^+)$ is affinoid. We can construct $Y_i$ analytically locally on $X$. Lemma~\ref{lemma:formal-models-of-prufer-domains} implies that we can freely replace $X$ with any non-empty open affinoid without changing the assumptions on $X$. Therefore, \cite[Lemma 2.2.8]{H3} implies that we can assume that $f\colon Y \to X$ factors as a composition of an open immersion $j\colon Y \to \ov{Y}^{/X}$ and a finite \'etale morphism $\ov{f}^{/X} \colon \ov{Y}^{/X} \to X$. Since $X$ is strongly noetherian, we conclude that the category $X_{\rm{fet}}$ of finite \'etale adic spaces over $X$ is equivalent to the category $C_{\rm{fet}}$ of finite \'etale $C$-algebras. Since $C$ is algebraically closed, we conclude that $\ov{Y}^{/X} = \sqcup_{i\in I} X_i$ is a disjoint union of a finite number of copies of $X$ ($X_i \simeq X$). Therefore, 
    \[
    \{j_i \colon Y_i\coloneqq X_i \cap Y \to Y\}_{i\in I}
    \]
    gives the desired covering of $Y$. \smallskip

    Now to conclude that $i_X^{-1} \colon \rm{Vect}(X_{\an}, \O_X^+/\varpi) \to \rm{Vect}(X_{\et}, \O_{X_{\et}}^+/\varpi)$ is an equivalence, it suffices to show that $i_X^{-1} \O_X^+/\varpi = \O_{X_\et}^+/\varpi$. Since $i_X^{-1}$ is exact, it suffices to show that $i_X^{-1} \O_X^+ = \O_{X_\et}^+$. For this, it suffices to show that $i_{X, *} \O_{X_\et}^+ = \O_X^+$, but this is evident from the definition. 
\end{proof}

\begin{lemma}\label{lemma:complete-prufer} Let $K$ be a non-archimedean field with an open and bounded valuation subring $K^+\subset K$ and a pseudo-uniformizer $\varpi \in K^+$. Let $K^{\rm{sep}}$ be a separable closure of $K$, and let $\{K_i\}_{i\in I}$ be a filtered system of finite subextensions $K\subset K_i \subset K^{\rm{sep}}$. For each $i\in I$, we put $K_i^+$ to be the integral closure of $K^+$ in $K_i$. Then the completed colimit 
\[
C^+ \coloneqq \big(\colim_{I} K_i^+\big)^{\wedge}_{(\varpi)}
\]
is a Pr\"ufer domain and $C\coloneqq C^+\big[\frac{1}{\varpi}\big]$ is an algebraically closed non-archimedean field. 
\end{lemma}
\begin{proof}
    First, we note that \cite[Lemma 1.6]{H0} implies that $C$ is the usual completion of the topological field $K^{\rm{sep}}$. Therefore, \cite[Proposition 3.4.1/3 and Proposition 3.4.1/6]{BGR} imply that $C$ is algebraically closed. So we only need to show that $C^+$ is a Pr\"ufer domain. \smallskip

    First, we note that \cite[Theorem 22.1]{Gilmer} ensures that $K^+$ is a Pr\"ufer domain. Then \cite[Theorem 22.3]{Gilmer} implies that each $K_i^+$ is a Pr\"ufer domain. Now we note that $\colim_I K_i^+$ is a domain, so \cite[Proposition 22.6]{Gilmer} ensures that it is a Pr\"ufer domain. Then \cite[Theorem 4]{Richman} implies that it suffices to show that every torsion-free $C^+$-module $M$ is flat. Clearly, $M\big[\frac{1}{\varpi}\big]$ is a flat $C=C^+\big[\frac{1}{\varpi}\big]$-module because $C$ is a field. Furthermore, \cite[Chap.~VII, Proposition 4.5]{Cartan-Eilenberg} applied to $A=M$ and $\Lambda=\colim_{I} K_i^+$ implies that $M$ is flat over $\colim_{I} K_i^+$. In particular, $M/\varpi M$ is flat over $C^+/\varpi \simeq (\colim_{I} K_i^+)/\varpi$. Therefore, \cite[Lemma 8.2/1]{B} concludes that $M$ is flat over $C^+$ and finishes the proof. 
\end{proof}

\begin{lemma}\label{lemma:trivialize-projective-modules} In the notation of Lemma~\ref{lemma:complete-prufer}, any finite projective $C^+/\varpi$-module $M$ is free.
\end{lemma}
\begin{proof}
    First, we note that $C^+/\varpi \simeq \colim_{I} (K_i^+/\varpi)$. Therefore, a standard approximation argument reduces the question to showing that every finite projective $K_i^+/\varpi$-module is finite free. Let us denote the residue field of (the rank-$1$ valuation ring) $K_i^\circ = \O_{K_i}$ by $k_i$. Then we observe that $K_i^{\circ\circ} = \rm{rad}(\varpi)$, and thus $K_i^+/\rm{rad}(\varpi) = K_i^+/K_i^{\circ\circ} \subset \O_{K_i}/K_i^{\circ\circ} = k_i$ is a domain. In particular, 
    \[
    |\Spec K^+_i/\varpi | = |\Spec K^+_i/K_i^{\circ\circ}| 
    \]
    is irreducible. Furthermore,  \cite[Ch.VI, \textsection 8.3, Thm.~1 and Ch.VI, \textsection 8.6, Prop.~6]{Bou} imply that each $K_i^+$ is semi-local. In particular, the ring $K^+_i/\varpi$ is semi-local as well. Therefore, \cite[\href{https://stacks.math.columbia.edu/tag/02M9}{Tag 02M9}]{stacks-project} and the above observation that $|\Spec K^+_i/\varpi|$ is irreducible guarantee that any finite projective $K_i^+/\varpi$-module is free.
\end{proof}

\begin{cor}\label{cor:prufer-trivial} In the notation of Lemma~\ref{lemma:complete-prufer}, put $X=\Spa(C, C^+)$. Then any $\O_{X_\et}^+/\varpi$-vector bundle is free.
\end{cor}
\begin{proof}
    Lemma~\ref{lemma:formal-models-of-prufer-domains}, Lemma~\ref{lemma:etale-analytic-prufer}, and Lemma~\ref{lemma:complete-prufer} imply that the category of $\O_{X_\et}^+/\varpi$-vector bundles is equivalent to the category of usual vector bundles on $\Spec C^+/\varpi$. Any such vector bundle is free due to Lemma~\ref{lemma:trivialize-projective-modules}. 
\end{proof}

Now we can prove the main result of this section: 

\begin{thm}\label{thm:precise-trivialization} Let $X$ be a strongly sheafy adic space over $\Spa(\Q_p, \Z_p)$, let $x\in X$ be a point, and let $\cal{E}$ be an $\O_{X^\diam}^+/p$-vector bundle. Then there is an affinoid open subset $x\in U_x \subset X$ and a finite \'etale surjective morphism $\widetilde{U}_x \to U_x$ such that $\cal{E}|_{\widetilde{U}_x} \simeq (\O_{\widetilde{U}_x^\diam}^+/p)^d$ for some integer $d$.
\end{thm}
\begin{proof}
        {\it Step~$1$. The space $X=\Spa(K, K^+)$ is a non-archimedean field $K$ and an open and bounded valuation subring $K^+\subset K$.} In this case, we fix a separable closure $K^{\rm{sep}}$ of $K$. We put $\{K_i\}_{i\in I}$ to be a filtered system of all finite sub-extensions $K\subset K_i \subset K^{\rm{sep}}$, we also put $K_i^+$ to be the integral closure of $K^+$ in $K_i$. Then Lemma~\ref{lemma:complete-prufer} ensures that $C^+\coloneqq \big(\colim_{I} K_i^+\big)^{\wedge}_{(p)}$ is a Pr\"ufer domain and $C\coloneqq C^+\big[\frac{1}{p}\big]$ is an algebraically closed non-archimedean field. Therefore, Corollary~\ref{cor:prufer-trivial} implies that $\cal{E}|_{\Spd(C, C^+)}$ is free. By construction, we have 
        \[
        \Spa(C, C^+) \approx  \lim_I \Spa(K_i, K_i^+) \to \Spa(K, K^+)
        \]
        is an affinoid strongly pro-\'etale morphism. Therefore, Lemma~\ref{lemma:spread-trivialization} implies that there is an index $i\in I$ such that $\cal{E}|_{\Spd(K_i, K_i^+)}$ is free. Now the result follows from the evident observation that $\Spa(K_i, K_i^+) \to \Spa(K, K^+)$ is a surjective finite \'etale morphism. \smallskip

        {\it Step~$2$. General $X$.} We consider the local ring $\O_{X, x}^+$. Then Step~$1$ products a finite separable extension $\wdh{k(x)} \subset K$ such that $\cal{E}|_{\Spd(K, K^+)}$ is free, where $K^+$ is the integral closure of $\wdh{k(x)}^+$ in $K$. \smallskip
        
        Now we note that \cite[Proposition 7.5.5(5)]{Bhatt-notes} implies that $\O_{X, x}^+$ is $p$-adically henselian, and there is a natural isomorphism $\big(\O_{X, x}^+\big)^{\wedge}_{(p)} \simeq \wdh{k(x)}^+$. Therefore, \cite[Proposition 5.4.54]{GR} implies that we can find a finite \'etale morphism $\O_{X, x} \to A$ such that $A \otimes_{\O_{X, x}} \wdh{k(x)} = K$. Since $\O_{X, x}$ is a local ring with residue field $k(x)$, we easily conclude that $\O_{X, x} \to A$ is also faithfully flat. Now we recall that $\O_{X, x} = \colim_{x\in V\subset X} \O_{X} (V)$, so a standard approximation argument implies that we can find an affinoid open $x\in V \subset X$ and a faithfully flat finite \'etale morphism $\O_{X}(V) \to A_V$ such that $A_{V} \otimes_{\O_{X}(V)} \O_{X, x} \simeq A$. \smallskip

        For each affinoid open subset $x\subset W\subset V$, we put $A_W \coloneqq A_V\otimes_{\O_X(V)} \O_X(W)$ and put $A_W^+$ to be the integral closure of $\O_X(W)$ in $A_W$. Then Lemma~\ref{lemma:finite-etale-topology} ensures that $(A_W, A_W^+)$ is a complete Tate-Huber pair for each open affinoid $x\subset W\subset V$. Furthermore, the corresponding morphism $f_W\colon \Spa(A_W, A_W^+) \to W$ is a finite \'etale surjection due to Lemma~\ref{lemma:finite-etale-surjective}. By construction we have that 
        \[
        \Spa(K, K^+) \approx \lim_{x\in W\subset V} \Spa(A_W, A_W^+) \to \Spa(A_V, A_V^+)
        \]
        is an affinoid strongly pro-\'etale morphism, and $\cal{E}|_{\Spd(K, K^+)}$ is free. Therefore, Lemma~\ref{lemma:spread-trivialization} implies that there is an open affinoid subspace $x\in U_x\subset V_x$ such that $\cal{E}|_{\Spd(A_U, A_U^+)}$ is free. Then $\widetilde{U}_x = \Spa(A_U, A_U^+)$ does the job.
\end{proof}

Now we summarize all the results about various $\O^+/p$-vector bundles below:

\begin{cor}\label{cor:different-vector-bundles-equivalent} Let $X$ be a strongly sheafy adic space over $\Spa(\Q_p, \Z_p)$. Then 
\begin{enumerate}[label=\textbf{(\arabic*)}]
    \item the categories $\rm{Vect}(X_{\et}; \O_{X_\et}^+/p)$, $\rm{Vect}(X^\diam_\qp; \O_{X^\diam_\qp}^+/p)$, and $\rm{Vect}(X^\diam_v; \O_{X^\diam}^+/p)$ are equivalent;
    \item These equivalences preserve cohomology groups;
    \item for any $\O_{X^\diam}^+/p$-vector bundle $\cal{E}$ and a point $x\in X$, there exists an open affinoid subspace $x\in U_x\subset X$ and a finite \'etale surjective morphism $\widetilde{U}_x \to U_x$ such that $\cal{E}|_{\widetilde{U}_x^\diam}$ is a free vector bundle.
\end{enumerate}
\end{cor}

\subsection{\'Etale coefficients}\label{section:etale-coefficients}

The main goal of this section is to relate the \'etale and $v$-cohomology groups of \'etale ``overconvergent '' $\O^+/p$-modules. \smallskip

We fix a strongly sheafy adic space over $\Spa(\Q_p, \Z_p)$. Then we note that any \'etale sheaf of $\bf{F}_p$-modules $\F$ on $X$ defines sheaves $\mu^{-1}\F$ and $\lambda^{-1}\mu^{-1}\F$ of $\bf{F}_p$-modules on $X^\diam_\qp$ and $X^\diam_v$ respectively (see Diagram~(\ref{eqn:many-morphisms}). In what follows, we abuse the notation and denote $(\lambda^{-1}\mu^{-1}\F)\otimes_{\bf{F}_p} \O_{X^\diam}^+/p$ simply by $\F\otimes \O_{X^\diam}^+/p$ for any $\F\in \bf{Shv}(X_\et; \bf{F}_p)$. Similarly, we denote by $(\mu^{-1}\F)\otimes_{\bf{F}_p} \O_{X^\diam_\qp}^+/p$ simply by $\F\otimes\O_{X^\diam_\qp}^+/p$. \smallskip

Before we go to the comparison results, we need to discuss some preliminary results on sheaves on pro-finite sets. They turn out to be tied up with overconvergent \'etale sheaves on strictly totally disconnected spaces. 

\begin{defn} Let $S$ be a pro-finite set, a sheaf of $\bf{F}_p$-modules $\F$ is {\it constructible} if there exists a finite decomposition of $S$ into a disjoint union of clopen subsets $S=\bigsqcup_{i=1}^n S_i$ such that $\F|_{S_i}$ is a constant sheaf of finite rank. 
\end{defn}

\begin{lemma}\label{lemma:constructible-serre} Let $S$ be a pro-finite set, and let $f\colon \F \to \G$ be a morphism of constructible sheaves of $\bf{F}_p$-modules. Then $\ker f$ and $\coker f$ are constructible.
\end{lemma}
\begin{proof}
    Since $S$ is pro-finite, each point $s\in S$ admits a clopen subset $s\in U_s\subset S$ such that both $\F|_{U_s}$ and $\G|_{U_s}$ are constant. Since $S$ is quasi-compact, we can find a finite disjoint union decomposition $S=\sqcup_{i=1}^n U_i$ such that both $\F|_{U_i}$ and $\G|_{U_i}$ are constant. So we can assume that both $\F$ and $\G$ are constant. Then it is easy to see that kernel and cokernel are constant as well. 
\end{proof}

\begin{lemma}\label{lemma:approximate-constructible-coherent} Let $S$ be a pro-finite set, and let $\F$ be a sheaf of $\bf{F}_p$-vector spaces. Then $\F \simeq\colim_I \F_i$ for a filtered system of constructible sheaves $\F_i$.
\end{lemma}
\begin{proof}
    We use \cite[\href{https://stacks.math.columbia.edu/tag/093C}{Tag 093C}]{stacks-project} with $\cal{B}$ being the collection of clopen subsets of $S$ to write $\F$ is a filtered colimit of the form
    \[
    \F \simeq \colim_I \coker\left(\bigoplus_{j=1}^m \ud{\bf{F}}_{p,V_j} \to  \bigoplus_{i=1}^n \ud{\bf{F}}_{p,U_i}\right).
    \]
    Now Lemma~\ref{lemma:constructible-serre} implies that each cokernel is constructible finishing the proof. 
\end{proof}

\begin{defn} A sheaf of $\bf{F}_p$-modules $\F$ on $X_\et$ is {\it overconvergent} if, for every specialization $\ov{\eta} \to \ov{s}$ of geometric points of $X$, the specialization map $\F_{\ov{s}} \to \F_{\ov{\eta}}$ is an isomorphism.
\end{defn}

\begin{defn} An \'etale sheaf of $\bf{F}_p$-modules $\F$ on a strictly totally disconnected perfectoid space $X$ is {\it special} if there exists a finite decomposition of $X$ into a disjoint union of clopen subsets $X=\bigsqcup_{i=1}^n X_i$ such that $\F|_{X_i}$ is a constant sheaf of finite rank.
\end{defn}

\begin{lemma}\label{lemma:approximate-by-special} Let $X$ be a strictly totally disconnected perfectoid space, and $\F$ an overconvergent \'etale sheaf of $\bf{F}_p$-modules. Then $\F \simeq\colim_I \F_i$ for a filtered system of special sheaves $\F_i$ of $\bf{F}_p$-modules.
\end{lemma}
\begin{proof}
    Since $X$ is strictly totally disconnected, the \'etale and analytic sites of $X$ are equivalent. So we can argue on the analytic site of $X$. By \cite[Lemma 7.3]{Sch2}, there is a continuous surjection $\pi\colon X \to \pi_0(X)$ onto a pro-finite set $\pi_0(X)$ of connected components. \smallskip
    
    {\it Step~$1$. The natural map $\pi^*\pi_*\F \to \F$ is an isomorphism:} It suffices to check that it is an isomorphism on stalks. Pick a point $x\in X$, then \cite[Lemma 7.3]{Sch2} implies that the connected component of $x$ has a unique closed point $s$. Then after unravelling all definitions, one gets that the map $(\pi^*\pi_*\F)_x \to \F_x$ is naturally identified with the specialization map $\F_s \to \F_x$ that is an isomorphism by the overconvergent assumption. \smallskip
    
    {\it Step~$2$. Finish the proof:} Lemma~\ref{lemma:approximate-constructible-coherent} ensures that $\pi_*\F \simeq \colim_I \G'_i$ is a filtered colimit of constructible sheaves. Since pullback commutes with all colimits, we get $\F \simeq \pi^*\pi_*\F \simeq \colim_I \pi^*\G'_i$. This finishes the proof since each $\G_i\coloneqq \pi^*\G'_i$ is special. 
\end{proof}

\begin{lemma}\label{lemma:quasiproetale-v} Let $X$ be a strongly sheafy adic space over $\Spa(\Q_p, \Z_p)$, and let $\F$ be an overconvergent \'etale sheaf of $\bf{F}_p$-modules. Then the natural morphism
\[
\O_{X^\diam_\qp}^+/p \otimes \F \to \bf{R}\lambda_*(\O_{X^\diam}^+/p\otimes \F)
\]
is an isomorphism.
\end{lemma}
\begin{proof}
    Since strictly totally disconnected spaces form a basis for the quasi-pro\'etale topology on $X^\diam$, it suffices to show that $a$ is an isomorphism on such spaces. Then we can write $\F \simeq \colim_I\F_i$ as a filtered colimit of special sheaves by Lemma~\ref{lemma:approximate-by-special}. One easily checks that $\a$ is a coherent morphism of algebraic topoi, so each $\rm{R}^i\lambda_*(\O_{X^\diam}^+/p \otimes -)$ commutes with filtered colimits by \cite[Exp. VI, Theoreme 5.1]{SGA4_2}. Thus, it suffices to prove the claim for a special $\F$. By definition of a special sheaf, there exists a disjoint decomposition $X=\bigsqcup_{i=1}^n X_i$ into clopen subsets such that $\F|_{X_i}$ is constant of finite rank. Since the question is local on $X^\diam_\qproet$, we can replace $X$ with each $X_i$ to assume that $\F$ is constant. In this case, the claim follows from Corollary~\ref{cor:v-coh-qpro-coh}. 
\end{proof}

\begin{rmk} We do not know if Lemma~\ref{lemma:quasiproetale-v} holds for non overconvergent \'etale sheaves $\F$.
\end{rmk}

Now we discuss the relation between \'etale and quasi-pro\'etale topology. \smallskip

\begin{lemma}\label{lemma:proetale-quasi-coefficients} Let $X$ be a strongly sheafy adic space over $\Spa(\Q_p, \Z_p)$, and let be $\F$ an overconvergent \'etale sheaf of $\bf{F}_p$-modules.  Then the natural morphism
\[
\O_{X_\et}^+/p\otimes \F \to \bf{R}\mu_*(\O_{X^\diam_\qp}^+/p \otimes \F)
\]
is an isomorphism.
\end{lemma}
\begin{proof}
    By Lemma~\ref{lemma:v-qproet-et}, the right hand side is canonically isomorphism to 
    \[
    \bf{R}\mu_*\mu^{-1}\left(\O_{X_\et}^+/p\otimes \F\right).
    \]
    So the result follows from \cite[Proposition 14.8]{Sch2}.
\end{proof}

Now we combine all these results together: 

\begin{lemma}\label{lemma:et-qp-v-overconv-coeff} Let $X$ be a strongly sheafy adic space over $\Spa(\Q_p, \Z_p)$, and $\F$ an overconvergent \'etale sheaf of $\bf{F}_p$-modules on $X$.  Then the natural morphisms
\[
\O_{X_\et}^+/p\otimes \F \to \bf{R}\mu_*\left(\O_{X^\diam_\qp}^+/p \otimes \F\right),
\]
\[
\O_{X^\diam_\qp}^+/p \otimes \F \to \bf{R}\lambda_*\left(\O_{X^\diam}^+/p\otimes \F\right)
\]
are isomorphisms.
\end{lemma}

\subsection{Application: $\O^+$ and $\O$ vector bundles}

In this section, we discuss the relation between $\O^+_{X^\diam}$ and $\O_{X^\diam}$ vector bundles in different topologies. As an application of the methods developed in this section, we reprove a theorem of Kedlaya--Liu saying that, for a perfectoid space $X$, the categories of $\O_{X^\diam}$-vector bundles in the analytic, \'etale, quasi-pro\'etale, and $v$-topologies are all equivalent. To achieve this result, we prove an intermediate claim that the categories of $\O^+_{X^\diam}$-vector bundles in the \'etale, quasi-pro\'etale, and $v$-topologies are equivalent. The results of this section will not be used in the rest of the paper. \smallskip

We define the categories of $v$, quasi-pro\'etale, and \'etale $\O^+$-vector bundles on $X$ (resp. $\O$-vector bundles on $X$) similarly to Definition~\ref{defn:vect-bundles}. \smallskip

We start by understanding the category of $\O^+_{X^\diam}$-torsors on an affinoid perfectoid space $X$.

\begin{lemma}\label{lemma:mod-top-nilpotents} Let $(R, R^+)$ be a perfectoid pair, and let $f\colon (R^+)^d \to (R^+)^d$ be an $R^+$-linear homomorphism such that $\ov{f} \colon (R^+/R^{\circ\circ})^d \to (R^+/R^{\circ\circ})^d$ is an isomorphism. Then $f$ is an isomorphism.
\end{lemma}
\begin{proof}
    Lemma~\ref{lemma:roots-of-pseudounformizer}\ref{lemma:roots-of-pseudounformizer-2} and a standard approximation argument imply that $f\rm{ mod } \varpi \colon (R^+/\varpi)^d \to (R^+/\varpi)^d$ is an isomorphism. Then \cite[\href{https://stacks.math.columbia.edu/tag/0315}{Tag 0315}]{stacks-project} implies that $f$ is surjective, put $K=\ker f$. We note that $K$ is derived $\varpi$-adically complete due to \cite[\href{https://stacks.math.columbia.edu/tag/091U}{Tag 091U}]{stacks-project}. Furthermore, our assumption implies that $K/\varpi K=0$, so \cite[\href{https://stacks.math.columbia.edu/tag/09B9}{Tag 09B9}]{stacks-project} ensures that $K=0$. In particular, $f$ is an isomorphism.
\end{proof}

\begin{lemma}\label{lemma:trivialize-integrally-perfectoid} Let $X=\Spa(R, R^+)$ be an affinoid perfectoid space over $\Spa(\Q_p, \Z_p)$, and let $\cal{E}$ be an $\O_{X^\diam}^+$-vector bundle. If $\cal{E}/p$ is a free $\O_{X^\diam}^+/p$-vector bundle, then $\cal{E}$ is a free $\O_{X^\diam}^+$-vector bundle.
\end{lemma}
\begin{proof}
    In this proof, we put $\m=R^{\circ\circ}$ and always do almost mathematics with respect to this ideal (see Lemma~\ref{lemma:almost-setup}). \smallskip
    
    Now Lemma~\ref{lemma:first-properties-structure-sheaves}\ref{lemma:first-properties-structure-sheaves-1} implies that $\bf{R}\Gamma(X_v^\diam, \cal{E}/p)$ is almost concentrated in degree $0$. Then Lemma~\ref{lemma:first-properties-structure-sheaves}\ref{lemma:first-properties-structure-sheaves-3}, \cite[\href{https://stacks.math.columbia.edu/tag/0A0G}{Tag 0A0G}]{stacks-project}, and Lemma~\ref{lemma:derived-complete-almost-ampl-mod-p} imply that $\bf{R}\Gamma(X_v^\diam, \cal{E})$ is almost concentrated in degree $0$. This implies that 
    \begin{equation}\label{eqn:mod-p}
    \m \otimes_{R^+} \bf{R}\Gamma(X_v^\diam, \cal{E}) = \bf{R}\Gamma(X_v^\diam, \m\otimes_{R^+} \cal{E}) = \bf{R}\Gamma(X_v^\diam, \m \cal{E}) \text{ and }
    \end{equation}
    \begin{equation}\label{eqn:integral} 
    \m \otimes_{R^+} \bf{R}\Gamma(X_v^\diam, \cal{E}/p) = \bf{R}\Gamma(X_v^\diam, \m\otimes_{R^+} \cal{E}/p) = \bf{R}\Gamma(X_v^\diam, \m\cal{E}/p\m\cal{E})
    \end{equation}
    are concentrated in degree $0$. Since $\cal{E}/p$ is trivial, we conclude that $\cal{E}/\m\cal{E}$ is a trivial $\O_{X^\diam}^+/\m$-vector bundle. We choose an isomorphism $\cal{E}/\m\cal{E} \simeq (\O_{X^\diam}^+/\m)^d$ and define a basis
    \[
    e''_1, \dots, e''_r\in \rm{H}^0(X^\diam_v, \cal{E}/\m\cal{E}).
    \]
    Then we consider the short exact sequence 
    \[
    0 \to \m \otimes_{R^+} \big(\cal{E}/p\big) \to \cal{E}/p\m\cal{E} \to \cal{E}/\m\cal{E} \to 0.
    \]
    Now (\ref{eqn:mod-p}) implies that we can lift $e''_i$ to elements $e'_i\in \rm{H}^0(X_v^\diam, \cal{E}/p\m\cal{E})$. Then we use the commutative diagram 
    \[
    \begin{tikzcd}
    0 \arrow{r}& \m\otimes_{R^+} p\cal{E} \arrow{r} \arrow{d} & \cal{E} \arrow{d} \arrow{r} & \cal{E}/p\m \cal{E} \arrow{d} \arrow{r} & 0 \\
    0 \arrow{r} & \cal{E} \arrow{r}{\cdot p} & \cal{E} \arrow{r} & \cal{E}/p \arrow{r} & 0
    \end{tikzcd}
    \]
    and (\ref{eqn:integral}) to conclude that the natural morphism $\rm{H}^0(X^\diam_v, \cal{E}/p\m\cal{E}) \to \rm{H}^0(X^\diam_v, \cal{E}/p)$ factors through $\rm{H}^0(X_v^\diam, \cal{E})/p \subset \rm{H}^0(X^\diam_v, \cal{E}/p)$. This implies that we can lift $e'_i$ to elements $e_i\in \rm{H}^0(X^\diam_v, \cal{E})$. This defines a morphism 
    \[
    \varphi \colon (\O_{X^\diam}^+)^d \to \cal{E}.
    \]
    By construction, $\varphi \text{ mod } \m$ becomes an isomorphism. We wish to show that this implies that $\varphi$ is an isomorphism. This can be checked $v$-locally on $X$, so we can assume that $\cal{E} \simeq (\O_{X^\diam}^+)^d$ and we need to check that $\varphi(X')$ is an isomorphism for any affinoid perfectoid $X' \to X$. Then the result follows directly from Lemma~\ref{lemma:first-properties-structure-sheaves} and Lemma~\ref{lemma:mod-top-nilpotents}. 
\end{proof}

\begin{cor}\label{cor:etale-locally-trivial} Let $X$ be a perfectoid space over $\Spa(\Q_p, \Z_p)$, and let $\cal{E}$ be an $\O_{X^\diam}^+$-vector bundle. Then, for each $x\in X$, there is an open subspace $x\in U_x \subset X$ and a finite \'etale surjective morphism $\widetilde{U}_x \to U_x$ such that $\cal{E}|_{\widetilde{U}_x}$ is trivial.
\end{cor}
\begin{proof}
    This formally follows from Corollary~\ref{cor:different-vector-bundles-equivalent} and Lemma~\ref{lemma:trivialize-integrally-perfectoid}. 
\end{proof}

Now we denote by 
\[
\mu^* = \mu^{-1} \otimes_{\mu^{-1} \O_{X_\et}^+} \O_{X^\diam_\qp}^+ \colon \rm{Vect}(X_\et; \O_{X_\et}^+) \to \rm{Vect}(X_\qp; \O_{X^\diam_\qp}^+) \text{ and}
\]
\[\lambda^* = \lambda^{-1} \otimes_{\lambda^{-1} \O_{X^\diam_\qp}^+} \O_{X^\diam}^+ \colon \rm{Vect}(X^\diam_\qp; \O_{X^\diam_\qp}^+) \to \rm{Vect}(X^\diam_v; \O_{X^\diam}^+)
\]
the pullback functors. 

Now we can show that the categories of $\O^+$-vector bundles in the \'etale, quasi-pro\'etale, and $v$ topologies are all equivalent: 

\begin{thm}\label{thm:integral-Kedlaya-Liu} Let $X$ be a pre-adic space over $\Spa(\Q_p, \Z_p)$. 
\begin{enumerate}[label=\textbf{(\arabic*)}]
    \item Then the functor $\lambda^*\colon \rm{Vect}(X^\diam_\qp; \O_{X^\diam_\qp}^+) \to \rm{Vect}(X^\diam_v; \O_{X^\diam}^+)$ is an equivalence. Furthermore, for any $\O_{X^\diam_\qp}^+$-vector bundle $\cal{V}$, the natural morphism 
    \[
    \cal{V} \to \bf{R}\lambda_*\lambda^* \cal{V}
    \]
    is an isomorphism;
    \item If $X$ is perfectoid, then the functor $\mu^*\colon \rm{Vect}(X_\et; \O_{X_\et}^+) \to \rm{Vect}(X^\diam_\qp; \O_{X^\diam_\qp}^+)$ is an equivalence. Furthermore, for any $\O_{X_\et}^+$-vector bundle $\cal{E}$, the natural morphism
    \[
    \cal{E} \to \bf{R}\mu_*\mu^* \cal{E}
    \]
    is an isomorphism.
\end{enumerate}
\end{thm}
\begin{proof}
    First, we note that the second claim is quasi-pro\'etale local on $X$, so we can assume that $X$ is a perfectoid space. Then the proof is very similar to that of Theorem~\ref{thm:v-quasi-proetale}. We spell out the main steps. \smallskip

    We first show that the natural morphisms 
    \[
    \alpha\colon \O_{X_\et}^+ \to \bf{R}\mu_* \O_{X^\diam_\qp}^+
    \]
    \[
    \beta \colon \O_{X^\diam_\qp}^+ \to \bf{R}\lambda_* \O_{X^\diam}^+
    \]
    are isomorphisms. For this, we note that Remark~\ref{rmk:etale-derived-complete} implies that $\O_{X_\et}^+$, $\O_{X^\diam_\qp}^+$, and $\O_{X^\diam}^+$ are derived $p$-adically complete and $p$-torsionfree. Therefore, we can check that $\alpha$ and $\beta$ are isomorphisms modulo $p$ (in the derived sense). This follows from Theorem~\ref{thm:quasi-proetale-etale} and Theorem~\ref{thm:v-quasi-proetale}. \smallskip

    This formally implies that the maps $\cal{V} \to \bf{R}\lambda_*\lambda^* \cal{V}$ and $\cal{E} \to \bf{R}\mu_*\mu^* \cal{E}$ are isomorphisms. This, in turn, formally implies that $\lambda^*$ and $\mu^*$ are fully faithful. Essential surjectivity of both functors follows from Corollary~\ref{cor:etale-locally-trivial}. 
\end{proof}

\begin{rmk} We note that \cite[Theorem 4.27]{Heuer-G-torsor} gives a much more general version of Theorem~\ref{thm:integral-Kedlaya-Liu}. However, Corollary~\ref{cor:etale-locally-trivial} does not seem to addressed in \cite{Heuer-G-torsor}. 
\end{rmk}

Now we discuss the case of $\O_X$-vector bundles. We first wish to show that any $\O_{X^\diam}$-vector bundle (in the $v$-topology) admits an $\O_{X^\diam}^+$-lattice \'etale locally on $X$. This will be our key tool to reduce questions about $\O$-vector bundles to the case of $\O^+$-vector bundles. For this, we will need a number of preliminary lemmas:

\begin{lemma}\label{lemma:formula-GL_n} Let $A$ be an $f$-henselian ring for some regular element $f\in A$, and let $\wdh{A}$ be its $f$-adic completion. Then the natural morphism 
\[
\rm{GL}_n(A\big[1/f\big])/\rm{GL}_n(A) \to \rm{GL}_n(\wdh{A}\big[1/f\big])/\rm{GL}_n(\wdh{A})
\]
is a bijection. 
\end{lemma}
\begin{proof}
    In this proof, we denote by $\ud{\rm{Vect}}_n(R)$ the groupoid of finite projective $R$-modules of rank-$n$, and by $\rm{Vect}_n(R)$ the set of isomorphism classes of finite projective $R$-modules of rank-$n$.\smallskip

    Now we start the proof. First, \cite[\href{https://stacks.math.columbia.edu/tag/0BNS}{Tag 0BNS}]{stacks-project} ensures that $(A \to \wdh{A}, f)$ is a gluing data. Second, \cite[\href{https://stacks.math.columbia.edu/tag/0BNW}{Tag 0BNW}]{stacks-project} ensures that any finite projective $A$-module is glueable. Therefore, \cite[\href{https://stacks.math.columbia.edu/tag/0BP2}{Tag 0BP2}]{stacks-project} and \cite[\href{https://stacks.math.columbia.edu/tag/0BP6}{Tag 0BP6}]{stacks-project} imply that the following diagram of groupoids
    \[
    \begin{tikzcd}[column sep = 4em]
        \ud{\rm{Vect}}_n(A) \arrow{r}{-\otimes_A \wdh{A}} \arrow{d}{-\otimes_A A_f} & \ud{\rm{Vect}}_n(\wdh{A}) \arrow{d}{-\otimes_{\wdh{A}} \wdh{A}\big[1/f\big]} \\
        \ud{\rm{Vect}}_n(A_f) \arrow{r}{-\otimes_{A_f} \wdh{A}\big[1/f\big]} & \ud{\rm{Vect}}_n(\wdh{A}\big[1/f\big])
    \end{tikzcd}
    \]
    is cartesian. Therefore, we can pass to homotopy groups at the free module $A^n$ to get a long exact sequence of pointed sets:
    \[
    \begin{tikzcd}
    0 \arrow{r} & \rm{GL}_n(A) \arrow{r} & \rm{GL}_n(A\big[1/f\big])\times \rm{GL}_n(\wdh{A}) \arrow{r} & \rm{GL}_n(\wdh{A}\big[1/f\big]) \arrow[dlll] \\
    \rm{Vect}_n(A) \arrow{r}{\alpha} & \rm{Vect}_n(\wdh{A}) \times \rm{Vect}_n(A\big[1/f\big])\arrow{r} & \rm{Vect}_n(\wdh{A}\big[1/f\big]) \arrow{r} & 0.
    \end{tikzcd}
    \]
    To prove the claim, it suffices to show that the fiber of $\alpha$ over the pair of trivial rank-$n$ modules is just a point. This follows from \cite[\href{https://stacks.math.columbia.edu/tag/0D4A}{Tag 0D4A}]{stacks-project} which even implies that the map $\rm{Vect}_n(A) \to \rm{Vect}_n(\wdh{A})$ is a bijection. 
\end{proof}

\begin{defn} Let $X$ be a pre-adic space over $\Spa(\Q_p, \Z_p)$. We define {\it the (pre-)sheaf of invertible matrices} $\rm{GL}_{n, X^\diam}$ on $X_v^\diam$ via the rule
\[
(S \to X^\diam) \mapsto \rm{GL}_{n}\big(\O_{S^\sharp}(S^\sharp)\big)
\]
for any affinoid perfectoid space $S$ over $X^\diam$. 

We define the {\it the (pre-)sheaf of integral invertible matrices} $\rm{GL}^+_{n, X^\diam}$ on $X_v^\diam$ via the rule
\[
(S \to X^\diam) \mapsto \rm{GL}_n\big(\O^+_{S^\sharp}(S^\sharp)\big)
\]
for any affinoid perfectoid space $S$ over $X^\diam$. 
\end{defn}

One easily checks that a $\rm{GL}_{n, X^\diam}$ is a $v$-sheaf since it is isomorphic to the diamond associated to the classical (pre)-adic space $\rm{GL}_{n, \Q_p}\times_{\Q_p} X$. Similarly, $\rm{GL}^+_{n, X^\diam}$ is a $v$-sheaf since it is isomorphic to the diamond associated to the (pre-)adic space $\big(\rm{GL}_{n, \Q_p}\cap \bf{D}^{n^2}_{\Q_p}\big)\times_{\Q_p} X$. \smallskip

For the next definition, we fix a pre-adic spce $X$ over $\Spa(\Q_p, \Z_p)$ and an $\O_{X^\diam}$-vector bundle $\cal{E}$. 

\begin{defn} The {\it sheaf of lattices} $\rm{Latt}_X(\cal{E})$ is the $v$-sheaf defined by the formula
\[
(S \to X^\diam) \mapsto \Big\{\cal{E}^+\in \rm{Vect}(S^{\sharp, \diam}_v; \O^+_{S^{\sharp, \diam}}), \varphi \colon \cal{E}^+\big[\frac{1}{p}\big] \xr{\sim} \cal{E}|_S\Big\}/\rm{isom}
\]
for each affinoid perfectoid $S \to X^\diam$ over $X^\diam$. 
\end{defn}

\begin{lemma}\label{lemma:formula-for-lattices} Let $X$ be a pre-adic space over $\Spa(\Q_p, \Z_p)$, and let $\cal{E}$ be an $\O_{X^\diam}$-vector bundle. Then, $v$-locally on $X^\diam$, the sheaf $\rm{Latt}_X(\cal{E})$ is isomorphic to $\rm{GL}_{n, X^\diam}/\rm{GL}^+_{n, X^\diam}$. 
\end{lemma}
\begin{proof}
    The claim is $v$-local on $X$ by design, so we can assume that $\cal{E}\simeq \O_{X^\diam}^d$. Then we note that $\rm{GL}_{n, X^\diam}$ acts $\cal{E}$, i.e., for any $g\in \rm{GL}_{n, X^\diam}(S)$ we have an isomorphism $g^*\colon \cal{E}_S \to \cal{E}_S$. Therefore, it also acts on $\rm{Latt}_X(\cal{E})$ via the rule
    \[
    g(\cal{E}^+, \varphi\colon \cal{E}^+ \xr{\sim} \cal{E}) = (\cal{E}^+, g^*\circ \varphi).
    \]
    Now let $\cal{E}^+_0 \subset \cal{E}$ be the trivial lattice $(\O_{X^\diam}^+)^d \subset \O_{X^\diam}^d = \cal{E}$, this defines a point $\ell_0 \in \rm{Latt}_X(\cal{E})$. Then the orbit map defines a morphism of sheaves $\alpha\colon \rm{GL}_{n, X^\diam} \to \rm{Latt}_X(\cal{E})$ via the rule
    \[
    g \mapsto g(\ell_0).
    \]
    The stabilizer of $\ell_0$ is equal to $\rm{GL}^+_{n, X^\diam}$, so $\alpha$ factors through an injective morphism
    \[
    \beta \colon \rm{GL}_{n, X^\diam}/\rm{GL}^+_{n, X^\diam} \hookrightarrow \rm{Latt}_X(\cal{E}).
    \]
    So we are only left to show that it is surjective. Let $S\to \rm{Latt}_X(\cal{E})$ be a point corresponding to a lattice $(\cal{E}^+, \varphi)$. We need to show that this point lies in the image of $\beta$ locally in the $v$-topology. By definition, there is a $v$-covering $S' \to S$ such that $\cal{E}^+|_{S'}$ becomes a free $\O_{X^\diam}^+$-vector bundle. But then there is an element $g\in \rm{GL}_{n, X^\diam}(S')$ such that $g(\cal{E}^+)=\cal{E}_0^+|_{S'}$. In particular, $(\cal{E}^+|_{S'}, \varphi|_{S'})$ lies in the image of $\beta(S')$.
\end{proof}

\begin{cor}\label{cor:lattices-are-continuous} Let $X$ be a pre-adic space over $\Spa(\Q_p, \Z_p)$, let $\cal{E}$ be an $\O_{X^\diam}$-vector bundle, and let $Z=\lim Z_i$ be a cofiltered limit of affinoid perfectoid spaces over $X$. Then the natural morphism
\[
\colim_I \rm{Latt}_X(\cal{E})(Z_i^\diam) \to \rm{Latt}_X(\cal{E})(Z^\diam)
\]
is a bijection.
\end{cor}
\begin{proof}
    Let $Z_i=\Spa(R_i, R_i^+)$, we put $R_\infty^+ \coloneqq \colim_I R_i^+$ and denote by $\wdh{R}^+_\infty$ its $p$-adic completion. Then we note that the claim is $v$-local on $Z$, so we can assume that $\cal{E}$ is a free $\O_{X^\diam}$-vector bundle. Then Lemma~\ref{lemma:formula-for-lattices} implies that it suffices to show that
    \[
    \rm{GL}_n(R_\infty^+\big[1/p\big])/\rm{GL}_n(R^+_\infty) \to \rm{GL}_n(\wdh{R}_\infty^+\big[1/p\big])/\rm{GL}_n(\wdh{R}^+_\infty) 
    \]
    is a bijection. This follows directly from Lemma~\ref{lemma:formula-GL_n}, \cite[\href{https://stacks.math.columbia.edu/tag/0ALJ}{Tag 0ALJ}]{stacks-project}, and \cite[\href{https://stacks.math.columbia.edu/tag/0FWT}{Tag 0FWT}]{stacks-project}. 
\end{proof}

\begin{cor}\label{cor:lattices-are-continuous-2} Let $X$ be a pre-adic space over $\Spa(\Q_p, \Z_p)$, let $\cal{E}$ be an $\O_{X^\diam}$-vector bundle, and let $Z \approx \lim_I Z_i \to Z_0$ be an affinoid strongly pro-\'etale morphism of strongly sheafy Tate-affinoid adic spaces over $X$. Then the natural morphism
\[
\colim_I \rm{Latt}_X(\cal{E})(Z_i^\diam) \to \rm{Latt}_X(\cal{E})(Z^\diam)
\]
is a bijection.
\end{cor}
\begin{proof}
This is a direct consequence of Corollary~\ref{cor:lattices-are-continuous} (one can argue as in Corollary~\ref{cor:cohomology-of-bundles-commute-tilde-limit}).
\end{proof}

\begin{cor}\label{cor:lattices-exist-etale-locally} Let $X$ be a strongly sheafy adic space over $\Spa(\Q_p, \Z_p)$, let $\cal{E}$ be an $\O_{X^\diam}$-vector bundle. Then there is an \'etale covering $X' \to X$, an $\O_{X'^\diam}^+$-vector bundle $\cal{E}^+$, and an isomorphism $\cal{E}^+\big[\frac{1}{p}\big]\simeq \cal{E}|_{X'}$. 
\end{cor}
\begin{proof}
    First, we note that we want to show that $\rm{Latt}_{X}(\cal{E})$ admits a section for some \'etale covering $X' \to X$. For this, we can assume that $X$ is an affinoid space. \smallskip
    
    If $X$ is strictly totally disconnected, the result follows from the observation that such a section exists after a $v$-surjection, Lemma~\ref{lemma:approximate-section}, and Corollary~\ref{cor:lattices-are-continuous}. Then the result follows from Lemma~\ref{lemma:str-tot-disc-covering} and Corollary~\ref{cor:lattices-are-continuous-2}. 
\end{proof}

We will later be able to prove a more precise version of Corollary~\ref{cor:lattices-exist-etale-locally}. But, before that, we show that all possible versions of $\O_X$-vector bundles coincide on perfectoid spaces. For this, we denote by $\pi\colon (X_\et, \O_{X_\et}) \to (X_{\an}, \O_{X})$ the natural morphism of ringed sites. We also denote by 
\[
\pi^* = \pi^{-1} \otimes_{\pi^{-1} \O_{X_{\an}}} \O_{X_{\et}} \colon \rm{Vect}(X_{\an}; \O_{X_{\an}}) \to \rm{Vect}(X_{\et}; \O_{X_{\et}}),
\]
\[
\mu^* = \mu^{-1} \otimes_{\mu^{-1} \O_{X_\et}} \O_{X^\diam_\qp} \colon \rm{Vect}(X_\et; \O_{X_\et}) \to \rm{Vect}(X^\diam_\qp; \O_{X^\diam_\qp}), 
\]
\[\lambda^* = \lambda^{-1} \otimes_{\lambda^{-1} \O_{X^\diam_\qp}} \O_{X^\diam} \colon \rm{Vect}(X^\diam_\qp; \O_{X^\diam_\qp}) \to \rm{Vect}(X^\diam_v; \O_{X^\diam})
\]
the natural pullback functors.

\begin{thm}[{\cite[Theorem 3.5.8]{KedLiu2}, \cite[Lemma 17.1.8]{Berkeley-notes}, \cite[Theorem 4.27]{Heuer-G-torsor}}]\label{thm:rational-vector-bundles-on-perfectoids} Let $X$ be a pre-adic space over $\Spa(\Q_p, \Z_p)$.
\begin{enumerate}[label=\textbf{(\arabic*)}]
    \item\label{thm:rational-vector-bundles-on-perfectoids-1} If $X$ is strongly sheafy, then $\pi^*\colon \rm{Vect}(X_{\an}, \O_{X})\to \rm{Vect}(X_{\et}, \O_{X_\et})$ is an equivalence. Furthermore, the natural morphism
    \[
    \cal{L} \to \bf{R}\pi_*\pi^*\cal{L} 
    \]
    is an isomorphism for any $\O_X$-vector bundle $\cal{L}$. Moreover, if $X=\Spa(A, A^+)$ for a strongly sheafy Tate ring $A$, then $\rm{Vect}(X_{\an}, \O_X)$ is equivalent to the category of finitely generated projective $R$-modules;
    \item\label{thm:rational-vector-bundles-on-perfectoids-2} If $X$ is perfectoid, then $\mu^*\colon \rm{Vect}(X_\et; \O_{X_\et}) \to \rm{Vect}(X_\qp; \O_{X^\diam_\qp})$ is an equivalence. Furthermore, the natural morphism
    \[
    \cal{E} \to \bf{R}\mu_*\mu^* \cal{E}
    \]
    is an isomorphism for any $\O_{X_\et}$-vector bundle $\cal{E}$;
    \item\label{thm:rational-vector-bundles-on-perfectoids-3} the functor $\lambda^*\colon \rm{Vect}(X^\diam_\qp; \O_{X^\diam_\qp}) \to \rm{Vect}(X^\diam_v; \O_{X^\diam})$ is an equivalence. Furthermore, the natural morphism
    \[
    \cal{V} \to \bf{R}\lambda_*\lambda^* \cal{V}
    \]
    is an equivalence for any $\O_{X^\diam_\qp}$-vector bundle $\cal{V}$.
\end{enumerate}
\end{thm}
\begin{proof}
    \ref{thm:rational-vector-bundles-on-perfectoids-1} follows from \cite[Theorem 8.2.22(c), (d)]{KedLiu1}. \smallskip

    Part~\ref{thm:rational-vector-bundles-on-perfectoids-3} is quasi-pro\'etale local on $X$, so we can assume that $X$ is an affinoid perfectoid for the purpose of proving \ref{thm:rational-vector-bundles-on-perfectoids-2} and \ref{thm:rational-vector-bundles-on-perfectoids-3}. \smallskip

    Then Theorem~\ref{thm:integral-Kedlaya-Liu} implies that the natural maps $\O_{X_\et}^+ \to \bf{R}\mu_* \O_{X^\diam_\qp}$ and $\O_{X^\diam_\qp} \to \bf{R}\lambda_* \O_{X^\diam}$ are isomorphisms. This formally implies that the maps $\cal{E} \to \bf{R}\mu_*\mu^* \cal{E}$ and $\cal{V} \to \bf{R}\lambda_*\lambda^* \cal{V}$ are isomorphisms. This, in turn, formally implies that $\mu^*$ and $\lambda^*$ are fully faithful. To show essential surjectivity, it suffices to show that any $\O_{X^\diam}$-vector bundle can be trivialized \'etale locally on $X$ (for a perfectoid space $X$). This follows from the combination of Corollary~\ref{cor:lattices-exist-etale-locally} and Corollary~\ref{cor:etale-locally-trivial}.
\end{proof}

Finally, we give a more refined version of Corollary~\ref{cor:lattices-exist-etale-locally}: 

\begin{cor}\label{cor:etale-locally-trivial-precise} Let $X$ be a strongly sheafy adic space over $\Spa(\Q_p, \Z_p)$, and let $\cal{E}$ be an $\O_{X^\diam}$-vector bundle. Then, for each $x\in X$, there is an open subspace $x\in U_x \subset X$, a finite \'etale surjective morphism $\widetilde{U}_x \to U_x$, and an $\O_{\widetilde{U}^\diam_x}^+$-vector bundle $\cal{E}_x^+$ such that $\cal{E}_x^+\big[\frac{1}{p}\big] \simeq \cal{E}|_{\widetilde{U}_x}$.
\end{cor}
\begin{proof}
    Using Corollary~\ref{cor:lattices-are-continuous-2} in place of Lemma~\ref{lemma:spread-trivialization}, we can repeat the argument of Theorem~\ref{thm:precise-trivialization} once we know that  $\cal{E}|_{\Spd(C, C^+)}$ admits a lattice for any morphism $\Spa(C, C^+) \to X$ such that $C$ is an algebraically closed non-archimedean field (and any open, integrally closed, bounded subring $C^+\subset C$).\footnote{The proof of Theorem~\ref{thm:precise-trivialization} ensures that it suffices to prove this claim for a very restrictive class of such pairs $(C, C^+)$, but this is irrelevent for the current proof.} For this, we note that $(C, C^+)$ is a perfectoid pair, so Theorem~\ref{thm:rational-vector-bundles-on-perfectoids} implies that the category of $\O_{\Spd(C, C^+)}$-vector bundles is equivalent to the category of finite dimensional $C$-vector spaces. In particular, any $\cal{E}|_{\Spd(C, C^+)}$ is a free bundle, so it clearly admits an $\O_{X^\diam}^+$-lattice. This finishes the proof.
\end{proof}

\section{Almost coherence of ``$p$-adic nearby cycles''}\label{section:almost-coherent-nearby-cycles}
\subsection{Introduction}

The main goal of this section is to study the ``$p$-adic Nearby Cycles'' sheaves $\mathbf{R}\nu_*{\O}_{X^\diam}^+$ and $\mathbf{R}\nu_*{\O}_{X^\diam}^+/p$ for a rigid-analytic variety $X$. We also study other versions with more general ``coefficients'' including $\O^+/p$-vector bundles in the $v$-topology, and sheaves of the form $\O^+_{X^\diam}/p\otimes \F$ for a Zariski-constructible sheaf $\F$ (see Definition~\ref{defn:zc}). These complexes turn out to be almost coherent; this makes it possible to study {\'etale} cohomology groups of rigid-analytic varieties using (almost) coherent methods on the special fiber.  \smallskip

Before giving precise definitions, let us explain the main motivation to study these sheaves and their relation with \'etale cohomology of rigid-analytic varieties in the simplest case of the ``nearby cycles'' of the sheaf $\O_{X_\et}^+/p$. In \cite{Sch1}, P.\,Scholze proved (\cite[Theorem 5.1]{Sch1}) that  \'etale cohomology groups $\rm{H}^i(X, \mathbf F_p)$ are finite for any smooth, proper rigid-analytic variety $X$ over an algebraically closed $p$-adic non-archimedean field $C$. There are two important ingredients: the almost primitive comparison theorem that says that $\rm{H}^i(X, \O_{X_\et}^+/p)$ are {\it almost} isomorphic to $\rm{H}^i(X, \bf{F}_p) \otimes \O_C/p$, and the almost finiteness of $\rm{H}^i(X,\O_{X_\et}^+/p)$. \smallskip

The proof of the almost finiteness result in \cite{Sch1} uses properness of $X$ in a very elaborate way; first, the proof constructs some specific ``good covering'' of $X$ by affinoids and then shows that there is enough cancelation in the \v{C}ech-to-Derived spectral sequence associated with that covering. We note that all terms of the second page of this spectral sequence are not almost finitely generated, but mysteriously there is enough cancelation so that the terms of the $\infty$-page become almost finitely generated. We refer to \cite[\textsection 5]{Sch1} for details. \smallskip

Our main goal is to give a more geometric and conceptual way to prove this almost finiteness result. Instead of constructing an explicit ``nice'' covering of $X$, we separate the problem into two different problems. We choose an admissible formal $\O_C$-model $\X$ of $X$ and consider the associated morphism of ringed topoi
\[
t\colon (X_{\text{\'et}}, \O_{X_\et}^+) \to (\X_{\text{Zar}}, \O_\X)
\]
that induces the morphism
\[
t\colon (X_{\text{\'et}}, \O_{X_\et}^+/p) \to (\X_{\text{Zar}}, \O_\X/p)=(\X_0, \O_{X_0})
\]
where $\X_0\coloneqq \X\times_{\Spf \O_C} \Spec \O_C/p$ is the mod-$p$ fiber of $\X$. Then one can write 
\[
\bf{R}\Gamma(X, \O_{X_\et}^+/p)\simeq \bf{R}\Gamma\left(\X_0, \bf{R}t_*\O_{X_\et}^+/p\right)
\]
so one can separately study the ``nearby cycles'' complex $\bf{R}t_*\O_{X_\et}^+/p$ and its derived global sections on $\X_0$. \smallskip

The key is that now $\X$ is proper over $\Spf \O_K$ by \cite[Lemma 2.6]{Lutke-proper}\footnote{Strictly speaking, his proof is written under the assumption that $\O_K$ is discretely valued. However, it can be easily generalized to the of a general rank-$1$ complete valuation ring $\O_K$.} (or \cite[Corollary 4.4 and 4.5]{Temkin-proper}). Thus, the Almost Proper Mapping Theorem~\ref{almost-proper-mapping} tells us that, to prove the almost finiteness of $\bf{R}\Gamma(X, \O_{X_\et}^+/p)$, it is sufficient only to show that $\bf{R}t_*\O_{X_\et}^+/p \in \mathbf{D}^+_{acoh}(X)$ has almost coherent cohomology sheaves. \smallskip

The main advantage now is that we can study the ``nearby cycles'' $\bf{R}t_*\O_{X_\et}^+/p$ {\it locally} on the formal model $\X$. So this holds for any admissible formal model and not only for proper ones. Moreover, the only place where we use properness of $X$ in our proof is to get properness of the formal model $\X$ to be able to apply the Almost Proper Mapping Theorem~\ref{almost-proper-mapping}. This allows us to avoid all elaborate spectral sequence arguments while at the same time making the essential part of the proof local on $\X$.  \smallskip

Now we discuss how we prove almost coherentness of $\bf{R}t_*\O_{X_\et}^+/p$. We will actually prove a much stronger almost coherentness statement that holds for all $\O^+/p$-vector bundles in the $v$-topology. However, we find it instructive to discuss the simplest case first. \smallskip

The main idea of the proof is similar to the idea behind the proof \cite[Lemma 5.6]{Sch1}: we reduce the general case to the case of an affine $\X$ with ``nice'' coordinates, where everything can be reduced to almost coherentness of certain continuous group cohomology via perfectoid techniques. In order to make this work, we have to pass to a finer topology that allows towers of finite \'etale morphisms. There are different possible choices, but we find the $v$-topology on the associated diamond $X^\diam$ of $X$ (in the sense of \cite{Sch2}) to be the most convenient for our purposes (see Section~\ref{section:proetale-v-sites} for the detailed discussion). \smallskip

The case of a general $\O_{X^\diam}^+/p$-vector bundle (see Definition~\ref{defn:vect-bundles}) will cause us more trouble; we will use the structure results from Section~\ref{section:trivialize} to handle a general $\O_{X^\diam}^+/p$-vector bundle. The main crucial input that we are going to use is that the category of $\O_{X^\diam}^+/p$-vector bundles is equivalent to the category of \'etale $\O_{X_\et}^+/p$-vector bundles and that, locally, any $\O_{X^\diam}^+/p$-vector bundle can be trivialized by some very particular \'etale covering (see Corollary~\ref{cor:different-vector-bundles-equivalent}).\smallskip

That being said, we can move to the formulation of the main theorem of this section. We refer to Section~\ref{section:proetale-v-sites} for the definition of the quasi-pro\'etale and $v$-topologies on $X^\diam$ for a rigid-analytic variety over a non-archimedean field $K$. These sites come with their ``integral'' structure sheaves $\O^+_{X^\diam}$, $\O^+_{X^\diam_\qp}$, and $\O_{X_\et}^+$ (see Definition~\ref{defn:v-structure-sheaf}) and a diagram of morphisms of ringed sites (see Diagram~\ref{eqn:many-morphisms}~and ~\ref{eqn:many-morphisms-ringed}): 
\begin{equation}\label{eqn:many-morphisms-5}
\begin{tikzcd}
\left(X^\diam_v, \O_{X^\diam}^+\right) \arrow[rrr, bend right, "\nu"] \arrow{r}{\lambda} & \left(X^\diam_\qproet, \O_{X^\diam_\qp}^+\right) \arrow{r}{\mu} & \left(X_\et, \O_{X_\et}^+\right) \arrow{r}{t} & \left(\X_{\rm{Zar}}, \O_\X\right)
\end{tikzcd}
\end{equation}
and the mod-$p$ version
\begin{equation}\label{eqn:many-morphisms-5-2}
\begin{tikzcd}
\left(X^\diam_v, \O_{X^\diam}^+/p\right) \arrow[rrr, bend right, "\nu"] \arrow{r}{\lambda} & \left(X^\diam_\qproet, \O_{X^\diam_\qp}^+/p\right) \arrow{r}{\mu} & \left(X_\et, \O_{X_\et}^+/p\right) \arrow{r}{t} & \left(\X_{\rm{Zar}}, \O_{\X_0}\right)\end{tikzcd}
\end{equation}

If there is any ambiguity in the meaning of $\nu$, we then denote it by $\nu_{\X}$ to explicitly specify the formal model for these functors. \smallskip

Recall that for a perfectoid field $K$, Lemma~\ref{lemma:almost-setup} ensures that the maximal ideal $\m\subset \O_K$ is an ideal of almost mathematics with flat $\widetilde{\m} \simeq \m^2=\m$. For the rest of this section, we fix a $p$-adic perfectoid field $K$, and always do almost mathematics with respect to the ideal $\m$.\smallskip

We are ready to formulate our first main result. We thank B.\,Heuer for the suggestion of trying to prove Theorem~\ref{thm:main-thm-small} for all $\O_{X_\diam}^+/p$-vector bundles. 

\begin{defn}\label{defn:small-mod-p} An $\O^+_{X^{\diam}}/p$-module $\cal{E}$ is a {\it small $\O^+_{X^{\diam}}/p$-vector bundle} if there is a finite \'etale surjective morphism $V\to U$ such that $\cal{E}|_{V^\diam_v} \simeq (\O_{V^\diam}^+/p)^r$ for some integer $r$. 
\end{defn}

\begin{thm}\label{thm:main-thm-small} Let $\X$ be an admissible formal $\O_K$-scheme with adic generic fiber $X$ of dimension $d$ and mod-$p$ fiber $\X_0$, and let $\cal{E}$ be an $\O^+_{X^{\diam}}/p$-vector bundle. Then 
\begin{enumerate}[label=\textbf{(\arabic*)}]
    \item\label{thm:main-thm-small-1} $\bf{R}\nu_*\cal{E}\in \bf{D}^+_{qc, acoh}(\X_0)$ and $(\bf{R}\nu_*\cal{E})^a\in \bf{D}^{[0, 2d]}_{acoh}(\X_0)^a$;
    \item\label{thm:main-thm-small-22} if $\X=\Spf A$ is affine, then the natural map 
    \[
    \widetilde{\rm{H}^i\left(X^\diam_v, \cal{E} \right)} \to \rm{R}^i\nu_*\left(\cal{E}\right)
    \]
    is an isomorphism for every $i\geq 0$;
    \item\label{thm:main-thm-small-3} the formation of $\rm{R}^i\nu_*(\cal{E})$ commutes with \'etale base change, i.e., for any \'etale morphism $\mf \colon \Y \to \X$ with adic generic fiber $f\colon Y\to X$, the natural morphism 
    \[
    \mf^*_0 \left(\rm{R}^i\nu_{\X, *}(\cal{E}) \right)\to \rm{R}^i\nu_{\Y, *}\left(\cal{E}|_{Y^\diam}\right)
    \]
    is an isomorphism for any $i\geq 0$;
    \item\label{thm:main-thm-small-4} if $\X$ has an open affine covering $\X=\bigcup_{i\in I} \sU_i$ such that $\cal{E}|_{(\sU_{i, K})^\diam}$ is small, then
    \[
    \left(\bf{R}\nu_{*}\cal{E}\right)^a \in \bf{D}^{[0, d]}_{acoh}(\X_0)^a;
    \]
    \item\label{thm:main-thm-small-5} there is an admissible blow-up $\X'\to \X$ such that $\X'$ has an open affine covering $\X'=\bigcup_{i\in I} \sU_i$ such that $\cal{E}|_{(\sU_{i, K})^\diam}$ is small. 
    
    In particular, there is a cofinal family of admissible formal models $\{\X'_i
    \}_{i\in I}$ of $X$ such that 
    \[
    \left(\bf{R}\nu_{\X'_i, *}\cal{E}\right)^a\in \bf{D}^{[0, d]}_{acoh}(\X'_{i, 0})^a.
    \]
    for each $i\in I$. 
\end{enumerate}
\end{thm}

\begin{rmk} We refer to Definition~\ref{defn:derived-almost-category-1} and Definition~\ref{defn:derived-almost-category-2} for the precise definition of all derived categories appearing in Theorem~\ref{thm:main-thm-small}. In order to avoid any confusion, we explicate that the expression $\left(\bf{R}\nu_{*}\cal{E}\right)^a \in \bf{D}^{[0, d]}_{acoh}(\X_0)^a$ means that the complex $\left(\bf{R}\nu_{*}\cal{E}\right)$ is almost concentrated in degree $[0, d]$ and each of its cohomology sheaves is almost coherent.
\end{rmk}

\begin{rmk} We note that Theorem~\ref{thm:main-thm-small}~\ref{thm:main-thm-small-1} implies that the nearby cycles $\bf{R}\nu_*\cal{E}$ is quasi-coherent on the nose (as opposed to being almost quasi-coherent). This is quite unexpected to the author since all previous results on the cohomology groups of $\O^+/p$ were only available in the almost category. 
\end{rmk}

\begin{rmk} If $K=C$ is algebraically closed, the proof gives a non-almost version of cohomological bound. Namely, we see that
\[
\bf{R}\nu_*\cal{E} \in \bf{D}^{[0, 2d]}_{acoh}(\X_0).
\]
However, we do not know if $\bf{R}\nu_*\cal{E}$ is concentrated in degrees $[0, d]$ on the nose (for a cofinal family of formal models).  
\end{rmk}

\begin{rmk} Ofer Gabber has informed the author that he knows an example of a smooth affinoid rigid-analytic variety $X$, a formal model $\X$, and an $\O_{X^\diam}^+/p$-vector bundle $\cal{E}$ such that $\bf{R}\nu_{\X, *}\cal{E}$ is {\it not} almost concentrated in degrees $[0, d]$.
\end{rmk}


One can prove a slightly more precise version in case $\cal{E}$ is equal to the tensor product of a Zariski-constructible \'etale sheaf of $\bf{F}_p$-modules and $\O_{X^\diam}^+/p$. 

\begin{defn}\label{defn:zc}\cite{Hansen-vanishing} An \'etale sheaf $\F$ of $\bf{F}_p$-modules is a {\it local system} if it is a locally constant sheaf with finite stalks.

An \'etale sheaf $\F$ of $\bf{F}_p$-modules is {\it Zariski-constructible} if there is a locally finite stratification $X= \bigsqcup_{i\in I} Z_i$ into Zariski locally closed subspaces $Z_i$ such that $\F|_{Z_i}$ is a local system.

The category $\bf{D}_{zc}(X; \bf{F}_p)$ is a full subcategory of $\bf{D}(X_\et; \bf{F}_p)$ consisting of objects with Zariski-constructible cohomology sheaves.
\end{defn}

\begin{rmk}\label{rmk:zc-overconvergent} Any Zariski-constructible sheaf $\F$ is overconvergent, i.e., for any morphism $\ov{\eta} \to \ov{s}$ of geometric points in $X_\et$, the specialization map $\F_{\ov{s}} \to \F_{\ov{\eta}}$ is an isomorphism.
\end{rmk}

Note that any sheaf of $\bf{F}_p$-modules on $X_\et$ can be treated as a sheaf on any of the sites $X_v^\diam$, $X^\diam_\qproet$, or $X_\proet$ via the pullback functors along the morphisms in Diagram~\ref{eqn:many-morphisms-5}. In what follows, we abuse the notation and implicitly treat a sheaf $\F$ as a sheaf on any of those sites. We also denote the tensor product $\F\otimes_{\bf{F}_p} \O_X^+/p$ simply by $\F\otimes \O_X^+/p$ in what follows.  \smallskip

Now we discuss an integral version of Theorem~\ref{thm:main-thm-small}.

\begin{thm}\label{main-thm} Let $\X$ be an admissible formal $\O_K$-scheme with adic generic fiber $X$ of dimension $d$ and mod-$p$ fiber $\X_0$, and $\F \in \bf{D}^{[r, s]}_{zc}(X; \bf{F}_p)$. 
Then 
\begin{enumerate}[label=\textbf{(\arabic*)}]
    \item\label{thm:main-thm-1} there is an isomorphism $\bf{R}t_*\left(\F\otimes \O_{X_\et}^+/p\right) \simeq \bf{R}\nu_*\left(\F\otimes \O_{X^\diam}^+/p\right)$;
    \item\label{thm:main-thm-22} $\bf{R}\nu_*\left(\F\otimes \O_{X^\diam}^+/p\right) \in \bf{D}^+_{qc, acoh}(\X_0)$, and $\bf{R}\nu_*\left(\F\otimes \O_{X^\diam}^+/p\right)^a \in \bf{D}^{[r, s+d]}_{acoh}(\X_0)^a$;
    \item\label{thm:main-thm-3} if $\X=\Spf A$ is affine, then the natural map 
    \[
    \widetilde{\rm{H}^i\left(X^\diam_v, \F\otimes \O_{X^\diam}^+/p \right)} \to \rm{R}^i\nu_*\left(\F\otimes \O_{X^\diam}^+/p\right)
    \]
    is an isomorphism for every $i\geq 0$;
    \item\label{thm:main-thm-4} the formation of $\rm{R}^i\nu_*\left(\F\otimes \O_{X^\diam}^+/p\right)$ commutes with \'etale base change, i.e., for any \'etale morphism $\mf \colon \Y \to \X$ with adic generic fiber $f\colon Y\to X$, the natural morphism 
    \[
    \mf^*_0 \left(\rm{R}^i\nu_{\X, *}\left(\F\otimes \O_{X^\diam}^+/p\right) \right)\to \rm{R}^i\nu_{\Y, *}\left(f^{-1}\F\otimes \O_{Y^\diam}^+/p\right)
    \]
    is an isomorphism for any $i\geq 0$.
\end{enumerate}
\end{thm}

\begin{defn}\label{defn:small-integrally} An $\O^+_{X^{\diam}}$-vector bundle $\cal{E}$ is a {\it small $\O^+_{X^\diam}$-vector bundle} if $\cal{E}/p\cal{E}$ is a small $\O_{X^\diam}^+/p$-vector bundle (see Definition~\ref{defn:small-mod-p}). 
\end{defn}

\begin{thm}\label{thm:main-thm-integral} Let $\X$ be an admissible formal $\O_K$-scheme with adic generic fiber $X$ of dimension $d$, and let $\cal{E}$ be an $\O^+_{X^{\diam}}$-vector bundle. 
Then 
\begin{enumerate}[label=\textbf{(\arabic*)}]
    \item\label{thm:main-thm-integral-1} $\bf{R}\nu_*\cal{E}\in \bf{D}^+_{qc, acoh}(\X)$ and $(\bf{R}\nu_*\cal{E})^a\in \bf{D}^{[0, 2d]}_{acoh}(\X)^a$;
    \item\label{thm:main-thm-integral-22} if $\X=\Spf A$ is affine, then the natural map 
    \[
    \rm{H}^i\left(X^\diam_v, \cal{E} \right)^{\Updelta} \to \rm{R}^i\nu_*\left(\cal{E}\right)
    \]
    is an isomorphism for every $i\geq 0$;
    \item\label{thm:main-thm-integral-3} the formation of $\rm{R}^i\nu_*(\cal{E})$ commutes with \'etale base change, i.e., for any \'etale morphism $\mf \colon \Y \to \X$ with adic generic fiber $f\colon Y\to X$, the natural morphism 
    \[
    \mf^* \left(\rm{R}^i\nu_{\X, *}(\cal{E}) \right)\to \rm{R}^i\nu_{\Y, *}\left(\cal{E}|_{Y^\diam}\right)
    \]
    is an isomorphism for any $i\geq 0$;
    \item\label{thm:main-thm-integral-4} if $\X$ has an open affine covering $\X=\bigcup_{i\in I} \sU_i$ such that $\cal{E}|_{(\sU_{i, K})^\diam}$ is small, then
    \[
    \left(\bf{R}\nu_{*}\cal{E}\right)^a \in \bf{D}^{[0, d]}_{acoh}(\X)^a;
    \]
    \item\label{thm:main-thm-integral-5} there is an admissible blow-up $\X'\to \X$ such that $\X'$ has an open affine covering $\X'=\bigcup_{i\in I} \sU_i$ such that $\cal{E}|_{(\sU_{i, K})^\diam}$ is small. 
    
    In particular, there is a cofinal family of admissible formal models $\{\X'_i
    \}_{i\in I}$ of $X$ such that 
    \[
    (\bf{R}\nu_{\X'_i, *}\cal{E})^a\in \bf{D}^{[0, d]}_{acoh}(\X'_{i})^a.
    \]
    for each $i\in I$. 
\end{enumerate}
\end{thm}

\begin{rmk} We refer to Definition~\ref{defn:almost-category-formal}  for the precise definition of all derived categories appearing in Theorem~\ref{thm:main-thm-integral}. \end{rmk}

\begin{rmk} One can also prove a version of Theorem~\ref{thm:main-thm-integral} for Zariski-constructible $\Z_p$-sheaves in the sense of \cite[Definition 3.32]{Bhatt-Hansen}. However, we prefer not to do this here as it does not require new ideas but instead complicates the notation. 
\end{rmk}

For the version of Theorem~\ref{thm:main-thm-integral} with the pro-\'etale site $X_\proet$ as defined in \cite{Sch1} and \cite{Sch-err}, see Theorem~\ref{thm:integral-structure-sheaf}. \smallskip

The rest of the paper is devoted to proving Theorem~\ref{main-thm}, Theorem~\ref{thm:main-thm-small}, and Theorem~\ref{thm:main-thm-integral} and discussing their applications. We have decided to work entirely in the $v$-site of $X^\diam$ because it is quite flexible for different types of arguments (e.g. proper descent, torsors under pro-finite groups, etc.). However, most of the argument can be done using the more classical pro-\'etale site defined in \cite{Sch1}. However, it is crucial to use the theory of diamonds to get an almost cohomological bound on $\bf{R}\nu_* \cal{E}$ for non-smooth $X$, and it also seems difficult to justify that the sheaves $\rm{R}^i\nu_{*}\cal{E}$ are quasi-coherent without using (at least) quasi-pro\'etale topology.

\subsection{Digression: geometric points}

In this section, we discuss preliminary results that will be used both in the proof of Theorem~\ref{main-thm} and in deriving applications from it. \smallskip 

We start the section by recalling some definitions.

\begin{defn}\cite[2.1.4]{Temkin-degree} An extension of non-archimedean fields\footnote{Recall that non-archimedean fields are complete by our convention.} $K\subset L$ is {\it topologically algebraic} if the algebraic closure of $K$ in $L$ is dense in $L$. Equivalently, $K\subset L$ is topologically algebraic if $L$ is a non-archimedean subfield of $\wdh{\ov{K}}$.
\end{defn}

\begin{lemma}\label{lemma:extension-fields}
\begin{enumerate}[label=\textbf{(\arabic*)}]
    \item\label{lemma:extension-fields-1} Let $K \subset L$ and $L\subset M$ be two topologically algebraic extensions of non-archimedean fields. Then $K\subset M$ is also topologically algebraic.
    \item\label{lemma:extension-fields-2} Let 
    \[
    \begin{tikzcd}
    N &\arrow{l} L \\
    M \arrow{u} & \arrow{l} \arrow{u} K
    \end{tikzcd}
    \]
    be a commutative diagram of non-archimedean fields such that $LM$ is dense in $N$ and $K\subset L$ is topologically algebraic. Then $M \subset N$ is also a topologically algebraic extension. 
\end{enumerate}
\end{lemma}
\begin{proof}
    \ref{lemma:extension-fields-1}: We know that $L\subset \wdh{\ov{K}}$ and $M\subset \wdh{\ov{L}}$ since both extensions are topologically algebraic. Since $\wdh{\ov{L}}$ is already algebraically closed, we conclude that $M\subset \wdh{\ov{L}} \subset \wdh{\ov{K}}$. \smallskip
    
    \ref{lemma:extension-fields-2}: First, we note that 
    \[
    LM \subset \wdh{\ov{K}}M \subset \wdh{\ov{K}M} \subset \wdh{\ov{M}},
    \]
    where the composites are taken inside $\wdh{\ov{N}}$. Then we note that $LM \subset N$ is dense, so inclusion $LM \subset \wdh{\ov{M}}$ uniquely extends to an inclusion $N \subset \wdh{\ov{M}}$. This implies that $M\subset N$ is topologically algebraic. 
\end{proof}

\begin{defn}\label{defn:geometric-point} A {\it geometric point} above point $x\in X$ of an analytic adic space $X$ is a morphism $x\colon \Spa(C(x), C(x)^+) \to X$ such that $C(x)$ is an algebraically closed non-archimedean field, $C(x)^+$ is an open and bounded valuation subring of $C(x)$, and the corresponding extension of completed residue fields $\wdh{k(x)} \subset C(x)$ is a topologically algebraic extension.
\end{defn}

\begin{rmk} If $\Spa(C(x), C(x)^+) \to X$ is a geometric point, then $C(x)$ can be identified with the completed algebraic closure of $\wdh{k(x)}$ (or, equivalently, of $k(x)$) and $C(x)^+$ with a valuation ring extending $\wdh{k(x)}^+$ (or, equivalently, $k(x)^+$). Therefore, Definition~\ref{defn:geometric-point} is more restrictive than \cite[Definition 2.5.1]{H3}, but coincides with the subclass of geometric points constructed in \cite[(2.5.2)]{H3}.
\end{rmk}

\begin{lemma}\label{lemma:restriction-mod-p} Let $K$ be a non-archimedean field with an open and bounded valuation sub-ring $K^+\subset K$ and a pseudo-uniformizer $\varpi$. Let $f\colon X \to Y$ be a morphism of locally of finite type $(K, K^+)$-adic spaces, and $\ov{y}\colon \Spa(C(y), C(y)^+) \to Y$ be a geometric point above $y\in Y$. Then the natural morphism 
\[
a\colon i^{-1}(\O_{X_\et}^+/\varpi) \to \O_{X_{\ov{y}, \et}}^+/\varpi
\]
is an isomorphism where $i\colon X_{\ov{y}} \to X$ is the ``projection'' of the geometric fiber $X_{\ov{y}} \coloneqq X\times_Y \Spa (C(y), C(y)^+)$ back to $X$.
\end{lemma}
\begin{proof}
   \cite[Proposition 2.5.5]{H3} ensures that it suffices to show that $a$ is an isomorphism on stalks at geometric points of $X_{\ov{y}}$. Now note that Lemma~\ref{lemma:extension-fields} implies that any geometric point $\ov{x}\colon \Spa(C(x), C(x)^+) \to X_{\ov{y}}$ defines a geometric point $\ov{x}'\colon \Spa(C(x), C(x)^+) \to X$ of $X$ by taking the composition of $\ov{x}$ with $i$. So it is enough to show that the natural map
   \begin{equation}\label{eqn:different-stalks}
   (\O_{X_\et}^+/\varpi)_{\ov{x}'} \simeq \left(i^{-1}(\O_{X_\et}^+/\varpi)\right)_{\ov{x}} \to (\O_{X_{\ov{y}, \et}}^+/\varpi)_{\ov{x}}
   \end{equation}
   is an isomorphism. But \cite[Proposition 2.6.1]{H3} naturally identifies both sides of (\ref{eqn:different-stalks}) with $C(x)^+/\varpi C(x)^+$ finishing the proof. 
\end{proof}

\begin{rmk} Lemma~\ref{lemma:restriction-mod-p} is very specific to the adic geometry (and quite counter-intuitive from algebraic point of view). Its scheme-theoretic version with $\O^+/\varpi$ replaced by $\O$ is false. The main feature of analytic adic geometry (implicitly) used in the proof is that the morphism $\O_{X, x}^+ \to k(x)^+$ becomes an isomorphism after the $\varpi$-adic completion. 
\end{rmk}

\begin{lemma}\label{lemma:finite-isomorphism} Let $(C, C^+)$ be a Huber pair of an algebraically closed non-archimedean field $C$, an open and bounded valuation sub-ring $C^+ \subset C$ and a pseudo-uniformizer $\varpi \in C^+$. Suppose that $(C, C^+) \to (D, D^+)$ is a finite morphism of complete Huber pairs with a local ring $D$. Then the natural morphism
\[
C^+/\varpi C^+ \to D^+/\varpi D^+
\]
is an isomorphism.
\end{lemma}
\begin{proof}
    First, we show that $C^+/\varpi C^+ \to D^+/\varpi D^+$ is injective. Suppose that $\ov{c}\in C^+/\varpi C^+$ is an element in the kernel and lift it to $c\in C^+$. The assumption on $c$ implies that $c=\varpi d$ for some $d\in D^+$. Then $d=c/\varpi \in C\cap D^+ = C^+$. Therefore, $\ov{c}=0$ in $C^+/\varpi C^+$.\smallskip
    
    Now we check surjectivity. Since $D$ is a local ring that is finite over an algebraically closed field $C$, we conclude that $D$ is an Artin local ring and $D/\rm{nil}(D) \simeq C$. Therefore, for every $d\in D^+$, we can find $c\in C$ and $d'\in \rm{nil}(D)$ such that $d=c + d'$. Since $\rm{nil}(D) \subset D^{\circ\circ} \subset D^+$, we conclude that $c=d-d'\in D^+\cap C =C^+$. Now note that $d'/\varpi$ is still a nilpotent element of $D$, thus $d'/\varpi \in \rm{nil}(D) \subset D^+$. So we conclude that
    \[
    d=c+\varpi(d'/\varpi)
    \]
    proving that $C^+/\varpi C^+ \to D^+/\varpi D^+$ is surjective. 
\end{proof}

\begin{cor}\label{cor:primitive-finite} Let $K$ be a $p$-adic non-archimedean field, and $K^+$ an open and bounded valuation subring of $K$. Let $f\colon X \to Y$ be a finite morphism of locally finite type $(K, K^+)$-adic spaces. Then the natural morphism
\[
c\colon \left(f_*\ud{\bf{F}}_p\right) \otimes \O_{Y_\et}^+/p \to f_*(\O_{X_\et}^+/p)
\]
is an isomorphim on $Y_\et$.
\end{cor}
\begin{proof}
We use \cite[Proposition 2.5.5]{H3} to ensure that it suffices to show that $c$ is an isomorphism on stalks at geometric points. Thus, \cite[Proposition 2.6.1]{H3} and Lemma~\ref{lemma:restriction-mod-p} imply that it suffices to show that the natural map
\[
\rm{H}^0_\et(X, \bf{F}_p) \otimes C^+/p \to \rm{H}^0_\et(X, \O_{X_\et}^+/p)
\]
is an isomorphism when $Y=\Spa(C, C^+)$ for an algebraically closed $p$-adic non-archimedean field $C$ and an open and bounded valuation subring $C^+\subset C$. In this case, $X=\Spa(D, D^+)$ for some finite morphism of Huber pairs $(C,C^+) \to (D, D^+)$. In particular, $D$ is a finite $C$-algebra, so it is a finite direct product of local artinian $C$-algebra. By passing to a direct factor of $D$ (or, geometrically, to a connected component of $\Spa(D, D^+)$), we can assume that $D$ is local. In particular, $D$ does not have any idempotents, and therefore $\Spa(D, D^+)$ is connected. In this case, we have
\[
\rm{H}^0_\et(X, \bf{F}_p) \otimes C^+/pC^+ \simeq C^+/pC^+,
\]
since $\rm{H}^0_\et(X, \bf{F}_p) \simeq \bf{F}_p$ since $\Spa(D, D^+)$ is connected. \smallskip

Now we observe that $\Spa(D, D^+)_{\rm{red}} \simeq \Spa(C, C^+)$, so all \'etale sheaves on $\Spa(D, D^+)$ do not have higher cohomology groups. Thus, we have
\[
\rm{H}^0_\et(X, \O_{X_\et}^+/p) \simeq D^+/pD^+.
\]
In particular, the question boils down to showing that the natural map
\[
C^+/p C^+ \to D^+/p D^+
\]
is an isomorphism. This was already done in Lemma~\ref{lemma:finite-isomorphism}.
\end{proof}

\begin{cor}\label{cor:primitive-finite-rigid} Let $K$ be a $p$-adic non-archimedean field, $f\colon X\to Y$ a finite morphism of rigid-analytic varieties over $K$, and $\F\in \bf{D}^b_{zc}(X; \bf{F}_p)$. Then the natural morphism
\[
c\colon \left(f_*\F\right) \otimes \O_{Y_\et}^+/p \to f_*\left(\F\otimes \O_{X_\et}^+/p\right)
\]
is an isomorphim on $Y_\et$.
\end{cor}
\begin{proof}
We recall that \cite[Proposition 3.6]{Bhatt-Hansen} says that $\bf{D}^b_{zc}(X; \bf{F}_p)$ is a thick triangulated subcategory of $\bf{D}(X_\et; \bf{F}_p)$ generated by objects of the form $g_*\ud{\bf{F}}_p$ for finite morphisms $g\colon X' \to X$. Since both claims in the question satisfy the $2$-out-of-$3$ property and are preserved by passing to direct summands, it suffices to prove the claim only for $\F=g_*\ud{\bf{F}}_p$. In this situation, the claim follows from Corollary~\ref{cor:primitive-finite} by a sequence of isomorphisms
\begin{align*}
    f_*\left(g_*\left(\ud{\bf{F}}_p\right)\right)\otimes \O_{Y_\et}^+/p & \simeq \left(f\circ g\right)_*(\ud{\bf{F}}_p) \otimes \O_{Y_\et}^+/p \\
    & \simeq (f\circ g)_*\left(\O_{X'_\et}^+/p\right)  \\
    & \simeq f_*\left(g_*\O_{X'_\et}^+/p\right) \\
    & \simeq f_*\left( g_*\ud{\bf{F}}_p\otimes \O_{X_\et}^+/p \right). \qedhere
\end{align*}
\end{proof}

\subsection{Applications}

The main goal of this section is to discuss some applications of Theorem~\ref{main-thm}. In particular, we show that ``p-adic nearby cycles'' commute with proper pushfowards and prove finiteness of the usual \'etale cohomology of proper rigid-analytic varieties. \smallskip 

For the rest of the section, we fix a $p$-adic algebraically closed field $C$ with its rank-$1$ valuation ring $\O_C$, maximal ideal $\m\subset \O_C$, and a good pseudo-uniformizer $\varpi\in \O_C$ (see Definition~\ref{defn:good-unifor}). We always do almost mathematics with respect to the ideal $\m$ in this section. If we need to consider a more general non-archimedean field, we denote it by $K$. \smallskip

The first non-trivial consequence of Theorem~\ref{thm:main-thm-integral} is that the $v$-cohomology groups of ${\O}_{X^\diam}^+$-vector bundles have bounded $p$-torsion. 

\begin{lemma}\label{lemma:small-torsion} Let $K$ be a $p$-adic perfectoid field, let $\X=\Spf A_0$ be an affine admissible formal $\O_K$-scheme with adic generic fiber $X$, and let $\cal{E}$ be an $\O_{X^\diam}^+$-vector bundle. Then the cohomology groups $\rm{H}^i(X^\diam_v, \cal{E})$ are almost finitely presented over $A_0$. In particular, they are $p$-adically complete and have bounded torsion $p^{\infty}$-torsion. 
\end{lemma}
\begin{proof}
This is a straightforward consequence of Theorem~\ref{thm:main-thm-integral}, Lemma~\ref{bounded-torsion} and Lemma~\ref{completion-finitely-generated}. 
\end{proof}

\begin{rmk} Lemma~\ref{lemma:small-torsion} implies that the $v$ cohomology groups of $\O_{X^\diam}^+$ behave pretty differently from the analytic cohomology groups of $\O_X^+$. Indeed, see \cite[Remark 9.3.4]{Bhatt-notes} (that can be easily adapted to the $p$-adic situation) for an example of an affinoid rigid-analytic variety with unbounded $p$-torsion in $\rm{H}^1_{\rm{an}}(X, \O_X^+)$.
\end{rmk}

\begin{thm}\label{thm:finiteness-bundles} Let $K$ be a $p$-adic perfectoid field, let $X$ be a proper rigid-analytic $K$-variety of dimension $d$, and let $\cal{E}$ be an $\O_{X^\diam}^+$-vector bundle (resp. $\O_{X^\diam}^+/p$-vector bundle). Then 
\[
\bf{R}\Gamma(X_v^\diam, \cal{E}) \in \bf{D}^{[0, 2d]}_{acoh}(\O_K)^a.
\]
\end{thm}
\begin{proof}
We firstly deal with the case of an $\O_{X^\diam}^+/p$-vector bundle $\cal{E}$. We choose an admissible formal model $\X$ of $X$ as in Part~\ref{thm:main-thm-small-5} of Theorem~\ref{thm:main-thm-small}. This formal model is automatically proper by \cite[Lemma 2.6]{Lutke-proper} and \cite[Corollary 4.4 and 4.5]{Temkin-proper}.  Now Theorem~\ref{thm:main-thm-small} implies that 
\[
\bf{R}\nu_*\left(\cal{E}\right)^a \in \bf{D}_{acoh}^{[0,d]}(\X_0)^a.
\]
Recall that the underlying topological spaces of $\X_0$ and the special fiber $\ov{\X}\coloneqq \X\times_{\Spf \O_C} \Spec \O_C/\m$ are the same. Thus, \cite[Corollary II.10.1.11]{FujKato} implies that $\X_0$ has Krull dimension $d$. Therefore, Theorem~\ref{almost-proper-mapping}, \cite[\href{https://stacks.math.columbia.edu/tag/0A3G}{Tag 0A3G}]{stacks-project} and Lemma~\ref{projection} imply that
\[
\bf{R}\Gamma(X_v^\diam, \cal{E})^a \simeq \bf{R}\Gamma\left(\X_0, \bf{R}\nu_*\left(\cal{E}\right)\right)^a \in \bf{D}_{acoh}^{[0,2d]}(\O_K/p)^a.
\]
The case of an $\O_{X^\diam}^+$-vector bundle follows from the $\O_{X^\diam}^+/p$-case, Corollary~\ref{cor:check-almost-coh-mod-p-2}, and Lemma~\ref{lemma:first-properties-structure-sheaves}~\ref{lemma:first-properties-structure-sheaves-3}.
\end{proof}

\begin{lemma}\label{almost-scholze} Let $X$ be a proper rigid-analytic variety over $C$ of dimension $d$, and  $\F$ a Zariski-constructible sheaf of $\bf{F}_p$-modules on $X_\et$. Then $\bf{R}\Gamma(X_v^\diam, \F\otimes \O_{X^\diam}^+/p)^a \in \bf{D}^{[0, 2d]}_{acoh}(\O_C/p)^a$.
\end{lemma}
\begin{proof}
The proof is analogous to the proof of Theorem~\ref{thm:finiteness-bundles} using Theorem~\ref{main-thm} in place of Theorem~\ref{thm:main-thm-small}.
\end{proof}

Now we discuss finiteness of classical \'etale cohomology groups. Later, we will generalize it to Zariski-constructible coefficients. 

\begin{lemma}\label{lemma:finiteness-cohomology-proper} Let $X$ be a proper rigid-analytic variety over $C$ of dimension $d$. Then 
\[
\bf{R}\Gamma(X, \bf{F}_p) \in \bf{D}^{[0, 2d]}_{coh}(\bf{F}_p)
\]
and the natural morphism
\[
\bf{R}\Gamma(X, \bf{F}_p) \otimes \O_C/p \to \bf{R}\Gamma(X_v^\diam, \O_{X^\diam}^+/p)
\]
is an almost isomorphism. 
\end{lemma}
\begin{proof}
{\it Step~$1$. $\bf{R}\Gamma(X_v^\diam, \O_{X^\diam}^{\flat, +})^a \in \bf{D}^{[0, 2d]}_{acoh}(\O^\flat_C)^a$}. We consider the tilted integral structure sheaf $\O_{X^\diam}^{\flat, +}$ (see Definition~\ref{defn:v-structure-sheaf-tilted}). Lemma~\ref{lemma:first-properties-structure-sheaves}\ref{lemma:first-properties-structure-sheaves-4} ensures that $\O_{X^\diam}^{\flat, +}$ is derived $\varpi^\flat$-adically complete and Lemma~\ref{lemma:first-properties-structure-sheaves}\ref{lemma:first-properties-structure-sheaves-5} implies that 
\[
[\O_{X^\diam}^{\flat, +}/\varpi^\flat] \simeq [\O_{X^\diam}^+/p] \simeq \O_{X^\diam}^+/p.
\] 
Therefore, \cite[\href{https://stacks.math.columbia.edu/tag/0BLX}{Tag 0BLX}]{stacks-project} guarantees that $ \bf{R}\Gamma\left(X_v^\diam, \O_{X^\diam}^{\flat, +}\right) \in \bf{D}\left(\O^\flat_C\right)$ is derived $\varpi^\flat$-adically complete. Moreover, Lemma~\ref{almost-scholze} implies 
\[
\left[\bf{R}\Gamma\left(X_v^\diam, \O_{X^\diam}^{\flat, +}\right)^a/\varpi^\flat\right] \simeq \bf{R}\Gamma(X_v^\diam, \O_{X^\diam}^+/p)^a \in \bf{D}^{[0, 2d]}_{acoh}(\O_C/p)^a. 
\]
Thus, Corollary~\ref{cor:check-almost-coh-mod-p-2} applied to $R=C^+=\O_C^\flat$ implies that $\bf{R}\Gamma\left(X_v^\diam, \O_{X^\diam}^{\flat, +}\right)^a\in \bf{D}^{[0,2d]}_{acoh}(\O_C^\flat)^a$. \smallskip

{\it Step~$2$. $\bf{R}\Gamma(X, \bf{F}_p) \in \bf{D}^{[0, 2d]}_{coh}(\bf{F}_p)$ and the natural morphism $\bf{R}\Gamma(X, \bf{F}_p) \otimes C^\flat \to \bf{R}\Gamma(X_v^\diam, \O_{X^\diam}^{\flat})$ is an isomorphism.} After inverting $\varpi^\flat$, Step~$1$ implies that
\[
\bf{R}\Gamma(X_v^\diam, \O_{X^\diam}^{\flat}) \in \bf{D}^{[0, 2d]}_{coh}(C^\flat).
\]
Since $\O_{X^\diam}^{\flat}$ is a sheaf of $\bf{F}_p$-algebras, we have a natural Frobenius morphism 
\[
F\colon \O_{X^\diam}^{\flat} \xr{f\mapsto f^p} \O_{X^\diam}^{\flat}
\]
that can be easily seen to be an isomorphism by Lemma~\ref{lemma:first-properties-structure-sheaves}\ref{lemma:first-properties-structure-sheaves-2} (and Remark~\ref{rmk:perfect-perfectoid}). Now we use the Artin-Shreier short exact sequence
\[
0 \to \ud{\bf{F}}_p \to \O_{X^\diam}^{\flat} \xr{F - \rm{Id}} \O_{X^\diam}^{\flat} \to 0
\]
on the $v$-site $X^\diam_v$ to get the associated long exact sequence\footnote{We implicitly use that $\rm{H}^i(X, \bf{F}_p) \simeq \rm{H}^i(X_v^\diam, \bf{F}_p)$ by \cite[Proposition 14.7, 14.8, and Lemma 15.6]{Sch2}.}
\[
\dots \to \rm{H}^i(X, \bf{F}_p) \to \rm{H}^i(X_v^\diam, \O_{X^\diam}^{\flat}) \xr{\rm{H}^i(F)-\rm{Id}} \rm{H}^i(X_v^\diam, \O_{X^\diam}^{\flat}) \to \rm{H}^{i+1}(X, \bf{F}_p) \to \dots
\]
We already know that each group $\rm{H}^i(X_v^\diam, \O_{X^\diam}^{\flat})$ is a finitely generated $C^\flat$-vector space, each $\rm{H}^i(F)$ is a Frobenius-linear automorphism, and $C^\flat$ is an algebraically closed field of characteristic $p$ (see \cite[Theorem 3.7]{Sch0}). Thus (the proof of) \cite[\href{https://stacks.math.columbia.edu/tag/0A3L}{Tag 0A3L}]{stacks-project} ensures that $\rm{H}^i(F)-\rm{Id}$ is surjective for each $i\geq 0$ (so $\rm{H}^i(X, \bf{F}_p)\simeq \rm{H}^i(X_v^\diam, \O_{X^\diam}^{\flat})^{F=1}$) and the natural morphism 
\[
\rm{H}^i(X, \bf{F}_p) \otimes C^\flat  \to \rm{H}^i(X_v^\diam, \O_{X^\diam}^{\flat})
\]
is an isomoprhism. In particular, we have $\dim_{\bf{F}_p} \rm{H}^i(X, \bf{F}_p) = \dim_{C^\flat} \rm{H}^i(X_v^\diam, \O_{X^\diam}^{\flat})$, the natural morphism
\[
\bf{R}\Gamma(X, \bf{F}_p) \otimes C^\flat \to \bf{R}\Gamma(X^\diam_v, \O_{X^\diam}^{\flat})
\]
is an isomorphism, and $\bf{R}\Gamma(X, \bf{F}_p)\in \bf{D}^{[0, 2d]}_{coh}(\bf{F}_p)$.\smallskip

{\it Step~$3$. The natural morphism $\bf{R}\Gamma(X,\bf{F}_p) \otimes \O_C/p \to \bf{R}\Gamma(X_v^\diam, \O_{X^\diam}^+/p)$ is an almost isomorphism.} It suffices to show that 
\[
\bf{R}\Gamma(X, \bf{F}_p)\otimes \O_C^\flat \to \bf{R}\Gamma(X_v^\diam, \O_{X^\diam}^{\flat, +})
\]
is an almost isomorphism. The version with $\O_{X^\diam}^+/p$ would follow by taking the derived mod-$\varpi^\flat$ reduction. Therefore, it suffices to show that
\[
\rm{H}^i(X, \bf{F}_p)\otimes \O_C^\flat \to \rm{H}^i(X_v^\diam, \O_{X^\diam}^{\flat, +})
\]
is an almost isomorphism for each $i\geq 0$. We consider the following commutative diagram
\[
\begin{tikzcd}
\rm{H}^i(X, \bf{F}_p)\otimes \O_C^\flat \arrow{r}{\a} \arrow{d}{i} & \rm{H}^i(X_v^\diam, \O_{X^\diam}^{\flat, +}) \arrow{d} \\
\rm{H}^i(X, \bf{F}_p)\otimes C^\flat \arrow{r}{\beta} & \rm{H}^i(X_v^\diam, \O_{X^\diam}^{\flat}).
\end{tikzcd}
\]
By Step~$2$, we know that $\beta$ is an isomorphism. Since $i$ is injective, we conclude that $\a$ is injective as well. So it suffices to show that $\a$ is almost surjective. \smallskip

Now we note that Frobenius acts on $\rm{H}^i(X, \bf{F}_p)\otimes \O_C^\flat$ by acting on $\O_C^\flat$ and Frobenius acts on $\rm{H}^i(X_v^\diam, \O_{X^\diam}^{\flat, +})$ by acting on $\O_{X^\diam}^{\flat, +}$. Furthermore, the map $\a$ is Frobenius-equivariant. The action on $\rm{H}^i(X, \bf{F}_p)\otimes \O_C^\flat$ is an isomorphism because $\O_C^\flat$ is perfect, and the action on $\rm{H}^i(X_v^\diam, \O_{X^\diam}^{\flat, +})$ is an isomorphism because Frobenius is already an isomorphism on the sheaf $\O_{X^\diam}^{\flat, +}$ due to Lemma~\ref{lemma:first-properties-structure-sheaves}~\ref{lemma:first-properties-structure-sheaves-2} (and Remark~\ref{rmk:perfect-perfectoid}). Therefore, it makes sense to consider the inverse Frobenius action $F^{-1}$ on both modules and $\a$ commutes with this action. \smallskip

We pick an element $x\in \rm{H}^i(X_v^\diam, \O_{X^\diam}^{\flat, +})$. Since $F$ is an isomorphism on $\rm{H}^i(X_v^\diam, \O_{X^\diam}^{\flat, +})$, there is an $x'\in \rm{H}^i(X_v^\diam, \O_{X^\diam}^{\flat, +})$ such that $F^m(x')=x$. Since $\rm{H}^i(X^\diam_v, \O^{\flat, +}_{X^\diam})$ is almost coherent, Lemma~\ref{bounded-torsion} implies that it has bounded $(\varpi^\flat)^{\infty}$-torsion. Combining this with the fact that $\beta$ is an isomorphism, we conclude that there is an integer $N$ and an element $y'\in \rm{H}^i(X_v^\diam, \bf{F}_p)\otimes \O_C^\flat$ such that $\a(y')=(\varpi^\flat)^Nx'$. Therefore,
\[
\big(\varpi^\flat\big)^{N/p^m}x=F^{-m}\left(\big(\varpi^\flat\big)^Nx'\right)=F^{-m}\left(\a\left(y'\right)\right)=\a\left(F^{-m}\left(y'\right)\right).
\]
Thus $(\varpi^\flat)^{N/p^m}x=\a(y)$ where $y=F^{-m}(y') \in \rm{H}^i(X, \bf{F}_p)\otimes \O_C^\flat$. Since $N/p^m$ can be made arbitrary small by increasing $m$, we conclude that $\a$ is almost surjective. 
\end{proof}

\begin{lemma}\label{lemma:primitive-zc} Let $X$ be a proper rigid-analytic variety over $C$ of dimension $d$, and $\F\in \bf{D}^{[r,s]}_{zc}(X; \bf{F}_p)$ for some integers $[r, s]$. Then 
\[
\bf{R}\Gamma(X, \F) \in \bf{D}^{[r, s+2d]}_{coh}(\bf{F}_p).
\]
\end{lemma}  
\begin{proof}
    First, \cite[Corollary 2.8.3]{H3} implies that $\bf{R}\Gamma(X, \F) \in \bf{D}^{[r, s +2d]}(\bf{F}_p)$. Therefore, it suffices to show that $\bf{R}\Gamma(X, \F) \in \bf{D}_{coh}(\bf{F}_p)$. For this, we recall that \cite[Proposition 3.6]{Bhatt-Hansen} says that $\bf{D}^b_{zc}(X, \bf{F}_p)$ is a thick triangulated subcategory of $\bf{D}(X_\et; \bf{F}_p)$ generated by objects of the form $f_*(\ud{\bf{F}}_p)$ for finite morphisms $f\colon X' \to X$. Since the claim we try to prove satisfies the $2$-out-of-$3$ property and is preserved by passing to direct summands, it suffices to prove the claim only for $\F=f_*(\ud{\bf{F}}_p)$. Then  Lemma~\ref{lemma:finiteness-cohomology-proper} and \cite[Proposition 2.6.3]{H3} imply that
    \[
    \bf{R}\Gamma\left(X, f_*\left(\ud{\bf{F}}_p\right)\right)\simeq \bf{R}\Gamma(X', \bf{F}_p) \in \bf{D}^{[0,2d]}_{coh}(\bf{F}_p). \qedhere
    \]
\end{proof}  

The last thing we discuss is the behaviour of the ``$p$-adic nearby cycles'' under proper pushforwards. We start with the following lemma:

\begin{lemma}\label{lemma:relative-primitive} Let $K$ be a $p$-adic perfectoid field $K$, let $f\colon X \to Y$ be a proper morphism of rigid-analytic varieties over $K$, and let $\F\in \bf{D}^b_{zc}(X; \bf{F}_p)$. Then the natural morphism
\[
\bf{R}f_*\F \otimes \O_{Y_\et}^+/p \to \bf{R}f_*(\F \otimes \O_{X_\et}^+/p)
\]
is an almost isomorphism. 
\end{lemma}
\begin{proof}
    The claim is local on $Y$, so we can assume that $Y$ is affinoid. Then a similar argument as in the proof of Corollary~\ref{lemma:primitive-zc} allows us to reduce to the case when $\F=g_*\left(\ud{\bf{F}}_p\right)$ for a finite map $g\colon X' \to X$. Therefore, Corollary~\ref{cor:primitive-finite} implies that it suffices to prove the claim for the morphism $f\circ g\colon X' \to Y$ and $\F=\ud{\bf{F}}_p$. \smallskip
    
    Now \cite[Proposition 2.5.5]{H3} guarantees that it suffices to show the claim on stalks at geometric points. Therefore, by Lemma~\ref{lemma:restriction-mod-p} we reduce the question to showing that, for any proper adic space $X$ over a geometric point $\Spa(C, C^+)$, the natural morphism
    \[
    \bf{R}\Gamma(X, \bf{F}_p)\otimes C^+/p \to \bf{R}\Gamma(X, \O_{X_\et}^+/p).
    \]
    is an almost isomorphism. Denote by $X^\circ \coloneqq X\times_{\Spa(C, C^+)} \Spa(C, C^\circ)$. Now \cite[Proposition 8.2.3(ii)]{H3} implies that $\bf{R}\Gamma(X, \bf{F}_p) \simeq \bf{R}\Gamma(X^\circ, \bf{F}_p)$, Lemma~\ref{lemma:perfectoid-rank-1} implies that $C^+/pC^+ \simeq^a \O_C/p\O_C$, and Corollary~\ref{cor:O+/p-cohomology} and Corollary~\ref{cor:base-change-result} imply that 
    \[
    \bf{R}\Gamma(X, \O_{X_\et}^+/p) \simeq^a \bf{R}\Gamma(X^\circ, \O_{X^\circ_\et}^+/p).
    \]
    Combining these results, we may replace $(C, C^+)$ with $(C, \O_C)$ and $X$ with $X^\circ$ to achieve that $\Spa(C, \O_C)$ is a geometric point of rank-$1$. In this case, the claim was already proven in Lemma~\ref{lemma:primitive-zc}.
\end{proof}

Now we show that $p$-adic nearby cycles commute with proper morphisms.

\begin{cor}\label{lemma:nearby-pushforward} Let $K$ be a $p$-adic perfectoid field $K$, let $\mf\colon \X \to \Y$ be a proper morphism of admissible formal $\O_K$-schemes with adic generic fiber $f\colon X \to Y$, and let $\F\in \bf{D}^b_{zc}(X; \bf{F}_p)$. Then the natural morphism
\[
\bf{R}\nu_{\Y, *}\left(\bf{R}f_{*} \F\otimes \O^+_{Y^\diam}/p\right) \to \bf{R}\mf_{0, *}\left(\bf{R}\nu_{\X, *}\left(\F\otimes \O_{X^\diam}^+/p\right)\right)
\]
is an almost isomorphism. 
\end{cor}
\begin{proof}
First, note that $\bf{R}f_*\F$ has overconvergent cohomology sheaves by \cite[Proposition 8.2.3(ii)]{H3} and Remark~\ref{rmk:zc-overconvergent}. Therefore, Lemma~\ref{lemma:et-qp-v-overconv-coeff} implies that
\[
\bf{R}\nu_{\Y, *}\left(\bf{R}f_{*} \F\otimes \O^+_{Y^\diam}/p \right) \simeq \bf{R}t_{\Y, *}\left(\bf{R}f_{*} \F\otimes \O^+_{Y_\et}/p\right),
\]
where $t_\Y\colon \left(Y_\et, \O_{Y_\et}^+/p\right) \to (\Y_0, \O_{\Y_0})$ is the natural morphism of ringed sites. Similarly, we have an isomorphism
\[
\bf{R}\mf_{0, *}\left(\bf{R}\nu_{\X, *}\left(\F\otimes \O_{X^\diam}^+/p\right)\right) \simeq \bf{R}\mf_{0, *}\left(\bf{R}t_{\X, *}\left(\F\otimes \O_{X_\et}^+/p\right)\right).
\]
Therefore, it suffices to show that the natural morphism
\[
\bf{R}t_{\Y, *}\left(\bf{R}f_{*} \F\otimes \O^+_{Y_\et}/p\right) \to \bf{R}\mf_{0, *}\left(\bf{R}t_{\X, *}\left(\F\otimes \O_{X_\et}^+/p\right)\right)
\]
is an almost isomorphism. \smallskip

For this, we observe that the commutative diagram of ringed sites
\[
\begin{tikzcd}
\left(X_\et, \O_{X_\et}^+/p\right) \arrow{r}{t_{\X}} \arrow{d}{f} & \left(\X_0, \O_{\X_0}\right) \arrow{d}{\mf_0}\\
\left(Y_\et, \O_{Y_\et}^+/p\right) \arrow{r}{t_{\Y}} & \left(\Y_0, \O_{\Y_0}\right)
\end{tikzcd}
\]
implies that 
\[
\bf{R}\mf_{0, *}\left(\bf{R}t_{\X, *}\left(\F\otimes \O_{X_\et}^+/p\right)\right) \simeq \bf{R}t_{\Y, *}\left(\bf{R}f_{*}\left(\F\otimes \O^+_{X_\et}/p\right)\right). 
\]
Therefore, the morphism
\[
\bf{R}t_{\Y, *}\left(\bf{R}f_{*}\F \otimes \O^+_{Y_\et}/p\right) \to \bf{R}\mf_{0, *}\left(\bf{R}t_{\Y, *}\left(\F\otimes \O_{Y_\et}^+/p\right)\right)  \simeq \bf{R}t_{\Y, *}\left(\bf{R}f_{*}\left(\F\otimes \O^+_{X_\et}/p\right)\right)
\]
is an almost isomorphism due to Lemma~\ref{lemma:relative-primitive} and Proposition~\ref{derived-pushforward}.
\end{proof}

\subsection{Perfectoid covers of affinoids}\label{section:perfectoid-covers}  
  
The main goal of this section is to show almost vanishing of higher $v$-cohomology groups of a small $\O^+_{X^\diam}/p$-vector bundle on an affinoid perfectoid space. Later on, we will apply it to certain pro-\'etale coverings of $\Spa(A, A^+)$ to reduce the computation of $v$-cohomology groups to the computation of \v{C}ech cohomology groups.  

\begin{setup}\label{set-up:covers} We fix 
\begin{enumerate}[label=\textbf{(\arabic*)}]
    \item a $p$-adic perfectoid field $K$ with its rank-$1$ open and bounded valuation ring $\O_K$ and a good pseudo-uniformizer $\varpi\in \O_K$ as in Definition~\ref{defn:good-unifor} (we always do almost mathematics with respect to the ideal $\m= \bigcup_{n}\varpi^{1/p^n}\O_K = K^{\circ \circ}$),
    \item an affine admissible formal $\O_K$-scheme $\X=\Spf A_0$ with adic generic fiber $X=\Spa(A, A^+)$;
    \item and an affinoid perfectoid pair $(A_\infty, A_\infty^+)$ (see Definition~\ref{defn:perfectoid-pair}) with a morphism $(A, A^+) \to (A_\infty, A_\infty^+)$ such that $\Spd(A_\infty, A_\infty^+) \to \Spd(A, A^+)$ is a $v$-covering (see Definition~\ref{defn:v-topology} and Definition~\ref{defn:diamondification});
    \item a small $\O_{X^\diam}^+/p$-vector bundle $\cal{E}$ (see Definition~\ref{defn:small-mod-p}). 
\end{enumerate}
\end{setup}

\begin{defn}\label{defn:integral-perfectoid} We say that a $p$-torsionfree (equivalently, $\varpi$-torsionfree) $\O_K$-algebra $R$ is {\it integrally perfectoid} if the Frobenius homomorphism $R/\varpi R \xr{x\mapsto x^p} R/\varpi^p R = R/pR$ is an isomorphism. 
\end{defn}

\begin{rmk}\label{rmk:infty-integral-perfectoid} \cite[Lemma 3.10]{BMS1} implies that this definition coincides with \cite[Definition 3.5]{BMS1} for $p$-torsionfree $\O_K$-algebras. In particular, $A_\infty^+$ is an integral perfectoid $\O_K$-algebra by \cite[Lemma 3.20]{BMS1}.
\end{rmk}

\begin{lemma}\label{lemma:perf-base-change} Under the assumption of Set-up~\ref{set-up:covers}, let $\mf\colon \Spf B_0 \to \Spf A_0$ be an \'etale morphism of admissible affine formal $\O_K$-schemes. Then $B^+_{\infty}\coloneqq B_0\wdh{\otimes}_{A_0}A^+_{\infty}$ is $p$-torsionfree integrally perfectoid $\O_K$-algebra. 
\end{lemma}
\begin{proof}
Firstly, we note that $A_0 \to B_0$ is a flat morphism by \cite[Proposition I.4.8.1]{FujKato}, so $B_0\otimes_{A_0} A^+_{\infty}$ is $\varpi$-torsionfree. Since the $\varpi$-adic completion of a $\varpi$-torsionfree algebra is $\varpi$-torsionfree, we conclude that $B_\infty^+=B_0\wdh{\otimes}_{A_0} A_{\infty}^+$ is $\varpi$-torsionfree. We see that the only thing we are left to show is that the Frobenius morphism
\[
B^+_{\infty}/\varpi B^+_{\infty} \to B^+_{\infty}/\varpi^p B^+_{\infty}
\]
is an isomorphism. We consider the commutative diagram
\[
\begin{tikzcd}[]
\Spec B^+_{\infty}/\varpi\arrow{rd}{F} \arrow{rrd}{\Phi_B^*} \arrow{rdd}[swap]{\mf_{\infty}/\varpi}& & \\
& \Spec \left(B^+_{\infty}/\varpi^p \otimes_{A^+_{\infty}/\varpi^p} A^+_{\infty}/\varpi\right) \arrow{r}[swap]{\Phi_A^*\times B_0} \arrow{d}{\alpha} & \Spec B^+_{\infty}/\varpi^p \arrow{d}{\mf_{\infty}/\varpi^p} \\
& \Spec A^+_{\infty}/\varpi \arrow{r}{\Phi_A^*} & \Spec A^+_{\infty}/\varpi^p.
 \end{tikzcd}
\]
We need to show that $\Phi^*_B$ is an isomorphism. We know that $\mf_{\infty}/\varpi^p$ and $\mf_{\infty}/\varpi$ are \'etale morphisms since $\mf$ is, and the Frobenious $\Phi_A^*$ an isomorphism by Remark~\ref{rmk:infty-integral-perfectoid}. Therefore, the morphism 
\[
\alpha \colon \Spec \left(B^+_{\infty}/\varpi^p \otimes_{A^+_{\infty}/\varpi^p} A^+_{\infty}/\varpi\right) \to \Spec A_{\infty}^+/\varpi
\]
is \'etale as a base change of the \'etale morphism $\mf_\infty/\varpi^p$. Thus, we conclude that $F$ is an \'etale morphism as a morphism between \'etale $A_\infty^+/\varpi$-schemes. Now we note that $\Phi^*_A\times B_0$ is an isomorphism since $\Phi_A^*$ is an isomorphism. Therefore, $\Phi^*_B$ is an \'etale morphism as a composition of an \'etale morphism and an isomorphism. However, $\Phi_B^*$ is a bijective radiciel morphism since it is the absolute Frobenius morphism. Thus, we conclude that it must be an isomorphism as any \'etale, bijective radiciel morphism is an isomorphism by \cite[Th\'eor\`eme 5.1]{SGA1}.  
\end{proof}

\begin{cor}\label{cor:restrict} Under the assumption of Set-up~\ref{set-up:covers}, let $\mf\colon \Spf B_0 \to \Spf A_0$ be an \'etale morphism of admissible affine formal $\O_K$-schemes.  Then 
\[
(B_\infty, B_\infty^+)\coloneqq \left(\left(B_0\wdh{\otimes}_{A_0} A_\infty^+\right)[1/p], B_0\wdh{\otimes}_{A_0}A_\infty^+\right)
\]
is a perfectoid pair. 
\end{cor} 
\begin{proof}
Lemma~\ref{lemma:perf-base-change} states that $B_\infty^+=B_0\wdh{\otimes}_{A_0} A^+_\infty$ is a $p$-torsionfree integral perfectoid. Now $B_0\otimes_{A_0} A^+_\infty$ is integrally closed in $B_0\otimes_{A_0} A^+_\infty[1/p]$ because $A^+$ is integrally closed in $A$ and $B_0$ is \'etale over $A_0$. Therefore, \cite[Lemma 5.1.2]{Bhatt-notes} ensures that the same holds after completion, i.e. $B_\infty^+$ is integrally closed in $B_\infty$. Thus \cite[Lemma 3.20]{BMS1} guarantees that $(B_\infty, B_\infty^+)$ is a perfectoid pair.
\end{proof}

\begin{lemma}\label{lemma:no-integral-closure} Under the assumption of Set-up~\ref{set-up:covers}, let $\mf\colon \Spf B_0 \to \Spf A_0$ be an \'etale morphism of admissible affine formal $\O_K$-schemes with adic generic fiber $\Spa(B, B^+) \to \Spa(A, A^+)$. Then the natural morphism
\[
\left(\left(B_0\wdh{\otimes}_{A_0} A_\infty^+\right)[1/p], B_0\wdh{\otimes}_{A_0}A_\infty^+\right) \to \left(B\wdh{\otimes}_A A_\infty, (B\wdh{\otimes}_A A_\infty)^+\right)
\]
is an isomorphism of Tate-Huber pairs.
\end{lemma}
\begin{proof}
    By \cite[Lemma 1.6]{H0}, $B\wdh{\otimes}_A A_\infty \simeq \left(B_0\wdh{\otimes}_{A_0} A_\infty^+\right)[1/p]$. Now $(B\wdh{\otimes}_A A_\infty)^+$ is defined to be the integral closure of the image of the map
    \[
    B^+\wdh{\otimes}_{A^+}A_\infty^+ \to B\wdh{\otimes}_A A_\infty.
    \]
    By \cite[Lemma 1.6]{H0}, we also have 
    \[
    B^+\wdh{\otimes}_{A^+}A_\infty^+ \simeq \left(B^+\otimes_{A^+}A_\infty^+\right)\otimes_{B_0\otimes_{A_0} A^+_\infty} \left(B_0\wdh{\otimes}_{A_0} A^+_\infty \right).
    \]
    Since $B^+$ is integral over $B_0$, we conclude that $B^+\wdh{\otimes}_{A^+}A_\infty^+$ is integral over $B_0\wdh{\otimes}_{A_0} A^+_\infty$. In particular, we see that $(B\wdh{\otimes}_A A_\infty)^+$ is integral over $B_0\wdh{\otimes}_{A_0} A^+_\infty$. However, Corollary~\ref{cor:restrict} implies that $B_0\wdh{\otimes}_{A_0}A_\infty^+$ is a sub-algebra of $B\wdh{\otimes}_A A_\infty$ that is integrally closed in $B\wdh{\otimes}_A A_\infty$. Thus we must have an equality
    \[
    B_0\wdh{\otimes}_{A_0}A_\infty^+ \simeq (B\wdh{\otimes}_A A_\infty)^+. \qedhere
    \]
\end{proof}

\begin{rmk} It will be crucial for our arguments later that $(B\wdh{\otimes}_A A_\infty)^+$ is {\it equal} to $B_0\wdh{\otimes}_{A_0}A_\infty^+$ and not simply to its integral closure. 
\end{rmk}

\begin{lemma}\label{lemma:almost-no-higher-coh-small} Under the assumption of Set-up~\ref{set-up:covers}, we put
\[
M_\cal{E}\coloneqq \rm{H}^0\left(\Spd(A_\infty, A_\infty^+)_v, \cal{E}\right).
\]
Then $M_\cal{E}$ is an almost faithfully flat, almost finitely presented $A_\infty^+/p$-module, and for every morphism $\Spa(D, D^+) \to \Spa(A_\infty, A_\infty^+)$ of affinoid perfectoid spaces, the natural morphism
\[
M_\cal{E}\otimes_{A_\infty^+/p} D^+/p \to \rm{H}^0\left(\Spd(D, D^+)_v, \cal{E}\right)
\]
is an almost isomorphism\footnote{We note that $\cal{E}$ is a sheaf on a (big) $v$-site of $\Spd(A, A^+)$, so it makes sense to evaluate $\cal{E}$ on $\Spd(D, D^+)$.}. Moreover, 
\[
\rm{H}^i\left(\Spd(A_\infty, A_\infty^+)_v, \cal{E}\right)\simeq^a 0
\]
for $i>0$.
\end{lemma}
\begin{proof}
    {\it Step~$1$. $\rm{H}^0\left(\rm{Spd}(A_\infty, A_\infty^+)_v, \cal{E}\right)$ is almost flat and almost finitely presented:} The smallness assumption implies that there is a finite \'etale surjection $\Spa(B, B^+) \to \Spa(A_\infty, A_\infty^+)$ such that $\cal{E}|_{\Spd(B, B^+)} \simeq (\O_{X^\diam}^+/p)^r$ for some integer $r\geq 0$. The adic space $\Spa(B, B^+)$ is affinoid perfectoid by \cite[Theorem 7.9]{Sch1}. \smallskip
    
    The natural morphism $A^+_\infty \to B^+$ is almost finitely presented and almost faithfully flat by \cite[Theorem 7.9]{Sch1} (see also \cite[Theorem 10.0.9]{Bhatt-notes} for the almost {\it faithfully} flat part). Since $\cal{E}|_{\Spd(B, B^+)}$ is trivial, Lemma~\ref{lemma:first-properties-structure-sheaves}~\ref{lemma:first-properties-structure-sheaves-1} implies that 
    \[
    \rm{H}^0\left(\Spd(B, B^+)_v, \cal{E}\right) \simeq^a (B^+/pB^+)^r.
    \]
    In particular, it is almost flat and almost finitely presented. We now want to descend these properties to $\rm{H}^0\left(\Spd(A_\infty, A_\infty^+)_v, \cal{E}\right)$. For this, we use Proposition~\ref{prop:properties-of-diamond} to recall that diamondification commutes with fiber products, and so  
    \begin{align*}
    \Spd(B, B^+)\times_{\Spd(A_\infty, A_\infty^+)} \Spd(B, B^+) &\simeq \left(\Spa(B, B^+)\times_{\Spa(A_\infty, A_\infty^+)} \Spa(B, B^+)\right)^{\diam} \\
    & \simeq \Spd\left(B \wdh{\otimes}_{A_\infty} B, (B \wdh{\otimes}_{A_\infty} B)^+\right).
    \end{align*}
    By the proof of \cite[Proposition 6.18]{Sch0} (and Lemma~\ref{lemma:perfectoid-almost}), we see that $B^+\wdh{\otimes}_{A_\infty^+} B^+ \to (B \wdh{\otimes}_{A_\infty} B)^+$ is an almost isomorphism (while, a priori, the latter group is the integral closure of the former one inside $B \wdh{\otimes}_{A_\infty} B$). In particular, 
    \[
    B^+/p \otimes_{A^+_\infty/p} B^+/p \simeq^a  (B \wdh{\otimes}_{A_\infty} B)^+/p(B \wdh{\otimes}_{A_\infty} B)^+.
    \]
    Thus 
    \[
    \rm{H}^0\left(\Spd\left(B \wdh{\otimes}_{A_\infty} B, (B \wdh{\otimes}_{A_\infty} B)^+\right)_v, \cal{E}\right)\simeq^a \left((B^+/p)^{\otimes^2_{A^+_\infty/p}}\right)^r
    \]
    and the two natural morphisms
    \[
    \rm{H}^0\left(\Spd(B, B^+)_v, \cal{E}\right) \otimes_{B^+/p} (B^+/p)^{\otimes^2_{A^+_\infty/p}} \to \rm{H}^0\left(\Spd\left(B \wdh{\otimes}_{A_\infty} B, (B \wdh{\otimes}_{A_\infty} B)^+\right)_v, \cal{E}\right)
    \]
    are almost isomorphisms. We use the sheaf condition and the previous discussion to get the following almost exact sequence
    \[
    0\to \rm{H}^0(\Spd\left(A_\infty, A^+_\infty\right)_v, \cal{E}) \to \rm{H}^0\left(\Spd(B, B^+)_v, \cal{E}\right) \to \rm{H}^0\left(\Spd(B, B^+)_v, \cal{E}\right) \otimes_{B^+/p} \left((B^+/p)^{\otimes 2}\right).
    \]

    Theorem~\ref{thm:faithuflly-flat descent} applied to the almost faithfully flat morphism $A^+_\infty/pA^+_\infty \to B^+/pB^+$ implies that the natural morphism
    \begin{equation}\label{eqn:almost-iso-1}
    \rm{H}^0\left(\Spd\left(A_\infty, A^+_\infty\right)_v, \cal{E}\right)\otimes_{A_\infty^+/p} B^+/p \to \rm{H}^0\left(\Spd(B, B^+)_v, \cal{E}\right)
    \end{equation}
    is an almost isomorphism. By the computation above, we know that $\rm{H}^0\left(\Spd(B, B^+)_v, \cal{E}\right)$ is almost faithfully flat and almost finitely presented over $B^+/pB^+$. Thus, the faithfully flat descent for flatness and almost finitely presented modules (see Lemma~\ref{almost-flat-descent} and Lemma~\ref{lemma:almost-flat-descent-flat}) implies that $\rm{H}^0\left(\Spd\left(A_\infty, A^+_\infty\right)_v, \cal{E}\right)$ is almost faithfully flat and almost finitely presented over $A^+_\infty/pA^+_\infty$. \smallskip
    
    {\it Step~$2$. $\rm{H}^0\left(\rm{Spd}(A_\infty, A_\infty^+)_v, \cal{E}\right)$ almost commutes with base change:} By the proof of \cite[Proposition 6.18]{Sch0} (and Lemma~\ref{lemma:perfectoid-almost}), we know that $\Spa(B, B^+)\times_{\Spa(A_\infty, A_\infty^+)} \Spa(D, D^+)$ exists as an adic space and is represented by $\Spa(R, R^+)$ for a perfectoid pair $(R, R^+)$ such that 
    \begin{equation}\label{eqn:almost-base-change}
        B^+/p \otimes_{A^+_\infty/p} D^+/p \to R^+/p
    \end{equation} 
    is an almost isomorphism. Thus, the proof of Step~$1$ and (\ref{eqn:almost-base-change}) imply that 
    \[
    \rm{H}^0\left(\Spd(D, D^+)_v, \cal{E}\right) \otimes_{A_\infty^+/p} B^+/p \to \rm{H}^0\left(\Spd(R, R^+)_v,\cal{E}\right)
    \]
    is an almost isomorphism. Now we wish to show that the natural morphism
    \[
    \rm{H}^0\left(\Spd(A_\infty, A^+_\infty)_v, \cal{E} \right) \otimes_{A_\infty^+/p} D^+/p \to \rm{H}^0\left(\Spd(D, D^+)_v, \cal{E} \right)
    \]
    is an almost isomorphism. By the faithfully flat descent, it suffices to check after tensoring against $B^+/p$ over $A_\infty^+/p$. Therefore, we use (\ref{eqn:almost-iso-1}) and (\ref{eqn:almost-base-change}) to see that it suffices to show that 
    \[
    \rm{H}^0\left(\Spd(B, B^+)_v, \cal{E}\right) \otimes_{B^+/p} R^+/p \to \rm{H}^0\left(\Spd(R, R^+)_v, \cal{E}\right)
    \]
    is an almost isomorphism. Now Lemma~\ref{lemma:first-properties-structure-sheaves}~\ref{lemma:first-properties-structure-sheaves-1} almost identifies (in the technical sense) this morphism with the identity moprhism
    \[
    (B^+/pB^+)^r \otimes_{B^+/p} R^+/p \to (R^+/pR^+)^r
    \]
    since $\cal{E}|_{\Spd(B, B^+)}$ is a trivial $\O^+/p$-vector bundle of rank $r$. This map is clearly an isomorphism. \smallskip
    
    {\it Step~$3$. $\rm{H}^i\left(\rm{Spd}(A_\infty, A_\infty^+)_v, \cal{E}\right)$ is almost zero for $i> 0$:} As in Step~$1$, we use that 
    \[
        \Spa(B, B^+) \to \Spa(A_\infty, A_\infty^+)
    \]
    is a finite \'etale morphism of affinoid perfectoid spaces to conclude that all fiber products 
    \[
    \Spa(B, B^+)^{j/\Spd(A_\infty, A_\infty^+)}
    \]
    are represented by affinoid perfectoid spaces $\Spa(B_j, B_j^+)$ and the natural morphisms 
    \[
    (B^+/pB^+)^{\otimes^j_{A_\infty^+/pA_\infty^+}} \to B_j^+/pB_j^+
    \]
    are almost isomorphisms. Since each restriction $\cal{E}|_{\Spd(B_j, B_j^+)}$ is trivial, Lemma~\ref{lemma:first-properties-structure-sheaves}~\ref{lemma:first-properties-structure-sheaves-1} ensures that higher cohomology of $\cal{E}$ on $\Spd(B_j, B_j^+)$ almost vanish. Thus, $\bf{R}\Gamma\left(\Spd(A_\infty, A_\infty^+)_v, \cal{E}\right)$ is almost isomorphic to the \v{C}ech complex associated to the covering $\Spd(B, B^+) \to \Spd(A_\infty, A_\infty^+)$. Step~$2$ implies that this complex is almost isomorphic to the standard Amitsur complex 
    \[
    0 \to M_{\cal{E}} \to M_{\cal{E}}\otimes_{A_\infty^+/p} B^+/p \to  M_{\cal{E}}\otimes_{A_\infty^+/p} B^+/p \otimes_{A_\infty^+/p} B^+/p \to \dots 
    \]
    Almost exactness of this complex follows from Lemma~\ref{lemma:amitsur-acyclic}. 
\end{proof}

\subsection{Strictly totally disconnected covers of affinoids}\label{section:strictly-totally-disconnected-covers}  
  
The main goal of this section is to eliminate almost mathematics in Lemma~\ref{lemma:almost-no-higher-coh-small} under some stronger assumptions on $A_\infty$. \smallskip

\begin{setup}\label{set-up:covers-strong} We fix 
\begin{enumerate}[label=\textbf{(\arabic*)}]
    \item a $p$-adic perfectoid field $K$ with its rank-$1$ open and bounded valuation ring $\O_K$ and a good pseudo-uniformizer $\varpi\in \O_K$ (we always do almost mathematics with respect to the ideal $\m= \bigcup_{n}\varpi^{1/p^n}\O_K = K^{\circ \circ}$),
    \item an affine admissible formal $\O_K$-scheme $\X=\Spf A_0$ with adic generic fiber $X=\Spa(A, A^+)$;
    \item a strictly totally disconnected affinoid perfectoid space  $\Spa(A_\infty, A_\infty^+)$ (see Definition~\ref{defn:totally-disconnected}) with a morphism 
    \[
    \Spa(A_\infty, A_\infty^+) \to \Spa(A, A^+)
    \]
    such that $\Spd(A_\infty, A_\infty^+) \to \Spd(A, A^+)$ is a $v$-covering and all fiber products
    \[
    \Spd(A_\infty, A_\infty^+)^{j/\Spd(A, A^+)}
    \]
    are strictly totally disconnected affinoid perfectoid spaces.
\end{enumerate}
\end{setup}

\begin{cor}\label{cor:restrict-strong} Under the assumption of Set-up~\ref{set-up:covers-strong}, let $\mf\colon \Spf B_0 \to \Spf A_0$ be an \'etale morphism of admissible affine formal $\O_K$-schemes.  Then 
\[
(B_\infty, B_\infty^+)\coloneqq \left(\left(B_0\wdh{\otimes}_{A_0} A_\infty^+\right)[1/p], B_0\wdh{\otimes}_{A_0}A_\infty^+\right)
\]
is a perfectoid pair and $\Spa(B_\infty, B_\infty^+)$ is a strictly totally disconnected (affinoid) perfectoid space.
\end{cor}
\begin{proof}
    Corollary~\ref{cor:restrict} already implies that $\Spa(B_\infty, B_\infty^+)$ is an affinoid perfectoid space. Moreover, Lemma~\ref{lemma:no-integral-closure} implies that
    \[
    \Spa(B_\infty, B_\infty^+) \simeq \Spa(B, B^+) \times_{\Spa(A, A^+)} \Spa(A_\infty, A_\infty^+),
    \]
    where $\Spa(B, B^+)$ is the generic fiber of $\Spf B_0$. So $\Spa(B_\infty, B_\infty^+) \to \Spa(A_\infty, A_\infty^+)$ is an \'etale morphism, thus the claim follows from \cite[Lemma 7.19]{Sch2}.
\end{proof}

\begin{lemma}\label{lemma:no-higher-coh-small-strong} Under the assumption of Set-up~\ref{set-up:covers-strong}, let $M_\cal{E}$ be an $A_\infty^+/pA_\infty^+$-module 
\[
M_\cal{E}\coloneqq \rm{H}^0\left(\Spd(A_\infty, A_\infty^+)_v, \cal{E}\right).
\]
Then $M_\cal{E}$ is a finite projective $(A_\infty^+/p)^r$-module, and for every morphism $\Spa(D, D^+) \to \Spa(A_\infty, A_\infty^+)$ of strictly totally disconnected affinoid perfectoid spaces, the natural morphism
\[
M_\cal{E}\otimes_{A_\infty^+/p} D^+/p \to \rm{H}^0\left(\Spd(D, D^+)_v, \cal{E}\right)
\]
is an isomorphism. Moreover, 
\[
\rm{H}^i\left(\Spd(A_\infty, A_\infty^+)^{j/\Spd(A, A^+)}_v, \cal{E}\right)\simeq 0
\]
for $i, j \geq 1$.
\end{lemma}
\begin{proof}
    Lemma~\ref{lemma:trivial-on-strictly-totally-disconnected} implies that we can replace $\Spa(A_\infty, A_\infty^+)$ by a finite clopen decomposition to assume\footnote{At this step, the map $\Spd(A_\infty, A_\infty^+) \to \Spd(A, A^+)$ might not be a $v$-covering anymore. But this will not matter for the rest of the proof.} that $\cal{E}|_{\Spd(A_\infty, A_\infty^+)} \simeq (\O_{\Spd(A_\infty, A^+_\infty)}^+/p)^r$ for some integer $r$. Then Corollary~\ref{cor:no-v-cohomology} implies that $M_{\cal{E}} \simeq (A_\infty^+/p)^r$. The same applies to $\cal{E}|_{\Spd(D, D^+)}$, therefore the natural morphism
    \[
    M_{\cal{E}} \otimes_{A_\infty^+/p} D^+/p = (A^+_\infty/p)^r \otimes_{A_\infty^+/p} D^+/p \to (D^+/p)^r
    \]
    is clearly an isomorphism. Furthermore, Corollary~\ref{cor:no-v-cohomology} implies that
    \[
    \rm{H}^i\left(\Spd(A_\infty, A_\infty^+)^{j/\Spd(A, A^+)}_v, \cal{E}\right) \simeq 0
    \]
    for $i, j\geq 1$ because we assume that all fiber products 
    \[
    \Spd(A_\infty, A_\infty^+)^{j/\Spd(A, A^+)}
    \]
    are representable by strictly totally disconnected (affinoid) perfectoid spaces.
\end{proof}

\begin{cor}\label{cor:restrict-3} Under the assumption of Set-up~\ref{set-up:covers-strong}, let $\mf\colon \Spf B_0\to \Spf A_0$ be an \'etale morphism, and let $(B_\infty, B^+_\infty)$ be the perfectoid pair from Corollary~\ref{cor:restrict}. Then the natural morphism 
\[
\Gamma\left(\Spd(A_\infty, A^+_\infty)^{j/\Spd(A, A^+)}_{v}, \cal{E}\right) \otimes_{A_0/p A_0} B_0/p B_0 \to \Gamma\left(\Spd(B_\infty, B^+_\infty)^{j/\Spd(B, B^+)}_{v}, \cal{E}\right).
\]
is an isomorphism for $j\geq 1$.
\end{cor} 
\begin{proof}
    For $j=1$, the result follows from Lemma~\ref{lemma:no-higher-coh-small-strong} and Corollary~\ref{cor:restrict-strong}. For $j>1$, we know that $X_j\coloneqq \Spd(A_\infty, A_\infty^+)^{j/\Spd(A, A^+)}$ is represented by a strictly totally disconnected perfectoid space. The morphism $X_j \to \Spd(A, A^+)$ defines a strictly totally disconnected perfectoid space $X_j^\sharp$ with a morphism $X_j^\sharp \to \Spa(A, A^+)$. One checks that $X_j^\sharp \to \Spa(A, A^+)$ satisfies the assumptions of Set-up~\ref{set-up:covers-strong}, so we can replace $\Spa(A_\infty, A_\infty^+)$ with $X^\sharp_j$ to reduce the case of $j>1$ to the case $j=1$. 
\end{proof}

\begin{cor}\label{restrict-2} Under the assumption of Set-up~\ref{set-up:covers-strong}, let $\mf\colon \Spf B_0\to \Spf A_0$ be an \'etale morphism, and let $\Spa(B, B^+)$ be the adic generic fiber of $\Spf(B_0)$. Then the natural morphism 
\[
\rm{H}^i\left(\Spd(A, A^+)_{v}, \cal{E}\right) \otimes_{A_0/p A_0} B_0/p B_0 \to \rm{H}^i\left(\Spd(B, B^+)_{v}, \cal{E}\right).
\]
is an isomorphism for $i\geq 0$.
\end{cor} 
\begin{proof}
Arguing as in the proof of Corollary~\ref{cor:restrict-3}, we see that Lemma~\ref{lemma:no-higher-coh-small-strong} implies that 
    \[
    \rm{H}^i\left(\Spd(A_\infty, A_\infty^+)^{j/\Spd(A, A^+)}_v, \cal{E} \right)\simeq 0
    \]
for $i, j \geq 1$. Therefore, cohomology groups $\rm{H}^i(\Spd(A, A^+)_v, \cal{E})$ can be computed via cohomology of the \v{C}ech complex associated to the covering $\Spd(A_\infty, A_\infty^+) \to \Spd(A, A^+)$. By Corollary~\ref{cor:restrict-strong}, the same applies to $\Spa(B, B^+)$ and the \v{C}ech complex associated to the covering $\Spd(B_\infty, B_\infty^+) \to \Spd(B, B^+)$. Therefore, the claim follows from Corollary~\ref{cor:restrict-3}. 
\end{proof}

\begin{cor}\label{cor:field-extension} Under the assumption of Set-up~\ref{set-up:covers}, let $K\subset C$ be a completed algebraic closure of $K$, and $\Spa(A_C, A_C^+)=\Spa(A, A^+)\times_{\Spa(K, \O_K)} \Spa(C, \O_C)$. Then the natural morphism 
\[
\rm{H}^i\left(\Spd(A, A^+)_{v}, \cal{E}\right) \otimes_{\O_K/p} \O_C/p \to \rm{H}^i\left(\Spd(A_C, A_C^+)_{v}, \cal{E}\right).
\]
is an almost isomorphism.
\end{cor} 
\begin{proof}
    The proof is similar to that of Corollary~\ref{cor:restrict-3} and Corollary~\ref{restrict-2}. The only change we need to make is that the fiber product
    \[
    \Spa(A_\infty, A_\infty^+)\times_{\Spa(K, \O_K)} \Spa(L, \O_L)
    \]
    is a strictly totally disconnected affinoid perfectoid space with the $+$-ring \emph{almost} isomorphic to $A_\infty^+\wdh{\otimes}_{\O_K} \O_L$. The strictly totally disconnected claim follows from \cite[Lemma 7.19]{Sch2} and the almost computation of the $+$-ring follows from the proof of \cite[Proposition 6.18]{Sch1}.  
\end{proof}

\subsection{Perfectoid torsors}\label{section:torsors}

We apply the results of Section~\ref{section:perfectoid-covers} to certain pro-\'etale covers of $\Spa(A, A^+)$ to see that the computation of $v$-cohomology groups can often be reduced to the computation of certain continuous cohomology groups. To make this precise, we need to define the notion of a $G$-torsor under a pro-finite group $G$.

\begin{defn} A {\it $v$-sheaf $\ud{G}$ associated to a pro-finite group} $G$ is a $v$-sheaf $\ud{G}\colon \rm{Perf}^{\rm{op}} \to \rm{Sets}$ such that $\ud{G}(S)=\rm{Hom}_{\rm{cont}}(|S|, G)$. \smallskip

A morphism of $v$-sheaves $X \to Y$ is a {\it $\ud{G}$-torsor} if it is a $v$-surjection and there is an action $a\colon \ud{G}\times X \to X$ over $Y$ such that the morphism $a\times_Y p_2\colon \ud{G}\times X \to X\times_Y X$ is an isomorphism, where $p_2\colon \ud{G}\times X \to X$ is the canonical projection.
\end{defn}

\begin{rmk} If a pro-finite group $G$ is a cofiltered limit of finite groups $G\simeq \lim_I G_i$, then $\ud{G} \simeq \lim_I \ud{G}_i$.
\end{rmk}

Now we can formulate the precise set-up we are going to work in. 

\begin{setup}\label{set-up:torsors} We fix
\begin{enumerate}[label=\textbf{(\arabic*)}]
    \item a $p$-adic perfectoid field $K$ with its rank-$1$ open and bounded valuation ring $\O_K$ and a good pseudo-uniformizer $\varpi\in \O_K$ (we always do almost mathematics with respect to the ideal $\m= \bigcup_{n}\varpi^{1/p^n}\O_K = K^{\circ \circ}$);
    \item an admissible formal $\O_K$-scheme $\X = \Spf A_0$ with adic generic fiber $X= \Spa(A, A^+)$; 
    \item a morphism $(A, A^+) \to (A_\infty, A^+_\infty)$ such $(A_\infty, A^+_\infty)$ is a perfectoid pair and $\Spd(A_\infty, A_\infty^+) \to \Spd(A, A^+)$ is a $\ud{\Delta}_\infty$-torsor under a pro-finite group $\Delta_\infty$;
    \item a small $\O_{X^\diam}^+/p$-vector bundle $\cal{E}$.
\end{enumerate}
\end{setup}

We start the section by studying the structure of the fiber products $\Spd(A_\infty, A_\infty^+)^{j/\Spd(A, A^+)}$ for $j\geq 1$. For a general $v$-cover, we cannot say much about these fiber products. However, we have much more control in the case of $\ud{G}$-torsors.

\begin{lemma}\label{lemma:torsor-over-perfectoid} Under the assumption of Set-up~\ref{set-up:torsors}, the fiber product $\Spd(A_\infty, A_\infty^+)^{j/\Spd(A, A^+)}$ is represented by an affinoid perfectoid space\footnote{Recall that $\Spd(A_\infty, A_\infty^+)$ is itself represented by an affinoid perfectoid $\Spa(A_\infty^\flat, A_\infty^{\flat, +})$.} $\Spa(T_j, T_j^+)$ for every $j\geq 0$. Moreover, for every $j\geq 0$,  
\[
\left(T_j, T_j^+\right)  \simeq \left(\rm{Map}_{\rm{cont}}(\Delta_\infty^{j-1}, A_\infty^\flat), \rm{Map}_{\rm{cont}}(\Delta_\infty^{j-1}, A_\infty^{\flat, +}) \right)
\]
and $T_j^{\sharp, +}/p T_{j}^{\sharp, +} \simeq T_j^+/\varpi^\flat T_j^+ \simeq \rm{Map}_{\rm{cont}}(\Delta_\infty^{j-1}, A_\infty^{+}/pA_\infty^+)$.
\end{lemma}
\begin{proof}
We first show that $\Spd(A_\infty, A_\infty^+)^{j/\Spd(A, A^+)}$ are representable by affinoid perfectoid spaces. We write a presentation of $\Delta_\infty = \lim_I \Delta_i$ as a cofiltered limit of finite groups. Since $\Spd(A_\infty, A_\infty^+) \to \Spd(A, A^+)$ is a $\ud{\Delta}_\infty$-torsor, we get
    \begin{align*}
    \Spd(A_\infty, A_\infty^+)^{j/\Spd(A, A^+)} & \simeq \Spd(A, A^+)\times \ud{\Delta}_\infty^{j-1} \\
    &\simeq \lim_I\left(\Spa(A_\infty^\flat, A_\infty^{\flat, +})\times \ud{\Delta}_i^{j-1} \right) \\
    & \simeq \lim_I \left(\Spa\left(\rm{Map}(\Delta^{j-1}_i, A_\infty^\flat), \rm{Map}(\Delta^{j-1}_i, A_\infty^{\flat, +})\right)\right) 
    \end{align*}
    is a cofiltered limit of affinoid perfectoid spaces, so it is an affinoid perfectoid space $\Spa(T_j, T_j^+)$ by \cite[Proposition 6.5]{Sch2}. Moreover, {\it loc. cit.} implies that $T_j^+$ is equal to the $\varpi^\flat$-adic completion of the filtered colimit $\colim_I \rm{Map}(\Delta^{j-1}_i, A_\infty^{\flat, +})$ and $T_j=T_j^+[\frac{1}{\varpi^\flat}]$. In particular, we already see that
    \begin{align*}
    T_j^{\sharp, +}/p T_{j}^{\sharp, +} \simeq T^+_j/(\varpi)^\flat T^+_j & \simeq \left(\colim_I \rm{Map}(\Delta^{j-1}_i, A_\infty^{\flat, +})\right)/(\varpi)^\flat \\
    & \simeq \colim_I \rm{Map}(\Delta^{j-1}_i, A_\infty^{\flat, +}/(\varpi)^\flat A_\infty^{\flat, +}) \\
    & \simeq \colim_I \rm{Map}(\Delta^{j-1}_i, A_\infty^{+}/\varpi A_\infty^{+})\\
    & \simeq \colim_I \rm{Map}(\Delta^{j-1}_i, A_\infty^{+}/p A_\infty^{+})\\
    & \simeq \rm{Map}_{\rm{cont}}(\Delta^{j-1}_\infty, A_\infty^{+}/p A_\infty^{+}).
    \end{align*}
    Now we compute $T_j^+$ and $T_j$. We start with $T_j^+$:
    \begin{align*}
    T_j^+ &\simeq \lim_n \left(\colim_I \rm{Map}(\Delta^{j-1}_i, A_\infty^{\flat, +})/(\varpi^\flat)^n\right) \\
    & \simeq \lim_n \left(\colim_I \rm{Map}(\Delta^{j-1}_i, A_\infty^{\flat, +}/(\varpi^\flat)^n A_\infty^{\flat, +})\right) \\
    & \simeq \lim_n \rm{Map}\left(\Delta_\infty^{j-1}, A_\infty^{\flat, +}/(\varpi^\flat)^n A_\infty^{\flat, +}\right) \\
    & \simeq \rm{Map}_{\rm{cont}}\left(\Delta_\infty^{j-1}, \lim_n  A_\infty^{\flat, +}/(\varpi^\flat)^n A_\infty^{\flat, +}\right)\\
    & \simeq \rm{Map}_{\rm{cont}}\left(\Delta_\infty^{j-1},  A_\infty^{\flat, +}\right). \end{align*}
    Since $\Delta_\infty$ is compact and $ A_\infty^{\flat}\simeq  A_\infty^{\flat, +}[\frac{1}{\varpi^\flat}]$, we also have
    \begin{align*}
        T_j & \simeq T_j^+[1/\varpi^\flat] \\
            & \simeq \colim_{\times \varpi^\flat, n} \rm{Map}_{\rm{cont}}\left(\Delta_\infty^{j-1},  A_\infty^{\flat, +}\right)\\
            & \simeq \rm{Map}_{\rm{cont}}\left(\Delta_\infty^{j-1},  \colim_{\times \varpi^\flat} A_\infty^{\flat, +}\right)\\
            & \simeq \rm{Map}_{\rm{cont}}\left(\Delta_\infty^{j-1}, A_\infty^{\flat}\right)
    \end{align*}
    finishing the proof. 
\end{proof}

\begin{warning}\label{warning:perfectoid-only-after-diamondification} The fiber product $\Spa(A_\infty, A_\infty^+) \times_{\Spa(A, A^+)} \Spa(A_\infty, A_\infty^+)$ in the category of adic spaces is often not a perfectoid space. This already happens for $\Spa(A, A^+)=\Spa(\Q_p, \Z_p)$ and $\Spa(A_\infty, A_\infty^+) = \Spa(\Q_p(\mu_{p^{\infty}})^{\wedge}, \Z_p[\mu_{p^{\infty}}]^{\wedge})$. However, Lemma~\ref{lemma:torsor-over-perfectoid} implies that the diamond $\Spd(A_\infty, A_\infty^+) \times_{\Spd(A, A^+)} \Spd(A_\infty, A_\infty^+)$ is always represented by an affinoid perfectoid space. 
\end{warning}

Note that since $\Spd(A_\infty, A_\infty^+) \to \Spd(A, A^+)$ is a $\ud{\Delta}_\infty$-torsor, there is a canonical continuous $A^+$-linear action of $\Delta_\infty$ on $A_\infty^+$. Now we want to relate $v$-cohomology groups of $\cal{E}$ to the continuous group cohomology of $\Delta_\infty$. This is done in the following lemmas:

\begin{lemma}\label{lemma:no-higher-coh-torsor} Under the assumption of Set-up~\ref{set-up:torsors}, we define $M_{\cal{E}}$ to be the $A_\infty^+/p$-module $\rm{H}^0\left(\Spd(A_\infty, A_\infty^+)_v, \cal{E}\right)$. Then $M_{\cal{E}}$ is almost faithfully flat, almost finitely presented $A_\infty^+/p$-module, and \[
\rm{H}^0(\Spd(A_\infty, A_\infty^+)^{j/\Spd(A, A^+)}_v, \cal{E}) \simeq^a \rm{Map}_{\rm{cont}}(\Delta_\infty^{j-1}, M_{\cal{E}}) \simeq^a \rm{Map}_{\rm{cont}}(\Delta_\infty^{j-1}, (M_{\cal{E}}^a)_!),
\]
\[
\rm{H}^i(\Spd(A_\infty, A_\infty^+)^{j/\Spd(A, A^+)}_v, \cal{E}) \simeq^a 0
\]
for every $i, j\geq 1$.
\end{lemma}
\begin{proof}
    Lemma~\ref{lemma:torsor-over-perfectoid} implies that all fiber products $\Spd(A_\infty, A_{\infty}^+)^{j/\Spd(A, A^+)}$ satisfy the assumptions of Lemma~\ref{lemma:almost-no-higher-coh-small}. Thus, Lemma~\ref{lemma:almost-no-higher-coh-small} and the computation of fiber products in Lemma~\ref{lemma:torsor-over-perfectoid} imply that 
    \[
    \rm{H}^i\left(\Spd(A_\infty, A_\infty^+)^{j/\Spd(A, A^+)}_v, \cal{E}\right) \simeq^a 0
    \]
    for every $i, j\geq 1$, and the natural morphism
    \[
    M_{\cal{E}}\otimes_{A_\infty^+/p} \rm{Map}_{\rm{cont}}\left(\Delta_\infty^{j-1}, A_\infty^+/p\right) \to \rm{H}^0\left(\Spd(A_\infty, A_\infty^+)^{j/\Spd(A, A^+)}_v, \cal{E}\right)
    \]
    is an almost isomorphism for every $j\geq 1$. Thus, it suffices to show that the natural morphism
    \[
    M_{\cal{E}}\otimes_{A_\infty^+/p} \rm{Map}_{\rm{cont}}\left(\Delta_\infty^{j-1}, A_\infty^+/p\right) \to \rm{Map}_{\rm{cont}}(\Delta_\infty^{j-1}, M_{\cal{E}})
    \]
    is an isomorphism. This can be done by writing $\Delta_\infty = \lim_I \Delta_i$ and reducing to the case of a finite group similarly to the proof of Lemma~\ref{lemma:torsor-over-perfectoid}. The almost isomorphism 
    \[
    \rm{Map}_{\rm{cont}}(\Delta_\infty^{j-1}, M_{\cal{E}}) \simeq^a \rm{Map}_{\rm{cont}}(\Delta_\infty^{j-1}, (M_{\cal{E}}^a)_!)
    \]
    is achieved similarly using that $(-)_!$ commutes with colimits being a left adjoint functor. 
\end{proof}

\begin{lemma} Under the assumption of Set-up~\ref{set-up:torsors}, we define $M_{\cal{E}}$ to be the $A_\infty^+/pA_\infty^+$-module $\rm{H}^0\left(\Spd(A_\infty, A_\infty^+)_v, \cal{E}\right)$. Then there is a canonical continuous action of $\Delta_\infty$ on $(M_{\cal{E}}^a)_!$ compatible with the action of $\Delta_\infty$ on $A_\infty^+/p$, i.e. $g(am)=g(a)g(m)$ for any $a\in A_\infty^+/p$ and $m\in M_\cal{E}$.  
\end{lemma}
\begin{proof}
    Lemma~\ref{lemma:torsor-over-perfectoid} ensures that the fiber product $\Spd(A_\infty, A_\infty^+)\times_{\Spd(A, A^+)} \Spd(A_\infty, A_\infty^+)$ is represented by an affinoid perfectoid space $\Spa(T_2, T_2^+)$ of characteristic $p$. Therefore, we can uniquely write it as $\Spd(S, S^+)$ for an untilt of $(T_2, T_2^+)$ corresponding to the morphism $\Spa(T_2, T_2^+)\to \Spd(A, A^+)\to \Spd(\Q_p, \Z_p)$.\smallskip
    
    Lemma~\ref{lemma:almost-no-higher-coh-small} implies that the descent data for the sheaf $\cal{E}$ provide us with an $(S^+/pS^+)^a$-isomorphism
    \[
    \left(S^+/p\right)^a \otimes_{\left(A_\infty^+/p\right)^a} \left(M_\cal{E}\right)^a \to \left(M_\cal{E}\right)^a \otimes_{\left(A_\infty^+/p\right)^a} \left(S^+/p\right)^a
    \]
    satisfying the cocycle condition. By Corollary~\ref{inner-hom-algebra}~\ref{inner-hom-algebra-2}, this defines an $\left(A_\infty^+/p\right)^a$-linear morphism
    \[
    \left(M_\cal{E}\right)^a \to \left(M_\cal{E}\right)^a \otimes_{\left(A_\infty^+/p\right)^a} \left(S^+/p\right)^a.
    \]
    By Lemma~\ref{lemma:almost-no-higher-coh-small} and Lemma~\ref{lemma:no-higher-coh-torsor}, this is equivalent to an $\left(A_\infty^+/p\right)^a$-linear morphism
    \[
    \left(M_\cal{E}\right)^a \to \rm{Map}_{\rm{cont}}\left(\Delta_\infty, (M^a_\cal{E})_!\right)^a.
    \]
    By Lemma~\ref{lemma:adjoint-almost}~\ref{adjoint-almost-3}, this is the same as an $\left(A_\infty^+/pA_\infty^+\right)$-linear morphism
    \[
    \varphi\colon \left(M^a_\cal{E}\right)_! \to \rm{Map}_{\rm{cont}}\left(\Delta_\infty, \left(M^a_\cal{E}\right)_!\right).
    \]
    This defines a morphism
    \[
    \gamma\colon \Delta_\infty \to \rm{Hom}_{A_\infty^+/p}\left(\left(M_\cal{E}\right)_!, \left(M_\cal{E}\right)_!\right)
    \]
    by the rule 
    \[
    \gamma(g)(m)=(\phi(m))(g).
    \]
    One checks that the cocycle condition translates into the statement that $\gamma$ is a group homomorphism, i.e. it defines an action of $\Delta_\infty$. Similarly, one checks that $A_\infty^+/p$-linearity of $\phi$ translates into the fact that this action is compatible with the action on $A_\infty^+/p$. And continuity of $\phi$ translates into the fact that $\gamma$ defines a continuous action, i.e., the natural morphism
    \[
    \colim_{U_i\triangleleft \Delta_\infty, \text{open}} (M^a_{\cal{E}})_!^{U_i} \to (M^a_{\cal{E}})_!^{\Delta_\infty}
    \]
    is an isomorphism. 
\end{proof}

\begin{cor}\label{cor:v-cohomology-group-cohomology} Under the assumption of Set-up~\ref{set-up:torsors}, we define $M_{\cal{E}}$ to be an $A_\infty^+/p$-module $\rm{H}^0\left(\Spd(A_\infty, A_\infty^+)_v, \cal{E}\right)$. Then 
\[
\rm{H}^i(\Spd(A, A^+)_v, \cal{E})\simeq^a \rm{H}^i_{\rm{cont}}(\Delta_\infty, (M_{\cal{E}}^a)_!).
\]
\end{cor}
\begin{proof}
    Lemma~\ref{lemma:no-higher-coh-torsor} implies that 
    \[
    \rm{H}^i\left(\Spd(A_\infty, A_\infty^+)^{j/\Spd(A, A^+)}_v, \cal{E} \right)\simeq^a 0
    \]
    for $i, j \geq 1$. Therefore, the cohomology groups $\rm{H}^i(\Spd(A, A^+)_v, \cal{E})$ can be almost computed via cohomology of the \v{C}ech complex associated to the covering $\Spd(A_\infty, A_\infty^+) \to \Spd(A, A^+)$. Moreover, Lemma~\ref{lemma:no-higher-coh-torsor} also implies that the terms of this complex can be almost identified with the bar complex computing continuous cohomology of the pro-finite group $\Delta_\infty$ with coefficients in the discrete module $(M^a_{\cal{E}})_!$. We leave it to the reader to verify that the differentials in the \v{C}ech complex coincide with the differentials in the bar complex computing continuous cohomology. 
\end{proof}

For future reference, we also discuss the following base change result:

\begin{lemma}\label{coh-flat-base-change} Let $G$ be a pro-finite group, and let $M$ be a discrete $R$-module that has a continuous $R$-linear action of $G$. Suppose that $R\to A$ is a flat homomorphism of rings. Then the canonical morphism $\rm{H}^i_{\rm{cont}}(G, M)\otimes_R A \to \rm{H}^i_{\rm{cont}}(G, M\otimes_R A)$ is an isomorphism for $i\geq 0$.
\end{lemma}
\begin{proof}
We first prove the claim for $\rm{H}^0$. Since $G$ acts on $M$ continuously, we can write $M=\colim_I M_i$ as a filtered colimit of $G$-stable $R$-submodules of $M$ such that the action of $G$ on $M_i$ factors through a finite group $G_i$. Since both $\rm{H}^0_{\rm{cont}}(G, -) \otimes_R A$ and $\rm{H}^0_{\rm{cont}}(G, - \otimes_R A)$ commute with filtered colimits, we can reduce to the case when the action of $G$ factors through a finite group quotient. In this case, the result is classical (see, for example, \cite[Exp.~V, Proposition 1.9]{SGA1}). \smallskip

In general, the result follows from the following sequence of isomorphisms
\begin{align*}
\rm{H}^i_{\rm{cont}}(G, M)\otimes_R A &\cong(\text{colim}_{H\triangleleft G, \text{open}}\rm{H}^i(G/H, M^H))\otimes_R A\\
& \simeq \text{colim}_{H\triangleleft G, \text{open}}(\rm{H}^i(G/H, M^H)\otimes_R A) \\
& \simeq \text{colim}_{H\triangleleft G, \text{open}}\rm{H}^i(G/H, M^H\otimes_R A)\\
& \simeq \text{colim}_{H\triangleleft G, \text{open}}\rm{H}^i(G/H, (M\otimes_R A)^H)\\
& \simeq \rm{H}^i_{\rm{cont}}(G, M\otimes_R A) \qedhere
\end{align*}
\end{proof}

\subsection{Nearby cycles are quasi-coherent}\label{universal}

We start the proof of Theorem~\ref{main-thm} and Theorem~\ref{thm:main-thm-small} in this Section. Namely, we show that the complex $\bf{R}\nu_*\cal{E}$ is quasi-coherent and commutes with \'etale base change for an $\O_{X^\diam}^+/p$-vector bundle $\cal{E}$.  The main idea is to apply the results of Section~\ref{section:perfectoid-covers} to a particular perfectoid covering of $X$. \smallskip

For the rest of this section, we fix a perfectoid $p$-adic field $K$ with a good pseudo-uniformizer $\varpi \in \O_K$ (see Definition~\ref{defn:good-unifor}). We always do almost mathematics with respect to the ideal $\m= \bigcup_{n}\varpi^{1/p^n}\O_K$. \smallskip 


\begin{lemma}\label{lemma:quasi-coherentness-easy} Let $\X=\Spf A_0$ an admissible affine formal $\O_K$-scheme with an affinoid generic fiber $X=\Spa(A, A^+)$, and let $\cal{E}$ be an $\O_{X^\diam}^+/p$-vector bundle. Then  $\rm{R}^i\nu_*\cal{E}$ is quasi-coherent for $i\geq 0$. More precisely, the natural morphism 
\[
\widetilde{\rm{H}^i(X_v^\diam, \cal{E})} \to \rm{R}^i\nu_*\cal{E}
\]
is an isomorphism for any $i\geq 0$.
\end{lemma}
\begin{proof}
    The universal property of the tilde-construction implies that we do have a natural morphism 
\[
c \colon \widetilde{\rm{H}^i(X^\diam_v, \cal{E})} \to \rm{R}^i\nu_*\cal{E}.
\]
Recall that $\rm{R}^i\nu_*\cal{E}$ is the sheafification of the presheaf defined by the rule
\[
\sU \mapsto \rm{H}^i(\sU_{K, v}^\diam, \cal{E}).
\]
Thus, in order to show that $c$ is an isomorphism, it suffices to show that the natural morphism
\[
\rm{H}^i(X^\diam_v,\cal{E}) \otimes_{A_0/p} (A_0/p)_f\to \rm{H}^i(\sU^\diam_{K, v}, \cal{E})
\]
is an isomorphism for any open formal subscheme $\Spf (A_0)_{\{f\}} \subset \Spf A_0$. We choose a covering $\Spa(A_\infty, A_\infty) \to \Spa(A, A^+)$ from Lemma~\ref{lemma:str-tot-disc-covering}. Then the result follows from Corollary~\ref{restrict-2} since $(A, A^+) \to (A_\infty, A_\infty^+)$ fits into Set-up~\ref{set-up:covers-strong}. 
\end{proof}

\begin{thm}\label{almost-quasi} Let $\X$ an admissible formal $\O_K$-scheme with adic generic fiber $X=\X_K$, and let $\cal{E}$ be an $\O_{X^\diam}^+/p$-vector bundle. Then  $\rm{R}^i\nu_*\cal{E}$ is quasi-coherent for $i\geq 0$. Furthermore, if $\mf\colon \Y \to \X$ is an \'etale morphism with generic fiber $f\colon Y\to X$, then the natural morphism
\[
\mf_0^*\left(\rm{R}^i\nu_{\X, *}\cal{E}\right) \to \rm{R}^i\nu_{\Y, *}\left(\cal{E}|_{Y^\diam_{v}}\right)
\]
is an isomorphism for any $i\geq 0$.
\end{thm}
\begin{proof}
    Both claims are local on $\X$ and $\Y$, so we can assume that $\X=\Spf A_0$ and $\Y=\Spf B_0$ are affine. Then quasi-coherence of $\rm{R}^i\nu_*(\cal{E})$ directly follows from Lemma~\ref{lemma:quasi-coherentness-easy}. In order to show that 
    \[
    \mf_0^*\left(\rm{R}^i\nu_{\X, *}\cal{E}\right) \to \rm{R}^i\nu_{\Y, *}\left(\cal{E}|_{Y^\diam_{v}}\right),
    \]
    it suffices to show that the natural morphism
    \[
    \rm{H}^i(X_v^\diam, \cal{E}) \otimes_{A_0/p} B_0/p \to \rm{H}^i(Y_v^\diam, \cal{E})
    \]
    is an isomorphism. This follows from Corollary~\ref{restrict-2} using the covering $\Spa(A_\infty, A_\infty^+) \to \Spa(A, A^+)$ from Lemma~\ref{lemma:str-tot-disc-covering}.
\end{proof}

For future reference, we also prove the following result:

\begin{lemma}\label{lemma:affinoid-base-change} Let $X=\Spa(A, A^+)$ be an affinoid rigid-analytic space over $K$, let $\cal{E}$ be an $\O_{X^\diam}^+/p$-vector bundle, and let $K\subset C$ be a completed algebraic closure of $K$. Then
\[
\rm{H}^i(X_v^\diam, \cal{E}) \otimes_{\O_K/p} \O_C/p \to \rm{H}^i(X_{C, v}^\diam, \cal{E})
\]
is an almost isomorphism.
\end{lemma}
\begin{proof}
    This follows directly from Corollary~\ref{cor:field-extension} using the covering $\Spa(A_\infty, A_\infty^+) \to \Spa(A, A^+)$ from Lemma~\ref{lemma:str-tot-disc-covering}.
\end{proof}

\subsection{Nearby cycles are almost coherent for smooth $X$ and small $\cal{E}$}\label{almost-coherent-complex}

The main goal of this section is to show that the complex $\bf{R}\nu_* \cal{E}$ has almost coherent cohomology sheaves for an admissible formal $\O_K$-scheme with {\it smooth} generic fiber. The main idea is to apply the results of Section~\ref{section:torsors} to a particular ``small'' perfectoid torsor cover of $X$, where one has good control over the structure group $\Delta_\infty$. \smallskip

For the rest of the section, we fix a $p$-adic perfectoid field $K$ with a good pseudo-uniformizer $\varpi \in \O_K$. We always do almost mathematics with respect to the ideal $\m= \bigcup_{n}\varpi^{1/p^n}\O_K$. \smallskip

Before we embark on the proof, we discuss the overall strategy of the proof. We proceed in four steps: first, we show the result for $\wdh{\bf{G}}_m^n$ and $\cal{E}=\O_{X^\diam}^+/p$, then we deduce the result for affine formal schemes such that the adic generic fiber admits a map to a torus $\bf{T}^n_{C}$ that is a composition of finite \'etale maps and rational embeddings. After that, we finish the proof for $\cal{E} = \O_{X^\diam}^+/p$ and a general smooth $X$ by choosing a ``good'' covering of $\X$, possibly after an admissible blow-up of $\X$. We reduce the general case to the case $\cal{E}=\O_{X^\diam}^+/p$ via Corollary~\ref{cor:different-vector-bundles-equivalent}. \smallskip

The main ingredient for the third step is Achinger's result (\cite[Proposition 6.6.1]{Ach1}) that any \'etale morphism $g\colon \Spa(A, A^+) \to \bf{D}^n_{K}$ can be replaced with a finite \'etale morphism 
\[
g'\colon \Spa(A, A^+) \to \bf{D}^n_{K}.
\]
The proof of this result in \cite{Ach1} is given only for rigid-analytic varieties over discretely valued non-archimedean fields, but we need to apply it in the perfectoid situation that is never discretely valued. So Appendix~\ref{achinger} provides the reader with a detailed proof of this result without any discreteness assumptions. \smallskip

Now we begin to realize the strategy sketched above. We consider $\X=\Spf \O_K\langle T_1^{\pm 1}, \dots, T_n^{\pm 1}\rangle$, and set $R^+\coloneqq \O_K\langle T_1^{\pm 1}, \dots, T_n^{\pm 1}\rangle$ and $R^+_m\coloneqq \O_K\langle T_1^{\pm 1/p^m}, \dots, T_n^{\pm 1/p^m}\rangle$. We note that the map $\Spf R^+_m \to \Spf R^+$ defines a $\mu_{p^m}^n$-torsor, thus $\mu_{p^m}^n$ continuously acts on $R^+_m$ by $R^+$-linear automorphisms. \smallskip

Now we consider the $R^+$-algebra 
\[
R^+_{\infty}=\O_K\langle T_1^{\pm 1/p^{\infty}}, \dots, T_n^{\pm 1/p^{\infty}}\rangle=\left(\text{colim}_n \ R^+_m \right)\widehat{}
\]
where $\ \widehat{} $ stands for the $p$-adic completion. It comes with a continuous $R^+$-linear action of the profinite group $\Delta_{\infty}\coloneqq \Z_p(1)^{n}=T_p(\mu_{p^{\infty}})$ on $R^+_{\infty}$. We trivialize $\Z_p(1)$ by choosing some compatible system of $p^i$-th roots of unity $(\zeta_p, \zeta_{p^2}, \zeta_{p^3}, \dots )$. To describe the action of $\Delta_{\infty}$ on $R^+_{\infty}$ we need the following definition:

\begin{defn} For any $a\in \Z[1/p]$, we define $\zeta^a$ as $\zeta^{ap^l}_{p^l}$ whenever $ap^l\in \Z$. It is clear to see that this definition does not depend on the choice of $l$.
\end{defn}

Essentially by definition, the $k$-th basis vector $\gamma_k\in \Delta_{\infty} \simeq \Z_p^{n}$ acts on $R^+_{\infty}$ as 
\[
\gamma_k(T_1^{a_1}\dots T_n^{a_n})=\zeta^{a_k}T_1^{a_1}\dots T_n^{a_n}. 
\]

\begin{lemma}\label{first-computation}\cite[Lemma 5.5]{Sch1} Let $R^+$, $R^+_{\infty}$ and $\Delta_{\infty}$ be as above. Then the cohomology groups $\rm{H}^i_{\rm{cont}}(\Delta_{\infty}, R^+_{\infty}/p)$ are almost coherent $R^+/p$-modules. And the natural map
\[
\rm{H}^i_{\rm{cont}}(\Delta_{\infty}, R^+_{\infty}/p) \otimes_{R^+/p} A^+/p\to \rm{H}^i_{\rm{cont}}(\Delta_{\infty}, R^+_{\infty}/p\otimes_{R^+/p} A^+/p) 
\]
is an isomorphism for a $p$-torsionfree $R^+$-algebra $A^+$ and $i\geq 0$.
\end{lemma}
\begin{proof}
We note that $R^+/p$ is an almost noetherian ring due Theorem~\ref{thm:top-ft-almost-noetherian}. Thus, Theorem~\ref{cor:almost-noetherian-almost-fg=almost-coh} implies that $\rm{H}^i_{\rm{cont}}(\Delta_{\infty}, R^+_{\infty}/p R^+_{\infty})$ is almost coherent if it is almost finitely generated. \smallskip

Now \cite[Lemma 7.3]{BMS1} says that $\bf{R}\Gamma_{\rm{cont}}(\Delta_{\infty}, R^+_{\infty}/p)$ is computed via the Koszul complex $K\left(R^+_{\infty}/p; \gamma_1 -1, \dots, \gamma_n -1\right)$. Then, similarly to \cite[Lemma 4.6]{Bhatt-spec}, we can write 
\[
K\left(R^+_{\infty}/p ; \gamma_1 -1, \dots, \gamma_n -1\right) = K\left(R^+/p; 0, 0, \dots, 0\right) \oplus \bigoplus_{\left(a_1, \dots, a_n\right)\in \left(\Z[1/p] \cap \left(0,1\right) \right)^n} K\left(R^+/p; \zeta^{a_1} -1, \dots, \zeta^{a_n} -1\right) 
\]

We observe that
\[
\rm{H}^i\left(K\left(R^+/p; 0, 0, \dots, 0\right)\right)=\wedge^i\left(R^+/p\right)
\]
is a free finitely presented $R^+/p$-module. For each $(a_1, \dots, a_n) \in \left(\Z[1/p] \cap \left(0,1\right)\right)^n$, we can assume that $a_1$ has the minimal $p$-adic valuation for the purpose of proving that 
\[
    K(R^+_{\infty}/p; \gamma_1 -1, \dots, \gamma_n -1)
\]
has almost finitely finitely generated cohomology groups. Then \cite[Lemma 7.10]{BMS1} implies that $\rm{H}^i\left(K\left(R^+/p; \zeta^{a_1} -1, \dots, \zeta^{a_n} -1\right)\right)$ is finitely presented over $R^+/p$ and $\zeta^{a_1}-1$-torsion module. Note that
\[
v_p(\zeta^{a_1}-1)=v_p(\zeta_{p^l}-1)=\frac{v(p)}{p^l-p^{l-1}} \to 0
\]
where $a_1=b/p^l$ with $\text{gcd}(b, p)=1$. Furthermore, for any $h\in \Z$, there are only finitely many indexes $(a_1, \dots, a_n) \in (\Z[1/p] \cap (0,1) )^n$ with $v_p(a_j)\geq h$. This implies that 
\[
\rm{H}^i_{\rm{cont}}(\Delta_{\infty}, R^+_{\infty}/p)=\rm{H}^i \left(K\left(R^+_{\infty}/p; \gamma_1 -1, \dots, \gamma_n -1\right)\right)
\]
is a finitely presented $R^+/p$-module up to any $\varpi^{1/p^n}$-torsion. In particular, this module is almost finitely presented. \smallskip

Now we show that $\rm{H}^i_{\rm{cont}}(\Delta_{\infty}, R^+_{\infty}/p)$ commutes with base change for any $\O_K$-flat algebra $A^+$. In order to show this, we observe that the $(R^+/p)[\Delta_{\infty}]$-module $R^+_{\infty}/p$ comes as a tensor product $M\otimes_{\O_K/p} R^+/p$ for the $(\O_K/p)[\Delta_{\infty}]$-module
\[
M\coloneqq \bigoplus_{(a_1, \dots, a_n) \in (\Z[1/p] \cap [0,1) )^n} (\O_K/p\O_K) T_1^{a_1} \dots T_n^{a_n}
\]
where the basis element $\gamma_k$ acts by 
\[
\gamma_k(T_1^{a_1}\dots T_n^{a_n})=\zeta^{a_k}T_1^{a_1}\dots T_n^{a_n}. 
\]
Therefore, the desired claim follows from a sequence of isomorphisms
\begin{align*}
\rm{H}^i_{\rm{cont}}(\Delta_{\infty}, R^+_{\infty}/p)\otimes_{R^+/p} A^+/p &\simeq \left( \rm{H}^i_{\rm{cont}}(\Delta_{\infty}, M)\otimes_{\O_K/p} R^+/p \right) \otimes_{R^+/p}  A^+/p \\
& \simeq \rm{H}^i_{\rm{cont}}(\Delta_{\infty}, M)\otimes_{\O_K/p} A^+/p \\
& \simeq \rm{H}^i_{\rm{cont}}(\Delta_{\infty}, M\otimes_{\O_K/p} A^+/p) \\
&  \simeq \rm{H}^i_{\rm{cont}}(\Delta_{\infty}, R^+_{\infty}/p  \otimes_{R^+/p} A^+/p ),
\end{align*}
where the third isomorphism uses Lemma~\ref{coh-flat-base-change}.
\end{proof}

Lemma~\ref{first-computation} combined with Corollary~\ref{cor:v-cohomology-group-cohomology} essentially settles the first step of our strategy. Now we move to the second step. We start with the following preliminary result:

\begin{lemma}\label{lemma:topologically-free-flat} Let $A_0$ be a topologically finitely presented $\O_K$-algebra, and $P$ a topologically free $A_0$-module, i.e. $P=\wdh{\bigoplus}_I A_0$ for some set $I$. Then $M$ is $A_0$-flat.
\end{lemma}
\begin{proof}
We start the proof by noting that \cite[\href{https://stacks.math.columbia.edu/tag/00M5}{Tag 00M5}]{stacks-project} guarantees that it suffices to show that $\rm{Tor}_1^{A_0}(P, M)=0$ for any finitely presented $A_0$-module $M$. We choose a presentation
\[
0 \to Q \to A_0^n \to M \to 0
\]
and observe that $Q$ is finitely presented because $A_0$ is coherent. So vanishing of $\rm{Tor}_1$ is equivalent to showing that
\[
P \otimes_{A_0} Q \to P\otimes_{A_0} A_0^n
\]
is injective.\smallskip

Now note that $Q[p^\infty]$, $A_0^n[p^\infty]$, and $M[p^\infty]$ are bounded by \cite[Lemma 7.3/7]{B}, so the same holds for $\bigoplus_I Q$, $\bigoplus_I A_0^n$, and $\bigoplus_I M$. Therefore, the usual $p$-adic completions of $\bigoplus_I Q$, $\bigoplus_I A_0^n$ and $\bigoplus_I M$ coincide with their derived $p$-adic completions. Since derived $p$-adic completion is exact (in the sense of triangulated categories) and coincides with the usual one on these modules, we get that the sequence
\[
0 \to \wdh{\bigoplus}_I Q \to \wdh{\bigoplus}_I A_0^n \to \wdh{\bigoplus}_I M \to 0
\]
is exact. \smallskip

Now we want to show that this short exact sequence is the
same as the sequence
\[
P \otimes_{A_0} Q \to P\otimes_{A_0} A_0^n \to P\otimes_{A_0} M \to 0.
\]
As a consequence, this will show that $P \otimes_{A_0} Q \to P\otimes_{A_0} A_0^n$ is injective. \smallskip

For each $A_0$-module $N$, there is a canonical map 
\[
P\otimes_{A_0} N \to \wdh{\bigoplus}_I N.
\]
So we have a morphism of sequences:
\[
\begin{tikzcd}
& P \otimes_{A_0} Q \arrow{r}\arrow{d} &P\otimes_{A_0} A_0^n \arrow{r} \arrow{d} & P\otimes_{A_0} M \arrow{r} \arrow{d} & 0\\
0 \arrow{r} & \wdh{\bigoplus}_I Q \arrow{r} & \wdh{\bigoplus}_I A_0^n \arrow{r} & \wdh{\bigoplus}_I M \arrow{r} & 0.
\end{tikzcd}
\]

The map $A_0^n \otimes_{A_0} P \to \wdh{\bigoplus}_I A_0^n$ is an isomorphism because $A_0^n \otimes_{A_0} P=P^n$ is already $p$-adically complete. This implies that the arrow
\[
M\otimes_{A_0} P\to \wdh{\bigoplus}_I M
\]
is surjective. But then 
\[
P\otimes_{A_0} Q \to \wdh{\bigoplus}_I Q
\]
is surjective since $M$ was an arbitrary finitely presented $A$-module. Now a diagram chase implies that
\[
M\otimes_{A_0} P \to \wdh{\bigoplus}_I M
\]
is also injective. And, therefore, it is an isomorphism. So 
\[
P\otimes_{A_0} Q \to \wdh{\bigoplus}_I Q
\]
is also an isomorphism. Therefore, these two
sequences are the same. In particular, 
\[
    P \otimes_{A_0} Q \to P \otimes_{A_0} A_0^n
\]
is injective.
\end{proof}

To establish the second part of our strategy, we will also need a slightly refined version of \cite[Lemma 4.5]{Sch1} specific to the situation of an \'etale morphism to a torus. \smallskip

We recall that we have defined 
\[
R^+\coloneqq \O_K\langle T_1^{\pm 1}, \dots, T_n^{\pm 1}\rangle,
\]
\[
R^+_m\coloneqq \O_K\langle T_1^{\pm 1/p^m}, \dots, T_n^{\pm 1/p^m}\rangle, \text{ and}
\]
\[
R^+_{\infty}=\O_K\langle T_1^{\pm 1/p^{\infty}}, \dots, T_n^{\pm 1/p^{\infty}}\rangle=\left(\text{colim}_n \ R^+_m \right)\widehat{}.
\] 
and a group $\Delta_{\infty}\simeq \bf{Z}_p^m$ continuously acts on $R_\infty^+$. We also define $R$ (resp. $R_m$, $R_{\infty}$) as $R^+[1/p]$ (resp. $R^+_m[1/p]$, $R^+_{\infty}[1/p]$). For an \'etale morphism $\Spa(A, A^+) \to \Spa(R, R^+)=\bf{T}^n$ we define a Huber pair 
\[
\left(A_m, A_m^+\right)\coloneqq \left(R_m \otimes_R A, (R_m\otimes_R A)^+\right) = \left(R_m \wdh{\otimes}_R A, (R_m\wdh{\otimes}_R A)^+\right),
\]
where $(R_m\wdh{\otimes}_R A)^+$ is the {\it integral closure} of the image of $R^+_m\wdh{\otimes}_{R^+} A^+$ in $R_m \wdh{\otimes}_R A$. Similarly, we define 
\[
A^+_\infty \coloneqq \left(\text{colim}_n \ A^+_m \right)\widehat{}
\]
and $A_\infty\coloneqq A^+_\infty[1/p]$.

\begin{lemma}\label{zamena-bazi}\cite[Lemma 4.5]{Sch1} Let $\Spa(A, A^+) \to \Spa(R, R^+)=\bf{T}^n$ be a morphism that is a composition of a finite \'etale maps and rational embeddings. Then $(A_\infty, A_\infty^+)$ is an affinoid perfectoid pair, $\Spd(A_\infty, A_\infty^+) \to \Spd(A, A^+)$ is a $\ud{\Delta}_\infty$-torsor, and, for any $n\in \bf{Z}$, there exists $m$ such that the morphism 
\[
A_{m}^+\wdh{\otimes}_{R_m^+} R_{\infty}^+ \to A_{\infty}^+
\]
is injective with cokernel annihilated by $\varpi^{1/p^n}$.
\end{lemma}
\begin{proof}
We note that \cite[Lemma 4.5]{Sch1} proves that $(A_\infty, A_\infty^+)$ is an affinoid perfectoid space (denoted by $(S_\infty, S_\infty^+)$ there). By construction (and Proposition~\ref{prop:properties-of-diamond}~\ref{prop:properties-of-diamond-6}), $\Spd(A_m, A_m^+) \to \Spd(A, A^+)$ is a $\ud{(\Z/p^m\Z)}^n$-torsor. So $\Spd(A_\infty, A_\infty^+) \simeq \lim_m \Spd(A_m, A_m^+)$ (see Proposition~\ref{prop:properties-of-diamond}~\ref{prop:properties-of-diamond-5}) is a $\ud{\Delta}_\infty\simeq \lim_m \ud{(\Z/p^n\Z)}^n$-torsor. Therefore, we are only left to show that, for any $n\in \bf{Z}$, there exists $m$ such that the morphism 
\[
A_{m}^+\wdh{\otimes}_{R_m^+} R_{\infty}^+ \to A_{\infty}^+
\]
is injective with the cokernel annihilated by $\varpi^{1/p^n}$. \smallskip

We denote by $\widetilde{A}_m$ the $p$-adic completion of $p$-torsionfree quotient of $A_{m}^+\otimes_{R_m^+} R_{\infty}^+$ ($\widetilde{A}_m$ is denoted by $A_m$ in \cite[Lemma 4.5]{Sch1}). Then \cite[Lemma 4.5]{Sch1} shows that, for any $n\in \Z$, there exists $m$ such that the map $\widetilde{A}_m \to A_{\infty}^+$ has the cokernel annihilated by $\varpi^{1/p^n}$. Moreover, the map becomes an isomorphism after inverting $p$. We observe that this implies that $\widetilde{A}_m \to A^+_{\infty}$ is injective as the kernel should be $p^{\infty}$-torsion, but the $p$-adic completion of a $p$-torsionfree ring is $p$-torsionfree. Thus, the only thing we need to show is that $A_{m}^+ \otimes_{R_m^+} R_{\infty}^+$ is already $p$-torsionfree for any $m$. We note that $R_{\infty}^+$ is topologically free as an $R_m^+$-module because

\begin{align*}
R_\infty^+ =\O_K\langle T_1^{\pm 1/p^\infty}, \dots, T_n^{\pm 1/p^\infty}\rangle & = \wdh{\bigoplus}_{(b_1, \dots, b_n) \in \bf{Z}^n \setminus m\bf{Z}^n} \O_K\langle T_1^{\pm 1/p^m}, \dots, T_n^{\pm 1/p^m}\rangle T_1^{1/p^{b_1}}\dots T_n^{1/p^{b_n}}\\
 & =\wdh{\bigoplus}_{(b_1, \dots, b_n) \in \bf{Z}^n \setminus m\bf{Z}^n} R_m^+ \cdot T_1^{1/p^{b_1}}\dots T_n^{1/p^{b_n}}.
\end{align*}
Thus, $R_\infty^+$ is $R_m^+$-flat for any $m$ due to Lemma~\ref{lemma:topologically-free-flat}. Therefore, $A_{m}^+ \otimes_{R_m^+} R_{\infty}^+$ is flat over $A_m^+$, so it is, in particular, $\O_K$-flat. As a consequence, it does not have any non-zero $p$-torsion. This finishes the proof.
\end{proof}

\begin{lemma}\label{second-calculation} Let $\X=\Spf A_0$ be an affine admissible formal $\O_K$-scheme with generic fiber $X=\Spa(A, A^+)$ that admits a map $f\colon X\to \bf{T}^n=\Spa(R, R^+)$ that factors as a composition of finite \'etale morphisms and rational embeddings. Then the cohomology groups
\[
\rm{H}^i(X^\diam_v, \O_{X^\diam}^+/p)
\]
are almost coherent $A_0/p$-modules for $i\geq 0$.
\end{lemma}
\begin{proof}
We denote the completed algebraic closure of $K$ by $C$. Then we note that Lemma~\ref{lemma:affinoid-base-change} implies that
\[
\rm{H}^i(X^\diam_v, \O_{X^\diam}^+/p)\otimes_{\O_K/p} \O_C/p \to \rm{H}^i(X^\diam_{C, v}, \O_{X_C^\diam}^+/p)
\]
is an almost isomorphism for any $i\geq 0$. Therefore, faithful flatness of the morphism $\O_K/p \to \O_C/p$ and Lemma~\ref{almost-flat-descent} imply that it suffices to prove the claim under the additional assumption that $K=C$ is algebraically closed.\smallskip

Theorem~\ref{thm:top-ft-almost-noetherian} ensures that $A_0$ is an almost noetherian ring, thus it suffices to show that $\rm{H}^i(X^\diam_v, \O_{X^\diam}^+/p)$ are almost finitely generated $A_0/p$-modules. \smallskip

Now the generic fiber $X$ is smooth over $C$, so \cite[Corollary 6.4.1/5]{BGR} implies that $A^+=A^{\circ}$ is a flat, topologically finitely type $\O_C$-algebra that is finite over $A_0$. Thus Lemma~\ref{extension} ensures that it suffices to show that $\rm{H}^i(X^\diam_v, \O_{X^\diam}^+/p)$ is almost finitely generated $A^+/pA^+$-modules for $i\geq 0$. We note that $A^+$ is almost noetherian as a topologically finitely generated $\O_C$-algebra, so almost coherent and almost finitely generated $A^+$-modules coincide.\smallskip

We consider a $\ud{\Delta}_\infty$-torsor $\Spd(A_\infty, A_\infty^+) \to \Spd(A, A^+)$ that is constructed in Lemma~\ref{zamena-bazi}. Thus, Corollary~\ref{cor:v-cohomology-group-cohomology} ensures that
\[
\bf{R}\Gamma(X_v^\diam, \O_{X^\diam}^+/p)\simeq^a \bf{R}\Gamma_{\rm{cont}}(\Delta_\infty, A_\infty^+/p).
\]
So we reduce the problem to showing that the complex $\bf{R}\Gamma_{\rm{cont}}(\Delta_{\infty}, A_\infty^+/p)$ has almost finitely generated cohomology modules.\smallskip

We pick any $\e\in \Q_{>0}$ and use Lemma~\ref{zamena-bazi} to find $m$ such that the map
\[
A_m^+\wdh{\otimes}_{R_m^+} R_{\infty}^+ \to A_{\infty}^+
\]
is injective with cokernel killed by $p^{\e}$. Thus, we conclude that the map
\[
A_m^+/p \otimes_{R_m^+/p} R_{\infty}^+/p \to A_{\infty}^+/p
\]
has kernel and cokernel annihilated by $p^{\e}$. Then it is clear that the induced map
\[
\rm{H}^i_{\rm{cont}}(\Delta_{\infty}, A_m^+/p \otimes_{R_m^+/p} R_{\infty}^+/p) \to \rm{H}^i_{\rm{cont}}(\Delta_{\infty}, A_{\infty}^+/p)
\]
has kernel and cokernel annihilated by $p^{2\e}$ for any $i\geq 0$. Therefore, Lemma~\ref{almost-finitely-presented} implies that it is sufficient to show that $\rm{H}^i_{\rm{cont}}(\Delta_{\infty}, A_m^+/p \otimes_{R_m^+/p} R_{\infty}^+/p)$ is almost finitely generated over $A^+/p$ for any $m\geq 0$ and any $i\geq 0$.\smallskip

The trick now is to consider the subgroup $p^m\Delta_{\infty}$ that acts trivially on $A_m^+/p$ to pull it out of the cohomology group by Lemma~\ref{first-computation}. More precisely, we consider the Hochschild--Serre spectral sequence
\[
\rm{E}^{i,j}_2=\rm{H}^i\left(\Delta_{\infty}/p^m\Delta_{\infty}, \rm{H}^j_{\rm{cont}}(p^m\Delta_{\infty}, A_m^+/p \otimes_{R_m^+/p} R_{\infty}^+/p)\right) \Rightarrow \rm{H}^{i+j}_{\rm{cont}}(\Delta_{\infty}, A_m^+/p \otimes_{R_m^+/p} R_{\infty}^+/p)
\]

We recall that the group cohomology of any finite group $G$ can be computed via an explicit bar complex. Namely, for a $G$-module $M$, the complex looks like
\[
\rm{C}^0(G, M) \xr{\rm{d^0}} \rm{C}^1(G, M) \xr{\rm{d}^1} \dots
\]
where 
\[
\rm{C}^i(G, M)=\left\{f:G^i \to M \right\} \simeq M^{\oplus i\cdot \#G}
\]
and 
\begin{eqnarray*}
d^i(f)(g_0, g_1, \dots, g_i)&=&g_0\cdot f(g_1, \dots, g_i)+ \\
&& \sum_{j=1}^i (-1)^j f(g_0, \dots, g_{j-2},g_{j-1}g_j, g_{j+1}, \dots, g_i)+(-1)^{i+1}f(g_0, \dots, g_{i-1}).
\end{eqnarray*}
In case $M$ is an $A^+/p$-module and $G$ acts $A^+/p$-linearly on $M$, all terms $\rm{C}^i(G, M)$ have a natural structure of an $A^+/p$-module, and the differentials are $A^+/p$-linear. Moreover, the terms $\rm{C}^i(G, M)$ are finite direct sums of $M$ as an $A^+/p$-module. In particular, they are almost coherent, if so is $M$. Thus, Lemma~\ref{main-coh} guarantees that all cohomology groups $\rm{H}^i(G, M)$ are almost coherent over $A^+/p$ if $M$ is almost coherent (equivalently, almost finitely generated) over $A^+/p$. \smallskip

We now apply this observation (together with Lemma~\ref{main-coh}) for 
\[
G=\Delta_{\infty}/p^m\Delta_{\infty} \text{ and } M=\rm{H}^j_{\rm{cont}}(p^m\Delta_{\infty}, A_m^+/p \otimes_{R_m^+/p} R_{\infty}^+/p)
\]
to conclude that it suffices to show that $\rm{H}^j_{\rm{cont}}(p^m\Delta_{\infty}, A_m^+/p \otimes_{R_m^+/p} R_{\infty}^+/p)$ is almost coherent (equivalently, almost finitely generated) over $A^+/p$ for any $j\geq 0$, $m\geq 0$. We note that $A_m^+$ is finite over $A^+$ by 
\cite[Corollary 6.4.1/5]{BGR}. Thus, Lemma~\ref{extension} implies that it is enough to show that $\rm{H}^j_{\rm{cont}}(p^m\Delta_{\infty}, A_m^+/p \otimes_{R_m^+/p} R_{\infty}^+/p)$ is almost finitely generated over $A_m^+/p$ for $i\geq 0$ and $m\geq 0$. Now we use Lemma~\ref{first-computation} to write
\[
\rm{H}^j_{\rm{cont}}(p^m\Delta_{\infty}, A_m^+/p \otimes_{R_m^+/p} R_{\infty}^+/p)\simeq \rm{H}^j_{\rm{cont}}(p^m\Delta_{\infty}, R_{\infty}^+/p) \otimes_{R_m^+/p} A_m^+/p
\]
Moreover, Lemma~\ref{first-computation} guarantees that $\rm{H}^j_{\rm{cont}}(p^m\Delta_{\infty}, R_{\infty}^+/p)$ is almost finitely generated over $R_m^+/p$. Thus $\rm{H}^j_{\rm{cont}}(p^m\Delta_{\infty}, R_{\infty}^+/p) \otimes_{R_m^+/p} A_m^+/p$ is almost finitely generated over $A_m^+/p$ by Lemma~\ref{trivial-base-change}.
\end{proof}

\begin{cor}\label{second-calculation-2} Let $\X=\Spf A_0$ and $X=\Spa(A, A^+)$ be as in Lemma~\ref{second-calculation}, and let $\cal{E}$ be a small $\O_{X^\diam}^+/p$-vector bundle. Then the cohomology group $\rm{H}^i(X^\diam_v, \cal{E})$ is almost coherent over $A_0/pA_0$ for any $i\geq 0$.
\end{cor}
\begin{proof}
Similarly to the proof of Lemma~\ref{second-calculation}, we can assume that $K=C$ is algebraically closed and $A_0=A^\circ=A^+$ is almost noetherian. \smallskip

By assumption, we can find a finite \'etale surjection $Y \to X$ that splits $\cal{E}$. Since $X$ is noetherian, we can dominate it by a Galois cover to assume that $Y \to X$ is a $G$-torsor for a finite group $G$ such that $\cal{E}|_{Y^\diam_v}\simeq (\O_{Y^\diam}^+/p)^r$ for some $r$. Then we have the Hochschild--Serre spectral sequence
\[
\rm{E}^{i,j}_2=\rm{H}^i\left(G, \rm{H}^j\left(Y^\diam_v, \left(\O_{Y^\diam}^+/p\right)^r\right)\right) \Rightarrow \rm{H}^{i+j}(X^\diam_v, \cal{E})
\]
Now note that \cite[Corollary 6.4/5]{BGR} implies that $\O_X^+(X) \to \O_Y^+(Y)$ is a finite morphism. Therefore, similarly to the proof of Lemma~\ref{second-calculation}, the argument with the explicit bar complex computing $\rm{H}^i(G, -)$ implies that it is sufficient to show that $\rm{H}^j\left(Y^\diam_v, \left(\O_{Y^\diam}^+/p\right)^r\right)$ is almost coherent over $\O_{Y^\diam}^+(Y^\diam)/p$ for $j\geq 0$. But this is done in Lemma~\ref{second-calculation}.
\end{proof}

\begin{lemma}\label{lemma:reasonable-cover-2} Let $K$ be a $p$-adic perfectoid field, let $\X$ be an admissible formal $\O_K$-scheme with adic generic fiber $X=\X_K$, and let $\cal{E}$ be an $\O_{X^\diam}^+/p$-vector bundle on $X^\diam_v$. Then there is a collection of
\begin{enumerate}[label=\textbf{(\arabic*)}]
    \item an admissible blow-up $\X'\to \X$,
    \item a finite open affine cover $\X'=\bigcup_{i\in I} \sU_i$,
\end{enumerate}
such that, for every $i\in I$, the restriction $\cal{E}|_{(\sU_{i, K})^\diam_v}$ is small.
\end{lemma}
\begin{proof}
    Corollary~\ref{cor:different-vector-bundles-equivalent} ensures that there is a finite open cover $X=\bigcup_{i\in I} U_i$ such that $\cal{E}|_{(U_{i, K})_v^\diam}$ can be trivialized by a finite \'etale surjection. Therefore, \cite[Lemma 8.4/5]{B} implies that there is an admissible blow-up $\X'\to \X$ with a covering $\X'=\bigcup_{i\in I} \sU_{i}$ such that $\sU_{i, K}=U_i$. We can then refine $\sU$  to assume that each $\sU_i=\Spf A_{i, 0}$ is affine. 
\end{proof}

\begin{thm}\label{almost-coh} Let $\X$ be an admissible formal $\O_K$-scheme with smooth adic generic fiber $X$ and mod-$p$ fiber $\X_0$. Then 
\[
\bf{R}\nu_*(\cal{E})^a\in \mathbf{D}^{+}_{acoh}(\X_0)^a
\] 
for any $\O_{X^\diam}^+/p$-vector bundle $\cal{E}$.
\end{thm}
\begin{proof}
First, we note that the claim is clearly Zariski-local on $\X$ and descends through rig-isomorphisms by the Almost Proper Mapping Theorem~\ref{almost-proper-mapping}. Thus Lemma~\ref{lemma:reasonable-cover-2} implies that it suffices to prove the theorem for $\X=\Spf A_0$ an affine formal $\O_K$-scheme and a small $\cal{E}$. \smallskip

Now we note that $\X$ is rig-smooth in the terminology of \cite[\textsection 3]{BLR3}. Thus, \cite[Proposition 3.7]{BLR3} states that there is an admissible blow-up $\pi\colon \X' \to \X$ and a covering of $\X'$ by open affine formal subschemes $\sU'_i$ with rig-\'etale morphisms $\mf'_i\colon \sU'_i \to \wdh{\bf{A}}_{\O_K}^{n_i}$, i.e. the adic generic fibers $\mf'_{i, K}\colon \sU'_{i, K} \to \bf{D}^{n_i}_K$ are \'etale. We apply the Almost Proper Mapping Theorem~\ref{almost-proper-mapping} again to conclude that it suffices to show the theorem for $\X'$. Moreover, since the claim is Zariski-local on $\X$, we can even pass to each $\sU'_i$ separately. So we reduce to the case where $\X=\Spf A_0$ is affine, admits a rig-\'etale morphism $\mf'\colon \colon \X \to \wdh{\bf{A}}^d_{\O_K}$, and $\cal{E}$ is small. \smallskip

We wish to reduce the question to the situation of Corollary~\ref{second-calculation-2}, though we are still not quite there. The key trick now is to use Theorem~\ref{etale-finite-etale-formal} to find a finite rig-\'etale morphism $\mf\colon \X \to \wdh{\bf{A}}^d_{\O_K}$. In particular, the generic fiber $\mf_K\colon X \to \bf{D}^d_K$ is a finite \'etale morphism. So the only thing we are left to do is to embedd $\bf{D}^d_K$ into $\bf{T}^d_K$ as a rational subset. This is done by observing that 
\[
\bf{D}^d_K \simeq \bf{T}^d_K\left(\frac{T_1-1}{p}, \dots, \frac{T_d-1}{p}\right) \subset \bf{T}^d_K.
\]
In particular, $X$ admits an \'etale morphism to a torus that is a composition of a finite \'etale morphism and a rational embedding. Therefore, Corollary~\ref{second-calculation-2} implies that 
\[
\bf{R}\Gamma(X^\diam_v, \cal{E})^a \in \bf{D}^+_{acoh}(A_0/pA_0)^a.
\]
Finally, we note that Lemma~\ref{almost-quasi} ensures that $
\widetilde{\bf{R}\Gamma(X^\diam_v, \cal{E})} \simeq \bf{R}\nu_*\cal{E}$, so
\[
\left(\bf{R}\nu_*\cal{E}\right)^a \in \bf{D}^+_{acoh}(\X_0)^a
\]
by Theorem~\ref{derived-schemes-2}. 
\smallskip
\end{proof}

\subsection{Nearby cycles are almost coherent for general $X$ and $\cal{E}$} The main goal of this section is to generalize Theorem~\ref{almost-coh} to the case of a general generic fiber $X$. The idea is to reduce the general case to the smooth case by means of Lemma~\ref{lemma:different-top}, resolution of singularities, and proper hyperdescent.\smallskip

For the rest of this section, we fix a perfectoid $p$-adic field $K$ with a good pseudo-uniformizer $\varpi \in \O_K$ (see Definition~\ref{defn:good-unifor}). We always do almost mathematics with respect to the ideal $\m= \bigcup_{n}\varpi^{1/p^n}\O_K$. \smallskip

\begin{lemma}\label{conclusion-finite} Let $\Spf A_0$ be an admissible affine formal $\O_K$-scheme with adic generic fiber $\Spa(A, A^+)$. Let $f\colon X \to \Spa(A, A^+)$ be a proper morphism with smooth $X$, and let $\cal{E}$ be an $\O^+_{\Spd(A, A^+)}/p$-vector bundle. Then $\rm{H}^i(X^\diam_v, \cal{E})$ is an almost coherent $A_0/p$-module for any $i\geq 0$.
\end{lemma}
\begin{proof}
First, \cite[Assertion (c) on p.307]{BL1} implies that we can choose an admissible formal $\O_K$-model $\X$ of $X$ with a morphism $\mf\colon \X \to \Spa A_0$ such that $\mf_K=f$. The map $f$ is proper by \cite[Lemma 2.6]{Lutke-proper} (or \cite[Corollary 4.4 and 4.5]{Temkin-proper}). Now we can compute 
\[
\bf{R}\Gamma(X^\diam_v, \cal{E}) \simeq \bf{R}\Gamma(\X_0, \bf{R}\nu_*\left(\cal{E})\right)
\]
Theorem~\ref{almost-coh} implies that $\bf{R}\nu_*\left(\cal{E}\right)\in \bf{D}^+_{acoh}(\X_0)$ as $X$ is smooth. Thus, Theorem~\ref{almost-proper-mapping} implies that 
\[
\bf{R}\Gamma(X^\diam_v, \cal{E}) \simeq \bf{R}\Gamma(\X_0, \bf{R}\nu_*\left(\cal{E})\right) \in \bf{D}^+_{acoh}(A_0/p). \qedhere
\]
\end{proof}

Now we recall the notion of a hypercovering that will be crucial for our proof. We refer to \cite[\href{https://stacks.math.columbia.edu/tag/01FX}{Tag 01FX}]{stacks-project} and \cite{conrad-hyper} for more detail.

\begin{defn} Let $\mathcal{C}$ be a category admitting finite limits. Let $\bf{P}$ be a class of morphisms in $\mathcal C$ which is stable under base change, preserved under
composition (hence under products), and contains all isomorphisms. A simplicial object $X_{\bullet}$ in $\mathcal C$ is said to be a {\it $\bf{P}$-hypercovering} if, for all $n \geq 0$, the natural adjunction map\footnote{See \cite[\textsection 3]{conrad-hyper} (or \cite[\href{https://stacks.math.columbia.edu/tag/0AMA}{Tag 0AMA}]{stacks-project}) for the definition of the coskeleton functor.}
\[
X_{\bullet} \to \text{cosk}_n(\text{sk}_n(X_{\bullet}))
\]
induces a map $X_{n+1} \to (\text{cosk}_n(\text{sk}_n(X_{\bullet})))_{n+1}$ in degree $n+1$ which is in $\bf{P}$. If $X_{\bullet}$ is an augmented simplicial complex, we make a similar definition but also require the case $n = -1$ (and then we say $X_{\bullet}$ is a $\bf{P}$-hypercovering of $X_{-1}$).
\end{defn}

\begin{lemma}\label{lemma:smooth-hypercovering} Let $X$ be a quasi-compact, quasi-separated rigid-analytic variety over $K$. Then there is a proper hypercovering $a\colon X_\bullet \to X$ such that all $X_i$ are smooth over $K$. 
\end{lemma}
\begin{proof}
    First, we note that quasi-compact rigid-analytic varieties over $\Spa(K, \O_K)$ admit resolution of singularities by \cite[Theorem 5.2.2]{Temkin-resolution}. Thus, the proof of \cite[Theorem 4.16]{conrad-hyper} (or \cite[\href{https://stacks.math.columbia.edu/tag/0DAX}{Tag 0DAX}]{stacks-project}) carries over to show that there is a proper hypercovering $a\colon X_{\bullet} \to X$ such that all $X_i$ are smooth over $\Spa(K, \O_K)$.
\end{proof}

\begin{lemma}\label{lemma:proper-v-hypercovering} Let $a\colon X_\bullet \to X$ be a proper hypercovering of a rigid-analytic variety $X$. Then $a^\diam \colon X_\bullet^\diam \to X^\diam$ is a $v$-hypercovering of $X^\diam$. 
\end{lemma}
\begin{proof}
    The functor $(-)^\diam$ commutes with fiber products by Proposition~\ref{prop:properties-of-diamond}~\ref{prop:properties-of-diamond-6}. So 
    \[
    ((\text{cosk}_n(\text{sk}_n X_{\bullet}))_{n+1})^\diam \simeq (\text{cosk}_n(\text{sk}_n X_{\bullet}^\diam))_{n+1}.
    \]
    Therefore, the only thing we need to show is that $(-)^\diam$ sends proper coverings to $v$-coverings. This follows from Lemma~\ref{lemma:naive-normal-v-covering} and Example~\ref{exmpl:naive-v-covering}. 
\end{proof}

\begin{thm}\label{acoh-nosmooth} Let $\X$ be an admissible formal $\O_K$-scheme with adic generic fiber $X$ and mod-$p$ fiber $\X_0\coloneqq \X\times_{\Spf\O_K}\Spec \O_K/p$. 
Then 
\[
\bf{R}\nu_*\cal{E} \in \mathbf{D}^{+}_{acoh}(\X_0)
\] 
for any $\O_{X^\diam}^+/p$-vector bundle $\cal{E}$.
\end{thm}
\begin{proof}
    The claim is Zariski-local on $\X$, so we can assume that $\X=\Spf A_0$ is affine. Thus, Lemma~\ref{almost-quasi} and Theorem~\ref{derived-schemes-2} ensure that it suffices to show that 
    \[
    \bf{R}\Gamma(X_v^\diam, \cal{E})\in \bf{D}^{+}_{acoh}(A_0/p).
    \]
    Lemma~\ref{lemma:smooth-hypercovering} shows that there is a proper hypercovering $a\colon X_\bullet \to X$ with smooth $X_i$, and Lemma~\ref{lemma:proper-v-hypercovering} implies that $a\colon X^\diam_\bullet\to X^\diam$ is then a $v$-hypercovering. \smallskip
    
    The proof of \cite[\href{https://stacks.math.columbia.edu/tag/01GY}{Tag 01GY}]{stacks-project} implies that there is a spectral sequence
    \[
    \rm{E}_1^{i,j}=\rm{H}^j\left(X^\diam_{i, v}, \cal{E}\right) \Rightarrow \rm{H}^{i+j}(X^\diam_v, \cal{E}).
    \]
    Lemma~\ref{conclusion-finite} guarantees that $\rm{H}^j(X^\diam_{i, v}, \cal{E})$ is almost coherent over $A_0/p$ for every $i, j \geq 0$. Therefore, Lemma~\ref{main-coh} guarantees that $\rm{H}^{i+j}(X^\diam_v, \cal{E})$ is almost coherent $A_0/p$ for every $i+j\geq 0$.  
\end{proof}

\subsection{Cohomological bound on nearby cycles}\label{section:coh-bound} The main goal of this section is to show that $\bf{R}\nu_*\cal{E}$ is almost concentrated in degrees $[0, d]$ for a small vector bundle $\cal{E}$. This claim turns out to be pretty hard. To achieve this result, we have to use a recent notion of perfectoidization developed in \cite{BS3} that gives a stronger version of the almost purity theorem in the world of diamonds. Our approach is strongly motivated by the proof of \cite[Proposition 7.5.2]{Guo}. \smallskip

For the rest of this section, we fix a perfectoid $p$-adic field $K$ with a good pseudo-uniformizer $\varpi \in \O_K$. We always do almost mathematics with respect to the ideal $\m= \bigcup_{n}\varpi^{1/p^n}\O_K$. \smallskip 

In this section, it is crucial that we work on the level of diamonds. The main observation is that the functor 
\[
(-)^\diamondsuit \colon \{\text{(Pre-)Adic Analytic Spaces}\} \to \{\text{Diamonds}\}
\]
is not fully faithful, so it is possible that a non-perfectoid (pre)-adic space becomes representable by an affinoid perfectoid space after diamondification (we already saw this phenomenon in Warning~\ref{warning:perfectoid-only-after-diamondification}). An explicit construction of such examples is the crux of our argument in this section. To construct such spaces, we need the following theorem of B.\,Bhatt and P.\,Scholze:

\begin{thm}\label{perfectoidization}\cite[Theorem 10.11]{BS3} Let $R$ be an integral perfectoid ring\footnote{We use \cite[Definition 3.5]{BMS1} as the definition for integral perfectoid rings here. This definition coincides with Definition~\ref{defn:integral-perfectoid} in the $p$-torsionfree case, but it is less restrictive in general.}. Let $R \to S$ be the $p$-adic completion of an integral map. Then there exists an integral perfectoid ring $S_{\perfd}$ together with a map of $R$-algebras $S \to S_{\perfd}$, such that it is
initial among all of the $R$-algebra maps $S \to S'$ for $S'$ being integral perfectoid.
\end{thm}

Now we show how this result can be used to obtain a cohomological bound on $\bf{R}\nu_*\left(\cal{E}\right)$. We recall that a torus 
\[
\bf{T}^d=\Spa\left(K\langle T_1^{\pm 1}, \dots, T_d^{\pm 1}\rangle, \O_K\langle T_1^{\pm 1}, \dots, T_d^{\pm 1}\rangle\right) = \Spa(R, R^+)
\]
admits a map 
\[
\bf{T}^d_\infty = \Spa\left(K\langle T_1^{\pm 1/p^{\infty}}, \dots, T_d^{\pm 1/p^{\infty}}\rangle, \O_K\langle T_1^{\pm 1/p^{\infty}}, \dots, T_d^{\pm 1/p^{\infty}}\rangle\right) \to \bf{T}^d
\]
such that $\bf{T}^d_\infty$ is an affinoid perfectoid space, and the map becomes a $\ud{\Delta}_\infty=\ud{\bf{Z}_p(1)}^d$-torsor after applying the diamondification functor.\smallskip

Now we can embed a $d$-dimensional disk $\bf{D}^d$ as a rational subdomain 
\[
\bf{D}^d=\bf{T}^d\left(\frac{T_1-1}{p}, \dots, \frac{T_n-1}{p}\right) \subset \bf{T}^d, 
\]
so the fiber product
\[
\bf{D}^d_\infty =  \bf{D}^d\times_{\bf{T}^d} \bf{T}^d_\infty \to \bf{D}^d
\]
is again an affinoid perfectoid covering of $\bf{D}^d$ by Lemma~\ref{zamena-bazi}. \smallskip

If $X=\Spa(A, A^+) \to \bf{D}^d$ is an arbitrary finite morphism, then the fiber product $X\times_{\bf{D}^d} \bf{D}^d_\infty$ may not be an affinoid perfectoid space (or even an adic space). However, it turns out that the associated diamond is always representable by an affinoid perfectoid space.

\begin{lemma}\label{123} Let $f\colon X=\Spa(A, A^+) \to \bf{D}^d$ be a finite morphism of rigid-analytic $K$-spaces. Then the fiber product $X^\diam_{\infty} \coloneqq X^{\diamondsuit}\times_{\bf{D}^{d, \diam}} \bf{D}^{d, \diam}_{\infty}$ is representable to an affinoid perfectoid space (of characteristic $p$). 
\end{lemma}
\begin{proof}
Let us say that $\bf{D}^d=\Spa(S, S^+)$ and $\wdh{\bf{D}}^d_{\infty}=\Spa(S_{\infty}, S^+_{\infty})$. The map $f$ defines an integral morphism $S^+ \to A^+$, we define 
\[
A_\infty^\dagger\coloneqq S_\infty^+\wdh{\otimes}_{S^+} A^+.
\]
This is a $p$-adic completion of an integral morphism over an integral perfectoid ring $S_\infty^+$ (see \cite[Lemma 3.20]{BMS1}), so there is a map
\[
A_\infty^\dagger \to (A_\infty^\dagger)_{\rm{perfd}}
\]
initial to an integral perfectoid ring. We define $A_\infty$ as $A_\infty^\dagger[1/p]$ and $A_\infty^+$ as the integral closure of $A_\infty^\dagger$ in $A_\infty$. Then $(A_\infty, A_\infty^+)$ is an affinoid perfectoid pair by \cite[Lemma 3.21]{BMS1}. Therefore, it suffices to show that the natural morphism
\[
\Spd(A_\infty, A_\infty^+) \to \Spd(A, A^+)\times_{\Spd(S, S^+)} \Spd(S_\infty, S_\infty^+)
\]
is an isomorphism. This can be easily checked on the level of rational point by the universal property of $(A_\infty^\dagger)_{\rm{perfd}}$ and the construction of the diamondification functor in Definition~\ref{defn:diamondification} (and \cite[Lemma 3.20]{BMS1} that relates affinoid perfectoid pairs and integral affinoid rings).
\end{proof}

\begin{thm}\label{coh-dimension} Let $\X=\Spf A_0$ be an admissible formal $\O_K$-scheme with adic generic fiber $X=\Spa(A, A^+)$ of dimension $d$, and let $\cal{E}$ be a small $\O_{X^\diam}^+/p$-vector bundle. Then 
\[
\bf{R}\Gamma(X^\diam_v, \cal{E})^a\in \bf{D}^{[0, d]}_{acoh}(A_0/p)^a.
\]
\end{thm}
\begin{proof}
    Lemma~\ref{conclusion-finite} ensures that $\bf{R}\Gamma(X_v^\diam, \cal{E})\in \bf{D}_{acoh}(A_0/p)$, so it suffices to show that 
    \[
    \rm{H}^i(X_v^\diam, \cal{E}) \simeq^a 0
    \]
    for $i>d$. The Noether Normalization Theorem (see \cite[Proposition 3.1.3]{B}) implies that there is a finite morphism $f\colon X\to \bf{D}^d$. We consider the $\ud{\Delta}_\infty \simeq \ud{\Z_p(1)}^d$-torsor
    \[
    X_\infty^\diam \simeq X^\diam\times_{\bf{D}^{d, \diam}} \bf{D}^{d,\diam}_\infty\to X^\diam.
    \]
    By Lemma~\ref{123}, $X_\infty^\diam$ is represented by an affinoid perfectoid space $\Spd(A_\infty, A_\infty^+) = \Spa(A_\infty^\flat, A_\infty^{\flat, +})$. Thus, we are in the situation of Set-up~\ref{set-up:torsors}. So Corollary~\ref{cor:v-cohomology-group-cohomology} implies that
    \[
    \rm{H}^i(X_v^\diam, \cal{E})\simeq^a \rm{H}^i_{\rm{cont}}(\Delta_\infty, (M^a_{\cal{E}})_!),
    \]
    where $M_{\cal{E}}\simeq \rm{H}^0(X^\diam_{\infty, v}, \cal{E})$. Therefore, the claim follows from the observation that the cohomological dimension of $\Delta_\infty\simeq \bf{Z}_p(1)^d\simeq \bf{Z}_p^d$ is $d$ due to \cite[Lemma 7.3]{BMS1}.
\end{proof} 

\subsection{Proof of Theorem~\ref{thm:main-thm-small}}

The main goal of this section is to give a full proof of Theorem~\ref{thm:main-thm-small}. Most of the hard work was already done in the previous sections. \smallskip

For the rest of this section, we fix a perfectoid $p$-adic field $K$ with a pseudo-uniformizer $\varpi \in \O_K$ as in Remark~\ref{rmk:roots-of-pseudounformizer}. We always do almost mathematics with respect to the ideal $\m= \bigcup_{n}\varpi^{1/p^n}\O_K$. \smallskip 

\begin{thm}\label{thm:main-thm-small-2} Let $\X$ be an admissible formal $\O_K$-scheme with adic generic fiber $X$ of dimension $d$ and mod-$p$ fiber $\X_0$, and let $\cal{E}$ be an $\O^+_{X^{\diam}}/p$-vector bundle. Then 
\begin{enumerate}[label=\textbf{(\arabic*)}]
    \item\label{thm:main-thm-small-2-1} $\bf{R}\nu_*\cal{E}\in \bf{D}^+_{qc, acoh}(\X_0)$ and $(\bf{R}\nu_*\cal{E})^a\in \bf{D}^{[0, 2d]}_{acoh}(\X_0)^a$;
    \item\label{thm:main-thm-small-2-2} if $\X=\Spf A$ is affine, then the natural map 
    \[
    \widetilde{\rm{H}^i\left(X^\diam_v, \cal{E} \right)} \to \rm{R}^i\nu_*\left(\cal{E}\right)
    \]
    is an isomorphism for every $i\geq 0$;
    \item\label{thm:main-thm-small-2-3} the formation of $\rm{R}^i\nu_*(\cal{E})$ commutes with \'etale base change, i.e., for any \'etale morphism $\mf \colon \Y \to \X$ with adic generic fiber $f\colon Y\to X$, the natural morphism 
    \[
    \mf^*_0 \left(\rm{R}^i\nu_{\X, *}(\cal{E}) \right)\to \rm{R}^i\nu_{\Y, *}\left(\cal{E}|_{Y^\diam}\right)
    \]
    is an isomorphism for any $i\geq 0$;
    \item\label{thm:main-thm-small-2-4} if $\X$ has an open affine covering $\X=\bigcup_{i\in I} \sU_i$ such that $\cal{E}|_{(\sU_{i, K})^\diam}$ is small, then
    \[
    \left(\bf{R}\nu_{*}\cal{E}\right)^a \in \bf{D}^{[0, d]}_{acoh}(\X_0)^a;
    \]
    \item\label{thm:main-thm-small-2-5} there is an admissible blow-up $\X'\to \X$ such that $\X'$ has an open affine covering $\X'=\bigcup_{i\in I} \sU_i$ such that $\cal{E}|_{(\sU_{i, K})^\diam}$ is small. 
    
    In particular, there is a cofinal family of admissible formal models $\{\X'_i
    \}_{i\in I}$ of $X$ such that 
    \[
    \left(\bf{R}\nu_{\X'_i, *}\cal{E}\right)^a\in \bf{D}^{[0, d]}_{acoh}(\X'_{i, 0})^a.
    \]
    for each $i\in I$. 
\end{enumerate}
\end{thm}
\begin{proof}
    The first part of \ref{thm:main-thm-small-2-1},~\ref{thm:main-thm-small-2-2}, and~\ref{thm:main-thm-small-2-3} follow from Theorem~\ref{almost-quasi} and Theorem~\ref{acoh-nosmooth}.  Now to show that $\bf{R}\nu_*\cal{E}$ is almost concentrated in degrees $[0, 2d]$, it suffices to show that, for every affine $\sU=\Spf A_0\subset \X$, the complex $\bf{R}\Gamma(\sU^\diam_{K, v}, \cal{E})^a$ (almost) lies
    in $\bf{D}^{[0, 2d]}(A_0/p)^a$. By Lemma~\ref{lemma:affinoid-base-change} and full faithful flatness of $\O_K/p \to \O_C/p$, it is sufficient to proof under the additional assumption that $K=C$ is algebraically closed. Then Theorem~\ref{thm:v-quasi-proetale} and Theorem~\ref{thm:quasi-proetale-etale} imply that 
    \[
    \cal{E}' \coloneqq \bf{R}\mu_*\bf{R}\lambda_*\cal{E}
    \]
    is an $\O_{X_\et}^+/p$-vector bundle concentrated in degree $0$. Therefore, 
    \[
    \bf{R}\Gamma(\sU^\diam_{C, v}, \cal{E}) \simeq \bf{R}\Gamma(\sU_{C, \et}, \cal{E}'),
    \]
    and 
    \[
    \bf{R}\Gamma(\sU_{C, \et}, \cal{E}') \in \bf{D}^{[0, 2d]}(A_0/p)
    \]
    due to \cite[Corollary 2.8.3 and Corollary 1.8.8]{H3}. \smallskip
    
    To show \ref{thm:main-thm-small-2-4}, we consider an open affine covering $\X=\bigcup_{i\in I} \sU_i$ and denote $\sU_i=\Spf A_{i}$. Then Part~\ref{thm:main-thm-small-2-2} implies that it suffices to show that 
    \[
    \bf{R}\Gamma((\sU_{i, K})^\diam_v, \cal{E})^a\in \bf{D}^{[0, d]}_{acoh}(A_i/p)^a
    \]
    for each $i\in I$. This follows from Theorem~\ref{coh-dimension} and the assumption that $\cal{E}|_{(\sU_{i, K})^\diam}$ is small. \smallskip
    
    \ref{thm:main-thm-small-2-5} now follows from Lemma~\ref{lemma:reasonable-cover-2}.
\end{proof}

\subsection{Proof of Theorem~\ref{main-thm}}

The main goal of this section is to prove Theorem~\ref{main-thm}. The idea is to reduce to the case of the constant Zariski-constructible sheaf through a sequence of reductions; in this case, the result follows directly from Theorem~\ref{thm:main-thm-small}.\smallskip

For the rest of this section, we fix a perfectoid $p$-adic field $K$ with a pseudo-uniformizer $\varpi \in \O_K$ as in Remark~\ref{rmk:roots-of-pseudounformizer}. We always do almost mathematics with respect to the ideal $\m= \bigcup_{n}\varpi^{1/p^n}\O_K$. \smallskip 

We consider the following diagram of morphisms of ringed sites: 

\begin{equation*}
\begin{tikzcd}
\left(X^\diam_v, \O_{X^\diam}^+/p\right) \arrow[rrr, bend right, "\nu"] \arrow{r}{\lambda} & \left(X^\diam_\qproet, \O_{X^\diam_\qp}^+/p\right) \arrow{r}{\mu} & \left(X_\et, \O_{X_\et}^+/p\right) \arrow{r}{t} & \left(\X_{\rm{Zar}}, \O_{\X_0}\right).
\end{tikzcd}
\end{equation*}

Both $\nu_*$ and $t_*$ will play an important role in the proof. 

\begin{lemma}\label{lemma:nearby-pushforward-finite} Let $\mf\colon \X \to \Y$ a finite morphism of admissible formal $\O_K$-schemes with adic generic fiber $f\colon X \to Y$, and $\F\in \bf{D}^b_{zc}(X; \bf{F}_p)$. Then the natural morphism
\[
\bf{R}\nu_{\Y, *}\left(f_{*} \F\otimes \O^+_{Y^\diam}/p\right) \to \bf{R}\mf_{0, *}\left(\bf{R}\nu_{\X, *}\left(\F\otimes \O_{X^\diam}^+/p\right)\right)
\]
is an isomorphism in $\bf{D}(\Y_0)$. 
\end{lemma}
\begin{proof}
    First, we note that $f$ is finite, and so $f_*\simeq \bf{R}f_{*}$ due to \cite[Proposition 2.6.3]{H3}. Now the proof of Corollary~\ref{lemma:nearby-pushforward} just goes through using Corollary~\ref{cor:primitive-finite-rigid} (that does not use Theorem~\ref{main-thm} as in input) in place of Lemma~\ref{lemma:relative-primitive}. 
\end{proof}

\begin{lemma}\label{lemma:descend-through-finite-morphisms} Let $f\colon X \to Y$ be a finite morphism of quasi-compact, quasi-separated rigid-analytic varieties over $K$, and $\F\in \bf{D}^{[r,s]}_{zc}(X; \bf{F}_p)$ such that 
\[
\bf{R}\nu_{\X, *}\left(\F\otimes \O_{X^\diam}^+/p\right)^a \in \bf{D}_{acoh}^{[r, s+d]}(\X_0)^a \text{ (resp. } \bf{R}\nu_{\X, *}\left(\F\otimes \O_{X^\diam}^+/p\right) \in \bf{D}^+_{qc, acoh}(\X_0))
\]
for any formal $\O_K$-model $\X$ of $X$. Then, for any formal $\O_K$-model $\Y$ of $Y$,
\[
\bf{R}\nu_{\Y, *}\left(f_*\F \otimes \O_{Y^\diam}^+/p\right)^a \in \bf{D}^{[r, s+d]}_{acoh}(\Y_0)^a \text{ (resp. } \bf{R}\nu_{\Y, *}\left(f_*\F \otimes \O_{Y^\diam}^+/p\right) \in \bf{D}^+_{qc, acoh}(\Y_0)). 
\]
\end{lemma}
\begin{proof}
    First, we note that we can choose a finite morphism $\mf\colon \X \to \Y$ such that its generic fiber $\mf_K$ is equal to $f$ (for example, this follows from \cite[Corollary II.5.3.3, II.5.3.4]{FujKato}).  \smallskip
    
    Now Lemma~\ref{lemma:nearby-pushforward-finite} ensures that the natural morhism
    \[
    \bf{R}\nu_{\Y, *}\left(f_{*} \F\otimes \O^+_{Y^\diam}/p\right) \to \bf{R}\mf_{0, *}\left(\bf{R}\nu_{\X, *}\left(\F\otimes \O_{X^\diam}^+/p\right)\right)
    \]
    is an isomorphism. Therefore, $\bf{R}\nu_{\Y, *}\left(f_* \F\otimes \O^+_{Y^\diam}/p\right)$ already lies in $\bf{D}_{acoh}(\Y_0)^a$ (resp. $\bf{D}_{qc, acoh}(\Y_0)$) by Theorem~\ref{almost-proper-mapping}. The cohomological bound follows from Proposition~\ref{derived-pushforward} and the fact that the finite morphism $\mf_0$ is (almost) exact on (almost) quasi-coherent sheaves.
\end{proof}

\begin{lemma}\label{lemma:bound-zc} Let $\X$ be an admissible formal $\O_K$-scheme with adic generic fiber $X$ of dimension $d$ and mod-$p$ fiber $\X_0$, and let $\F \in \bf{D}^{[r, s]}_{zc}(X; \bf{F}_p)$. Then 
\[
\bf{R}t_*\left(\F\otimes\O_{X_\et}^+/p\right) \simeq \bf{R}\nu_*\left(\F\otimes\O_{X^\diam}^+/p\right) \in \bf{D}^+_{qc, acoh}(\X_0), \text{ and}
\] 
\[
\bf{R}\nu_*\left(\F\otimes\O_{X^\diam}^+/p\right)^a \in \bf{D}^{[r, s+d]}_{qc, acoh}(\X_0)^a
\]
\end{lemma}
\begin{proof}
Lemma~\ref{lemma:et-qp-v-overconv-coeff} and Remark~\ref{rmk:zc-overconvergent} imply that \[
\bf{R}t_*\left(\F\otimes\O_{X_\et}^+/p\right) \simeq \bf{R}\nu_*\left(\F\otimes\O_{X^\diam}^+/p\right).
\]
In what follows, we will freely identify these sheaves. Also, we can assume that $\F$ is concentrated in degree $0$, i.e., $\F$ is a usual Zariski-constructible sheaf. \smallskip

{\it Step~$1$: The case of a local system $\F$.} In this case, $\cal{E}\coloneqq \F\otimes \O_{X^\diam}^+/p$ fits the assumption of Theorem~\ref{thm:main-thm-small-2}. Since an $\bf{F}_p$-local system on any rigid-analytic variety $Y$ splits by a finite \'etale cover, so $\F\otimes \O_{X^\diam}^+/p$ is small for any open affinoid $U\subset X$. Thus, the desired claim follows from Theorem~\ref{thm:main-thm-small-2}. \smallskip

{\it Step~$2$: Case of a zero-dimensional $X$.} If $X$ is of dimension $0$, then any Zariski-constructible sheaf on $X$ is a local system. So the claim follows from Step~$1$. \smallskip

Now we argue by induction on $\dim X$. We suppose the claim is known for every rigid-analytic variety of dimension less than $d$ (and any Zariski constructible $\F$) and wish to prove the claim for $X$ of dimensiond $d$. \smallskip

{\it Step~$3$: Reduction to the case of a reduced $X$.} Consider the reduction morphism $i\colon X_{\rm{red}} \to X$. Then $i_\et$ is an equivalence of \'etale topoi, we see that 
\[
i_{*} i^{-1} \F \to \F
\]
is an isomorphism. Thus the claim follows from Lemma~\ref{lemma:descend-through-finite-morphisms}. \smallskip

{\it Step~$4$: Reduction to the case of a normal $X$.} Consider the normalization morphism $f\colon X' \to X$. It is finite by \cite[Theorem 2.1.2]{C} and an isomorphism outside of a nowhere dense Zariski-closed subset $Z$. We use \cite[Proposition 2.6.3]{H3} and argue on stalks to conclude that the natural morphism $\F \to f_*f^{-1}\F$ is injective. Therefore, there is an exact sequence
\[
0\to \F \to f_*f^{-1}\F \to i_* \G \to 0
\]
where $i\colon Z\to X$ is a Zariski-closed immersion with $\dim Z<\dim X$ and $\G$ is a Zariski-constructible sheaf on $Z$. Now the induction hypothesis and Lemma~\ref{lemma:descend-through-finite-morphisms} ensure that
\[
\bf{R}\nu_*\left(i_*\G \otimes \O_{X^\diam}^+/p\right)\in \bf{D}^{+}_{qc, acoh}\left(\X_0\right),
\]
\[
\bf{R}\nu_*\left(i_*\G \otimes \O_{X^\diam}^+/p\right)^a \in \bf{D}^{[0, d-1]}_{acoh}\left(\X_0\right)^a.
\]
Therefore, it suffices to show the claim for $f_*f^{-1} \F$. Thus,Lemma~\ref{lemma:descend-through-finite-morphisms} guarantees that it suffices to show that
\[
\bf{R}\nu_{\X',*}\left(f^{-1} \F\otimes \O_{{X'}^\diam}^+/p\right)\in \bf{D}^{+}_{qc, acoh}\left(\X'_0\right),
\]
\[
\bf{R}\nu_{\X',*}\left(f^{-1} \F\otimes \O_{{X'}^\diam}^+/p\right)^a \in \bf{D}^{[r, s+d]}_{acoh}\left(\X'_0\right)^a
\]
for any admissible formal $\O_K$-model $\X'$ of $X'$. So we may and do assume that $X$ is normal. \smallskip

{\it Step~$5$: Reduction to the case $\F=\ud{\bf{F}}_p$.} By definition of a Zariski-constructible sheaf, there is a nowhere dense Zariski-closed subset $i\colon Z\to X$ with the open complement $j\colon U\to X$ and an $\bf{F}_p$-local system $\bf{L}$ on $U$ such that $\F|_{U} \simeq \bf{L}$. In particular, there is a short exact sequence
\[
0 \to j_!\bf{L} \to \F\to i_*\F|_{Z} \to 0.
\]
Similarly to the argument in Step~$4$, it suffices to prove the claim for $\F=j_!\bf{L}$. \smallskip

Then ``m\'ethode de la trace'' (see \cite[\href{https://stacks.math.columbia.edu/tag/03SH}{Tag 03SH}]{stacks-project}) implies that there is a finite \'etale covering $g\colon U' \to U$ of degree prime-to-$p$ such that $\bf{L}'\coloneqq \bf{L}|_{U'}$ is an iterated extension of constant sheaves $\ud{\bf{F}}_p$. Then $\bf{L}$ is a direct summand of $g_*\left(\bf{L}'\right)$. Thus, it is enough to prove the claim for 
\[
\F = j_!\left(g_*\bf{L}'\right).
\]
Moreover, it suffices to prove the claim for $\F=j_! \left(g_* \ud{\bf{F}}_p\right)$ because the claim of Lemma~\ref{lemma:bound-zc} satisfies the $(2)$-out-of-$(3)$ property, and both functors $g_*$ and $j_!$ are exact. \smallskip

Now we use \cite[Theorem 1.6]{Hansen-vanishing} to extend $g$ to a finite morphism $g'\colon X' \to X$. Then a similar reduction shows that it suffices to prove the claim for $\F=g'_*\left(\ud{\bf{F}}_p\right)$. This case follows from Step~$1$ and Lemma~\ref{lemma:descend-through-finite-morphisms}.
\end{proof}

\begin{thm}\label{main-thm-2} Let $\X$ be an admissible formal $\O_K$-scheme with adic generic fiber $X$ of dimension $d$ and mod-$p$ fiber $\X_0$, and $\F \in \bf{D}^{[r, s]}_{zc}(X; \bf{F}_p)$. 
Then 
\begin{enumerate}[label=\textbf{(\arabic*)}]
    \item\label{thm:main-thm-2-1} there is an isomorphism $\bf{R}t_*\left(\F\otimes \O_{X_\et}^+/p\right) \simeq \bf{R}\nu_*\left(\F\otimes \O_{X^\diam}^+/p\right)$;
    \item\label{thm:main-thm-2-2} $\bf{R}\nu_*\left(\F\otimes \O_{X^\diam}^+/p\right) \in \bf{D}^+_{qc, acoh}(\X_0)$, and $\bf{R}\nu_*\left(\F\otimes \O_{X^\diam}^+/p\right)^a \in \bf{D}^{[r, s+d]}_{acoh}(\X_0)^a$;
    \item\label{thm:main-thm-2-3} if $\X=\Spf A$ is affine, then the natural map 
    \[
    \widetilde{\rm{H}^i\left(X^\diam_v, \F\otimes \O_{X^\diam}^+/p \right)} \to \rm{R}^i\nu_*\left(\F\otimes \O_{X^\diam}^+/p\right)
    \]
    is an isomorphism for every $i\geq 0$;
    \item\label{thm:main-thm-2-4} the formation of $\rm{R}^i\nu_*\left(\F\otimes \O_{X^\diam}^+/p\right)$ commutes with \'etale base change, i.e., for any \'etale morphism $\mf \colon \Y \to \X$ with adic generic fiber $f\colon Y\to X$, the natural morphism 
    \[
    \mf^*_0 \left(\rm{R}^i\nu_{\X, *}\left(\F\otimes \O_{X^\diam}^+/p\right) \right)\to \rm{R}^i\nu_{\Y, *}\left(f^{-1}\F\otimes \O_{Y^\diam}^+/p\right)
    \]
    is an isomorphism for any $i\geq 0$.
\end{enumerate}
\end{thm}
\begin{proof}
    \ref{thm:main-thm-2-1} and \ref{thm:main-thm-2-2} follow from Lemma~\ref{lemma:bound-zc}. Now \ref{thm:main-thm-2-3} follows from Lemma~\ref{derived-schemes} and the isomorphism
    \[
    \bf{R}\Gamma\left(\X_0, \bf{R}\nu_*\left(\F\otimes \O_{X^\diam}^+/p\right)\right) \simeq \bf{R}\Gamma(X_v^\diam, \F\otimes \O_{X^\diam}^+/p).
    \]
    
    We are left to show \ref{thm:main-thm-2-4}. By \ref{thm:main-thm-2-1}, it suffices to show that the natural morphism
    \[
    \mf_0^*\left(\rm{R}^it_{\X, *}\left(\F\otimes \O_{X_\et}^+/p\right)\right) \to \rm{R}^it_{\Y, *}\left(f^{-1}\F\otimes \O_{Y_\et}^+/p\right)
    \]
    Moreover, \cite[Proposition 3.6]{Bhatt-Hansen} ensures that it suffices to prove the claim for $\F=g_*\left(\ud{\bf{F}}_p\right)$ for some finite morphism $g\colon X'\to X$. Then we can lift it to a finite morphism $\mathfrak{g}\colon \X' \to \X$ as in the proof of Lemma~\ref{lemma:descend-through-finite-morphisms}. Then we have a commutative diagram
    \begin{equation}\label{diagram:big}
    \begin{tikzcd}
    & (Y'_\et, \O_{Y'_\et}^+/p) \arrow{ld}{f'} \arrow{dd}{g'}\arrow{rrr}{t_{\Y'}} & & & (\Y'_0, \O_{\Y'_0})\arrow{ld}{\mf'_0} \arrow{dd}{\mathfrak{g}'_0}\\
    (X'_\et, \O_{X'_\et}^+/p)\arrow{dd}{g} \arrow{rrr}{t_{\X'}} & & & (\X'_0, \O_{\X'_0})\arrow{dd}{\mathfrak{g}_0} & \\
    & (Y_\et, \O_{Y_\et}^+/p)\arrow{ld}{f}\arrow{rrr}{t_{\Y}} & & & (\Y_0, \O_{\Y_0})\arrow{ld}{\mf_0} \\
    (X, \O_{X_\et}^+/p) \arrow{rrr}{t_{\X}}& & & (\X_0, \O_{\X_0})
    \end{tikzcd}
    \end{equation}
    with $\Y'=\Y\times_{\X}\X'$ and $Y'$ being its adic generic fiber. Then we have a sequence of isomorphisms:
    \begin{align*}
        \mf^*_0 \left(\bf{R}t_{\X, *}\left(g_*\left(\bf{\ud{F}}_p\right) \otimes\O_{X_\et}^+/p\right) \right) & \simeq \mf_0^*\left(\bf{R}t_{\X,  *} \left( \bf{R}g_*\O_{X'_\et}^+/p\right)\right)\\
        & \simeq \mf^*_0\left( \bf{R}\mathfrak{g}_{0, *}\left(\bf{R}t_{\X', *}\O_{X'_\et}^+/p\right) \right) \\
        & \simeq \bf{R}\mathfrak{g}'_{0, *}\left( {\mf'_0}^*\left(\bf{R}t_{\X',*}\O_{X'_\et}^+/p\right) \right)  \\
        & \simeq \bf{R}\mathfrak{g}'_{0, *}\left(\bf{R}t_{\Y', *}\left(\O_{Y'_\et}^+/p\right)\right)\\
        & \simeq \bf{R}t_{\Y, *}\left(\bf{R}g'_*\O_{Y'_\et}^+/p\right) \\
        & \simeq \bf{R}t_{\Y, *}\left(g'_*\left(\ud{\bf{F}}_p\right) \otimes \O_{Y_\et}^+/p\right) \\
        & \simeq \bf{R}t_{\Y, *}\left(f^{-1}\left(g_*\ud{\bf{F}}_p\right) \otimes\O_{Y_\et}^+/p\right)
    \end{align*}
    The first isomorphism holds by (the proof of) Corollary~\ref{cor:primitive-finite-rigid}. The second isomorphism is formal and follows from Diagram~\ref{diagram:big}. The third isomorphism holds by flat base change applied $\mf_0$. The fourth isomorphism follows from Theorem~\ref{thm:main-thm-small-2} applied to $\cal{E}=\O_{{X'}^\diam}^+/p$ and the \'etale morphism $\Y' \to \X'$. The fifth isomorphism is formal again. The sixth isomorphism follows from (the proof of) Corollary~\ref{cor:primitive-finite-rigid}. Finally, the last isomorphism follows from \cite[Theorem 4.3.1]{H3}. 
\end{proof}

\subsection{Proof of Theorem~\ref{thm:main-thm-integral}} The main goal of this section is to prove Theorem~\ref{thm:main-thm-integral}. The proof is a formal reduction to the case of $\O_{X^\diam}^+/p$-vector bundles. After that, we also discuss a version of this theorem for the classical pro-\'etale site from \cite{Sch1}. \smallskip

For the rest of this section, we fix a perfectoid $p$-adic field $K$ with a good pseudo-uniformizer $\varpi \in \O_K$. We always do almost mathematics with respect to the ideal $\m= \bigcup_{n}\varpi^{1/p^n}\O_K$. \smallskip 

\begin{lemma}\label{lemma:vector-bundle-complete} Let $X$ be a rigid-analytic variety over $K$, and let $\cal{E}$ be an $\O_{X^\diam}^+$-vector bundle on $X$. Then $\cal{E}$ is derived $p$-adically complete. 
\end{lemma}
\begin{proof}
    Lemma~\ref{thm:main-thm-integral} implies that it suffices to prove the claim $v$-locally on $X^\diam_v$. Therefore, we may and do assume that $\cal{E}=\left(\O_{X^\diam}^+\right)^r$ for some integer $r$. Then the claim follows from Lemma~\ref{lemma:first-properties-structure-sheaves}~\ref{lemma:first-properties-structure-sheaves-3}. 
\end{proof}

\begin{lemma}\label{lemma:almost-coherent-integral} Let $\X=\Spf A_0$ be an affine admissible formal $\O_K$-scheme with adic generic fiber $X=\Spa(A, A^+)$ of dimension $d$, and let $\cal{E}$ be an $\O_{X^\diam}^+$-vector bundle. Then 
\[
\bf{R}\Gamma\left(X^\diam_v, \cal{E}\right)^a \in \bf{D}^{[0, 2d]}_{acoh}(A_0).
\]
Moreover, 
\[
\bf{R}\Gamma\left(X^\diam_v, \cal{E}\right)^a \in \bf{D}^{[0, d]}_{acoh}(A_0)
\]
if $\cal{E}$ is small (see Definition~\ref{defn:small-integrally}). 
\end{lemma}
\begin{proof}
    Lemma~\ref{lemma:vector-bundle-complete} implies that $\cal{E}$ is derived $p$-adically complete. Thus, the result follows from Theorem~\ref{thm:main-thm-small}, \cite[\href{https://stacks.math.columbia.edu/tag/0A0G}{Tag 0A0G}]{stacks-project}, and Corollary~\ref{cor:check-almost-coh-mod-p-2}.
\end{proof}

\begin{lemma}\label{lemma:restrict-integral} Let $\X=\Spf A_0$ be an admissible affine formal $\O_K$-scheme with adic generic fiber $X=\Spa(A, A^+)$, and $\mf\colon \Spf B_0 \to \Spf A_0$ an \'etale morphism with adic generic fiber $f\colon Y \to X$, and $\cal{E}$ an $\O_{X^\diam}^+$-vector bundle on $X$. Then the natural morphism 
\[
r\colon \bf{R}\Gamma\left(X^\diam_{v}, \cal{E}\right) \otimes_{A_0} B_0 \to \bf{R}\Gamma\left(Y^\diam_v, \cal{E}\right).
\]
is an isomorphism.
\end{lemma} 
\begin{proof}
    The morphism $A_0 \to B_0$ is flat since $\mf$ is \'etale. Now Lemma~\ref{lemma:almost-coherent-integral} and Lemma~\ref{completion-finitely-generated} ensure that the cohomology groups of both $\bf{R}\Gamma\left(X^\diam_{v}, \cal{E}\right) \otimes_{A_0} B_0$ and $\bf{R}\Gamma\left(Y^\diam_v, \cal{E}\right)$ are (classically) $p$-adically complete. In particular, both complexes are derived $p$-adically complete. So it suffices to show that $r$ is an isomorphism after taking derived mod-$p$ fiber (see \cite[\href{https://stacks.math.columbia.edu/tag/0G1U}{Tag 0G1U}]{stacks-project}). Then the claim follows from Theorem~\ref{thm:main-thm-small-2}~\ref{thm:main-thm-small-2-3}~\ref{thm:main-thm-small-2-4}. 
\end{proof}

\begin{thm}\label{thm:main-thm-integral-2} Let $\X$ be an admissible formal $\O_K$-scheme with adic generic fiber $X$ of dimension $d$, and let $\cal{E}$ be an $\O^+_{X^{\diam}}$-vector bundle. 
Then 
\begin{enumerate}[label=\textbf{(\arabic*)}]
    \item\label{thm:main-thm-integral-2-1} $\bf{R}\nu_*\cal{E}\in \bf{D}^+_{qc, acoh}(\X)$ and $(\bf{R}\nu_*\cal{E})^a\in \bf{D}^{[0, 2d]}_{acoh}(\X)^a$;
    \item\label{thm:main-thm-integral-2-2} if $\X=\Spf A$ is affine, then the natural map 
    \[
    \rm{H}^i\left(X^\diam_v, \cal{E} \right)^{\Updelta} \to \rm{R}^i\nu_*\left(\cal{E}\right)
    \]
    is an isomorphism for every $i\geq 0$;
    \item\label{thm:main-thm-integral-2-3} the formation of $\rm{R}^i\nu_*(\cal{E})$ commutes with \'etale base change, i.e., for any \'etale morphism $\mf \colon \Y \to \X$ with adic generic fiber $f\colon Y\to X$, the natural morphism 
    \[
    \mf^* \left(\rm{R}^i\nu_{\X, *}(\cal{E}) \right)\to \rm{R}^i\nu_{\Y, *}\left(\cal{E}|_{Y^\diam}\right)
    \]
    is an isomorphism for any $i\geq 0$;
    \item\label{thm:main-thm-integral-2-4} if $\X$ has an open affine covering $\X=\bigcup_{i\in I} \sU_i$ such that $\cal{E}|_{(\sU_{i, K})^\diam}$ is small, then
    \[
    \left(\bf{R}\nu_{*}\cal{E}\right)^a \in \bf{D}^{[0, d]}_{acoh}(\X)^a;
    \]
    \item\label{thm:main-thm-integral-2-5} there is an admissible blow-up $\X'\to \X$ such that $\X'$ has an open affine covering $\X'=\bigcup_{i\in I} \sU_i$ such that $\cal{E}|_{(\sU_{i, K})^\diam}$ is small. 
    
    In particular, there is a cofinal family of admissible formal models $\{\X'_i
    \}_{i\in I}$ of $X$ such that 
    \[
    (\bf{R}\nu_{\X'_i, *}\cal{E})^a\in \bf{D}^{[0, d]}_{acoh}(\X'_{i})^a.
    \]
    for each $i\in I$. 
\end{enumerate}
\end{thm}
\begin{proof}
First, \ref{thm:main-thm-integral-2-5} follows directly from Lemma~\ref{lemma:reasonable-cover-2}. Therefore, we only need to prove \ref{thm:main-thm-integral-2-1}-\ref{thm:main-thm-integral-2-4}. These claims are local on $\X$, so we can assume that $\X=\Spf A$ is affine. Then it suffices to show that, for every \'etale morphism $\Spf B_0\to \Spf A_0$ with adic generic fiber $Y\to X$, 
\[
\rm{H}^i(Y_v^\diam, \cal{E}|_{Y^\diam})
\]
is almost coherent for $i\geq 0$, 
\[
\rm{H}^i(Y_v^\diam, \cal{E}|_{Y^\diam}) \simeq^a 0
\]
for $i>2d$ (resp. for $i>d$ if $\cal{E}$ is small), and the natural morphism
\[
\rm{H}^i(X_v^\diam, \cal{E}) \otimes_{A_0} B_0 \to \rm{H}^i(Y_v^\diam, \cal{E}|_{Y^\diam})
\]
is an isomorphism (see Lemma~\ref{useful} and its proof). The first two claims follow from Lemma~\ref{lemma:almost-coherent-integral}, while the last one follows from Lemma~\ref{lemma:restrict-integral} (and $A_0$-flatness of $B_0$). 
\end{proof}

Let us also mention a version of Theorem~\ref{thm:main-thm-integral} for the pro-\'etale site of $X$ as defined in \cite{Sch1} and \cite{Sch-err}. It will be convenient to have this reference in our future work. In what follows, $\wdh{\O}_X^+$ is the completed integral structure sheaf on $X_\proet$ (see \cite[Definition 4.1]{Sch1}), and 
\[
\nu'\colon (X_\proet, \wdh{\O}_X^+) \to (\X_{\rm{Zar}}, \O_{\X})
\]
is the evident morphism of ringed sites.

\begin{thm}\label{main-thm-proetale} Let $\X$ be an admissible formal $\O_K$-scheme with adic generic fiber $X$ of dimension $d$ and mod-$p$ fiber $\X_0$. 
Then 
\begin{enumerate}[label=\textbf{(\arabic*)}]
    \item $\bf{R}\nu'_*\left(\O_{X}^+/p\right) \in \bf{D}^+_{qc, acoh}(\X_0)$ and $\bf{R}\nu'_*\left(\O_{X}^+/p\right)^a \in \bf{D}^{[0, d]}_{acoh}(\X_0)^a$;
    \item if $\X=\Spf A$ is affine, then the natural map 
    \[
    \widetilde{\rm{H}^i\left(X_\proet, \O_X^+/p \right)} \to \rm{R}^i\nu'_*\left(\O_X^+/p\right)
    \]
    is an isomorphism for every $i\geq 0$;
    \item the formation of $\rm{R}^i\nu'_*\left(\O_X^+/p\right)$ commutes with \'etale base change, i.e., for any \'etale morphism $\mf \colon \Y \to \X$ with adic generic fiber $f\colon Y\to X$, the natural morphism 
    \[
    \mf^*_0 \left(\rm{R}^i\nu'_{\X, *}\left(\O_X^+/p\right) \right)\to \rm{R}^i\nu'_{\Y, *}\left(\O_{Y}^+/p\right)
    \]
    is an isomorphism for any $i\geq 0$;
\end{enumerate}
\end{thm}
\begin{proof}
    By \cite[Corollary 3.17]{Sch1}, $\bf{R}\nu'_*\left(\O_X^+/p\right) \simeq \bf{R}t_*\left(\O_{X_\et}^+/p\right)$. So the results follow formally from Theorem~\ref{main-thm-2}. 
\end{proof}

\begin{thm}\label{thm:integral-structure-sheaf} Let $\X$ be an admissible formal $\O_K$-scheme with adic generic fiber $X$ of dimension $d$. Then 
\begin{enumerate}[label=\textbf{(\arabic*)}]
    \item\label{integral-structure-sheaf-1} $\bf{R}\nu'_*\wdh{\O}_X^+\in \bf{D}^+_{qc, acoh}(\X)$ and $(\bf{R}\nu'_*\wdh{\O}_X^+)^a\in \bf{D}^{[0, d]}_{acoh}(\X)^a$;
    \item\label{integral-structure-sheaf-2} if $\X=\Spf A$ is affine, then the natural map 
    \[
    \rm{H}^i\left(X_\proet, \wdh{\O}_X^+ \right)^{\Updelta} \to \rm{R}^i\nu'_*\wdh{\O}_X^+
    \]
    is an isomorphism for every $i\geq 0$;
    \item\label{integral-structure-sheaf-3} the formation of $\rm{R}^i\nu'_*(\cal{E})$ commutes with \'etale base change, i.e., for any \'etale morphism $\mf \colon \Y \to \X$ with adic generic fiber $f\colon Y\to X$, the natural morphism 
    \[
    \mf^* \left(\rm{R}^i\nu'_{\X, *}\left(\wdh{\O}_X^+\right) \right)\to \rm{R}^i\nu'_{\Y, *}\left(\wdh{\O}_Y^+\right)
    \]
    is an isomorphism for any $i\geq 0$;
\end{enumerate}
\end{thm}
\begin{proof}
    The proof is identical to the proof of Theorem~\ref{thm:main-thm-integral-2} once one establishes that the sheaf $\wdh{\O}_X^+$ is $p$-adically derived complete. For this, see \cite[Remark 5.5]{BMS1}.
\end{proof}

\appendix

\section*{Appendix}

\section{Derived complete modules}\label{appendix:derived-complete}

The main goal of this section is to collect some standard results on derived complete modules that seem difficult to find in the literature. \smallskip

For the rest of the section, we fix a ring $R$ with an element $\varpi\in R$. 

\begin{defn1}\label{defn:derived-complete} A complex $M\in \bf{D}(R)$ is {\it $\varpi$-adically derived complete} (or just derived complete) if the natural morphism $M\to \bf{R}\lim_n [M/\varpi^n]$ is an isomorphism. 
\end{defn1}

\begin{rmk1}\label{rmk:our-defn-stacks-project} This definition coincides with \cite[\href{https://stacks.math.columbia.edu/tag/091S}{Tag 091S}]{stacks-project} due to \cite[\href{https://stacks.math.columbia.edu/tag/091Z}{Tag 091Z}]{stacks-project}.
\end{rmk1}

\begin{lemma1}\label{lemma:derived-complete-ampl-mod-p} Let $M\in \bf{D}(R)$ be a derived complete complex. Then
\begin{enumerate}[label=\textbf{(\arabic*)}]
    \item\label{lemma:derived-complete-ampl-mod-p-1} $M\in \bf{D}^{\geq d}(R)$ if $[M/\varpi] \in \bf{D}^{\geq d}(R/\varpi)$.
    \item\label{lemma:derived-complete-ampl-mod-p-2} $M\in \bf{D}^{\leq d}(R)$ if $[M/\varpi] \in \bf{D}^{\leq d}(R/\varpi)$;
\end{enumerate}
\end{lemma1}
\begin{proof}
    \ref{lemma:derived-complete-ampl-mod-p-1}: By shifting, we can assume that $d=0$. Now suppose that $[M/\varpi] \in \bf{D}^{\geq 0}(R/\varpi)$. Then we use an exact triangles
    \[
    [M/\varpi]\to [M/\varpi^n] \to [M/\varpi^{n-1}]
    \]
    to ensure that $[M/\varpi^n]\in \bf{D}^{\geq 0}(R/\varpi^n)$ for every $n\geq 0$. Now we use that $M$ is derived complete to see that the natural morphism
    \[
    M \to \bf{R}\lim_n [M/\varpi^n M]
    \]
    is an isomophism. By passing to cohomology groups (and using that $\lim$ has cohomological dimension $1$), we see that 
    \[
    0 \to \rm{R}^1\lim_n\rm{H}^{i-1}([M/\varpi^n]) \to \rm{H}^i(M) \to \lim_n\rm{H}^i([M/\varpi^n]) \to 0
    \]
    are exact for any integer $i$. This implies that $\rm{H}^i(M) = 0$ for $i\leq 0$, i.e. $M\in \bf{D}^{\geq 0}(R)$. \smallskip
    
    \ref{lemma:derived-complete-ampl-mod-p-2}: Similarly, we can assume that $d=0$. Then the same inductive argument shows that $[M/\varpi^{n}]\in \bf{D}^{\leq 0}(R/\varpi^n)$ and we have short exact sequences
    \[
    0 \to \rm{R}^1\lim_n\rm{H}^{i-1}([M/\varpi^n]) \to \rm{H}^i(M) \to \lim_n\rm{H}^i([M/\varpi^n]) \to 0.
    \]
    This implies that $M\in \bf{D}^{\leq 1}(R)$ and $\rm{H}^1(M)=\rm{R}^1\lim_n \rm{H}^0([M/\varpi^n])$. Now note that the exact triangle
    \[
    [M/\varpi] \to [M/\varpi^n] \to [M/\varpi^{n-1}]
    \]
    and the fact that $[M/\varpi] \in \bf{D}^{\leq 0}(R/\varpi)$ imply that $\rm{H}^0([M/\varpi^n]) \to \rm{H}^0([M/\varpi^{n-1}])$ is surjective, so $\rm{R}^1\lim_n \rm{H}^0([M/\varpi^n])=0$ by the Mittag-Leffler criterion.  
\end{proof}

\begin{lemma1}\label{lemma:derived-complete-after-!} Let $R$ be a ring with an ideal of almost mathematics $\m$ and an element $\varpi \in \m$. Let $M\in \bf{D}(R)$ be a derived $\varpi$-adically complete complex. Then $\widetilde{\m}\otimes M$ is also derived $\varpi$-adically complete complex. 
\end{lemma1}
\begin{proof}
    Consider the exact triangle
    \[
    \widetilde{\m} \otimes M \to M \to Q.
    \]
    Since $\widetilde{\m}\otimes M \to M$ is an almost isomorphism, we see that cohomology groups of $Q$ are almost zero. In particular, they are $\varpi$-torsion, so derived complete. Therefore, $Q$ is derived complete (for example, by \cite[\href{https://stacks.math.columbia.edu/tag/091P}{Tag 091P}]{stacks-project} and \cite[\href{https://stacks.math.columbia.edu/tag/091S}{Tag 091S}]{stacks-project}). Now derived completeness of $M$ and $Q$ implies derived completeness of $\widetilde{\m}\otimes M$. 
\end{proof}

\begin{lemma1}\label{lemma:derived-complete-almost-ampl-mod-p} Let $R$ be a ring with an ideal of almost mathematics $\m$ and an element $\varpi \in \m$. Let $M\in \bf{D}(R)$ be a $\varpi$-adically derived complete complex. Then
\begin{enumerate}[label=\textbf{(\arabic*)}]
    \item $M^a\in \bf{D}^{\geq d}(R)^a$ if $[M^a/\varpi] \in \bf{D}^{\geq a}(R/\varpi R)^a$.
    \item $M^a\in \bf{D}^{\leq d}(R)^a$ if $[M^a/\varpi] \in \bf{D}^{\leq a}(R/\varpi R)^a$.
\end{enumerate}
\end{lemma1}
\begin{proof}
    Lemma~\ref{lemma:derived-complete-after-!} guarantees that $\widetilde{\m}\otimes M$ is derived $\varpi$-adically complete. Therefore, the claim follows from Lemma~\ref{lemma:derived-complete-ampl-mod-p} applied to $\widetilde{\m}\otimes M$. 
\end{proof}

Now we fix an $R$-ringed site $(X, \O_X)$. 

\begin{defn1}\label{defn:derived-complete-sheaf} A complex $M\in \bf{D}(X; \O_X)$ is {\it $\varpi$-adically derived complete} (or just derived complete) if the natural morphism $M\to \bf{R}\lim_n [M/\varpi^n]$ is an isomorphism. 
\end{defn1}

\begin{rmk1} This definition coincides with \cite[\href{https://stacks.math.columbia.edu/tag/0999}{Tag 0999}]{stacks-project} by \cite[\href{https://stacks.math.columbia.edu/tag/0A0E}{Tag 0A0E}]{stacks-project}.
\end{rmk1}

\begin{lemma1}\label{lemma:derived-complete-global-sections} Let $\cal{B} \subset \rm{Ob}(X)$ be a basis in a ringed site $(X, \O_X)$, and let $M\in \bf{D}(X; \O_X)$. Then $M$ is $\varpi$-adically derived complete if and only if $\bf{R}\Gamma(U, M)$ is $\varpi$-adically derived complete for any $U\in \cal{B}$.
\end{lemma1}
\begin{proof}
Suppose that $M$ is $\varpi$-adically derived complete. Then $\bf{R}\Gamma(U, M)$ is derived $\varpi$-adically complete for any $U\in \rm{Ob}(X)$ by \cite[\href{https://stacks.math.columbia.edu/tag/0BLX}{Tag 0BLX}]{stacks-project}.  \smallskip

Now suppose that $\bf{R}\Gamma(U, M)$ is $\varpi$-adically derived complete for any $U\in \cal{B}$, and consider the derived $\varpi$-adic completion $M \to \wdh{M}$ with the associated distinguished triangle:
\[
M \to \wdh{M} \to Q.
\]
We wish to show that $Q\simeq 0$. In order to show it, it suffices to establish that $\bf{R}\Gamma(U, Q)\simeq 0$ for any $U \in \cal{B}$. Now we use \cite[\href{https://stacks.math.columbia.edu/tag/0BLX}{Tag 0BLX}]{stacks-project} to conclude that 
\[
\bf{R}\Gamma(U, \wdh{M})\simeq \wdh{\bf{R}\Gamma(U, M)},
\]
so we get the distinguished triangle
\[
\bf{R}\Gamma(U, M) \to \wdh{\bf{R}\Gamma(U, M)} \to \bf{R}\Gamma(U, M).
\]
Since $\bf{R}\Gamma(U, M)$ is derived $\varpi$-adically complete by the assumption, so we see that the morphism 
\[
\bf{R}\Gamma(U, M) \to \wdh{\bf{R}\Gamma(U, Q)}
\]
is an isomorphism. Therefore, we conclude that $\bf{R}\Gamma(U, Q)\simeq 0$. This finishes the proof. 
\end{proof}

\section{Perfectoid rings}\label{Section:perfectoid-valuation}
The main goal of this Appendix is to recall the main structural results about perfectoid rings.

\begin{defn1}\label{defn:perfectoid-field}\cite[Definition 3.6]{Sch2} A non-archimedean field $(K, |\ .\ |_K)$ is a {\it perfectoid field} if there is a pseudo-uniformizer $\varpi \in K$ such that $\varpi^p \ | \ p$ in $\O_K=\{x\in K \ | \ |x|\leq 1\}$ and the $p$-th power Frobenius map
\[
\Phi\colon \O_K/\varpi \O_K \to \O_K/\varpi^p\O_K
\]
is an isomorphism. 
\end{defn1}

\begin{defn1}\label{defn:valuation-perfectoid} A complete valuation ring $K^+$ is a {\it perfectoid valuation ring} if $K\coloneqq \rm{Frac}(K^+)$ is a perfectoid field with its valuation topology. \smallskip

A Huber pair $(K, K^+)$ is a {\it perfectoid field pair} if $K$ is a perfectoid field and $K^+$ is an open and bounded valuation subring. 
\end{defn1}

\begin{rmk1}\label{rmk:perfectoid-microbial} Any perfectoid valuation ring $K^+$ is automatically microbial (see \cite[L9, Proposition 9.1.3 and Definition 9.1.4]{Seminar}). Any rank-$1$ valuation ring $K^+ \subset K^{++} \subset K$ defines the same topology on $K$ by \cite[Ch. VI, \textsection 7.2, Prop. 3]{Bou}. Therefore, $K^{++}$ must be equal to $K^\circ$ the set of powerbounded elements. In particular, there is a unique rank-$1$ valuation ring between $K^+$ and $K$ that we denote by $\O_K$, and the associated rank-$1$ valuation on $K$ by $|\ .\ |_K\colon K \to \bf{R}_{\geq 0}$. 
\end{rmk1}

\begin{lemma1}\label{lemma:equivalent-perfectoid-field}\cite[Proposition 3.8]{Sch2} Let $K$ be a non-archimedean field. Then $K$ is a perfectoid field if and only if the following conditions hold: 
\begin{enumerate}[label=\textbf{(\arabic*)}]
    \item $K$ is not discretely valued,
    \item $|p|_K < 1$,
    \item the Frobenius morphism $\Phi\colon \O_K/p\O_K \to \O_K/p\O_K$ is surjective. 
\end{enumerate}

\end{lemma1}


We wish to show that the ideal $\m=K^{\circ\circ} \subset K^+$ defines an ideal of almost mathematics in $K^+$. For the future reference, it will be convenient to do in a more general set-up of perfectoid pairs.

\begin{defn1}\label{defn:perfectoid-pair}\cite[Definition 3.1]{Sch2}  A complete Tate-Huber pair $(R, R^+)$ is a {\it perfectoid pair} if $R$ is a uniform Tate ring containing a pseudo-uniformizer $\varpi_R \in R^\circ$ such that $\varpi_R^p \ | \ p$ in $R^\circ$ and the Frobenius homomorphism $R^{\circ}/\varpi_R R^\circ \xr{x\mapsto x^p} R^{\circ}/\varpi_R^p R^\circ$ is an isomorphism. \smallskip

A Tate-Huber pair $(R, R^+)$ is {\it $p$-adic perfectoid pair} if it is a Huber pair, and $R$ is $p\neq 0$ in $R$. \smallskip

A Tate ring $R$ is a {\it perfectoid ring} if $(R, R^\circ)$ is a perfectoid pair. 
\end{defn1}

\begin{rmk1} It is not, a priori, clear that a perfectoid ring $R$ that is a field is a perfectoid field (in the sense of Definition~\ref{defn:perfectoid-field}). The problem is to verify that $R$ has a non-archimedean topology on it. This turned out to be always true by \cite{Kedlaya-Banach}.
\end{rmk1}

\begin{rmk1}\label{rmk:perfect-perfectoid} By \cite[Proposition 3.5]{Sch2}, a complete Tate ring $R$ of characteristic $p$ is perfectoid if and only if $R$ is perfect as a ring, i.e. the Frobenius morphism is an isomorphism. 
\end{rmk1}

\begin{rmk1}\label{rmk:only-surjectivity} In the definition of a perfectoid pair, it suffices to require $R^{\circ}/\varpi_R R^\circ \xr{x\mapsto x^p} R^{\circ}/\varpi_R^p R^\circ$ to be surjective. This map actually turns out to be always injective. Moreover, this condition turns out to be equivalent to the surjectivity of the Frobenius map
\[
R^\circ/pR^\circ \to R^\circ/pR^\circ.
\]
In particular, it is independent of a choice of a pseudo-uniformizer $\varpi_R^p \ | \ p$, see \cite[Remark 3.2]{Sch2} for more detail. Therefore, if $R$ is an algebra over a perfectoid field $K$ with a pseudo-uniformizer $\varpi_K\in \O_K$, one can always take $\varpi_R = \varpi_K$. In particular, every perfectoid ring in the sense of \cite[Definition 5.1]{Sch0} is a perfectoid ring in the sense of Definition~\ref{defn:perfectoid-pair}.
\end{rmk1}

\begin{lemma1}\cite[Lemma 3.10]{Sch2}\label{lemma:roots-of-pseudounformizer} Let $(R, R^+)$ be a perfectoid pair. Then there is a pseudo-uniformizer $\varpi\in R^{\circ\circ}$ such that 
\begin{enumerate}[label=\textbf{(\arabic*)}]
    \item\label{lemma:roots-of-pseudounformizer-1} $\varpi^p \ |\ p$ in $R^\circ$;
    \item\label{lemma:roots-of-pseudounformizer-2} $\varpi$ admits a compatible sequence of $p^n$-th roots $\varpi^{1/p^n}\in R^+$ for $n\geq 0$.
\end{enumerate}
In this case, $R^{\circ \circ} = \bigcup_{n\geq 0} \varpi^{1/p^n}R^+$.
\end{lemma1}
\begin{proof}
    \cite[Lemma 3.10]{Sch2} says that there is a pseudo-uniformizer $\varpi \in R^{\circ\circ}\subset R^+$ such that $\varpi^p \ | \ p$ in $R^\circ$, and there is a compatible sequence of the $p^n$-th roots $\varpi^{1/p^n}\in R^\circ$ for $n\geq 0$. Since $R^+$ is integrally closed, we conclude that all $\varpi^{1/p^n}$ must lie in $R^+$. Since $R^{\circ\circ}$ is a radical ideal $R^+$ and contains $\varpi$, it clearly contains $\bigcup_{n\geq 0} \varpi^{1/p^n}R^+$. \smallskip
    
    Now we pick an element $x\in R^{\circ\circ}$, and wish to show that $x\in \bigcup_{n\geq 0} \varpi^{1/p^n}R^+$. Since $x$ is topologically nilpotent, we can find an integer $m$ such that \[
    x^{p^m} \in \varpi R^+
    \]
    Therefore, $x^{p^m}=\varpi a$ for $a\in R^+$. Thus
    \[
    \left(\frac{x}{\varpi^{1/p^m}}\right)^{p^m} = a\in R^+.
    \]
    Therefore, $\frac{x}{\varpi^{1/p^m}} \in R^+$ because $R^+$ is integrally closed in $R$. So $x\in \varpi^{1/p^m}R^+$.
\end{proof}

\begin{rmk1}\label{rmk:roots-of-pseudounformizer} If $(R, R^+)$ is a $p$-adic perfectoid pair, then one can choose $\varpi$ such that $\varpi^p R^+ = p R^+$. Indeed, \cite[Lemma 3.20]{BMS1} implies that $R^+$ is perfectoid in the sense of \cite[Definition 3.5]{BMS1}. Thus the desired $\varpi$ exists by \cite[Lemma 3.9]{BMS1}.
\end{rmk1}

\begin{defn1}\label{defn:good-unifor} A pseudo-uniformozer $\varpi \in R^+$ of a $p$-adic perfectoid pair $(R, R^+)$ is {\it good} if $\varpi R^+ = pR^+$ and $\varpi$ admits a compatible sequence of $p$-power roots.
\end{defn1}

For the rest of the section, we fix a perfectoid pair $(R, R^+)$ and an ideal $\m=R^{\circ \circ}$. Our goal is to show that $\m$ defines a set-up for almost mathematics, i.e. $\widetilde{\m}=\m\otimes_{R^+} \m$ is $R^+$-flat and $\m^2=\m$.

\begin{lemma1}\label{lemma:almost-setup} Let $(R, R^+)$ be a perfectoid pair, and $\m=R^{\circ\circ}$ the associated ideal of topologically nilpotent elements. Then $\m$ is flat over $R^+$ and $\widetilde{\m} \simeq \m^2=\m$.
\end{lemma1}
\begin{proof}
    Lemma~\ref{lemma:roots-of-pseudounformizer} implies that $\m$ is flat as a colimit of free modules of rank-$1$. \smallskip
    
    Now we wish to show that $\m^2=\m$. We take any element $x\in \m$, by Lemma~\ref{lemma:roots-of-pseudounformizer} we know that $x=\varpi^{1/p^n}a$ for some integer $n$ and $a\in R^+$. Therefore, 
    \[
    x=\left(\varpi^{1/p^{n+1}}\right)^{p-1}\left(\varpi^{1/p^{n+1}}a\right) \in \m^2.
    \]
    
    Now we consider the short exact sequence
    \[
    0 \to \m \to R^+ \to R^+/\m \to 0.
    \]
    By flatness of $\m$, we conclude that it remains exact after applying the functor $-\otimes_{R^+} \m$. Therefore, the sequence
    \[
    0 \to \widetilde{\m}\to \m \to \m/\m^2 \to 0
    \]
    is exact. Since $\m^2=\m$, we conclude that 
    \[
    \widetilde{\m} \simeq \m^2=\m. \qedhere
    \]
\end{proof}

\begin{lemma1}\label{lemma:perfectoid-almost} Let $(R, R^+)$ be a perfectoid pair. Then the natural inclusion $\iota\colon R^+ \to R^{\circ}$ is an almost isomorphism.
\end{lemma1}
\begin{proof}
    Clearly, the map $\iota \colon R^+ \to R^\circ$ is injective, so it suffices to show that its cokernel is almost zero, i.e., annihilated by any $\e\in \m$. Pick an element $x\in R^\circ$, then $\e x\in R^{\circ\circ} \subset R^+$. Therefore, we conclude that $\e(\coker \iota)=0$ finishing the proof. 
\end{proof}

\section{Strongly sheafy adic spaces}\label{section:appendix-etale}

In this appendix, we discuss the notion of strongly sheafy spaces following \cite{Hansen-Kedlaya} and \cite{KedLiu1}.

\subsection{Preliminary results}

In this section, we discuss some results about general Tate-Huber pairs. 

\begin{lemma}\label{lemma:finite-etale-topology} Let $(A, A^+)$ be a complete Tate-Huber pair with a pair of definition $(A_0 \subset A^+, \varpi)$, and $A \to B$ is a finite \'etale morphism. Topologize $B$ using its natural $A$-module topology (see \cite[Appendix B.3]{ZavQuotient}). Then $(B, B^+)$ is a complete Tate-Huber pair where $B^+$ is the integral closure of $A^+$ in $B$. 
\end{lemma}
\begin{proof}
    {\it Step $1$: $B$ is complete in its natural topology}. Since $B$ is finite \'etale, $B$ is a projective $A$-module of finite rank. Then there is another finite $A$-module $M$ such that $B\oplus M \simeq A^{\oplus n}$. Consider the projection $p\colon A^{\oplus n} \to B$, the natural topology on $B$ coincides with the quotient topology (see \cite[Lemma B.3.2]{ZavQuotient}). Using the fact that $A$ is a Huber ring, it is not difficult to show that the quotient topology on $B$ should coincide with the subspace topology. Since $A^{\oplus n}$ is complete, we conclude that the natural topology on $B$ is separated. Therefore, the same applies to $M$ since we never used the ring structure on $B$. Then $B$ is closed in $A$ as a kernel of a continuous homomorphism with a separated target. In particular, $B$ is complete in its subspace (equivalently, quotient) topology, and as discussed above, this topology coincides with the natural topology. So it is complete in its natural topology. \smallskip
    
    {\it Step $2$: $B$ admits a finite set of $A$-module generators $x_1, \dots, x_n$ that are integral over $A_0$}. Pick any finite set $x'_1, \dots, x'_n \in B$ of $A$-module generators. It suffices to show that $x_i=\varpi^c x'_i \in B$ are integral over $A_0$ for some integer $c$. So it is enough to show that, for any $b\in B$, there is an integer $c$ such that $\varpi^c b$ is integral over $A_0$. \smallskip
    
    By definition, $b$ is integral over $A$. So we can find a monic equation
    \[
    b^n+a_{n-1}b^{n-1}+ \dots +a_0=0
    \]
    with $a_k\in A$ for $k=0, \dots, n-1$. Then there is an integer $c$ such that $\varpi^c a_k\in A_0$ for $k=0, \dots, n-1$. Thus, the equation
    \[
    (\varpi^cb)^n + a_{n-1}\varpi^c(\varpi^cb)^{n-1} +\dots + a_0\varpi^{cn}=0
    \]
    shows that $\varpi^c b$ is integral over $A_0$. \smallskip
    
    {\it Step~$3$: An $A_0$-subalgebra $B_0$ of $B$ generated by $x_1, \dots, x_n$ is finite as an $A_0$-module.} Clearly, this algebra is finitely generated over $A_0$ as an algebra and every element is integral. Therefore, it is finite. \smallskip
    
    {\it Step~$4$: $B_0$ is open in $B$ and the induced topology coincides with the $\varpi$-adic one.} Choose some $A_0$-module generators $b_1, \dots, b_m\in B_0$. Clearly, $B_0\left[\frac{1}{\varpi}\right]=B$, so the $A$-linear morphism
    \[
    q\colon \bigoplus_{i=1}^m Ae_i \to B
    \]
    sending $e_i$ to $b_i$ is surjective. By \cite[Lemma 2.4(i)]{H1}, $q$ is open. In particular, the topology on $B$ is the quotient topology along $q$. Therefore, $B_0$ is open in $B$ as
    $q^{-1}(B_0)$ is a subgroup that contains an open subgroup $\oplus_{i=1}^m A_0e_i$. Moreover, the topology on $B_0$ is $\varpi$-adic since $B_0=q(\oplus_{i=1}^m A_0e_i)$ with the quotient topology, and the topology on $\oplus_{i=1}^m A_0e_i$ is already $\varpi$-adic. \smallskip
    
    {\it Step~$5$. $(B, B^+)$ is a complete Huber pair:} We have already shown that $B$ is complete in its natural topology and $(B_0, \varpi)$ is a pair of definition for this topology. Therefore, $B$ is a Huber ring. It suffices to show that $B^+$ is open, integrally closed, and lies in $B^\circ$. Openness is clear since $B_0 \subset B^+$, and $B^+$ is integrally closed by definition. One also easily shows that $B^+ \subset B^\circ$ because $B^+$ is integral over $A^+ \subset A^\circ$. 
\end{proof}

\begin{lemma}\label{lemma:finite-etale-surjective} Let $(A, A^+)$ and $(B, B^+)$ be as in Lemma~\ref{lemma:finite-etale-topology}. Then $\Spec B \to \Spec A$ is surjective if and only if $\Spa(B, B^+) \to \Spa(A, A^+)$ is surjective. 
\end{lemma}
\begin{proof}
    First, we assume that $\Spec B \to \Spec A$ is surjective. To show that $\Spa(B, B^+) \to \Spa(A, A^+)$ is surjective, we need to show that $B \wdh{\otimes}_{A} \wdh{k(x)} \neq 0$ for any $x\in \Spa(A, A^+)$. Now \cite[Lemma B.3.5]{ZavQuotient} and Lemma~\ref{lemma:finite-etale-topology} ensure that $B \wdh{\otimes}_{A} \wdh{k(x)} = B\otimes_A \wdh{k(x)}$. To finish the proof, we note that $B\otimes_A k \neq 0$ for any field $k$ and a homomorphism $A \to k$ (in particular, this holds for $A \to \wdh{k(x)}$). \smallskip

    Now we assume that $\Spa(B, B^+) \to \Spa(A, A^+)$ is surjective. Then we note that \cite[Lemma 1.4]{H1} implies that every maximal ideal $\m\subset A$
    admits a valuation $v\in \Spa(A, A^+)$ such that $\rm{supp}(v) = \m$. This implies that the image of the morphism $\Spec B \to \Spec A$ contains all closed points of $\Spec A$. Since \'etale morphisms are open, we conclude that $\Spec B \to \Spec A$ must be surjective. 
\end{proof}

Now we discuss the notion of semi-uniform Tate-Huber pairs.

\begin{lemma}\label{lemma:uniform-+} Let $(A, A^+)$ be a (possibly noncomplete) Tate-Huber pair. Then $A^+$ is bounded if and only if $A$ is uniform (i.e., $A^\circ$ is bounded).
\end{lemma}
\begin{proof}   
    Clearly, $A^+$ is bounded if $A^\circ$ is bounded. So we assume that $A^+$ is bounded and we wish to show that $A^\circ$ is bounded as well. Choose a ring of definition $A_0\subset A^+$ and a pseudo-uniformizer $\varpi\in A_0$. Since $A^+$ is bounded, we conclude that there is an integer $N$ such that $A^+\subset \frac{1}{\varpi^N} A_0$. Now we note that $A^{\circ\circ} \subset A^+$ since $A^+$ is integrally closed and open. Since $\varpi$ is topologically nilpotent and any element $a\in A^\circ$ is powerbounded, we conclude that $\varpi A^{\circ} \subset A^{\circ \circ} \subset A^+$. Therefore, $A^{\circ} \subset \frac{1}{\varpi^{N+1}} A_0$, i.e., $A^\circ$ is bounded. 
\end{proof}

The above lemma motivates the following definition: 

\begin{defn} A (possibly noncomplete) Tate-Huber pair $(A, A^+)$ is {\it uniform} if $A^+ \subset A$ is bounded. 
\end{defn}

\begin{rmk} \cite[Proposition 1]{H0} implies that $(A, A^+)$ is uniform if and only if the subspace topology on $A^+$ coincides with the $\varpi$-adic topology for a (equivalently, any) choice of a pseudo-uniformizer $\varpi\in A^+$. Lemma~\ref{lemma:uniform-+} implies that it is equivalent to asking that the subspace topology on $A^\circ$ coincides with the $\varpi$-adic topology. 
\end{rmk}

\begin{lemma}\label{lemma:uniform-reduced} Let $A$ be a (possibly noncomplete) Tate ring. If $A$ is Hausdorff, then $A$ is reduced. In particular, any complete uniform Tate ring is reduced.
\end{lemma}
\begin{proof}
    Let $a\in A$ be a nilpotent element. We choose a pseudo-uniformizer $\varpi\in A^{\circ}$. Then $\frac{a}{\varpi^n}$ is nilpotent for any $n\geq 0$. In particular, it is bounded, so $\frac{a}{\varpi^n} \in A^{\circ}$. Thus, $a\in \varpi^n A^{\circ}$ for any $n\geq 0$. Since $A$ is uniform, we conclude that the topology on $A^{\circ}$ coincides with the $\varpi$-adic topology. Since $A$ is Hausdorff, we conclude that $\cap_{n\geq 0} \varpi^n A^{\circ}=0$. Thus, $a=0$ finishing the proof.
\end{proof}

\begin{defn}\label{defn:semi-uniform-completion} Let $(A, A^+)$ be a (possibly noncomplete) Tate-Huber pair with a pseudo-uniformizer $\varpi\in A^+$. The {\it uniformization} of $(A, A^+)$ is the Tate-Huber pair $(A_{u}, A_{u}^+)$, where $A_{u}^+=A^+$, $A_{u}=A$, and the topology on $A_{u}$ is induced from the $\varpi$-adic topology on $A_{u}^+$.

The {\it uniform completion} of $(A, A^+)$ is the Tate-Huber pair $(\wdh{A}_{u}, \wdh{A}^+_{u})$ obtained as the completion of $(A_{u}, A_{u}^+)$ (see \cite[Lemma 1.6]{H0}).
\end{defn}

\begin{rmk} We leave it to the reader to check that the uniformization is indeed a Tate-Huber pair and that it is independent of the choice of a pseudo-uniformizer $\varpi\in A^+$. In fact, uniformization is a left adjoint functor to the inclusion of uniform Tate-Huber pairs into the category of all Tate-Huber pairs. Likewise, uniform completion is a left adjoint functor to the inclusion of complete uniform Tate-Huber pairs into all Tate-Huber pairs. 
\end{rmk}

Now we discuss the relation between the topology of $\Spa(A, A^+)$ and its uniform completion.

\begin{lemma}\label{lemma:subdomains-completion} Let $(A, A^+)$ be a (possibly noncomplete) Tate-Huber pair. Then the natural morphisms
\[
\Spa(A_{u}, A_{u}^+) \to \Spa(A, A^+),~ \Spa(\wdh{A}, \wdh{A}^+) \to \Spa(A, A^+),\text{ and } \Spa(\wdh{A}_{u}, \wdh{A}^+_{u}) \to \Spa(A, A^+)
\]
are homeomorphisms that induce bijections on the sets of rational subdomains.
\end{lemma}
\begin{proof}
    First, \cite[Proposition 3.9]{H0} implies that the natural morphism $\Spa(\wdh{A}, \wdh{A}^+) \to \Spa(A, A^+)$ is a homeomorphism that induces a bijection on the sets of rational subdomains. Applying the same result to $(A_{u}, A_{u}^+)$, we see that it suffices to show the claim for $\Spa(A_{u}, A_{u}^+) \to \Spa(A, A^+)$. \smallskip
    
    For this, we note that rational subdomains on both sides are indexed by tuples $(f_1, \dots, f_n, g)$ of non-zero elements in $A$ generating the unit ideal, we conclude that it suffices to show that $\Spa(A, A^+) \to \Spa(A_{u}, A^+_{u})$ is a bijection. After unraveling the definition, we see that it suffices to show that any continuous (in the usual topology) valuation $v\colon A \to \Gamma_v\cup\{0\}$ is continuous in the topology induced from the $\varpi$-adic topology on $A^+$. Since $v$ is continuous, \cite[Corollary 9.3.3]{Seminar} implies that $v(\varpi)\in \Gamma_v$ is cofinal and $v(\varpi a)<1$ for any ring of definition $\varpi\in A_0$. Likewise, {\it loc.~cit.} implies that it suffices to show that $v(\varpi a)<1$ for any $a\in A^+$. This follows from \cite[Corollary 1.3]{H0}, which ensures that we can always find a ring of definition $A_0\subset A^+$ which contains both $a$ and $\varpi$. 
\end{proof}

\begin{lemma}\label{lemma:complete-+-henselian} Let $(A, A^+)$ be a complete Tate-Huber pair, and let $\varpi\in A^+$ be a pseudo-uniformizer. Then $A^+$ is $\varpi$-adically henselian. 
\end{lemma}
\begin{proof}
    \cite[Corollary 1.3]{H0} ensures that $A^+$ is a filtered colimit of its subrings of definition $A_0\subset A^+$. Therefore, the result follows from \cite[\href{https://stacks.math.columbia.edu/tag/0ALJ}{Tag 0ALJ}]{stacks-project} and \cite[\href{https://stacks.math.columbia.edu/tag/0FWT}{Tag 0FWT}]{stacks-project}. 
\end{proof}

\begin{lemma}\label{lemma:finite-etale-completion} Let $(A, A^+)$ be a (possibly noncomplete) Tate-Huber pair with a pseudo-uniformizer $\varpi\in A^+$. Suppose that $A^+$ is $\varpi$-adically henselian, then the natural functors 
\[
-\otimes_A \wdh{A} \colon A_{\rm{f\et}} \to \wdh{A}_{\rm{f\et}}
\]
\[
-\otimes_A \wdh{A}_{u} \colon A_{\rm{f\et}} \to \wdh{A}_{u, \rm{f\et}}
\]
are equivalences. Furthermore, the natural maps $\rm{Idem}(A) \to \rm{Idem}(\wdh{A})$ and $\rm{Idem}(A) \to \rm{Idem}(\wdh{A}_{u})$ are bijections. 
\end{lemma}
\begin{proof}
    The semi-uniform completions of $(A, A^+)$ and $(\wdh{A}, \wdh{A}^+)$ coincide, so Lemma~\ref{lemma:complete-+-henselian} implies that it suffices to prove both claims for $\wdh{A}_{u}$. \smallskip

    The claim about finite \'etale algebras follows immediately from \cite[Proposition 5.4.54]{GR} and the observation that $\varpi\in A^+$ is a regular element. Then claim about idempotents follows from \cite[\href{https://stacks.math.columbia.edu/tag/09XI}{Tag 09XI}]{stacks-project} and the observation that any idempotent $e\in A$ must lie in $A^+$ because it is integral over $\Z$. 
\end{proof}

\subsection{Noetherian approximation}

The main goal of this section is to show a version of noetherian approximation for complete Tate-Huber pairs. The main result of this section was originally shown in \cite[Proposition 2.6.2]{KedLiu1} in a slightly different language. \smallskip

To motivate the definition below, we want to mention one important subtlety of working with complete Tate-Huber pairs: this category does not admit filtered colimits. However, this issue can be remedied by considering Tate-Huber pairs together with the choice of a ring of definition. 

\begin{defn} A {\it Tate-Huber quadruple} $(A, A^+, A_0, \varpi)$ is a quadruple of a Tate-Huber $(A, A^+)$, a ring of definition $A_0$, and a pseudo-uniformizer $\varpi \in A_0$. A {\it morphism of Tate-Huber quadruples} 
\[
f\colon (A, A^+, A_0, \varpi) \to (B, B^+, B_0, \pi)
\]
is a continuous ring homomorphism $f\colon A \to B$ such that $f(A^+)\subset B^+$, $f(A_0) \subset B_0$, and $f(\varpi) = \pi$.
\end{defn}

For the next definition, we fix a filtered system $\{(A_i, A^+_i, A_{i, 0}, \varpi), f_{i,j}\}_{i\in I}$ of Tate-Huber quadruples\footnote{We slightly abuse notation and denote pseudo-unifomizer in $A_{i, 0}$ by the same letter $\varpi$.}.

\begin{defn} The {\it filtered colimit} of $\{(A_i, A^+_i, A_{i, 0}, \varpi), f_{i,j}\}_{i,j\in I}$ is the Tate-Huber quadruple
\[
(\colim_I A_i, \colim_I A_i^+, \colim_I A_{i, 0}, \varpi),
\]
where we topologize $A_i$ by requiring $\colim_I A_{i, 0} \subset \colim_I A_i$ to be a ring of definition with a pseudo-uniformizer $\varpi$.\footnote{We note that this implies that the subspace topology on $\colim_I A_{i, 0}$ is equal to the $\varpi$-adic topology. We warn the reader that the colimit topology on $A_{i, 0}$ is usually {\it different} from the $\varpi$-adic one.} 

The {\it completed filtered colimit} of $\{(A_i, A^+_i, A_{i, 0}, \varpi), f_{i,j}\}_{i\in I}$ is the Tate-Huber quadruple
\[
(A_\infty, A_\infty^+, A_{\infty, 0}, \varpi) \coloneqq (\wdh{\colim_I A_i}, \wdh{\colim_I A_i^+}, \wdh{\colim_I A_{i, 0}}, \varpi),
\]
where $\wdh{-}$ stands for the topological completion in the sense of \cite[Chapitre III, \textsection 3, n.4]{Topologie-1} (see also \cite[Chapitre III, \textsection 6, n.5 and Chapitre III, \textsection 7, n.2]{Topologie-1})
\end{defn}

\begin{rmk}\label{rmk:explicit-completion} \cite[Lemma 1.6]{H0} gives a very explicit description of completed filtered colimits. Namely, the ring $A_{\infty, 0}$ is equal to $(\colim_I A_{i, 0})^{\wedge}_{\varpi}$, the usual $\varpi$-adic completion of $\colim_{I} A_{i, 0}$. The ring $A_{\infty}^+ = (\colim_I A_{i}^+) \otimes_{\colim_I A_{i, 0}} A_{\infty, 0}$, and $A_\infty = A_{\infty}^+\big[\frac{1}{\varpi}\big] = A_{\infty, 0}\big[\frac{1}{\varpi}\big]$.
\end{rmk}

For the next definition, we fix a filtered system of $\{(A_i, A_i^+), f_{i, j}\}_{i,j\in I}$ of (not necessarily uniform) Tate-Huber pairs with a compatible choice of pseudo-uniformizers $\varpi \in A_i^+$. 

\begin{defn}\label{defn:uniform-colimit} The {\it uniform filtered colimit} of $\{(A_i, A^+_i), f_{i,j}\}_{i,j\in I}$ is the filtered colimit of Tate-Huber quadruple $\{A_{i, u}, A^+_{i, u}, A^+_{i, u}, \varpi\}_{i\in I}$ (see Definition~\ref{defn:semi-uniform-completion}).

The {\it completed uniform filtered colimit} of $\{(A_i, A^+_i), f_{i,j}\}_{i,j\in I}$ is the completed filtered colimit of Tate-Huber quadruple $\{A_{i, u}, A^+_{i, u}, A^+_{i, u}, \varpi\}_{i\in I}$.
\end{defn}

\begin{rmk}\label{rmk:explicit-uniform-completion} Remark~\ref{rmk:explicit-completion} implies that the completed uniform filtered colimit is explicitly given by the pair $\Big(\big(\wdh{\colim_I A^+_i}\big)\big[\frac{1}{\varpi}\big], \wdh{\colim_I A^+_i}\Big)$, where $\wdh{(-)}$ stands for the $\varpi$-adic completion.
\end{rmk}

Now we wish to prove a version of noetherian approximation for complete (uniform) Tate-Huber pairs. Before we do this, we need to invoke the following basic fact: 

\begin{lemma}\label{lemma:quotient-closed-ideal} Let $(A,A^+)$ be a complete Tate-Huber pair, and let $I\subset A$ be a closed ideal. Then $(A/I, (A/I)^+)$ is a complete Tate-Huber pair, where $(A/I)^+$ is the integral closure of $A^+/(I\cap A^+)$ in $A/I$.
\end{lemma}
\begin{proof}
    First, we choose a ring of definition $A_0$ and a pseudo-uniformizer $\varpi \in A_0$. Then $A/I$ is complete in the quotient topology due to \cite[Chapitre IX, \textsection 3, Proposition 1]{topologie} and \cite[Chapitre IX, \textsection 9, Proposition 4]{topologie}. Then \cite[Proposition 2.4(ii)]{H1} ensures that the natural morphism $\pi\colon A \to A/I$ is open. Therefore, $\pi(A_0) \subset A/I$ is an open subset such that the subset topology coincides with the $\varpi$-adic topology. Furthermore, (the image of) $\varpi$ is clearly a topologically nilpotent unit in $A/I$. Therefore, we conclude that $A/I$ is a complete Tate ring. Thus, we only need to show that $(A/I)^+$ is open, integrally closed, and lies in $(A/I)^{\circ}$. It is closed integrally by construction and is open because it contains $\pi(A_0)$. Finally, we note that by construction we have $(A/I)^+ \subset \pi(A^\circ)^+ \subset (A/I)^{\circ}$, where $\pi(A^\circ)^+$ is the integral closure of $\pi(A^\circ)$ in $A/I$.
\end{proof}

\begin{example} Let $A$ be a Tate ring and $e\in A$ is an idempotent element. Then the ideal $eA$ is closed in $A$ since it is equal to the kernel of the continuous multiplication by $(1-e)$ map.
\end{example}

Finally, we are ready to prove the main results of this section.                                                                            

\begin{lemma}[{cf., \cite[Proposition 2.6.2]{KedLiu1}}]\label{lemma:noetherian-approximation} Let $(A, A^+)$ be a complete Tate-Huber pair, let $A_0\in A^+$ be a ring of defintion, and let $\varpi\in A_0$ be a pseudo-uniformizer. Then there is a filtered system of Tate-Huber quadruples $\{(A_i, A_i^+, A_{i, 0}, \varpi)\}_{i\in I}$ such that 
\begin{enumerate}
    \item each $A_i$ is a strongly noetherian complete Tate algebra;
    \item the completed filtered colimit $(A_\infty, A_\infty^+, A_{\infty, 0}, \varpi)$ of $\{(A_i, A_i^+, A_{i, 0}, \varpi)\}_{i\in I}$ is isomorphic to $(A, A^+, A_0, \varpi)$. 
\end{enumerate}
\end{lemma}
\begin{proof}
    The choice of a pseudo-uniformizer $\varpi \in A^+$ defines a map $(\Z((t)), \Z\llbracket t\rrbracket) \to (A, A^+)$. Then we put $I$ to be the filtered poset of all finite subsets $S\subset A^+$. For each $S\in I$, we consider the unique $(\Z((t)), \Z\llbracket t\rrbracket)$-linear continuous morphism
    \[
    \alpha_S \colon \Big(\Z((t))\langle X_f\rangle_{f\in S}, \Z\llbracket t\rrbracket\langle X_f\rangle_{f\in S}\Big) \to \Big(A, A^+\Big)
    \]
    that sends $X_f$ to $f$. We put $I_S=\ker \alpha_S$ and $I_S^+ = I_S \cap \Z\llbracket t\rrbracket\langle X_f\rangle_{f\in S}$. The ideal $I_S$ is closed because it is the kernel of a continuous morphism. Therefore, Lemma~\ref{lemma:quotient-closed-ideal} gives us a complete Tate-Huber pair
    \[
    (A_S, A_S^+) \coloneqq \big(\Z((t))\langle X_f\rangle_{f\in S}/I_S, (\Z((t))\langle X_f\rangle_{f\in S}/I_S)^+\big)
    \]
    that admits an injective continuous morphism $\ov{\alpha}_S\colon (A_S, A_S^+) \to (A, A^+)$. We finally define 
    \[
    A_{S, 0} \coloneqq \big(\Z\llbracket t\rrbracket \langle X_f\rangle_{f\in S}/I_S^+\big) \cap A_0. 
    \]
    This subring is clearly open and bounded, so it is a ring of definition due to \cite[Proposition 1]{H0}. Finally, we put $\varpi\in A_{S, 0}$ to be the image of $t$. Therefore, we note that 
    \[
    \{(A_S, A_S^+, A_{S, 0}, \varpi)\}_{S\in I}
    \]
    with natural (injective) transition maps is a filtered system of Tate-Huber quadruples. Using the explicit description of completed filtered colimits from Remark~\ref{rmk:explicit-completion}, we note that the (uncompleted) filtered colimit of $\{(A_S, A_S^+, A_{S, 0}, \varpi)\}_{S\in I}$ coincides with the completed filtered colimit, and it is isomorphic to $(A, A^+, A_0, \varpi)$. To finish the proof, we only need to show that each $A_S$ is strongly noetherian. This follows from the fact that $\Z((t))$ admits a noetherian ring of definition and $A_S$ is topologically finite type over $\Z((t))$. 
\end{proof}

\begin{lemma}\label{lemma:uniform-noetherian-approximation} Let $(A, A^+)$ be a complete uniform Tate-Huber pair, and let $\varpi\in A_0$ be a pseudo-uniformizer. Then there is a filtered system of complete uniform Tate-Huber pairs $\{(A_i, A_i^+)\}_{i\in I}$ such that 
\begin{enumerate}
    \item each $A_i$ is strongly noetherian;
    \item the completed uniform filtered colimit $(A_\infty, A_\infty^+)$ of $\{(A_i, A_i^+)\}_{i\in I}$ is isomorphic to $(A, A^+)$. 
\end{enumerate}
\end{lemma}
\begin{proof}
    The proof is similar to that of Lemma~\ref{lemma:noetherian-approximation}. We define the index set $I$ as in the proof of Lemma~\ref{lemma:noetherian-approximation}. Likewise, for any $S\in I$, we define $A_S$, $A_S^+$, and $A_{S, 0}$ as in the proof of Lemma~\ref{lemma:noetherian-approximation} as well. Then we wish to show that the Huber-Tate pair $(A_S, A_S^+)$ is uniform. Once we know this fact, the rest of the argument is the same. \smallskip
    
    Now we show that $(A_S, A_S^+)$ is uniform. we note that $\Z\llbracket t\rrbracket$ is excellent due to \cite[\href{https://stacks.math.columbia.edu/tag/07QW}{Tag 07QW}]{stacks-project}. Therefore, \cite[Main Theorem 2]{excellent-completion} implies that $A_{S, 0}$ is also excellent. Then we recall that $A_{S}^+$ was defined as the integral closure of $A_{S, 0}$ inside $A$. Since $A$ is reduced due to Lemma~\ref{lemma:uniform-reduced}, \cite[\href{https://stacks.math.columbia.edu/tag/03GH}{Tag 03GH}]{stacks-project} and \cite[\href{https://stacks.math.columbia.edu/tag/07QV}{Tag 07QV}]{stacks-project} imply that $A_{S}^+$ is a finite $A_{S, 0}$-module. Therefore, there is an integer $n$ such that $A_{S}^+\subset \frac{1}{\varpi^n} A_{S, 0}$, i.e., $A_S^+$ is bounded. This finishes the proof. 
\end{proof}

\subsection{\'Etale maps}\label{section:etale-maps}

In this section, we discuss (strongly) \'etale maps of general complete Tate-Huber pairs. We also show that strongly \'etale morphisms satisfy approximation along completed (uniform) filtered colimits of complete Tate-Huber pairs. 

\begin{defn}\label{defn:strongly-etale} A morphism $(A, A^+) \to (B, B^+)$ of complete Tate-Huber pairs is a {\it rational subdomain} if there is a finite set of non-zero elements $f_1, \dots, f_n, g\in A$ which generates the unit ideal in $A$ and $(B, B^+)=\big(A\langle \frac{f_i}{g}\rangle, A\langle \frac{f_i}{g}\rangle^+\big)$ as an $(A, A^+)$-algebra. We denote by $(A, A^+)_{rsd}$ the poset\footnote{\cite[Proposition 1.3]{H1} implies that there is at most unique continuous $(A, A^+)$-linear morphism between two rational subdomains over $(A, A^+)$.} of rational subdomains of $(A, A^+)$ (it coincides with the poset of rational subdomains of $\Spa(A, A^+)$)

A morphism $(A, A^+) \to (B, B^+)$ of complete Tate-Huber pairs is {\it strongly finite \'etale} if $A\to B$ is finite \'etale and $B^+$ is the integral closure of $A^+$ in $B$. We denote by $(A, A^+)_{\rm{s}\et}$ the category of finite \'etale $(A, A^+)$-pairs and all $(A, A^+)$-linear morphisms between them. 

A morphism $(A, A^+) \to (B, B^+)$ of complete Tate-Huber pairs is {\it strongly \'etale} if it can be written as a finite composition of finite \'etale morphisms and rational subdomains. We denote by $(A, A^+)_{\rm{s}\et}$ the category of \'etale $(A, A^+)$-pairs and all $(A, A^+)$-linear morphisms between them. 
\end{defn}

\begin{rmk} In what follows, we will freely use the fact that the category of complete Tate-Huber pairs admits pushouts. Explicitly, the pushout $(B, B^+)\otimes_{A, A^+} (C, C^+)$ is given by
\[
\big( B\wdh{\otimes}_A C, (B\wdh{\otimes}_A C)^+\big),
\]
where $B\wdh{\otimes}_A C$ is the completed tensor product, and $(B\wdh{\otimes}_A C)^+$ is the integral closure of (the image of) $B^+\wdh{\otimes}_{A^+} C^+$ in $B\wdh{\otimes}_A C$.
\end{rmk}

\begin{rmk} By definition, strongly \'etale maps are closed under composition. Lemma~\ref{lemma:finite-etale-topology} implies that strongly finite \'etale maps are closed under pushouts of complete Tate-Huber pairs (since completion is not needed). Therefore, all strongly \'etale maps are also closed under pushouts in the category of complete Tate-Huber pairs.
\end{rmk}

\begin{rmk}\label{rmk:strongly-finite-etale-algebraic} Lemma~\ref{lemma:finite-etale-topology} ensures that, for any complete Tate-Huber pair $(A, A^+)$, there is an equivalence $(A, A^+)_{\sfet} \simeq A_{\fet}$.
\end{rmk}

Now we wish to show that the category of strongly \'etale $(A, A^+)$-pairs satisfies approximation with respect to completed filtered colimits. It will be convenient to first prove a version of this result for completed uniform filtered colimits. For this, we need a number of preliminary lemmas. 

\begin{defn}\label{defn:clopen-immersion} A morphism $f\colon (A, A^+) \to (B, B^+)$ of complete Tate-Huber pairs is a {\it clopen immersion} if $A\to B$ is a topological quotient morphism, $B^+$ is equal to the integral closure of $A^+$, and $\ker f$ is generated by an idempotent element $e_f\in A$. 
\end{defn}

\begin{rmk} If $f\colon (A, A^+) \to (B, B^+)$ is a clopen immersion, then  Lemma~\ref{lemma:quotient-closed-ideal} implies that 
\[
(B, B^+)= \Big(A/e_fA, \big(A^+/(e_fA\cap A^+)\big)^+\Big)
\]
as complete Tate-Huber pairs. 
\end{rmk}

\begin{rmk}\label{rmk:clopen-idempotent} We note that the idempotent $e_f\in A$ in Definition~\ref{defn:clopen-immersion} is unique if exists. In particular, two clopen immersions $f\colon (A, A^+) \to (B, B^+)$ and $g\colon (A, A^+) \to (B, B^+)$ coincide if and only if $e_f=e_g$. 
\end{rmk}

For the next definition, we fix a complete Tate-Huber pair $(A, A^+)$ and an $(A, A^+)$-linear morphism $f\colon (B, B^+) \to (C, C^+)$ of complete Tate-Huber pairs. 

\begin{defn}\label{defn:graph-diagonal} The {\it graph} of $f$ is the unique continuous $(A, A^+)$-linear morphism 
\[
    \Gamma_f \colon \big(B\wdh{\otimes}_A C, (B\wdh{\otimes}_A C)^+\big) \to (C, C^+)
\] which sends $b\otimes c$ to $f(b)c$. 

The diagonal of $f$ is the morphism $\Delta_f\colon \big(C\wdh{\otimes}_B C, (C\wdh{\otimes}_B C)^+\big) \to (C, C^+)$ that sends $c\otimes c'$ to $cc'$. 
\end{defn}

\begin{lemma}\label{lemma:clopen-diagonal} Let $f\colon (A, A^+) \to (B, B^+)$ and $g\colon (B, B^+) \to (C, C^+)$ be morphisms of complete Tate-Huber pairs such that $g$ and $h\coloneqq g\circ f$ are strongly \'etale. Then the morphisms $\Delta_g$ and $\Gamma_g$ are clopen immersions.
\end{lemma}
\begin{proof}
    First, we note that it suffices to prove the claim for $\Delta_g$ (for all $g$). This follows from the following pushout square (for simplicity, we suppress the $+$-rings in the diagram below)
    \[
    \begin{tikzcd}
        B\wdh{\otimes}_A B \arrow{r}{\Delta_f}\arrow{d}{\rm{id}\otimes g} & B\arrow{d}{g} \\
        B \wdh{\otimes}_A C \arrow{r}{\Gamma_g} & C
    \end{tikzcd}
    \]
    and the observation that clopen immersions are preserved by pushouts. Now we show that it suffices to show that $\Delta_g$ is a clopen immersion, it suffices to show that $\Delta_f$ and $\Delta_h$ are clopen immersions. For this, we consider the following diagram (for simplicity, we suppress the $+$-rings in the diagram below)
    \begin{equation}\label{eqn:composition-clopen}
    \begin{tikzcd}
    B\wdh{\otimes}_A B \arrow{r}{\Delta_f} \arrow{d}{g\otimes g} & B \arrow{d} & \\
    C\wdh{\otimes}_A C \arrow{r}{\alpha}\arrow[rr, bend right, "\Delta_{h}"] & C\wdh{\otimes}_B C \arrow{r}{\Delta_g} & C,
    \end{tikzcd}
    \end{equation}
    where the left square is a pushout square. Now, if $\Delta_f$ is a clopen immersion, then $\alpha$ is a clopen immersion as well. Now since $\alpha$ and $\Delta_{h}$ are clopen immersions, then $\Delta_g$ is a clopen immersion as well (with $e_{\Delta_g}=\alpha(e_{\Delta_{h}})$). Therefore, we reduce the question to showing that $f$ and $h$ are clopen immersions. In other words, we can assume that $g$ is a strongly \'etale morphism. \smallskip

    Now we use Diagram~(\ref{eqn:composition-clopen}) and the observation that clopen immersions are preserved by compositions to conclude that it suffices to prove the result separately for strongly finite \'etale morphisms and rational subdomains. If $g$ is a rational subdomain, then $\Delta_g$ is clearly an isomorphism. Therefore, it suffices to assume that $g$ is a strongly finite \'etale map. In this case, Lemma~\ref{lemma:finite-etale-topology} and \cite[Lemma B.3.5]{ZavQuotient} imply that $B\wdh{\otimes}_A B = B\otimes_A B$. Therefore, the result follows from the algebro-geometric claim that $\Spec B \to \Spec (B\otimes_A B)$ is a clopen immersion for a finite \'etale $A \to B$. 
\end{proof}

\begin{thm}\label{thm:approximation-etale} Let $\{(A_i, A_i^+)\}_{i\in I}$ be a filtered system of complete uniform Tate-Huber pairs, and let $(A_\infty, A_\infty^+)=(\wdh{\colim_I A^+_i}\big[\frac{1}{\varpi}\big], \wdh{\colim_I A_i^+})$ be its completed uniform filtered colimit. Then
\begin{enumerate}[label=\textbf{(\arabic*)}]
    \item\label{thm:approximation-etale-11} the natural map $|\Spa(A_\infty, A_\infty^+)| \to \lim_I |\Spa(A_i, A_i^+)|$ is a homeomorphism of spectral spaces;
    \item\label{thm:approximation-etale-2} the natural map $\colim_I (A_i, A_i^+)_{\rm{rsd}} \to (A_\infty, A_\infty^+)_{\rm{rsd}}$ is a bijection;
    \item\label{thm:approximation-etale-3} the natural functor $2\text{-}\colim_I (A_i, A_i^+)_\sfet \to (A_\infty, A_\infty^+)_{\sfet}$ is an equivalence;
    \item\label{thm:approximation-etale-4}  the natural functor $2\text{-}\colim_I (A_i, A_i^+)_\set \to (A_\infty, A_\infty^+)_{\set}$ is an equivalence.
\end{enumerate}
\end{thm}
\begin{proof}
    Let us denote by $\big(\ov{A}, \ov{A}^+\big)=\big(\colim_I A^+_i\big[\frac{1}{\varpi}\big], \colim_I A_i^+\big)$ the uncompleted uniform filtered colimit of $(A_i, A_i^+)$. Then we easily see that the natural morphism $|\Spa(\ov{A}, \ov{A}^+)| \to \lim_I |\Spa(A_i, A_i^+)|$ is a homeomorphism. Now \cite[Proposition 3.9]{H0} implies that the natural map $|\Spa(A_\infty, A_\infty^+)| \to |\Spa(\ov{A}, \ov{A}^+)|$ is a homeomorphism that induces bijection on rational subdomains. This already proves \ref{thm:approximation-etale-11}. To see \ref{thm:approximation-etale-2}, we use that the homeomorphism $|\Spa(\ov{A}, \ov{A}^+)| \simeq |\Spa(A_\infty, A_\infty^+)|$ induces a bijection on rational subdomains. Since every rational subdomain of $|\Spa(\ov{A}, \ov{A}^+)|$ is defined at a finite level, we conclude that the natural morphism
    \[
    \colim_I (A_i, A_i^+)_{\rm{rsd}} \to (A_\infty, A_\infty^+)_{\rm{rsd}}
    \]
    is surjective. The map is also injective due to \label{thm:approximation-etale-1}  and \cite[\href{https://stacks.math.columbia.edu/tag/0A30}{Tag 0A30}]{stacks-project}. This finishes the proof of \ref{thm:approximation-etale-2}. \smallskip

    \ref{thm:approximation-etale-3} follows from Lemma~\ref{lemma:finite-etale-completion}, Remark~\ref{rmk:strongly-finite-etale-algebraic}, and a standard (algebraic) approximation for finite \'etale algebras. \smallskip
    
    Now we show \ref{thm:approximation-etale-4}. First, we set up some notation. For any complete Tate-Huber $(A_i, A_i^+)$-pair $(B_i, B_i^+)$ and $i'>i\in I\sqcup \{\infty\}$, we put
    \[
    (B'_i, B'^+_i)\coloneqq \big(B_i\wdh{\otimes}_{A_i} A_{i'}, (B_i\wdh{\otimes}_{A_i} A_{i'})^+\big).
    \]
    {\it Observation.} For any compatible sequence of complete Tate-Huber $(A_i, A_i^+)$-pairs $(B_i, B_i^+)$, the uniform completion of $(B_\infty, B_\infty^+)$ and the compeleted uniform filtered colimit of $\{(B_i, B'_i)\}_{i\in I}$ are isomorphic as $(A, A^+)$-pairs. 

    In what follows, we will freely use this observation. Finally, we are ready to start the proof. \smallskip
    
    {\it Step~$0$. Essential surjectivity} Using {\it Observation}, Lemma~\ref{lemma:subdomains-completion}, and Lemma~\ref{lemma:finite-etale-completion}, we can inductively reduce the question to showing that any rational subdomain (resp. finite \'etale pair) over $(A_\infty, A^+_\infty)$ comes from a finite level. This follows directly from \ref{thm:approximation-etale-11} (resp. \ref{thm:approximation-etale-2}).  \smallskip

    {\it Step~$1$. Faithfulness} We fix two systems of compatible morphisms $f_i, g_i \colon (B_i, B_i^+) \to (C_i, C_i^+)$ in $(A_i, A_i^+)_{\set}$ and wish to show that, if their pushouts to $(A_\infty, A_\infty^+)$ coincide, then $f_i=g_i$ for some $i\gg 0$. For this, we set $f_\infty, g_\infty \colon (B_\infty, B_\infty^+) \to (C_\infty, C_\infty^+)$ to be the morphisms induced by $f_i$ and $g_i$ respectively. \smallskip
    
    The graphs $\Gamma_{f_i}$, $\Gamma_{g_i}$, $\Gamma_{f_\infty}$, and $\Gamma_{g_\infty}$ are clopen immersions due to Lemma~\ref{lemma:clopen-diagonal}. We also notice that $f_i=g_i$ (resp. $f_{\infty}=g_\infty$) if and only if $\Gamma_{f_i}=\Gamma_{g_i}$ (resp. $\Gamma_{f_\infty}=\Gamma_{g_\infty}$). Furthermore, Remark~\ref{rmk:clopen-idempotent} implies that $\Gamma_{f_i}=\Gamma_{g_i}$ (resp. $\Gamma_{f_\infty}=\Gamma_{g_\infty}$) if and only if $e_{\Gamma_{f_i}}=e_{\Gamma_{g_i}}$ (resp. $e_{\Gamma_{f_\infty}}=e_{\Gamma_{g_\infty}}$). Thus, we reduce the question to showing that if two idempotents $e, e'\in A_i$ become equal in $A_\infty$, they are already equal in $A_j$ for some $j>i$. This follows from Lemma~\ref{lemma:finite-etale-completion} and usual properties of filtered colimits. This finishes the proof of faithfulness. \smallskip

    {\it Step~$2$. Fullness.} We start with two compatible sequences $(B_i, B_i^+)$ and $(C_i, C_i^+)$ of elements in $(A_i, A_i^+)_\set$ and a continuous $(A_\infty, A_\infty^+)$-linear morphism $f_\infty\colon (B_\infty, B_\infty^+) \to (C_\infty, C_\infty^+)$, we wish to show that it is defined over $(A_i, A_i^+)$ for some $i\in I$. For this, we can freely replace $I$ with $I_{\geq i_0}$ for some $i_0$ to assume that $I$ has a minimal element $i_0\in I$.\smallskip

    We write the morphism $g_{i_0}\colon (A_{i_0}, A_{i_0}^+) \to (B_{i_0}, B_{i_0}^+)$ is a composition of $n$ morphisms, each of which is finite \'etale or a rational subdomain. We argue by induction on $n$. If $n=0$, then $(A_{i_0}, A_{i_0}^+) = (B_{i_0}, B_{i_0}^+)$ and then the result is obvious (the morphism $f_\infty$ must be the structure morphism $(A_\infty, A_\infty^+) \to (C_\infty, C_\infty^+)$, so it descends to any $i\in I$). \smallskip

    Now we do the induction step. We write $g_{i_0}$ as a composition 
    \[
    (A_{i_0}, A_{i_0}^+) \xr{h_{i_0}} (B'_{i_0}, B'^+_{i_0}) \xr{g'_{i_0}} (B_{i_0}, B_{i_0}),
    \]
    where $g'_{i_0}$ is either strongly finite \'etale or a rational subdomain, and $h_{i_0}$ is a composition of at most $n-1$ finite \'etale and rational subdomain morphisms. By induction, we know that there is $i\in I$ such that the morphism $f'_{\infty}\coloneqq f_{\infty} \circ g'_{\infty} \colon (B'_{\infty}, B'^+_{\infty}) \to (C_\infty, C^+_\infty)$ is defined over $i\in I$. So we can replace $i_0$ with $i$ to assume that there is a morphism
    \[
    f'_{i_0} \colon (B'_{i_0}, B'^+_{i_0}) \to (C_{i_0}, C^+_{i_0})
    \]
    such that its pushout to $(A_\infty, A_\infty^+)$ is equal to $f'_\infty$. \smallskip

    Consider the following diagram
    \begin{equation}\label{eqn:doing-induction}
    \begin{tikzcd}
    & (C_{i}, C_{i}^+) & \\
    (A_{i}, A_{i}) \arrow{ru} \arrow{r} & (B'_{i}, B'^+_{i}) \arrow{u}{f'_{i}} \arrow{r}{g'_{i}}& (B_{i}, B^+_{i}) \arrow[ul, swap, dotted, "f_{i}"].
    \end{tikzcd}
    \end{equation}
    for $i\geq i_0\in I$. The proof of faithfulness boils down to constructing a morphism $f_{i}$ such that Diagram~(\ref{eqn:doing-induction}) commutes and the pushout of $f_{i}$ to $(A_\infty, A_\infty^+)$ is equal to $f_\infty$. For this, we consider two cases.

    {\it Case~$1$. $g'_{i_0}$ is a rational subdomain.} In this case, \cite[Proposition 1.3]{H1} implies that, for each $i\in I\sqcup \infty$, there is at most one $f_{i}$ which makes Diagram~(\ref{eqn:doing-induction}) commute. Furthermore, it exists if and only if 
    \[
    |\Spa(f'_{i})| \colon |\Spa(C_{i}, C^+_{i})| \to |\Spa(B'_{i}, B'^+_{i})|
    \]
    factors through $|\Spa(B_{i}, B^+_{i})| \subset |\Spa(B'_{i}, B'^+_{i})|$. {\it Loc.~cit.} implies that $|\Spa(f'_\infty)|$ factors through $|\Spa(B_\infty, B^+_{\infty})| \subset |\Spa(B'_\infty, B'^+_{\infty})|$. Therefore, {\it Observation}, Lemma~\ref{lemma:subdomains-completion}, Part~\ref{thm:approximation-etale-11}, \cite[\href{https://stacks.math.columbia.edu/tag/0A2S}{Tag 0A2S}]{stacks-project}, and \cite[\href{https://stacks.math.columbia.edu/tag/0A2X}{Tag 0A2X}]{stacks-project} imply that there is an index $i\in I$ such that 
    \[
    |\Spa(f'_{i})| \colon |\Spa(C_{i}, C^+_{i})| \to |\Spa(B'_{i}, B'^+_{i})|
    \]
    factors through $|\Spa(B_{i}, B^+_{i})| \subset |\Spa(B'_{i}, B'^+_{i})|$. This defines a morphism $f_{i}\colon (B_i, B_i^+) \to (C_i, C_i^+)$ which makes Diagram~(\ref{eqn:doing-induction}). Furthermore, its pushout to $(A_\infty, A_\infty^+)$ equals to $f_\infty$ due to its uniqueness. \smallskip

    {\it Case~$2$. $g'_{i_0}$ is a finite \'etale morphism.} Then we consider a pushout diagram
    \[
    \begin{tikzcd}
        (C_i, C_i^+) \arrow{r}{\alpha_i} & (D_i, D_i^+)\coloneqq (B_i\wdh{\otimes}_{B'_i} C_i, (B_i\wdh{\otimes}_{B'_i} C_i)^+)  \\
        (B'_i, B'^+_i) \arrow{u}{f'_i} \arrow{r}{g'_i} & \arrow{u}{f''_i} (B_i, B_i^+).
    \end{tikzcd}
    \]
    and notice that morphisms $f_i$ that makes Diagram~(\ref{eqn:doing-induction}) commute are in bijection with morphisms $\beta_i\colon (D_i, D_i^+) \to (C_i, C_i^+)$ such that $\beta_i\circ \alpha_i = \rm{id}$, i.e., they are in bijection with sections of $\alpha_i$. Now we note that $\alpha_i$ are finite \'etale as pushouts of finite \'etale morphisms. Therefore, we note that the question boils down to showing that any section of a finite \'etale morphism $\alpha_\infty\colon (C_\infty, C_\infty^+) \to (D_\infty, D^+_\infty)$ comes from a finite level. This follows from {\it Observation}, Lemma~\ref{lemma:finite-etale-completion}, and Part~\ref{thm:approximation-etale-2} (applied to the filtered system $\{(C_{i}, C_{i}^+)\}_{i\in I}$). 
\end{proof}

\begin{cor}\label{cor:completion-strongly-etale} Let $(A, A^+)$ be a complete Tate-Huber pair with the uniform completion $(\wdh{A}_{u}, \wdh{A}^+_{u})$. Then the natural functor
\[
(A, A^+)_{\set} \to (\wdh{A}_{u}, \wdh{A}^+_{u})_{\set}
\]
is an equivalence. 
\end{cor}
\begin{proof}
    This follows directly from Theorem~\ref{thm:approximation-etale} applied to the constant filtered system $\{(A, A^+)\}$. 
\end{proof}

\begin{cor}\label{cor:approximation-etale-general} Let $\{(A_i, A_i^+, A_{i, 0}, \varpi)\}_{i\in I}$ be a filtered system of complete Tate-Huber quadruples, and let $(A_\infty, A_\infty^+, A_{\infty, 0}, \varpi) \coloneqq (\wdh{\colim_I A_i}, \wdh{\colim_I A_i^+}, \wdh{\colim_I A_{i, 0}}, \varpi)$ be its completed filtered colimit. Then
\begin{enumerate}[label=\textbf{(\arabic*)}]
    \item\label{cor:approximation-etale-general-11} the natural map $|\Spa(A_\infty, A_\infty^+)| \to \lim_I |\Spa(A_i, A_i^+)|$ is a homeomorphism of spectral spaces;
    \item\label{cor:approximation-etale-general-2} the natural map $\colim_I (A_i, A_i^+)_{\rm{rsd}} \to (A_\infty, A_\infty^+)_{\rm{rsd}}$ is a bijection;
    \item\label{cor:approximation-etale-general-3} the natural functor $2\text{-}\colim_I (A_i, A_i^+)_\sfet \to (A_\infty, A_\infty^+)_{\sfet}$ is an equivalence;
    \item\label{cor:approximation-etale-general-4}  the natural functor $2\text{-}\colim_I (A_i, A_i^+)_\set \to (A_\infty, A_\infty^+)_{\set}$ is an equivalence.
\end{enumerate}
\end{cor}
\begin{proof}
    This follows directly from Theorem~\ref{thm:approximation-etale} and Corollary~\ref{cor:completion-strongly-etale} by replacing $A_i$, $A_i^+$, $A_\infty$, and $A_\infty^+$ by their uniform completions. 
\end{proof}

\subsection{Strongly sheafy adic spaces}\label{section:strongly-sheafy-spaces}

In this section, we define the notion of a strongly sheafy adic space. We also define the \'etale structure sheaves on such spaces. \smallskip

\begin{defn}\label{defn:strongly-sheafy} \cite[Definition 4.1]{Hansen-Kedlaya} A complete Tate ring $A$ is {\it strongly sheafy} if $A\langle T_1, \dots, T_d\rangle$ is sheafy for any integer $d\geq 0$. 

A Tate-affinoid (pre-)adic space $X=\Spa(A, A^+)$ is {\it strongly sheafy} if $A$ is strongly sheafy. 

An adic space $X$ is {\it strongly sheafy} if there is an open covering of $X$ by strongly sheafy Tate-affinoids. 
\end{defn}

\begin{example} A strongly noetherian Tate ring $A$ is strongly sheafy (see \cite[Theorem 2.2]{H1}). Likewise, a sousperfectoid Tate ring $A$ is strongly sheafy as well (see \cite[Definition 7.1 and Corollary 7.4]{Hansen-Kedlaya}).   
\end{example}

\begin{rmk} \cite[Proposition 5.5]{Hansen-Kedlaya} and (the proof of) \cite[Theorem 5.6]{Hansen-Kedlaya} imply that, if $(A, A^+)$ is a sheafy complete Tate-Huber pair and $\Spa(A, A^+)$ is a strongly sheafy adic space, then $A$ is a strongly sheafy Tate ring.
\end{rmk}

\begin{rmk}\label{rmk:etale-sheafy} \cite[Theorem 5.6 and Definition 5.4]{Hansen-Kedlaya} imply that, if $(A, A^+)$ is a strongly sheafy Tate-Huber ring and $(A, A^+) \to (B, B^+)$ is a strongly \'etale morphism, then $(B, B^+)$ is strongly sheafy as well.
\end{rmk}

The above remark allows us to make the following definition: 

\begin{defn}\label{defn:strongly-etale-spaces} A morphism of strongly sheafy Tate-affinoids $\Spa(B, B^+) \to \Spa(A, A^+)$ is {\it strongly finite \'etale} if $(A, A^+) \to (B, B^+)$ is strongly finite \'etale (in the sense of Definition~\ref{defn:strongly-etale}). 

A morphism of strongly sheafy Tate-affinoids $\Spa(B, B^+) \to \Spa(A, A^+)$ is {\it (affinoid) strongly \'etale} if $(A, A^+) \to (B, B^+)$ is strongly \'etale (in the sense of Definition~\ref{defn:strongly-etale}).
\end{defn}

\begin{rmk} We note that finite disjoint unions of rational subdomains $\sqcup_{i\in I} X_i \to X=\Spa(A, A^+)$ are strongly \'etale as they can be decomposed as a composition $\sqcup_{i\in I} X_i \to \sqcup_{i\in I}  X \to X$, where the first morphism is a rational subdomain and the second morpism is strongly finite \'etale. 
\end{rmk}

\begin{rmk}\label{rmk:union-strongly-etale} More generally, if $\{X_i \to X\}_{i\in I}$ is a finite family of strongly \'etale morphisms, then $\sqcup_{i\in I} X_i \to X$ is a strongly \'etale morphism as well.
\end{rmk}


\begin{lemma}\label{lemma:etale-sheafy} Let $Y$ be a strongly sheafy adic space, let $X$ be a pre-adic space (in the sense of \cite[Definition 8.2.3]{KedLiu1}), and let $f\colon X \to Y$ be an \'etale morphism (in the sense of \cite[Definition 8.2.16]{KedLiu1}). Then $X$ is a strongly sheafy adic space.
\end{lemma}
\begin{proof}
    The claim is local in the analytic topology on both $X$ and $Y$. Therefore, we can assume that $X=\widetilde{\Spa(B, B^+)}$, $Y=\widetilde{\Spa(A, A^+)}$ for a strongly sheafy Tate ring $A$, and the morphism $(A, A^+) \to (B, B^+)$ is strongly \'etale. Then $(B, B^+)$ is strongly sheafy due to Remark~\ref{rmk:etale-sheafy}. 
\end{proof}

Finally, we can define the \'etale integral and structure (pre-)sheaves on strongly sheafy spaces: 

\begin{defn}\label{defn:etale-structure-pre-sheaves} Let $X$ be a strongly sheafy adic space. The {\it \'etale structure pre-sheaf} $\O_{X_\et}$ is a pre-sheaf of rings on $X_\et$ defined via the assignment
\[
(Y \xr{\text{\'etale}} X) \mapsto \O_Y(Y)
\]
with evident transition maps. The {\it integral \'etale structure pre-sheaf} $\O^+_{X_\et}$ is a pre-sheaf of rings on $X_\et$ defined via the assignment
\[
(Y \xr{\text{\'etale}} X) \mapsto \O_Y^+(Y)
\]
with evident transition maps.
\end{defn}

Before we show that $\O_{X_\et}$ and $\O_{X_\et}^+$ are sheaves, we need to prove the following basic lemma:

\begin{lemma}\label{lemma:check-+-ring-after-surjection} Let $\{\varphi_i \colon (A, A^+) \to (B_i, B_i^+)\}_{i\in I}$ be a family of morphisms of complete Tate-Huber pairs such that $\cup_{i\in I} |\Spa(\varphi_i)|\big(|\Spa(B_i, B^+_i)|\big) = |\Spa(A, A^+)|$, and let $a\in A$. Then $a\in A^+$ if and only if $\varphi_i(a)\in B_i^+$.
\end{lemma}
\begin{proof}
    If $a\in A^+$, then clearly $\varphi_i(a)\in B_i^+$ for every $i\in I$. Now we assume that $\varphi_i(a)\in B_i^+$ for all $i\in I$ and wish to show that $a\in A^+$. First, \cite[Lemma 3.3(i)]{H0} (or \cite[Proposition 1.6(iv)]{H1}) implies that 
    \[
    A^+=\{f\in A\ | \ v(f)\leq 1~ \forall v\in \Spa(A, A^+)\}.
    \]
    Therefore, we wish to show that $v(a)\leq 1$ for any $v\in \Spa(A, A^+)$. For this, we choose $i\in I$ and a $w_i\in \Spa(B_i, B_i^+)$ such that $\Spa(\varphi_i)(w)=v$. Then we know that 
    \[
    v(a) = w(\varphi_i(a)) \leq 1.
    \]
    This finishes the proof. 
\end{proof}

\begin{lemma}\label{lemma:etale-structure-sheaves} Let $X$ be a strongly sheafy adic space. Then the \'etale pre-sheaves $\O_{X_\et}$ and $\O_{X_\et}^+$ are sheaves. 
\end{lemma}
\begin{proof}
    Let $\{Y_i \to Y\}_{i\in I}$ be a covering in $X_\et$. We wish to verify the sheaf axiom for $\O_{X_\et}$ and $\O_{X_\et}^+$ with respect to this covering. Lemma~\ref{lemma:etale-sheafy} implies that the usual (analytic) pre-sheaves $\O_{Y}^+$, $\O_{Y_i}^+$, $\O_Y$, and $\O_{Y_i}$ are sheaves (in the analytic topology). Therefore, we can verify the sheaf condition analytically locally on $Y$ and $Y_i$. Therefore, we can assume that all spaces involved are strongly sheafy Tate-affinoids and all morphisms are strongly \'etale. In this case, sheafinness of $\O_{X_\et}$ follows from the last sentence of \cite[Proposition 5.5]{Hansen-Kedlaya} and (the proof of) \cite[Theorem 5.6]{Hansen-Kedlaya}. Then sheafinness of $\O_{X_\et}^+$ follows from sheafinness of $\O_{X_\et}$ and Lemma~\ref{lemma:check-+-ring-after-surjection}. 
\end{proof}

\section{Achinger's result in the non-noetherian case}\label{achinger}

Recall that P.\,Achinger proved a remarkable result \cite[Proposition 6.6.1]{Ach1} that says that an affinoid rigid-analytic variety $X=\Spa(A, A^+)$ that admits an \'etale map to a closed unit disc $\mathbf{D}^n_{K}$ also admits a {\it finite} \'etale map to $\mathbf{D}^n_K$ provided that $K$ is the fraction field of a complete DVR $R$ with residue field of characteristic $p$. This result is an analytic analogue of a more classical result of Kedlaya (\cite{Ked-finite} and \cite[Proposition 5.2.1]{Ach1}) that an affine $k$-scheme $X=\Spec A$ that admits an \'etale map to an affine space $\mathbf{A}^n_k$ also admits a {\it finite} \'etale to $\mathbf{A}^n_k$ provided that $k$ has characteristic $p$. 

We generalize P.\,Achinger's result to the non-noetherian setting. The proof essentially follows the ideas of \cite{Ach1}, we only need to be slightly more careful at some places due to non-noetherian issues. We also show a version of this result for formal schemes.


\begin{lemma1}\label{noether-normalization} Let $k$ be a field of characteristic $p$, and let $A$ be a finite type $k$-algebra such that $\dim A \leq d$ for some integer $d$. Suppose that $x_1, \dots, x_d \in A$ some elements of $A$, and $m$ is any integer $m\geq 0$. Then there exist elements $y_1, \dots, y_d \in A$ such that the map $f\colon k[T_1, \dots, T_d] \to A$, defined as $f(T_i)=x_i+y_i^{p^m}$ is finite. 
\end{lemma1}
\begin{proof}
We extend the set $x_1, \dots, x_d$ to some set of generators $x_1, \dots, x_d, \dots, x_n$ of $A$ as a $k$-algebra. This defines a presentation $ A=k[T_1, \dots, T_d, \dots, T_n]/I $
for some ideal $I\subset k[T_1, \dots, T_r, \dots, T_n]$. We prove the claim by induction on $n-d$. 

The case of $n-d=0$ is trivial as then the map $f\colon k[T_1, \dots, T_d] \to A$, defined by $f(T_i)=x_i$, is surjective. Therefore, it is finite. 

Now we do the induction argument, so we suppose that $n-d\geq 1$. We consider the elements
\[
x'_i=x_i -x_n^{p^{im'}}, \ i=1, \dots, n-1
\]
for some integer $m'\geq m$. Now the assumption $n\geq d +1$ and Krull's principal ideal theorem imply that we can choose some non-zero element $g\in I$, thus we have an expression
\[
g(x'_1+x_n^{p^{m'}}, x'_2+x_n^{p^{2m'}}, \dots, x'_{n-1}+x_n^{p^{(n-1)m'}}, x_n)=0
\]
Now \cite[\textsection 1]{redbook} implies that there is some large $m'$ such that this expression is a polynomial in $x_n$ with coefficients in $k[x'_1, \dots, x'_{n-1}]$ and a non-zero leading term. We may and do assume that this leading term is $1$. So $x_n$ is integral over a subring of $R$ generated by $x'_1, \dots, x'_{n-1}$, we denote this ring by $R'$. Since $x_i=x'_i+x_n^{p^{im'}}$, we conclude that $R$ is integral over $R'$. Moreover, $R$ is finite over $R'$ because it is finite type over $k$. Now we note that \cite[Theorem 9.3]{M1} implies that $\dim R' \leq \dim R \leq d$, and $R'$ is generated by $x'_1, \dots, x'_{n-1}$ as a $k$-algebra. So we can use the induction hypothesis to find some elements
\[
y'_1, \dots, y'_{d}\in R'
\]
such that the morphism $f'\colon k[T_1, \dots, T_d] \to R'$, defined as $f'(T_i)=x'_i+(y'_i)^{p^m}$, is finite. Therefore, the composite morphism
\[
f\colon k[T_1, \dots, T_d] \to R
\]
is also finite. We now observe that
\[
f(T_i)=x'_i+(y'_i)^{p^m}=x_i+x_n^{p^{im'}}+(y'_i)^{p^m}=x_i+(x_n^{p^{im'-m}}+y'_i)^{p^m}
\]
Therefore, the set $(y_i\coloneqq x_n^{p^{im'-m}}+y'_i)_{i=1,\dots, d}$ does the job. 
\end{proof}

\begin{lemma1}\label{finiteness-mod-ideal} Let $\O$ be a complete valuation ring of rank-$1$ with maximal ideal $\m$ and residue field $k$. Suppose that $f\colon A \to B$ is a morphism of topologically finitely generated $\O_K$-algebras. Then $f$ is finite if and only if $f\otimes_{\O} k\colon A\otimes_{\O} k \to B\otimes_{\O} k$ is finite.
\end{lemma1}
\begin{proof}
The ``only if'' part is clear, so we only need to deal with the ``if'' part. We recall that \cite[Lemma (28.P), p. 212]{M} says that $A\to B$ is finite if and only if $A/\pi \to B/\pi$ is finite for some pseudo-uniformizer $\pi \in \O$. So we only need to show that finiteness of $A\otimes_{\O} k \to B\otimes_{\O} k$ implies that there is a pseudo-uniformizer $\pi\in \O$ such that $A/\pi \to B/\pi$ is finite. Then we note that the maximal ideal $\m$ is a filtered colimit of its finitely generated subideal $\{I_{j}\}_{j\in J}$. Moreover, the valuation property of the ring $\O$ implies that this colimit is actually direct and that $I_{j}=(\pi_j)$ is principal for any $j\in J$. We also observe that each $\pi_j$ is a pseudo-uniformizer since $\O$ is of rank-$1$. Thus we see that 
\[
A\otimes_\O k \to B\otimes_\O k = \text{colim}_{j\in J} \left(A/\pi_j \to B/\pi_j\right)
\]
and $A/\pi_j \to B/\pi_j$ is a finite type morphism by the assumption that both $A$ and $B$ are topologically finitely generated. Then \cite[\href{https://stacks.math.columbia.edu/tag/07RG}{Tag 07RG}]{stacks-project} implies that there is $j\in J$ such that $A/\pi_j \to B/\pi_j$ is finite. Therefore, $A\to B$ is finite as well.
\end{proof}

Before going to the proof of Theorem~\ref{etale-finite-etale-formal}, we need to show a result on the dimension theory of rigid-analytic varieties spaces that seem to be missing in the literature. It seems that there is no generally accepted definition of a dimension of adic spaces. We define the dimension as $\dim X= \sup_{x\in X} \dim \O_{X,x}$, this is consistent with the definition of dimension in \cite[Definition II.10.1.1]{FujKato}. We denote by $X^{\text{cl}} \subset X$ the set of all classical points of $X$.

\begin{lemma1}\label{dimension} Let $f\colon X=\Spa(B, B^+) \to Y=\Spa(A, A^+)$ be an \'etale morphism of rigid-analytic varieties over a complete rank-$1$ field $K$, then $\dim B\geq \dim A$. If $Y$ is equidimensional, i.e. $\dim \O_{Y, y}=\dim Y$ for any classical point $y\in Y^{\text{cl}}$, then we have an equality $\dim B=\dim A$. In particular, if $f\colon \Spa(A, A^+) \to \mathbf{D}^d_{K}$ is \'etale, then $\dim A =d$.
\end{lemma1}
\begin{proof}
We note that \cite[Proposition II.10.1.9 and Corollary II.10.1.10]{FujKato} imply that 
\[
\dim X=\dim B=\sup_{x\in X^{\text{cl}}}(\dim \O_{X,x}), \ \text{and} \ \dim Y=\dim A=\sup_{y\in Y^{\text{cl}}}(\dim \O_{Y,y}). 
\] 
Since $f$ is topologically finite type, it sends classical points to classical points. Therefore, \cite[Lemma 1.6.4, Corollary 1.7.4, and Proposition 1.7.9]{H3} imply that the map $\O_{Y, f(x)} \to \O_{X, x}$ is finite \'etale for any $x\in X^{\text{cl}}$. Thus, we see that
\[
\dim B = \sup_{x\in X^{\text{cl}}} (\dim \O_{X, x}) = \sup_{x\in X^{\text{cl}}} (\dim \O_{Y, f(x)}) \leq \dim Y
\]
It is also clear that this inequality becomes an equality, if $Y$ is equidimensional. 

Finally, we claim that $\mathbf{D}^d_{K}=\Spa(K\langle T_1, \dots, T_d\rangle, \O_K\langle T_1, \dots, T_d\rangle)=\Spa(A, A^+)$ is equidimensional. Pick any classical point $x\in (\mathbf{D}^d_{K})^{\text{cl}}$ and a corresponding maximal ideal $\m_x\in K\langle T_1, \dots, T_d\rangle$. Then we know that  $A_{\m_x}$ and $\O_{\mathbf{D}^d_{K}, x}$ are noetherian by \cite[Proposition 0.9.3.9, Theorem II.8.3.6]{FujKato}, and $\wdh{\O_{\mathbf{D}^d_{K}, x}} \simeq \wdh{A_{\m_x}}$ by \cite[Proposition II.8.3.1]{FujKato}. Therefore, we get
\[
\dim \O_{\mathbf{D}^d_{K}, x}=\dim \wdh{\O_{\mathbf{D}^d_{K}, x}}=\dim \wdh{A_{\m_x}} =\dim A_{\m_x}=d
\]
where the last equality comes from \cite[Proposition 0.9.3.9]{FujKato}.
\end{proof}

For the rest of the section we fix a complete rank-$1$ valuation ring $\O$ with the fraction field $K$ and the characteristic $p$ residue field $k$. We refer to \cite[\textsection 1.9]{H3} for the construction of the adic generic fiber of a topologically finitely generated formal $\O$-scheme. The only thing we mention here is that it sends an affine formal scheme $\Spf A$ to the affinoid adic space $\Spa(A\otimes_{\O}K, A^+)$, where $A^+$ is the integral closure of the image $\Imm(A \to A\otimes_{\O}K)$.

\begin{thm1}\label{etale-finite-etale-formal} In the notation as above, let $g\colon \Spf A \to \widehat{\mathbf{A}}^d_{\O}$ be a morphism of flat, topologically finitely generated formal $\O$-schemes such that the adic generic fiber $g_K\colon \Spa(A\otimes_{\O}K, A^+) \to \mathbf{D}^d_{K}$ is \'etale. Then there is a {\it finite} morphism $f\colon \Spf A \to \widehat{\mathbf{A}}^d_{\O}$ that is \'etale on adic generic fibers.
\end{thm1}
\begin{proof}
First of all, we note that Lemma~\ref{dimension} implies that $\dim A\otimes_{\O}K = d$. Now \cite[Theorem 9.2.10]{FujKato} says that there exists a finite injective morphism
\[
\varphi\colon \O\langle T_1, \dots, T_d \rangle \to A 
\]
with an $\O_K$-flat cokernel. This implies that $K\langle T_1, \dots, T_d \rangle \to A\otimes_{\O} K$ is finite and injective. Flatness of $\coker \varphi$ implies that the map
\[
k[T_1, \dots, T_d] \to A\otimes_\O k
\]
is also finite and injective, so $\dim A \otimes_{\O} k = d$. Now we finish the proof in two slightly different ways depending on $\operatorname{char} K$. \smallskip

{\it Case 1, $\operatorname{char} K=p$}:  We consider the morphism $g^\#\colon \O\langle T_1, \dots, T_d \rangle \to A$ induced by $g$. We define $x_i\coloneqq g^{\#}(T_i)$ for $i=1, \dots, d$. Since $\dim A\otimes_\O k=d$ we can apply Lemma~\ref{noether-normalization} for the residue classes $\ov{x_1}, \dots, \ov{x_d}$ and $m=1$ to get elements $\ov{y_1}, \dots, \ov{y_d} \in A\otimes_\O k$ such that the map
\[
\ov{f^\#}\colon k[T_1, \dots, T_d] \to A\otimes_\O k, \text{ defined as } \ov{f^\#}(T_i)=\ov{x_i}+\ov{y_i}^p \text{ for } i=1, \dots, d
\]
is finite. We lift $\ov{y_i}$ in an arbitrary way to elements $y_i\in A$, and define 
\[
f^\#\colon \O\langle T_1, \dots, T_d \rangle \to A
\]
as $f^\#(T_i)=x_i+y_i^p$ for any $i=1, \dots, d$. This map is finite by Lemma~\ref{finiteness-mod-ideal}. 

Now we note that $X\coloneqq \Spa(A\otimes_{\O}K, A^+)$ is smooth over $K$, so \cite[Proposition 2.6]{BLR3} says that \'etaleness of $f_K\colon X \to \mathbf{D}^d_{K}$ is equivalent to the bijectivity of the map
\[
f_K^*\Omega^1_{\mathbf{D}^d_{K}/K} \to \Omega^1_{X/K}
\]
This easily follows from \'etaleness of $g_K$ and the fact that $d(x_i+y_i^p)=d(x_i)$ in characteristic $p$.\\

{\it Case 2, $\operatorname{char} K=0$}: We denote $\Spf A$ by $\X$ and its adic generic fiber $\Spa(A\otimes_{\O}K, A^+)$ by $X$. Then we use \cite[Proposition 2.6]{BLR3} once again to see that the map
\[
g_K^*\Omega^1_{\mathbf{D}^d_{K}/K} \to \Omega^1_{X/K}
\]
is an isomorphism. Since $(\widehat{\Omega^1}_{\X/\O})_K \simeq \Omega^1_{X/K}$ and the same for $\widehat{\mathbf{A}}^d_{\O}$ and $\mathbf{D}^d_K$, we conclude that the fundamental short exact sequence (\cite[Proposition I.3.6.3, Proposition I.5.2.5 and Theorem I.5.2.6]{FujKato})
\[
g^*\widehat{\Omega}^1_{\widehat{\mathbf{A}}^d_{\O}/\O} \to \widehat{\Omega}^1_{\X/\O} \to \widehat{\Omega}^1_{\X/\widehat{\mathbf{A}}^d_{\O}} \to 0
\]
implies that $\left(\widehat{\Omega}^1_{\X/\widehat{\mathbf{A}}^d_{\O}}\right)_K=0$. Furthermore, we know that  
\[
\widehat{\Omega}^1_{\X/\widehat{\mathbf{A}}^d_{\O}} \cong\left(\widehat{\Omega}^1_{A/\O \langle T_1, \dots, T_d \rangle}\right)^{\Updelta}
\]
for a finite $A$-module $\widehat{\Omega}^1_{A/\O \langle T_1, \dots, T_d \rangle}$ (see \cite[Corollary I.5.1.11]{FujKato}). We denote this module by $\wdh{\Omega}^1_g$ for the rest of the proof, and recall that condition $\left(\wdh{\Omega}^1_{\X/\widehat{\mathbf{A}}^d_{\O}}\right)_K=0$ is equivalent to $\widehat{\Omega}^1_{g} \otimes_\O K =0$. Using finiteness of $\widehat{\Omega}^1_{g}$ and adhesiveness of $A$, we conclude that there is an integer $k$ such that
\[
p^k\widehat{\Omega}^1_{g}=0
\]
as $p$ is a pseudo-uniformizer in $\O$. Now, similarly to the case of $\text{char} K=p$, we consider the morphism
\[
g^\#\colon \O\langle T_1, \dots, T_d \rangle \to A
\]
and define $x_i\coloneqq g^\#(T_i)$ for $i=1, \dots, d$. Again, using that $\dim A\otimes_\O k=d$ we can apply Lemma~\ref{noether-normalization} for the residue classes $\ov{x_1}, \dots, \ov{x_d}$ and $m=k+1$ to get elements $\ov{y_1}, \dots, \ov{y_d} \in A\otimes_\O k$ such that the map
\[
\ov{f^\#}\colon k[T_1, \dots, T_d] \to A\otimes_\O k, \text{ defined as } \ov{f^\#}(T_i)=\ov{x_i}+\ov{y_i}^{p^{k+1}} \text{ for } i=1, \dots, d,
\]
is finite. We lift $\ov{y_i}$ to some elements $y_i\in A$ and define
\[
f^\#\colon \O\langle T_1, \dots, T_d\rangle \to A
\]
by $f^{\#}(T_i)=x_i+y_i^{p^{k+1}}$. The map $f^{\#}$ is finite by Lemma~\ref{finiteness-mod-ideal}. \smallskip

We are only left to show that the induced map
\[
f\colon X \to \wdh{\mathbf{A}}^d_{\O}
\]
is \'etale on adic generic fibers. We claim that $p^k(\wdh{\Omega}^1_f)=0$. Indeed, we use \cite[Proposition I.5.1.10]{FujKato} to trivialize $\wdh{\Omega}^1_{\O \langle T_1, \dots, T_d \rangle/\O}\simeq \oplus_{i=1}^d dT_i  \O \langle T_1, \dots, T_d \rangle$, so we have the fundamental exact sequence
\[
\bigoplus_{i=1}^d AdT_i \xr{dT_i \mapsto d(x_i+y_i^{p^{k+1}})} \wdh{\Omega}^1_{A/\O} \to \wdh{\Omega}^1_f \to 0
\]
As $d(y_i^{p^{k+1}})$ is divisible by $p^{k+1}$. Therefore, we see that modulo $p^{k+1}$ this sequence is equal to
\[
\bigoplus_{i=1}^d A/p^{k+1}dT_i \xr{dT_i \to d(x_i)} \wdh{\Omega}^1_{A/\O}/p^{k+1} \to \wdh{\Omega}^1_f/p^{k+1} \to 0
\]
Thus, we see that $\wdh{\Omega}^1_f/p^{k+1}\simeq \wdh{\Omega}^1_g/p^{k+1}$. In particular, 
\[
\left(p^k\wdh{\Omega}^1_f\right)/p\left(p^k\wdh{\Omega}^1_f\right)=\left(p^k\wdh{\Omega}^1_g\right)/p\left(p^k\wdh{\Omega}^1_g\right)=0
\]
by the choice of $k$. Therefore, $p^k\wdh{\Omega}^1_f=0$ by \cite[Lemma (28.P), p. 212]{M}. By passing to the adic generic fiber we get the map $f_K\colon X \to \mathbf{D}^d_K$ such that 
\[
d(f_K)\colon f_K^*\Omega^1_{\mathbf{D}^d_K/K} \to \Omega^1_{X/K}
\]
is surjective. However, we recall that $X$ and $\mathbf{D}^d_K$ are both smooth rigid-analytic varieties of (pure) dimension $d$. Thus $d_{f_K^*}$ is a surjective map of vector bundles of the same dimension $d$, so it must be an isomorphism. Finally, \cite[Proposition 2.6]{BLR3} implies that $f_K$ is \'etale. 
\end{proof} 

\begin{cor1}\label{etale-finite-etale-rigid} Let $K$ be a complete rank-$1$ valuation field with a valuation ring $\O_K$, and the residue field $k$ of characteristic $p$. Suppose that $g\colon X=\Spa(A, A^+) \to \mathbf{D}^d_K$ is an \'etale morphism of  affinoid rigid-analytic $K$-varieties. Then there exists a {\it finite} \'etale morphism $f\colon X \to \mathbf{D}^d_K$.
\end{cor1}
\begin{proof}
We note that \cite[Lemma 4.4]{H1} implies that $A^+=A^{\circ}$, so the map $g$ corresponds to the map 
\[
g^\#\colon (K\langle T_1, \dots, T_d \rangle, \O_K\langle T_1, \dots, T_d \rangle) \to (A, A^{\circ})
\]
of Tate-Huber pairs. Theorem~\ref{etale-finite-etale-formal} implies that it suffices to show that the image of $\O_K\langle T_1, \dots, T_d \rangle$ inside some ring of definition $A_0\subset A$. \smallskip

Since $A$ is topologically finitely generated, we can extend $g^\#$ to a surjection 
\[
\varphi\colon K\langle T_1, \dots, T_d, X_1, \dots, X_n \rangle \twoheadrightarrow A.
\]
Then  
\[
A_0\coloneqq \varphi(\O_K\langle T_1, \dots, T_d, X_1, \dots, X_n \rangle)
\] 
is bounded and it is open due to the Banach Open Mapping Theorem (\cite[Lemma 2.4 (i)]{H1}). Thus, it is a ring of definition containing $g^\sharp\left(\O_K\langle T_1, \dots, T_d\rangle\right)$.
\end{proof}


\bibliography{biblio}

\end{document}